\documentclass[11pt,reqno]{amsart}
\usepackage{hyperref}
\usepackage{amsmath}
\usepackage{amsthm}
\usepackage{amssymb}
\usepackage[dvipsnames]{xcolor}  
\usepackage{graphicx}
\usepackage{pdfsync}     
\usepackage{cleveref}
\usepackage{tikz}
\usepackage{mathrsfs,amscd}
\usepackage{braket}
\usepackage{psfrag} 
\usepackage{mathtools} 
\usepackage{esint} 
\usepackage{nccbbb} 
\usepackage[utf8]{inputenc} 
\usetikzlibrary{decorations.pathmorphing}
\usepackage{verbatim}
\usetikzlibrary{patterns}
\usepackage{subcaption}
\usepackage{pgfplots}
\pgfplotsset{compat=1.18} 

\usepgfplotslibrary{fillbetween}

\usetikzlibrary{decorations.pathreplacing, calligraphy}
\usetikzlibrary{arrows.meta}

\def\nc{\newcommand}
\nc{\red}[1]{{\color{red} #1}}
\nc{\green}[1]{{\color{green} #1}}
\nc{\cyan}[1]{{\color{cyan} #1}}
\definecolor{purple}{rgb}{0.9,0,0.8}
\nc{\purple}[1]{{\color{purple} #1}}
\definecolor{gray}{rgb}{0.5,0.5,0.5}
\nc{\gray}[1]{{\color{gray} #1}}
\nc{\Blue}[1]{{\color{blue} #1}}

\def\blue{}

\newtheorem{theorem}{Theorem}[section]
\newtheorem{lemma}[theorem]{Lemma}
\newtheorem{corollary}[theorem]{Corollary}
\newtheorem{proposition}[theorem]{Proposition}
\newtheorem{remark}[theorem]{Remark}
\newtheorem{definition}[theorem]{Definition}
\newtheorem{cns}[theorem]{Construction}

\newcommand{\eqnsection}{
   \renewcommand{\theequation}{\thesection.\arabic{equation}}
\makeatletter
 \csname @addtoreset\endcsname{equation}{section} 
\makeatother}

\def \be{\begin{equation}}
\def \ee{\end{equation}}
\def \bt{\begin{theorem}} 
\def \et{\end{theorem}}
\def \bl{\begin{lemma}} 
\def \el{\end{lemma}}
\def \bea{\begin{eqnarray}}
\def \eea{\end{eqnarray}}
\def \bas{\begin{eqnarray*}}
\def \eas{\end{eqnarray*}}

\renewcommand{\eqref}[1]{(\ref{#1})} 
\newcommand{\abbr}[1]{{\sc\lowercase{#1}}}

\def \al{\alpha}
\def \bb{\beta}
\def \ga{\gamma}

\def \de{\delta}
\def \De{\Delta}

\def \vep{\varepsilon}
\def \la{\lambda}

\def\btet{\boldsymbol{\theta}}
\def\bpsi{\boldsymbol{\psi}}
\def\bmu{\boldsymbol{\bar \mu}}
\def\boeta{\boldsymbol{\eta}}
\def\bsig{\boldsymbol{\sigma}}
\def \th{\theta}

\def \ze{\zeta}

\def \({\left(}
\def \){\right)}

\def \nn{\nonumber}
\def \bc{\begin{center} }
\def \ec{\end{center} }
\def \bs{\begin{slide} }
\def \es{\end{slide} }

\def \un{\underline}
\newcommand{\uu}{\un}
\def \wh{\widehat}
\newcommand{\wt}{\widetilde}

\def \ff{\infty}
\newcommand{\one}{\bf{1}}
\def\square{{\vcenter{\vbox{\hrule height.3pt
        \hbox{\vrule width.3pt height5pt \kern5pt
           \vrule width.3pt}
        \hrule height.3pt}}}}
\def\qed{{\hfill $\square$ \bigskip}}

\DeclareMathOperator*{\avsum}{\mathrlap{\hspace{.2ex}\rotatebox[origin=c]{20}{\bf---}}{\sum}}
\DeclareMathOperator*{\avvsum}{\mathrlap{\hspace{.2ex}\rotatebox[origin=c]{20}{\bf--}}{\sum}}
\def\deq{ \overset{\rm dist}{=}}
\def\dconv{\overset{\rm dist}{\Longrightarrow}}
\newcommand{\Pbm}{\mathbb{P}}
\newcommand{\Ebm}{\mathbb{E}}

\newcommand{\bA}{\mathbb{A}}
\newcommand{\bB}{\mathbb{B}}
\newcommand{\bC}{\mathbb{C}}

\newcommand{\bE}{\mathbb{E}}
\newcommand{\bF}{\mathbb{F}}
\newcommand{\bG}{\mathbb{G}}

\newcommand{\bI}{\mathbb{I}}

\newcommand{\bM}{\mathbb{M}}
\newcommand{\bN}{\mathbb{N}}

\newcommand{\bP}{\mathbb{P}}
\newcommand{\bQ}{\mathbb{Q}}
\newcommand{\bR}{\mathbb{R}}
\newcommand{\bS}{\mathbb{S}}
\newcommand{\Ipi}{[-\pi,\pi)}
\newcommand{\bT}{\mathbb{T}}

\newcommand{\bZ}{\mathbb{Z}}

\renewcommand{\AA}{\mathcal{A}}
\newcommand{\BB}{\mathcal{B}}
\newcommand{\CC}{\mathcal{C}}
\newcommand{\DD}{\mathcal{D}}
\newcommand{\EE}{\mathcal{E}}
\newcommand{\FF}{\mathcal{F}}
\newcommand{\GG}{\mathcal{G}}
\newcommand{\HH}{\mathcal{H}}
\newcommand{\II}{\mathcal{I}}
\newcommand{\JJ}{\mathcal{J}}

\newcommand{\LL}{\mathcal{L}}

\newcommand{\NN}{\mathcal{N}}
\newcommand{\OO}{\mathcal{O}}
\newcommand{\PP}{\mathcal{P}}
\newcommand{\QQ}{\mathcal{Q}}
\newcommand{\RR}{\mathcal{R}}

\newcommand{\TT}{\mathcal{T}}
\newcommand{\UU}{\mathcal{U}}

\newcommand{\WW}{\mathcal{W}}

\newcommand{\YY}{\mathcal{Y}}

\newcommand{\sfA}{\mathsf {A}}
\newcommand{\sfB}{\mathsf {B}}
\newcommand{\sfb}{\mathsf {b}}

\newcommand{\sfG}{\mathsf {G}}

\newcommand{\sfJ}{\mathsf {J}}
\newcommand{\sfK}{\mathsf {K}}
\newcommand{\sfL}{\mathsf {L}}
\newcommand{\sfM}{\mathsf {M}}
\newcommand{\sfMo}{\mathsf {M_\star}}

\newcommand{\sfT}{\mathsf {T}}

\newcommand{\sfX}{\mathsf {X}}

\newcommand{\sfa}{\mathsf{a}}
\newcommand{\sfc}{\mathsf{c}}
\newcommand{\sfd}{\mathsf{d}}
\newcommand{\sfe}{\mathsf{e}}
\newcommand{\sfk}{\mathsf{k}}
\newcommand{\sfl}{\mathsf{L}}
\newcommand{\sfm}{\mathsf{m}}
\newcommand{\sfo}{\mathsf{o}}
\newcommand{\sfs}{\mathsf{s}}
\newcommand{\sft}{\mathsf{t}}
\newcommand{\sfu}{\mathsf{u}}
\newcommand{\sfv}{\mathsf{v}}
\newcommand{\sfw}{\mathsf{w}}
\newcommand{\sfx}{\mathsf{x}}
\newcommand{\sfxc}{\sfx_{\sfc}}
\newcommand{\sfxb}{\sfx_{\sfb}}
\newcommand{\sfy}{\mathsf{y}}
\newcommand{\sfz}{\mathsf{z}}

\def \E{{\mathbf E}}

\newcommand{\barP}{\wh\bP}
\newcommand{\sfI}{\bar{\mathsf I}}
\newcommand{\sfIp}{\mathsf{I}'}
\newcommand{\tsfI}{\mathsf{I}^{\sfe}_{\sfKo}}
\newcommand{\sfII}{\mathsf{I}^{\sfT}}
\newcommand{\sfIw}{\mathsf{I}_{\sfwm}}
\newcommand{\sfIwe}{\mathsf{I}^{\sfe}_{\sfwm}}
\newcommand{\sfIs}{\mathsf{I}^{\star}}
\newcommand{\sfIos}{\mathsf{I}_\star}
\newcommand{\sfJo}{\sfJ^o}
\newcommand{\sfJm}{\partial \sfJ}
\newcommand{\sfJw}{\sfJ_\sfa}
\newcommand{\sfJwp}{\sfJ'_\sfa}
\newcommand{\sfJs}{\sfJ_\star}
\newcommand{\sfJus}{\sfJ^\star}
\newcommand{\sfJupl}{\sfJ}
\newcommand{\sfJup}{\sfJ_\sfe}
\newcommand{\sfJl}{\sfJ_\ell}

\newcommand{\osfm}{\overline{\sfm}}
\newcommand{\sfrho}{\boldsymbol{\rho}}

\newcommand{\sfco}{2}

\def \npp{{\bf 0}}
\newcommand{\bTs}{\partial \bC}

\newcommand{\BMS}{X}

\def\GJ{\sfG}
\def\GJo{\sfG^o}

\newcommand{\PRzu}{\bP_R^{\sfz,\sfu}}
\newcommand{\PRzv}{\bP_R^{\sfz,\sfv}}
\newcommand{\PRsu}{\bQ_{R_1}^{\sfs,\sfu}}
\newcommand{\PRsb}{\bQ_{R_1}^{\sfs,\sfv_\bb}}
\newcommand{\PRso}{\bQ_R^{m,\npp}}

\newcommand{\FFuz}{\FF^{\vep}_{u,z}}
\newcommand{\HHk}{\HH_{\sfk_\vep(z)}}
\newcommand{\HHu}{\HH^\star_u}
\newcommand{\HHs}{\HH^\star_\infty}

\newcommand{\srs}{{\sf r}_\star}
\newcommand{\ello}{\ell_o}
\newcommand{\CCT}{\CC^\star_{\vep,\bT^2}}
\newcommand{\CCdd}{{\CC}^{(\vep)}_{\bT_{\delta}}}
\newcommand{\CCd}{{\CC}^{(\vep)}_{\dlc,\srs}}
\newcommand{\CCdz}{{\CC}^{(\vep)}_{0,\srs}}
\newcommand{\CCz}{{\CC}^{(\vep)}_{\dlc,0}}
\newcommand{\CCex}{\CC^{(\vep),\sfe\sfx}_{\dlc,\srs}}
\newcommand{\CCp}{{\CC\PP}^{(\vep)}_{\srs}}
\newcommand{\CCpp}{{\CC\PP}^{(\vep),\tau_0}_{\dlc,\srs}}

\newcommand{\KQ}{\mathsf Q}

\def\sfKz{\sfK_{0}}
\def\sfKo{\sfK_{1}}
\def\sfKop{\sfK'_{1}}
\def\sfKt{\sfK_{2}}
\def\sfKth{\sfK_{3}}
\def\sfKtp{\sfK^\star_{2}}
\def\bKz{\bar b_{\sfKz}}
\def\sfae{\varpi_{\vep}}
\def\sfKf{\sfK_{4}}
\def\sfKfv{\sfK_{5}}
\def\tsfKfv{\widetilde{\sfK}_{5}}
\def\sfKp{\sfK_{6}}
\def\sfKw{\sfK_{\sfa}}
\def\nvf{n_\vep^{(4)}}
\def\nvfv{n_\vep^{(5)}}
\def\wnvfv{\tilde{n}_\vep^{(5)}}
\def\nvp{n_\vep^{(6)}}
\def\bst{\sfK'_{6}}
\def\sfLw{\sfL} 
\def\sfLo{\sfL_o}
\def\sfLp{\sfLo'}
\def\sfLwp{\sfLw'}
\def\sfLn{\sfL_1}
\def \dlc{\mathsf{\delta_\star}}
\def \cst{\kappa_\star}
\def \cstp{\wt{\kappa}_\star}
\def \ka{\kappa}
\def \gak{{\kappa_o}}
\def \gakm{{\kappa_o^-}}
\def \gakp{{\kappa_o^+}}
\def \kapo{{\kappa_1}}
\def \kapom{{\kappa_1^-}}

\def \kapp{{\kappa'}}
\def \kappp{{\kappa'_+}}
\def \kappm{{\kappa'_-}}
\def\okappa{{\kappa''}}
\def \okappam{{\kappa''_-}}
\def\kab{\kappa_{a,b}}

\def \boetak{\boeta^{(\vep)}_k}
\def \boetaj{\boeta^{(\vep)}_j}
\def \boetasj{\boeta^{\star,(\vep)}_j}
\def \boetasl{\boeta^{\star,(\vep)}_\ell}
\def \boetaKo{\boeta^{(\vep)}_{\sfKo}}
\def \boetasKo{\boeta^{\star,(\vep)}_{\sfKo}}
\def\boetasKop{\boeta^{\star,(\vep)}_{\sfKop}}
\def \Istar{\bI_\star^{(\vep)}} 

\newcommand{\xs}{x_\star}
\newcommand{\xz}{x_o}

\def \UUk{\UU^{(\vep)}_{\bb,k}} 
\def \UUz{\UU^{(\vep)}_{\bb,[z]}}
\def\UUv{\UU^{(\vep)}_{\sfv,[z]}}
\def\whUUp{\wh{\UU}^{(\vep),\langle z \rangle}_{\bb}}
\def\whUUpk{\wh{\UU}^{(\vep),\langle z \rangle}_{\sfk,\bb}}

\def\whUUvk{\UU^{(\vep),\bb,\langle z \rangle}_{\sfk,\sfv}}
\def\whUUvkp{\UU^{(\vep),\bb,\langle z' \rangle}_{\sfk,\sfv}}
\def \whUUz{\wh{\UU}^{(\vep)}_{\bb,[z]}}
\def \UUp{\UU^{(+,\vep)}_{\bb,[z]}}

\def\UUpb{\UU^{(+,\vep)}_{\bb',[z]}}

\def\UUy{\UU^{(\vep)}_{\bb}}
\def\UUmy{\UU^{(\star,\vep)}_{\bb,\sfKo}}
\def\UUmp{\UU^{(\star,\vep)}_{\bb,\sfKop}}
\def\UUpy{\UU^{(+,\vep)}_{\bb}}
\def\UUx{\UU^{(\star,\vep)}_{\bb,\sfx}}
\def\UUxp{\UU^{(\star,\vep)}_{\bb',\sfx'}}
\def\UUpx{\UU^{(\star,\vep)}_{\bb,\sfx}}
\def\UUxx{\UU^{(\star,\vep)}_{\bb,\sfx,\sfx'}}
\def\UUpl{\UU^{(+,\vep)}_{\bb,\ell}}
\def\UUpp{\UU^{(+,\vep)}_{\bb'}}
\def\UUplp{\UU^{(+,\vep)}_{\bb',\ell}}

\def\sfwm{\measuredangle}

\def \Yxz{Y^{(\vep)}_{\sfx,[z]}}
\def\Yxiz{Y^{(\vep)}_{\sfx_i,[z]}}
\def \Yvz{Y^{(\vep)}_{\sfv,[z]}}
\def\Yyz{Y^{(\vep)}_{\sfy,[z]}}
\def\Vbz{V^{[z]}_{\bb}}
\def\Vbp{V^{\langle z' \rangle}_{\bb}}
\def\Vb{h^{\sfwm}_{\vep}}
\def \Vxz{V^{(\vep)}_{\sfx,[z]}}
\def\Vxiz{V^{(\vep)}_{\sfx_i,[z]}}
\def\Vyz{V^{(\vep)}_{\sfy,[z]}}
\def\Vvz{V^{(\vep)}_{\sfv,[z]}}

\def\wX{\wh{\sfX}}

\def\Tvjz{T_{\sfv,j}^{[z]}}
\def\Tv1z{T_{\sfv,1}^{[z]}}
\def\TvKoz{T_{\sfv,\sfKo}^{[z]}}
\def\TxKoz{T_{\sfx,\sfKo}^{[z]}}
\def\TxKo2z{T_{\sfx,\sfKo-2}^{[z]}}
\def\Tjz{T_{\bb,j}^{[z]}}
\def\TKoz{T_{\bb,\sfKo}^{[z]}}
\def\Tvjy{T_{\sfv,j}^{\bb,(y)}}
\def\Tjy{T_{\bb,j}^{(y)}}

\newcommand{\trav}[7]{\ensuremath{T_{#5,#6\to #7}^{#2,#4\overset{#1}{\to}#3} }}
\newcommand{\travshort}[7]{\trav{#1}{#2}{#3}{#4}{#5}{#6}{#7}}

\def\etavjz{\eta_{\sfv,j}^{[z]}}

\def\etajz{\eta_{\bb,j}^{[z]}}
\def\etav1z{\eta_{\sfv,1}^{[z]}}
\def\whetavjz{\wh{\eta}_{\sfv,j}^{\,[z]}}
\def\whetajz{\wh{\eta}_{\bb,j}^{\,[z]}}

\def\etavjy{\eta_{\sfv,j}^{\bb,(y)}}
\def\etaxny{\eta_{\sfx,n_\vep+1}^{\bb,(y)}}
\def\etaxnyp{\eta_{\sfx',n_\vep+1}^{\bb,(y)}}
\def\etajy{\eta_{\bb,j}^{(y)}}
\def\whetajy{\wh{\eta}_{\bb,j}^{\,(y)}}
\def\whetavjy{\wh{\eta}_{\sfv,j}^{\bb,(y)}}
\def\whetaky{\wh{\eta}_{\bb,k}^{\,(y)}}
\def\why{\wh{\eta}^{\,(y)}}

\def\whetaynp{\wh{\eta}_{\bb,n'}^{\,(y)}}
\def\whetaynpp{\wh{\eta}_{\bb,n''}^{\,(y)}}
\def\etaxvjp{\eta_{\sfv,j}^{\sfx,\langle z \rangle}}
\def\etabvjp{\eta_{\sfv,j}^{\bb,\langle z \rangle}}
\def\etavjp{\eta_{\sfv,j}^{\langle z \rangle}}
\def\etajp{\eta_{\bb,j}^{\langle z \rangle}}
\def\whetavjp{\wh{\eta}_{\sfv,j}^{\,\langle z \rangle}}
\def\etap{\eta^{\,\langle z \rangle}}

\def\etappp{\eta^{\sfx', \langle z' \rangle}}
\def\etabp{\eta^{\,\sfx,\langle z'' \rangle}}
\def\etaps{\eta^{\langle u \rangle,s}}
\def\etapsm{\eta^{\langle u \rangle,s_-}}
\def\whetapsm{\wh{\eta}^{\langle u \rangle,s_-}}
\def\whetap{\wh{\eta}^{\,\langle z \rangle}}
\def\whetajp{\wh{\eta}_{\bb,j}^{\,\langle z \rangle}}
\def\whetanp{\wh{\eta}_{\bb,n'}^{\,\langle z \rangle}}

\def\whetajpp{\wh{\eta}_{\bb,j}^{\,\langle z' \rangle}}
\def\whetajppp{\wh{\eta}_{\bb,j}^{\,\langle z'' \rangle}}
\def\whetaKpp{\wh{\eta}_{\bb,\sfKo}^{\,\langle z' \rangle}}
\def\whetaKppp{\wh{\eta}_{\bb,\sfKo}^{\,\langle z'' \rangle}}
\def \Ybz{\wh{\eta}^{\,[z]}_{\bb,\sfKo}}
\def\YbKp{\wh{\eta}^{\,[z]}_{\bb,\sfKop}}
\def\Ybzp{\wh{\eta}^{\,\langle z' \rangle}_{\bb,\sfKop}}
\def \Ybpz{\wh{\eta}^{\,[z]}_{\bb',\sfKo}}
\def\Ybpp{\wh{\eta}^{\,[z]}_{\bb',\sfKop}}
\def\whes{\wh{\eta}_\star}

\def\Ljz{L_{\bb,j}^{[z]}}

\def\LKoz{L_{\bb,\sfKo}^{[z]}}

\def\Ljp{L_{\bb,j}^{\langle z \rangle}}
\def\Ljpp{L_{\bb,j}^{\langle z' \rangle}}
\def\Lkp{L_{\bb,k}^{\langle z \rangle}}
\def\LKop{L_{\bb,\sfKo}^{\langle z \rangle}}
\def\LKopp{L_{\bb,\sfKo}^{\langle z' \rangle}}
\def\Ljy{L_{\bb,j}^{(y)}}
\def\LKoy{L_{\bb,\sfKo}^{(y)}}
\def\LKop{L_{\bb,\sfKop}^{(y)}}

\def\thbj{\theta_{q,\bb}^{(j)}}
\def\thbl{\theta_{q,\bb}^{(\ell)}}

\def\thbKo{\theta_{q,\bb}^{(\sfKo)}}

\def\Akz{\AA_{\sfk}^{[z]}}
\def\Akbz{\AA_{\sfk,\bb}^{[z]}}
\def\Akzm{\AA_{\sfk}^{[z-2\log \sfk]}}
\def\AUvz{\AA\, \UU_{\sfk,\sfv}^{(\vep),\bb,[z]}}

\def \BBb{\BB^{(\vep)}_{\bb,[z]}}
\def \BBy{\BB^{(\vep)}_{\bb}(y)}
\def \BBby{\BB^{(\vep)}_{\bb}}
\def \WWy{\WW^{(\vep)}_{\bb}(y)}
\def \WWb{\WW^{(\vep)}_{\bb,[z]}}
\def\BW{\BB\WW}
\def\RU{\RR\UU}
\def\CCe{\CC^{(\vep)}}
\def\wCe{\wh{\CC}^{(\vep)}}
\def\QQe{\RU^{(\vep)}}
\def\QQs{\RU^{(\star)}_\sfv}
\def\wQs{\wh{\QQ}}
\def\whcI{\AA}
\def\wtcI{\bar \AA}
\def\RUb{\RU^{(\vep)}_{\bb,[z]}}
\def\BWz{\BW_{\bb}^{[z]}}
\def\BWs{\RR\BW_{\bb}^{[z]}}
\def\BWp{\BW_{\bb}^{\,\langle z \rangle}}
\def\BWpp{\BW_{\bb}^{\,\langle z' \rangle}}
\def\BWb{\BW^{(\vep)}_{[\sfKo,n']}}
\def\BWlb{\BW^{(\vep),(y)}_{\bb,[\sfKo,\sfKt]}}
\def\BWl{\BW^{(\vep),(y),\Cap}_{[\sfKo,\sfKt]}}
\def\Bla{\BB^{(\vep),(a)}_{[\ell,\sfKt]}} 
\def\BWlo{\BW^{(\vep),(y_1),\Cap}_{[\sfKo,\sfKt]}}
\def\Bs{\BB^{(\vep),\Cap}_{[\sfKth,n']}}
\def\Bl{\BB^{(\vep),\smallfrown}_{[\sfKo,\sfKt]}}
\def\Bll{\BB^{(\vep)}_\ell}
\def\Bsl{\BB^{(\vep)}_{\star,\ell}}
\def\BKo{\BB^{(\vep)}_{\sfKo}}
\def\wtBl{\BB^{(\vep)}_{\pm}}
\def\wtBKo{\BB^{(\vep)}_{+,\sfKo}}
\def\wtBm{\BB^{(\vep)}_{-}}
\def\wtBp{\BB^{(\vep)}_{+}}
\def\wtBlp{\BB^{(\vep)}_{+,\ell}}
\def\wtBz{\BB^{(\vep)}_{o}}
\def\BBl{\BB^{(\vep)}_{\bb,\ell}}
\def\BBpl{\BB^{(\vep)}_{\bb',\ell}}
\def\BBKth{\BB^{(\vep)}_{[\sfKth,n'']}}
\def\BBKthp{\BB^{(\vep)}_{[\sfKth',n'']}}
\def\whBBl{\wh{\BB}^{(\vep)}_{\bb,\ell}}
\def\whBBlp{\wh{\BB}^{(\vep)}_{\bb,\ell'}}
\def\Rkp{\RR_{\sfk}^{\langle z \rangle}}
\def\Rkpm{\RR_{\sfk}^{\langle z - 3 \log \sfk \rangle}}
\def\Rkz{\wt{\RR}_{\sfk}^{\langle z \rangle}}
\def\Rkpp{\wt{\RR}_{\sfk}^{\langle z' \rangle}}
\def\Rkvp{\RR^{\langle z \rangle}_{\sfk,\sfv}}
\def\Rkbp{\RR^{o,\langle z \rangle}_{\sfk,\bb}}

\def\Rkvpp{\RR^{\bb,\langle z \rangle}_{\sfk,\sfv}}
\def\Rkvppp{\RR^{\bb,\langle z' \rangle}_{\sfk,\sfv}}
\def\RRv{\RR^{(\vep)}_{\sfv,[z]}}
\def\RRy{\RR^{(\vep)}_{\sfv}(y)}
\def\RUy{\RR\UU^{(\vep)}_{\bb}(y)}
\def\RUyz{\RR\UU^{(\vep)}_{\bb}}
\def\RRyz{\RR^{(\vep)}_{\sfv}}
\def\gst{g_\star}
\def\gstb{g_o}
\def\gon{g_1}
\def\gtw{g_2}
\def\gth{g_3}
\def\ggi{g_i}
\def\gfo{g_{\sfKf}}
\def\hg{g}
\def\gKo{g_{\sfKo}^{(\vep)}}
\def\whgo{\wh{g}_{\sfKo}^{(\vep)}}
\def\whgl{\wh{g}_{\ell}^{(\vep)}}
\def\whgC{\wh{g}_{\sfKo,j}^{(\vep),\Cap}}
\def\whgw{\wh{g}_{\sfKo,j}^{(\vep),\sfwm}}
\def\whgk{\wh{g}_{\sfKo,k}^{(\vep),\sfwm}}
\def\wtgo{\wh{g}_{\pm,\sfKo}^{(\vep)}}

\def\wtgm{\wh{g}_{-,\sfKo}^{(\vep)}}

\def\wtgp{\wh{g}_{+,\sfKo}^{(\vep)}}
\def\wtglp{\wh{g}_{+,\ell}^{(\vep)}}
\def\gKth{g_{\sfKth}^{(\vep)}}

\def\gfl{\wh{g}_{\ell}^{(\vep)}}
\def\whgth{\wh{g}_{\sfKth}^{(\vep)}}

\def\htq{h^{(\vep)}_{q}(\theta_q)}
\def\bgB{\bar g^{(\vep)}
}
\def\Delj{\De^{(\vep)}_{\btet,\ell}}
\def\Deljj{\De^{(\vep)}_{\btet,j}}
\def\DelKth{\De^{(\vep)}_{\btet,\sfKth}}

\def\DelKo{\De^{(\vep)}_{\btet,\sfKo}}
\def\bDl{\bar \De^{(\vep)}_{\btet}
}
\def\Pii{P^{(i)}}
\def\Pon{P^{(1)}}
\def\Ptw{P^{(2)}}
\def\Ptwl{P^{(2,{\sf Late})}}
\def\Pm{P^{(\star)}}
\def\dst{\xi}
\def\wbG{\widehat{\bG}}

\eqnsection

\begin{document}
\title[Brownian cover time for 2D torus]{Limit law for Brownian cover time\\
of the two-dimensional torus}
\author[Amir Dembo\;\; Jay Rosen\;\; Ofer Zeitouni]
{Amir Dembo\;\;  Jay Rosen\;\; 
Ofer Zeitouni}

\date{Oct 18,  2026}

\thanks{Amir Dembo was partially supported by NSF grant DMS-2348142.}
\thanks{Jay Rosen was partially supported by the Simons Foundation.}
\thanks{Ofer Zeitouni was partially supported by the Israel Science Foundation grant \#615$\backslash$24}


\begin{abstract}
Let $\CC_{\vep}$ denote the time it takes a 
Brownian motion to come within distance $\vep$ of every point in 
the two-dimensional torus $\bT^2$ of unit area.  We prove that for an
explicit non-random centering $\osfm_\vep$, the variables 
$\{\sqrt{\CC_{\vep}} - \osfm_\vep \}$ 
converge in distribution as $\vep \to 0$, to a randomly shifted Gumbel law.
\end{abstract}
\maketitle

\section{Introduction}
Let $\bM$ be a smooth, compact, connected,
two-dimensional Riemannian manifold without boundary. For each $\sfx\in \bM$ let
$\CC_{\sfx,\vep,\bM}$ be  the amount of time needed for a  
 Brownian motion (started uniformly on $\bM$) to come within (Riemannian) distance
$\vep$ of $\sfx$. 
There has been much interest in the precise asymptotics, as $\vep \to 0$, of the
$\vep$-cover time 
$\widetilde{\CC}^\star_{\vep,\bM} := \sup_\sfx \{ \CC_{\sfx,\vep,\bM} \}$. 
In this context, \cite[Theorem 1.3]{DPRZ} shows that a.s.
\begin{equation}\label{eq:dprz}
\lim_{\vep \to 0} \frac{\sqrt{\CC^\star_{\vep,\bM}}}{|\log \vep^2|} = 1 \,, \qquad
\text{where} \qquad \CC^\star_{\vep,\bM} := \frac{2 \pi}{\sfA_{\bM}} \widetilde{\CC}^\star_{\vep,\bM} 
\end{equation}
and $\bA_{\bM}$ denote the Riemannian area of $\bM$. Setting hereafter
\begin{align}\label{def:mn}
\osfm_\vep & :=\sfrho_\vep |\log  \vep|  \,,  \quad
\sfrho_\vep := \sfco - \frac{\log |\log \vep|}{\sfco |\log \vep|} 
\end{align}
\cite[(1.2)]{BK} shows that for the two-dimensional torus 
$\bT^2$ with $\sfA_{\bT^2}=1$, 
\begin{equation}
  \label{eq-BK1.2}
\sqrt{\CCT} - \osfm_\vep = o(\log |\log \vep|) \quad \mbox{\textrm{with high probability (\abbr{whp})}}\,.
\end{equation}
\blue{We refer to $\CCT$ as the \emph{normalized} $\vep$-cover time of $\bT^2$.  Likewise, 
for various subsets of $\bT^2$, the normalized $\vep$-cover time of the subset is its $\vep$-cover time multiplied by 
$2\pi$. We also note in passing that}
an analogue of \eqref{eq-BK1.2} for the cover time of the discrete torus by simple random walk was 
obtained in \cite{Abe}.

It is widely expected, by analogy with the case of extremes
of the two-dimensional discrete Gaussian free field,
that the random variable 
$\sqrt{\smash[b]{\CC^\star_{\vep,\bM}}}$, properly centered, converges in distribution
(see Section \ref{sec-background} for more about this prediction).
An important step in this direction was taken in \cite{BRZ2}, where
it is shown that
for the two-dimensional unit sphere $\bS^2$ (where $\sfA_{\bS^2} = 4\pi$), 
the collection 
$\{ g_\vep:=\sqrt{ \smash[b]{\CC^\star_{\vep,\bS^2}}}- \osfm_\vep \}$ is tight, and it is 
conjectured that the sequence $g_\vep$ converges in distribution.

Our goal in this paper is to settle  this conjecture for the case of the
two-dimensional torus $\bT^2$.
Recall that a standard Gumbel random variable $W_\star$ has the distribution
\begin{equation}
  \label{eq-Gumble}
  \bP(W_\star\leq x)=e^{-e^{-x}}.
\end{equation}
\begin{theorem}
\label{theo-covertime}
For some finite random variable $X_\star>0$ and any fixed $z\in \bR$,
\begin{align}
\lim_{\vep \to 0} \bP \big( \sqrt{\CCT}  \le \osfm_\vep + z  \big)
=  \bE \(\exp\{-X_\star \mathrm{e}^{-\sfco z}\} \) :=  \bP(Y_\star \leq z) \,.
\label{ct2.25b}
\end{align}
where
$Y_\star\stackrel{d}{=} (W_\star + \log X_\star)/2$.
\end{theorem}
\noindent 
That is, the normalized cover times $\{  \sqrt{\smash[b]{\CCT}} - \osfm_\vep \}$ 
converge in distribution to $Y_\star$,  which is one half of
a standard Gumbel variable
shifted by $\log X_\star$.
(The random variable $X_\star$ is the limit in law of 
a "derivative martingale" type variable $e^{-4 \ello} \overline{X}$,  see \eqref{def:mn-alt},  \eqref{dfn:Xstar-K}).

\subsection{Background and discussion}
\label{sec-background}
The motivation for studying the  cover time of graphs and manifolds, besides their intrinsic interest, stems from applications in the sciences. Among such graphs, the $d$-dimensional torus is probably the simplest and most natural, and it was long recognized
that the two-dimensional case is the hardest, see the quote from the (unfinished and unpublished) book of Aldous and Fill appearing in \cite[Section 1.1]{DPRZ}. To our knowledge, it first appears in the lovely article \cite{Wilf}.  We refer to the introduction of \cite{DPRZ} for further 
pointers to the early (pre 2000) literature on the problem.

A line of research closely related to the current work deals with the
extrema of Gaussian logarithmically correlated fields, of which
a prime example is the Gaussian free field in dimension $d=2$, and its discrete counterpart, the \abbr{DGFF}. Recall that the \abbr{DGFF} is the centered Gaussian field
on a box $V_N\subset \bZ^2$ of side $N$, whose covariance is the Green function
of the discrete Laplacian in $V_N$ with Dirichlet boundary conditions. Indeed,
at about the same time that \cite{DPRZ} was written, analogous results
for the \abbr{DGFF} appeared in \cite{BDG}.  A series of works 
(\cite{BoDeZe, BrZe, DiZe})
improved the result for the \abbr{DGFF}, until finally \cite{BDZ} and \cite{BiLo,BiLo1} gave a complete description of the maximum of the \abbr{DGFF} and the landscape at
high values of the field; a theory for general Gaussian log-correlated fields has also been developed in \cite{DiRZ}. All these works build on the modified second moment method introduced by Bramson to study the maximum of branching Brownian motion \cite{BR78,Bramson1}.

The seminal work \cite{DLP},  using isomorphism theorems from \cite{EKMRS},
showed a direct link between the cover time problem and the maximum of the
Gaussian free field,  for general graphs.  Sharper versions \cite{Di, Lu, zhai} then established
that the ratio of the square-root of the cover time to
the maximum of the 
\abbr{DGFF} converges to a constant.  A consequence of the latter is that
the main results of \cite{BDG} and \cite{DPRZ} follow from each other.

Going beyond the leading term, and in particular proving
distributional limit theorems, proved challenging. 
The arguably simplest
setup is that of the square-root of the
cover time of the binary tree 
by simple random walk; in that setup, the maximum of the \abbr{GFF} is simply the maximum of a branching (Gaussian) random walk, and detailed information on the latter
is available since \cite{Aid}.  For the closely related model of Branching Brownian motions, this goes further back to \cite{Bramson1}.
It is shown in \cite{DingZeitouni-ASharpEstimateForCoverTimesOnBinaryTrees} that the logarithmic correction terms 
in these two problems differ; however, with a proper rescaling, the asymptotics
of the cover time were evaluated in \cite{CoLo, DRZ}, and the analogue 
of Theorem \ref{theo-covertime} proved.

The methods in the latter papers are quite different. Indeed, 
the approach of \cite{CoLo} uses  an isomorphism theorem to reduce
the study to that of a Gaussian free field, while \cite{DRZ} derives the 
limit law by directly using a second moment method and barrier estimates 
from \cite{BRZ}. 

Extending the isomorphism methods to obtain sharp results
on the cover time of other graphs seems challenging, mostly because the
isomorphism theorems lead naturally to random walks with an additional
renewal structure that is not present in the case of the simple random walk.
(See  \cite{AB1} and \cite{LS26} for some recent results along these lines.)\footnote{After this work was completed, we learnt of 
 \cite{arXiv:2610:01852}, which claims the cover time limit result for the discrete torus, based on the method of \cite{LS26}.} 
For this reason, in this work we expand on the second moment methods developed in \cite{BRZ2} and \cite{DRZ} (using barrier estimates from \cite{BRZ}). We next outline the main ideas.

\subsection{High level description of the proof}
We first recall the main steps in our proof of convergence of the cover
time for the binary tree.  Let  $\TT_n$ denote the binary tree of depth $n$, whose only vertex of degree $2$  is  attached to an extra vertex 
$\rho$, called the  \abbr{root}. The \abbr{cover time} $\CC_n$ of $\TT_n$ is
the number of steps taken by a (discrete time) simple random walk (\abbr{srw}) 
started at $\rho$, till visiting 
all vertices of $\TT_n$. The main result of \cite{DRZ} is that $\sqrt{\smash[b]{2^{-(n+1)} \CC_n}}$, centered around the
explicit constant $m_n=c_\ast n-c_\ast^{-1}\log n$ with
$c_\ast=\sqrt{2\log 2}$, converges in distribution to a randomly shifted 
Gumbel  law. 

The proof in \cite{DRZ} proceeds in several steps.   As in the main body of the current paper, 
$\sfKz$ and $\sfKfv$ are constants that eventually tend to infinity only after $n$ tends to infinity.
\newline
$\bullet$  {\bf Step I} One translates time to excursion counts: let  
$N_k(t)$ denote the total number of excursions made by the \abbr{srw} at a fixed level 
$k \ge 1$ by time $t$. Then,
one has that the variance of $t^{-1/2} 2^{-k} N_k(t)$ goes to zero 
as $t\to\infty$ followed by $k \to \infty$. Thus,  
time $\frac{1}{2}(m_n+z)^2$ can be replaced with the proxy \blue{of completing} 
$\blue{\sfk_{n}(z)}= \bE [N_{\sfKz}(\frac{1}{2}(m_n+z)^2)]$ \blue{excursions at level $\sfKz \gg 1$.}
\newline
$\bullet$ {\bf Step II} Let $\blue{u} \in D_{\sfKz}$ denote a vertex at level $\sfKz$.  Let 
$N(u)$ denote the number of excursions at $u$
  until a total of  $\sfk_{n}(z)$ excursions have accumulated at level $\sfKz$. 
  It is then not hard to check that $\{\sqrt{2N(u)},  u \in D_{\sfKz}\}$ 
  converges in distribution, as $n\to\infty$, to a Gaussian vector of mean $m_n+z$
  and variance of order $\sfKz$ (this step in \cite{DRZ}
  is borrowed from \cite{CoLo}).
\newline 
$\bullet$
{\bf Step III} Conditioning at a fixed level $\sfKz$ and using that the subtree 
  rooted at each $u \in D_{\sfKz}$ is just $\TT_{n-\sfKz}$,  one reduces the problem to proving that a finite,   positive limit
\begin{equation}
\label{righttail.1}
\alpha_{\ast}:=
\lim_{y \to \infty}\lim_{n\to \infty} 
y^{-1} \mathrm{e}^{c_\ast y} 
\bP (\sqrt{2t^{\ast}_n} -m_n >y)
\end{equation}
exists,  where
$t_n^*$ denotes the number of \abbr{srw} excursions to the \abbr{root} before all
vertices of $\TT_n$ have been visited.
(In reality, one shows \eqref{righttail.1} with $\limsup$ and $\liminf$
replacing the limit in $n$, but we gloss over such details in this high level description.)
\newline 
$\bullet$ {\bf Step IV} To prove \eqref{righttail.1}, one uses the modified second moment method, with a twist. For a fixed vertex $\blue{\sfv} \in D_n$ \blue{and $y \gg 1$,  let $T^{(y)}_{\sfv, j}$} denote the number of excursions at the ancestor 
of $\sfv$ at level $j$ during the first $\frac{1}{2}(m_n+\blue{y})^2$ excursions of the \abbr{srw} at the root.  As in \cite{BRZ}, one notes that $\{ \blue{\eta^{(y)}_{\sfv,j}} := 
\sqrt{\smash[b]{2 T^{(y)}_{\sfv,j}}},  j\geq 1 \}$ is a (discrete) Bessel process. 
One then introduces the event
 $\blue{\wh{\UU}^{(+,n)}_\sfv (y)} :=\{\eta^{(y)}_{\sfv,j}\geq \gamma(j), j < n, \eta^{(y)}_{\sfv,n}=0\}$ 
 for a well-chosen barrier function \blue{$\gamma(j)$} and
 shows that the probability that the event in \eqref{righttail.1} occurs
 together with the complement of $\wh{\UU}^{(+,n)}_\sfv(y)$ for some $\sfv \in D_n$, is
  of lower order than our target as  $n\to\infty$ followed by $y \to\infty$.
\newline
$\bullet$ {\bf Step V} It is natural to expect that the probability in \eqref{righttail.1} is asymptotic to 
\blue{$\sum_{\sfv \in D_n} \bP\big(\wh{\UU}^{(+,n)}_\sfv(y)\big)$ when  $n \gg y \gg 1$.} 
  Unfortunately, this is in fact false. Instead, following
  a strategy from \cite{BDZ},
  one fixes a (large) constant
  $\sfKfv$ and considers the events 
  \begin{align*}
    &\blue{\UU^{(+,n)}_\bb(y)} =\{\eta^{(y)}_{\bb,j} \geq \gamma(j),  j \le n-\sfKfv,
  \exists \text{ descendant $\sfv$ of $\bb$ with $\eta^{(y)}_{\sfv,n}=0$}\}, \qquad  \bb \in D_{n-\sfKfv}. 
\end{align*}
Then, the second moment method shows that 
$\sum_{\bb \in D_{n-\sfKfv}} \bP\big(\UU^{(+,n)}_\bb (y)\big)
$ is for large $\sfKfv$,  a good proxy for the probability in \eqref{righttail.1}.
\newline
$\bullet$
 {\bf Step VI} Finally, one evaluates $\bP\big(\UU^{(+,n)}_\bb(y)
 \big)$
using extensively the barrier estimates from \cite{BRZ}.
  
All steps in this sketch require a significant modification in order
to work in the setup of $\bT^2$. We discuss these modifications next.

  For  Step I,  it is convenient to partition the torus $\bT^2$ into 
  boxes of side length about $e^{-\sfKz}$ which we call \textit{cities},
  separated by thin channels.
  The variable $N_k(t)$ is then replaced by the total number of excursions 
  across all channels. Note that the information contained in the 
  excursions across channels is considerably  more detailed than the information 
  contained in $N_k(t)$ on the tree. 
 This step  is carried out in \Cref{sec-exc-count}.

  Next, motivated by Step V in the preceding sketch, we tile cities by boxes of side length 
  $\vep  e^{\sfKfv}$ (also separated by thin channels), which we call \textit{blocks}.
With 
$n=\log(1/\vep)$ we replace the excursion counts from the tree by 
  the Brownian excursion counts between consecutive circles $\partial \bB_\sfd(\sfv_\bb,r_j)$ 
  of radii $r_j=\vep e^{n-j}$ around the center $\blue{\sfv_\bb}$ of \blue{the corresponding block $\bb$} 
  (see \Cref{subsec-circle}).  
  Per block this is again a discrete Bessel process, but in contrast with the situation for the tree,
  here it is no longer true that for two blocks which are, say,  of distance $r_{\ell+3}$ of each other,
  the corresponding Bessel processes coincide up to scale \blue{$\ell$}. 
  This is, technically, a major source of difficulty, which was addressed
  in \cite{BRZ2} by certain \textit{continuity estimates}.  Unfortunately, 
  the continuity argument in \cite{BRZ2} used in a crucial way the symmetries of
  the sphere $\bS^2$, and a significant effort is needed here to bypass this.
  These \blue{refined continuity} estimates form the core of \Cref{sec:cont}.

Another source of difficulty when dealing with $\bT^2$ instead of the tree is
that when conditioning on the excursions across channels, one obtains location
information that modifies the law of the excursions. To keep track of that, we
introduce,  as in \cite{BRZ2}, for each ball \blue{$\bB_\sfd(\sfv_\bb,r_j)$,}  the empirical measure of angle differences
\blue{of the excursions at level $j$ and control its}
distance from the ergodic law of such differences.  In particular, at a level $\sfKo > \sfKz$, 
we keep track, for each block,  of the number of excursions and their  angular ``defect''.  The knowledge of the empirical measure
of angle differences  at scale $j$,  together with the excursion counts at that level,  
are crucial for providing \blue{enough decoupling to make} the second moment method work.

With these preliminaries, we can now describe the modifications in our proof when moving from
$\TT_n$ to $\bT^2$. First, in lieu of Step II,  we prove  in  \Cref{sec:end-game}
that for $\sfKz,\sfKo$ large enough (with $\sfKo-\sfKz=c \log \sfKz$), 
the field of excursion counts and angular  defects at level $\sfKo$ converges in distribution 
to a Gaussian field.  Steps III and IV then proceed as in the tree,  except that one now 
has to consider the information provided by the field of excursion counts and 
angular defects. \blue{This leads to our more complex truncated events $\UUp$ and hence to 
considerably more effort in the reduction} (cf.  \Cref{sec-barrier} and \Cref{sec:apriori}).

Step V is similar to the analogous one in \cite{BRZ2}, except that 
\blue{lacking here} the symmetries of the sphere,  makes the 
control of second moments significantly more delicate, as 
the decoupling (which is essentially trivial on the tree), and a major 
part of \cite{BRZ2}, requires new arguments, see \Cref{sec:2nd-mom-tight}.

The heart of Step VI,  where we significantly depart from previous work,
lies in an \emph{exact} computation,  for a single block,
of the (conditional) probability of the event $\UUp$. 
Such a computation appears in \cite{BRZ2},  but with multiplicative constant
errors,  which are too big for our needs here.   Nevertheless,  using in a crucial
way Lemma \ref{recursionends2}
allows us in \Cref{sec:angle-diff} to decouple 
the angular information at level $\sfKo$ from \blue{uncovering the block
and from the excursion counts beyond a certain large level $\sfKth$,  whereupon the 
contribution from this part of our event is computed by utilizing results from \cite{DRZ}.}
We also use in a crucial way a decoupling idea in \cite{BRZ2} that controls the effect of angular defects on barrier probabilities. However, to get the level of precision needed,
we have to consider certain conditional local central limit theorems (\abbr{CLT}s), and obtain exact probabilities without constant multiplicative factors. This computation is contained in \Cref{sec:good-cover}.

\blue{We note in passing that} the only place in the argument where we crucially use the fact that we work on $\bT^2$ (and not on another compact connected manifold) is in Steps I and II:  we \blue{tile} $\bT^2$ by cities and in this tiling all cities are identical (up to a shift).  This allows us to employ a simple ergodic argument for controlling the convergence of the excursion counts and the angular  defects to Gaussian fields.
We do not know how to carry this step \blue{out} for other manifolds,  \blue{not even} for the two-dimensional sphere.

  \subsection{Structure of the paper}
\medskip
\Cref{sec:outline} \blue{carries out} the proof of \Cref{theo-covertime}
\blue{by establishing our analog of Step I and reducing the remaining tasks outlined above,  to:}
\newline
$\bullet$ \blue{The a priori} estimates of \Cref{prop-apriori},  derived in \Cref{sec:apriori},
\blue{which form our modification of Steps III and IV}. 
\\
$\bullet$ \blue{For Step VI,  the relevant} precise first moment asymptotic \blue{of \Cref{thm:end-game}}
is evaluated via Propositions \ref{prop:angle-diff} and 
\ref{prop:proxy-good-cover},  which \blue{in turn} are \blue{proven}
 in Sections \ref{sec:angle-diff} and \ref{sec:good-cover},
respectively.
\\
$\bullet$ \blue{Our version of Step V which is essentially the}
second moment bound of \Cref{prop:2nd-moment} \blue{and is proven} in \Cref{sec:2nd-mom-tight}
(where we also prove tightness,  see \Cref{prop:tight}).
\\
$\bullet$ \blue{Our replacement of Step II, which is the} conditional limit 
of \Cref{prop:Gauss-field-continuity},  and is established in \Cref{sec:end-game}.

After this reduction,  we proceed to prove in \Cref{sec:prelim} 
various preliminary facts that are needed for tackling those steps,
supplemented by the various coupling and continuity estimates of \Cref{sec:cont}
\blue{(which we shall utilize in Sections \ref{sec:apriori} and \ref{sec:end-game}).}
 
\subsection{Notation}
Throughout we order intersections of multiple events about excursions,
from the lowest level (largest radii),  to the highest (smallest radii),  
except for expectations of a non-indicator variable $X$ and events $A_i$,  written 
as $\bE[X;A_1,A_2,\cdots]$ even if $X$ involves the highest excursion level
(e.g.  see \eqref{eq:good-cover}). 
We also use $\displaystyle{\varlimsup_{\dlc,\sfKfv,\sfKz,\vep}}\!\! f(\cdot)$
and $\displaystyle{\varliminf_{\dlc,\sfKfv,\sfKz,\vep}} \!\! f(\cdot)$ for 
the limsup and liminf, of $f(\cdot)$ when first $\vep \to 0$, then $\sfKz \to \infty$ 
followed by $\sfKfv \to \infty$ and finally $\dlc \to 0$, using  
$\displaystyle{\mathop{\underline{\overline{\lim }}}_{\dlc,\sfKfv,\sfKz,\vep}} \!\! f(\cdot)$
to indicate that these two limits coincide 
\blue{(we similarly consider certain limits as $\sfKp \to \infty$ and then $\srs \to 0$,
between $\sfKfv \to \infty$ and $\dlc \to 0$).}
For a measure $\mu$ on a subset of a Euclidean space,  the notation $d\mu/dx$ denotes the Radon-Nikodym derivative of $\mu$ with respect to Lebesgue measure,  and when we use this notation 
it is implicitly assumed that $\mu$ is absolutely continuous with respect to the Lebesgue measure.  
We use the Landau notation $\Omega(\cdot)$,  $O(\cdot)$,  $o(\cdot)$,  and write e.g.
$o_{\sfKo}$ to emphasize the parameter in which the asymptotics is taken.

\section{Proof of \Cref{theo-covertime}: the decomposition to key steps}\label{sec:outline}
This section provides the key definitions and statements of the main intermediate steps, with pointers to subsequent sections where those are proved.  \blue{After introducing in} \Cref{subsec-cities} the notion of cities,  blocks and boundary zones,  we show that we can replace the cover time by
\blue{that}
 of the torus minus the boundary zones. In \Cref{sec-exc-count} we show that we can replace the cover time by a proxy, the count of inter-city excursions. This builds on a \abbr{CLT} for (Doeblin) Markov chains (see \Cref{lem:conc}) and is essentially the main place  where we use the structure of $\bT^2$. In \Cref{subsec-circle} we introduce excursion counts between concentric circles around a block,  \blue{which we use} 
 in \Cref{sec-barrier} \blue{to} introduce truncated events (including barriers and restrictions on the empirical measure of angular increments of excursions),  \blue{and to} state \Cref{prop-apriori},  an a priori \blue{control} (proved in \Cref{sec:apriori}) on the difference between the truncated events
and the original covering event.  \blue{\Cref{subsec-exact} establishes \Cref{thm:end-game}, 
our first moment estimate on the truncated events,  as a consequence of Propositions 
\ref{prop:angle-diff}, 
\ref{prop:proxy-good-cover},
\ref{lem:good-cover-end}
and 
\ref{prop:Gauss-field-continuity},
whose proofs are deferred to Sections \ref{sec:angle-diff},
\ref{sec:good-cover},  \ref{sec:good-cover-end} and \ref{sec:end-game},
respectively.}  \Cref{subsec-pftheo-covertime} \blue{concludes} the proof of \Cref{theo-covertime},  
by combining \Cref{thm:end-game},  the second moment computation of \Cref{prop:2nd-moment}
(which is proved in \Cref{sec:2nd-mom-tight}),  and tightness 
(\Cref{prop:tight}, again proved in \Cref{sec:2nd-mom-tight}).
\subsection{Cities, city blocks 
and boundary zones}
\label{subsec-cities}

For proving Theorem \ref{theo-covertime} we fix a small constant $\dlc>0$  
and large  constants $\sfKfv$, $\sfKz$ (we shall eventually take,  after $\varepsilon\to 0$,  limits in the constants so that
 $\sfKz \to \infty$ is
followed by $\sfKfv \to \infty$ and then $\dlc \to 0$).

At the core of our argument lies a tiling of the torus $\bT^2$ by $\sfMo=\sfMo(\sfKz)$ 
identical squares $\{\bA^\sfe_u\}$ of disjoint interiors. We denote by
$e^{-\sfKz}$ the edge-length of such a tiling 
square $\bA^\sfe_u$ which we call the $u$-th \emph{extended city} (with
$\sqrt{\sfMo}=e^{\sfKz}$ integer, see \Cref{rem:int-effects} for further restrictions and conventions). 
We 
further associate with each $\bA^\sfe_u$ an inner concentric smaller square $\bA_u$
of edge-length $(1-2\dlc)e^{-\sfKz}$ to which we refer as the $u$-th \textit{city}
and a further concentric sub-square of edge-length 
$(1-6\dlc)e^{-\sfKz}$ that we denote by $\bA_u^\sfo$ and 
call the $u$-th \textit{inner city}, while calling $\bA_u^\sfe \setminus \bA_u^\sfo$
the $u$-th city \textit{boundary zone} \blue{(see \Cref{fig:cities}(B))}.
Of course
our $\sfMo$ pairs 
$(\bA_u,\bA^\sfo_u)$ of city and inner-city,  are $\bT^2$-shifts 
of each other. 
 
Next, we decompose the Brownian path by recording its successive alternate 
hitting locations of the sets $\bTs^\sfe:=\cup_u \partial \bA^\sfe_u$ and 
$\bTs:=\cup_u \partial \bA_u$. 
Specifically, denote by $\sfz_{01}$ the  starting 
point of the Brownian path,  with $\sfz_{02}$ the location of the path at 
\blue{the stopping} time $T_0$,  
which is 
its first hitting time of $\bTs$.
We then set $\sfz_{j1} \in \bTs^\sfe$ and $\sfz_{j2} \in \bTs$,
$j \ge 1$,  as the alternating, successive  hitting locations of $\bTs^\sfe$ followed by 
$\bTs$ and  then $\bTs^\sfe$, etc.  We let   $T_j$, $j \ge 1$,  denote the time elapsed
between the visits to $\bTs$ recorded by the points 
$\sfz_{(j-1)2}$ and $\sfz_{j2}$. 

With these notations in place,  we
denote by $\{\HH_k, k \ge 0\}$ the filtration generated by the Brownian \textit{inter-city} 
excursions, with $\HH_0$ generated by the path up to time $T_0$ while 
$\HH_k$ for $k \ge 1$ is generated,  beyond the path in $[0,T_0]$, by  the part of 
the Brownian path during the first $k$ successive excursions from $\bTs^\sfe$ to $\bTs$,
which we recall start and end at
the (random) locations $\sfz_{j1}$ and $\sfz_{j2}$, respectively, 
with $j=1,\ldots,k$.
We further partition 
the Brownian segments which do not participate in the definition of
$\HH_\infty$, into the 
$\sfKz$-excursions of each city.
These are, for the $u$-th city, 
the successive Brownian 
excursions starting at $\partial \bA_u$
and ending upon hitting  $\bA^\sfe_u$ \blue{(See \Cref{fig:cities}(A))}.
Fixing $u$,  denote by $\HHu$
the $\sigma$-algebra generated by 
the whole Brownian path,  excluding only the $\sfKz$-excursions of the city $u$
and for each $k$ let $\tau_u(k)$ be
the $\HH_k$-measurable number of $\sfKz$-excursions of $u$
which are completed 
during the first $k$ inter-city excursions.  \blue{We emphasize that here and throughout the paper,
whenever we refer to a $\sigma$-algebra generated by a collection of paths (or excursions),  
it is} up to time parametrization,  \blue{excluding} the information on time in the definition of each path.  
That is,  each path is a (directed) curve up to time re-parameterization \blue{and by a collection 
of path we mean an ordered list of such directed curves}.

Denoting by
$\sfx^{(u)}_i := (\sfx^{(u)}_{i1},\sfx^{(u)}_{i2})$ the starting and ending points of the $i$-th such 
$\sfKz$-excursion, we consider the $\sigma$-algebras $\FF_{u,s} := \sigma(\sfx^{(u)}_i,i=1,\ldots,s)$ 
generated by the collection of 
such pairs that correspond to the first $s$ of the $\sfKz$-excursions 
of the $u$-th city.  In particular,  
\blue{
\begin{equation}\label{eq:FF-to-HH}
\FF_{u,\tau_u(k)} \subset \HH_k \subset \HH_\infty \subset \HHu  \qquad  \forall u \in [1,\sfMo],  \quad  \forall k \ge 1 \,.
\end{equation}}
Further,   if
$\sfz_{01} \notin \bA_u \setminus \partial \bA_u$ then
$\sfx^{(u)}_{i} = (\sfz_{(j-1)2},\sfz_{j1})$ for 
some $j \in [1,k]$ and all $i \le \tau_u(k)$.

We likewise further tile each inner city $\bA_u^\sfo$
by disjoint squares of edge-length $\vep e^{\sfKfv}$, 
which we call the \textit{extended (city) blocks} within that city, referring 
to the corresponding concentric inner squares of edge-length 
$\vep (1-6 \dlc) e^{\sfKfv}$, as the  blocks within that city, and calling 
the set difference between a block and the extended-block containing it,
the \emph{block boundary zone}.
The open ball of radius $r$ centered at $\sfv \in \bT^2$, corresponds to 
\[
\bB_\sfd (\sfv,r) := \{x \in  \bT^2 : \sfd(\sfv, x) < r \} \,, 
\]
where $\sfd(\cdot,\cdot)$ stands for the usual (Riemannian) distance on $\bT^2$. 

\begin{figure}[ht]
\centering

\begin{subfigure}[t]{0.48\textwidth}
\centering

\begin{tikzpicture}[scale=1.25]

\path[use as bounding box] (0,-0.9) rectangle (6,6);

\fill[gray!25]
(0,0) rectangle (6,6);

\foreach \i in {0,1,2,3}{
    \foreach \j in {0,1,2,3}{

        \fill[white]
            ({1.5*\i+0.15},{1.5*\j+0.15})
            rectangle
            ({1.5*\i+1.35},{1.5*\j+1.35});

        \draw[line width=0.7pt]
            ({1.5*\i+0.15},{1.5*\j+0.15})
            rectangle
            ({1.5*\i+1.35},{1.5*\j+1.35});
    }
}

\draw[line width=0.8pt] (0,0) rectangle (6,6);

\foreach \x in {1.5,3,4.5}
    \draw[line width=0.8pt] (\x,0) -- (\x,6);

\foreach \y in {1.5,3,4.5}
    \draw[line width=0.8pt] (0,\y) -- (6,\y);
\draw[black,thick, dashed] plot [smooth,tension=1] coordinates {(0.1,2.4) (0.8,2.5)  (1.2,1.8) (1.5,2)};
  \draw [red,line width=1.5pt] plot [smooth, tension=1] coordinates {(1.5,2) (1.4,2.5) (1.6,3) (1.65,3.3)};
 \draw[black,thick,dashed] plot [smooth,tension=1] coordinates {(1.65,3.3) (1.85,3.6) (2.1,3.2) (2.5,4.5)};
  \draw [red,line width=1.5pt] plot [smooth, tension=1] coordinates {(2.5,4.5) (3,4.6) (3.6, 4.4) (4,4.4)};
\draw[black,thick,dashed] plot [smooth,tension=1] coordinates {(4,4.4) (3.6,3.6) (4.1,3.8) (4.5,4)};
  \draw [red,line width=1.5pt] plot [smooth, tension=1] coordinates {(4.5,4) (4.6,3.5) (4.4,3.2) (4.5,2.2) (4.65,2)};
\draw[black,thick,dashed] plot [smooth,tension=1] coordinates {(4.65,2) (5,2.82) (5.5,2.75) (5.8, 2.1) (5.3,1.5)};
  \draw [red,{-Stealth},line width=1.5pt] plot [smooth, tension=1] coordinates {(5.3,1.5) (5,1.4) (4.5,1.5) (4.5,1) (4.5,0.8) (4.6,0)};
  
  \node at (3.85,0.85){\small $\bA_u^\sfo$};

\node at (0.1,2.4) [left] {\small $\sfz_{02}$};
\node at (1.6,2) [right] {\small $\sfz_{11}$};
\node at (1.65,3.6) [right] {\small $\sfz_{12}$};
\node at (2.5,4.6) [above] {\small $\sfz_{21}$};
\node at (4,4.6) [above] {\small $\sfz_{22}$};
\node at (4.6,4) [right] {\small $\sfz_{31}$};
\node at (4.65,2) [right] {\small $\sfz_{32}$};
\node at (5.3,1.4) [below] {\small $\sfz_{41}$};

\draw (0.1,2.4) circle (0.5pt);
\draw (1.5,2)  circle (0.5pt);
\draw(1.65,3.3)  circle (0.5pt);
\draw (2.5,4.5)  circle (0.5pt);
\draw (4,4.4)  circle (0.5pt);
\draw (4.5,4)  circle (0.5pt);
\draw (4.65,2)  circle (0.5pt);
\draw (5.3,1.5)  circle (0.5pt);

\end{tikzpicture}

\caption{}
\label{fig:blocks}
\end{subfigure}
\hfill
\begin{subfigure}[t]{0.48\textwidth}
\centering

\begin{tikzpicture}[scale=1.25]

\path[use as bounding box] (0,-0.9) rectangle (6,6);
\fill[gray!25]
 (0,0) rectangle (6,0.25);
\fill[gray!25]
 (5.75,0) rectangle (6,6);
\fill[gray!25]
 (0,5.75) rectangle (6,6);
\fill[gray!25]
 (0,0) rectangle (0.25,6);
\draw (0,0) rectangle (6,6);
\draw (0.25,0.25) rectangle (5.75,5.75);
\draw (0.65,0.65) rectangle (5.35,5.35);

\draw (2,4) circle(0.70);
\draw (2,4) circle(0.85);

\draw[dashed] (3.85,3.85) rectangle (4.15,4.15);
\draw (1.8,3.8) rectangle (2.2,4.2);
\draw [dashed] (1.85,3.85) rectangle (2.15,4.15);
\draw (1.8,2.1) rectangle (2.2,2.5);
\draw [dashed] (1.85,2.15) rectangle (2.15,2.45);


\draw[thick,<->] (0,3) -- (0.25,3);
\draw[thick,<->] (0,1) -- (0.65,1);
\draw[thick,<->] (0,-0.2) -- (6,-0.2);
\draw[thick,<->] (1.8,2.5)--(2.2,2.5);
\draw[thick,<->] (3.85,4.15)--(4.15,4.15);
\node at (-0.5,2.95)
    {\footnotesize $\dlc e^{-\sfKz}$};

\node at (2.0,2.75) 
    {\footnotesize $\vep e^{\sfKfv}$};

\node at (2.1,1.9)
    {\footnotesize Extended block};

\node at (1.35,1)
    {\footnotesize $3\dlc e^{-\sfKz}$};
    
\node at (4.1,4.4)
    {\footnotesize $(1-6\dlc) \vep e^{\sfKfv}$};    

\node at (3,-0.5)
    {\footnotesize $e^{-\sfKz}$, Extended city, $\bA_u^{\sfe}$};

\node at (4,1)
    {\footnotesize Inner city, $\bA_u^\sfo$};

\node at (4.1,3.6)
    {\footnotesize City block};

\end{tikzpicture}

\caption{}
\label{fig:outer}
\end{subfigure}
\caption{(A) Inter-city excursions (solid red) and excursions within a city (dashed black); (B) A city and its blocks.}
\label{fig:cities}

\end{figure}

\begin{remark}\label{rem:int-effects}
The attentive reader must have realized that our tiling of the inner cities by extended blocks
requires in addition to $\sqrt{\sfMo}=e^{\sfKz}$ integer, also
that $(1-6 \dlc) e^{-\sfKz}$ be an integer multiple of $\vep e^{\sfKfv}$. We take 
care of this integer effect by adjusting $\dlc$ to some $\vep$-dependent 
$\dlc(\vep)$ which is within $\frac{1}{6} \vep e^{\sfKz+\sfKfv}$ of $\dlc$.  As this does
not impact any of our arguments, we stay with $\dlc$ for the sake of keeping compact notation.
\end{remark}

Our first step is to replace $\CCT$ by certain modified normalized $\vep$-cover times.  These will achieve two goals: first, allow for measuring time in terms of city excursions, and second, allow us to study the cover time of $\bT^2$ except for a patch near the starting point. This will be useful in deriving certain a priori barrier estimates. Readers primarily interested in the main argument  may on a first reading disregard these variations of the normalized cover time.
\begin{lemma} 
Let $\CCdd$ denote
the normalized $\vep$-cover time of certain sets $\bT_\delta \subset \bT^2$ chosen independently of the 
path on $\bT^2$.  If the area of the closure of 
$(\bT_\delta)^\complement$ decays to zero with $\delta \to 0$,  then starting the path uniformly on $\bT^2$,
\label{lem-claim1} 
\begin{equation}\label{eq:removal}
\lim_{\delta \to 0} \sup_{\vep \le 1/2} \bP(\CCdd \ne \CCT ) = 0 \,.
\end{equation}
\end{lemma}
\begin{proof} Clearly,  $\CCdd \le \CCT$ with equality \blue{unless the 
last $\vep$-ball within $\bT^2$ to intersect the Brownian path,  is centered at $(\bT_\delta)^\complement$.}
Starting our Brownian motion uniformly over $\bT^2$ results 
by symmetry,  with  a uniform distribution for the center of this (random) $\vep$-ball.  Thus,  per 
choice of $\delta$ and $\vep$,  the probability in \eqref{eq:removal}
 is precisely the 
area of the closure of 
$(\bT_\delta)^\complement$,  which by assumption goes to zero as $\delta \to 0$.
\end{proof}

\blue{
\begin{definition}\label{dfn:CCd}
Denote by $\CCd=\CCd (\sfKz,\sfKfv)$ the normalized $\vep$-cover time,  excluding the 
open balls $\bB_\sfd(\cdot,\vep)$ that intersect some city boundary zone or some block 
boundary zone,  and those balls which are within an extended city that intersects the 
square of edge-length $\srs \le e^{-1}$ around the Brownian motion starting point.
\end{definition}}

Hereafter we use $\displaystyle{\varlimsup_{\dlc,\sfKfv,\sfKz,\vep}}\!\! f(\cdot)$
and $\displaystyle{\varliminf_{\dlc,\sfKfv,\sfKz,\vep}} \!\! f(\cdot)$ for 
the limsup, respectively the liminf, of $f(\cdot)$ when first $\vep \to 0$, then $\sfKz \to \infty$ 
followed by $\sfKfv \to \infty$,  $\dlc \to 0$,  using  
$\displaystyle{\mathop{\underline{\overline{\lim }}}_{\dlc,\sfKfv,\sfKz,\vep}} \!\! f(\cdot)
$
to indicate that these two limits coincide (we similarly consider certain limits as $\sfKp \to \infty$ 
\blue{and then $\srs \to 0$,}  between $\sfKfv \to \infty$ and $\dlc \to 0$).

\begin{remark}\label{rem:good-ic}
Since uniformly in $\sfKz$,  the aggregate area of the \blue{$\vep$-blow-up of} 
boundary zones of all blocks and cities
goes to zero \blue{as $\vep \to 0$ then $\sfKfv \to \infty$ and} $\dlc \to 0$,  it follows from  \Cref{lem-claim1} that
\begin{align}
\varlimsup_{\dlc,\sfKfv,\sfKz,\vep} \bP \big( \CCz \ne \CCT  \big) = 0.
\label{ct2.25cj}
\end{align}
Note that the law of $\CCT$ is independent of the initial law of the Brownian motion
and \eqref{ct2.25cj} applies for any initial law which is uniform within each $\bA^\sfe_u$
(or for any initial law,  if we also uniformly shift,  independently of the Brownian path,
our  tiling of $\bT^2$ by cities).  Moreover,  in \Cref{prop:tight} we will show that 
\begin{align}
\lim_{\dlc,\srs} \varlimsup_{\sfKfv,\sfKz,\vep} \bP \big( \CCd \ne \CCT \big) = 0 \,,
\label{ct2.25-new}
\end{align}
so \Cref{theo-covertime} is a consequence of having that at any fixed $z \in \bR$,
\begin{align}
{\mathop{\underline{\overline{\lim }}}_{\dlc,\srs,\sfKfv,\sfKz,\vep}}
\bP \big( \CCd \le  (\osfm_\vep + z )^2  \big) = \bP(Y_\star  \le z) \,,
\label{ct2.25c}
\end{align}
where per $\dlc$ and $\srs$ we consider a non-random subset of $\sfM = (1-\srs^2+o_{\sfKz} (1)  )\sfMo$ 
cities and start the path at a point on $\bTs=\cup_u \partial \bA_u$
of distance at least $\srs/e$ from all blocks in those cities.
\end{remark}

\subsection{Inter-city excursion count as proxy of time}\label{sec-exc-count}
$~$
Per \Cref{rem:good-ic} 
we start \abbr{wlog} at some $\sfz_{01}=\sfz_{02} \in \bTs$
(namely with $T_0=0$).
Recall the preceding notation of $\{T_j\}$ and note that the law of the Markov pairs 
$\{(\sfz_{(j-1)2},T_j), j \ge 1\}$
depends only on $\sfKz$, $\dlc$ and the tiling $\{(\bA^\sfe_u,\bA_u)\}$. 
Thanks to continuity and positivity of the Poisson kernel, and the compactness of $\bTs$, the 
conditional law of $\sfz_{12}$ given $\sfz_{02}$ 
has a uniformly bounded away from zero density with respect to Lebesgue measure on $\bTs$.
The Markov transition probability  of $(\sfz_{12},T_2)$ given $(\sfz_{02},T_1)$ is composed 
of the law of $\sfz_{12}$ given $\sfz_{02}$ and the conditional law of $T_2$ given $\sfz_{12}$,
whose \abbr{mgf} is continuous in $\sfz_{12}$ and finite in some 
neighborhood of $0$. Utilizing again the compactness of $\bTs$ we 
deduce that these Markov transition probabilities 
form a Doeblin aperiodic chain, satisfying 
\cite[(16.9)]{Meyn-Tweedie}
(with respect to the product 
of Lebesgue measure on $\bTs$ and an Exponential($c$) density
on $\bR_+$). As such, this uniformly ergodic 
chain has a unique invariant measure $\mu_\sfe$ 
on $\bTs \times \bR_+$ to which its transition semi-group 
converges uniformly with respect to the initial state
(see \cite[Theorem 16.0.2(ii)]{Meyn-Tweedie}). This yields
the following sharp concentration for the time 
\begin{equation}\label{dfn:Sk}
S_k := \sum_{j=1}^k T_j \,,
\end{equation}
it takes the Brownian path to complete $k$ inter-city excursions.
\begin{lemma}\label{lem:conc}
Fixing $\dlc>0$ and $\sfKz$,  let $\sft_1 := \sft_1(\sfKz,\dlc) =  \bE_{\mu_\sfe} [T_1]$.
The sequence $\big\{ \sqrt{S_k} - \sqrt{k\, \sft_1}\big\}$ converges in distribution as $k \to \infty$
to a normal variable of zero mean,  whose finite variance $\sfv_\sfe/(4 \sft_1) \to 0$
when $\sfKz \to \infty$.
\end{lemma}
\begin{proof} By the preceding reasoning, and in particular, the uniform bound on
$\bE[e^{c T_j}|\sfz_{(j-1)2}]$, the ergodic 
chain $\{(\sfz_{(j-1)2},T_j), j \ge 1\}$, is further $V$-uniformly ergodic 
for $V(z,t)=e^{ct}$ and some $c>0$ (see \cite[(16.2)]{Meyn-Tweedie}).
By the \abbr{clt} for such chains (e.g. \cite[Theorem 17.0.1]{Meyn-Tweedie}), 
we have that, regardless of $\sfz_{02} \in \bTs$,
the sequence $(S_k - k \sft_1)/\sqrt{k}$ converges in law to a normal variable 
of zero mean and the finite variance 
\begin{equation}\label{eq:var-erg}
\sfv_\sfe :=
\lim_{k \to \infty} \frac{1}{k} \bE_{\mu_{\sfe}} \Big[ (S_k - k \sft_1)^2 \Big] = 
{\sf var}_{\mu_\sfe} (T_1) + 2 \sum_{j=2}^\infty {\sf cov}_{\mu_\sfe} (T_j,T_1)\,,
\end{equation}
where ${\sf cov}_{\mu_\sfe}(\cdot,\cdot)$ denotes the covariance under 
$\mu_\sfe$.  It thus remains only to show that 
\begin{equation}\label{eq-amirswriteup}
\sfv_\sfe = o_{\sfKz}( \bE_{\mu_\sfe} [T_1] ) \,.
\end{equation}
To see \eqref{eq-amirswriteup}, let $W_t$ denote the Brownian motion on $\bT^2$ 
equipped with a coordinate system so $(0,0)$ is a left-lower corner of the 
extended city $\bA^\sfe_{1}$.  Recall from our tiling that $e^{\sfKz} \in \bZ_+$,
hence the process 
$Y_t = W_t \!\!\! \mod \!\! (e^{-\sfKz})$ is a Brownian motion on 
a torus $
\bA^\sfe_1$ of edge-length $e^{-\sfKz}$. By  
the invariance of $\bTs$ and $\bTs^\sfe$ to $(e^{-\sfKz} \bZ)^2$-translations 
on $\bT^2$, the sequence $\{T_j$, $j \ge 1\}$ merely records the successive 
times it takes $\{Y_t\}$ to get from $\partial \bA_1$ to $\partial \bA^\sfe_1$ 
and back to $\partial \bA_1$, while our folding of $\bT^2$ 
maps $\mu_\sfe$ to the joint distribution of such $T_1$ (for $Y_t$), and the unique
invariant hitting measure of $\partial \bA_1$ 
associated with the preceding excursions. By the scale invariance of Brownian motion, 
it thus follows that $\sfv_\sfe$ of \eqref{eq:var-erg} scales as
$e^{-4\sfKz}$, whereas $\sft_1 = \bE_{\mu_\sfe} [T_1]$ scales as $e^{-2\sfKz}$, 
thereby establishing \eqref{eq-amirswriteup}.
\end{proof}

Next,  let $\UUk$ denote the event 
of having an unhit $\vep$-ball within the block $\bb$ by 
the end of the $k$-th inter-city excursion,  with the corresponding 
$\vep$-cover times \blue{measured in terms of inter-city excursions}, 
\begin{equation}\label{cov-time-exc}
\CCex := \inf\{ k \ge 1  :  \bigcup_\beta \,  \UUk  = \emptyset \} \,,
\end{equation}
where $\beta$ runs over our $\vep$-dependent set of city blocks within the $\sfM$
inner cities \blue{of sufficient distance from the starting point}).

In view of Lemma \ref{lem:conc} we further replace 
$\CCd$ in \eqref{ct2.25c} by $(2 \pi) \sft_1 \CCex$
(with the factor $2\pi$ due to the re-scaling on the right hand side of \eqref{eq:dprz} when
$\bM=\bT^2$).
Specifically,  setting for integer $\ello > 15$, 
\begin{align}\label{def:mn-alt}
r_0:=e^{-\ello} \quad {\rm and} \quad 
\sfm_\vep =\blue{ \sfm_{\vep,r_0}} :=  \sfrho_\vep \log (r_0/ \vep) = \osfm_\vep - 2 \ello +o_\vep(1) \,,
\end{align}
\begin{equation}\label{dfn:bar-svep}
\sfk_\vep (z) := \bar \sfs_\vep(z) / (\pi \sft_1)\,, \qquad 
\bar \sfs_\vep(z)  := \frac{1}{2} \big(\sfm_\vep + z)^2 \,,\qquad  z \in \bR \,,
\end{equation}
\blue{we see that for small $\dlc$,  $\srs$ and any $z$,  $\sfKfv$, }
\begin{align}
\varlimsup_{\sfKz,\vep} \Big| \bP \big( \CCd \le  (\sfm_\vep + z)^2  \big) -  \, \bP \big(  \, \CCex \le  \sfk_\vep(z)  \big) \Big| = 0 \,.
\label{ct2.25cjj}
\end{align}
Taking $[\sfk_\vep(z)]$ instead of $\sfk_\vep(z)$ corresponds to $o_\vep(1)$ change in $z$,  so 
\abbr{wlog} we assume hereafter that $\sfk_\vep(z) \in \bN$.
Further, as $k \mapsto \UUk$ are non-increasing,  partitioning our 
blocks per inner city $\bA_u^o$ address,  results with
\begin{equation}\label{dfn:UUz}
\{ \CCex \le \sfk_\vep(z) \} = \{ \forall u  \in [1,\sfM] :  \bigcup_{\beta \in u} \UUz = \emptyset \}  \,,\quad
{\rm where} 
\quad \UUz:=\UU^{(\vep)}_{\bb, \sfk_\vep(z)} \,.
\end{equation}
In particular,  setting $\FFuz := \FF_{u,\tau_u(\sfk_\vep(z))}$ and the $[0,1]$-valued random variables
\begin{equation}\label{def:barQu}
Q^{(\vep)}_{u,z} =\bP \big[    \Lambda^{(\vep)}_{u,z} \ge 1 |  \FFuz \big], \qquad 
\Lambda^{(\vep)}_{u,z} := \sum_{\beta \in u} {\bf 1}_{\UUz} \,,
\end{equation}
we have that for any $z \in \bR$,
\begin{equation}\label{eq:CC-Quz}
\bP(\CCex \le \sfk_\vep(z) ) = \bE \Big[   \prod_{u=1}^{\sfM}  (1- Q^{(\vep)}_{u,z} ) \Big] \,.
\end{equation}
\blue{In view of \eqref{ct2.25c},  \eqref{def:mn-alt} and \eqref{ct2.25cjj},  we will estimate the right hand side for
$\srs=r_0$ and $Q^{(\vep)}$ evaluated at $z+2\ello$}.

\subsection{Excursions between concentric circles}\label{subsec-circle} Set  
\blue{$\srs=r_0$ of} \eqref{def:mn-alt} and 
\begin{equation}\label{def.nvep}
\log (r_0/\vep)=n_\vep+\sfae,  \quad {\rm for} \quad n_\vep \in \bN \quad {\rm and}  \quad \sfae \in [0,1) \,.
\end{equation}
In the sequel we introduce integers $\sfK_i$,  $1 \le i \le 6$,
setting
\begin{equation}\label{def.n-vep-i} 
n_\vep^{(i)}:=n_\vep-\sfK_i  \quad \blue{{\rm for} \quad 4 \le i \le 6 \quad {\rm and} \quad \wnvfv := \nvfv - \log \dlc \,.}
\end{equation}
To assist us in dealing with \blue{the right hand side of \eqref{eq:CC-Quz}} let $\sfKtp:=(\sfKt+\sfKth)/2$ and 
introduce for each block $\beta$ the collection of concentric circles $\partial \bB_\sfd(\sfv_\beta,r_j)$ 
centered at its center $\sfv_\beta$,  with radii 
\begin{equation}\label{dr.1}
r_j=r_0 e^{-j},  \quad j \in \{0,\ldots, \sfKtp\},  
\;\;
r_j=\vep e^{n_\vep-j},  \;\;  j \in \{\sfKtp+1, \ldots,  n_\vep\} \,.
\end{equation}
In particular,  $r_{n_\vep}=\vep$ with ratio $e$ between consecutive radii,  
except for the larger ratio $r_{\sfKtp}/ r_{\sfKtp+1}= e^{1+\sfae}$ to bridge from 
independent of $\vep$,  city-like,  large scales,  to small,  block-like 
scales,  proportional to $\vep$.  Consider for $j \in [0,\nvfv]$
the various Brownian excursions which start at $\partial \bB_\sfd(\sfv_\beta,r_j)$
and end upon first exiting $\bB_\sfd(\sfv_\beta,r_{j-1})$, which we call the 
\textit{$j$-level excursions} of block $\beta$
(all such circles are inscribed within the
fundamental domain of $\bT^2$ around $\sfv_\beta$). 
Note that for such a choice of 
radii,  except at $\sfKtp$,  from $\partial \bB_\sfd(\sfv,r_j)$ the Brownian path is 
equally likely to hit $\partial \bB_\sfd(\sfv,r_{j-1})$ prior to $\bB_\sfd(\sfv,r_{j+1})$, 
or vice versa.  Further,  fixing 
$\sfKf=\sfKfv+5 \log \sfKfv$ and $\sfKo=\sfKz+5 \log \sfKz$ 
(in particular,   \abbr{wlog} take $\sfKz$ large enough so $e^{-\sfKo} \le \dlc e^{-\sfKz}$
and see \Cref{rem:more-int-effects} on why we ignore various integer effects),
we shall focus on block excursions at levels $j \in [\sfKo,\nvf]$,  see \Cref{fig:village}.
\begin{figure}[h]
\begin{center}
\begin{tikzpicture}
\draw (0,0) circle (3); 
\draw (-0.1,-0.1) rectangle (0.1,0.1);
\draw (0,0) circle (3.2); 
\draw (0,0) circle (1); 
\draw (0,0) circle (1.15); 
\draw[thick,->] (0,0) -- (0,3.2);
\node at (0.5, 2) {\small $r_{\sfKo-1}$};
\draw[thick,->] (0,0) -- (-0.3,3);
\node at (-0.5, 2.5) {\small $r_{\sfKo}$};
\draw[thick,->] (0,0) -- (-1.1,-0.34);
\node at (0.05, -0.36) {\tiny $r_{\nvf-1}$};
\draw[thick,->] (0,0) -- (0.60,0.81);
\node at (0.6, 0.1) {\tiny $r_{\nvf}$};
\end{tikzpicture}
\end{center}
\caption{Excursion radii}
\label{fig:village}
\end{figure}

An advantage over \cite{BRZ,BRZ2}, where the  
excursions are from $\partial \bB_\sfd(\sfv,r_{j-1})$ to 
first entering $\bB_\sfd(\sfv,r_j)$, is that our excursions, per fixed $j$, $\sfv$,
and conditional on their starting and ending positions, are
independent of the whole Brownian path outside $\bB_\sfd(\sfv,r_{j-1})$. Nevertheless, 
we freely use estimates from \cite{BRZ,BRZ2} for excursion counts such as
\begin{align}
T_{j}^{\sfv,\sfKo,s} &:= \#\{  
\mbox{\rm  
excursions from } \partial \bB_\sfd(\sfv,r_j) \;\;
 \mbox{\rm  to } \partial \bB_\sfd(\sfv,r_{j-1})  
\; \mbox{\rm during} 
\nonumber \\
& \quad [s] \; \; \mbox{\rm first 
 excursions from }
 \partial \bB_\sfd(\sfv,r_{\sfKo-1}) \;\;
 \mbox{\rm  to } \partial \bB_\sfd(\sfv,r_{\sfKo}) 
\} 
 \label{Tjsdef}
\end{align}
(where $s \ge 1$,  $j \in  [\sfKo,n_\vep]$ and $\sfv \in \bT^2$),  since within each 
$\sfKo$-level excursion, the number of contained,  $j$-level excursions 
will be the same whether one counts their inward parts (towards $\sfv$),
as in \cite{BRZ,BRZ2},
or their outward/return parts,
as we do here.  The collection of 
counts $\{T^{\sfv_\beta,\sfKo,s}_j\}$ 
is invariant under rotations of the Brownian path
around $\sfv_\beta$, thereby allowing for a relatively precise
computation of probabilities of events involving only these variables, via the
\abbr{srw} on $\bZ$ induced by the order in which the path hits the circles $\partial \bB_\sfd(\sfv_\beta,r_j)$.
In contrast, due to the rectangular shape of blocks and cities, certain angles affect the probabilities
of events within  the blocks, and outside the cities (in particular affecting the law of $Q^{(\vep)}_{u,z}$). 
As shown in the sequel,  the preceding choice of $\sfKo-\sfKz$ and $\sfKf-\sfKfv$ allows us to effectively 
decouple the block's interior and city's exterior from the intermediate excursions tower profile 
(namely, $\{T^{\sfv_\beta,\sfKo,s}_j\}$).

\subsection{Truncation: barriers,  angles and effective $\vep$-cover}\label{sec-barrier}
Fixing $\dlc>0$,  $\srs=r_0$ and parameters 
\begin{equation}\label{dfn:aux-par}
0< \wt{c},  \gak,  \quad 3\gak \vee \frac{1}{6} < \kapo<\frac{1}{4}<\frac{1}{2}< \kapp < \okappa < \frac{3}{4},  \quad \cst<\infty,  
\end{equation}
similarly to \cite{BRZ},  we replace the event $\UUz$ that block $\bb$ is
$\vep$-uncovered,  by truncated events that are more convenient to handle.  
\blue{Specifically,   here}
they consist of the restricted \blue{$\vep$-uncovering} event $\RUb$ of \eqref{def:res-u-z},  
which we further intersect with the barrier restrictions $\BBb$ of \eqref{eq:barrier-event-z} and 
the angle restrictions $\WWb$ of  \eqref{eq:was-event-z},  to form the truncated events
\begin{equation}\label{eq-good}
  \UUp := \RUb \cap \BBb \cap \WWb \subset\UUz \,.
\end{equation}
We will establish in  \Cref{sec:apriori} the following a priori estimate for our truncation.
\begin{proposition}\label{prop-apriori}
For large $\cst$ and small $\wt{c}$,  $\srs=r_0$,  fixing $z \in \bR$,  $\dlc>0$ and $\gak,\kapo,\kapp,\okappa$
as in \eqref{dfn:aux-par},
\begin{align}\label{eq:U-minus}
 \varlimsup_{\sfKp,\sfKfv,\sfKz,\vep} & \; \bP \Big(
\bigcup_\beta \big\{ \UUz \setminus \UUp  \big\} 
 \Big) = 0 \,.
\end{align}
\end{proposition} 
\blue{Setting} analogously to \eqref{def:barQu},  the counts and conditional \blue{probabilities}
\begin{align}\label{def:barQub}
Q^{(+,\vep)}_{u,z} &:= \bP \big[    \Lambda^{(+,\vep)}_{u,z} \ge 1 |  \FFuz \big], \qquad 
\Lambda^{(+,\vep)}_{u,z} := \sum_{\beta \in u} {\bf 1}_{\UUp} \,,
\end{align}
\blue{we deduce from \eqref{ct2.25cjj},  \eqref{eq:CC-Quz},  \eqref{eq-good} and \eqref{eq:U-minus} that}
\begin{equation}\label{eq:CC-Quz-plus}
\varlimsup_{\sfKp,\sfKfv,\sfKz,\vep} \Big|  \bP \big( \CCd \le  (\sfm_\vep + z)^2  \big)
- \bE \big[   \prod_{u=1}^{\sfM}  (1- Q^{(+,\vep)}_{u,z} ) \big] \Big| = 0 \,.
\end{equation}
\blue{To recap,  recalling \eqref{def:mn-alt} we get \Cref{theo-covertime} upon} showing that for $\srs=r_0=e^{-\ello}$,  
\begin{equation}\label{eq:ct-exc-plus}
{\mathop{\underline{\overline{\lim }}}_{\dlc,\srs,\sfKp,\sfKfv,\sfKz,\vep}} \,
\bE \Big[   \prod_{u=1}^{\sfM}  (1- Q^{(+,\vep)}_{u,z+2\ello} ) \Big] 
 = \bP(Y_\star \le z)\,.
\end{equation}
\begin{remark}\label{rem:more-int-effects} 
To define $\UUp$ we set $\vep$-dependent levels 
$\sfKw:=(21 \log n_\vep)^{1/\gak}$,  $\sfKt:=\sfKo + 2 (21 \log n_\vep)^{1/\gak}$
and $\sfKth:=\sfKt+ 5 \log n_\vep$.  Our proofs are essentially unaffected 
by inserting an integer part in our specifications of $\sfKo$,  $\sfKw-\sfKo$,  
$\sfKt-\sfKo$,  $\sfKth-\sfKt$ and $\sfKf-\sfKfv$ (thereby yielding integer
$\sfKw$ and $\sfK_i$, $1 \le i \le 5$),  and for a lighter notation we proceed without doing so. 
Indeed,  all our arguments rely only on having gaps of suitable order between 
$\sfKz<\sfKo < \sfKw<\sfKt<\sfKth$ and between $\sfKf >\sfKfv > \sfKp$, while taking first 
$\sfKz \to \infty$ before $\sfKfv \to \infty$, and only at the very end considering $\sfKp \to \infty$.
\end{remark}

\blue{We next}  
specify the events $\UUp$ 
in \eqref{eq-good} \blue{via}
the following three steps.\\
{\bf Step 1 - nets,  barriers and excursion counts}.\\
$\bullet$ \blue{{\bf Nets}}:
Fix $\widetilde{c} \in (0,\frac{1}{2})$ 
and \abbr{wlog} assume that $e^{\sfKfv}/(\sqrt{2} \wt{c})$ is a dyadic integer.  With
\begin{align}\label{dfn:F-star}
\bF^\star& :=  \sqrt{2} \,  \wt{c} \,  \bZ^2 \cap [-\sfL,\sfL]^2  \quad \text{for} \quad 2 \sfL := (1-6\dlc) e^{\sfKfv}
\end{align}
a minimal $\wt{c}$-net for 
$[-\sfL,\sfL]^2$,  denote 
the corresponding $\widetilde{c} \, \vep$-nets for each block $\bb$ and for their union
by
\begin{align}\label{dfn:F-beta}
\bF^\beta := \sfv_\beta + \vep \bF^\star \quad \text{and} \quad \bF_{n_\vep} &:= \bigcup_{\beta} \bF^\beta\, .
\end{align}
Clearly,  $|\bF^\bb|=|\bF^\star|=c(\wt{c},\dlc) e^{2\sfKfv}$ 
and the dyadic dilutions $\{ 2^\ell \bF_{n_\vep},  \ell \ge 0 \}$ of $\bF_{n_\vep}$ 
form nested sub-lattices $\bF_j$ which are $\wt{c}\,  r_j$-nets for the union of all blocks in $\bT^2$
(by taking $\ell=\lfloor \log_2 (r_j/\vep) \rfloor$),  such that 
for some $c=c(\wt c,r_0)$ finite, 
\begin{equation}
\left|\bF_{j}\right| \le c \,  e^{2j},  \qquad \forall j \in [0,n_\vep] \,.
\label{dr.2}
\end{equation}
$\bullet$ \blue{{\bf Barriers}: Our barriers for excursion counts are based on} the curves
\begin{align}\label{eq:norm-exc}
\gamma_{b} (j+1):= \sfrho_\vep \log(r_j/\vep)+ b j_{n_\vep}^\gak\;\;
\hbox{where}  \;\;\;  j_\ell := j \wedge (\ell-j), \quad j \in [0,\ell], \quad \ell\in \bN \,.
\end{align}
In particular,  for any $b \in \bR$, 
\begin{align}
\gamma_b(1)&=\sfm_\vep,  \;\; \text{with} \;\; \gamma_0(\sfK_i)=\sfm_\vep - \sfrho_\vep (\sfK_i-1),  \;\; 1 \le i \le 2,  \nn \\
\gamma_b(n_\vep+1)&=0,\;\;\;\; \text{with} \;\;\; \gamma_0(n_\vep^{(i)})=\sfrho_\vep (\sfK_i+1),   \qquad \;  4 \le i \le 6  \,.
\label{eq:special-norm}
\end{align}
\begin{definition}\label{rem:nominal}
Using the notation of \eqref{eq:norm-exc},  we set for $j \ge 0$ and $y \in \bR$, 
\begin{equation}\label{eq:nominal}
 \sfs_j (y)  :=  \frac{1}{2} ( \gamma_0(j)  + y)^2
\end{equation}
(so $\sfs_1(y)=\bar \sfs_\vep(y)$ of \eqref{dfn:bar-svep}).  
\end{definition}
\noindent $\bullet$ \blue{{\bf Excursion counts}}: \blue{We set similarly to \eqref{Tjsdef},  for 
$j \in [\sfKo,n_\vep]$ the following excursion related counts.
\begin{definition}\label{def:ex-cnt}
For $z \in \bR$ and $\sfv \in \bT^2$,  let $\etavjz := \sqrt{\smash[b]{2\Tvjz}}$ for the count 
\begin{align}
\Tvjz := \#\{ &  \mbox{\rm excursions from } \partial \bB_\sfd(\sfv,r_j) \;\; \mbox{\rm to } 
\partial \bB_\sfd(\sfv,r_{j-1})  \; 
\mbox{\rm during the first } \; \sfk_\vep(z) \mbox{\rm inter-city excursions} \}.
\end{align}
Similarly,  for any block $\bb$,  $y \in \bR$ 
and $\sfv \in \bF^\beta$ set $\etavjy := \sqrt{\smash[b]{2 \Tvjy}}$ for the count
\begin{align}
\Tvjy := \#\{ & 
\mbox{\rm
excursions from } \partial \bB_\sfd(\sfv,r_j) \;\;
\mbox{\rm  to } \partial \bB_\sfd(\sfv,r_{j-1})  
\nonumber \\
\; & \mbox{\rm during  the}
 \;\; [\sfs_{\sfKo}(y)] \; \mbox{\rm first 
$\sfKo$-level excursions of the block $\beta$} \}.
 \label{Tjvdef}
\end{align}
We further consider the corresponding "excess excursion counts"
\begin{align}\label{eq:eta-centered}
\whetavjz := \etavjz - \gamma_0(j) \quad \rm{and} \quad \whetavjy=\etavjy-\gamma_0(j) \,.
\end{align} 
In case $\sfv=\sfv_\bb$,  we shall use the shorthand $\Tjz$, $\etajz$ and $\whetajz$ 
for  $\Tvjz$,  $\etavjz$ and $\whetavjz$,  likewise using the notations $\Tjy=\Tvjy$,
$\etajy=\etavjy$ and $\whetajy=\whetavjy$ when $\sfv=\sfv_\bb$.
\end{definition}
We note in passing $\Tjy=T_{j}^{\sfv_\bb,\sfKo,s}$ of \eqref{Tjsdef},  at $s=\sfs_{\sfKo}(y)$
and make the following two useful observations.
\begin{remark}\label{rem:y-to-z}
In view of  \eqref{eq:nominal} and \eqref{eq:eta-centered}
the path completes within the first $\sfk_\vep(z)$ inter-city excursions
precisely $\sfs_{\sfKo} (\Ybz)$ of the $\sfKo$-level excursions to block $\beta$.
Hence,  for $\sfv \in \bF^\bb$ each $\Tvjz$ count equals $\Tvjy$ at $y=\Ybz$.
\end{remark}}
\begin{remark}\label{rem:alt-UU}
 Note that for any $z$ and $\bb$,  
\begin{equation}\label{dfn:UU-alt}
\UUz \subseteq \Big\{ \min_{\sfv \in \bF^\bb} \{ T_{\sfv,n_\vep+1}^{[z]} \} = 0  \Big\} := \whUUz \,.
\end{equation}
Thus,  $\bP(\UUz \cap \AA) \le \bP(\whUUz \cap \AA)$ for any event $\AA$ of interest.  In 
particular,  we can and shall replace $\UUz$ by $\whUUz$ throughout our proof of \eqref{eq:U-minus}.
\end{remark}
\noindent
{\bf Step 2 - barrier events}.   
\blue{Consider} for the level intervals 
\begin{equation}\label{dfn:J-intervals}
\sfJupl:=\sfJo \cup \sfJm \cup \{\nvfv\}\, , \;\; 
\sfJo:=(\sfKo,\sfKt) \cup (\sfKth,\nvf) \,,\;\;
\sfJm:=\{\sfK_i, i \le 3, \nvf \}\,,
\end{equation}
the excursion count windows  
\begin{equation}\label{dfn:windows}
\sfI_j := \left\{\!\!\!\begin{array}{cll}
& (j_{n_\vep}^{\gak}, j_{n_\vep}^{\okappa}),  & \qquad j \in \sfJo \,,\\
& \,\,(j_{n_\vep}^{\kapo},  j_{n_\vep}^{\kapp}),  & \qquad j \in \sfJm,   \\
& (-\infty,\sfKfv), & \qquad j = \nvfv \,,
\end{array}
\right.
\end{equation}
setting the outer-block double barrier events 
\begin{align}\label{eq:barrier-event-z}
\!\!\!\! \BBy \! := & \bigcap_{j \in \sfJupl}  \{ \whetajy \in \sfI_j  \} \quad {\rm and} \quad 
\blue{ \BBb \!  :=  \BBby (\Ybz)}  =  \bigcap_{j \in \sfJupl} \{ \whetajz \in \sfI_j  \} 
\end{align}
\blue{(see \Cref{rem:y-to-z})}.
Next,  considering the events $\UUy(y)$ \blue{and $\UU^{(\vep)}_\sfv(y)$}, 
of having an unhit $\vep$-ball within the block $\beta$,  \blue{or
centered within the (smaller) square of edge length $\sqrt{2} \wt{c} \,  \vep$ 
around $\sfv \in \bF^\bb$,  respectively},  during the 
$[\sfs_{\sfKo}(y)]$ first $\sfKo$-level excursions of block $\bb$,
one has \blue{from \Cref{rem:y-to-z}} that for any $\vep$ and $z$, 
\begin{equation}\label{eq:time-ident}
\UUz = \UU^{(\vep)}_{\bb}(\Ybz) = \blue{\bigcup_{\sfv \in \bF^\bb} \UUv} \quad {\rm where} \quad 
\UUv := \UU^{(\vep)}_{\sfv}(\Ybz) \,.
\end{equation}
Setting $\sfJs :=[\wnvfv,\nvp]$,  our restricted-uncovering events are then
\begin{align}\label{def:res-u-z}
\!\!\!\! \RUb \!\! :=\blue{\RUyz (\Ybz)} = \! \bigcup_{\sfv \in \bF^\beta}  \big( \UUv \cap \RRv \big),  \quad &\text{for} \qquad 
\!\!\!\! \RRv \!\!  := \blue{\RRyz (\Ybz)} = \bigcap_{j \in \sfJs} \{ \etavjz > \gamma_{-1}(j) \} \,,\\
\label{eq:barrier-event}
 \RUy \!\! := \!\!\! \bigcup_{\sfv \in \bF^\beta}  \big( \UU^{(\vep)}_\sfv(y) \cap \RRy \big)
\;\; &\text{and} \;\; 
 \RRy \!\! := \bigcap_{j \in \sfJs} \!\!\!\
\{ \etavjy > \gamma_{-1}(j) \} \,.
\end{align}
\begin{figure}[h]
\begin{tikzpicture}
    \begin{axis}[
width=18cm, height=10cm,
        axis lines = left,
xticklabels={},
yticklabels={},
scaled y ticks = false,
scaled x ticks = false,
        domain = 0:100000,
        samples = 100,
        ymin =0, ymax = 220000,
        xmin =0, xmax = 110000,
xtick = {1000,5000,13500,20000,25000,75000,80000,95000,100000},
xticklabels = {$\sfKz$,$\sfKo$,$\sfKw$,$\sfKt$,$\sfKth$,$\quad\sfKf$,$\sfKfv$,$\sfKp$,$n_\varepsilon$}
 ]
\draw[<-, thick, blue] (axis cs:95100, 3500) -- (axis cs:100000, 3500);
\draw[<-, thick, brown] (axis cs:80000, 2000) -- (axis cs:100000, 2000);
\draw[<-,thick,green] (axis cs:75000, 500) -- (axis cs:100000, 500);
\draw[dashed, gray, thick] (axis cs:25000,0) -- (axis cs:25001,155000);
\draw[dashed, gray, thick] (axis cs:1000,0) -- (axis cs:1000,198000);
\draw[dashed, gray, thick] (axis cs:5000,0) -- (axis cs:5000,192000);
\draw[dashed, gray, thick] (axis cs:20010,0) -- (axis cs:20010,165000);
\draw[dashed, gray, thick] (axis cs:75000,0) -- (axis cs:75010,55000);
\draw[dashed, gray, thick] (axis cs:80000,0) -- (axis cs:80000,45000);
\draw[dashed, gray, thick] (axis cs:94990,0) -- (axis cs:94990,20000);
\draw[dashed, gray, thick] (axis cs:13500, 0) -- (axis cs:13500,20000);
\draw[dashed, gray, thick] (axis cs:99950,0) -- (axis cs:99950,20000);
\draw[-{Stealth[scale=0.5]}, blue,thick] (axis cs:4900,190000) -- (axis cs:4900,195000);

\draw[line width=5pt, black](axis cs:20000,170000) -- (axis cs:20000, 175000);
\draw[line width=4pt, black](axis cs:25100,162000) -- (axis cs:25100, 168000);
\draw[line width=4pt, black](axis cs:75000,64000) -- (axis cs:75000, 68000);
\draw[line width=4pt, black](axis cs:80000,19500) -- (axis cs:80000, 68000);

\node[above right] at (axis cs:20000,163000) {$\sfI_{\sfKt}$};
\node[above right] at (axis cs:2400,190000){$y$};
\node[above right] at (axis cs:25100,155000) {$\sfI_{\sfKth}$};
\node[above right] at (axis cs:75000,55000) {$\sfI_{n''}$};
\node[above right] at (axis cs:100000,100) {$\UUy(y)$};

\draw[
            decorate, 
            decoration={brace, amplitude=6pt,mirror}, 
            thick, 
            black
        ] (axis cs:5000,22000) -- (axis cs:13500,22000) 
        node[midway, yshift=-15pt, font=\small] {$\WWy$}; 

\draw[
            decorate, 
            decoration={brace, amplitude=6pt,mirror}, 
            thick, 
            black
        ] (axis cs:80000,48000) -- (axis cs:100000,48000) 
        node[midway, yshift=-15pt, font=\small] {$\RUy$}; 

\draw[
            decorate, 
            decoration={brace, amplitude=6pt,mirror}, 
            thick, 
            black
        ] (axis cs:5000,48000) -- (axis cs:20000,48000) 
        node[midway, yshift=-15pt, font=\small] {$\BWlb$}; 

\draw[
            decorate, 
            decoration={brace, amplitude=6pt,mirror}, 
            thick, 
            black
        ] (axis cs:25000,48000) -- (axis cs:75000,48000) 
        node[midway, yshift=-15pt, font=\small] {$\BBKth(y)$}; 


        \addplot[red,dashed] {200000-2*x};
 \addplot[black,dashed,domain=0:5000] {200000-2*x+40*x^(0.48)};
\addplot[black,dashed,domain=0:5000] {200000-2*x+20*x^(0.7)};
 \addplot[black,thick,domain=5000:20000,name path=up] {200000-2*x+40*x^(0.48)};
 \addplot[black,thick,domain=5000:20000,name path=down] {200000-2*x+20*x^(0.7)};
\addplot[gray!30] fill between[
            of=up and down,
            soft clip={domain=5000:20000} 
        ];
\addplot[black,dashed,domain=20000:25000] {200000-2*x+40*x^(0.48)};
\addplot[black,dashed,domain=20000:25000] {200000-2*x+20*x^(0.7)};

\addplot[black,thick,domain=25000:50050,name path=up1] {200000-2*x+40*x^(0.48)};
 \addplot[black,thick,domain=25000:50050,name path=down1] {200000-2*x+20*x^(0.7)};
\addplot[gray!30] fill between[
            of=up1 and down1,
            soft clip={domain=25000:50050} 
        ];

\addplot[black,thick,domain=50000:75000,name path=up2] {200000-2*x+40*(100000-x)^(0.48)};
 \addplot[black,thick,domain=50000:75000,name path=down2] {200000-2*x+20*(100000-x)^(0.7)};
\addplot[gray!30] fill between[
            of=up2 and down2,
            soft clip={domain=50000:75000} 
        ];

\addplot[black,dashed,domain=75000:80000] {200000-2*x+40*(100000-x)^(0.48)};
\addplot[black,dashed,domain=75000:80000] {200000-2*x+20*(100000-x)^(0.7)};

\addplot[black,dashed,domain=80000:95000,name path=up3] {100000-5000*sin((x-80000)/40)};
 \addplot[black,thick,domain=80000:95000,name path=down3] {200000-2*x-20*(100000-x)^(0.7)};
\addplot[gray!30] fill between[
            of=up3 and down3,
            soft clip={domain=80000:95000} 
        ];

\addplot[black,dashed,domain=95000:100000] {200000-2*x-20*(100000-x)^(0.7)};

    \end{axis}
\end{tikzpicture}
\caption{\footnotesize{Scales and restrictions.  Given excess $y$ at scale $\sfKo$,  we impose the angular restriction 
event $\WWy$ at scales between $\sfKo$ and $\sfK_a$ and various barrier restrictions 
(the event $\BWlb$ corresponds to  $\WWy$ and the shaded area between 
$\sfKo$ and $\sfKt$ with $\BBKth(y)$ for the shaded area between
$\sfKth$ and $n''=\nvf$,  having stronger restrictions $\sfI$ at boundary scales $\sfKt$,  $\sfKth$,  $n''$ and
the event $\RUy$ of restricted block uncovering).  The scales prior to $\sfKz$ 
correspond to what happens outside the city in which block $\beta$ lies.}}
\label{fig-scales}
\end{figure}
\noindent{\bf Step 3 - angle restrictions.} For $j \ge \sfKo$, $q \ge 1$,  we denote by $\thbj \in \Ipi 
$,  the angle-difference
between the start point on $\partial \bB_\sfd(\sfv_\bb,r_j)$ of the $q$-th $j$-level excursion of block $\bb$
and its end-point on $\partial \bB_\sfd(\sfv_\bb,r_{j-1})$
and for $y,z \in \bR$,  consider the corresponding empirical measures 
\begin{equation}\label{dfn:emp-ang}
\Ljy := \frac{1}{\Tjy} \sum_{q=1}^{\Tjy} 
 \delta_{\thbj} \qquad {\rm and} \qquad \Ljz := \frac{1}{\Tjz} \sum_{q=1}^{\Tjz} 
 \delta_{\thbj} \,, 
\end{equation}
during the $[\sfs_{\sfKo}(y)]$ first $\sfKo$-level excursions of $\bb$ 
and the first $\sfk_\vep(z)$ inter-city excursions,  respectively.  
Absent further conditioning,  for any $\beta$ and $j \in [\sfKo, \sfKt]$,  the variables $\thbj$ 
are \abbr{iid} of marginal law $\bar \mu$ matching that of the angle-difference 
between a Brownian path starting point on $\partial \bB_\sfd(\npp,r_1)$ and its first
hitting point of $\partial \bB_\sfd(\npp,r_0)$.
In particular, $\bar \mu$ has a continuous, positive density with respect to the uniform 
measure on $\Ipi$,  given by the Poisson kernel,
\begin{equation}\label{dfn:Poiss}
\blue{2 \pi}\,  \frac{d\bar \mu}{d\theta} = \frac{r_{0}^2-r_{1}^2}{|r_{0} e^{i\theta}-r_{1}|^2} = 
\frac{1}{\sqrt{1+\bar c^2} -\bar c \cos \theta}\,, 
\qquad \bar c := \frac{2e}{e^2-1} \blue{<} 1.
\end{equation}
Next,  recall the \blue{Kantorovich-Wasserstein $W_\infty$-distance} between probability measures 
$\mu_i \in M_1(\Ipi)$, 
\begin{align}
\sfd_\infty(\mu_1,\mu_2) &:= \sup_{t \in [0,1]} \{ |\mu_1^{-1}(t) - \mu_2^{-1} (t) | \} \,,
\label{dfn:dsup}
\end{align}
where \blue{$\mu^{-1} (t):=\inf\{\theta: \mu([-\pi,\theta]) \ge t \}$
(see \cite[(2.3)]{BL19})},  denoting by 
$\bB^\infty_\mu(c)$ the $\sfd_\infty$ ball of radius $c$ and center $\mu \in M_1(\Ipi)$.  
Setting $\boetaj := \frac{\cst}{\sfm_\vep} \log j$ for some large constant $\cst > 1$ and 
\begin{equation}\label{dfn:angle-width}
\boetasj := \left\{\!\!\!\begin{array}{cll}
& 2 \boetaj,  & \qquad j = \sfKo \,,\\
& 2 \cst \boetaj,  & \qquad j = \sfKo+3,   \\
& 2 \cst^2 \boetaj, & \qquad \text{otherwise} \,,
\end{array}
\right.
\end{equation}
we impose on the angle-differences 
empirical measures of \eqref{dfn:emp-ang} the controls  
\begin{align}\label{eq:was-event}
\WWy &:=
\! \bigcap_{j \in \sfJwp} \{ \Ljy \in \bB^\infty_{\bar \mu} (\boetasj) \}  \,,
 \quad {\rm with} \quad 
\sfJwp := \{\sfKo\} \cup \sfJw \quad {\rm and} \quad \sfJw :=[\sfKo+3,\sfKw] \,.
\end{align}
\blue{In view of \Cref{rem:y-to-z}} we have that
\begin{align}
\label{eq:was-event-z}
\WWb &:=\blue{\WW^{(\vep)}_{\bb}(\Ybz) =} 
\! \bigcap_{j \in \sfJwp} \! \{ \Ljz \in \bB^\infty_{\bar \mu} (\boetasj) \}  
\end{align}
\blue{and consequently also} the identities 
\begin{equation}
\label{eq:good-cover-y}
\Big\{\UUp \Big\}_\bb = \Big\{ \UUpy(\Ybz) \Big\}_\bb 
\quad {\rm for} \quad \UUpy(y) := \RUy \cap \BBy \cap \WWy 
 \end{equation}
(compare with \eqref{eq-good}).

\subsection{First moment computation for events $\UUp$}\label{subsec-exact}
Setting \blue{in terms of 
$\UUp$ and $\Lambda^{(+,\vep)}_{u,z}$ of  \eqref{def:barQub}},  the conditional expectations
\begin{align}\label{dfn:wt-Q}
\wt{Q}^{(1,\vep)}_{u,z} &:= \bE \big[ \,\Lambda^{(+,\vep)}_{u,z} \,|\, \FFuz  \big] = 
\sum_{\beta \in u} \bP \big[ \,\UUp \,| \, \FFuz  \, \big] \,,\\
\label{dfn:wt-Q2}
\wt Q^{(2,\vep)}_{u,z} &:= \bE\big[ \Lambda^{(+,\vep)}_{u,z} (\Lambda^{(+,\vep)}_{u,z}-1) 
| \FFuz \big] = \sum_{\bb' \ne \bb \in u}  \bP \big[ \,\UUpb \cap \UUp \,| \, \FFuz \, \big] \,,
\end{align}
\blue{our next result about the limiting first moment is key to proving \eqref{eq:ct-exc-plus}}.
\begin{theorem}\label{thm:end-game}
Fix $z$,  $\srs=r_0$,  $\dlc$ and $\sfKp$.  For some $\overline X$ 
whose law depends also on $\sfKz$,  $\sfKfv$,  but not on $\vep$ or $z$,
\begin{align}
 \bE \Big[   \prod_{u=1}^{\sfM}  (1-  \widetilde{Q}^{(1,\vep)}_{u,z})_+ \Big] 
= 
\bE \Big[  \exp\big(-\sum_{u=1}^{\sfM} \wt{Q}^{(1,\vep)}_{u,z}\big)  \Big] + o_{\sfKfv,\sfKz,\vep} (1) 
= 
\bE \Big[ \exp( - \overline X e^{-\sfco z} ) \Big] + o_{\sfKfv,\sfKz,\vep}(1)  \,.
\label{eq:end-game}
\end{align}
 \end{theorem}
\blue{The proof of \Cref{thm:end-game},  given at the end of this section is
based on four propositions whose proofs are deferred to later sections.  
For stating those propositions,  we first} introduce a few convenient notations.
\begin{definition}\label{dfn:ggi-barsi}
For $i=1,2,3$ let $\sfIs_i$ denote
the subset of $y \in \sfI_{\sfK_i}$ for which $\sfs_{\sfK_i}(y) \in \bN$ and
 \begin{equation}\label{dfn:ggi}
\ggi(y):={\bf 1}_{\sfIs_i} (y) y e^{-\sfco y},   \qquad i=1,2,3.
\end{equation}
\blue{We similarly} denote by $\sfIs_i$,  $i=4,5$,  those $y \in  \sfI_{n^{(i)}_\vep}$ for which
\begin{equation}\label{eq:nominalno-vep}
 \bar \sfs_{\sfK_i} (y)  :=  \frac{1}{2} (\sfco\,  \sfK_i  + y)^2 \in \bN \,.
\end{equation}
We denote by $\underline{\sfx} := \{ (\sfx_{q1},\sfx_{q2}), q=1,\ldots,[\sfs_{\sfKo} (y)] \}$,
the ordered collection of start and end points of the first $\sfKo$-level excursions 
of block $\beta$ \blue{and isolate}
the `large scale' constraints part of \eqref{eq:good-cover-y},  \blue{see \Cref{fig-scales},}
\begin{align}
\label{eq:BW-K12}
\BWlb & :=  \WWy 
\bigcap_{j \in [\sfKo,\sfKt]} \big\{ \whetajy \in \sfI_j  \big\}  \,.
\end{align}
\end{definition}
The law of $\BWlb$ 
given $\LKoy$ determines $\bP(\UUpy(y)| \underline{\sfx})$,
as \blue{stated next and proved} in \Cref{sec:angle-diff}.
\begin{proposition}\label{prop:angle-diff}
For 
some $0<\uu \alpha(\sfKp)<\bar \alpha(\sfKp)<\infty$
and $o_{\sfKfv,\vep} (1) \to 0$ (as $\vep \to 0$,  then $\sfKfv \to \infty$),  uniformly in $\beta$, $\underline{\sfx}$ 
and $y$,  
\begin{equation}\label{eq:good-cover}
e^{2 (\nvfv-\sfKth)} \, \bP(\UUpy(y) | \underline{\sfx}) = \big( \alpha(\sfKfv,\sfKp)+o_{\sfKfv,\vep} (1) \big)
\, \bE \big[ \gth(\wh{\eta}^{\,(y)}_{\bb,\sfKth});  \BWlb  \,| \LKoy\big]  \,,
\end{equation}
where $\alpha(\sfKfv,\sfKp) \in [\uu \alpha(\sfKp), \bar \alpha(\sfKp)]$ for all $\sfKp$ large enough.
\end{proposition}
\begin{remark}\label{rmk:no-vep-inJs}
We assume hereafter that $\vep$ is small enough 
so for all $j \in [\nvf,n_\vep]$ replacing $\sfrho_\vep$ by its limit does not change  
$\{ \sfs_{j} (y) : y \in \sfI_{j} \} \cap \bZ$,  nor which integers exceed $\frac{1}{2} \gamma_{-1}(j)^2$.
Since $\sfKp \gg 1$ we also freely re-adjust \eqref{eq:norm-exc} and \eqref{eq:eta-centered} to 
have $\gamma_{b}(n_\vep-\ell):=\sfco \ell + b \ell^{\gak}$ whenever $n_\vep - \ell \in \sfJs$ 
(by slightly changing $\gak>0$).
Similarly,  
we may \abbr{wlog} change to non-strict inequalities in \eqref{eq:barrier-event-z} and \eqref{def:res-u-z}
(e.g.,  when referring to \cite{BRZ,BRZ2}).
\end{remark}

Next,  for $\sfI_{\sfKo}$ of \eqref{dfn:windows} denote by 
\begin{equation}\label{dfn:I-star-ext}
\Istar  \subset \sfI_{\sfKo} \times \bB^\infty_{\bar \mu}(\boetasKo)
\end{equation}
the collection of $y \in \sfIs_1$ and all empirical measures of 
$\sfs_{\sfKo}(y) \in \bN$ angle-differences in $\Ipi$ which are 
in $\bB^\infty_{\bar \mu}(\boetasKo)$.  As $\widehat{\eta}^{(y)}_{\bb,\sfKo}=y$,  the barrier events 
$\BBy$  \emph{are empty whenever $y \notin \sfIs_1$}.  Likewise, 
$\WWy$ are empty whenever $\LKoy \notin \bB^\infty_{\bar \mu}(\boetasKo)$,  so 
$\Istar$ is the collection of pairs $(y,\LKoy)$ for which  the event  $\BWlb$ may occur.  
Setting for probability measures $\nu$ on $\Ipi$, 
\begin{equation}\label{dfn:hat-L-alpha}
\Vb (\nu) :=  \frac{\bar c}{2} \sfm_\vep \big[
\langle \cos(\cdot), \bar \mu \rangle
-  \langle \cos (\cdot), \nu \rangle
\big] \,,
\end{equation}
we establish in 
\Cref{sec:good-cover} the following proxy for the right hand side of \eqref{eq:good-cover}.
\begin{proposition}\label{prop:proxy-good-cover}
Let
\begin{equation}\label{eq:g-def}
\gstb (y,v) := y e^{-\sfco y} e^{- \sfco v} 
\,,
\qquad \sfb := \frac{1}{2} (\sfco \bar{c})^2 {\rm var}_{\bar\mu}(\theta) \,,
\end{equation}
for $\bar{c}$ and $\bar\mu$ of \eqref{dfn:Poiss}.  Then, 
for large $\cst$ and some $o_{\sfKo,\vep} (1) \to 0$ (as $\vep \to 0$ followed by 
$\sfKo \to \infty$),  
\begin{align}\label{eq:K1-to-K3}
[r_{\sfKo}/r_{\sfKth}]^2\,  
\, \bE \big[ \gth(\wh{\eta}^{\,(y)}_{\bb,\sfKth});  \BWlb  \,| \LKoy\big]  
&=  (1+o_{\sfKo,\vep} (1)) e^{-\sfb} \gstb (y,\Vb(\LKoy))\, ,
\end{align}
uniformly in $\beta$ and uniformly 
with respect to all $(y,\LKoy) \in \Istar$.
\end{proposition}

\blue{Combining Propositions \ref{prop:angle-diff} and \ref{prop:proxy-good-cover} we arrive at the following explicit
formula for $\bP(\UUpy(y) | \underline{\sfx})$.}
\begin{corollary}\label{cor:UUpy}
For some non-random 
$o_{\sfKfv,\sfKz,\vep}(1) \to 0$ as $\vep \to 0$ then $\sfKz \to \infty$ followed by $\sfKfv \to \infty$,
uniformly in $\bb$ and $\underline{\sfx}$ such that $(y,\LKoy) \in \Istar$
\begin{equation}\label{UUpy-asymp}
\bP(\UUpy(y) | \underline{\sfx}) =  \big(1+o_{\sfKfv,\sfKz,\vep} (1) \big) 
\big[\vep e^{\sfKfv}/r_{\sfKo}\big]^2
\alpha(\sfKfv,\sfKp) e^{-\sfb} \gstb (y,\Vb(\LKoy)) \,.
\end{equation}
\end{corollary}

\subsubsection*{Limiting variables and proof of  \Cref{thm:end-game}}
\label{subsec-thm2.8proof}
\blue{Before proceeding with the proof,  we introduce certain} $\sigma$-algebras 
which will be often used in the paper and explain the extra freedom within our barrier events.
\begin{definition}\label{def:Fvjs}
Let $\GG^\tau$ denote the $\sigma$-algebra generated by the Brownian motion up to the stopping time $\tau$,  
denoting by  $\GG^{\sfv, \tau}_j$ the $\sigma$-algebra $\GG^\tau$
augmented by all the Brownian excursions from $\partial \bB_\sfd (\sfv, r_{j-1})$ 
to $\partial \bB_\sfd(\sfv,r_j)$.  In particular,  $\GG^{\sfv, \tau}_j$ contains the start and end points of all 
excursions from $\partial \bB_{\sfd}(\sfv,r_j)$ to  $\partial\bB_{\sfd}(\sfv,r_{j-1})$.
We further use $\GG^{\sfv}_j$ for $\GG^{\sfv,\tau}_j$ with $\tau$ the Brownian entrance time of 
$\bB_\sfd (\sfv,r_{j-1})$. 
\end{definition}

\begin{remark}\label{rmk:z-a-vep}
The value of each $Q^{(+,\vep)}_{u,z}$,  $\wt{Q}^{(1,\vep)}_{u,z}$ or $\wt{Q}^{(2,\vep)}_{u,z}$ is 
unaffected if one replaces in their definition $\FFuz$ with the  $\sigma$-algebra $\HHk$ 
or even with the larger $\sigma$-algebra $\HHu$ \blue{(see \eqref{eq:FF-to-HH})}.  
Similarly,  the joint law of $\{ (\Ybz,\LKoz) \}_{\bb \in u}$ 
conditional on $\FFuz$ coincides with their joint law conditional on $\HHu$.
\end{remark}

Note that $|\Vb(\nu)| \le \rho$ for $\Vb(\cdot)$ of \eqref{dfn:hat-L-alpha} and 
any $\nu \in \bB^\infty_{\bar \mu} (2 \rho/\sfm_\vep)$
(this follows from the definition \eqref{dfn:dsup},  see \eqref{dfn:dw-orig},  \eqref{dfn:dw} for the justification).
Thus,  
\begin{equation}\label{eq:Ball-to-angle}
(y,\LKoy) \in \Istar \quad \Longrightarrow \quad \Vb(\LKoy) \in [-\cst \log \sfKo, \cst \log \sfKo] := \sfIw \,.
\end{equation}
\blue{Recall that \eqref{UUpy-asymp} applies only for $(y,\LKoy) \in \Istar$  (as its \abbr{lhs} is otherwise zero).
Further,  for  $y=\Ybz$ and $\uu{\sfx}$ the vector of start and end points of the first $\TKoz$ 
excursions of block $\bb$ at level  $\sfKo$,  the conditional probabilities on the \abbr{lhs} 
of \eqref{UUpy-asymp} are precisely $\bP(\UUp\,|\, \GG^{\sfv_\bb}_{\sfKo})$,  while as in 
 \Cref{rem:y-to-z},  on its right hand side we  replace $\LKoy$ by the $\GG^{\sfv_\bb}_{\sfKo}$-measurable $\LKoz$
 to arrive at 
\begin{equation}\label{UUpz-asymp}
\bP(\UUp\,|\, \GG^{\sfv_\bb}_{\sfKo}) =  \big(1+o_{\sfKfv,\sfKz,\vep} (1) \big) 
\Big[\frac{\vep e^{\sfKfv}}{r_{\sfKo}}\Big]^2
\alpha(\sfKfv,\sfKp) e^{-\sfb} \gstb (\Ybz,\Vb(\LKoz)) {\bf 1}_{(\Ybz,\LKoz) \in \Istar} \,.
\end{equation} 
Next,  setting  
\begin{equation}\label{def:bar-bKz}
\bKz := (1-6 \dlc)^2 r_0^{-2} e^{2 (\sfKo-\sfKz)} e^{-\sfb} \quad {\rm and} \quad \Vbz :=\Vb (\LKoz) \,,
\end{equation}
as $\FFuz \subset \HHk \subset \GG^{\sfv_\beta}_{\sfKo}$,
summing over $\bb \in u$ the expectation of \eqref{UUpz-asymp}
given $\HHk$ it follows from \eqref{dfn:wt-Q},  \eqref{def:bar-bKz} and \Cref{rmk:z-a-vep},  that 
\begin{align}\label{eq:wtQ1-bd}
\wt{Q}^{(1,\vep)}_{u,z} &=(1+o_{\sfKfv,\sfKz,\vep}(1)) \alpha (\sfKfv,\sfKp) \bKz
\avsum_{\bb  \in u} \bE \Big[ \gstb \big( \Ybz, \Vbz \big) ; \Ybz \in \sfI_{\sfKo},
\LKoz \in \bB^\infty_{\bar \mu} (\boetasKo) \,\big|\,  \HHk  \Big]
\end{align}}
where $\avvsum_{\bb \in u}$ denotes the average over the $[(1-6\dlc) e^{-\sfKz}/(\vep e^{\sfKfv})]^2$
blocks in city $u$.
\begin{remark}\label{rmk:free-kappas}
For convenience we have set in \eqref{dfn:windows}
the same parameters $\kapo$ and $\kapp$ for all
our narrow windows (namely,  throughout $\sfJm$),  but our arguments 
allow for different choices at each such window (provided they all satisfy \eqref{dfn:aux-par}).
For example,  we ease our task of proving Propositions \ref{prop:angle-diff} and \ref{prop:proxy-good-cover} 
by reducing $\kapo$ and increasing $\kapp$ at $j=\sfKth$ (thereby making
$\sfI_{\sfKth}$ much wider than $\sfI_{\sfKt}$),  and apply similar devices also when proving 
Propositions \ref{lem:good-cover-end} and \ref{prop:2nd-moment}.
\end{remark}
\begin{remark}\label{rmk:first-mom}
\blue{\Cref{prop-apriori},  \Cref{cor:UUpy} and their proofs
apply} for any parameters satisfying \eqref{dfn:aux-par} with $\cst$ large.  Also,  
while we have set $\sfKth = \sfKt + 5 \log n_\vep$,   $\sfKt = \sfKo + 2 \sfKw$ and 
$\sfKw = (21 \log n_\vep)^{1/\gak}$,  
our proof of \eqref{UUpy-asymp} works verbatim as soon as for some $C<\infty$,
\begin{equation}\label{eq:flex-Ki}
C (\log n_\vep)^{1/\gak} \ge \sfKth \ge  \sfKt + 2 \log (5 n_\vep^2),  \quad 
  \sfKt \ge \frac{1}{C} (\log n_\vep)^{1/\gak},  \quad \sfKt \ge 2 \sfKw \,.
\end{equation}
The same applies for \eqref{eq:wtQ1-bd}.
\end{remark}
\blue{The proof of the following approximation of the right hand side of \eqref{eq:wtQ1-bd} is}
given in  \Cref{sec:good-cover-end}. 
\begin{proposition}\label{lem:good-cover-end}
\blue{For $\bKz$ of \eqref{def:bar-bKz} and $\alpha(\sfKfv,\sfKp)$ of \Cref{prop:angle-diff} consider}
the random variables 
\begin{align}
\sfX^{(\vep)}_{u,z}  &:=  \alpha (\sfKfv,\sfKp) \bKz 
\,  \avsum_{\bb \in u} \bE \Big[ \gst \big( \Ybz -z,\Vbz \big)  
\,\big|\,  \HHk\Big]\,,
\label{def:X-eps}
\end{align}
where \blue{
 \begin{align}\label{eq:g-def-alt}
\gst (y,v) &:=  \gstb (y,v) \big( 1- \inf\{  |y-y'|+|v-v'| : (y',v') \in \sfI_{\sfKo} \times \sfIw \} \big)_+ \,.
\end{align}
Then,  there exist random variables $0 \le \bar Q^{(1,\vep)}_{u,z} \le \sfX^{(\vep)}_{u,z} e^{-\sfco z}$ such that 
\begin{equation}\label{eq:wX-small}
\varlimsup_{\sfKz,\vep} \bP \Big[ \sum_{u=1}^{\sfM}  \big(\sfX^{(\vep)}_{u,z} e^{-\sfco z} 
- \bar Q^{(1,\vep)}_{u,z} \big)  > \delta \Big] = 0 \qquad \forall \delta >0,
\end{equation}
and} for $\wt{Q}_{u,z}^{(1,\vep)}$ of \eqref{dfn:wt-Q},  we have that uniformly in $u$,
\begin{align}\label{eq:eqv-wt-Q}
\wt{Q}_{u,z}^{(1,\vep)} = \big(1+o_{\sfKfv,\sfKz,\vep}(1)\big) \bar Q^{(1,\vep)}_{u,z} 
\,.
\end{align}
\end{proposition}
\begin{remark}\label{rmk:Xu-cnv-zero}
Clearly,  $\| \gst \|_\infty \le e^{-c\,  \sfKo^{\kapo}}$ for some universal $c>0$ \blue{(as 
$y \ge \sfKo^{\kapo}-1$ and $v \ge -\cst \log \sfKo-1$} throughout the support of $\gst(\cdot)$).
Hence,  our choice of 
$\sfKo-\sfKz = 5 \log \sfKz$ implies that for any $z$ fixed,  
\begin{equation}\label{eq:Xu-small}
\lim_{\sfKz \to \infty} \sup_{\vep, u} \{ \bar Q^{(1,\vep)}_{u,z}
 \} = 0  \,.
\end{equation}
\end{remark}
In \Cref{sec:end-game} we prove the convergence in law,  as $\vep \to 0$, 
of $\{\sfX^{(\vep)}_{u,z},u \in [1,\sfM]\}$ of \eqref{def:X-eps},  which will allow us in turn 
to take full advantage of \eqref{eq:wX-small}-\eqref{eq:Xu-small}.
\begin{proposition}\label{prop:Gauss-field-continuity}
Fix $z$,  
$\dlc$, $\srs=r_0$  and $\sfKz$ large.  As $\vep \to 0$
the random vectors $\{ \sfX^{(\vep)}_{u,z},  u \in [1,\sfM] \}$
of \eqref{def:X-eps} converge in law to some $\{ \sfX_{u},  u \in [1,\sfM] \}$ which is
independent of $z$.
\end{proposition}

\begin{proof}[Proof of \Cref{thm:end-game}] Fixing $z$,  $\dlc$,  $\srs$ and $\sfKp$,  consider the variable  
\begin{equation}\label{dfn:Xstar-K}
\overline X :=  \sum_{u=1}^{\sfM}  \sfX_u 
\end{equation}
whose law depends only on $\dlc$,  $\srs$, 
$\sfKz$, $\sfKfv$ and $\sfKp$.
By \Cref{prop:Gauss-field-continuity},  we have that \blue{as $\vep \to 0$, 
\begin{align*}
\sum_{u=1}^{\sfM} X^{(\vep)}_{u,z} e^{-\sfco z} 
 \dconv  \overline X e^{-\sfco z} \,.
\end{align*}
In particular,  by \eqref{eq:wX-small},  taking $\vep \to 0$ then $\sfKz \to \infty$, } 
\begin{align}
\varlimsup_{\sfKfv,\sfKz,\vep}\,\, \Big|
\bE \big[  \exp\big(-\sum_{u=1}^{\sfM} \bar Q^{(1,\vep)}_{u,z}\big)  \big] 
- 
\bE \big[ \exp( - \overline X e^{-\sfco z} ) \big]  \Big| = 0 \,.
\label{eq:end-game1}
\end{align}
With $\bar Q^{(1,\vep)}_{u,z} \ge 0$,  by the uniform in $u$ and $\vep$ convergence of  \eqref{eq:Xu-small},   
\begin{equation}\label{eq:end-game2}
\varlimsup_{\sfKfv,\sfKz,\vep}\,\,  \bE \Big[  \exp\big(-\sum_{u=1}^{\sfM} \bar Q^{(1,\vep)}_{u,z}\big)  \Big] -
 \bE \Big[   \prod_{u=1}^{\sfM}  (1-  \bar Q^{(1,\vep)}_{u,z})_+ \Big] = 0 \,.
\end{equation}
Thanks to the uniform in $u$ multiplicative bound of \eqref{eq:eqv-wt-Q} \blue{we deduce that
\eqref{eq:end-game1} and \eqref{eq:end-game2} hold with}
$\bar Q^{(1,\vep)}_{u,z}$ replaced by  $\wt{Q}_{u,z}^{(1,\vep)}$,
thereby concluding our proof of \eqref{eq:end-game}.
\end{proof}
\subsection{Second moment,  tightness and the proof of \Cref{theo-covertime}}
\label{subsec-pftheo-covertime}
In \Cref{sec:2nd-mom-tight} we show that the events $\{ \UUp\}_{\beta \in u}$ 
are sufficiently uncorrelated to yield the following second moment control.  
\begin{proposition}\label{prop:2nd-moment}
Fix $\srs$,  $\dlc$,  $\sfKp$ and parameters as in \eqref{dfn:aux-par}. 
There exists $\dst=\dst(\sfKz,\sfKfv) \ge 0$ such that 
\begin{equation}\label{eq:C-star-good}   
\varlimsup_{\sfKfv,\sfKz} \dst(\sfKz,\sfKfv) = 0 
\end{equation}
and for any fixed $z$,
\begin{equation}
\varlimsup_{\sfKfv,\sfKz,\vep} \bP \Big( \sum_{u=1}^{\sfM} \wt Q^{(2,\vep)}_{u,z} 
\ge \dst + \dst \sum_{u=1}^{\sfM} \wt Q^{(1,\vep)}_{u,z}   \Big) = 0 \,.
\label{eq:sec-mom-bd}
\end{equation}
\end{proposition}
\blue{Utilizing this,  we also} prove in  \Cref{sec:tight} (see also \Cref{subsec:tools}),
the $\bT^2$-analog of the main result of \cite{BRZ2}.
\begin{proposition}\label{prop:tight}
The collection $\Big\{\sqrt{\CCT} -\osfm_\vep, \; \vep>0\Big\}$ is 
tight \blue{and \eqref{ct2.25-new} holds.}
\end{proposition} 
Now,  by definition $\wt{Q}^{(1,\vep)}_{u,z} \ge Q^{(+,\vep)}_{u,z} \ge 
( \wt Q^{(1,\vep)}_{u,z} -\wt Q^{(2,\vep)}_{u,z})_+$ 
(see \eqref{def:barQub},  \eqref{dfn:wt-Q} and \eqref{dfn:wt-Q2}),  so 
for any $u$,  
\begin{equation}\label{eq:CS-u}
(1-\wt Q^{(1,\vep)}_{u.z})_+ \le 1-Q_{u,z}^{(+,\vep)} \le 
e^{-Q^{(+,\vep)}_{u,z}} \le
e^{-(\wt Q^{(1,\vep)}_{u,z}-\wt Q^{(2,\vep)}_{u.z} )_+} 
\end{equation}
Considering the expectation of the product over $u$ in \eqref{eq:CS-u},  then 
utilizing \eqref{eq:sec-mom-bd} and bounded convergence,  yields that 
\[
\bE \Big[   \prod_{u=1}^{\sfM}  (1-  \wt{Q}^{(1,\vep)}_{u,z})_+ \Big] 
\le 
\bE \Big[   \prod_{u=1}^{\sfM}  \big(1- Q^{(+,\vep)}_{u,z} \big) \Big] 
\le 
e^{\dst} \,  \bE \Big[  \exp\big(-(1-\dst) \sum_{u=1}^{\sfM} \wt{Q}^{(1,\vep)}_{u,z}\big)  \Big] + o_{\sfKfv,\sfKz,\vep} (1) \,. 
\]
\blue{We combine these bounds with \eqref{eq:end-game} of \Cref{thm:end-game},  
both at $z+2\ello$},  and with $\dst(\sfKz,\sfKfv) \to 0$,  arrive at
\begin{equation}\label{eq:eff-ct-exc}
\varlimsup_{\sfKfv,\sfKz,\vep} \,\Big| 
\bE \big[   \prod_{u=1}^{\sfM}  \big(1- Q^{(+,\vep)}_{u,z+2\ello} \big) \big] 
-
\bE \big[ \exp( - \blue{\overline X_\star} e^{-\sfco z} ) \big]  \Big| = 0 \,,
\end{equation}
\blue{where $\overline{X}_\star:=\overline{X} e^{-4 \ello}$ for $\overline X$ of \eqref{dfn:Xstar-K}.
Taking now the limit in $\sfKp$ then $\srs=e^{- \ello}$ and finally $\dlc$,  we deduce upon 
combining \eqref{eq:eff-ct-exc} with \eqref{ct2.25-new}, 
and \eqref{eq:CC-Quz-plus} at $z+2\ello$,  that 
\[
\varlimsup_{\dlc,\srs,\sfKp,\sfKfv,\sfKz,\vep} \,\Big| \bP \big( \CCT \le (\osfm_\vep + z)^2  \big) 
-
\bE \big[ \exp( - \blue{\overline X_\star} e^{-\sfco z} ) \big]  \Big| = 0 \,.
\]
Now,  as $\CCT$ has nothing to do with $\dlc,\srs,\sfKp,\sfKfv,\sfKz$ 
on which the law of $\overline X_\star$ may depend,  both terms 
above must converge for any fixed $z \in \bR$,  with}
\begin{align}
\lim_{\vep \to 0} \bP \big( \sqrt{\CCT}  \le  \osfm_\vep + z  \big)
= {\mathop{\underline{\overline{\lim }}}_{\dlc,\srs,\sfKp,\sfKfv,\sfKz,\vep}}
\bE \(\exp\{-\overline X_\star  \mathrm{e}^{-\sfco z}\} \)
\label{ct2.25b-new}
\end{align}
and the limit on the right hand side taken in the order as above.  
We can now provide the 
\begin{proof}[Proof of \Cref{theo-covertime}]
In view of \Cref{prop:tight},  the right hand side of \eqref{ct2.25b-new} coincides with the \abbr{cdf} 
of a finite-valued random variable $Y_\star$.  \blue{By Helly's selection theorem} 
the laws of $\overline X_\star$ converge vaguely along some subsequence 
(of $\sfKz$, $\sfKfv$, $\sfKp$,$\srs$,
$\dlc$), 
\blue{hence} such convergence must also hold weakly \blue{and} 
to some $(0,\infty)$-valued limit in law $X_\star$.
In particular, the Laplace transforms of $\overline X_\star$ 
converge pointwise along such subsequence to the Laplace transform of $X_\star$. In view
of \eqref{ct2.25b-new} the latter does not depend on the chosen subsequence, thereby 
establishing \eqref{ct2.25b}.
\end{proof}

\section{Preliminaries: barriers,  off-centering,  quantiles and coupling}\label{sec:prelim}

This section contains several preliminary results.  \Cref{subsec:bend} is devoted to the statement and calculus proof of  
\Cref{lem-bend1}, which provides a 
 (deterministic) domination of the curve
$\gamma_{b}(j)$, $j\in J:=[s,t]$ by curves for which the sharp barrier estimates from \cite{BRZ,DRZ} apply. 
In \Cref{subsec:off-center} \blue{we state and prove}
\Cref{lem-off-center-applies},  \blue{comparing for small $\sfd(\sfu,\sfv)/R$} 
the joint law of excursion counts for concentric circles
$\partial \bB_\sfd(\sfv,r e^{-j})$,  \blue{$j \ge 0$,  during the first $m$ excursions 
from $\partial \bB_{\sfd}(\sfu,R)$ to $\partial \bB_\sfd(\sfu,e R)$,  possibly conditional 
on counts of excursions around $\sfu$ at scales $R_i \ge e R$,
with the $m$-th product law for independent,  unconditional paths at $\sfu=\sfv$
(which depends only on $r/R$ and $m$).}
These estimates 
on the off-centering effect allow us 
to control in \Cref{lem-xyN} certain hitting probabilities 
of two balls by multiple excursions and to
bound in  \Cref{sec:g4-bds} the effective contribution 
to the 
uncovering probability from the part of the path within the block.  Finally, 
\Cref{subsec:quant}  collects various 
facts about quantiles of empirical measures such as those of angle-difference that 
appear in the events $\WWy$,  while \Cref{subsec:cond-exc}
provides a coupling for \abbr{iid} variables conditional on their sum, 
which is key for relating conditional and unconditional Brownian excursions, 
as in \Cref{prop:coupling}.

\subsection{Curved versus relatively curved lines}\label{subsec:bend}
We begin with the calculus lemma alluded to above.
\bl\label{lem-bend1} For $0 \le s < t \le n$ and $\gak \in (0,\frac{2}{3}]$,  
\begin{equation}
s^{ \gak}_{n}+\(t_n^{\gak}-s^{ \gak}_{n}\)\frac{l-s}{t-s}+\(l-s\)_{t-s}^{\gak} 
\ge l_{n}^{\gak},\qquad 
\forall l \in [s, t].\label{bent.1}
\ee
Further,  if $t \le n/2$ or $s \ge n/2$,  then the inequality \eqref{bent.1} holds for all $\gak \in (0,1]$.
\el
Before bringing the proof of \Cref{lem-bend1},  we state a useful corollary \blue{for
$\gak < \frac{2}{3}$,  where we multiplied \eqref{bent.1} at $n=n_\vep$ by $b$ and 
added on both sides the straight line $\gamma_0(l)$ of \eqref{eq:norm-exc}}.
\begin{corollary}\label{cor:bend1}
Recall the curved lines $\gamma_{b}(\cdot)$ of \eqref{eq:norm-exc}.
For any $J=[s,t]\subset [0,n_\vep]$,  introduce the curves
\begin{equation}\label{dfn:gamma-J}
\gamma_{b}^{J} (l \blue{+1}) := \gamma_{b} (s \blue{+1}) +
(s-l) \sfrho_\vep^{J,b} + b (l-s)_{t-s}^\gak \,, \quad 
\sfrho_\vep^{J,b} := \sfrho_\vep + b \frac{s_{n_\vep}^\gak-t_{n_\vep}^\gak}{t-s} \,.
\end{equation}
Then,  
\begin{equation}\label{gamma-J-dom}
\mbox{\rm sign}(b)\gamma_{b}^J(\cdot) \ge \mbox{\rm sign}(b) \gamma_{b}(\cdot) \quad \mbox{\rm on $J$}.
\end{equation}
Further,  \eqref{gamma-J-dom} remains true if $\sfrho_\vep$ in  both  definitions of $\gamma_b$ and $\gamma_b^J$ is replaced by 
an arbitrary slope $\sfrho$.
\end{corollary}
\begin{proof}[Proof of \Cref{lem-bend1}] Let $a:=s/n$, $b:=(t-s)/n>0$ and $c:=(n- t)/n < 1-a$.  Dividing 
\eqref{bent.1}  by $n^\gak$ and expressing it via $a, b, c$ and 
$x:=(l-s)/(t-s) \in [0,1]$,  we see that  \eqref{bent.1} is a consequence of 
\begin{equation}
f (x) :=  (c)_1^\gak x + (a)_1^\gak (1-x) + b^\gak  (x)_1^\gak - (a+b x)_1^\gak \ge 0 \,, \quad \forall x \in [0,1]\,.
\label{bent.1-sc}
\end{equation}
We proceed to prove \eqref{bent.1-sc},  first for $0 \le c \le a < \xs := \frac{1}{2}$.  In that case,  setting  
$\xz \in [0,\xs]$ such that $b \xz = \xs - a$,  
it is not hard to see that for $\gak \in (0,1]$,  the continuous function $f(\cdot)$ is concave in each 
of the intervals $I_1=(0,\xz)$,  $I_2=(\xz,\xs)$ and $I_3=(\xs,1)$ (where we rely on $b \le 1-a$ 
and on having also $2 b x \le b
$ throughout $I_2$).
Thus,  $\min_{I_j} f(x)$ must be at one of the ends of $I_j$,  for $j=1,2,3$.  Noting that 
$f (0)=f (1)=0$,  it remains only to check that 
\begin{align*}
f(\xz) =   a^\gak + ( b  \xz)^\gak - (a+b \xz)^\gak - (a^\gak - c^\gak) \xz &\ge 0 \,,  \\
f(\xs)= c^\gak + (b \xs)^\gak - (1-a-b \xs)^\gak  + (a^\gak-c^\gak) \xs &\ge 0 \,.
\end{align*}
Having $\gak \le 1$,  $a \ge c$ and $1-a-b \xs = c + b \xs$,  implies that $f(\xs) \ge 0$,  leaving only the task 
of verifying that $f(\xz) \ge 0$,  or equivalently,  to show that for each $a \in [0, \xs]$,  the concave function 
\[
g_a(c) :=  (1-a-c) [ a^\gak + (\xs - a)^\gak - \xs^\gak] - (a^\gak - c^\gak) (\xs-a) \,,
\]
is non-negative on $[0,a]$.  
With $g_a (a) = (1-2a) [ a^\gak + (\xs - a)^\gak - \xs^\gak] \ge 0$,  it thus 
suffices to check that $g_a (0) \ge 0$.  Setting $v:=(1-2a)^{1/3}$ and recalling that $2\xs =1$,  this is equivalent to 
showing that
\begin{equation}\label{bent.2}
h_\gak(v) := (1-v^3)^{\gak} - (1+v^3) (1-v^{3\gak}) \ge 0 \,,  \qquad \forall v \in [0,1]\,.
\end{equation}
As $\gak \mapsto h_\gak(v)$ decreases monotonically (for $v \in [0,1]$),  it suffices to establish \eqref{bent.2}
at $\gak=\frac{2}{3}$.  To this end,  it is not hard to verify that the polynomial 
$p(v):=[(1-v^3)^2-(1+v^3)^3(1-v^2)^3]/v^2(1-v)^2$
is non-negative on $[0,1]$.  Thus,  $(1-v^3)^{\frac{2}{3}} \ge (1+v^3)(1-v^2)$ for any $v \in [0,1]$.
This is equivalent to \eqref{bent.2} at $\gak=\frac{2}{3}$,  thereby 
completing the proof of \eqref{bent.1} in case $n \le s+t$,  $s < n/2$.

Note that $a+bx+c+b(1-x)=1$,  hence the value of $f(\cdot)$ is unaffected by 
exchanging both $a \leftrightarrow c$ and $x \leftrightarrow (1-x)$.  The preceding proof
thus also establishes \eqref{bent.1-sc} for  $0 \le a \le c < \xs$ and $\gak \in (0,\frac{2}{3}]$,
thereby extending the scope of \eqref{bent.1} to all $t>n/2>s$.

We next move to $\xs \le c < 1-a$,  in which case 
$a+bx \le \frac{1}{2}$ for any $x \in [0,1]$.  Yet again $f(0)=f(1)=0$ and $f(\cdot)$ is now concave
in $I_1=(0,\xs)$ and in $I_2=(\xs,1)$.  With $2\xs=1$,  it thus suffices in this case to check that  
\begin{equation}\label{bent.2-sc}
f(\xs) = [(1-c)^\gak- a^\gak] \xs + [(2a)^\gak + b^\gak - (b+2a)^\gak] \xs^\gak \ge 0 \,.
\end{equation}
To establish \eqref{bent.2-sc} and thereby also \eqref{bent.1},  note that for any 
$\gak \in [0,1]$,  the left bracket in \eqref{bent.2-sc} is non-negative since $1-c > a$,  
with its right bracket non-negative by the concavity of $x \mapsto x^\gak$.  We are done,  as 
exchanging $a \leftrightarrow c$ and $x \leftrightarrow (1-x)$ extends  the latter reasoning 
to apply also for $\xs \le a < 1-c$.
\end{proof}

\subsection{The joint law of off-centered Brownian excursion counts}\label{subsec:off-center}

We start by estimating the effect of a small off-centering
on certain probabilities 
of Brownian excursion counts.
See  \Cref{angle-off} for an illustration of our setting.

Throughout this section, we let $\NN_{\sfv,r,\ell'}$ denote the $\sigma$-algebra
generated by the counts of excursions between the concentric circles 
$\partial \bB(\sfv,  r e^{-i})$,  $0 \le i \le \ell'$ (or $\partial \bB_{\sfd}(\sfv,r e^{-i})$ when $\sfv \in \bT^2$).
Viewing these excursion counts as points in $\mathbb{Z}_+^{\ell'}$, 
elements of $\NN_{\sfv,r,\ell'}$ are thus countable unions of elements of $\mathbb{Z}_+^{\ell'}$,
\blue{but we often use $\NN_{\sfv,r}$ when the value of} $\ell'$ is inconsequential.  We let $\PRzu(\cdot)$ denote the law of  a planar Brownian path  started from $\sfz$ and 
stopped upon exiting $\bB(\sfu,e R)$.
\blue{When applying $\PRzv(\cdot)$ at  $\sfz \in \partial \bB(\sfv, d_\star R)$} to events in $\NN_{\sfv,r,\ell'}$ 
we use the short-hand notation 
$\bP_R^{d_\star}(\cdot)$,  since the probabilities involved do not depend on $\sfz$,  see Remark \ref{rmk:off-center}.  With some abuse of notation, we continue to use $\PRzu(\cdot)$  and $\bP_R^{d_\star}(\cdot)$ for the law of Brownian motion on  $\bT^2$, if $R\leq e^{-2}$.
\begin{lemma}\label{lem:easy-off}
Fix non-negative $\rho$,  $d_\star$,  $r$,  $R$,  $\ell' < \infty$ and $\sfz, \sfu,\sfv \in \bR^2$ such that 
\begin{align}\label{eq:par-rt}
\frac{r}{R} & \le d_\star (1-\rho) \le 1- \rho \,,  \quad  \rho \le e^{-1} \\
|\sfu-\sfv| &\le \rho R \quad {\rm and} \quad |\,   |\sfz-\sfv| - d_\star R | \le \rho d_\star R \,.
\label{eq:par-uy}
\end{align}
Denoting by
$\wt{A}_k \in \NN_{\sfv,r,\ell'}$ the event of having $k$ excursions from $\partial \bB(\sfv, r)$ 
to $\partial \bB(\sfv, r/e)$,  for some $c<\infty$,  any $k \ge 0$ and such parameters,
\begin{equation}\label{eq:off1-center-bd}
 (1+c \rho)^{-(k+1)} \, \bP_R^{d_\star} (A) \le 
\PRzu(A) \le (1+c \rho)^{(k+1)} \, \bP_{R}^{d_\star}(A) \qquad \forall A \in \NN_{\sfv, r,\ell'},   A \subseteq \wt{A}_k \,.
\end{equation}
\end{lemma}

\begin{remark}\label{rmk:off-center} If $\rho=0$ then $\sfu=\sfv$ where by radial symmetry
$\sfz \mapsto \PRzv(A)$ depends 
only on $|\sfz-\sfv|=d_\star R$ when $A \in \NN_{\sfv,r,\ell'}$.  It is further
also invariant under translation $(\sfz,\sfv) \to  (\sfz-\sfv,\npp)$ 
sending $\NN_{\sfv,r,\ell'} \mapsto \NN_{\npp,r,\ell'}$,  motivating 
our short-hand notation $\bP^{d_\star}_{R}$.  By the same reasoning,  it suffices to prove
\Cref{lem:easy-off} for $\sfv=\npp$ (translating also $\sfu \to \sfu-\sfv$).  Having done so,  since for 
any $a>0$ the law $\bP_a^{a \sfz,a \sfu}$ on $\NN_{\npp,a r,\ell'}$ matches 
that of 
$\bP_1^{\sfz,\sfu}$ on $\NN_{\npp,r,\ell'}$,  we can also fix $R=1$ when proving \Cref{lem:easy-off}.
\end{remark}

\begin{figure}[h]
\begin{center}
\begin{tikzpicture} 
\draw [thin,dashed](0.3,0.3) circle(1.9);
\draw [very thick](0,0) circle(3);
 \draw[thick, dotted] (0,0) -- (0.4,0.4);
 \draw[thick, dotted] (0.3,0.3) -- (1.1,0.9);
 \draw(0.3,0.3) circle (0.7);
 \draw[very thick] (0.3,0.3) circle (1);
 \draw[very thick] (0.3,0.3) circle (1.3);
 \draw [thick, dotted] (0.3,0.3) -- (1,-1.5);
 \draw [thick, dotted] (0.3,0.3) -- (1.6,0);
 \draw [thick, dotted] (0,0) -- (-1,-2.8);
  \draw [red,{-Stealth},thick] plot [smooth, tension=0.5] coordinates {(-1.35,1) (-1,0.9)};
  \draw [red,{-Stealth},thick] plot [smooth, tension=0.5] coordinates {(-1,0.9) (0.2,0.9) (0.5,2) (0.9,1) (1.1,1.5) (1.4,-0.6) (2.9,0.65)};
  \fill (0,0) circle (1.5pt) node[left] {$\sfu$};
  \fill (0.3,0.3) circle (1.5pt) node[right] {$\npp$};
   \fill (-1.35,1) circle (1.5pt) node[below] {$\sfz$};
   \fill (.25,.25)  node[left] {$\rho$};
    \fill (0,0.7)  node[right] {$r/e$};
    \fill (1.5,0) node[above] {$r$};
    \fill (0.8, -1.1)  node[right] {$d_\star (1+\rho)$};
     \fill (-0.8, -2.3)  node[left] {$e$};
     \fill (0,-3) node[above] {$\bB(\sfu,e)$};
\end{tikzpicture}
\end{center}
\caption{Off centered excursions:
$\sfv=\npp$,  $R=1$ and a path satisfying $\wt{A}_3$.}
\label{angle-off}
\end{figure}

Restating \Cref{lem:easy-off} for off-centered Brownian excursion counts in $\bT^2$ yields the following.
\begin{corollary}\label{cor:easy-off}
Fix parameters as in \eqref{eq:par-rt} for $R \le e^{-2}$ and $\sfz,\sfu,\sfv \in \bT^2$ such that 
\begin{equation}\label{eq:z1-v-cond}
\sfd(\sfu,\sfv) \le \rho R
 \quad {\rm and}  \quad 
 |\sfd(\sfz,\sfv) - d_\star R| \le \rho \,d_\star  R \,.
\end{equation}
Then,  \eqref{eq:off1-center-bd} holds for the Brownian motion on $\bT^2$.
\end{corollary}

\begin{proof}[Proof of \Cref{lem:easy-off}]  Recalling \Cref{rmk:off-center},  we set $\sfv=\npp$,  $R=1$ and 
fix $\rho$,  $r$ and $d_\star$ satisfying \eqref{eq:par-rt}.  Considering  
a Brownian path from $\sfz=t \sfy$ with $|\sfy|=d_\star$,  stopped upon first exiting $\bB(\sfu,e)$, 
we are to show that for some $c$ finite,  any $\rho \le e^{-1}$ and $k \ge 0$,  
\begin{equation}\label{eq:off-center1}
(1  + c \rho)^{-(k+1)} \bP_1^{\sfy,\npp} (A) \le \bP_1^{\,  t \,  \sfy,\sfu} (A) \le 
(1+ c \rho)^{k+1} \bP_1^{\sfy,\npp}(A)  \qquad \forall A \in \NN_{\npp,r},  A \subseteq \wt{A}_k,
 |\sfu| \le \rho,  |\sfy| = d_\star \,,
\end{equation}
which in view of \eqref{eq:par-rt}-\eqref{eq:par-uy} should hold for any choice of
\begin{equation}\label{eq:par-rho}
|t-1| \le \rho,  \quad  
r \le d_\star (1-\rho) \le 1-\rho \,.
\end{equation}
Setting $B_1 := \bB(\npp,  r/e) \subset \bB(\npp,r):= B_0$ and $B^{\pm}:=\bB(\npp,e \pm \rho)$ 
we have that for any such $\sfu$,  $\sfy$,  $\rho$,  $t$,  $d_\star$,  $r$,  
\[
 B_1 \subset B_0 \subset B^- 
 \subset \bB(\sfu,e)  \subset B^+ 
\quad {\rm and} \quad  
t \,  \sfy \in B^- \setminus B_0 \,.
\]
We define $\bP_{+,k}^{d}(A)$ and $\bP_{-,k}^{d}(A)$ where 
 $\bP_{\pm, k}^{d}(A)$ stand for the probability that starting at $\partial \bB(\npp,d)$,
the excursion counts specified by $A \in \NN_{\npp,r}$
occur during $k$ excursions from $\partial B_0$ to $\partial B_1$ 
that do not exit $B^{\pm}$ and upon returning to $\partial B_0$ at the end of the $k$-th such 
excursion,  the path exits $B^{\mp}$ before re-entering $B_1$. 
Since any $A \in \NN_{\npp,r}$,  $A \subseteq \wt{A}_k$ specifies having 
exactly $k$ excursions from $\partial B_0$ to $\partial B_1$,  we have that
\begin{equation}\label{eq:cont-off}
\bP_{-,k}^{ t \,  d_\star}(A) \le \bP_1^{\,  t \,  \sfy,\sfu} (A) \le  \bP_{+,k}^{t \,  d_\star}(A) \,,
\end{equation}
In particular,  
\eqref{eq:cont-off} holds for $\bP_1^{\sfy,\npp}(A)$ and 
we get \eqref{eq:off-center1} upon showing that for $c$ finite and all relevant $r$,  $d_\star$,  $\rho$,  $k$,
\begin{equation}\label{eq:cnt-sym}
\sup_{|t-1| \le \rho} \{  \bP_{+,k}^{ t \,  d_\star}(A) \} \le (1+c \rho)^{k+1} \inf_{|t-1| \le \rho} \{ \bP_{-,k}^{t \,  d_\star}(A) \} 
\qquad \forall A \in \NN_{\npp,r},  A \subseteq \wt{A}_k \,.
\end{equation}
With $\bP_{\pm,k}^{d}(A)$ concerning only the joint law of excursion counts between concentric circles,  we have 
by the Markov property of such counts that $\bP_{+,k}^{t \, d_\star}(A|\wt{A}_k)=\bP_{-,k}^{t \, d_\star}(A|\wt{A}_k)$
is further independent of $t \in [1-\rho,1+\rho]$.  It thus suffices to establish \eqref{eq:cnt-sym} only for 
$A=\wt{A}_k$.  To this end,  recall that a Brownian path starting at $\partial \bB(\npp,d)$ 
with  $\frac{r}{e} \le d \le e-\rho$,  enters $B_1$ before exiting $B^\pm$,  with probability
\begin{equation}\label{eq:nom-q}
q_{\pm \rho}(d)=\frac{\log(e \pm \rho)-\log (d)}{\log (e \pm \rho )-\log (r /e)}
\end{equation}
(see e.g. \cite[(2.1)]{BRZ2}),  hence by definition 
\begin{equation}\label{eq:off-bd}
\bP_{\pm,0}^{d}(\wt{A}_0) = 1-q_{\mp \rho}(d), \qquad
\bP_{\pm,k}^{d}(\wt{A}_{k}) = q_{\pm\rho} (d) q_{\pm \rho} (r)^{k-1} (1-q_{\mp \rho} (r)), \quad k \ge 1 \,.
\end{equation}
\blue{Utilizing this,  one can check that}
the extremes in \eqref{eq:cnt-sym} occur at $t = 1 \pm \rho$.
Setting $p(\pm \rho):=q_{\pm \rho}(r)$,  $s(\pm \rho):=q_{\pm \rho} (d_\star (1 \mp \rho))$,
it \blue{then} suffices to show that for some $c$,  
all $\rho \le  e^{-1}$,  $d_\star \le 1$ and $r$ as specified,   
\begin{equation}\label{eq:calc-bd}
\max\Big\{  \frac{1-s(-\rho)}{1-s(\rho)},  \frac{s(\rho)}{s(-\rho)},   \frac{p(\rho)}{p(-\rho)},  \frac{1-p(-\rho)}{1-p(\rho)} 
\Big\} \le 1+c\rho \,.
\end{equation}
Setting $b:=\log(e/r) \ge \log(e/d_\star) := b_\star \ge 1$, 
note that $1-p(0) = \frac{1}{b+1} \le \frac{1}{2}$ while $s(0)=\frac{b_\star}{b+1} \in [1-p(0),p(0)]$, 
with \eqref{eq:calc-bd} holding if $b |p'(x)|$ and $b |s'(x)|$ are bounded,  uniformly in $|x| \le e^{-1}$ 
and $b \ge b_\star \ge 1$.  To verify the latter property,  check that 
\[
p(x)=1-\frac{1}{1+b+f_1(x)}\,,   \qquad s(x)=1 + \frac{f_2(x)+b_\star-1-b}{1+b+f_1(x)} \,,
\]
where $f_1(x):=\log(1+x/e)$ and $f_2(x):=-\log(1-x)$.
Consequently,  with $f_i(0)=0$,  $f_1(x) \ge -1$ and $f_i'(x)$,  $i=1,2$ uniformly bounded on $|x| \le e^{-1}$, 
the same applies for $b|p'(x)|$ and $b|s'(x)|$,  which by the preceding completes the proof of the lemma.
\end{proof}

Analogously to \eqref{eq:cont-off},  the probability $q^{\sfz,\sfu}_{r,\sfv}$
of entering $\bB_\sfd(\sfv,r R/e)$ before exiting $\bB_\sfd(\sfu,e R)$ is sandwiched between 
$\blue{s(\pm \rho)=}1-\bP_{\pm,0}^{t\,d_\star}({\wt A}_0)$ \blue{for $t=d_\star (1 \mp \rho)$.  
As before $s(0)=\frac{b_\star}{b+1}=\frac{\log(e/d_\star)}{1+\log (e/r)}$,  so
the uniform boundedness of $b |s'(x)|$ implies that}
for some $c$ finite and all $\rho \le \rho_o$
\begin{equation}\label{eq:pv-u}
(1- c \rho) \frac{\log(e/d_\star)}{1+\log(e/r)} \le q^{\sfz,\sfu}_{r,\sfv} \le (1+c \rho) \frac{\log(e/d_\star)}{1+\log(e/r)} \,.
\end{equation}

We proceed with a generalization of \Cref{cor:easy-off} to $m \ge 1$ paths in $\bT^2$, 
where $\NN_{\sfv,r,m,k}$ denotes the collection of events in the corresponding
product $\sigma$-algebra $(\NN_{\sfv,r})^m$ 
where the aggregate number of excursions from $\partial \bB_\sfd(\sfv,  r)$
to $\partial \bB_\sfd(\sfv,  r/e)$ is at most $k \ge 0$ and by $(\bP_R^{1})^m$ the $m$-th product law 
of  $\bP_R^{d_\star}$ with $d_\star=1$
(recall that $\NN_{\sfv,r}=\NN_{\sfv,r,\ell'}$ for some $\ell'$,  hence
events in $\NN_{\sfv,r,m,k}$ are countable unions of elements of $\mathbb{Z}_+^{m \ell'}$).
We also let $\PRsu(\cdot)$ denote the law of a Brownian motion 
on $\bT^2$,  
starting outside $\bB_\sfd(\sfu,R_1)$ and stopped at the end of its 
$\sfs$-th excursion from  
$\partial \bB_{\sfd} (\sfu,R_1)$ to $\partial \bB_{\sfd} (\sfu,\blue{e R_1})$.
\bl\label{lem-off-center-applies} 
Fix $\sfu \in \bT^2$,  $m,\sfs,\sfl \ge 1$, 
strictly decreasing $R_i$,  $0 \le i \le \sfl$ with \blue{$R_0=e R_1 \le e^{-1}$},  $R:=R_{\sfl}=R_{\sfl-1}/e$
and $\sfv \in \bB_{\sfd}(\sfu,\rho R)$,  $\rho \le e^{-1}$,  $r \le (1-\rho)R$.
Let
 $\AA^{(m)}$ denote an event involving only the excursion 
counts between $\partial \bB_{\sfd} (\sfu,R_i)$ and
$\partial \bB_{\sfd} (\sfu,R_{i-1})$,  $2 \le i \le \sfl$,  with precisely $m$ excursions \blue{specified for}
$i=\sfl$.  Then,  for $c$ of \Cref{lem:easy-off} and any $k \ge 0$, 
\begin{equation}\label{eq:off-center-bd} 
 (1+c \rho)^{-(k+m)} \, (\bP_R^{1})^m (\Gamma) \le 
\PRsu(\Gamma|\AA^{(m)}) 
\le (1+c \rho)^{(k+m)} \,  (\bP_R^{1})^m (\Gamma) \qquad \forall \Gamma \in \NN_{\sfv,r, m,k} \,,
\end{equation}
where in the term $\PRsu(\Gamma|\AA^{(m)}) $, we consider for the event $\Gamma$ the excursion counts induced by the  $m$ excursions  between $\partial \bB_{\sfd} (\sfu,R_\sfl)$ and
$\partial \bB_{\sfd} (\sfu,R_{\sfl-1})$.
\el

\begin{figure}[h]
\begin{center}
\begin{tikzpicture}
\draw (0,0) circle (3); 
\draw (0,0) circle (3.7); 
\draw (0,0) circle (1); 
\draw (0,0) circle (1.3); 
\draw (-0.3,0.05) circle (0.15);
\draw [dashed] (0,0) circle (0.5);
%
%
\draw [red,{-Stealth}, dashed,thick] plot [smooth, tension=2] coordinates {(-0.4,4.2) (0.2,3.5) (0.8,4) (0.7,2.92)};
\draw [blue, {-Stealth},dashed,thick]  plot [smooth, tension=2] coordinates { (0.7,2.92) (1, 2.6) (1.5,3.4)};
\draw [red,{-Stealth},thick] plot [smooth, tension=1] coordinates {  (1.5,3.4) (1.9, 3.7) (2.3,2.7) (2.8, 3.0) (2.5,1.7)};
\draw [blue,{-Stealth},thick] plot [smooth, tension=1] coordinates {  (2.5,1.7) (2,1.7) (0.5, 1.8) (0,1)};
\draw [black,{-Stealth},thick] plot [smooth, tension=1] coordinates {(0,1) (0.3,0.4) (0.9,0.96)};
\draw [green,{-Stealth},thick] plot [smooth, tension=2] coordinates { (0.9,0.96) (1.3,1.5) (1,0.1)};
\draw [black,{-Stealth},,thick] plot [smooth, tension=1] coordinates { (1,0.1) (0.4,-0.2) (0.8,-1.05)};
\draw [green,{-Stealth},thick] plot [smooth, tension=2] coordinates {  (0.8,-1.05) (1.1, -1.7) (0,-1)};
\draw [black,{-Stealth},thick] plot [smooth, tension=1] coordinates { (0,-1) (-0.4,0) (-1.3,0)};
\draw [green,{-Stealth},thick,dashed] plot [smooth, tension=1.5] coordinates {  (-1.3,0) (-2,-0.7) (-2.7,-0.6) (-3.5,-1.2) };
\draw [red,{-Stealth},thick] plot [smooth, tension=1.5] coordinates {   (-3.5,-1.2)  (-4,-2) (-1.1,-2.8)};
\fill (-.3,0.05) circle (1pt) node [above] {$\sfv$};
\fill (0,0) circle (1pt) node [below] {$\sfu$};
\fill (0,1) circle (1.5pt) node [left] {\small $\sfz_1$};
\fill (1,0.1) circle (1.5pt) node [below] {\small $\sfz_2$};
\fill (0,-1) circle (1.5pt) node [right] {\small $\sfz_3$};
\node at  (-0.6,-0.53) {$R$}; 
\node at (-1.3,-0.9) {$e R$};
\node at (0,-3.2) {$R_1$};
\node at (0,-3.9) {$R_0$};
\fill (-0.4,4.2) [red] circle (1.5pt);
\fill (0.9,0.96) [green] circle (1.5pt) node [above] {\small $\sfy_1$}; 
\fill (0.8,-1.05) [green] circle (1.5pt) node [right] {\small $\sfy_2$}; 
\fill (-1.3,0) [green] circle (1.5pt) node [left] {\small $\sfy_3$}; 
\fill (1.5,3.4) [red] circle (1.5pt) node [right] {\small $\sfw_1$}; 
\fill (-3.5,-1.2) [red] circle (1.5pt) node [above left] {\small $\sfw_2$}; 
\fill (0.7,2.92) [blue] circle (1.5pt) node [left] {\small $\sfx_1$}; 
\fill (2.5,1.7) [blue] circle (1.5pt) node [right] {\small $\sfx_2$}; 
\fill (-1.1,-2.8) [blue] circle (1.5pt) node [below] {\small $\sfx_3$}; 
%
%
%
%
%
\end{tikzpicture}
\end{center}
\caption{Excursions decomposition: $(\NN_{\sfv,r})^3$ generated by $m=3$ inner {\bf black} paths
(from $\sfz_j$ to $\sfy_j$,  $j \le 3$),   within 
$\sfs=2$ \red{outer excursions} (from $\sfw_i$ to $\sfx_{i+1}$,  $i \le 2$);
Event $\AA^{(m)}$ is determined by the
ordered list of \textcolor{blue}{blue} and \green{green} paths
(where the dashed \textcolor{blue}{blue} and \green{green} avoid $\bB_{\sfd}(\sfu,R)$).
}
\label{fig:village1}
\end{figure}
\begin{proof}  It suffices to consider $\sfl \ge 2$ and assume \abbr{wlog} 
that for some non-random 
\begin{equation}\label{eq:mi-s}
0 = m_0  \le m_1 \le \cdots \le m_{\sfs} := m,
\end{equation}
the event $\AA^{(m)}$ dictates having precisely $m_i$ of our inner excursions completed by the end of the $i$-th 
outer excursion (e.g.  \Cref{fig:village1} has $\sfs=2$,  $m_1=0$ and $m=m_2=3$).
The proof is based on a path decomposition.
Using throughout our running convention,  of viewing 
the path as a directed curve up to time re-parameterization,  we proceed to decompose the first $\sfs$ 
outer excursions,  between circles of radii $R_1$ and $R_0$ into outer segments 
$\RR_0$ till first entering $\bB_{\sfd}(\sfu,R_1)$ at some $\sfx_1$ (the 
dashed red curve in \Cref{fig:village1}),  followed by 
$\RR_i$ from some $\sfw_i \in \partial \bB_{\sfd}(\sfu,R_0)$ to some 
$\sfx_{i+1} \in \partial \bB_{\sfd}(\sfu,R_1)$ (marked there by solid red),  interlaced 
between out-connecting segments $\CC_i$ (within $\bB_{\sfd}(\sfu,R_0)$,  marked in blue),
that lead from $\sfx_i$ to $\sfw_i$ if $m_{i}=m_{i-1}$,  or otherwise
end at $\sfz_{m_{i-1}+1} \in \partial \bB_{\sfd}(\sfu,R_{\sfl})$
(with a dashed line in former case and a solid one in the latter).  In case $m_{i+1}>m_i$ we further decompose
the remainder of the $i$-th outer excursion (between $\CC_i$ and $\RR_{i}$),  to its inner parts 
$\BB_{j}$,   $j \in (m_{i-1},m_{i}]$,  from $\sfz_j \in \partial \bB_{\sfd}(\sfu,R_{\sfl})$ to 
$\sfy_j \in \partial \bB_{\sfd}(\sfu,R_{\sfl-1})$ staying within $\bB_{\sfd}(\sfu,R_{\sfl-1})$
(marked black in \Cref{fig:village1}),  interlaced between the in-connecting paths 
$\GG_{j}$,  $j \in (m_{i-1}, m_i]$ which stay outside $\bB_{\sfd}(\sfu,R_{\sfl})$,
leading from $\sfy_j$ to $\sfz_{j+1}$ (as solid green lines),  except for $\GG_{m_i}$
which goes out,  ending at $\sfw_i$ (marked by a green dashed line).  This decomposition is uniquely 
specified by the parameters of \eqref{eq:mi-s},  where for illustration,  the case
depicted in  \Cref{fig:village1} corresponds to having 
\[
\RR_0,  \CC_1,\RR_1,\CC_2,\BB_1,\GG_1,\BB_2,\GG_2,\BB_3,\GG_3,\RR_2\,.
\] 
Note that $\Gamma \in \sigma(\BB_j,j \le m)$ and $\AA^{(m)} \in \sigma(\CC_i,\GG_j, j \le m, i \le \sfs)$.  In particular, 
$\Gamma \in \NN_{\sfv,r,m,k}$ is a countable union of disjoint sets of the form $\Gamma_1 \times \cdots \times \Gamma_m$,  for $\Gamma_j \in \NN_{\sfv,r}$ such that $\Gamma_j \subset \wt{A}_{k_j}$,  $\sum_j k_j \le k$ 
specifying the relevant excursion counts (around $\sfv$),  during the segment $\BB_j$.  Likewise,  
$\AA^{(m)}$ is a countable disjoint union of product sets $\AA_0 \times \AA_1 \cdots \times \AA_m$ 
where $\AA_0$ specifies the excursion counts between balls around $\sfu$,
of radii $R_i$ and $R_{i-1}$,  $2 \le i < \sfl$,  
up to the start of the segment $\BB_1$,  with $\AA_j$,  $j \ge 1$ specifying such excursion counts 
that occur from the end of $\BB_{j}$ to the start of $\BB_{j+1}$.
By the law of total probability,  it thus suffices to establish 
\eqref{eq:off-center-bd} for $\Gamma=\Gamma_1 \times \cdots \times \Gamma_m$ 
and $\AA^{(m)}=\AA_0 \times \cdots \times \AA_m$
(excluding all choices $\AA_j$ and $\Gamma_j$ of zero probability,  or which are inconsistent 
with the given $\{m_i\}$).
Fixing such $\{\Gamma_j\}$ and $\{\AA_j\}$,  the stated bound follows from having for any $1 \le j \le m$,
\begin{equation}\label{eq:off-center-prod}
(1+c \rho)^{-(k_j+1)} \, \bP_R^{1} (\Gamma_j)  \le 
\frac{\PRsu(\Gamma_1 \times \cdots \times \Gamma_j\,|\, \AA_0 \times\cdots \times \AA_j)}
{\PRsu(\Gamma_1 \times \cdots \times \Gamma_{j-1}\,|\, \AA_0 \times\cdots \times \AA_{j-1})} \le
 (1+c \rho)^{(k_j+1)} \bP_R^{1} (\Gamma_j)  \,.
\end{equation}
Indeed,  applying \eqref{eq:off-center-prod} backwards from $j=m$ to $j=1$ establishes \eqref{eq:off-center-bd}
for such product sets.  

We thus turn to proving \eqref{eq:off-center-prod}. The ratio of probabilities in the middle of
\eqref{eq:off-center-prod} equals
\begin{equation}
\label{eq-ofernew300926a}
\PRsu (\Gamma_j \, | \, \AA_0 \times\cdots \times \AA_{j-1}, \Gamma_1 \times \cdots \times \Gamma_{j-1})
\frac{\PRsu(\AA_j \,|\,  \Gamma_1 \times \cdots \times \Gamma_{j},\AA_0 \times\cdots \times \AA_{j-1})}
{\PRsu(\AA_j \,|\, \AA_0 \times\cdots \times \AA_{j-1})} \,.
\end{equation}
\blue{Denoting by} $\PP_j$ the $\sigma$-algebra generated by the path up to the start of the segment 
$\BB_{j}$, 
 and by 
$\PP_{j}^+$ the $\sigma$-algebra generated by the path up to the end of the segment $\BB_j$
for $1 \le j \le m$, we \blue{will show} that
\begin{align}\label{eq:radial-eq}
& \PRsu(\AA_j \,|\,  \PP_j^+)= \PRsu(\AA_j ),\\
(1+c \rho)^{-(k_j+1)} \, \bP_R^{1} (\Gamma_j) \le
& \PRsu(\Gamma_j \, | \, \PP_j)
 \le(1+c \rho)^{(k_j+1)} \bP_R^{1} (\Gamma_j)  \,.
\label{eq:bd-Gammaj}
\end{align}
\blue{Since the} conditioning events in the ratio of probabilities in \eqref{eq-ofernew300926a} are $\PP_j^+$ measurable,  by \eqref{eq:radial-eq} that ratio must equal one,  whereby \eqref{eq:bd-Gammaj} together 
with the $\PP_j$-measurability of $ \AA_0 \times\cdots \times \AA_{j-1}, \Gamma_1 \times \cdots \times \Gamma_{j-1}$ would yield \eqref{eq:off-center-prod}.
Turning to show \eqref{eq:radial-eq}, 
note that by the strong Markov property,  
the conditioning on $\PP_j^+$ can only affect the probability of the event $\AA_j$ 
via the $\PP_j^+$-measurable value of $\sfy_j$. 
But $\AA_j$ involves only 
excursion counts between concentric circles centered at $\sfu$,  so its
conditional 
probability 
does not depend on the specific value of $\sfy_j \in \partial \bB_\sfd(\sfu,R_{\sfl-1})$,
yielding
the identity \eqref{eq:radial-eq}.  In \eqref{eq:bd-Gammaj} we have 
the conditional probability of the event $\Gamma_j \in \NN_{\sfv,r}$ occurring 
during our $j$-th inner part $\BB_j$ given $\PP_j$, 
which  
by the above reasoning is between the minimum and maximum,
over $\sfd(\sfz,\sfu)=R$,
of $\PRzu(\Gamma_j)$ of \Cref{cor:easy-off}.
Setting $d_\star=1$,  our assumptions on $r$ and $\rho$ match \eqref{eq:par-rt}, 
with \eqref{eq:z1-v-cond} holding by the triangle inequality.  Thus,    
applying \eqref{eq:off1-center-bd} for $A=\Gamma_{j} \subset \wt{A}_{k_j}$ we
verify \eqref{eq:bd-Gammaj} and thereupon complete our proof of \eqref{eq:off-center-bd}.
\end{proof}

\begin{remark}\label{rmk:off-center-applies}
Our proof of \Cref{lem-off-center-applies} is easily adapted to handle 
$\sfs$ independent planar Brownian paths from $\partial \bB(\sfu,R_1)$ to  
$\partial \bB (\sfu,R_0)$ (conditioned on $\AA^{(m)}$ of the form specified in the lemma).
In particular,  taking $\sfl=1$ and $\sfs=m$ 
then corresponds to $m$ \abbr{iid} Brownian paths,  each stopped at its first 
exit from $\partial \bB_{\sfd}(\sfu,e R_\sfl)$.
\end{remark}

Building on \eqref{eq:pv-u} we next control,  in preparation for \eqref{eq:sec-mom-bd}, 
the conditional probabilities
that in the setting of \Cref{lem-off-center-applies} 
all $m$ excursions from $\partial \bB_{\sfd}(\sfu,e^{\ell} \vep)$ to $\partial \bB_{\sfd}(\sfu,e^{\ell+1} \vep)$
avoid $\bB_\sfd(\sfx,\vep/e) \cup \bB_\sfd(\sfx',\vep/e)$.
 \begin{figure}[h]
\begin{center}
\begin{tikzpicture}
\draw (0,0) circle(3);
\draw (0,0) circle(3.5);
\draw (0,0) circle(1.5);
 \draw [thick, dashed](0,0) -- (.5,.5);
  \draw [thick, dashed](0,0) -- (0,-1.5);
  \draw (.5,.5) circle(0.25);
  \draw(1.6,1.6) circle(0.25);
 \draw [thick,dashed] (3,0) -- (0,0);
  \draw [thick,dotted](-1.4,2.65) -- (.5,.5);
  \draw [thick,dashed](3.2, -1.5) -- (0,0);
    \draw [red,thick] (-1.4,2.65) -- (-1.2,2.85); 
  \draw [red,thick,{Stealth[reversed]}-Stealth,smooth, tension=1.3] plot coordinates {(-1.2,2.85) 
  (0.5, 3.1) (1.5,1.9) (0.7,1.6) (.7,.65) (2.9,2)};
  \fill (0,0) circle (1.5pt) node[below left] {$\sfu$};
  \fill (.5,.5) circle (1.5pt) node[below] {$\sfx$};
    \fill (.7,.65) circle (1.5pt) node[above] {$\sfz'$};
     \fill (1.6,1.6) circle (1.5pt) node[right] {$\sfx'$};
   \fill (-1.4,2.65) circle (1.5pt) node[above left] {$\sfz$};
   \fill (.3,.3)  node[left] {$\rho e^{-n'}$};
    \fill ( 2.3, 0)  node[above] {$e^{-n'}$};
     \fill ( 2.45, -.3)  node[below] {$ e^{-n'+1}$};
       \fill (-0.85, 1.8)  node[left] {$\sfd(\sfz,\sfx)$};
        \fill ( 0, -.8)  node[right] {$e^{-n''}$};
          \draw[thick, dotted] (0.7,0.65) -- (1.6,1.6);
\end{tikzpicture}
\end{center}
\caption{Bounding $q_{\ell,\sfx,\sfx'}^{\sfz,\sfu}$: Avoiding the $\sfx'$-ball after hitting the $\sfx$-ball ($\ell=\sfKf, j=\sfKfv,  \rho\leq \sfKfv^{-5}$).}
\label{angle-off-m1}
\end{figure}
\bl \label{lem-xyN}
For some $C'$, 
all $\sfKfv$ large and $\sfu \in \bT^2$, 
if $j \ge \sfKfv+ \log \dlc + 2$,  $\sfd(\sfx,\sfx') \in [e^{j-1}\vep ,e^{j}\vep]$,  $\sfd(\sfx,\sfu) \le e^{\sfKfv} \vep$, 
$\ell=j+6\log j$ and $y \in [0, \ell]$,  then for $R=e^{\ell} \vep$ and  
$m=\frac{1}{2} (\sfco \ell+y)^2$ paths as 
in  \Cref{lem-off-center-applies},  
\begin{align}
\sup_{\vep}  \{ \PRsu\big(\Gamma^{(\vep)}_{\sfx,\sfx',1} \times \cdots \times \Gamma^{(\vep)}_{\sfx,\sfx',m}
\,|\,\AA^{(m)}\big) \}
\leq & C' e^{-4j} e^{-3y} 
\,, 
\label{xyN.2}
\end{align}
where 
$\Gamma^{(\vep)}_{\sfx,\sfx',j}$ is the event of the $j$-th path
avoiding $\bB_\sfd(\sfx,\vep/e) \cup \bB_\sfd(\sfx',\vep/e)$.  
\el
\begin{proof} The path decomposition in our proof of \Cref{lem-off-center-applies} yields 
that 
\[
\PRsu\big(\Gamma^{(\vep)}_{\sfx,\sfx',1} \times \cdots \times \Gamma^{(\vep)}_{\sfx,\sfx',m}
\,|\,\AA^{(m)}\big) 
\le \sup_{\sfz \in \partial \bB_{\sfd}(\sfu,R)} \big( 1- q^{\sfz,\sfu}_{\ell,\sfx,\sfx'} \big)^m
 \qquad {\rm for} \qquad 
q^{\sfz,\sfu}_{\ell,\sfx,\sfx'} := \PRzu \Big(({\Gamma^{(\vep)}_{\sfx,\sfx'}})^\complement\Big) \,,
\]
\blue{where $\Gamma^{(\vep)}_{\sfx,\sfx'}$ is the event that our (stopped) path
avoids $\bB_\sfd(\sfx,\vep/e) \cup \bB_\sfd(\sfx',\vep/e)$}.
With $4 \ge \frac{m}{\ell^{2}}  \ge 2 + 2 \frac{y}{\ell}
$ 
and $4 j \ge 3 \ell$ for all $\sfKfv$ large enough,  utilizing that $1-q \le e^{-q}$
we get \eqref{xyN.2} upon showing 
that for $C$ finite and all $\sfz \in \partial \bB_{\sfd}(\sfu,R)$,
\bea\label{ver.16-new}
q^{\sfz,\sfu}_{\ell,\sfx,\sfx'} \ge 2 \,  \frac{(j-C)}{\ell^2} \,.
\eea
To this end,  note that the event $({\Gamma^{(\vep)}_{\sfx,\sfx'}})^\complement$ of our 
(stopped) path 
entering $\bB_\sfd(\sfx,\vep/e) \cup \bB_\sfd(\sfx',\vep/e)$
contains the event of hitting one of these (small) balls 
and thereafter 
missing the other ball.   Setting $r=e^{-\ell}$ consider
the latter two disjoint events,  where by the strong Markov property upon first hitting 
the relevant ball of radius $\vep/e$ (say,  at location $\sfz'$),  
\begin{equation}\label{eq:lbd-union-2-balls}
q^{\sfz,\sfu}_{\ell,\sfx,\sfx'} \ge q^{\sfz,\sfu}_{r,\sfx} \inf_{\sfd(\sfz',\sfx)=\vep/e} \{ 1- q^{\sfz',\sfu}_{r,\sfx'} \} 
+ q^{\sfz,\sfu}_{r,\sfx'}  \inf_{\sfd(\sfz',\sfx')=\vep/e} \{ 1- q^{\sfz',\sfu}_{r,\sfx}\}  
\end{equation}
(see \Cref{angle-off-m1} \blue{depicting a sample from the first event considered on the right hand side)}.
Setting $d_\star=1$ for the first term in each product in \eqref{eq:lbd-union-2-balls} 
while $d_\star =\sfd(\sfx,\sfx')/R$ for the second term in each product,  we get from 
\eqref{eq:pv-u} the lower bound $\frac{1-c\rho}{\ell+2} (\frac{j}{\ell+2}-c\rho)$ for each product.
In our setting $\rho \to 0$ as $\sfKfv \to \infty$,  with uniformly bounded
$$
(\ell+2) \rho \le (\ell+2) \max(e^{-j}, e^{-\ell} (e^j + e^{\sfKfv}))\,.
$$
Thus,  each of the two terms we add on the right hand side of \eqref{eq:lbd-union-2-balls} is at least
$\frac{j-(C-8)}{(\ell+2)^{2}} \ge \frac{j-C}{\ell^2}$ and \eqref{ver.16-new} follows. 
\end{proof}

\subsection{Contribution to uncovering from within each block}\label{sec:g4-bds}
Our goal  \blue{here is to approximate,  within} a universal multiplicative factor,
the probability of the event 
$\RUy$  of restricted block uncovering,  see \eqref{eq:barrier-event} and  Figure
\ref{fig-scales}. (This probability eventually contributes the term $\alpha(\sfKfv,\sfKp)$ to $\bP(\UUpy(y) | \underline{\sfx})$ 
in \Cref{prop:angle-diff}.) To do so, it is convenient to re-scale \blue{space} 
so the $\varepsilon$ dependence disappears \blue{upon blowing up by factor 
$e^{-2}/(e^{\sfKf} \vep)$}. 
Instead of keeping track of excursion start points according to the actual law, and keeping in mind the off-center results of Section \ref{subsec:off-center}, we will work in this section  with the proxy law $\bar \bP^{(\star)}_{z}$ of 
$m:=\bar \sfs_{\sfKf} (z)$ \abbr{iid} Brownian excursions from $\partial \bB(\npp,e^{-2})$
to $\partial \bB(\npp,e^{-1})$, starting uniformly on $\partial \bB (\npp,e^{-2})$, 
later \blue{showing} in Section \ref{sec:angle-diff} that the term $\alpha(\sfKfv,\sfKp)$ is described,  
\blue{see \eqref{eq:alpha-tot}},  in terms of a function  $\gfo(z)$ that we introduce \blue{in} \eqref{dfn:g4-alt}.
\blue{Specifically},  set the random variables
\begin{align}\label{def:why}
\eta_{\sfv,j} &:= \sqrt{2\,  \trav{m}{\npp}{e^{-1}}{e^{-2}}{\sfv}{r_{j}}{r_{j-1}}} \quad 
\text{where} \quad r_j :=e^{-(j+2)}\,
\quad \text{and} \quad  \| \sfv \| \le  r_{\sfLw}  \quad \text{for} \quad \sfLw := \sfKf-\sfKfv 
\end{align}
(in terms of traversals $\partial \bB (\sfv, r_j)\to\partial \bB (\sfv,r_{j-1})$ 
within our $m$ excursions,  see also \Cref{def-trav}),  and
\begin{align}
\label{def:eta-sharp}
\eta^\sharp_{\sfv,\sfKf} &:= \min \{ \eta_{\sfu,\sfKf} : \|\sfu-\sfv\|_\infty \le \tfrac{\wt{c}}{\sqrt{2}} r_{\sfKf}\}
\quad \text{where} \quad  \wt{c} \in (0, \tfrac{1}{2}).
\end{align}
(\blue{These
$r_j$ correspond} to the change of scale discussed above,  and should not be \blue{confused} with
$r_j$ of \eqref{dr.1}. The scales involved indeed correspond to the area marked as $\RUy$ in Figure \ref{fig-scales}.)

Introduce now the function
\begin{equation}\label{dfn:g4-alt}
\gfo(z) := e^{2 \sfLw} {\bf 1}_{\sfIs_4}(z) \bar \bP^{(\star)}_{z}(\bigcup_{\sfv \in \bF^{\npp}} \QQs 
\bigcap \{\eta_{\npp,\sfLw} < 3 \sfKfv \}) \,,
\end{equation}
for the net $\bF^\npp := r_{\sfKf} \bF^\star$ (see \eqref{dfn:F-star}),  and events 
\begin{equation}\label{dfn:QQs}
\QQs := 
 \Big\{  \eta_{\sfv,j} > \gamma_{-1} (j),  j \in \sfJus, \;\; \eta^\sharp_{\sfv,\sfKf} = 0 
\Big\} \,,
\end{equation}
where $\Delta:=-\log \dlc$ and the barrier parameters are (see \Cref{fig:EnlandFnl}), 
\begin{align}\label{dfn:tsfKfv}
\tsfKfv &:=\sfKfv-\Delta,  
\quad \sfLo:=\sfKf-\tsfKfv,  \quad  \sfLn:= \tsfKfv-\sfKp,  \quad \sfJus:=[\sfLo,\sfLo+\sfLn] \,,\\
\gamma_{-1}(\sfLo+\ell) &:=\sfco (\tsfKfv-\ell) - (\tsfKfv-\ell)^{\gak}  \,,  \qquad \ell \in [0,\sfLn] \,.
\label{dfn:gam} 
\end{align}
The next lemma gives sharp bounds on $\gfo(z)$ \blue{that we use} in \Cref{lem:alpha-unif-bd} 
to \blue{uniformly bound} $\alpha(\sfKfv,\sfKp)$.
\bl\label{lem:g4-bds} 
For some $0<c(\sfKp)<C(\sfKp)<\infty$ and all $\sfKfv$ large enough 
\begin{align}\label{ubd2-g4}
\gfo(z) &\le  C(\sfKp)\frac{z}{\sfKf}   e^{- \sfco z}  e^{-\frac{z^2}{8 \sfKf}} \,,\\
\label{eq:lbd-alpha}
\gfo(z) &\ge c(\sfKp) \frac{z}{\sfKf} e^{-2 z}\,, 
\qquad \forall z \in \sfIs_4 \cap  [2 \sqrt{\sfKf}, 4 \sqrt{\sfKf}]\,.
\end{align}
\el
Before proving this lemma at the end of the section,  we will estimate the 
$\bar \bP^{(\star)}_{z}$-probability of the events 
\begin{align}\label{def:QQs}
\wQs_{\sfv} &:= \{\eta_{\sfv,j} > \gamma^{\sfJus}_{-1}(j) \, ,  \;\; j \in \sfJus, \;\;  \eta_{\sfv,\sfKf+1}=0 \}\, ,\\
\label{def:whcI}
\whcI_{\sfv}  &:= \{ \eta_{\sfv,j} \in \sfIp_j,  \; j \in \sfJus,  \;\;
\eta_{\sfv,\sfKf} = 0 \},  
\end{align}
using for each $\ell \in [0,\sfLn]$ the windows and barrier (see \Cref{fig:EnlandFnl2}),  
\begin{figure}
\begin{center}
\begin{tikzpicture}[scale=0.5]

\begin{scope}[shift={(-5,0)}]
    \draw   (-5,0) -- (14.4,0) ;
     \draw   (-5,0) -- (-5,9.2) ;
    \draw  (0,0) -- (0,9.2) ;
  \draw [dotted] (13,0) -- (13,9);
    \draw[-] (0,7.1) .. controls (6,3) .. (13,0.3);
\draw[dotted,thick] (0,8.9) .. controls (2,5) and  (4,6.9) .. (6,5);
\draw[dotted,thick] (6,5) .. controls (8,3.2) and  (10,3) .. (13,1);
  \end{scope}
   \begin{scope}[shift={(-5,0)}]
\draw[-] (0,8.3) -- (14,0);
 \node at (-5,-0.4) {\small $0$};
 \node at (0,-0.4) {\small $\sfLo$};
  \node at (5,-0.4) {\small $\sfLo+\ell$};
   \node at (4.7,2.6) {\small $\gamma_{-1}(\sfLo+\ell)$};
 \node at (14,-0.4) {\small $\sfKf$};
 \node at (12,-0.4) {\small $\sfLo+\sfLn$};
  \node at (-1.6,7.0) {\small $2\tsfKfv-\tsfKfv^\gak$};
    \node at (-0.7,8.2) {\small $2\tsfKfv$};
   \node at (-.4,9.0) {\small $x$};
     \node at (13.2,1) {\small $y$};
  \end{scope}
\end{tikzpicture}
\end{center}
\caption{Depiction of $\gamma_{-1}(\sfLo+\ell)$ with a sample from $\wQs_\npp$ (dotted,  see \eqref{eq:decom-p0}
for $x$ and $y$).}
\label{fig:EnlandFnl}
\end{figure}
\begin{align}
\sfIp_{\,\sfLo+\ell} &:=  [ \sfco (\tsfKfv -\ell) + (\tsfKfv-\ell)^{\gak} ,\sfco (\tsfKfv -\ell) + (\tsfKfv-\ell)^{\okappa}]\,,
\label{dfn:Ip-sfJs}\\
\gamma_{-1}^{\sfJus} (\sfLo+\ell) &:=\gamma_{-1} (\sfLo) +   
\frac{\ell}{\sfLn}  \big[ \gamma_{-1} (\sfLo+\sfLn) -  \gamma_{-1} (\sfLo)\big]
 - (\ell_{\sfLn})^\gak \,.
 \label{dfn:gam-Js}
\end{align}
Now,  since $\bB(\sfv,r_{\sfKf+1}) \subset \bB(\sfu,r_{\sfKf})$ whenever
$\|\sfu-\sfv\|_\infty \le \tfrac{\wt{c}}{\sqrt{2}} r_{\sfKf}$,  we have that 
\begin{equation}\label{eq:sharp-cont}
\{ \eta_{\sfv,\sfKf}=0 \} \subseteq \{ \eta^\sharp_{\sfv,\sfKf} = 0 \} \subseteq \{ \eta_{\sfv,\sfKf+1}=0 \}\,.
\end{equation}
By \Cref{cor:bend1},  see \eqref{gamma-J-dom}, 
we also have that 
$\sfIp_j \subset (\gamma_{-1}(j),\infty) \subset (\gamma_{-1}^{\sfJus}(j),\infty)$, 
so
in view of \eqref{eq:sharp-cont},
\begin{equation}\label{eq:A-RU-Q-hatQ}
\whcI_{\sfv}\subset \QQs \subset 
\wQs_{\sfv}  \,.
\end{equation}
We will thus establish \eqref{ubd2-g4} by bounding above $\sum_{\sfv} \bar \bP^{(\star)}_z (\wQs_{\sfv})$ 
and \eqref{eq:lbd-alpha} by employing the second moment method for
$\bigcup_{\sfv} \whcI_{\sfv}$,  where our lower bounds on $\bar \bP^{(\star)}_z (\whcI_{\sfv})$
match such upper bounds for the larger events
\begin{equation}
\label{def:Aell}
\whcI_{\sfv,\ell} := \{ \eta_{\sfv,j} \in \sfIp_j,  \; j \in 
\sfJus \setminus \{ \ell-2,\ell-1 \},  \;\; \eta_{\sfv,\sfKf} = 0 \}.
\end{equation}
Since $\|\sfv\|$ is small,  we shall see in the proof of  \Cref{lem:g4-cnc}(d) that \Cref{lem-off-center-applies} 
reduces the relevant tasks to bounds for 
$\sfv=\npp$,  where the key is to control the following Markov transition kernels.
\begin{definition}\label{dfn:Markov}
Under $\bar \bP^{(\star)}_z$,  let $q^{(u)}_{\ell',\ell}(x',x)=q^{(u)}_{\ell-\ell'}(x',x)$ for  $u \ge 0$,  $\ell \ge \ell'$, 
be the probability that during $(x')^2/2$ excursions at level $\ell'$ around $\npp$,
the path completed exactly $x^2/2$ excursions at level $\ell$ and at least 
$(x' \wedge x - u)^2/2$ such excursions at levels $j \in (\ell,\ell')$,  using
$q_{\ell-\ell'}(x',x)$ when $u=x' \wedge x$.  Similarly,  denote by
$q^{(\ge)}_{\sfJus}(x',x)$ the probability that in addition to completing exactly $\frac{x^2}{2}$ 
excursions at level $\sfLo+\sfLn$ during $(x')^2/2$ such excursions at level $\sfLo$,  for all
$j \in \sfJus$ the $j$-th level excursion count around $\npp$ exceeds $\gamma^{\sfJus}_{-1}(j)^2/2$
and by $q^{\Cap}_{\sfJus}(x',x)$ the probability that for each $j \in \sfJus$
the number of such excursions 
is further $y^2/2$ for some $y \in \sfIp_j$.
\end{definition}
Indeed,  
for $\sfv=\npp$ the counts of concentric nested excursions 
form a Markov chain and summing over $x^2, y^2 \in 2 \bZ_+$ we have the decompositions 
\begin{align}\label{eq:decom-p0}
\bar \bP^{(\star)}_{z} ( \wQs_{\npp}) &=
\sum_{x,y} q_{0,\sfLo} (2 \sfKf+z,x)q^{(\ge)}_{\sfJus}(x,y) q_{\sfLo+\sfLn, \sfKf+1}(y,0),
\\
\bar \bP^{(\star)}_z (\whcI_{\npp}) &=
\label{eq:new-p0}
\sum_{x,y} q_{0,\sfLo} (2 \sfKf+z,x)q^{\Cap}_{\sfJus}(x,y) q_{\sfLo+\sfLn,\sfKf}(y,0) \,,
\end{align}
with a similar decomposition for $\bar \bP^{(\star)}_z(\whcI_{\npp,\ell})$ of \eqref{def:Aell}, 
except for unrestricted counts at two consecutive levels.  
Utilizing those decompositions yield the following bounds on the relevant events.

\begin{lemma}\label{lem:g4-cnc}
Fix $0 < \kappa \le \gak < \frac{1}{3} \wedge \kapo < \frac{1}{2} < \kapp < \okappa < \frac{3}{4}$.
\newline
(a).  There exists $u_o<\infty$ such that for  any  $(x \wedge y)^\kappa \ge u \ge u_o$ 
and $|x-y| \le (x \wedge y)^{\okappa}$, 
\begin{equation}\label{eq:q3-bd}
q_3^{(u)}(x,[y,y+2)) \ge \frac{1}{4}
 q_3(x,[y,y+2)) \,.
\end{equation}
(b). Fixing $\Delta \ge 0$ and $\sfKp$ large,  for some $C_{\sfKp} < \infty$,  any $\sfKfv$ large and 
$z \in (\sfKf^{\kapo},  \sfKf^{\kapp})$ with $\bar \sfs_{\sfKf}(z) \in \bN$,
\begin{align}\label{ubd-wQ}
 \bar \bP^{(\star)}_{z} ( \wQs_{\npp}) &\le C_{\sfKp} \frac{z}{\sfKf} e^{-2 \sfKf} e^{- 2z}  e^{-z^2/(8 \sfKf)} \,,\\
 \label{ubd-Aell}
 \bar \bP^{(\star)}_{z} (\whcI_{\npp,\ell}) &\le C_{\sfKp} \frac{z}{\sfKf} e^{-2 \sfKf} e^{- 2z}  e^{-z^2/(8 \sfKf)}\,, 
 \qquad \quad \forall \ell \in \sfJus \,.
\end{align}
(c).  There exist $c_{\sfKp} >0$ so if in addition $z \in [2 \sqrt{\sfKf}, 4\sqrt{\sfKf}]$,  then for $\whcI_{\npp}$ of \eqref{def:whcI},
\begin{align}\label{eq:mom-bd-cent}
\bar \bP^{(\star)}_z(\whcI_{\npp}) &\ge c_{\sfKp}  \frac{z}{\sfKf}  e^{-2 \sfKf} e^{-2 z} \,.
\end{align}
(d).  The
\blue{bounds \eqref{ubd-wQ}-\eqref{eq:mom-bd-cent} apply} 
also for   
$\wQs_{\sfv}$,  $\whcI_{\sfv}$ and $\whcI_{\sfv,\ell}$,  uniformly in $\| \sfv\| \le \sfKfv^{-5}$.
\end{lemma}

\begin{proof} (a).  We show \eqref{eq:q3-bd} by bounding above the value of 
$q_3(x,[y,y+2))-q_3^{(u)}(x,[y,y+2))$.  Specifically,  setting $z:=x \wedge y - u$,  
by our definition of $q_3^{(u)}(x,\cdot)$,  the Markov property and the monotonicity 
of $z' \mapsto q_j(z',[y,\infty))$,  it suffices to show that 
\[
q_1(x,(0, z]) q_2 (z,[y,\infty)) + q_2(x,(0,z]) q_1(z,[y,\infty)) \le \frac{3}{4} q_3(x,[y,y+2)) \,.
\]
Setting $m:={x^2}/{2}$,  $\wt{m}:={z^2}/{2}$,  $m':={y^2}/{2}$ and 
$I_{m'}:=[0, 2(\sqrt{2m'}+1))$,  this amounts in turn to
\begin{equation}\label{eq:q3-bd-alt}
\bP(S^{(1)}_m \le \wt{m}) \bP(S^{(2)}_{\wt{m}} \ge m') + 
\bP(S^{(2)}_m \le \wt{m}) \bP(S^{(1)}_{\wt{m}} \ge m')  \le \frac{3}{4} \bP(S^{(3)}_m -m' \in I_{m'})\,,
\end{equation}
where $S^{(j)}_m$ denote the sum of $m$ 
\abbr{iid} non-negative integer variables $\{\GJ'_i\}$ of the form 
\begin{equation}\label{dfn:GJp}
\bP(\GJ'_i=0)=1-q  \quad \text{and} \quad \bP(\GJ'_i=k)=q^2 (1-q)^{k-1},  \qquad \forall k \ge 1
\end{equation}
for the parameter $q:={1}/{(j+1)}$ (see \cite[(2.1)]{BRZ2}). 
\blue{One can then} verify that 
\begin{equation}\label{def:q-cgf}
\Lambda_q(\lambda) :=\log \bE[e^{-\lambda \GJ_i'}] = 
\log \Big( 1 - \frac{q (e^{\lambda}-1)}{q+e^{\lambda}-1} \Big)  \qquad \forall \lambda \ge 0 \,,
\end{equation}
whereby $\GJ'_i$ has finite \abbr{mgf} up to $\lambda=\log(1-q)$,  with   
\begin{equation}\label{eq:GJp-var}
\bE[\GJ'_i]=-\Lambda_q'(0)=1 \qquad \hbox{and} \qquad \hbox{var}(\GJ_i')=\Lambda_q''(0)=\frac{2}{q} (1-q) \,.
\end{equation}
In particular,  
$\wh{S}^{(j)}_m :=(S^{(j)}_m-m)/\sqrt{2jm}$ has mean zero and variance one.  Under our
assumptions on $x$ and $y$,  \abbr{wlog} $I_{m'} \supset [0,\sqrt{6 m})$.  Moreover,  for 
$m':=m+t \sqrt{m}$,  $t_+:=t \vee 0$ and any $u \in [u_o, (2 (m \wedge m'))^{\kappa/2}]$,
\[
\wt{m} - m  = - (-t)_+ \sqrt{m} - u \sqrt{2(m \wedge m')} + \frac{1}{2} u^{2}  \le - (u+ (-t)_+) \sqrt{m} 
\]
provided $u_o$ is large enough \blue{(the last inequality holds since 
$u \le (x \wedge y)^{\kappa} \le x \wedge y - \frac{x}{\sqrt{2}}$ for 
$(x-y)_+ \le (x \wedge y)^{\okappa}$ and $x \wedge y$ large enough)}.
We thus replace \eqref{eq:q3-bd-alt} by the stronger requirement 
\begin{align}\label{eq:q3-new}
\bP\Big(\wh{S}^{(1)}_m \le -\frac{u+(-t)_+}{\sqrt{2}}\Big) \bP\Big(\wh{S}^{(2)}_{\wt{m}} \ge \frac{u+t_+}{2}\Big) &+ 
\bP\Big(\wh{S}^{(2)}_m \le -\frac{u+(-t)_+}{2}\Big) \bP\Big( \wh{S}^{(1)}_{\wt{m}} \ge \frac{u+t_+}{\sqrt{2}}\Big) 
 \nonumber \\
&\le \frac{3}{4} \Big[ \bP\big(\wh{S}^{(3)}_m \ge \frac{t}{\sqrt 6} +1 \big) -
\bP\big(\wh{S}^{(3)}_m \ge \frac{t}{\sqrt 6} \big) \Big] \,.
\end{align}
Denoting by $\Phi(\cdot)$ the standard normal \abbr{cdf},  we have that for some $u_o$ finite, 
any $u \ge u_o$ and all $t$,
\begin{align}\label{eq:q3-Gauss}
e^{\frac{t^2}{25}} \Big[ \Phi\big(-\frac{u+(-t)_+}{\sqrt{2}}\big) \Big(1-\Phi\big(\frac{u+t_+}{2}\big)\Big) &+ \Phi\big( -\frac{u+(-t)_+}{2}\big) \Big(1-\Phi\big(\frac{u+t_+}{\sqrt{2}}\big)\Big) \Big] \nonumber \\
&\le \frac{3}{5} \Big(\Phi\big(\frac{t}{\sqrt{6}}+1 \big) - \Phi\big(\frac{t}{\sqrt{6}}\big) \Big) \,.
\end{align}
Indeed,  as \eqref{eq:q3-Gauss} is of the form $f_2(t,u) \le f_1(t)$ for continuous $f_1,f_2 >0$, 
with $\sup_u f_2(t,u)/f_1(t) \to 0$ as $|t| \to \infty$ and for fixed $b$,  
$\sup_{|t| \le b} \{f_2(t,u)\} \to 0$ as $u \to \infty$,  
it holds for all $t$ and $u \ge u_0$,  provided $u_o$ is large enough.

Now,  fixing $\xi<1/6$ and increasing $u_o$ as needed,  the normal approximation of 
\cite[Thm. ~VIII.2.1]{Petrov} yields \eqref{eq:q3-new} for $u_o \le u \le 2 m^{\kappa/2}$,  
$|t| 
\le m^{\xi}$.  To deal with larger $|t|$ first 
note that our assumptions on $x$ and $y$ imply for large $u_o$ that $|t| \le 2 m^{\okappa/2} \le m^{3/8}$,  
with the \abbr{lhs} of \eqref{eq:q3-new} still uniformly bounded by the \abbr{lhs} of \eqref{eq:q3-Gauss} 
(increasing  $u_o$,  hence $m \wedge m'$,  so the factors $(\pm t^3/\sqrt m) \lambda(\pm t/\sqrt m)$ 
in \cite[Thm.~VIII.2.1]{Petrov} are accommodated by the left-most
exponential factor in \eqref{eq:q3-Gauss}).  Moreover,  mimicking the derivation of \cite[Thm. ~VIII.2.8]{Petrov},  
apply \cite[Thm.~VIII.2.1]{Petrov} for each term on the right hand side of 
\eqref{eq:q3-new} then plug the tail asymptotic of $1-\Phi(\cdot)$ to deduce that their difference,  
for any  $|t| \in [m^{\xi},m^{3/8}]$ and $m \ge m_0$,
is at least $\frac{2}{3} \bP\big(\wh{S}^{(3)}_m \ge t' \big)$ with
$t':=|{t}/{\sqrt 6}| \wedge |{t}/\sqrt{6}+1|$.  
To complete the proof \blue{note} that by \cite[Thm.~VIII.2.1]{Petrov},
 the latter quantity is bounded below by the right hand side of \eqref{eq:q3-Gauss},  up to another 
multiplicative factor,  embedded within the factor $\exp({t^2}/{25})$ of \eqref{eq:q3-Gauss}.
\newline
(b). Clearly,  the probability $q_\ell(y,0)$ of missing a concentric inner ball during 
${y^2}/{2}$ excursions is 
\begin{equation}\label{eq:miss-ball}
q_{\ell}(y,0) =\Big(1-\frac{1}{\ell+1}\Big)^{y^2/2} \le e^{-\frac{y^2}{2 \ell +2}} \,.
\end{equation}
Considering $\ell=\sfKp+1$,  we \blue{can} ignore hereafter the negligible contribution to $\bar \bP^{(\star)}_{z} ( \wQs_{\npp})$
\blue{in \eqref{eq:decom-p0}}
from $y \ge \sfKp \sqrt{6 \sfKf}$.  Now,  bounding $q^{(\ge)}_{\sfJus}(\cdot)$ by \blue{$q_{\sfLo,\sfLo+\sfLn}(\cdot)$}
and summing over $y$ \blue{in \eqref{eq:decom-p0}}
leads to 
\[
\bar \bP^{(\star)}_{z} ( \wQs_{\npp}) \le \sum_x q_{\tsfKfv+1}(x,0),
\]
so by the preceding reasoning we further restrict the sum in \eqref{eq:decom-p0} to 
$x \le a':= 2 \tsfKfv + (\tsfKfv)^{\okappa}$.  Note that by \cite[Theorem 1.1(a)]{BRZ} we have 
for some $C<\infty$ and any $x \ge a := 2\tsfKfv-(\tsfKfv)^\gak$ and $y \ge b := 2\sfKp-\sfKp^\gak$ as above,  
\begin{equation}\label{eq:BRZ-1.1a}
q^{(\ge)}_{\sfJus}(x,[y,y+2)) \le C  \frac{(1+x-a)(1+y-b)}{\sfLn}\sqrt{\frac{x}{y\sfLn}} e^{-\frac{(x-y)^2}{2\sfLn}} \,.
\end{equation}
One may further adapt the proof of  \cite[(1.10)]{BRZ} (taking there $L=\sfLo$),  to have for
some $C$, $c$ finite,  upon setting $\sfLp:=\sfLo+c$,  that for $z \in (\sfKf^{\kapo},  \sfKf^{\kapp})$ 
and any such $x$,  
\begin{equation}\label{eq:BRZ-1.10m}
q_{\, \sfLo} (2\sfKf+z,[x,x+2)) \le \frac{C}{\sqrt{4\pi \sfLp}} e^{-\frac{(2\sfKf+z-x)^2}{2\sfLp}} \,.
\end{equation}
Indeed,  \cite[(1.10)]{BRZ} bounds $q_{\ell}(x',[x,x+2))$ in case $\max(x',x) \le c' \ell$; That restriction is only used in evaluating the last expression in the display after \cite[(4.1),  page 140]{BRZ}.  But its conclusion 
applies as long as $|x'-x| \le c' \ell$,  since combining the sub-Gaussian tails of $w=\sqrt{2\mathcal{L}_1}-x'$ and 
the expansion of \cite{BRZ},  yields the bound on \cite[(1.10)]{BRZ},  up to the uniformly bounded multiplicative factor    
$\int_\bR dw \sqrt{1+|w|/x'} e^{-c w^2} \int_0^2 e^{(x-x')(w-t)/\ell} dt$.
When $|x-x'|>c' \ell$ one gets,  by comparing the law of $(\sqrt{2\mathcal{L}_1}-x')/\sqrt{2}$ 
there to the standard Gaussian law,  an extra error term of the type 
$e^{c(|x'-x|/\ell^2}$ (for some $c \ge {1}/{2}$),  which is absorbed by changing to $2(\ell+c)$ in the exponent in the bound,  
as we have done. 

To help appreciate the following,  we point out that for large $\sfKfv$ we have $\sfLn=(1+o(1)) \sfKfv$
and $\sfLp = O(\log \sfKfv)$.  Plugging the bounds \eqref{eq:miss-ball}-\eqref{eq:BRZ-1.10m} into \eqref{eq:decom-p0} 
we thus find that for some $C'$ finite and any $\sfKf$ large enough, 
\[  
\bar \bP^{(\star)}_{z} ( \wQs_{\npp}) 
\le \frac{C'}{\sfKf} \sum_{x=a}^{a'}  \frac{(1+x-a)}{\sqrt{4\pi \sfLp}} e^{-\frac{(2\sfKf+z-x)^2}{2\sfLp}}
\sum_{y=0}^\infty e^{-\frac{(x-y)^2}{2\sfLn}} \sqrt{y} e^{-y^2/(2 \sfKp+4)},
\]
where now the sums in $x,y$ go over even integers,  and we have used the fact that 
per value of $\sfKp$,  the ratio $x/\sfLn \le a'/\sfLn$ is bounded 
uniformly in $\tsfKfv$.  Using this again,  we see that the preceding sum over $y$ is at most 
$C'_{\sfKp} e^{-x^2/(2\sfLn)}$ for 
some $C'_{\sfKp}$ finite,  resulting with 
\begin{equation}\label{eq:temp-bd}
\bar \bP^{(\star)}_{z} ( \wQs_{\npp}) 
\le \frac{C_{\sfKp}}{3 \sfKf} \sum_{x=a}^{a'}  \frac{(2c+x-a)}{\sqrt{4\pi \sfLp}} e^{-\frac{(2\sfKf+z-x)^2}{2\sfLp}} e^{-\frac{x^2}{2\sfLn}} \,.
\end{equation}
Now using $(A+B)^{2}\geq A^2/2-B^2$ twice,  we have for any $z,w,v$ and $\sfL \ge \sfKp - c \ge  0$,  that 
 \[
 2 ( \sfKp - c) + \frac{(z+w-v+2\sfKp-2c)^2}{2\sfL} \ge \frac{(z+w-v)^2}{4\sfL} 
 \ge  \frac{z^2}{8\sfL} - \frac{(w-v)^2}{4\sfL} \,.
\]
Thus,  \blue{since} $\sfKf-\sfLp
=\sfLn+\sfKp-c$,  \blue{upon} setting $x=2\sfKf - 2\sfLp -v+w+z$  we get that for $\sfLn \ge \sfLp$,  
\begin{align}
\frac{(2\sfKf+z-x)^2}{2\sfLp} + \frac{x^2}{2\sfLn}  
=& \frac{(2\sfLp+v-w)^2}{2\sfLp} + \frac{(2\sfLn+2\sfKp-2c-v+w+z)^2}{2\sfLn} \nonumber \\
 = & 2\sfKf + 2z + 2 (\sfKp - c) + \frac{(v-w)^2}{2\sfLp}  + \frac{(z+w-v+2\sfKp-2c)^2}{2\sfLn} 
 \nonumber \\
 \ge & 2\sfKf + 2 z + \frac{(w-v)^2}{4\sfLp} + \frac{z^2}{8\sfLn} \,.
\label{eq:quad-ineq}
\end{align}
Now,  setting $v:=(\tsfKfv)^\gak$,  or alternatively $x=a-2c+w+z$,  combining \eqref{eq:temp-bd} with \eqref{eq:quad-ineq},  results with
\begin{equation}\label{eq:riemann-sum}
\bar \bP^{(\star)}_{z} ( \wQs_{\npp}) 
\le \frac{C_{\sfKp}}{3 \sfKf} e^{-2\sfKf} e^{-2z} e^{-\frac{z^2}{8\sfKf}}
 \sum_{w} \frac{e^{-\frac{(w-v)^2}{4\sfLp}}}{\sqrt{4\pi \sfLp}} (w+z)  \,,
\end{equation}
where the sum is over all $w+z \in 2 \bZ_+$.   For $\sfLp \gg 1$,  the sum in \eqref{eq:riemann-sum} is
the Riemann sum approximation of
the expectation of $(W+z)_+$ for $W$ a $N(v,2\sfLp)$ variable.  Since for $v,z \ge 0$ the latter expectation 
is at most $v+z+\sqrt{2\sfLp} \le 2z$ (since $z \ge \sfKf^{\kapo} \ge \sqrt{2\sfLp} + (\tsfKfv)^\gak$),   
this completes the proof of \eqref{ubd-wQ}.

Turning to \eqref{ubd-Aell},  increasing $\Delta$ by two has no effect on the bound \eqref{ubd-wQ},  so 
we may and shall assume that  
$\ell \ge \sfLo+\Delta+3$.  Then,  \eqref{ubd-wQ} holds,  uniformly in $\ell$,  when replacing 
$\wQs_{\npp}$ with 
\[
\whcI_{\npp,\ell} \cap \{\eta_{\npp,\ell-i} > \gamma_{-1}(\ell-i), i=1,2 \}
\subset \wQs_{\npp} \,,
\] 
so by the Markov property of $j \mapsto \eta_{\npp,j}$,  the bound \eqref{ubd-Aell} follows upon showing 
that for $\sfKp$ large,  all $\ell \in [\sfLo+3,\sfLo+\sfLn]$ and $x^2, y^2 \in 2\bN$,
\begin{align}
 \bar \bP^{(\star)} \big( \eta_{\npp,\ell-i} > \gamma_{-1}(\ell-i), i=1,2,
\eta_{\npp,\ell}=y\,\big|\,
\eta_{\npp,\ell-3}=x\big)
& \ge \frac{1}{4}  \,  q_3 (x,y) \, \quad  \forall x \in \sfIp_{\ell-3},  y \in \sfIp_\ell \,.
\label{eq:drop2-bar}
\end{align}
With $\sfKf-\ell \ge \sfKp$ large,   as in our proof of \eqref{ubd-wQ} we can replace on 
both sides the point estimates at $y$,  by the interval $[y,y+2)$ and restrict to $y \in 2 \bN$.   
Having done so,  here $\gamma_{-1}(\ell-1) \le \gamma_{-1}(\ell-2) < x \wedge y- u$
for $u=(x \wedge y)^{\gak}$ (compare \eqref{dfn:gam} and \eqref{dfn:Ip-sfJs}),  with the 
conditional probability on the \abbr{lhs} 
thus being at least $q_3^{(u)}(x,[y,y+2))$.  With $\sfKf-\ell \ge \sfKp$ large,  having both
$x \in \sfIp_{\ell-3}$ and $y \in \sfIp_\ell$ implies that $|x-y| \le (x \wedge y)^{\okappa}$ and 
$u \ge u_o$ of part (a),  hence \eqref{eq:q3-bd} yields the desired conclusion.

\noindent
(c).  Aiming at \eqref{eq:mom-bd-cent},  we 
restrict \eqref{eq:new-p0} hereafter to $x^2,y^2 \in 2 \bZ_+$ such that 
$\frac{x-2 \tsfKfv}{\sqrt{\sfLn}} \in [1,6]$ and
$\frac{y - 2\sfKp}{\sqrt{\sfKp}} \in [1,2]$,  whereby 
$q_{\sfKp}(y,0) \ge c'_{\sfKp} >0$ (see \eqref{eq:miss-ball}).
Recall from \eqref{dfn:Ip-sfJs} and \Cref{dfn:Markov} that $q^{\Cap}_{\sfJus}(x,y)$ specifies having 
\begin{equation}\label{dfn:capq}
\sfco (\tsfKfv -\ell) + (\tsfKfv-\ell)^{\gak} \le \eta_{\npp,\sfLo+\ell} \le
\sfco (\tsfKfv -\ell) + (\tsfKfv-\ell)^{\okappa} \qquad \forall \ell \in [0,\sfLn] \,.
\end{equation}
In the sequel we bound $q^{\Cap}_{\sfJus}(x,y)$ below,  using \cite[Theorem 1.1]{BRZ} 
after showing that
\begin{align}\label{eq:fab-bd}
2 ( \tsfKfv-\ell) + (\tsfKfv-\ell)^{\gak} &\le f_{a,b}(\ell;\sfLn) +  \ell_{\sfLn}^{\gak}   \qquad \qquad \quad \forall \ell \in [0,\sfLn] \,,
\end{align}
where $f_{a,b}(\ell;\sfLn):=a+(b-a)\ell/\sfLn$
for $a:=2\tsfKfv+\tsfKfv^{\gak}$,  $b:=2\sfKp + \sfKp^{\gak} \gg 1$,  and further
(see \Cref{fig:EnlandFnl2}),  showing that for any $x,y$ as above,  
\begin{align}
 f_{x,y}(\ell;\sfLn) + \wt{C} \ell_{\sfLn}^{\kapp} &\le 2(\tsfKfv-\ell)+(\tsfKfv-\ell)^{\okappa}  \qquad \forall \ell \in [0,\sfLn] \,.
\label{eq:fxy-bd}
\end{align}
To get \eqref{eq:fab-bd}-\eqref{eq:fxy-bd},  first note that for $0<\gak<\frac{1}{2}<\okappa<1$ and all $\lambda,\varphi\in [0,1]$,  we have
\begin{align}\label{eq:fab-norm-bd}
(\varphi+\lambda)^{\gak} &\le \lambda (\varphi+1)^{\gak} + (1-\lambda) \varphi^{\gak} +  (\lambda \wedge (1-\lambda))^{\gak}\,,\\
\lambda + (1-\lambda) \varphi^{\okappa} &\le (\varphi+\lambda)^{\okappa} 
\label{eq:fxy-norm-bd}
\end{align}
(indeed \blue{\eqref{eq:fxy-norm-bd} holds by the concavity of $\lambda \mapsto (\varphi+\lambda)^{\okappa}$,
while} $\lambda ( (\varphi+1)^{\gak} - \varphi^{\gak}) \ge 0$ and 
$(\varphi+\lambda)^{\gak} \le \varphi^\gak + \lambda^\gak$ settle \eqref{eq:fab-norm-bd}
when $\lambda \le {1}/{2}$; otherwise increase its \abbr{lhs} to
$(\varphi+1)^{\gak}$ and note that $(1-\lambda) ( (\varphi+1)^{\gak} - \varphi^{\gak}) 
\le (1-\lambda) \le (1-\lambda)^{\gak})$.  Now,  fixing ${1}/{2} < \kapp < (1-\gak) \wedge \okappa$,  
we get \eqref{eq:fab-bd} after substituting $\varphi=\sfKp/\sfLn$,  $\lambda = 1-\ell/\sfLn$ in \eqref{eq:fab-norm-bd}. 
Next,  
\begin{figure}
\begin{center}
\begin{tikzpicture}[scale=0.5]

\begin{scope}[shift={(-5,0)}]
    \draw   (-5,0) -- (14.4,0) ;
     \draw   (-5,0) -- (-5,13.5) ;
    \draw  (0,0) -- (0,13.5) ;
  \draw [dotted] (13,0) -- (13,13);
    \draw[red] (0,6.8) .. controls (6,7) .. (13,1);
      \draw[-] (0,6.8) .. controls (6,5) .. (13,1);
    \draw[-] (0,11.5) .. controls (1,11) and (3,10.5).. (13,4);
  \draw[red] (0,9.7) .. controls (6,8.3) .. (13,2.5);
\draw[dotted,thick] (0,9.7) .. controls (3,7.3 ) and  (5,7.3) .. (6,7);
\draw[dotted,thick] (6,7) .. controls (8,5.5) and  (10,3.8) .. (13,2.5);
  \end{scope}
   \begin{scope}[shift={(-5,0)}]
\draw[blue] (0,5.6) -- (14,0);
 \node at (-5,-0.4) {\small $0$};
 \node at (0,-0.4) {\small $\sfLo$};
  \node at (5,-0.4) {\small $\sfLo+\ell$};
 \node at (14,-0.4) {\small $\sfKf$};
 \node at (12,-0.4) {\small $\sfLo+\sfLn$};
  \node at (-2.3,6.7) {\small $a=2\tsfKfv+\tsfKfv^\gak$};
  \node at (-1.6,11.6) {\small $2\tsfKfv+\tsfKfv^\okappa$};
   \node at (-.8,5.5) {\small $2\tsfKfv$};
   \node at (-.6,9.7) {\small $x$};
     \node at (13.2,2.5) {\small $y$};
         \node at (15.3,1.3) {\small $b=2\sfKp + \sfKp^{\gak}$};
  \end{scope}
\end{tikzpicture}
\end{center}
\caption{\small 
Depiction of $\sfIp_{ \sfLo+\ell}$ with \red{$f_{a,b}(\ell;\sfLn)+\ell_{\sfLn}^\gak$,  $f_{x,y}(\ell;\sfLn)+\wt{C}\ell_{\sfLn}^\kapp$} and a sample from $\whcI_\npp$ (dotted).}
\label{fig:EnlandFnl2}
\end{figure}
consider \eqref{eq:fxy-bd} for our largest
$x=2 \tsfKfv + 6 \sqrt{\sfLn}$,  $y=2 \sfKp + 2 \sqrt{\sfKp}$ and
assuming $\sfKp^{\okappa-\kapp} \ge \wt{C}+6$,  increase its \abbr{lhs} by noting that 
$\wt{C} \ell_{\sfLn}^{\kapp} \le \xi (\tsfKfv-\ell)^{\okappa}$ for $\xi := \wt{C} \sfKp^{\kapp-\okappa} \in (0,1)$. 
Subtract $\xi (\tsfKfv-\ell)^{\okappa}$ from both sides,  using $2 \sqrt{\sfKp} 
\le 
(1-\xi) \sfKp^{\okappa}$ and $6 \sqrt{\sfLn} 
\le (1-\xi) \sfLn^{\okappa}$ 
and some calculations to reduce  \eqref{eq:fxy-bd} to \eqref{eq:fxy-norm-bd} 
with $\varphi=\sfKp/\sfLn$,  $\lambda = 1-\ell/\sfLn$.

In view of \eqref{dfn:capq}-\eqref{eq:fxy-bd},  
\[
q^{\Cap}_{\sfJus}(x,y) \ge \bar \bP^{(\star)} \Big( \{ \eta_{\npp,\sfLo+\sfLn}=y \}
\bigcap_{\ell \in [0,\sfLn]} \{ 
 f_{a,b}(\ell;\sfLn) +  \ell_{\sfLn}^{\gak} 
\le \eta_{\npp,\sfLo+\ell} \le
 f_{x,y}(\ell;\sfLn) + \wt{C} \ell_{\sfLn}^{\kapp} \}
  \,  \big| \, \eta_{\npp,\sfLo}=x\Big) 
\]
and summing over the possible values of $y$,  we deduce from 
\cite[(1.2)]{BRZ} (at $\eta=3$, $C=\delta=1$, $\vep=\kapp-{1}/{2}$ and $\wt{C}$ as needed),
that for some $c>0$,
\begin{equation}\label{eq:BRZ-1.2}
\sum_y q^{\Cap}_{\sfJus}(x,y) \ge
6 \,  c \sum_y \frac{(1+x-a)(1+y-b)}{\sfLn}\sqrt{\frac{x}{y \sfLn}} e^{-\frac{(x-y)^2}{2\sfLn}} 
\ge \frac{c\,  z}{\sfKf} e^{-\frac{x^2}{2\sfLn}} 
\end{equation}
(indeed $\sfLn^{-1/2} \ge \sfKf^{-1/2} \ge {{z}/{(3 \sfKf)}}$ and 
for $\sfKp$ large,  our $x,y$ satisfy
$x \in [2\sfLn,3\sfLn]$,  $y -b \ge {\sqrt{y/3}}$,  $x -a \ge {\sqrt{\sfLn}/{2}}$).  
Similarly to \eqref{eq:q3-bd-alt},  note that $q_{\, \sfLo}(2\sfKf+z,x)= \bP(S_m = x^2/2)$
for the sum $S_m=S^{(\sfLo)}_m$ of $m:={(2\sfKf+z)^2/2}$ \abbr{iid} variables $\{\GJ'_i\}$
as in \eqref{dfn:GJp},  having here $q:={{1}/{(\sfLo+1)}}$.  It thus follows from 
\eqref{eq:new-p0} and the range of $x$ for which we derived \eqref{eq:BRZ-1.2} that 
\begin{equation}\label{eq:pre-pz}
\bar \bP^{(\star)}_z(\whcI_{\npp}) \ge c'_{\sfKp}  \frac{c z}{\sfKf}  
\bE \Big[ e^{-S_m/\sfLn}; \sqrt{\sfLn} \le \sqrt{2 S_m} - 2\tsfKfv \le 6 \sqrt{\sfLn} \,  \Big] \,.
\end{equation}
Now,  for all large $\sfKfv$ and $2 \sqrt{\sfKf} \le z \le 4\sqrt{\sfKf}$,
\[
\frac{m}{\sfLn+\sfLo} = \frac{(2\sfKf+z)^2}{2(\sfKf-\sfKp)} 
= 2 \sfKf + 2z + 2 \sfKp + \frac{(z+2\sfKp)^2}{2(\sfKf-\sfKp)}\,,
\]
with the right-most term uniformly bounded,  and in addition
\[
\Big[ \sqrt{2m - m^{3/4}},\sqrt{2m+m^{3/4}}\,\Big] \subset 
\Big[\sqrt{2m} - \frac{m^{1/4}}{\sqrt{2}},\sqrt{2m} + \frac{m^{1/4}}{\sqrt{2}}\Big]
\subset \big[2 \tsfKfv + \sqrt{\sfLn}, 2\tsfKfv + 6 \sqrt{\sfLn} \,  \big] \,.
\]
Thus,  in view of \eqref{eq:pre-pz},  we get \eqref{eq:mom-bd-cent} upon showing that for some  
$c>0$ and all $\sfKfv$ large,
\begin{equation}\label{eq:pz-alt}
\bE \Big[ e^{-S_m/\sfLn}; |S_m - m| \le  \frac{1}{2} m^{3/4} \Big] \ge c \,  e^{-m/(\sfLn+\sfLo)} \,.
\end{equation}
Turning to \blue{show}  \eqref{eq:pz-alt},  with $\Lambda_q(\lambda)$ as in \eqref{def:q-cgf} and 
$\wt{S}_m$ the partial sum of the \abbr{iid} variables $\{\wt{\GJ}_i\}$ such that 
\[
\bP(\wt{\GJ}_i=k) = e^{-k/\sfLn} e^{-\Lambda_q(1/\sfLn)} \bP(\GJ'_i = k),  \qquad \forall k \ge 0 \,,
\]
{we have 
that}
\begin{equation}\label{eq:rnd}
\bE \Big[ e^{-S_m/\sfLn}; |S_m - m| \le  \frac{1}{2} m^{3/4} \Big] = 
\bP\Big( |\wt{S}_m - m| \le  \frac{1}{2} m^{3/4} \Big) e^{m \Lambda_q(1/\sfLn)} 
\end{equation}
With $\wt{\sfL}:=(1-e^{-1/\sfLn})^{-1} \ge \sfLn$,  substituting $q^{-1}=\sfLo+1$
in \eqref{def:q-cgf},  we find that for $\sfLn,\sfLo \ge 1$, 
\[
m \Lambda_q(1/\sfLn) = m \log \Big( 1 - \frac{1}{\wt{\sfL}+ \sfLo} \Big)  \ge  m \log \Big( 1 - \frac{1}{\sfLn + \sfLo} \Big) \ge 
-\frac{m}{\sfLn+\sfLo} - \frac{m}{(\sfLn+\sfLo)^2} \,.
\]
The right-most term is at least $-8$ when $\sfKf \ge 4+ 8 \sfKp$,  so by 
\eqref{eq:rnd} it suffices for \eqref{eq:pz-alt} to show that 
\[
\bP\Big( |\wt{S}_m - m| >  \frac{1}{2} m^{3/4} \Big) \le \frac{1}{2} \,,
\]
which by Markov's inequality is a consequence of having
\begin{equation}\label{eq:var-mean-wtG}
{\rm var}[\wt{\GJ}_i]  
= \Lambda_q''(\frac{1}{\sfLn}) 
= O(\sfLo) 
\le \frac{1}{16} \sqrt{m} \quad \text{and} \quad 
1- \bE [\wt{\GJ}_i] = \Lambda'_q(1/\sfLn) - \Lambda_q'(0) = O(\sfLo/\sfLn) 
\le \frac{1}{4} m^{-1/4} \,.
\end{equation}
(d).  With $I_x:=[x,x+1)$,  we have that 
\begin{equation}\label{eq:wQs-off}
\bar \bP^{(\star)}_{z} ( \wQs_{\sfv}) = \sum_{x=0}^\infty
\bar \bP^{(\star)}_{z} ( \wQs_{\sfv},  \eta_{\sfv,\sfLo} \in I_x) 
:= \sum_{x=0}^{\infty} p_\sfv(z,x).
\end{equation} 
By rotation symmetry, 
$p_{\sfv}(z,x) = p_{\rho}(z,x)$ depends only on 
$\rho:= \|\sfv\| \le \sfKfv^{-5}$.  
We are in the unconditional setting of \Cref{rmk:off-center-applies} for $R_1=R=e^{-2}$ and 
$\sfs=m=\bar \sfs_{\sfKf}(z) \le 3 \sfKf^2$ independent planar Brownian excursions.
The event we consider 
$\Gamma:=\wQs_{\sfv} \cap \{ \eta_{\sfv,\sfLo} \in I_x \}$ is in $\NN_{\sfv,r,m,k}$ 
for $r=e^{-(\sfLo+3)}$ and $k=\frac{1}{2} (x+1)^2$.  
We thus have by 
\eqref{eq:off-center-bd},  that for all large enough $\sfKfv$ and $x \le \sfKfv^{2}$ 
the joint probability $\PRso (\Gamma)$ 
is at most $2 (\bP_R^{1})^m(\Gamma)$
and hence $p_\rho(z,x) \le 2 p_0(z,x)$ for any such $x$.
Further,  the probability of hitting a ball of radius $r_{\sfKf+1}$ 
during one Brownian excursion between concentric balls of radii 
$r_{\sfLo}$ and $r_{\sfLo-1}$ is,  uniformly in the starting and ending points of such excursion,  
at least $2c/\sfKfv$ for some $c=c(\Delta)>0$.
The probability that at least ${ (\sfKfv^2)^2/2}$ such excursions miss this small ball
is thus at most $e^{-c \,  \sfKfv^3}$,  resulting in view of \eqref{eq:wQs-off} with 
\begin{equation}\label{eq:bd-hQv}
\bar \bP^{(\star)}_{z} ( \wQs_{\sfv}) \le  2 \sum_{x=0}^{\infty} p_0 (z,x) + e^{-c \,  \sfKfv^3} = 
2 \bar \bP^{(\star)}_{z} ( \wQs_{\npp}) + e^{-c \,  \sfKfv^3} \,,
\end{equation}
thereby extending the validity of \eqref{ubd-wQ} to all $\|\sfv\| \le \sfKfv^{-5}$.  Applying such 
reasoning with $\whcI_{\sfv,\ell}$ instead of $\wQs_{\sfv}$
extends \eqref{ubd-Aell} from $\whcI_{\npp,\ell}$ 
to $\whcI_{\sfv,\ell}$ (uniformly in $\sfv$ and $\ell \in \sfJus$),  while considering 
instead $\whcI_{\sfv}$ shows that the much larger than 
$e^{-c \sfKfv^3}$ lower bound of \eqref{eq:mom-bd-cent} extends to  
$\whcI_{\sfv}$,  uniformly over $\|\sfv\| \le \sfKfv^{-5}$.
\end{proof}

\begin{proof}[Proof of \Cref{lem:g4-bds}]
Here $\|\sfv\| \le r_{\sfLw} \le \sfKfv^{-5}$,  so with $|\bF^\npp| = |\bF^{\star}| \le e^{2\sfKfv}/(2 \wt{c}^2)$,  it
follows from \eqref{dfn:g4-alt},  \eqref{eq:A-RU-Q-hatQ},
the union bound and \Cref{lem:g4-cnc},  that for any $z \in \sfIs_4$
\[
\gfo(z) \le  \frac{e^{2 \sfKf}}{2 \, \wt{c}^2}  \max\big\{ \bar \bP^{(\star)}_{z} ( \wQs_{\sfv}) \,:\,
\| \sfv \|  \le \sfKfv^{-5} \big\} 
\le   C(\sfKp) \frac{z}{\sfKf} e^{- 2 z}  e^{-\frac{z^2}{8 \sfKf}} 
\]
as stated in \eqref{ubd2-g4}.  

Next,  consider the sub-lattice $\bF' \subset \bF^{\npp}$ obtained when up-scaling the cell-size parameter $\wt{c}$ of
$\bF^\star$,  see \eqref{dfn:F-star},  to be $\wt{c} := \sqrt{2} e e^{\bst}$ for some $\bst  \ge \sfKp + 6 \log \sfKp$,
denoting by $N_z$ the number of events $\{\whcI_{\sfv},  \sfv \in \bF'\}$ of \eqref{def:whcI} that occur 
under $\bar \bP^{(\star)}_z$ when defining $\sfJus$ with $\Delta=0$.  
Recall \eqref{eq:A-RU-Q-hatQ} that $\whcI_{\sfv}\subset \QQs$ so fixing $z \in \sfIs_4$ 
with $z \in [2 \sqrt{\sfKf},4 \sqrt{\sfKf}]$,  it suffices for \eqref{eq:lbd-alpha},  see \eqref{dfn:g4-alt},  to show that both
\begin{equation}\label{eq:lbd-Av}
\bP(N_z \ge 1) \ge c(\sfKp) \frac{z}{\sfKf} e^{-2 \sfLw} e^{-2 z} 
\end{equation}
and for some $c>0$, 
\begin{equation}\label{eq:ubd-eta-L}
\bar \bP^{(\star)}_{z} (\eta_{\npp,\sfLw} \ge 3 \sfKfv) \le e^{-c \,  \sfKfv} 
\end{equation}
(uniformly in $\sfKfv$ and $z$ as above).  To this end,  first note that the \abbr{lhs} of \eqref{eq:ubd-eta-L} is 
precisely the probability that the sum $S_m^{(\sfLw)}$ of $m \le 3 \sfKfv^2$ \abbr{iid} variables $\GJ'_i$
of the form \eqref{dfn:GJp} at $q={{1}/{(\sfLw+1)}}$,  exceeds ${9 \sfKfv^2/2}$.  Now, 
from  \eqref{def:q-cgf} we get by elementary algebra that for any $\sfLw \ge 1$ and 
$0 \le \lambda \le {1}/(2 \sfLw)$,  
\[
\log \bE\Big[e^{\lambda \big(S_m^{(\sfLw)}-m\big)}\Big] 
= m \big( \Lambda_q(-\lambda) - \lambda \big) \le 3 \lambda^2 \sfLw  m \,.
\]
The probability in \eqref{eq:ubd-eta-L} is thus at most $\exp(-{c\sfKfv^2/{\sfLw}})$ for some $c>0$
($\lambda={1/(12 \sfLw)}$ gives $c={{1}/{16}}$),   as claimed.  

Turning to {prove} \eqref{eq:lbd-Av},  compare $\whcI_{\sfv,\ell}$ of \eqref{def:Aell} with 
\begin{align}\label{def:A-ell-ellp} 
\wtcI_{\sfv,\ell'} &:= \{ \eta_{\sfv,\ell'} \in \sfIp_{\ell'},  \eta_{\sfv,\sfKf} = 0 \} \,,   \qquad \ell':=\ell+5 \log (\sfKf-\ell) \,, 
\end{align}
noting that 
\[
\ell \in [\sfLw+2,\sfKf-\bst] \quad \Longrightarrow  \quad \ell' \in \sfJus,  \quad \text{hence}
\quad \whcI_{\sfv} \subset \whcI_{\sfv,\ell} \subset  \wtcI_{\sfv,\ell'}  \quad \forall \sfv \in \bF' \,.
\]
Further,  by construction
\begin{equation}\label{def:Gell}
\bG'_{\ell}  :=  \left\{ (\sfv,\sfv') \in (\bF')^2 : 0 < 2 r_{\ell-1} \le \sfd(\sfv,\sfv') < 2 r_{\ell-2}\right\}
\end{equation}
is empty except for $\ell \in [\sfLw +2, \sfKf-\bst]$,  so we establish 
\eqref{eq:lbd-Av} by utilizing the bound
\begin{align}
\bP(N_z \ge 1) &\ge \;\;\;\;  \bE[N_z] \quad  \; - \;  \bE[N_z (N_z-1)] \nonumber \\
 \ge |\bF'| & \inf_{\sfv \in  \bF'} \{ \bar \bP^{(\star)}_{z} (\whcI_{\sfv}) \}
- \sum_{\ell=\sfLw+2}^{\sfKf-\bst} |\bG'_\ell| \sup_{(\sfv,\sfv') \in \bG_\ell} \{ \bar \bP^{(\star)}_z 
(\wtcI_{\sfv,\ell'} \cap \whcI_{\sfv',\ell}) \} 
\label{eq:lbd-Ny}
\end{align}
(with $\bG_\ell$ as in \eqref{def:Gell},  except for considering all $\sfv,\sfv' \in \bF^\npp$).
Indeed,  in the sequel we show that for large $\sfKp$,  some $C_1(\sfKp)<\infty$,  all $\sfKfv \gg \sfKp$ 
and $z$,  $\ell$ as above,  
\begin{align}\label{eq:2mom-ubd}
e^{2 (\sfKf-\ell)} 
\sup_{(\sfv,  \sfv') \in \bG_\ell} \big\{ \bar \bP_z^{(\star)} (\whcI_{\sfv',\ell} \cap  \wtcI_{\sfv,\ell'} ) \big\} & \le 
C_1(\sfKp) \frac{z}{\sfKf} e^{-2 \sfKf} e^{-2 z}  e^{-\frac{1}{4}(\sfKf-\ell)^{\gak}}\,.
\end{align}
Further,  $|\bF'| \ge c\,  e^{2 (\sfKfv-\bst)}$ and $|\bG'_{\ell}| \le c^{-1} e^{2 (\sfKfv-\bst)} e^{2 (\sfKf-\ell-\bst)}$ 
for some $c>0$ (with at most $c' e^{2 (\sfKf-\bst)}$ points of $\bF'$ per unit area,  clearly 
$\sup_{\sfv} | \bF' \cap \bB(\sfv,  2 r_{\ell-2})| \le c^{-1}  e^{2 (\sfKf- \bst)} r_{\ell-2}^{2}$).  
Plugging in \eqref{eq:lbd-Ny} those estimates, 
the lower bound \eqref{eq:mom-bd-cent} (which holds throuhgout $\bF'$,  see \Cref{lem:g4-cnc}(d)),  
and the upper bound of \eqref{eq:2mom-ubd},  we get \eqref{eq:lbd-Av} 
upon considering $\bst$ large enough.

Turning to show \eqref{eq:2mom-ubd},  recalling the $\sigma$-algebra $\GG^{\sfv}_\ell$ 
which excludes only the 
excursions from $\partial \bB(\sfv,r_\ell)$ to $\partial\bB (\sfv,r_{\ell-1})$ (see \Cref{def:Fvjs}),
we claim that for some $C<\infty$,  any relevant $\ell$ and $\sfv \in \bF^\npp$, 
\begin{equation}\label{eq:Jay-alpha}
\bP(\wtcI_{\sfv,\ell'} \, |\, \GG^{\sfv}_\ell) \le C e^{-2 (\sfKf-\ell')} e^{-(\sfKf-\ell')^{\gak}} 
\le C e^{-2 (\sfKf-\ell)} e^{-\frac{1}{2}(\sfKf-\ell)^{\gak}} \,.
\end{equation}
Indeed,  for $\sfKp \ge c_{\gak}$ the second inequality is a trivial consequence of having 
$\ell'-\ell=5 \log (\sfKf-\ell)$ and $\sfKf -\ell \ge \sfKf-\ell' \ge \sfKp$.  For the first inequality,  
we set $T_{\ell'}^{\sfv,\ell,\wt{m}}$
for the number of traversals $\partial \bB (\sfv, r_{\ell'})\to\partial \bB (\sfv,r_{\ell'-1})$ 
accumulated within the first $\wt{m}$ such traversals $\partial \bB (\sfv, r_{\ell-1})\to\partial \bB (\sfv,r_{\ell})$,
with $T_{\sfKf}^{\sfv,\ell',m'}$ defined analogously (see also \Cref{def-trav}).   Equipped with this 
notation and that of \eqref{def:why},  we re-write the events {of \eqref{def:A-ell-ellp}} as
\begin{align*}
\wtcI_{\sfv,\ell'}
& =
\! \bigcup_{\wt{m}} \! \bigcup_{m' \in \sfII_{\ell'}} \! \{ \trav{m}{\npp}{e^{-1}}{e^{-2}}{\sfv}{r_{\ell}}{r_{\ell-1}}
 = \wt{m}, 
  T_{\ell'}^{\sfv,\ell,\wt{m}} = m',   T_{\sfKf}^{\sfv,\ell',m'} = 0  \},
\end{align*} 
where $\sfII_j :=\{ m' \in \bN: \sqrt{2m'} \in \sfIp_j \}$.  Here
$(1+ 50 (\ell'-\ell) e^{-(\ell'-\ell)})^{m'}
\le C$ for some $C<\infty$ and all $m' \in \sfII_{\ell'}$ (by our choice of $\ell$,  $\ell'$).  Noting that 
$\trav{m}{\npp}{e^{-1}}{e^{-2}}{\sfv}{r_{\ell}}{r_{\ell-1}}$ and the start and end points of all excursions from 
$\partial \bB(\sfv,r_{\ell})$ to $\partial \bB(\sfv,r_{\ell-1})$,  are 
measurable on $\GG^{\sfv}_\ell$,  we thus deduce in view of Lemma \ref{recursionends} that
\begin{align*}
\bP(\wtcI_{\sfv,\ell'} \, |\, \GG^{\sfv}_\ell)
& \le C 
\sum_{\wt{m}} {\bf 1}_{\{\trav{m}{\npp}{e^{-1}}{e^{-2}}{\sfv}{r_{\ell}}{r_{\ell-1}}=\wt{m}\}} 
\sum_{m' \in \sfII_{\ell'}} \bP(T_{\ell'}^{\sfv,\ell,\wt{m}} = m'  \, |\, \GG^{\sfv}_\ell) \bP(T_{\sfKf}^{\sfv,\ell',m'} = 0 ) \\
&\le C \sup_{m' \in \sfII_{\ell'}} \{ \bP(T_{\sfKf}^{\sfv,\ell',m'} = 0 ) \} = C \sup_{m' \in \sfII_{\ell'}}
\big(1-\frac{1}{\sfKf-\ell'+1}\big)^{m'}
\end{align*}
\blue{(using also that $\{\trav{m}{\npp}{e^{-1}}{e^{-2}}{\sfv}{r_{\ell}}{r_{\ell-1}}=\wt{m}, 
T_{\ell'}^{\sfv,\ell,\wt{m}} = m' \}_{\wt{m},m'}$ are disjoint events)}. 
Since any $m' \in \sfII_{\ell'}$ is at least $\frac{1}{2} (2 (\sfKf-\ell') + (\sfKf-\ell')^{\gak})^2$,  
the first inequality of \eqref{eq:Jay-alpha} follows by elementary algebra.
Further,  for any $(\sfv,\sfv') \in \bG_\ell$, 
\[
\bB(\sfv,r_{\ell-1}) \subset \bB(\sfv',r_{\ell-3}) 
\setminus \bB(\sfv',r_{\ell-1}) \subset \bB(\npp,1),
\] 
hence $\whcI_{\sfv',\ell}$ is measurable on $\GG^{\sfv}_{\ell}$.  In particular,  
\[
\bar \bP_z^{(\star)} (\whcI_{\sfv',\ell} \cap  \wtcI_{\sfv,\ell'} ) = 
\bar \bE_z^{(\star)} \big[ \bP( \wtcI_{\sfv,\ell'} | \GG^{\sfv}_{\ell} ) ; \whcI_{\sfv',\ell} \big] \le 
\| \bP( \wtcI_{\sfv,\ell'} | \GG^{\sfv}_{\ell} ) \|_\infty \bar \bP_z^{(\star)} (\whcI_{\sfv',\ell} )
\]
and we get \eqref{eq:2mom-ubd} by combining \eqref{eq:Jay-alpha} with 
our uniform bound on $\bar \bP_z^{(\star)} (\whcI_{\sfv',\ell} )$ (see \Cref{lem:g4-cnc}(d)).
\end{proof}


\subsection{Controlling quantiles of 
empirical measures}\label{subsec:quant}

The \blue{Kantorovich-Wasserstein} $W_1$-distance 
\begin{align}
\sfd_\sfw(\mu_1,\mu_2) & := \inf_{\stackrel{\xi \in M_1(\Ipi^2)}{\xi_1= \mu_1,\xi_2=\mu_2}} \Big\{ \int_{\Ipi^2} |\th-\wt{\th}| {\rm d} \xi (\th,\wt{\th}) \Big\} 
= \sup \{| \int_{-\pi}^{\pi} h {\rm d} \mu_1 - \int_{-\pi}^{\pi} h {\rm d} \mu_2| : \|h\|_{\sf Lip} \le 1 \}  
 \label{dfn:dw-orig}
\end{align}
\blue{on $M_1(\Ipi)$}
(see \cite[Def. 2.1 \& Thm. 2.5]{BL19}),
which is used by \cite{BRZ2},  is dominated by \blue{the 
Kantorovich-Wasserstein
$W_\infty$-distance}
$\sfd_\infty(\cdot,\cdot)$ of \eqref{dfn:dsup}. 
Indeed,  we have for any $\mu_i \in M_1(\Ipi)$, 
\begin{equation}
\label{dfn:dw}
\sfd_\sfw(\mu_1,\mu_2) = \| \mu_1([-\pi,\theta])-\mu_2([-\pi,\theta]) \|_1 
 = \int_0^1 | \mu_1^{-1} (t) - \mu_2^{-1}(t) | {\rm d} t \le \sfd_\infty(\mu_1,\mu_2).
\end{equation}
(see \cite[Thm.  2.9 \& 2.10]{BL19}).
We proceed with 
certain additional useful facts about 
\blue{the $W_\infty$-distance}.
\begin{lemma}\label{lem:angle-diff} 
Fix $\mu \in M_1(\Ipi)$ such that $\frac{d \mu}{d\theta} \ge c^{-1}$ for some $c$ finite.  Then
for any $\nu \in M_1(\Ipi)$,
\begin{equation}\label{eq:quantiles}
\sfd_\infty(\nu,\mu) \le c\,  \sfd_{\rm KS}(\nu,\mu),  \quad {\rm where}
 \quad \sfd_{\rm KS}(\nu,\mu):= \|\nu([-\pi,\cdot])-\mu([-\pi,\cdot])\|_\infty 
\end{equation}
\blue{denotes the Kolmogorov-Smirnov distance}.  In particular,  \blue{for discrete quantile measures}
\begin{equation}\label{eq:ds-emp}
L_m^\mu := \frac{1}{m} \sum_{q=1}^m \delta_{\mu^{-1}(q/m)} \in \bB^{\infty}_{\mu}(c/m) ,
\end{equation}
and we have that for the empirical measure $L_{\mu}(\sfs)$ of $\sfs \ge 1$ \abbr{iid} samples from $\mu$,
\begin{equation}\label{lbl:d-inf-bd}
\bP \big( \sfd_\infty (\mu, L_{\mu}(\sfs)) > \frac{c u}{\sqrt{\sfs}} \big) \le 2 e^{-2 u^2} \,, 
\qquad \forall u > 0 \,.
\end{equation}
Moreover,  for such $\mu$,  any $\boeta,\wt{\boeta} \ge 0$,  $\wt{m} -m  \in [0, \wt{\boeta} \,  m]$ and
$\{\th_q,  \blue{1 \le q \le \wt{m}} \} \subset \Ipi$, 
\begin{equation}\label{eq:dinf-comp}
L_{m} := \frac{1}{m} \sum_{q=1}^m \delta_{\theta_q} 
 \in  \bB^{\infty}_{\mu}(\boeta)
\;\; \Longrightarrow \;\; L_{\wt{m}} \in \bB^{\infty}_{\mu}\big(\boeta+ c \,  \wt{\boeta} \big) \,.
\end{equation}
\end{lemma}
\begin{proof} With $\frac{d \mu}{d\theta} \ge c^{-1}$,  setting $h =c \,  \sfd_{\rm KS}(\nu,\mu)$ yields by the definition 
of $\sfd_{\rm KS}(\cdot)$,  that for any $\theta \in [-\pi,\pi)$,
\begin{align*}
 \mu([-\pi,\theta-h])  
 \le  \mu([-\pi,\theta]) - \sfd_{\rm KS}(\nu,\mu)  \le \nu([-\pi,\theta]) 
 &\le  \mu([-\pi,\theta]) + \sfd_{\rm KS}(\nu,\mu) 
 \le \mu([-\pi,\theta+h]) \,.
\end{align*}
This implies that $\sfd_\infty(\nu,\mu) \le h$  (see  \cite[Thm. 2.12]{BL19}),  as stated in \eqref{eq:quantiles}.
By definition,  $\sfd_{\rm KS}(L_m^\mu,\mu) \le \frac{1}{m}$,  so plugging 
$\nu = L_m^\mu$ in \eqref{eq:quantiles} results with \eqref{eq:ds-emp}.  Moreover, 
applying \eqref{eq:quantiles} for $\nu=L_{\mu}(\sfs)$,  we see that \eqref{lbl:d-inf-bd} 
follows from the Dvoretsky-Kiefer-Wolfowitz inequality
\begin{equation}\label{eq:dkw}
\bP \big( \sfd_{\rm KS} (\mu, L_{\mu}(\sfs)) > \frac{u}{\sqrt{\sfs}} \big) \le 2 e^{-2 u^2} 
\qquad \forall u > 0 
\end{equation}
(see \cite{Massart} for a proof of \eqref{eq:dkw}). 

Next,  with $\wt{m}\geq m$, let $\theta_{(i)}$, $i=1,\ldots,m$  and $\wt{\th}_{(i)}$, $i=1,\ldots,\wt{m}$
denote the order statistics of
$\{\theta_q\}_{q\in [m]}$ and  $\{\theta_q\}_{q\in[\wt{m}]}$,  that is $\{\theta_q\}_{q\in[m]}=\{\th_{(q)}\}_{q\in [m]}$
and
 $\{\theta_q\}_{q\in[\wt{m}]}=\{\wt{\th}_{(q)}\}_{q\in [\wt{m}]}$,  with
\[
 \{ \th_{(0)} \le \cdots \le \th_{(m)} \} \,,  \quad
\{ \wt{\th}_{(0)} \le  \cdots \le \wt{\th}_{(\wt{m})} \} \,,
\]
where for convenience we set  $\th_{(0)}=\wt{\th}_{(0)}=-\pi$.  \blue{Adopting the convention 
of $\th_{(q)}=\pm \,  \pi$ when $q>m$ or $q<0$,  respectively,  necessarily}
$\th_{(q-(\wt{m}-m))} \le \wt{\th}_{(q)} \le \th_{(q)}$ for all $q \in [\wt{m}]$.  Taking
$q = [t \wt{m}]$ and noting that $1 \le \frac{\wt{m}}{m} \le  1+ \wt{\boeta}$,  \blue{this} implies that for all $t \in [0,1]$,
\begin{equation}\label{eq:order-th}
(L_{m})^{-1}(t-\wt{\boeta}) \le (L_m)^{-1}(\frac{q}{m} - \wt{\boeta}) \le (L_{\wt m})^{-1}(t) \le (L_m)^{-1}(\frac{q}{m})
\le (L_{m})^{-1}(t+\wt{\boeta}) \,.
\end{equation}
Since $\| \mu^{-1} \|_{\rm Lip} \le c$,  see \cite[\blue{display following} (5.7) \& Cor.5.12(d)]{BL19},  
combine \eqref{eq:order-th} with our assumption \blue{in \eqref{eq:dinf-comp}} that $\sfd_\infty(L_{m},\mu) \le \boeta$,  to get
in view of  \eqref{dfn:dsup},  that for all $t \in [0,1]$
\[
-\boeta - c\,  \wt{\boeta} + \mu^{-1} (t) \le 
-\boeta + \mu^{-1}(t-\wt{\boeta}) \le (L_{\wt m})^{-1}(t) \le \boeta + \mu^{-1}(t+\wt{\boeta}) \le 
\boeta + c \,  \wt{\boeta} + \mu^{-1}(t)
\,,
\]
as claimed in \eqref{eq:dinf-comp}.
\end{proof}

Our next lemma is utilized in \Cref{sec:good-cover} for deducing from the control gained 
by \eqref{lbl:d-inf-bd} on angle-differences for two independent Brownian excursions,
a uniform total-variation bound between their first exit points on an intermediate ball.
\begin{lemma}\label{lem-coupb}
For some $C<\infty$,  any $r<e^{-1}$,  $|\sfu|=e r$,  $\sfz'=r e^{i \theta'}$ and 
$\sfz = r e^{i \theta} \in \bB_\sfd(\sfz',r/(2e))$,
\[
\sfd_{\rm TV}(\zeta_{\sfz},\zeta_{\sfz'}) \le \frac{C}{r} |\sfz-\sfz'| \le C |\theta - \theta'| \,,
\]
where $\zeta_{\sfz}$ denotes the first hitting point on $\partial \bB_\sfd(\sfz',r/e)$ 
by a Brownian path started at $\sfz$ and conditioned to exit $\bB_\sfd({\bf 0},e \, r)$ 
at the given point $\sfu$ (with $\zeta_{\sfz'}$ corresponding to $\sfz=\sfz'$).
\end{lemma}
\begin{proof} The right-inequality being trivial,
we only prove the left one.  Further,
being confined to $\bB_\sfd({\bf 0},1)$, we can replace \abbr{wlog} the unit torus 
$\bT^2$ by $\bR^2$ and thereby, using Brownian scaling, fix \abbr{wlog} $r=1$. 
Having done so, we abbreviate $\bB:=\bB(\sfz',1/e)$ and
let $f_{\sfz}(\cdot)$ and $f_{\sfz'}(\cdot)$ denote the normalized densities 
of $\zeta_{\sfz}$ and $\zeta_{\sfz'}$ with respect to a unit surface measure on $\partial \bB$. 
Denoting by $q$ the Poisson kernel for the exit density of $\bB$ 
and by $p$ the corresponding kernel for
exit density of $\bB({\bf 0}, e)$, we have that   
\begin{equation}\label{10.11}
f_{\sfz'}(\sfx)=\frac{q (\sfz',\sfx)\, p (\sfx,\sfu)}{p (\sfz',\sfu)}\,, \qquad
f_{\sfz}(\sfx)=\frac{q(\sfz,\sfx) p (\sfx,\sfu)}{p(\sfz,\sfu)} \,.
\end{equation}
From these identities we get that
\begin{align}\label{10.12new}
 \sfd_{\rm TV}(\zeta_{\sfz},\zeta_{\sfz'})
 &= \frac{1}{2} \int_{\partial \bB} | f_{\sfz'}(\sfx) - f_{\sfz}(\sfx) | d\sfx \nonumber \\
& \le \frac{1}{2} \int_{\partial \bB} \big|  \frac{p (\sfz',\sfu)}{p (\sfz,\sfu)} - 1 \big| f_{\sfz'} (\sfx) d\sfx 
+ \int_{\partial \bB} \big| 1 - \frac{q(\sfz',\sfx)}{q(\sfz,\sfx)} \big|  f_{\sfz} (\sfx) d\sfx \nonumber \\
 & \le  \frac{1}{2} \Big|  \frac{p (\sfz',\sfu)}{p (\sfz,\sfu)} - 1 \Big| + \frac{1}{2} \sup_{\sfx \in \partial \bB}
\Big| 1 - \frac{q(\sfz',\sfx)}{q(\sfz,\sfx)} \Big| \,.
\end{align}
Now, since $|\sfz|=|\sfz'|$, we have from \eqref{dfn:Poiss} that 
\begin{equation}
\frac{p (\sfz',\sfu)}{p (\sfz,\sfu)} - 1 = \frac{ |\sfz-\sfu|^{2}}{|\sfz'-\sfu|^{2}} - 1 =
2 \Re\Big\{ \frac{\sfz-\sfz'}{|\sfz'-\sfu|}\frac{\bar\sfz'-\bar\sfu}{|\sfz'-\sfu|} \Big\}
+ \frac{ |\sfz-\sfz'|^{2}}{|\sfz'-\sfu|^{2}} \,.
 \label{10.12r2}
\end{equation}
Further, as $|\sfz-\sfz'| \le 2|\sfz'-\sfu|/(e-1)$ and $|\sfz'-\sfu| \ge e-1$,
we deduce from \eqref{10.12r2} that the 
first term on the right hand side of \eqref{10.12new} is at most
\[
\big(1 + \frac{1}{e-1} \big) \frac{|\sfz-\sfz'|}{|\sfz'-\sfu|} \le C_1 |\sfz-\sfz'| 
\]
(with $C_1
= e/(e-1)^2$). Next, note that by symmetry $q(\sfz',\cdot)=1$, with
\begin{equation}
1- \frac{1}{q (\sfz,\sfx)} = 1 - \frac{|\sfz-\sfx|^{2}}{|\sfz'-\sfx|^{2}-|\sfz'-\sfz|^{2}}
= \frac{2\Re\{(\sfz'-\sfz)(\bar \sfz - \bar \sfx)\}}{|\sfz'-\sfx|^{2}-|\sfz'-\sfz|^{2}} \,.
\label{10.12r4}
 \end{equation}
Since $|\sfz'-\sfz| \le |\sfz'-\sfx|/2=1/(2e)$, the second term on the right hand side of 
\eqref{10.12new} is bounded by 
\[
\frac{4}{3} \frac{|\sfz'-\sfz| |\sfz - \sfx|}{|\sfz'-\sfx|^2} \le 2\frac{|\sfz'-\sfz|}{|\sfz'-\sfx|} \le 
2 e  |\sfz-\sfz'|  \,,
\] 
 yielding the stated result.
\end{proof}

\subsection{Coupling \abbr{iid} variables given their sum}\label{subsec:cond-exc}

Suppose $\bZ_+$-valued variables $L^{(k)}$ satisfy 
for some $\mu \ge 1$, $c>0$ and all $k \ge 1$,
\begin{align}\label{eq:mean-var}
\bE L^{(k)} &=\mu \,, \quad \inf_k \{ k^{-1} \sfv_k \} \ge c  \quad {\rm for } \quad \sfv_k := {\sf var} (L^{(k)})\,, \\
\label{eq:Ber-Geo}
\bP(L^{(k)} \ne 0) &\le \frac{1}{ck}\,, \quad
\sup_k \Big\{ \bE\Big[\exp(c L^{(k)}/k) \,\big|\, L^{(k)} \ne 0 \Big] \Big\} < \infty \,,\\
\label{eq:qi-spread} 
\liminf_{k \to \infty} & \Big\{ k^{1/c} \inf_{i \ge 1} \{ \mu - i q_i^{(k)} \} \Big\}  > 0  
\quad {\rm for } \quad q_i^{(k)} := \bP( L^{(k)} = i) \,.
\end{align}
Restricting hereafter to $\mu m + b \in \bZ_+$, denote by $\bQ^{(k)}_{m,b}$ the 
joint law of \abbr{iid} copies $\{L_j^{(k)}, j \le m\}$ of $L^{(k)}$,  conditional on the event 
$\{S_m^{(k)}=\mu m +b\}$, where
\begin{equation}\label{dfn:Sl}
S^{(k)}_\ell := \sum_{i=1}^\ell L^{(k)}_i \,.
\end{equation}
We have the following coupling among  $\bQ^{(k)}_{m,b}$ with ``not too far'' values of $(m,b)$.
\begin{lemma} \label{lem:emp}
For any $\ka>0$ let
\begin{align}\label{dfn:Thetak}
\Theta_\ka := \big\{ m,k,r \in \bZ_+, b \in \bR : \; & \;  
k \ge \ka^{-3} \,, \;
\frac{|b|}{\sqrt{m}} \le  \frac{k \vee r}{\ka} \le \frac{\log m}{\ka^2}   \big\} \,.
\end{align}
Consider $\{ \wt{L}_j,  j \le m\}$ and $\{\wh{L}_j, j \le \wh{m}\}$ 
of laws $\bQ^{(k)}_{m,b}$ and $\bQ^{(k)}_{\wh{m},\wh{b}}$, respectively
(for the same unconditional marginal law of $L^{(k)}$ satisfying \eqref{eq:mean-var}--\eqref{eq:qi-spread}).
Extending with zeros the shorter among these two sequences, 
there exists a coupling of $\{ \wt{L}_j\}$ and $\{\wh{L}_j\}$
such that for small $\ka>0$ and any $C_1$ finite,
\begin{align}\label{cond-pre-a}
\sup_{\Theta_\ka} \sup_{(\wh{m},\wh{b}) \in \Gamma} \, \Big\{ \frac{e^r}{k} 
\bP \big( \, \bar{\Delta}_{(m,b),(\wh{m},\wh{b})} \ge 5  r \ka \sqrt{m}\, \big) \Big\} < \infty \,, 
\end{align}
where
\begin{align}
\!\!\!\Gamma\!=\!\Gamma^{\blue{(\mu)}}_{m,k,r,b} \!:= \!\{ (\wh{m},\wh{b}) : \wh{m}, \mu \wh{m} + \wh{b} \in \bZ_+,\,  | \mu m + b - \mu \wh{m} - \wh{b}| \le C_1 r k ,\,
|\mu \wh{m} - \mu m - b |  \le C_1 \sqrt{m  r \sfv_k}\}
\label{dfn:Gamma}
\end{align}
and
\begin{equation}\label{dfn:bar-Del}
\bar{\Delta}_{(m,b),(\wh{m},\wh{b})} := \sum_{j=1}^{m \vee \wh{m}} {\bf 1}_{\{\wt{L}_j \ne \wh{L}_j\}} \,.
\end{equation}
\end{lemma}
\begin{proof} 
Fixing $\ka>0$, as \eqref{cond-pre-a} trivially holds for $r \le r_o + \log k$ and any 
$r_o=r_o(\ka)$ finite, it must also hold when $m \le m_o$, for any finite $m_o(\ka)$.
We thus consider hereafter only the asymptotics,  within $\Theta_\ka$, of the
probability on the \abbr{lhs} of \eqref{cond-pre-a},  when $\wh{m} < 2m$ and
both $m \to \infty$ and $r - \log k \to \infty$. 
Further, setting $t_m:=C k \log m$ and for some $C=C(\ka)$ finite,
we can and shall intersect the event on the \abbr{lhs} of \eqref{cond-pre-a} with the event 
\begin{align}\label{eq:E0}
\EE_o &:= \{ \max_{j \le m} \{\wt{L}_j\} \vee \max_{j \le \wh{m}} \{\wh{L}_j\} \le t_m \}\,.
\end{align}
Indeed, $\bP(L_j > t_m,$ some $j \le 3m) \le 
C_0 m^{-C'}$,
for $C' = c \, C -1$ with $c$ as on the right hand side of \eqref{eq:Ber-Geo},  
some $C_0$ finite and all $k \ge 1$. In the sequel, after \eqref{eq:bd-whb},  we show that 
$\bP(S_{\wh{m}} = \mu \wh{m} + \wh{b}) \ge O(m^{-s/2})$ for some finite $s=s(\ka,c)$,
uniformly over $(\wh{m},\wh{b}) \in \Gamma$, as well as for $(m,b) \in \Theta_\ka$.
Consequently, $\bP(\EE_o^{\complement}) \le O(m^{-C'+s/2})$ which for $C=C(\ka)$ large enough, and 
$r \le \frac{1}{\ka} \log m$ is at most $e^{-r}$, as claimed.
We proceed to control the difference between the profiles 
\begin{equation}\label{dfn:Nmb}
N_i^{(m,b)} := \sum_{j=1}^m {\bf 1}_{\{L_j=i\}} \,, \qquad i \ge 1 \,,
\end{equation}
and $\{N^{(\wh{m},\wh{b})}_i \}$ induced by $\bQ_{m,b}$ and $\bQ_{\wh{m},\wh{b}}$ 
on the \emph{un-ordered counts of the corresponding} $\{L_j\}$, 
uniformly over $(m,b)$, $(\wh{m},\wh{b})$ as in the statement of the lemma.
Specifically, we claim that 
there exist $C=C(\ka)$ finite, sub-probability weights $\{\varphi_i\}$ and some 
$n_i=n_i^{(m,b)}$, $\wh{n}_i = n_i^{(\wh{m},\wh{b})}$ 
such that for all $(\wh{m},\wh{b}) \in \Gamma$ \blue{as well as for}
$(\wh{m},\wh{b})=(m,b)$,
\begin{equation}\label{eq:main-est}
\sum_{i=1}^{t_m} \bQ_{\wh{m},\wh{b}} ( N_i^{(\wh{m},\wh{b})} \notin \II_i ) \le C k e^{-r}  \,,
\qquad
\sum_{i=1}^{t_m} |n_i-\wh{n}_i| \le h \,,
\end{equation}
where $\II_i=\II_i^{(\wh{m},\wh{b})} :=(\wh{n}_i-\varphi_i \wh{h},\wh{n}_i+\varphi_i \wh{h})$,
$\wh{h}:=r \ka \sqrt{\wh{m}}
$ and $h:=r \ka \sqrt{m}$. To this end, first note that
for any $i,\wh{m} \ge 1$, $\ell \ge 0$ and $\wh{b} \in \bZ$,
\begin{equation}\label{eq:basic-bd}
\bQ_{\wh{m},\wh{b}} (N_i^{(\wh{m},\wh{b})}=\ell) = \bP(N_i^{(\wh{m})}=\ell) \frac{\bP(S^{[i]}_{\wh{m}-\ell}=
\mu \wh{m}+ \wh{b}-\ell i)}{\bP(S_{\wh{m}}=\mu \wh{m}+\wh{b})} \,,
\end{equation}
with $N_i^{(\wh{m})}$ a Binomial$(\wh{m},q_i)$ and
$S^{[i]}_{\wh{m}}$ denoting the partial sums of \abbr{iid} variables $L_j^{[i]} \deq (L_j | \{L_j \ne i\})$,
that is, $\bP(L_1^{[i]}=j)=\frac{q_j}{1-q_i} {\bf 1}_{j \ne i}$ for $j \in \bZ_+$.
Setting hereafter for any $i \ge 1$,
\[
\mu^{[i]} := \bE L_1^{[i]} \,, \qquad \sfv^{[i]} := {\sf var}(L_1^{[i]}) \,, \qquad u_i = q_i (1-q_i) \,,
\]
and recalling that $\bE L_1 =\mu$, it is easy to check that 
\begin{equation}\label{eq:mu-i}
(1-q_i) \mu^{[i]} + i q_i = \mu \qquad \& \qquad 
\sfv_k = (1-q_i) \sfv^{[i]} + (i-\mu^{[i]})^2 u_i \,.
\end{equation}
In choosing $\{\wh{n}_i\}$ we are guided by the normal approximation 
$\bar \bQ_{\wh{m},\wh{b}}$ for the conditional law of
 $N_i^{(\wh{m},\wh{b})}$ given $S_{\wh m}= \mu \wh{m}+\wh{b}$,
where we replace the Binomial $N_i^{(\wh{m})}$ with a gaussian $\bar N_i^{(\wh{m})}$ 
of mean $\wh{m} q_i$ and variance $\wh{m} u_i$, which is independent of
the gaussian $\bar S^{[i]}_{\wh{m} - \wh{m} q_i}$ of mean zero and
variance $\wh{m} (1-q_i) \sfv^{[i]}$. The gaussian law thus induced on 
\begin{equation}\label{eq:bar-Sn-def}
\bar S_{\wh m} := \, i \bar N^{(\wh{m})}_i + (\wh{m}-\bar N_i^{(\wh{m})}) \mu^{[i]} 
+ \bar S^{[i]}_{\wh{m} - \wh{m} q_i} 
\end{equation}
has thanks to \eqref{eq:mu-i},  the same mean $\mu \wh{m}$ and variance $\wh{m} \sfv_k$ as 
$S_{\wh m}$.
Further, 
with ${\sf cov}(\bar N_i^{(\wh{m})},\bar S_{\wh{m}})=(i-\mu^{[i]}) \wh{m}\, u_i$, we have that
under $\bar \bQ_{\wh{m},\wh{b}}$ the law of $\bar N_i^{(\wh{m})}$ 
given the event $\{\bar S_{\wh{m}}= \mu \wh{m} + \wh{b}\}$,
is gaussian of mean and variance
\begin{align}\label{eq:ci-def}
n_i^{(\wh{m},\wh{b})} := \bE[\bar N_i^{(\wh{m})}|\bar S_{\wh{m}}=\mu \wh{m}+\wh{b}]  &= \wh{m} q_i + 
u_i  (i-\mu^{[i]}) \wh{b}/\sfv_k \,, \\
{\sf var} (\bar N^{(\wh{m})}_i|\bar S_{\wh{m}}=\mu \wh{m}+\wh{b}) 
&\le  
{\sf var} (\bar N^{(\wh{m})}_i) \le \wh{m} q_i \,.
\label{eq:cond-var-bd}
\end{align}
Setting $\phi(x):=(2\pi)^{-1/2} e^{-x^2/2}$ and 
\begin{equation}\label{eq:ofer-xyz}
x_{\wh b}:=\frac{{\wh b}}{\sqrt{\sfv_k \wh{m}}}, \qquad 
z_\ell:=\frac{\ell-\wh{m}q_i}{\sqrt{u_i \wh{m}}}\,, \qquad 
y_\ell:=\frac{\wh{b} - (i-\mu^{[i]}) (\ell-\wh{m}q_i)}{\sqrt{\sfv^{[i]} (1-q_i) \wh{m}}} \,,
\end{equation}
it follows from \eqref{eq:bar-Sn-def} that for some $C=C(\ka)$ finite and all $r$ large enough 
\begin{align}\label{eq:gaussian1}
\bar \bQ_{\wh{m},\wh{b}} (\bar N_i^{(\wh{m})}  \notin \II_i \,  | \, \bar{S}_{\wh m} = \mu \wh m + \wh b ) &=
\frac{\sqrt{\sfv_k \wh{m}}}{\phi(x_{\wh{b}})} \int_{\II_i^{\complement}}
\frac{\phi(z_\ell)}{\sqrt{u_i \wh{m}}} \frac{\phi(y_\ell)}{\sqrt{\sfv^{[i]} (1-q_i)\wh{m}}} d\ell 
\nonumber \\
& \le C \exp\big(-(\varphi_i \wh{h})^2/(2 \wh{m}q_i) \big)
\end{align}
(with the  inequality due to the usual tail bounds for centered gaussian of variance at most
$\wh{m} q_i$, see \eqref{eq:cond-var-bd}). 
Note from \eqref{eq:Ber-Geo} that for $q_i=q_i^{(k)}$ of \eqref{eq:qi-spread},
\begin{equation}\label{eq:exp-tail}
\sup_k \Big\{ k \sum_{i=1}^\infty q_i e^{c i/k} \Big\} < \infty \,,
\end{equation}
whereby it follows that 
\begin{equation}\label{var:ubd}
k^{-1} \sfv_k \le  k \sum_{i=1}^\infty q_i (i/k)^2
\end{equation}
is also uniformly bounded above. 
Proceeding to take $\varphi_i := \delta e^{c i/(4k)} \sqrt{q_i}$, by Cauchy-Schwarz and \eqref{eq:exp-tail},
\[
\sum_{i=1}^\infty \varphi_i  \le \delta \Big(k \sum_{i=1}^\infty q_i e^{c i/k} \Big)^{1/2}
\Big( \sum_{i=1}^\infty \frac{1}{k} e^{-c i/(2k)} \Big)^{1/2} \le c' \delta \,,
\]
so some $\delta(c)>0$ results with $\sum_i \varphi_i \le 1$ 
for all values of $k$ considered here. With such a choice of $\varphi_i$ and having
$\wh{h}=r \ka \sqrt{\wh{m}}$, it follows from \eqref{eq:gaussian1} that
\begin{equation}\label{eq:gaussian2}
\bar \bQ_{\wh{m},\wh{b}} (\bar N_i^{(\wh{m})}  \notin \II_i \,  | \, \bar{S}_{\wh m} = \mu \wh m + \wh b ) 
\le C \exp(-(r \ka \delta)^2 e^{ci/(2k)}/2) \,.
\end{equation}
Consequently, up to the 
use of $\bar \bQ_{\wh{m},\wh{b}}$ instead of $\bQ_{\wh{m},\wh{b}}$,  the
\abbr{lhs} of \eqref{eq:main-est} follows upon summing 
\eqref{eq:gaussian2} over $i \le t_m$.

We next justify the preceding move from $\bQ_{\wh{m},\wh{b}}$ to $\bar \bQ_{\wh{m},\wh{b}}$ by 
using the local \abbr{CLT} of \cite[Thm. VII.3.16]{Petrov} for sums 
of integer valued \abbr{iid} in approximating the terms involving $S^{[i]}_{\wh{m}}$ and $S_{\wh{m}}$ in  
\eqref{eq:basic-bd}. 
Specifically, from the second part of \eqref{dfn:Gamma} we know that $\wh{m} \in [m/2,2m]$ whereas by
\eqref{dfn:Gamma} and our restrictions on $r$ and $k$, 
\begin{equation}\label{eq:bd-whb}
|\wh{b}|  \le C_1 (r k + \sqrt{m r \sfv_k}) \le 2 C_1 \sqrt{\wh{m} r \sfv_k} \le \frac{\sqrt{m r k}}{\ka}  
\end{equation}
(provided $\ka \le \ka_0(C_1)$ is small enough). In particular, the standardized argument 
$x_{\wh{b}}$ of \eqref{eq:ofer-xyz} in the normal approximation 
for $\bP(S_{\wh{m}}=\mu \wh{m}+\wh{b})$ is by \eqref{eq:bd-whb} 
such that $|x_{\wh b}| \le 2\sqrt{r}/(\ka \sqrt{c})$. Thanks to 
\eqref{eq:Ber-Geo}, the $\rho$-th cumulant of the law of $L_j = L_j^{(k)}$ 
scales at most as $k^\rho$, with the multiplicative correction terms 
in the normal approximation of \cite[Thm. VII.3.16]{Petrov} are thus $o(1)$ as long as 
$k |x_{\wh{b}}| = o(\wh{m}^{1/6})$, which applies under our restriction to $k \vee r=O(\log m)$. 
Further, in view of \eqref{var:ubd}, the leading term $(\sfv_k \wh{m})^{-1/2} \phi(x)$ is for $x=x_{\wh b}$
then at least $O(m^{-s/2})$ for some finite $s=s(\ka,c)$. Examining how 
the bounds on $I_1$-$I_4$ of \cite[(VII.3.7)]{Petrov} depend on the law of $L_1^{(k)}$, 
we thus deduce that upon taking the expansion in \cite[Thm. VII.3.16]{Petrov} up to the $s$-th level, 
the resulting additive error be negligible, 
uniformly in the range of $\wh{m},\wh{b},m,b,k$ we consider. 
Proceeding to show that 
\begin{equation}\label{eq:pre-gaussian}
\sum_{i=1}^{t_m} \sum_{\ell \notin \II_i}
\bP(N_i^{(\wh m)}=\ell) \bP(S^{[i]}_{\wh m-\ell}=\mu \wh m+\wh b-\ell i) \le C \frac{k e^{-r}}{\sqrt{\sfv_k \wh m}} \phi(x_{\wh b}) \,,
\end{equation}
recall that the right hand side is at least $O(m^{2-s'/2})$ for some finite $s'=s'(\ka,c)$.
Summing over at most $\wh{m} t_m = o(m^{3/2})$ terms on the \abbr{lhs} of \eqref{eq:pre-gaussian}, 
thanks to the sub-gaussian tails of the Binomial $N_i^{(\wh m)}$ of mean $\wh m q_i$ and variance 
$\wh m u_i$ we may and shall restrict, for some 
$c'=c'(s')$ finite, the inner sum in \eqref{eq:pre-gaussian} to $\II_i^{\complement} \cap \JJ_i$ for the interval 
$\JJ_i$ of width $2 c' \sqrt{\wh m q_i \log m}$ centered at $\wh m q_i$. In particular, 
thanks to \eqref{eq:Ber-Geo} we know that $1-q_i \ge q_0$ is bounded away from zero
uniformly in $k$ and $i \ge 1$, which within the interval $\JJ_i$ 
applies also for $(\wh m-\ell)/\wh m$. By the same reasoning, 
it follows from \eqref{eq:qi-spread} that 
\begin{equation}\label{eq:mu-i-bd}
\sfv^{[i]} \ge q_0 (\mu^{[i]})^2 \ge q_0 (\mu-i q_i)^2 \ge O(k^{-2/c}) \,.
\end{equation}
With $\sfv^{[i]}$ 
decaying at most polynomially in $k \le O(\log m)$, applying \cite[Thm. VII.3.16]{Petrov} 
to each of the terms $\bP(S^{[i]}_{\wh m-\ell}=\mu \wh m+ \wh b -\ell i)$ with an expansion up to the $s'$-th level,
results with a negligible additive error. Further, thanks to the \abbr{lhs} of 
\eqref{eq:mu-i}, the standardized argument in such a normal approximation 
is $\wt{y}_\ell = y_\ell \sqrt{(1-q_i)/(1-\ell/\wh m)}$ for $y_\ell$ of \eqref{eq:ofer-xyz},
where
for some $C=C(\ka)$, uniformly over $\ell \in \JJ_i$ 
and $\wh b,\wh m, k$ satisfying our restrictions, 
\[
|\wt{y}_\ell| \le C k^{1/c}   [ \sqrt{r k}  + c' \sqrt{\sfv_k \log m} ]  \le C^2 (\log m)^{1+1/c} \,.
\]
The right hand side of \eqref{eq:Ber-Geo} implies a similar bound on the \abbr{mgf} of 
$L_1^{[i]}/k$, uniformly in $i \ge 1$ and $k$, hence the same $k^{\rho}$ growth 
rate of the $\rho$-th cumulant of its law. This in turn yields a uniformly $o(1)$ 
multiplicative correction terms in the preceding normal approximation throughout 
the interval $\JJ_i$ of $\ell$ values of interest here (where $k |\wt{y}_\ell|$ stays uniformly
bounded by $(\log m)^{2+1/c}$). Finally,  
$|y_\ell^2-\wt{y}_\ell^2| 
= o(1)$, uniformly over $\JJ_i$, allowing us to replace $\phi(\wt y_\ell)$ by $\phi(y_\ell)$
in the leading term of the approximation. 
Upon utilizing also the normal approximation 
for the Binomial $N_i^{(\wh m)}$, at the standardized argument $z_\ell$ of \eqref{eq:ofer-xyz},
we deduce that it suffices for \eqref{eq:pre-gaussian}, to show that 
\[
\sum_{i=1}^{t_m} \frac{\sqrt{\sfv_k \wh m}}{\phi(x_{\wh b})} \sum_{\ell \notin \II_i} 
\frac{\phi(z_\ell)}{\sqrt{u_i \wh m}} \frac{\phi(y_\ell)}{\sqrt{\sfv^{[i]} (1-q_i) \wh m}} \le C k e^{-r} \,.
\]
As $m \to \infty$, the inner discrete sum over 
$\ell \notin \II_i$ is bounded, up to some universal constant $C=C(\ka)$ by the corresponding
integral over $d\ell$. With the latter integral being precisely the conditional gaussian computation
that led to \eqref{eq:gaussian1}, we thus arrive at the \abbr{lhs} of \eqref{eq:main-est}.

Next, recall the \abbr{lhs} of \eqref{eq:mu-i} that $u_i |i -\mu^{[i]}| = q_i |i-\mu|$. Since  
$\sum_{i} (i-\mu) q_i =0$ (see \eqref{eq:mean-var}), also
\begin{equation}\label{eq:ui-sum-bd}
\sum_{i \ne 0} u_i |i-\mu^{[i]}| \le 2 \sum_{i} q_i (\mu-i)_+ \le 2 \mu \,.
\end{equation}
Our choice \eqref{eq:ci-def} for $\{n_i, \wh{n}_i\}$ thus results for some $C=C(\mu)$ finite, with 
\begin{equation}\label{eq:mean-bd}
\sum_{i=1}^{\infty} |n_i -\wh{n}_i| \le |\wh{m}-m| (1-q_0) + 2 \mu \frac{|\wh{b}- b|}{\sfv_k} \le 
\frac{C \sqrt{m}}{\ka} ((r/k) \vee 1) \,,
\end{equation}
where we have used here the lower bound \eqref{eq:mean-var} on $k^{-1} \sfv_k$, the upper 
bounds \eqref{eq:bd-whb} on $|\wh{b}|, |b|$, and \eqref{eq:Ber-Geo} on $(1-q_0)$, as well as the bound
\begin{equation}\label{bd:m-diff}
|\wh{m}-m| \le \mu^{-1} [\, |b| + C_1 \sqrt{m r \sfv_k}\, ] \le \frac{2 (r \vee k) \sqrt{m}}{\mu \ka} \,,
\end{equation}
due to the second part of \eqref{dfn:Gamma} and the right-inequality in \eqref{eq:bd-whb}. Since 
$(C/\ka) \sqrt{m} ((r/k) \vee 1) \le h$ 
for all $r$ large enough, any $k \ge \ka^{-3}$ and $\ka \le \ka_0(C)$,
we thus arrive at the right hand side of \eqref{eq:main-est}. Now, setting 
\[
\wh{\Delta}_{(m,b),(\wh{m},\wh{b})} := \sum_{i \ge 1} | N_i^{(m,b)} - N^{(\wh{m},\wh{b})}_i | \,,
\]
for the independent collections $\{N_i^{(m,b)}\}$, $\{N^{(\wh{m},\wh{b})}_i \}$, note that under
the event $\EE_o$ of \eqref{eq:E0},
\[
\wh{\Delta}_{(m,b),(\wh{m},\wh{b})} 
\le \sum_{i=1}^{t_m} |N_i^{(m,b)} - n_i| 
 + \sum_{i=1}^{t_m} |\wh{N}_i^{(\wh{m},\wh{b})} -\wh{n}_i| + \sum_{i=1}^{t_m} |\wh{n}_i - n_i| 
\]
and therefore the bounds of \eqref{eq:main-est} yield that for any $\ka>0$ as above
\begin{align}\label{cond-pre-b}
\sup_{\Theta_\ka, \; (\wh{m},\wh{b}) \in \Gamma} \, \Big\{ \frac{e^r}{k} 
\bP \big( \wh{\Delta}_{(m,b),(\wh{m},\wh{b})} \ge 4 r \ka \sqrt{m} \big) \Big\} < \infty \,.
\end{align}
Assuming hereafter, for the sake of definiteness, that $\wh{m} \le m$, setting 
$\wh{L}_j \equiv 0$ for all $\wh{m} < j \le m$ results with  
\begin{equation}\label{bd:N0-diff}
| N_0^{(\wh{m},\wh{b})} - N_0^{(m,b)} | = \big| \sum_{i \ge 1} (N_i^{(m,b)} - N_i^{(\wh{m},\wh{b})} ) \big| 
\le \wh{\Delta}_{(m,b),(\wh{m},\wh{b})} \,.
\end{equation}
Now, conditional upon the independently generated $\{N^{(m,b)}_i\}$ and $\{N^{(\wh{m},\wh{b})}_i \}$, 
the pairs $\{(\wt{L}_j, \wh{L}_j), j \le m \}$, can and will be aligned so that $\wt{L}_j = \wh{L}_j$ for  
$\sum_{i \ge 0} (N_i^{(m,b)} \wedge N_i^{(\wh{m},\wh{b})} )$ indices. Applying a uniformly chosen 
permutation of $\{1,\ldots,m\}$ to these $m$ pairs, attains the marginal 
law $\bQ_{m,b}$ of $\{\wt{L}_j, j \le m \}$ and we
thereafter recover the marginal law $\bQ_{\wh{m},\wh{b}}$ of the shorter sequence by
exchanging the set $\{\wh{L}_j  \ne 0 : \wh{m} < j \le m \}$ whose size we denote by $V$,
with a uniformly chosen subset of the same size from $\{ \wh{L}_j =0 : j \le \wh{m} \}$.
In view of \eqref{bd:N0-diff}, this coupling guarantees having 
\[
\bar{\Delta}_{(m,b),(\wh{m},\wh{b})} \le V + 
\frac{1}{2} \sum_{i \ge 0} |N_i^{(m,b)} - N_i^{(\wh{m},\wh{b})}| \le  
V + \wh{\Delta}_{(m,b),(\wh{m},\wh{b})} \,,
\]
which in view of \eqref{cond-pre-b} results with the desired statement \eqref{cond-pre-a} provided that
\begin{equation}\label{bd:V}
\sup_{\Theta_\ka, \; (\wh{m},\wh{b}) \in \Gamma} \, \Big\{ \frac{e^r}{k} 
\bP \big( V \ge r \ka \sqrt{m} \big) \Big\} < \infty \,, 
\end{equation}
for any $\ka>0$ small enough. 
Turning to establish \eqref{bd:V}, note that given $\{N_i^{(\wh{m},\wh{b})} \}$, the variable 
$V$ counts how many of the $N_{\ge 1} =\sum_{i \ge 1} N_i^{(\wh{m},\wh{b})}$ items
which are uniformly distributed among the given $m$ bins, falls within the last $m-\wh{m}$ bins.
As such, conditioned on $N_{\ge 1}=\ell$, the variable $V$ is stochastically dominated by 
the variable $\bar V$ of Binomial($\ell,\frac{m-\wh{m}}{m-\ell}$) law. Turning to bound 
the upper tail of $\bar V$,
in view of \eqref{eq:main-est} we can \abbr{wlog} replace the random value of $N_{\ge 1}$ by 
\[
\ell := \wh{h}+ \sum_{i \ge 1} \wh{n}_i \le h + m (1-q_0) + 2 \mu |\wh{b}| / \sfv_k 
\le r \ka \sqrt{m} + \frac{m}{ck}  + 2 \mu \frac{(r/k) \vee 1)}{\ka c} \sqrt{m} \,,
\]
where we relied on having $\wh{m} \le m$ as well as 
\eqref{eq:ci-def} and \eqref{eq:ui-sum-bd} to arrive at the first inequality,
and on our choice of $h$, the \abbr{lhs} of \eqref{eq:Ber-Geo}, \eqref{eq:bd-whb} and 
\eqref{eq:mean-var} for the second inequality. The latter bound applies also when setting
$\ell = h + \sum_{i \ge 1} n_i$ in case $m \le \wh{m}$,
and assuming \abbr{wlog} that $\sqrt{m}/(ck) \ge r \ka + 2 \mu r/(\ka c)$,
we see that that regardless of the sign of $m-\wh{m}$,  we can replace $N_{\ge 1}$, 
by $\bar{\ell}=\lceil 2 m/(ck) \rceil$. Upon doing so, we have thanks to \eqref{bd:m-diff}, that 
for $\ka>0$ small enough and $r \ge r_o(\ka)$, 
\[
\bE [\bar V] = \frac{|m-\wh{m}| \bar{\ell}}{m-\bar{\ell}} \le \frac{8 ((r/k) \vee 1)}{c \mu \ka} \sqrt{m}
\le \frac{r \ka}{2} \sqrt{m}\,.
\]
We thus get that \eqref{bd:V} holds 
by the exponential upper tails of Binomial variables such as $\bar V$.
\end{proof}

Note that \eqref{eq:mean-var}--\eqref{eq:qi-spread} hold with $\mu=1$ and some $c>0$,  \blue{for 
$L^{(k)}$ of the form \eqref{dfn:GJp},  with
\begin{equation}\label{dfn:GJpL}
q_0^{(k)} := \bP(L^{(k)}=0)=1-\frac{1}{k+1}  \quad \text{and} \quad q^{(k)}_i := \bP(L^{(k)}=i)=\frac{k^{(i-1)}}{(k+1)^{i+1}} 
\,,  \qquad \forall i \ge 1 \,.
\end{equation}}
Then,  for $S_\ell=S_\ell^{(k)}$ of \eqref{dfn:Sl} and the stopping times
\begin{equation}\label{def:tau-h}
\tau_h := \inf\{ \ell \ge 0 :  S_\ell \ge h \} \,,
\end{equation}
denote by $\bQ^\star_{m'}$ the  joint law of $\{ L^\star_{i}, i \le \tau_{m'} \}$,  
where $L_i^\star=L_i$, $i < \tau_{m'}$ and $L^{\star}_{\tau_{m'}} \le L_{\tau_{m'}}$ is set as needed 
in order to force the event
\begin{equation}\label{eq:Amb-star-def}
\AA^\star_{m'} := \Big\{ \sum_{i=1}^{\tau_{m'}} L^\star_i = m' \Big\} \,.
\end{equation}
Now,  relying on Lemma \ref{lem:emp} we proceed to bound the number of positive $\wt{L}_i$-s and $L^\star_i$-s 
in unmatched \abbr{iid} Geometric blocks,  under the laws $\bQ_{m,b}$ and $\bQ^\star_{m+b}$.
\begin{lemma} \label{lem:coupling} 
For $L^{(k)}$ distributed as \eqref{dfn:GJpL},
consider
$\{ \wt{L}_i,  i \le m\}$ of law $\bQ_{m,b}$ and $\{ L_{i}^\star, i \le \tau_{m'} \}$ of law
$\bQ^\star_{m'}$, $m'=m+b$. Choose 
independently of both sequences,  \abbr{iid} Geometric($1/c')$ random variables $\{\chi_j, j \ge 0\}$   and 
construct the corresponding blocks $\{\wt{\LL}_j, j \in [0,M_m] \}$ and $\{\LL^\star_j, j \in [0,M_{\tau_{m'}}] \}$ by 
summing consecutive $\chi_j$ entries from the finite sequences $\{\wt{L}_i\}$ and $\{L^\star_i\}$ 
(cutting short the last sum as needed). For small $\ka>0$ and any $(m,b,k,r)$ in $\Theta_\ka$,
there exists some $\JJ \subseteq [1,M_{m \wedge \tau_{m'}})$ with 
$\wt{\LL}_j=\LL_j^\star$ when $j \in \JJ$, and
\begin{align}\label{cond-a}
\sup_{\Theta_\ka} \, \Big\{ \frac{e^r}{k} 
\bP \big( \bar \Delta_{m,b} \ge 6  r c' \ka \sqrt{m} \big) \Big\} < \infty \,,
\end{align}
where $\bar \Delta_{m,b}$ count the total number of $\wt{L}_i \ge 1$ and $L^\star_i \ge 1$ which are
from within the collection $\bar \UU_{m,b}$, of non-matched blocks $j \notin \JJ$ (including 
the first and last blocks in each sequence, which we do not aim to match).
\end{lemma}
\begin{proof}
Recall, as in the proof of Lemma \ref{lem:emp},
that it suffices to consider only the asymptotic, within $\Theta_\ka$, of the
probability on the \abbr{lhs} of \eqref{cond-a},  when  both $m \to \infty$, 
$r - \log k \to \infty$. 
Setting $\wh{M}:=\tau_{m'}-1$ and $\wh{B} := S_{\wh{M}}-\wh{M}$, 
note that $\{L_j^\star, \, j \le \wh{M} \}$ is distributed according to a mixture of $\bQ_{\wh{m},\wh{b}}$, 
for non-negative integers $\wh{m}+\wh{b} < m + b$ and mixture coefficients 
\begin{align}
\bQ^\star_{m'} (\, \wh M=\wh m, \, \wh B = \wh b\, )  
 =
\bP(S_{\wh{m}}= \wh{m}+\wh{b}\,)\, \bP(L_1 \ge  m+b- \wh{m}-\wh{b}) \,.
\label{eq:mixture}
\end{align} 
In view of \eqref{eq:mixture}, 
$\bQ_{m'}^\star\big(\, (\wh{M}, \wh{B}) \notin \Gamma\big) \to 0$ faster than $k e^{-r}$
for $\Gamma$ of \eqref{dfn:Gamma}, with $\mu=1$, and some $C_1$ finite
(indeed, this follows by the exponential tails as in the right hand side of \eqref{eq:Ber-Geo}  
for our $L_1$, and the sub-gaussian tail of $\wh{M}$ up to  $O(\log m)$ standard deviations). 
Thus, restricting hereafter the mixture in \eqref{eq:mixture} to $(\wh{m},\wh{b}) \in \Gamma$,
we deduce that for $\ka>0$ small enough as in this application of Lemma \ref{lem:emp},
\begin{align}\label{cond-star}
\sup_{\Theta_\ka} \, \Big\{ \frac{e^r}{k} 
\bP \big( \bar{\Delta}^\star \ge 5  r \ka \sqrt{m} \big) \Big\} < \infty \,, 
\end{align}
where upon extending with zeros the shorter of the two sequences, 
\[
\bar{\Delta}^\star := \sum_{i=1}^{m \vee \wh{M}} {\bf 1}_{\{\wt{L}_i \ne L^\star_i\}} \,.
\]
Note that each $i \ge 1$ for which $\wt{L}_i \ne L^\star_i$, contributes at most the 
corresponding block length $\chi_j$ to $\bar{\Delta}_{m,b}$.  Since the \abbr{iid} $\{\chi_j\}$ are independent 
of the collection $\{\wt{L}_i,L_i^\star, i \ge 1 \}$, this in turn results with the conditional law of
$\bar{\Delta}_{m,b}$ given $\bar{\Delta}^\star=\ell$ being stochastically dominated 
by $\sum_{j=1}^{\ell+2} \chi_j$. In view of \eqref{cond-star}, we thus get \eqref{cond-a} 
(utilizing the exponential tails of the \abbr{iid} $\chi_j$ of mean $c'$). 
\end{proof}

\section{Barriers and effective cover}
\label{sec:apriori}
The bulk of this long and technical section  (up to and including \Cref{subsec-clock-adjust}),  is devoted to the proof of \Cref{prop-apriori}. The estimates we derive in the course of that proof serve us in
\Cref{sec:good-cover-end} to \blue{prove} Proposition \ref{lem:good-cover-end} and in Section \ref{subsec:tools}
to derive estimates \blue{for proving our} second moment and tightness \blue{results} (namely,  in the proofs of
Propositions \ref{prop:2nd-moment} and  \ref{prop:tight}).  \blue{As many arguments} in this section are variants of \blue{those} in \cite{BRZ2},  we \blue{aim to} avoid repetitions \blue{while} providing enough details to follow the needed modifications.
The reader is advised to keep a copy of \cite{BRZ2} handy when following the arguments in this section.

The proof of  \Cref{prop-apriori} is through a series of reductions. There are two main issues at hand: moving from a centralized clock to excursion counts, and introducing some slack in the barriers and windows in various quantities, while controlling the variability in excursion counts around balls of different centers.  \blue{To this end,  we}
introduce \blue{some required} notation \blue{and} recalling $\{r_j\}$ of \eqref{dr.1} we \blue{start with the following} Brownian traversal times. 
\begin{definition}\label{def-trav} For any $0<r < R \le \tilde{r}<\tilde{R}$, integer $m \ge 1$
and $\sfv,\sfx \in \bT^2$ such that $\sfv \in \bB_\sfd(\sfx,\tilde{r}-R)$,
	let $\trav{m}{\sfx}{\tilde R}{\tilde r}{\sfv}{r}{R}$ be the number of
	traversals  
	$\partial \bB_{\sfd}(\sfv, r)\to\partial \bB_{\sfd}(\sfv,R)$
	during $m$ excursions $\partial \bB_{\sfd}(\sfx,\tilde{r})\to\partial \bB_{\sfd}(\sfx,\tilde{R})$.
These traversals are interlaced with the traversals $\partial \bB_{\sfd}(\sfv, R)\to\partial \bB_{\sfd}(\sfv,r)$,  so we set 
$\trav{m}{\sfx}{\tilde R}{\tilde r}{\sfv}{R}{r}=\trav{m}{\sfx}{\tilde R}{\tilde r}{\sfv}{r}{R}$.  
	For $j > k \ge 1$ we often use the abbreviations 
\begin{equation}
\label{eq-calgary}
 T_{\sfv,j}^{\sfx,m}:= \trav{m}{\sfx}{r_0}{r_1}{\sfv}{r_j}{r_{j-1}}, \quad
  T_{j}^{\sfx,k,m} 
  :=\trav{m}{\sfx}{r_{k-1}}{r_{k}}{\sfx}{r_{j}}{r_{j-1}}\,,
 \quad T^{\sfx,m}_j := T^{\sfx,1,m}_j = T^{\sfx,m}_{\sfx,j} 
\end{equation}
(also with the roles of $r_j$ and $r_{j-1}$ reversed).  
More generally,  for Brownian stopping times $\tau$ and $\tau'$,  denote by $T^{[\tau,\tau']}_{\sfv,j}$ the number 
of traversals $\partial \bB_{\sfd}(\sfv, r_j)\to\partial \bB_{\sfd}(\sfv,r_{j-1})$ during 
$[\tau, \tau']$,  with $T^{[\tau,\tau']}_{\sfv,j} \equiv 0$ if $\tau' \le \tau$. 

We also consider the following approximation of $\sfs_1(\cdot)$ of \eqref{eq:nominal},
\begin{equation}\label{clt.2}
\wt \sfs (y) = \wt \sfs_{n_\vep}(y) := n_\vep( 2 n_\vep  - \log n_\vep + y) = \sfs_1(y/2) +  O((\log n_\vep)^2) \,,
\end{equation}
and keep track of the time it takes the Brownian motion to complete a prescribed
number of excursions.  Specifically,  for $0<r<r'$, $m \ge 1$,  we set
\begin{align}
  \tau_{\sfv,r,r'} (m)=&
  \mbox{\rm 
  time it takes the Brownian motion to complete $m$ excursions} \nn 
  \\
  &\mbox{\rm from  
$\partial \bB_{\sfd}(\sfv,r)$ to $\partial \bB_\sfd(\sfv,r')$},
 \label{eq-tauab}
\\
\tau_{\sfv}(m) :=&\tau_{\sfv,r_1,r_0}(m) \,.
\label{eq-taux}
\end{align}
\end{definition}
\begin{remark}
Note the compatibility property 
\begin{equation}
\{ T_{j}^{\sfx,k,m}=T_{j}^{\sfx,t}, \;\; \forall  j > k\}, \mbox{\,\; on the event \; }\left\{ T_{k}^{\sfx,t} = m  \right\}.\label{eq:TkandTCompatability}
\end{equation}
Further,  with ${GW}_j(m)$ denoting the \blue{branching process of a}
Geometric(1/2) offspring distribution,  started at $GW_1(m)=m$,  we have the identity in law
\begin{equation}\label{eq:Jay1.12}
\{T_j^{\sfx,m}\}_j \deq \{{GW}_j(m)\}_j\,.
\end{equation}
\end{remark}

Fixing small $\wt{c}$,  $\srs$ and $\dlc$,  recall 
our $\wt{c} \, \vep$-nets $\bF^\bb$ of the blocks of $\bT^2$,  where
$\bF:=\bF_{n_\vep}=\bigcup_\bb \{ \bF^\bb \}$ is of size $c \, \vep^{-2}$ 
for some $c=c(\wt{c})$ finite  (see \eqref{dfn:F-star}-\eqref{dfn:F-beta}).
For any $\sfv \in \bF$, $j \in [1,n_\vep]$, $k \ge 1$, let 
\begin{align}
T_{\sfv,j}^{\partial \bC^\sfe (k)} := \#\{ & 
\mbox{\rm 
excursions from } \partial \bB_\sfd(\sfv,r_j) \;\;
\mbox{\rm  to } \partial \bB_\sfd(\sfv,r_{j-1})  \,
 \mbox{\rm during the first }
 k \;  \mbox{\rm inter-city excursions} \},
 \label{dfn:T-Cek}
\end{align}
so in particular $\Tvjz=T_{\sfv,j}^{\partial \bC^\sfe (\sfk_\vep(z))}$.  Recall $\etavjz:=\sqrt{2 \Tvjz}$
with $\etajz=\eta^{[z]}_{\sfv_\bb,j}$ and $\whetajz=\etajz-\gamma_0(j)$.

Fixing the various parameters of \eqref{dfn:aux-par},  we further enlarge $\sfJwp=[\sfKo,\sfKw]$ and introduce the events
\begin{align}\label{dfn:bw}
\!\!\!\! \BWz  \!\!  &:=   \!\ \bigcap_{j\in \sfJwp} \! \{ \Ljz \in \bB^\infty_{\bar \mu} (2 \boetaj) \}
\!
\bigcap_{j \in \sfJupl}  \{ \whetajz \in \sfI_j  \} ,\\
\BWs &:= 
\!\! \bigcap_{\sfv \in \bF^\beta} \!\!  \RRv \bigcap \BWz
\label{dfn:bws}
\end{align}

We now give a road-map to the proof of \Cref{prop-apriori}.  First, 
we recall  from  \Cref{rem:alt-UU},
  \eqref{eq:barrier-event-z},  \eqref{def:res-u-z}
and \eqref{eq:was-event-z}
that for proving  \Cref{prop-apriori} it suffices to show that for fixed $z$ 
\begin{equation}\label{eq:U-plus}
\varlimsup_{\sfKp,\sfKfv,\sfKz,\vep}
\bP \Big(\bigcup_{\bb} \Big\{ (\BWs)^{\complement} \bigcap \big\{
\min_{\sfv \in \bF^\bb}  \eta_{\sfv,n_\vep+1}^{[z]} = 0 \big\}  \Big\} \Big) = 0 
\end{equation}
(\blue{where the union is over $\bb$ in the $\sfM$ inner cities we consider in \eqref{dfn:UUz})}.
\blue{The proof of the latter relies on the a-prori controls on a rough barrier
that we establish in Proposition \ref{Thm4.1+4.2-Jay}} (whose proof is given in  \Cref{sec:rough-barrier}).  
\blue{Specifically},  Proposition \ref{Thm4.1+4.2-Jay} will allow us to reduce the proof of 
\Cref{eq:U-plus} to showing \Cref{eq:U-plus-a}. 
The proof of \Cref{eq:U-plus-a},  in turn,  uses Propositions \ref{theo-pre220} (angular and excess control) and \ref{Jay-4.3-4.3rev} (clock synchronization),  whose proofs are provided 
in Sections \ref{sec:Jay-26} and \ref{subsec-clock-adjust}, respectively.  Indeed, with Propositions \ref{theo-pre220} we reduce the proof of
\Cref{eq:U-plus-a} to proving \Cref{eq:U-plus-alt}, and the proof of the latter, which utilizes Proposition \ref{Jay-4.3-4.3rev}, is provided
at the end of the current section.

We now turn to the actual  proof of \eqref{eq:U-plus}. 
While the latter is 
stated in terms of inter-city excursions and working with such is crucial for us,  it is much easier 
to bound the relevant probabilities when events involving $\sfv \in \bF$ are considered 
within the first $\sfs_1(z)$ excursions of level one around $\sfv$ (or $\sfv_\bb$),  as we 
do in Sections \ref{sec:rough-barrier} and \ref{sec:Jay-26}.  Specifically,  
we set 
$\etaxvjp:= \sqrt{2T_{\sfv,j}^{\sfx, \sfs_1(z)} }$ for $j \ge 1$,  $\sfx,\sfv \in \bF$,  
with the shorthand $\etabvjp:=\eta^{\sfv_\bb,\langle z \rangle}_{\sfv,j}$,  $\etavjp:=\eta_{\sfv,j}^{\sfv,\langle z \rangle}$, 
$\etajp:=\etap_{\sfv_\bb,j}$,  and show  in \Cref{sec:rough-barrier}
that for large values of $\sfk$,  the barrier events
\begin{align}
 \label{def:Jay-Thm4.2}
\Rkz  &:= \bigcap_{\bb} \Big\{ \Rkbp \bigcap_{\sfv \in \bF^\beta} \Rkvpp \Big\} 
\end{align}
hold \abbr{whp},  where 
\begin{align}
\label{a216.0}
\Rkbp & := \!\!\!\! \bigcap_{j \in [\sfk,\nvf]} \!\!\! \Big \{ \etajp \ge \gamma_{-1}(j) \Big\},   \quad
\Rkvpp := \!\!\!\!\!\! \bigcap_{j \in [\wnvfv,n_\vep-\sfk]} \!\!\!\!\big\{  
\etabvjp
 \ge \gamma_{-1}(j) \big\} \,.
\end{align}
In words, for a given $\beta$, the event $\Rkbp$ is a lower barrier event for excursion counts around the center 
of the $\beta$-block (between levels $\sfk$ and $\nvf\blue{=n_\vep-\sfKf}$)
while the event $\Rkvpp$ gives,  for $\sfv\in \bF^\beta$,  a lower barrier event for excursion counts 
around $\sfv$ (between levels $\wnvfv\blue{=n_\vep - \sfKfv - \log \dlc}$ and $n_\vep-\sfk$).

\blue{
\begin{remark}\label{rmk:z-lz-adjust}
Recall \Cref{rem:good-ic} that we have the Brownian path starting at some point on $\bTs$
which for all blocks $\bb$ we consider when $\srs \ge r_0$ (e.g.  in \eqref{eq:U-plus}), 
is outside $\bB_{\sfd}(\sfv_\bb,r_1)$.  This aids in comparing 
$\{\etavjz\}_{j,\sfv \in \bF^\bb}$ to $\{\etabvjp\}_{j,\sfv \in \bF^\bb}$ (which we 
monitor only \emph{after first hitting} $\partial \bB_{\sfd}(\sfv_\bb,r_1)$,  see 
\eqref{eq-calgary}).
\end{remark}}
Hereafter $z':=z-\log \sfk$,  $z'':=z+\log \sfk$.  In \Cref{Thm4.1+4.2-Jay}
we also consider the clock events 
\begin{align}\label{def:Jay-Thm4.1}
\Akz & := \bigcap_\bb \Akbz  \quad \text{where} \quad
\Akbz := 
\bigcap_{\sfx',\sfx \in \bF^\bb}
\Big\{ 
 S_{\sfk_\vep(z)} \in [ \tau_{\sfx'}(\sfs_1(z')), \tau_{\sfx}(\sfs_1(z''))] \Big\} 
\end{align}
for $S_k$ of \eqref{dfn:Sk}.   This proposition will be used to transfer the bounds of 
Sections \ref{sec:rough-barrier} and \ref{sec:Jay-26} to bounds on the events in \eqref{eq:U-plus},
\blue{where thanks to \Cref{rmk:z-lz-adjust}},
\begin{align}
\label{def:Jay-Thm4.1-cons}
\Akbz & \;\; \Longrightarrow \;\;
 \Big\{ \etavjz \in [\etappp_{\sfv,j},\etabp_{\sfv,j}],   
\quad \forall  j,  \sfx',  \sfx,  \sfv \in \bF^\bb \Big\} \,. 
\end{align} 

\begin{proposition}\label{Thm4.1+4.2-Jay}
For $\Rkz$ of \eqref{def:Jay-Thm4.2},  $\Akz$ of \eqref{def:Jay-Thm4.1},  small $\wt{c}$, $\srs=r_0$,  there exists $C<\infty$ such that for
any $z\in \bR$ and $\sfk \ge 1$, 
\begin{align}
\varliminf_{\sfKz,\vep} \bP(\Akz) \ge 1-  C \sfk^{-4}  \quad \text{\rm and,  if $11 z \ge - \sfk^{\gak}$, then} \quad 
\varliminf_{\sfKz,\vep}
\bP\big(\Akz \cap \Rkpp  \big) &\ge 1-  C \sfk^{-4}\,.
\label{eq:Jay-Thm-4.2}
\end{align}
\end{proposition}
\noindent
By \eqref{def:Jay-Thm4.1-cons} we have under $\Akbz$ that 
$\eta_{\sfv,n_\vep+1}^{\bb,\langle z' \rangle} \le \eta_{\sfv,n_\vep+1}^{[z]}$ for any $\sfv \in \bF^\bb$.  Thus,  setting 
for $\sfv \in \bF^\bb$ the events
\begin{align}\label{def:UUvk}
\whUUvk &:= \Rkbp \cap \Rkvpp \cap \big\{ \eta_{\sfv,n_\vep+1}^{\bb,\langle z \rangle} = 0 \big\},  \\
\AUvz &:= \big\{ \Akbz,  \eta_{\sfv,n_\vep+1}^{[z]} = 0 \big\}  \setminus \BWz \,,
\label{def:AUvz}
\end{align}
we have in view of \eqref{eq:Jay-Thm-4.2} that it suffices for \eqref{eq:U-plus} to show that,  for any fixed $\sfk$, 
\begin{equation}\label{eq:U-plus-a}
\varlimsup_{\sfKfv,\sfKz,\vep}
\bP \Big( \bigcup_{\bb, \sfv \in \bF^\bb} \big\{ \whUUvkp \cap \AUvz \big\} \Big) = 0 
\end{equation}
(we have \abbr{wlog} that $\sfKp \ge \sfk$,  so \eqref{def:Jay-Thm4.1-cons} implies that 
$\Rkvppp \subset  \RRv
$,
\blue{see \Cref{rmk:z-lz-adjust}},
allowing us to replace hereafter
the complement of $\BWs$ of \eqref{dfn:bws} by the (smaller) complement of $\BWz$ of \eqref{dfn:bw}).

Next,  denoting by $\Ljp$ the angle-difference empirical measures of $j$-level excursions of $\bb$ during that 
block's first $\sfs_{1}(z)$ excursions of level one (that is,  replacing the count 
$\Tjy=T_j^{\sfv_\bb,\sfKo,\sfs_{\sfKo}(y)}$ in \eqref{dfn:emp-ang} by $T_j^{\sfv_\bb,\sfs_1(z)}$),
in preparation for \eqref{eq:U-plus-a} 
we derive in \Cref{sec:Jay-26} the following nearly optimal bounds.
\begin{proposition}\label{theo-pre220} Fix parameters as in \eqref{dfn:aux-par},  with 
$\wt{c},  \dlc$,  $r_0$  small and $\cst$ large.  Set
\begin{equation}\label{def:UUpk}
\whUUpk:= \bigcup_{\sfv \in \bF^\bb} \whUUvk = \Rkbp \cap \whUUp  \quad \text{where} \quad  
\whUUp :=\bigcup_{\sfv\in \bF^\bb} 
\Big\{ \Rkvpp,  \eta_{\sfv,n_\vep+1}^{\bb,\langle z \rangle} = 0\Big\}
\,.
 \end{equation}
Further,   \blue{recall $n_\vep^{(i)}=n_\vep-\sfK_i$ of \eqref{def.n-vep-i}, 
setting for} $\whetajp:=\etajp-\gamma_0(j)$ \blue{the events}
\begin{align}
\!\!\!\! \BWp  \!\!  &:= \! \bigcap_{j\in \sfJwp}  \!  \{ \Ljp \in \bB^\infty_{\bar \mu} (\boetaj) \}
\bigcap_{j \in \sfJup} \! \{ \whetajp \in \sfIp_j  \} \,,
 \label{dfn:bw-Jay}
\end{align} 
where $\sfJup:=[\sfKo-3,\nvf] \cup \{\nvfv\}$ and we have replaced \eqref{dfn:windows} by the narrower
restrictions of 
\begin{equation}\label{dfn:windows-alt}
\sfIp_j := \left\{\!\!\!\begin{array}{cll}
& (-\infty,\frac{1}{2} \sfKfv), & \qquad j = \nvfv \,, \\
& (\sfKf^{\kapo},  \frac{1}{2} \sfKf^{\kapp}),  & \qquad j=\nvf,   \\
& (2 j^{\kapo},  \frac{1}{2} j^{\kapp}),  & \qquad j \le \sfKo,  j = \sfKt,  j=\sfKth \,,  \\
&(2 j_{n_\vep}^\gak, \frac{1}{2} j_{n_\vep}^{\okappa}),  & \qquad  j \in (\sfKo,\nvf)\,,  j \ne \sfKt,\sfKth \,.
\end{array}
\right.
\end{equation}
Then,  for some $C_{\sfk}<\infty$,  all large enough $\sfKfv$ and any fixed $z$,  $\sfk$,
\begin{align}
\varlimsup_{\vep \to 0} \Big\{ e^{2 n_\vep} \sup_{\bb,\sfv \in \bF^{\bb}} 
\Pbm \big( \whUUvk,  \whetap_{\bb,\nvf} \in (\sfKf^{\kapo},\sfKf^{\kapp})
 \big)  \Big\} & \le C_{\sfk}  (1+z_{+}) e^{-2z},  \label{atot-new}
\\
\varlimsup_{\vep \to 0} \Big\{ e^{2 \nvfv} \sup_\bb \Pbm \big( \whUUpk \setminus \BWp
 \big)  \Big\} &=o_{\sfKfv,\sfKz}(1)  (1+z_{+}) e^{-2z}\,.
\label{atot}
\end{align}
\end{proposition}
\noindent
Clearly,  \eqref{eq:U-plus-a} is an immediate consequence of the following stronger statement
\begin{align}
\varlimsup_{\sfKfv,\sfKz,\vep} \Big\{ e^{2 \nvfv} \sup_\bb 
\bP \big(\bigcup_{\sfv \in \bF^\bb} \{ \whUUvkp \cap \AUvz \}  \big) \Big\} = 0 \,.
\label{atot2}
\end{align}
Further,  for any $\bb$,  $\sfk$,  $\vep$,  $z$ and $z'$,
$$
 \bigcup_{\sfv \in \bF^\bb} \{ \whUUvkp \cap \AUvz \} \subset 
 \bigcup_{\sfv \in \bF^\bb} \{ \whUUvkp \cap \AUvz \cap \BWpp \}
  \bigcup_{\sfv \in \bF^\bb} \{ \whUUvkp \setminus \BWpp \} \,.
$$
Hence,   utilizing \eqref{atot} at $z'$ we in turn get \eqref{atot2} upon showing that 
\begin{equation}\label{eq:U-plus-alt}
\varlimsup_{\sfKfv,\sfKz,\vep} \Big\{ e^{2 n_\vep} 
\sup_{\bb, \sfv \in \bF^\bb} \bP \big(\whUUvkp \cap \AUvz \cap \BWpp\big) \Big\} = 0 
\end{equation}
and then considering the union of the events in \eqref{eq:U-plus-alt} over the $c e^{2\sfKfv}$ 
possible  values of $\sfv \in \bF^\bb$. 

When deriving \eqref{eq:U-plus-alt} we bound probabilities of events involving the traversal times $\Tvjz$ 
by noting that \blue{for our cities which are far enough from the starting point},
up to the negligible events $(\Akz)^{\complement}$,  we may consider instead $T^{\sfv,m}_j$ 
for \emph{random} $m \in [a,b]$.  Our next result,  to which \Cref{subsec-clock-adjust} is devoted, 
allows us to mitigate the impact of this uncertainty in the value of $m$.
\begin{proposition}\label{Jay-4.3-4.3rev}
Recall the stopping times $\tau_{\sfx}(a)$ of \eqref{eq-taux} and the various $\sigma$-algebras of \Cref{def:Fvjs}.
For $\delta \in (0,1)$,  some $c=c(\delta)>0$,  any $b>a$ and all $\sfx$, 
\begin{equation}\label{A.2-n}
\bP \big( T_j^{\sfx, b}  - T_j^{\sfx,a} \ge \ell \,|\,  \GG^{\tau_{\sfx}(a)} \,\big) \le e^{-c \, \ell/j},  \qquad 
\forall j \ge 1,  \;\; \ell \ge (1+2 \delta) (b-a) \,.
\end{equation}
%
Further,  there exist $c>0$ so that  if $j \le  n_\vep$ and $T^{[\tau,\tau']}_{\sfx,j}$ is $\GG^{\sfx, \tau}_j$-measurable,
then
\begin{equation}\label{eq:Jay4.14-n}
\sup_{\sfv \in  \bB_\sfd(\sfx, r_j)}
\bP( T^{[\tau,\tau']}_{\sfv,n_\vep+1}=0 \;\big| \; \GG^{\sfx,\tau}_j,  T^{[\tau,\tau']}_{\sfx,j} \ge \ell )
 \le e^{-c \,  \ell/(n_\vep-j)} \,,   \qquad \forall \ell \ge 1 \,.
\end{equation}
\end{proposition}

Prior to applying \Cref{Jay-4.3-4.3rev},  we first modify the events in \eqref{eq:U-plus-alt} 
by utilizing the next lemma,  whose proof is deferred to the end of this sub-section.
\begin{lemma}\label{lem:adjust-AUz}
For some $c_\star=c_\star(\cst)>0$,  large enough $\sfKfv,\sfKz$,  small enough $\vep$,  and all
$\bb$,  $\sfv \in \bF^\bb$, 
\begin{equation}\label{eq:U-plus-event}
\BWpp \cap \AUvz \subset 
\{\tau' \ge \tau \} \cap  
\{ T_{\sfv,n_\vep+1}^{[\tau,\tau']} = 0 \big\} \cap \Big\{  \bigcup_{j \in \sfJwp} \WW^\star_{\sfx,j} \bigcup_{j \in \sfJupl} \BB^\star_{\sfx,j} \Big\} \,,
\end{equation}
where 
$\tau' := S_{\sfk_\vep(z)}$,
$\tau:=\tau_\sfx(a)$, $\sfx:=\sfv_\bb$,
$a:=\sfs_1(z')$,  $b:=\sfs_1(z'')$ and
\begin{equation}\label{dfn:bb-extra-i}
\WW^\star_{\sfx,j} :=   \{ T_{j}^{\sfx,b} - T_{j}^{\sfx,a} \ge c_\star n_\vep \log j \}\,, \qquad
\BB_{\sfx,j}^\star := \{  T_{\sfx,j}^{[\tau,\tau']} \ge \ell_j\}\,,\qquad \ell_j:=\frac{1}{3} (n_\vep-j) j_{n_\vep}^{\kapp}  \,.
\end{equation}
\end{lemma}

\begin{proof}[Proof of \eqref{eq:U-plus-alt}:] By \Cref{lem:adjust-AUz} and the union bound,  we get 
\eqref{eq:U-plus-alt} upon showing that for $\sfx=\sfv_\bb$,
\begin{equation}\label{eq:U-plus-123}
\varlimsup_{\vep \to 0}
\Big\{ e^{2 n_\vep} \sum_{j=\sfKo}^{\sfKw} \sup_{\bb,\sfv \in \bF^\bb} \bP \Big( \whUUvkp \cap \BWpp ,
T^{\sfx,b}_{j} -T^{\sfx,a}_{j} \ge c_\star n_\vep \log j
 \Big) \Big\} = 0 \,,
\end{equation}
as soon as $\sfKz$ is large enough (in terms of $\sfk$ and $z$),  while for $\tau=\tau_\sfx(a) \le \tau' = S_{\sfk_\vep (z)}$,
\begin{equation}\label{eq:U-plus-45}
\varlimsup_{\vep \to 0}
\Big\{ e^{2 n_\vep} \sum_{j \in \sfJupl}\sup_{\bb,\sfv \in \bF^\bb} \bP \Big( \whUUvkp \cap \BWpp ,  \tau' \ge \tau,
T_{\sfv,n_\vep+1}^{[\tau,\tau']} = 0,  T^{[\tau,\tau']}_{\sfx,j} \ge \ell_j \Big) \Big\}  = o_{\sfKz,\sfKfv} (1) \,.
\end{equation}
Starting with \eqref{eq:U-plus-123},  note that for some $C_{\sfk}$ finite,  any  $j \in \sfJwp=[\sfKo,\sfKw]$ and 
$z =o(n_\vep)$,  
\begin{equation}\label{eq:b-a-bd}
2(b-a)=2(\sfs_1(z'')-\sfs_1(z')) = 2 (z''-z') (\sfm_\vep+z)
\le C_{\sfk} n_\vep \le c_\star n_\vep \log j
\end{equation}
provided $\sfKz$ is large enough.  In that case we deduce from \eqref{A.2-n} that as $\sfKw = o(n_\vep)$,  
\begin{equation}\label{eq:step1-123}
\sum_{j=\sfKo}^{\sfKw} 
\sup_{\bb} \bP\Big( T^{\sfx,b}_{j} -T^{\sfx,a}_{j} \ge c_\star n_\vep \log j \,\big|\, \GG^{\tau} \Big) \le 
\sfKw e^{-c \, c_\star n_\vep \log \sfKw/\sfKw} = o_\vep(1) \,.
\end{equation}
Further,  in view of \eqref{atot-new},  for some other finite $C_{\sfk,z}$, 
\begin{equation}\label{eq:step2-123}
\varlimsup_{\vep \to 0} \Big\{ e^{2 n_\vep} \sup_{\bb,\sfv \in \bF^{\bb}} 
\Pbm \big( \whUUvkp \cap \BWpp \big)  \Big\} \le C_{\sfk,z} 
\end{equation}
and as $\whUUvkp \cap \BWpp$ is measurable on $\GG^{\tau}$ we get \eqref{eq:U-plus-123} 
upon combining \eqref{eq:step1-123} and \eqref{eq:step2-123}.

For \eqref{eq:U-plus-45} we use instead \eqref{eq:Jay4.14-n} at levels $j \in \sfJupl$. 
Indeed,  $\bF^\bb \subset \bB_\sfd(\sfx,r_{\nvfv})$ and since the part of the Brownian path 
(after $\tau$) that is not measurable on  $\GG^{\sfx, \tau}_j$  is contained within $\bB_\sfd(\sfx,r_{j-1})$,
the excursion counts $T^{[\tau,\tau']}_{\sfx,j}$ can be read off the \blue{path,  after $\tau$}, 
outside $\bB_\sfd(\sfx,r_{j-1})$ and therefore these counts are $\GG^{\sfx, \tau}_j$-measurable
(even though the time $\tau'=S_{\sfk_\vep(z)}$ is not,  since it depends on the time spent by the Brownian 
motion inside $\bB_\sfd(\sfx,r_{j-1})$).
Applying the tail bound \eqref{eq:Jay4.14-n} in this setting,  we thus deduce that 
\begin{equation}\label{eq:step1-45}
\sum_{j=\sfKo}^{n_\vep-\sfKfv}
\sup_{\bb,\sfv \in \bF^\bb} 
\bP\Big( T_{\sfv,n_\vep+1}^{[\tau,\tau']} = 0\,|\, \GG^{\sfx,\tau}_{j}, 
T^{[\tau,\tau']}_{j} \ge \ell_j  \Big) \le \sum_{j=\sfKo}^{n_\vep-\sfKfv} 
e^{- c \, j_{n_\vep}^{\kapp}/3} = o_{\sfKz,\sfKfv}(1) \,.
\end{equation}
As before,  combining \eqref{eq:step2-123} and  \eqref{eq:step1-45} yields \eqref{eq:U-plus-45},
thereby completing the proof.
\end{proof}

\begin{proof}[Proof of Lemma \ref{lem:adjust-AUz}] 
For $\tau'=S_{\sfk_\vep}(z)$ and  $\tau=\tau_{\sfv_\bb}(\sfs_1(z'))$,   the event $\{\tau' \ge \tau\}$ is part of  the
conditions in $\Akbz$ (see \eqref{def:Jay-Thm4.1}).  In that case we also have that 
$T_{\sfv,n_\vep+1}^{[\tau,\tau']} \le T_{\sfv,n_\vep+1}^{[z]}$ and recall  \eqref{def:AUvz} that 
$\{ T_{\sfv,n_\vep+1}^{[z]} = 0\}$ is part of $\AUvz$.  Next,  
recall \eqref{def:Jay-Thm4.1-cons} 
that $\Akbz$ implies having $\whetajz \in [\whetajpp,\whetajppp]$ for all $j$.  Thus,  setting 
\begin{align*}
\WW'_{\bb,j}:= \Big\{ \Akbz, \BWpp,  \Ljz \notin \bB^{\infty}_{\bar \mu} (2\boetaj) \Big\}\,,&
\qquad
\BB'_{\bb,j} := \Big\{  \Akbz, \BWpp,  \whetajz \ge j_{n_\vep}^{\kapp} \Big\}
\end{align*}
upon comparing \eqref{dfn:bw} at $z$ with \eqref{dfn:bw-Jay} at $z'$,   
and recalling \eqref{dfn:J-intervals},  \eqref{dfn:windows},  it follows that 
\begin{align}
 \Akbz \cap \BWpp \setminus  \BWz &\subset \bigcup_{j \in \sfJwp}  \WW'_{\bb,j} \bigcup_{j \in \sfJupl} 
 \BB'_{\bb,j} \,.
\label{dfn:bb-extra}
\end{align}
In view of \eqref{dfn:bb-extra} we conclude the proof upon showing that 
$\WW'_{\bb,j} \subset \WW^\star_{\sfv_\bb,j}$ and $\BB'_{\bb,j} \subset \BB^\star_{\sfv_\bb,j}$.  
To this end,   recall that $\Akbz$ implies that
$\tau' \in [\tau_{\sfv}(a),\tau_{\sfv}(b)]$ and $\BWpp$ 
implies that  $\Ljpp \in \bB^{\infty}_{\bar \mu} \big(\boetaj \big)$.  Hence by 
\eqref{dfn:emp-ang},  \eqref{dfn:Poiss} 
 and \eqref{eq:dinf-comp} 
(at $\wt{\boeta}=c^{-1} \,  \boetaj$,  $c := 2\pi (\sqrt{1+\bar{c}^2} + \bar c)$),  if in addition 
$\Ljz \notin \bB^{\infty}_{\bar \mu} (2 \boetaj)$ then necessarily 
$ \Tjz \ge (1 + c^{-1} \,  \boetaj) T_{j}^{\sfv_\bb,a}$,  which under $\Akbz \cap \BWpp$ yields
for $c_\star := {\cst}/{(2c)}$,  any $\vep$ small enough and all $j \in \big[\sfk,\sfKw]$,
\begin{equation}\label{eq:Tb-Ta-W}
T_{j}^{\sfv_\bb,b} - T_{j}^{\sfv_\bb,a} \ge \Tjz -  T_{j}^{\sfv_\bb,a} \ge  c^{-1} \,  \boetaj T_{j}^{\sfv_\bb,a} \ge 
\frac{\gamma_0(\sfKw)^2}{\sfm_\vep} c_\star  \log j \ge c_\star n_\vep \log j 
\end{equation}
(recall \eqref{def:mn-alt},  \eqref{def.nvep},   \eqref{eq:norm-exc},  \eqref{dfn:bw-Jay},
\eqref{dfn:windows-alt} and $\sfKw=o(n_\vep)$).  That is,  $\WW'_{\bb,j} \subset \WW^\star_{\sfv_\bb,j}$ as claimed.
Similarly,  for any  $j \in \sfJupl$,  having in $\BB'_{\bb,j}$ 
\begin{equation}\label{eq:diff-hateta}
\whetajz \ge j_{n_\vep}^{\kapp} \quad \text{and} \quad \whetajpp \le \frac{1}{2} j_{n_\vep}^{\kapp} \,,  
\end{equation}
implies that for $\sfv=\sfv_\bb$ and all $\vep$ small enough 
\begin{equation}\label{eq:Tb-Ta-nvf}
T_{\sfv,j}^{[\tau_{\sfv}(a),\tau']}=T_{\sfv,j}^{[z]} - T_{j}^{\sfv,a} \ge \frac{1}{2} \gamma_0(j)
(\whetajz-\whetajpp)
\ge \frac{1}{4} \gamma_0(j) j_{n_\vep}^{\kapp}  
\ge  \frac{1}{3} (n_\vep-j) j_{n_\vep}^{\kapp} \,.
\end{equation}
That is,  $\BB'_{\bb,j} \subset \BB^\star_{\sfv_\bb,j}$,  which concludes the proof of the lemma.
\end{proof}

\subsection{Rough barrier and continuity: Proof of  \Cref{Thm4.1+4.2-Jay}}\label{sec:rough-barrier}

\Cref{Thm4.1+4.2-Jay} is a direct consequence of Propositions \ref{theo-etrt} and \ref{gprop:MainGoalEasyRegime} below,  which are the analogs 
of \cite[Theorems   3.1 and  3.8]{BRZ2} and of \cite[Proposition 3.4]{BRZ2},  respectively,
with $n_\vep$ standing for what is called there $L$.
Here $z$ can be negative,  which we mitigate in 
\Cref{gprop:MainGoalEasyRegime} by taking $\sfk$ large,  where by restricting
to levels $[\sfk,n_\vep-\sfk]$ we also pick up the crucial extra term $\sfk^\gak/5$ 
in the exponent of \eqref{g17.0}.
\begin{proposition}
\label{theo-etrt} For small $\wt{c}, r_0>0$ there exist finite $
c=c(r_0,\wt{c})$ and 
$
d=d(r_0,\wt{c})$ such that with $\wt\sfs(\cdot)$ of \eqref{clt.2} and $\tau_{\sfv}(\cdot)$ of \eqref{eq-taux}, for any $\vep$ sufficiently small and all $z,z_0$ satisfying $|z| \vee z_0 \leq \sqrt{n_{\vep}}$, 
\begin{equation}
\Pbm\( \wt \sfs \left(z-d\sqrt{z_{0}}\right)     \leq \pi \tau_{\sfv}\(\wt \sfs \left(z \right) \)\leq  \wt \sfs \left(z+d\sqrt{z_{0}}\right),\,\forall \sfv \in \bF
\)\geq 1-ce^{-4z_{0}}.\label{r1.1}
\end{equation}
Moreover,  with $\sfs_1(\cdot)$ of \eqref{eq:nominal},  for any such $\vep$ and $z$,
\begin{equation}
\Pbm(\QQ)\geq 1-ce^{-n_{\vep}/2}, \label{g17.0b}
\end{equation} 
where 
\begin{equation*}
  \QQ:=\{ \tau_{\sfy}\(\sfs_1(z)-10\)\leq \tau_{\sfx}\(\sfs_1(z)\)\leq \tau_{\sfy}\(\sfs_1(z)+10\),\,\forall 
  \sfx,\sfy \in \bF
  , \, \sfd (\sfx,\sfy)\leq  r_{n_{\vep/2}}\}.
\end{equation*}
\end{proposition}

For  the next proposition we introduce the rough barrier events
\begin{align}
\label{eq:JayThm2.2eq}
\Rkp  &:= \bigcap_{\sfv \in \bF} \Rkvp \quad  {\rm with}  \quad  
\Rkvp := \!\!\! \bigcap_{j \in [\sfk,n_\vep-\sfk]} \!\!\!  \Big\{ \etavjp \ge \gamma_{-1}(j) \Big\} \,.
\end{align}
In comparison with (the intersection of) the barrier events $\Rkbp$ and $\Rkvpp$ of \eqref{a216.0}, the centering around which the excursions are measured is always $\sfv$, and no gaps are present.
\begin{proposition}
\label{gprop:MainGoalEasyRegime} 
For small  $\wt{c}$, $r_0$ and $c,d$ as in Proposition \ref{theo-etrt},  we have for any $z \in \bR$,  $\sfk \ge 1$,
\begin{align}
 \widetilde{q}_{\sfk}(z) :=  \varlimsup_{\vep \to 0} \bP \big( (\Rkp)^{\complement} \big)
& \leq 
c
e^{-2(z  - d  \sqrt{|z|})} e^{-\sfk^{\gak}/5}\,.
\label{g17.0}
\end{align}
\end{proposition}

\begin{proof}[Proof of \Cref{Thm4.1+4.2-Jay} (assuming Propositions \ref{theo-etrt}  and \ref{gprop:MainGoalEasyRegime})]
Taking $\wt{c}$, $r_0$ small enough,  we have by \eqref{r1.1} (replacing $d\sqrt{z_0}$ by $t$) that for fixed $t$,
\begin{align}
\label{eq-110926b}
1 - \varliminf_{\vep\to 0} \,  \inf_{|z| \le \sqrt{n_\vep}} \,
\bP \Big(\wt \sfs (z-t) \le \pi \tau_{\sfv} (\, \wt \sfs (z)) \le \wt \sfs (z+t),\, 
 \forall \sfv \in \bF \Big)
 \le c e^{-4 t^2/d} := q_{t}  \,.
\end{align}
Letting
\begin{equation}\label{eq:needed-Jay2.1}
C_1=C_1(\wt{c},r_0):=
 4 \sup_{\sfk} \big\{ \sfk^4 q_{\, \blue{2 \log \sfk-6}}
\big\}< \infty \,,
\end{equation}
we proceed to \blue{establish}  the first claim in \eqref{eq:Jay-Thm-4.2} \blue{by showing}
that, without restriction on $z,\sfk$,
\begin{equation}\label{eq:Jay-Thm-4.1}
\varliminf_{\sfKz,\vep}
\bP\big(\Akz \cap \Akzm \big) \ge 1 - C_1 \sfk^{-4} \;.
\end{equation}
To see \eqref{eq:Jay-Thm-4.1},  recall from \eqref{dfn:bar-svep} and the proof of \Cref{lem:conc} (see \eqref{eq-amirswriteup})
that the time $S_\ell$ it takes
to complete $\ell$ inter-city excursions is such that 
$(S_{\sfk_\vep(z)} -  \bar \sfs_\vep(z)/\pi)/\sfm_\vep \stackrel{p}{\to} 0$ 
as $\vep \to 0$ followed by $\sfKz \to \infty$.  In view of the right hand side of \eqref{def:Jay-Thm4.1}, 
it suffices for \eqref{eq:Jay-Thm-4.1} to show that for any $z$,
\begin{equation}
\label{eq-110926a}
\varliminf_{\vep\to 0} \bP \Big(
\tau_{\sfv} (\sfs_1 (z')) \le \frac{1}{\pi} \bar \sfs_\vep(z-1) \le \frac{1}{\pi} \bar \sfs_\vep(z+1) \le
\tau_{\sfv} (\sfs_1 (z'')),  \forall \sfv \in \bF \Big) \ge 1 -  \frac{1}{2} \blue{C_1} \sfk^{-4} \,.
\end{equation}
(Indeed, if the latter hold for $z$ then it also holds when $z$ is replaced by  $z-2\log \sfk$, 
with the appropriate changes in $z'$ and $z''$.)
Recall that $\bar \sfs_\vep(y)=\sfs_1(y)$ and  
$\wt \sfs (2y-2) \le \sfs_1(y) \le \wt \sfs(2y+2)$ for any $y$ and all $\vep \le \vep_o$  
(see \eqref{clt.2}).  It thus suffices for \eqref{eq-110926a} 
to show that for $t=2 \log \sfk-6$ and $\wh{z} = 2 z \blue{\pm} 2( \log \sfk -1)$, 
\[
\varliminf_{\vep\to 0} \bP \Big(
\wt \sfs (\hat{z}-t) \le \pi \tau_{\sfv} (\wt \sfs (\wh{z})) \le \wt \sfs(\wh{z}+t)  
, \;\; \forall \sfv \in \bF \Big) \ge 1 -  \frac{1}{4} C_1 \sfk^{-4} \,,
\]
which in turn follows from \eqref{eq-110926b} and \eqref{eq:needed-Jay2.1}.

Next,  
by \eqref{def:Jay-Thm4.1-cons},  the event  $\Akzm$ implies that $\Rkpm \subset \Rkpp$ 
and therefore 
\begin{equation}
\label{eq-120926c}
\big(\Akz\cap \Akzm \cap\Rkpm\big)\subset  \big(\Akz\cap \Akzm \cap\Rkpp\big)\subset \big(\Akz \cap\Rkpp\big).\end{equation}
We show below that,  with $\wt{q}_\sfk(\cdot)$ of \eqref{g17.0},  it holds that
\begin{equation}\label{eq:needed-Jay2.2}
C_2=C_2(\widetilde{c},r_0) := \sup_{\sfk \ge 1} \sup_{11 z \ge - \sfk^\gak} \{ \sfk^4 \wt{q}_\sfk (z-3 \log \sfk) \} < \infty \,.
\end{equation}
We then  get the second half of \eqref{eq:Jay-Thm-4.2} with $C=C_1+C_2$  by combining \eqref{eq:Jay-Thm-4.1},
\eqref{eq-120926c} and \Cref{gprop:MainGoalEasyRegime} with $z$ replaced by $z-3\log \sfk$.

It remains to check \eqref{eq:needed-Jay2.2}. Since
$2 d \sqrt{|z-3\log \sfk|} \le 10 d^2 + (|z| +3 \log \sfk)/10$,   we get from \eqref{g17.0} that 
\[
\sfk^4 \widetilde{q}_\sfk(z-3 \log \sfk) \le c e^{10\, d^2} \sfk^{11} e^{- \frac{1}{5} \sfk^{\gak} - 2z + \frac{|z|}{10}}
\]
which is more than enough for \eqref{eq:needed-Jay2.2}.
\end{proof}

In our proofs of Propositions \ref{theo-etrt} and \ref{gprop:MainGoalEasyRegime},
the function $h$ which is used in \cite{BRZ2} to modify radii 
from the sphere to planar coordinates,  is the identity,  while due to the different 
symmetries of $\bT^2$ compared to the sphere,  the continuity estimates 
in \cite[Section 5]{BRZ2},  which serve as input to various steps in the proof,  
must now be reworked,  as we do in \Cref{sec-continest} (and utilize via Lemmas \ref{rlem:Continuity lemma} 
and \ref{lem:4.58m}). 
Recall also the nested $\frac{\wt{c}}{2} r_j$-nets $\bF_{j}$ satisfying \eqref{dr.2},  for $r_j$ as in \eqref{dr.1},
using also the shorthand 
\begin{equation}\label{def:tz-alpha}
\blue{ t_z:=\sfs_1(z) \quad {\rm of} \quad \eqref{eq:nominal} \quad 
{\rm and} \quad \alpha(\cdot):=\gamma_{-1}(\cdot)
\quad {\rm of} \quad \eqref{eq:norm-exc}. }
\end{equation}
The first ingredient for \Cref{theo-etrt} is the following uniform tail estimate on
$\tau_{\sfv,a,b} (m)$ of \eqref{eq-tauab}.  The proof is given 
 in  \cite[Prop.  8.10]{BK},  see also \cite[Lemma 2.2]{DPRZ}, and is therefore not repeated here.
\bl\label{lem-clt}
Setting $\kab:=\frac{1}{\pi} \log(b/a)$,  for some $c'>0$ any $0<a<b<1/2$,  
$\delta \in (0,c')$,  $m \ge 2$ and $\sfv, \sfy\in\bT^{2}$,
\begin{align}
\Pbm^{\sfy}\Big( (1-\de) m \le 
\kab^{-1} \tau_{\sfv, a,b}(m) \le (1+\de) m \Big)& \geq 1 - e^{-c \de^{2} m},\label{clt.1}
\end{align}
where  $c= c(a,b):=c' b (1-a/b)^{6} | \log a|^{-2}$.
\el
Setting as in \cite{BRZ2},  
\begin{equation}\label{eq:tilde-ri}
\tilde{r}_1=r_1-n_\vep^{-1},\quad
\tilde{r}_0=r_0+n_\vep^{-1},\quad
\tilde{r}_{-1}=2 \tilde{r}_0,\quad
\tilde{r}_{-2}=\tilde{r}_{-1} \frac{\tilde{r_0}}{\tilde{r_1}},  
\end{equation}
the second ingredient for \Cref{theo-etrt} is the following continuity estimate on excursion counts,  an analog of \\cite[Lemma 6.2]{BRZ2}.  The proof is a simple adaptation of 
the proof of \cite[Lemma 6.2]{BRZ2},  
using \Cref{lem-cont1} in lieu of \cite[Lemma 5.1]{BRZ2} there,  and is therefore omitted.
\bl
\label{rlem:Continuity lemma}
Fixing $r_1$,  $\wt{c}>0$,  there exists $a_0$,  $c$ finite,  so for any {$|z| \vee z_0 \le \sqrt{n_\vep}$} 
and $n_{\vep}$ large enough, 
\begin{equation}
\Pbm \Big[\exists \sfy \in  \bF_{\frac{3}{2}\log n_{\vep}}\cap  \bB_{\sfd}\left({\bf 0},\wt{c}\, \tilde r_{0}\right)
\label{req: continuity goal}
 \mbox{ s.t. }\left|T_{{\bf 0},\tilde{r}_{1} }^{{\bf 0},\tilde{r}_{-1},\wt \sfs (z)}-
T_{\sfy,\tilde{r}_{1} }^{{\bf 0},\tilde{r}_{-1},\wt \sfs (z)}\right|\ge a_0 \sqrt{z_{0}} n_{\vep} \Big]\le ce^{-4z_{0}}.
\end{equation}
\el

\begin{proof}[Proof of \Cref{theo-etrt}]
Recall that $\wt \sfs (z) = 2 n_\vep^2 (1+o_\vep (1))$,  see \eqref{clt.2}.  Hence,  
by \Cref{lem-clt} with $m=\wt \sfs (z)$ and $\delta=d\sqrt{z_{0}} n_\vep/m$,  
we get that for all $\vep$ small enough,  uniformly in $\sfv, \sfy \in\bT^{2}$,
\bea
\Pbm^{\sfy} \Big( \wt \sfs (z-d\sqrt{z_0}) \leq  
\kab^{-1} \tau_{\sfv,a,b}(\wt \sfs(z)) \leq  \wt \sfs (z+d\sqrt{z_{0}}) \Big) 
\geq 1-e^{-c d^2 z_0/3}.
\label{clt.3a}
\eea
Taking $a=r_1$, $b=r_0$ and $d^2 c(a, b) \ge 30$ in \eqref{clt.3a} yields 
\begin{equation}\label{clt.3a-new}\blue{
\inf_{\sfv,\sfy \in \bT^2}
\Pbm^{\sfy} \big( \,  \wt \sfs \left(z-d\sqrt{z_{0}}\right)     \leq \pi \tau_{\sfv}\(\wt \sfs \left(z \right) \)
\leq  \wt \sfs \left(z+d\sqrt{z_{0}}\right) \,
\big) \geq 1-e^{-10 \, z_{0}} \,.}
\end{equation}
With $|\bF_{n_{\vep}}|$ diverging in $\vep$, 
we get \eqref{r1.1} out of \eqref{clt.3a-new} via a chaining argument,  as in the 
derivation of \cite[Theorem 3.1]{BRZ2} out of \cite[(6.4)]{BRZ2}.  Indeed,  we
replace the factors $4$ and $\sqrt{z}$ of
\cite[(6.4)]{BRZ2} by ${1}/{\pi}$ and $\sqrt{z_0}$ of \eqref{clt.3a-new},
then having here $h(r)=r$,  $L=n_\vep$ and $\tilde{r}_i$ of \eqref{eq:tilde-ri} set 
as in \cite{BRZ2},   we replace \cite[Lemma 6.2]{BRZ2}  by \Cref{rlem:Continuity lemma}.

Finally,  the direct 
proof of \cite[Theorem 3.8]{BRZ2} via a union bound (see \cite[Section 6]{BRZ2}),  applies
here with no change,  to yield \eqref{g17.0b}.
\end{proof}

Turning to the proof of \Cref{gprop:MainGoalEasyRegime},  \blue{recall 
$T_{\sfv,l}^{\npp,t}= \trav{t}{\npp}{r_0}{r_1}{\sfv}{r_l}{r_{l-1}}$ of \eqref{eq-calgary}
(with $T_l^{\npp,t}$ when $\sfv=\npp$),  and} $I_u=[u,u+1)$,  introducing,  as in \cite[(3.20) and (3.30)]{BRZ2},  the events 
\begin{align}
\mathcal{B}_{t,l}&:=\Big\{ \exists \sfv \in \bF_{n_{\vep}}\cap \bB_{\sfd}\left({\bf 0},\wt{c}\, r_{l}\right)
\quad \mbox{such that} \quad
\sqrt{2T_{\sfv, l}^{{\bf 0}, t}}\le \alpha(l)\Big\},\label{17.defb}\\
  \label{eq-Ddef}
\mathcal{D}_{t,l}(j)&:=\mathcal{D}_{{\bf 0},t,l}(j) = \Big\{ \sqrt{2T_{l-2}^{{\bf 0},t}}\in 
I_{\alpha(l)+j}\Big\}
\end{align}
(note that $\mathcal{B}_{t,l}$ involve traversal counts centered at $\sfv \ne {\bf 0}$).  Further 
noting that \blue{$l \ge \sfk + 4 \log (\alpha(l)+j+1)$ for  $\sfk$ fixed}, 
$n_\vep$ large and any $j \le 12 n_\vep$,  $l \ge 5 \log n_\vep$,  
we define for such $l$ and $j$ the events
\begin{equation}
  \label{wCdefr}
\wt{\CC}_{t,l}(j) :=\bigcap_{l' \in [\sfk,\blue{l-h_{l,j}}]}
\big\{ \sqrt{2T_{l'}^{{\bf 0},t}} \ge \alpha(l') \big\} \quad \blue{ {\rm where} \quad h_{l,j}:=4 \log (\alpha(l)+j+1) \,.}
\end{equation} 

The key ingredient is then our next lemma,  which replaces \cite[Lemmas 3.7 and 3.10]{BRZ2} \blue{and is proved}
after the statement of the continuity estimates in \Cref{lem:4.58m}.
\begin{lemma}
  \label{lem-sec3.3}
For $r_0, \wt{c}>0$ small enough,  some finite $C=C(r_0, \wt{c})$,  \blue{$t=t_{\wh z}$ of \eqref{def:tz-alpha}
and} any fixed
$\wh z \in \bR$,  all large enough $n_\vep$,  any $l \in  (5 \log n_\vep, n_{\vep}]$ and $j \in [0, 12 n_\vep]$,
\begin{align}\label{cb.15y}
\Pbm\left(\BB_{t,l}\,\big |\,\wt{\CC}_{t,l} (j) \cap\DD_{t,l} (j)\right)    &  \le Ce^{-4j } \,, \\
\Pbm\left(\mathcal{B}_{t,l}\,\big|\, \mathcal{D}_{t,l}(   j)\right) & \leq C e^{ -4j}\,.
\label{17.1sv}
\end{align}
Further,  \eqref{17.1sv} holds also for $l \in [\sfk,5 \log n_\vep]$ and $j \in [\frac{1}{4} l^{\gak},12l]$.
\end{lemma}

\blue{To match notations with \cite[(3.55) and  (3.58)]{BRZ2},  we set here
\begin{equation}\label{BRZ2-355}
\wt r_l:=(1-\beta_{l,j}) r_l \quad {\rm and} \quad \wt r^+_{l-1}:=(1+\beta_{l,j}) r_{l-1} \quad {\rm where} \quad 
\beta_{l,j}:=\frac{1}{\al(l)+j} \le \frac{1}{2}  
\end{equation}
(not to be confused with our use of $\bb$ as the block index).
Recalling} the number $\trav{t}{\sfx}{\tilde R}{\tilde r}{\sfv}{R}{r}$ of traversals  
$\partial \bB_{\sfd}(\sfv, R)\to\partial \bB_{\sfd}(\sfv,r)$ during $t$ excursions 
$\partial \bB_{\sfd}(\sfx,\tilde{r})\to\partial \bB_{\sfd}(\sfx,\tilde{R})$,  our
proof below of \Cref{lem-sec3.3} relies on the 
following continuity estimate,  whose proof is deferred to \Cref{sec:4.58m}.
\bl\label{lem:4.58m}  
For small enough $\wt{c}>0$,  any $d_0>0$,  some finite $j_o$,  any 
$\blue{t}, m \ge 1$,
 all $j_o \le j \le 12 n_\vep$,  \blue{$l \ge h_{l,j}+1$ of \eqref{wCdefr} and}
$\wt \sfy,  \sfy \in \bB_{\sfd}\left({\bf 0},\wt{c}\, r_{l}\right)$ 
with $\tilde q:=\sfd \left(\sfy, \wt \sfy\right)/r_{l} \ge \frac{\wt{c}}{4} \blue{\beta_{l,j}}$ 
\blue{for $\beta_{l,j}$ of \eqref{BRZ2-355}},
\begin{equation}
  \label{eq-thanks4lem}
 \sup_{u\in I_{\alpha(l)+j}}
 \Pbm\Big(
\big| 
\trav{u^2/2}{{\bf 0}}{r_{l-3}}{r_{l-2}}{\wt \sfy}{\wt r_{l-1}^+}{\wt r_l}
 -\trav{u^2/2}{{\bf 0}}{r_{l-3}}{r_{l-2}}{\sfy}{\wt r_{l-1}^+}{\wt r_l}
 \big| \ge \frac{d_{0}}{2} j  \, m \, 
u \sqrt{\tilde q} \,\Big|\,  \GG^{\bf 0}_{\blue{l-h_{i,j}}},
T_{l-2}^{{\bf 0},t} = \frac{u^2}{2} \Big)
\leq  e^{-4(j+ m)}\,.
\end{equation}
The same applies without conditioning,  for any $l \in [3,n_\vep]$ and $m \le 1/(j_o \beta_{l,j})$.
\el

\begin{proof}[Proof of \Cref{lem-sec3.3}] 
\blue{We start by showing that \eqref{cb.15y} holds for any $t \ge 1$}.  Recall that for 
the values of $l$,  $j$ considered here,  $\bar k := l - h_{l,j} \ge \sfk$.  Hence, 
upon replacing $\CC_{0,t,l}$ of \cite[(3.49)]{BRZ2} by $\wt{\CC}_{t,l}(j)$ 
(i.e.,  dropping part of the barrier in order to handle the conditioning without symmetry),  
and setting as in \cite[(3.58)]{BRZ2}
\begin{equation}\label{BRZ2-358}
\blue{ k_{l,j}:=l+\log(2/\beta_{l,j}), }
\end{equation}
restricted \abbr{wlog} to $\blue{k_{l,j}} \le n_\vep$,
the proof of \eqref{cb.15y} is identical to that of \cite[Lemma 3.10]{BRZ2},   apart from replacing 
\cite[Lemma 3.11]{BRZ2} by the bound of \eqref{eq-thanks4lem},  where $d_0$ and $\wt \sfy$
are chosen as in \cite[(3.67)]{BRZ2}.  \blue{Indeed,  since} $\sfy,  \wt \sfy \in \bF_{\blue{k_{l,j}}}$ they are either 
the same or of distance at least $\frac{\wt{c}}{2}  r_{k_{l,j}} \ge \frac{\wt{c}}{4} \beta_{l,j} r_l$,  \blue{while} 
the measurability of $\wt{\CC}_{t,l} (j)$ on $\GG^{\bf 0}_{\bar k}$ (of \Cref{def:Fvjs}),  allowed us to 
enlarge the conditioning in our analog of \cite[(3.72)]{BRZ2}.

Turning to prove \eqref{17.1sv} \blue{at $t=t_{\wh{z}}$ and fixed $\wh{z} \in \bR$},
we follow the proof of \cite[Lemma 3.7]{BRZ2},  noting that \blue{for any $t \ge 1$}
the bound of \cite[(3.66)]{BRZ2} 
applies for $\mathfrak{P}_{31}$ of \cite[(3.61)]{BRZ2} which involves concentric excursion counts.
Hence,  \blue{setting for $\wt r_{l}$ of \eqref{BRZ2-355} and $k_{l,j}$ of \eqref{BRZ2-358} the event
\[
 \mathcal{B}_{t,l,j}:=\Big\{ \exists \sfy\in \bF_{k_{l,j}}\cap \bB_{\sfd}\left(\npp,\tilde{c}\, r_{l} \right)
 \mbox{ such that }\sqrt{2T_{\sfy,\wt r_{l}}^{\npp, t}}\le \alpha(l)\Big\}
\]
(which is $\mathcal{B}_{t,l}^{\beta,k}$ of \cite[(3.57)]{BRZ2}), 
and using the notation of} \cite[(3.60)]{BRZ2} it suffices to bound 
\begin{align}\label{eq-P3}
\mathfrak{P}_{32} & := \Pbm \Big(\mathcal{B}_{t,l,j}\cap 
\Big\{\sqrt{2T_{\wt r_{l}}^{\npp, t }}>\alpha(l)+\frac{j}{2}\Big\}\,\big|\, \DD_{t,l}(j) \Big) \,.
\end{align}
To this end,  for $l > 5 \log n_{\vep}$ and $j \le 12 n_\vep$ we proceed as before, 
apart from omitting $\CC_{0,t,l}$,  up to the top line of \cite[(3.72)]{BRZ2}.  In this case
one has that $l > \sfk + h_{l,j}$, 
so we can again further condition on $\GG^{\bf 0}_{\bar k}$ and utilize  \eqref{eq-thanks4lem}.  

Fixing hereafter $\wh{z}$,  $l \in [\sfk,5 \log n_\vep]$ and $j \in [ l^{\gak}/4,12l]$,  
for \blue{$t=t_{\wh{z}}$ and} $\sqrt{2t'} \in I_{\al(l)+j}$ we
decompose the event $\DD_{t,l}(j)$ in \eqref{eq-P3} according to the 
value $t'$ of $T_{l-2}^{\npp,t}$. 
  Thereby,  we 
condition on having precisely $t'$ of the $(l-2)$-level excursions around $\npp$ 
completed during the first $t=(\sfm_\vep+\wh{z})^2/2$ (outer)
excursions from $\partial \bB_\sfd(\npp,r_1)$ to $\partial \bB_\sfd(\npp,r_0)$.  These inner excursions 
thus correspond to the measure $\Pbm_{t,b,b}$ of \Cref{prop:coupling} with
$R=r_1$ and at the values $k=l-3$,
$t=t_{\wh{z}}$ and $b:=t'-t$,  
which we then compare with the unconditional law
$\Pbm^\star_{t'}$ of the $t'$ first inner excursions.  To this end, 
we proceed to verify that $(t,l-3,5j,b) \in \Theta_{\ka_0}$ of \eqref{dfn:Thetak} for some universal $\ka_0>0$.  Indeed,
we assume \abbr{wlog} that $\sfk \ge 3 + \ka_0^{-3}$ and that $n_\vep$ is large enough so that
$\sqrt{t}  \in [n_\vep,2n_\vep]$.  Then,  for some $\ka_0>0$,  both 
$5 j \le 60 l \le 300 \log n_\vep \le \ka_0^{-1} \log t$ and
$|b| \le (\sqrt{2t} + j) | j-\sfrho_\vep l-l^{\gak}-\wh z\,| \le  \ka_0^{-1} (l \vee 5j) \sqrt{t}$.
With all the conditions of \Cref{prop:coupling} satisfied,  we see from  \eqref{eq:coupling} that with  
some universal $\bar c >1$ and a single event $\EE^{t'}_{t,l}(j)$ 
of probability at most $C l e^{-5j}$,  given that $\{T_{l-2}^{\npp,t}=t'\}$,  
we match all but at most 
\begin{equation}\label{eq:N-def}
N := 35 {\bar c} \ka_0 j \sqrt{t_{\wh{z}}} \,,
\end{equation}
of our conditional inner excursions with the first such $t'$ unconditional inner excursions. 
Noting that,
\[
\mathfrak{P}_{32} \le C l e^{-5j} + \sup_{\sqrt{2t'} \in I_{\al(l)+j}} \Pbm_{t,b,b} \Big(\mathcal{B}_{t,l,j}\cap 
\Big\{\sqrt{2T_{\wt r_{l}}^{\npp, t}}>\alpha(l)+\frac{j}{2}\Big\}\cap \EE^{t'}_{t,l}(j)^{\complement}\,  \Big) \,,
\]
we follow the argument leading to \cite[(3.71)]{BRZ2},  now with the law
$\Pbm_{t,b,b}$ instead of conditioning on $\CC_{0,t,l} \cap \DD_{t,l}(j)$
while intersecting each event there with $\EE^{t'}_{t,l}(j)^{\complement}$.  With $l e^{-j}$ bounded uniformly over
$j \ge  l^\gak/4$,
we arrive,  as in \cite[(3.74)]{BRZ2},  to $\mathfrak{P}_{32} \le C e^{-4j}$,  upon replacing 
\cite[Lemma 3.11]{BRZ2} by the following claim:\\
For small $\wt{c}>0$,  
any $m \le \log (2/\blue{\beta_{l,j}})$,  $\sqrt{2t'} \in  I_{\alpha(l)+j}$ and
$\wt \sfy,  \sfy \in \bB_{\sfd}\left(\npp,\wt{c} \,  r_{l}\right)$, 
with $\tilde q:=\sfd \left(\sfy, \wt \sfy\right)/r_{l} \ge \frac{\wt{c}}{4} \beta_{l,j}$,
\begin{equation}\label{eq-mod-458}
\Pbm_{t,b,b} \Big(
\big| 
\trav{t'}{\npp}{r_{l-3}}{r_{l-2}}{\wt \sfy}{\wt r_{l-1}^+}{\wt r_l}
-\trav{t'}{{\bf 0}}{r_{l-3}}{r_{l-2}}{\sfy}{\wt r_{l-1}^+}{\wt r_l}
 \big|\ge d_{0} j  \, m \, 
\sqrt{2 t' \tilde q} \,,\; \EE^{t'}_{t,l}(j)^{\complement} \Big) \leq  3 e^{-4(j+ m)}  \,.
\end{equation}
\blue{Turning to show \eqref{eq-mod-458},  we set}
$\tilde{R}=r_{l-3}$,  $\tilde{r} = r_{l-2}$,  $R=\wt r_{l-1}^+$ and $r=\wt r_l$
(which satisfy \eqref{eq-ratiobounds}),  \blue{noting}
that 
the tail bound \eqref{basicexp2} applies for the contribution to the difference 
within \eqref{eq-mod-458} from the $N$ possibly \emph{non-matched excursions} 
$\partial \bB_{\sfd}(\npp,\tilde{r})\to\partial \bB_{\sfd}(\npp,\tilde{R})$.  Thus,  using 
the matching to $\Pbm_{t'}^\star$ up to $N$ inner excursions
(which holds on $\EE^{t'}_{t,l}(j)^{\complement}$),  we get \eqref{eq-mod-458} by combining 
the unconditional version of \eqref{eq-thanks4lem} at $u=\sqrt{2t'}$,  with the bound of 
\eqref{basicexp2} at $\varphi= {d_{0}}{2} j  \, m \,  \sqrt{2 t' \tilde q}/2$ and $N$ of \eqref{eq:N-def}.
Indeed,  with $\tilde q \le 2 \wt c$ and $t'/t_{\wh{z}} \to 1$ for $\vep \to 0$,  clearly
$\varphi \ge 2 c_3 \tilde q N$ for small enough $\wt c (c_3 \bar c \ka_0/d_0)$.
Further,  since $j,  l = o(n_\vep)$ we have that $2 t' \tilde q \ge {\wt c}t' \blue{\beta_{l,j}} /2= \Omega(n_\vep)$,
resulting with $\varphi = \Omega (\sqrt{n_\vep})$.  Similarly $N = \Omega(n_\vep)$, 
while $j$ and $m$ are at most $O(\log n_\vep)$,  hence $\varphi \wedge N \ge 8 (j+m)$ and 
the proof \blue{of \eqref{eq-mod-458}} is thus complete. 
\end{proof}

\begin{proof}[Proof of \Cref{gprop:MainGoalEasyRegime}]
We assume \abbr{wlog} that $\sfk^\gak \ge 10 z_-$ (otherwise \eqref{g17.0} holds trivially).  Setting 
$L' := n_\vep - (10 \log n_\vep)^{1/\gak}$ we first establish \eqref{g17.0} while restricting in \eqref{eq:JayThm2.2eq}
to the ``not too large" levels $j \in [\sfk, L']$.  To this end,  
utilizing \eqref{r1.1} allows us,  as in \cite[(3.11)]{BRZ2},  to intersect 
the event in \eqref{eq:JayThm2.2eq} with 
\[
\AA_{z,d \sqrt{z_0}} := \Big\{ \tau_{\sfy} (t_{z-d\sqrt{z_0}}) \le \tau_{\sfx} (t_z) \le \tau_{\sfy}(t_{z+d\sqrt{z_0}}),  \;\;
\forall \sfx,\sfy \in \bF_{n_\vep} \Big\} \,,
\]
where taking here $z_0=|z|+\sfk^{\gak}$ (instead of $z_0=z$ of \cite[(3.11)]{BRZ2}),  guarantees having 
$\Pbm(\AA_{z,d\sqrt{z_0}}^{\complement})$ much smaller than the bound of  \eqref{g17.0}.  We thus
follow the proof in \cite[Section 3.1]{BRZ2},  replacing \cite[(3.15)]{BRZ2} by 
$\wh z := z - d \sqrt{z_0}$,  to deduce that it suffices to show  for some 
$\uu{\theta} =\uu{\theta} (l,n_\vep)$,  $\bar \theta=\bar \theta (l,n_\vep)$ and any $l \in [\sfk,L']$,  
\begin{equation}
\Pbm(\sqrt{2T_{l-2}^{{\bf 0},t_{\wh z}}} \le \alpha(l)+\uu{\theta}) + 
\sum_{j=\uu{\theta}}^{\bar \theta-1} \Pbm(\BB_{t_{\wh z},l} \cap \DD_{t_{\wh z},l}(j) )
+ \Pbm(\sqrt{2T_{l-2}^{{\bf 0},t_{\wh z}}} \ge \alpha(l)+\bar \theta) 
\label{eq:brz-331}
 \le 
c \,
e^{-2 \wh z  - 2 l - l_{n_\vep}^{\gak}}. 
\end{equation}
To this end,  we set $\uu{\theta}=l_{n_\vep}^{\gak}/4$ and $\bar{\theta}=12l$ 
up to $l=5\log n_\vep$,  taking thereafter $\bar{\theta}=12 n_\vep$.  Then,  using \eqref{17.1sv} and 
bounding $\Pbm(\DD_{t_{\wh z},l}(j))$ via \cite[(3.25)-(3.26)]{BRZ2},  
we arrive at the stated bound of \eqref{eq:brz-331},  for the sum over $j$ values.
Further,  $\sqrt{2 t_{\wh z}} \ge \al(l)+\frac{1}{4} l_{n_\vep}^\gak$ 
(since by assumption $\wh z \ge -\sfk \sfrho_\vep$),  with the
left-most term in \eqref{eq:brz-331} handled as in \cite[(3.27)]{BRZ2}.  
Similarly,  with $\sqrt{2 t_{\wh z}} \le \al(l) + 4 ( l \vee 5 \log n_\vep)$,  the right-most term of \eqref{eq:brz-331}
is also suitably controlled,  this time by \cite[(3.26)]{BRZ2} (upon noting that 
$ (\bar \theta-\wh z - l_{n_\vep}^{\gak})^2 /4\ge 2 \bar \theta + l$ whenever
$\bar \theta \ge 12 (l \vee \wh z)$).

Turning to deal with the ``large levels" beyond $L'$,  consider now  as in \cite[(3.38)]{BRZ2},  the events 
\begin{equation}
  \label{defG}
\blue{\GG_{z,l}} :=\bigcap_{l' \in [\sfk, l]} \Big\{ \min_{\sfv \in \bF_{n_{\vep}}\cap \bB_{\sfd} \left({\bf 0},\wt{c}r_0 \right)}
\big\{ \sqrt{2T^{\sfv, t_z}_{l'}} \big\} \ge \alpha(l')\Big\} 
\end{equation}
(not to be confused with the $\sigma$-algebras of \Cref{def:Fvjs}).
By a union bound over the finite cover $\bF_0$,  it suffices for \eqref{g17.0},  to show that when  $l=n_\vep-\sfk$,
\begin{equation}\label{eq-thanks1}
\Pbm(\GG_{z,l}^{\complement}) \le  c e^{-2(z-d \sqrt{|z|})} e^{-\frac{1}{5} \sfk^{\gak}} \,.
\end{equation}
We have already proved \eqref{eq-thanks1} for $l = L'$,  so by 
following the proof of \cite[Proposition 3.4]{BRZ2},  we complete the  proof of  \Cref{gprop:MainGoalEasyRegime},
upon establishing our next lemma (which replaces \cite[Lemma 3.9]{BRZ2}).
\end{proof}
\begin{lemma}
  \label{lem-westjet}
In the setting of \Cref{gprop:MainGoalEasyRegime},  for $\QQ$ of \eqref{g17.0b},  some finite $c$
and all $l \in (L', n_{\vep}-\sfk]$, 
\begin{equation}
\Pbm\left(\GG_{z,l}^\complement\cap\GG_{z,l-2}\cap\QQ \right)\le c (1+z_+) e^{-2z} e^{-l_{n_{\vep}}^{\gak}} .
\label{Gls}
\end{equation}
\end{lemma}
\begin{proof} Following the argument leading to \cite[(3.51)]{BRZ2},  it suffices to show that 
\begin{equation}\label{brz-351}
\sum_{j=0}^{8 n_\vep} \Pbm (\BB_{t_z,l} \cap \wt \CC_{t_z,l}(j) \cap \DD_{t_z,l}(j) ) 
\le c (1+z_+)  e^{-2 z} e^{-2l - l_{n_{\vep}}^{\gak}}.
\end{equation}
To this end,  recall that $l\geq L'\geq 5 \log n_\vep$, hence $h_{l,j}\leq ((l-\sfk)/2)\wedge(j+ l_{n\vep})\wedge 2\log n_\vep$.
Thus, for $\sfk^\gak \ge 10 z_-$ and  $\sfk \le l_{n_\vep}$,    the barrier estimate of \Cref{lem-bendd} yields that
\begin{align*}
\Pbm\left(\wt\CC_{t_{z},l} (j) \cap\DD_{t_{z},l}(j)\right) = \blue{\KQ_{[\sfk,l-h_{l,j}]}^{(1,l-2)}(z,j-(l-2)_{n_\vep}^\gak)} 
&
  \leq  C e^{-2l-2z-2l_{n_{\vep}}^{\gak}} e^{2j}  (1+\sqrt{\sfk}+z_{+} ) (1+j +\sqrt{l_{n_\vep}})\,.
\end{align*}
With $L' \ge 5 \log n_\vep$,  the proof of \eqref{brz-351} is completed
by combining the preceding bound with that of  \eqref{cb.15y} and summing over $j$.
 \end{proof}

\subsection{Angle and excursion windows: Proof of \Cref{theo-pre220}}\label{sec:Jay-26}
Writing for $\sfIp_j$ of \eqref{dfn:windows-alt},  
\begin{equation}
\boxminus_{\sfv,z, j}:=\{\whetavjp \in \sfIp_j \},  \qquad
\boxminus_{\bb,z; j_1,\ldots, j_m}:=\bigcap_{i=1}^m \boxminus_{\sfv_\bb,z, j_i},\label{14.4}
\end{equation}
recall (see \eqref{def.n-vep-i}, \eqref{def:UUpk} and \eqref{dfn:bw-Jay}) 
that for $\sfJup=[\sfKo-3,\nvf] \cup \{\nvfv\}$,
\begin{align} \label{a216.02}
\whUUpk = \Rkbp \cap \whUUp \quad \text{and} \quad 
\BWp = \boxminus_{\bb,z;  \sfJup} \bigcap_{k \in \sfJwp} \{ \Lkp \in \bB^\infty_{\bar \mu} (\boetak) \} 
 \,.
\end{align}
(\blue{where} throughout,  the notation $\boxminus_{\bb,z;  \sfJup}$ 
\blue{stands for} $\boxminus_{\bb,z;\sfKo-3,\ldots,\nvf,\nvfv}$).
\Cref{theo-pre220} is a direct consequence of our next lemma,  as \eqref{lem-first-moment} is a 
stronger bound than \eqref{atot-new},  while in view of \eqref{a216.02},
the bound \eqref{atot} follows from the other four statements of the lemma.  \blue{Recall 
\Cref{rmk:z-lz-adjust} that we consider here a path starting outside 
$\bB_{\sfd}(\sfv_\bb,r_1)$,  allowing the use of \Cref{lem-off-center-applies} and various 
barrier estimates from \Cref{sec:BoundaryCrossing}}.
\bl\label{theo-win220}
Set $n'=\nvf$,  $n''=\nvfv$ and fix $0<3 \gak < \kapo < \frac{1}{4}$,  $\okappa > \kapp > \frac{1}{2}$.  Then,  
for $\cst$ large enough,  $z$ fixed and $\hg(z)= C_{\sfk} (1+z_+) e^{-2z}$ for some $C_{\sfk}$ finite, 
 \begin{align}
\varlimsup_{\vep \to 0} \Big\{ e^{2 n''} \sup_{\bb} 
\Pbm \big( \whUUpk \cap \boxminus^{\complement}_{\sfv_\bb,z,n'} \big) \Big\}  
&= o_{\sfKfv}(1) \hg(z)\,,
\label{agar.67a}\\
  \varlimsup_{\vep \to 0} \Big\{ e^{2 n''} \sup_{\bb} 
 \Pbm \big(\Rkbp  \cap \boxminus_{\sfv_\bb,z,n'} \cap \boxminus^{\complement}_{\sfv_\bb,z,n''}\big)\Big\}
 &=o_{\sfKfv}(1) \hg(z)\,, \label{win5-good}\\
 \varlimsup_{\sfKfv, \vep} \Big\{ e^{2 n_\vep} \sup_{\bb,\sfv \in \bF^\bb} \Pbm \big(\whUUvk \cap 
\boxminus_{\bb,z,n'} \big) \Big\} & \le \hg(z) \,,
\label{lem-first-moment}\\
\label{lem-win1}
 \varlimsup_{\vep \to 0} \Big\{ e^{2 n''} \sum_{j=\sfKo-3}^{n'-1} \sup_{\bb} 
 \Pbm \big( \whUUpk  \cap \boxminus_{\bb,z,n'} \cap \boxminus^{\complement}_{\sfv_\bb,z,j}
   \big) \Big\} &=  o_{\sfKfv,\sfKz}(1) \hg(z)\,,\\
\varlimsup_{\vep \to 0}  \Big\{ e^{2n''} \sum_{k=\sfKo}^{\sfKw} \sup_{\bb} \Pbm \big(
\whUUp \cap \boxminus_{\bb,z; \sfJup^{(k)}} \cap  \{ \Lkp \notin \bB^\infty_{\bar \mu} (\boetak) \}
 \big) \Big\} & = o_{\sfKz}(1) \hg(z),
 \label{pawfud0.1}   
\end{align}
for the subset $\sfJup^{(k)}:=\{ k-3,[k,\sfKt],[\sfKth,n'],n''\}$ of $\sfJup$
(and with the larger $\sfI_{n'}$,  $\sfI_{n''}$ of \eqref{dfn:windows} in \eqref{lem-first-moment}).
\el 
\begin{remark} Recalling \eqref{def.nvep}
that there are at most $e^{2(n''+1)}/(r_0^2 \sfMo)$ blocks in a city,  it follows from \eqref{def:UUpk},
\eqref{agar.67a} and \eqref{lem-first-moment} that for some $C_{\sfk}$ finite,  any $\sfk \ge 1$,  $z \le 0$, 
$\dlc$,  $\srs$ and $\sfKfv$ large enough, 
\begin{align}\label{first-moment-patch}
\varlimsup_{\vep \to 0} 
\Pbm \big( \bigcup_{u=\sfM+1}^{\sfMo}\bigcup_{\bb \in u} \bigcup_{\sfv \in \bF^\bb}
\whUUvk
\big)  & \le  C_{\sfk} \,  \srs^2 \,  r_0^{-2} e^{-2z},
\end{align}
where the union is over all 
cities which are excluded from $\CCd$ (see 
\Cref{rem:good-ic}),
and we assume here re-starting the Brownian motion sufficiently further away for
\Cref{rmk:z-lz-adjust} to apply to all those cities.
\end{remark}

\begin{proof} We establish \eqref{agar.67a}--\eqref{pawfud0.1}, respectively,  in parts (a)--(e) below,
modifying \abbr{wlog} the barrier $\gamma_{-1}(j)$  
used in the events $\Rkvpp$ of \eqref{a216.0},  as explained in \Cref{rmk:no-vep-inJs}.\\
(a).  Recall \eqref{def:UUpk} that 
\[
\whUUpk = \bigcup_{\sfv\in \bF^\bb} \Big\{ \Rkvpp,  \eta_{\sfv,n_\vep+1}^{\bb,\langle z \rangle} = 0 \Big\}
\bigcap \Rkbp \subset \bigcup_{\sfv\in \bF^\bb} 
\Big\{ \eta_{\sfv,n_\vep+1}^{\bb,\langle z \rangle} = 0\Big\} \bigcap \Rkbp
\]
hence by a union bound over the $c e^{2 \sfKfv}$ possible values of $\sfv \in \bF^\bb$,  it suffices
for \eqref{agar.67a} to show that
\begin{equation}\label{Jay-agar.67a}
\varlimsup_{\vep \to 0} \Big\{ e^{2 n_\vep} \sup_{\bb,\sfv \in \bF^\bb} 
\Pbm \big( \{ \eta_{\sfv,n_\vep+1}^{\bb,\langle z \rangle} = 0 \} \cap \Rkbp \cap  \boxminus^{\complement}_{\sfv_\bb,z,n'} \big) 
\Big\}  = o_{\sfKfv}(1) \hg(z)\,.
\end{equation}
To this end,  \blue{setting $\AA^{\langle z \rangle}_{[\sfk,n'],\bb}(m) := \{ \Rkbp,  T^{\langle z \rangle}_{\bb,n'}=m\}$
we claim that for large $\sfKfv$ any $z$,  $\sfk$,  $\vep$ and $\bb$,
\begin{equation}\label{eq:pre-winx.2}
\sup_{\sfv \in \bF^\bb} \Pbm \(  \eta_{\sfv,n_\vep+1}^{\bb,\langle z \rangle} = 0  \, |\,   
\AA^{\langle z \rangle}_{[\sfk,n'],\bb}(m)  \) \leq 
 e^{-\frac{m}{\sfKf} (1-\frac{2}{\sfKf})} \,. 
\end{equation}
Indeed,  we are here in the setting of \Cref{lem-off-center-applies},  for $\sfu=\sfv_\bb$,  
$\sfl=n'$,  $R_i=r_i$,  $R=r_{n'}$,  $\sfs=\sfs_1(z)$,  
$\AA^{(m)}=\AA^{\langle z \rangle}_{[\sfk,n'],\bb}(m)$ and $\rho = r_{n''}/R = \sfKfv^{-5}$.  
The \abbr{lhs} of \eqref{eq:pre-winx.2} is of the form $\PRsb (\Gamma|\AA^{(m)})$ for
the event $\Gamma:=\{\eta_{\sfv,n_\vep+1}^{\bb,\langle z \rangle} = 0\} \in \NN_{\sfv,\vep,m,0}$
where our $m$ excursions from $\partial \bB_{\sfd}(\sfv_\bb,r_{n'})$ 
to $\partial \bB_{\sfd}(\sfv_\bb,r_{n'-1})$ avoid $\bB_{\sfd}(\sfv,\vep/e)$.
Hence,  by \eqref{eq:off-center-bd} followed by \eqref{eq:pv-u} at $d_\star=1$ and $r=\vep/R=e^{-\sfKf}$, 
the \abbr{lhs} of \eqref{eq:pre-winx.2} is at most 
\[
(1-\inf_{\sfz,\sfv} \{ q^{\sfz,\sfv_\bb}_{r,\sfv}\})^m 
\le  \Big(1- \frac{1- c \rho}{\sfKf+2}\Big)^m \le e^{-m\frac{1-c \rho}{\sfKf+2}} \le 
e^{-\frac{m}{\sfKf} (1-\frac{2}{\sfKf})} 
\]
(with the right hand side holding once $c \rho \,  \sfKf^2 \le 4$).  Next,  recall \Cref{rmk:no-vep-inJs} and \eqref{a216.0} that 
$m=\frac{1}{2}(2 \sfKf+y)^{2}$ for $y=\whetanp \ge -\sfKf^{\gak}$.  For all $\sfKfv$ large enough  
\[
\inf_{y \ge -\sfKf^{\gak}} \Big\{ \frac{(2\sfKf+y)^{2}}{2\sfKf} (1-\frac{2}{\sfKf}) 
- 2(\sfKf+y) - \frac{y^2}{4\sfKf} \Big\} 
\]
is uniformly bounded below,  hence \eqref{eq:pre-winx.2} implies that for} some 
universal $C$,  all $\vep$,  any $z$,  $\sfk$,  $\sfKfv$ and $\bb$,
 \begin{equation}
\sup_{\sfv \in \bF^\bb} \Pbm \(  \eta_{\sfv,n_\vep+1}^{\bb,\langle z \rangle} = 0  \, |\,   \Rkbp,   T^{\langle z \rangle}_{\bb,n'}  \) \leq 
C  e^{-2\sfKf}  e^{ -2 \whetanp} e^{-(\whetanp)^2/(4 \sfKf)} \, .
\label{winx.2}
\end{equation} 
\blue{Now,  in view of \eqref{14.4} at $j=n'$ and \eqref{winx.2},  by the union bound it suffices 
for \eqref{Jay-agar.67a} to show} that 
\begin{align}
e^{2 n'} \!\!\!\!\!\! \sum_{y\in [-\sfKf^\gak,\sfKf^{\kapo}] \cup [\frac{1}{2} \sfKf^{\kapp},\infty)} \!\! \sup_\bb \Pbm ( \Rkbp,  \whetanp \in I_y) e^{-2 y} e^{-y^2/(4 \sfKf)} &= o_{\sfKfv}(1) \hg(z)\,,
\label{winx.1}
\end{align}
where $I_y=[y,y+1)$ and throughout this proof summation over $y$ in some range means the summation over $y$ integer in that range.
To this end,  from \blue{\eqref{a216.0} and \Cref{lem-bendc}}
we have that
\begin{equation}\label{agar.5old}
\Pbm ( \Rkbp,  \whetanp \in I_y) \blue{= \KQ_{[\sfk,n']}^{(1,n')}(z,y)}
\le \hg(z) e^{-2 n'} \sfKf^{-1/2} (1+y+\blue{\sfKf^\gak})  e^{2y} \,,
\end{equation}
which yields \eqref{winx.1} since for $\gak \le \kapo < \frac{1}{4}$ and $\kapp>\frac{1}{2}$,
\begin{equation}\label{agar.66-67}
2 \sfKf^{-1/2} \sum_{j=-\sfKf^{\gak}}^{\sfKf^{\kapo}} |j|
+ \sum_{j= \frac{1}{2} \sfKf^{\kapp}}^\infty  j e^{-j^2/(4 \sfKf)} 
 = o_{\sfKf}(1)\,.
\end{equation}
(b).  Recall  that excursion counts of concentric balls with a constant ratio of radii 
form a Markov chain whose transition kernels are $q_{\ell}(x',x)$ of \Cref{dfn:Markov}.
Upon using \eqref{agar.5old} we thus arrive at \eqref{win5-good} once we show that 
for $\sfLw:=\sfKf-\sfKfv$,
\begin{equation}\label{win5-good-alt}
\sum_{y\in \sfIp_{n'}} (1+y)  e^{2y}  
\sum_{j\geq \sfKfv/2} \sup_{x' \in I_{2 \sfKf+y}} \{ q_{\, \sfLw}(x',I_{2\sfKfv+j-1}) \} =o_{\sfKfv}(1) e^{-2 \sfLw} \,. 
\end{equation}
As in \eqref{eq:BRZ-1.10m},  for some $C$ finite,  $c \ge 0$,  $\sfLwp=\sfLw+c$ and uniformly in 
$y,j$ appearing in the sums of \eqref{win5-good-alt}
\begin{equation}\label{eq:BRZ-1.10mod}
\sup_{x' \in I_{2 \sfKf+y}} \{ q_{\,\sfLw}(x',I_{2\sfKfv+j-1}) \}
\le \frac{C}{\sqrt{\sfLwp}} e^{-\frac{(j-y-2\sfLwp)^2}{2\sfLwp}} \le \frac{C}{\sqrt{\sfLwp}} e^{-2\sfLwp-2y} e^{-\frac{j^2}{16\sfLwp}} 
\end{equation}
(since with $\kapp<1$,  once $\sfKfv$ is large enough we have that  $j-y \ge \frac{1}{2} j \ge 16 \sfLwp$). 
Summing the right hand side of \eqref{eq:BRZ-1.10mod}
over all integer $j \ge \frac{1}{2} \sfKfv$ and bounding by the appropriate Riemann integral 
\blue{we get by standard tail} estimates on the Gaussian error function \blue{that for large enough $\sfKfv$,}
\begin{align}\label{win5-bad}
\blue{\sum_{j\geq \sfKfv/2} \sup_{x' \in I_{2 \sfKf+y}} \{ q_{\, \sfLw}(x',I_{2\sfKfv+j-1}) \} \le } 
8 C e^{-2\sfLw-2y} e^{-\frac{\sfKfv^2}{128 \sfLw}} \,.
\end{align}
Plugging 
\eqref{win5-bad} 
in \eqref{win5-good-alt} and 
\blue{having} $\sfLw \le O(\log \sfKfv)$,  \blue{yields,  as claimed,  that}
the \abbr{lhs} of \eqref{win5-good-alt} \blue{is at most}
\[
\blue{8 C e^{-2\sfLw} e^{-\frac{\sfKfv^2}{128 \sfLw}} \sum_{y\in \sfIp_{n'} \cap \bZ} (1+y) \le}   
8 C  \sfKf^{\,  2 \kapp} e^{-\frac{\sfKfv^2}{128 \sfLw}}  e^{-2 \sfLw} = o_{\sfKfv}(1) e^{-2 \sfLw} \,.
\]
\noindent
(c).  In view of \eqref{def:UUvk},  \eqref{dfn:windows-alt} and \eqref{14.4} we have that 
\[
\whUUvk \cap \boxminus_{\bb,z,n'} \subset 
\Big\{  \Rkvpp,  \eta_{\sfv,n_\vep+1}^{\bb,\langle z \rangle} = 0,
 \Rkbp,  \whetanp \in \sfI_{n'} \Big\} \,.
\]
By \eqref{agar.5old},  it thus suffices for \eqref{lem-first-moment},   
to show that
\begin{align}\label{first-moment-alt}
\varlimsup_{\sfKfv,\vep} \Big\{ \frac{e^{2\sfKf}}{\sqrt{\sfKf}} \sum_{y \in \sfI_{n'}\cap \bZ} y  e^{2y} \sup \{ f_m : \sqrt{2m} \in  I_{2\sfKf+y} \} \Big\}
& \le C_\sfk  \,,  
\end{align}
where \blue{for $\AA^{\langle z \rangle}_{[\sfk,n'],\bb} (m)$ of \eqref{eq:pre-winx.2}},
\begin{align}\label{dfn:f-y}
f_m & := \sup_{\sfv \in \bF^\bb}   
 \Pbm\Big( 
 \Rkvpp,  \eta_{\sfv,n_\vep+1}^{\bb,\langle z \rangle} = 0
   \,\big|\,  \AA^{\langle z \rangle}_{[\sfk,n'],\bb} (m) \Big)  \,.
\end{align}
\blue{Setting $\wt{n}:=\wnvfv$,  the conditional on $\AA^{(m)}:=\AA^{\langle z \rangle}_{[\sfk,n'],\bb}(m)$ probability 
in \eqref{dfn:f-y} is of the form $\PRsb (\Gamma|\AA^{(m)})$ of \Cref{lem-off-center-applies}, 
with same values of $\sfu$,  $\sfs$,  $\sfl$,  $R_i$,  $\rho$ and $\AA^{(m)}$ as in 
our derivation of \eqref{eq:pre-winx.2},  while having here  
$\Gamma:=\{\eta_{\sfv,n_\vep+1}^{\bb,\langle z \rangle} = 0\} \cap \Rkvpp \in (\NN_{\sfv,r})^m$
for $r=e \,  r_{\wt{n}}$ (see \eqref{a216.0}).  Setting $k= \tsfKfv^4$,  
similarly to our derivation of \eqref{eq:bd-hQv},  
we further split the event $\Gamma$ to the disjoint union of $\Gamma^{(k)} := \Gamma \cap \{ \eta^{\bb,\langle z \rangle}_{\sfv,\wt{n}} \le \sqrt{2k} \} \in \NN_{\sfv,r,m,k}$ 
and $\Gamma \setminus \Gamma^{(k)} \subset \{T^{\sfv,\wt{n},k}_{n_\vep+1} = 0 \}$.  From 
\eqref{eq:off-center-bd},  given our specification of $k$ and $\rho$ and having $m = O(\sfKf^2)$,  it follows that  
for $\sfKfv$ large enough
\begin{equation}\label{eq:rem-angles}
\sup_{\sfv \in \bF^\bb} \{ \PRsb (\Gamma^{(k)}|\AA^{(m)}) \} \le (1+c \rho)^{m+k} (\bP^{1}_R)^m(\Gamma^{(k)}) \le 
2 (\bP^{1}_R)^m (\Gamma) \,,
\end{equation}
where on the right hand side we have the probability of $\Gamma$ for $\sfv=\sfv_\bb$ and $m$ 
\abbr{iid} unconditioned,  stopped Brownian paths.  Setting \abbr{wlog} 
$\sfv=\npp$ and scaling space by $e^{-2}/R$ we deduce similarly to \eqref{eq:A-RU-Q-hatQ} that 
\[
2 (\bP^1_R)^m(\Gamma) \le 2 \bar \bP^{(\star)}_{u} ( \wQs_{\npp}) 
\]
for $u = \sqrt{2m}-2\sfKf$,  $\sfKp=\sfk$ and $\wQs_\npp$ of \eqref{def:QQs}.   
Further,  as in the derivation of \eqref{eq:bd-hQv},  
the probability of hitting $\bB_{\sfd} (\sfv,\vep/e)$ by a Brownian path
from $\partial \bB_{\sfd} (\sfv,r_{\wt{n}})$ to $\partial \bB_{\sfd} (\sfv, e \,  r_{\wt{n}})$ is,  
uniformly over its starting and end points,  at least $c/\tsfKfv$ for some 
$c>0$.  Since $\{ T^{\sfv,\wt{n},k}_{n_\vep+1} = 0 \}$ requires that $\tsfKfv^4$ such paths 
all miss this ball,  we deduce from \eqref{ubd-wQ},  that}
\begin{equation}\label{eq:sup-fm}
\sup \{ f_m : \sqrt{2m} \in  I_{2\sfKf+y} \} \le  e^{-c \,  \tsfKfv^3} + 2 
\sup_{u \in I_y} \bar \bP^{(\star)}_{u} ( \wQs_{\npp}) 
\le c_\sfk \frac{y}{\sfKf} e^{-2 \sfKf} e^{- 2y}  e^{-y^2/(8 \sfKf)} \,.
\end{equation}
By the preceding,  the \abbr{lhs} of \eqref{first-moment-alt} is therefore,  as claimed,  at most 
\[
c_\sfk \sup_{\sfKf \ge 1} \Big\{ \sum_{y \in \bZ} y^2 \sfKf^{-3/2}  e^{-y^2/(8 \sfKf)} \Big\} := C_\sfk < \infty \,.
\]
(d).  By our reasoning in part (c),  we get \eqref{lem-win1} upon showing that for 
all large $\sfKz$,  $\sfKfv$ and any $y \in \sfIp_{n'}$
\begin{equation}\label{eq:updn-alt}
\sup_{\bb,\vep} \Big\{
 e^{2 n'} \sum_{j=\sfKo-3}^{n'-1} \sum_{u \notin \sfIp_j} \Pbm ( \Rkbp,  \whetajp \in I_u,
\whetanp \in I_y) \Big\} 
= o_{\sfKfv,\sfKz}(1)  \hg(z) \sfKf^{-1/2} y e^{2y} \,,
\end{equation}
where throughout \blue{this} proof,  the summation is over \blue{integer values} of $u$ in the indicated range.
Starting with $j \in [\sfKo-3,\frac{1}{2} n_\vep]$,  since $z=o(\sqrt{j}\,)$
we deduce from the bound of 
\blue{\Cref{lem-bendc} for $\KQ_{[\sfk,n']}^{(1,j,n')}$ (see \Cref{rmk:straight-line}), }
that for any \blue{$u \in [- j^{\gak},20 \sqrt{n_\vep j}]$} and $y \gg \sfKf^{\gak}$, 
\begin{align}\label{agar.45p} 
e^{2 n'}  \Pbm ( \Rkbp,  \whetajp \in I_u,   \whetanp  \in I_y) 
\blue{=e^{2 n'} \KQ^{(1,j,n')}_{[\sfk,n']}(z,u,y)} 
& \leq \hg(z) e^{2y}\,  \frac{y}{\sqrt{\sfKf}} \,  \frac{(1+u+j^{\gak})^2}{j^{3/2}} e^{-\frac{(z-u)^{2}}{8j}} 
\nn\\ & 
\leq  \hg(z) \sfKf^{-1/2} y e^{2y} f_{j,0}(u) \,,
\end{align}
where we set for $s \in [0,j)$ and $u \in \bZ$, 
\begin{equation}\label{def:fjs}
f_{j,s}(u):= \frac{(j^{\gak} +u_+)^2}{(j-s)^{3/2}} e^{- \frac{u^{2}}{16(j-s)}}.  
\end{equation}
\blue{We ignore hereafter the aggregate contribution from larger values of $u$,  which is negligible
(see \eqref{eq:neglig-A}).  Moreover,  in view of \eqref{agar.45p},}
the contribution of $j \in [\sfKo-3,\sfKo] \cup \{\sfKt,\sfKth\}$ has a negligible effect on \eqref{eq:updn-alt},  
since similarly to \eqref{agar.66-67},  for $\gak \le \kapo < \frac{1}{2}<\kapp$ and any such $j$, 
\begin{equation}\label{agar.66-new}
\sum_{u=- j^{\kapo}}^{2 j^{\kapo}} f_{j,0}(u) + \sum_{u \ge \frac{1}{2} j^{\kapp}} f_{j,0}(u) \le 
\frac{2}{j^{3/2}} \Big[ \sum_{u=0}^{2 j^{\kapo}} (j^{\gak}+u)^2
+ \sum_{u = \frac{1}{2} j^{\kapp}}^{\infty} u^2 e^{-\frac{u^2}{16j}}  \Big]  = o_{\sfKo}(1)\,.
\end{equation}
The same applies to the contribution to \eqref{eq:updn-alt} by all $j \le \frac{1}{2} n_\vep$,  since for 
$\gak<\frac{1}{6}$ and $\okappa>\frac{1}{2}$, 
\begin{equation}\label{agar.66-ext}
\sum_{j=\sfKo}^{\frac{1}{2} n_\vep} \Big[
\sum_{u=- j^{\gak}}^{2 j^{\gak}} f_{j,0}(u) + \sum_{u \ge \frac{1}{2} j^{\okappa}} f_{j,0}(u) \Big] \le 
C_{\gak,\okappa}  \sum_{j=\sfKo}^{\infty} j^{3\gak-3/2} = o_{\sfKz}(1) \,.
\end{equation}
Since $\okappa >\kapp$ and $\kapo>\gak$ we have that $2 (u-y)^2 \ge u^2 + y^2$
for any $u \notin \sfIp_j$ and $y \in \sfIp_{n'}$,  hence for $j \in [\frac{1}{2} n_\vep,n')$ 
and any such integer $y$ and $u \ge - j_{n_\vep}^\gak$,  
we similarly \blue{deduce from \Cref{lem-bendc} and \Cref{rmk:straight-line}} that
\begin{align}\label{agar.45jj} 
e^{2 n'} \Pbm ( \Rkbp,  \whetajp \in I_u,   \whetanp  \in I_y)
\blue{=e^{2 n'} \KQ^{(1,j,n')}_{[\sfk,n']}(z,u,y)} 
& \leq 
\hg(z)   \,
\frac{(1+u+j_{n_\vep}^{\gak})^2}{(n'-j)^{3/2}}  e^{2y}\,  \frac{y}{\sqrt{\sfKf}} 
e^{-\frac{(u-y)^2}{8(n'-j)}} 
\nn\\
& \leq \hg(z) \sfKf^{-1/2} y e^{2y} f_{j_{n_\vep},\sfKf}(u) 
e^{-\frac{y^2}{16(n'-j)}} \,.
\end{align}
With $f_{j,s}(\cdot) \le 4 f_{j,0}(\cdot)$ when either $j \ge 2s$ or $u \ge C \sqrt{j}$,  the aggregate contribution to 
\eqref{eq:updn-alt} from $u \ge  j_{n_\vep}^{\okappa}/2$ and from $|u| \le 2 j_{n_\vep}^{\gak}$ 
and $j \leq n'-\sfKf$,  be as in the right hand side of \eqref{agar.66-ext} except for having $\sfKf$ instead of $\sfKo$.  
Moreover,  thanks to the term $\exp(-{y^2}/{(16(n'-j))})$,  the aggregate contribution from 
$|u| \le 2 j_{n_\vep}^{\gak}$ and $n'-j \in [1,\sfKf]$ is for $y \ge \sfKf^{\kapo}$ and $\kapo>3\gak$ at most
\[
(3 \sfKf^{\gak})^3 \sum_{\ell=1}^{\sfKf} \ell^{-3/2} e^{-\frac{\sfKf^{2\kapo}}{16\ell}} = o_{\sfKfv}(1), 
\]
which thereby completes the proof of  \eqref{eq:updn-alt}.

\noindent
(e).   Enumerating the angle differences $\theta_q \in \Ipi$ 
between the start and end-point of those excursions,  we set
\begin{equation}
\WW^{\ff,\in m}_{\bb,k}(\sfs)=\Big\{ \,  \sfd_\infty\Big(\frac{1}{\sfs} \sum_{q=1}^\sfs \delta_{\theta_q},\bar \mu\Big) 
\in   \frac{c_{0}}{2 \sqrt{\sfs}} I_m \Big\} \, .  \label{eq:LBTruncatedSumzf}
\end{equation}
Recall that $\boetak=({\cst}\log k)/{\sfm_\vep} $,   so for any $\cst \ge c_0$,  one has that  
\begin{equation}\label{eq:inclusion}
\Big\{ \frac{1}{\sfs} \sum_{q=1}^{\sfs} \delta_{\theta_q} \notin \bB^\infty_{\bar \mu} 
\Big(\frac{\sfm_\vep}{2 \sqrt{\sfs}} \boetak \Big)  \Big\} \subset \bigcup_{m \ge \log k} \WW^{\ff,\in m}_{\bb,k}(\sfs) \,,
\qquad  \forall k,  \sfs \,.
\end{equation}
Let $N_{k,a}=[\sfs_k(a+1)]$ denote the largest possible number of $k$-level excursions
in case $\whetap_{\bb,k} \in I_a$.  In \eqref{pawfud0.1},
$\Lkp$ are the empirical measures of $\{ \theta_q\}_{q=1}^{\sfs}$,  with  
$\sfs \in [\sfs_k(a),N_{k,a}]$ in case $\whetap_{\bb,k} \in I_a$.  Since $\sfs_k(a) \ge N_{k,a}-\sqrt{2 N_{k,a}}$ and 
\[
\lim_{\vep \to 0} \;\; \Big\{ \sup_{k \in [\sfKo,\sfKw],a \in \sfI_k}  \big|\sfm_\vep^{-2} N_{k,a} - \frac{1}{2}\big| \Big\}  = 0 \,,
\]
it thus follows from \eqref{eq:dinf-comp} and \eqref{eq:inclusion} that for all $\vep$ small enough 
\begin{equation}\label{eq:Lkp-WNa}
\Big\{ \Lkp \notin \bB^\infty_{\bar \mu} (\boetak) \Big\} 
\subset \bigcup_{m \ge \log k} \WW^{\ff,\in m}_{\bb,k}(N_{k,a}) , 
\qquad \forall k \in [\sfKo,\sfKw],  \quad \forall  a \in \sfIp_k \,.
\end{equation}
Consequently,  we control the event in \eqref{pawfud0.1} by 
adapting the proof of \cite[(4.42)]{BRZ2}
though we have the event $\whUUp  \cap \boxminus_{\bb,z;n',n''}$ at our finer 
scales $[n',n_\vep]$,  instead of the easier setting of \cite[(4.61)]{BRZ2},  where $n'=n_\vep$.
Due to the angle information in the events of \eqref{eq:LBTruncatedSumzf},
we lack
the radial symmetry of the events on scales $[k,n']$,
but $\boxminus_{\bb,z;\sfKt,\sfKth}$ implies at most $3 n_\vep^2$ 
excursions at levels $\sfKt$ and $\sfKth$.  \blue{So,  partitioning our events 
per $\{ \whetap_{\bb,\sfKt} \in I_{j'},\whetap_{\bb,\sfKth} \in I_{j''}\}$},  with 
$\sfKth-\sfKt=5 \log n_\vep$ and $\GG^{\sfv_\bb}_{\sfKt}$ 
the $\sigma$-algebra of excursions 
from $\partial \bB_\sfd(\sfv_\bb,r_{\sfKt-1})$ to $\partial \bB_\sfd(\sfv_\bb,r_{\sfKt})$,
we get as in \cite[(4.66)]{BRZ2},  that for some $c$ finite,  uniformly in $\bb$ and 
our range of $j'$, $j''$,  $y$,   
\begin{equation}\label{eq:BRZ-466}
\bP\Big(\boxminus_{\bb,z; [\sfKth,n']} \cap \{ \whetanp \in I_y \} \,\Big| \,  T^{\bb,\sfKt,s}_{\sfKth},
\GG_{\sfKt}^{\sfv_\bb}  \Big) \le c \;
\bP\Big(\boxminus_{\bb,z; [\sfKth,n']} \cap \{ \whetanp \in I_y \} \,\Big| \,  T^{\bb,\sfKt,s}_{\sfKth} \Big) \,.
\end{equation}
With the event on the right hand side of \eqref{eq:BRZ-466} having radial symmetry,  we thus 
combine the right hand side of \eqref{eq:sup-fm} with a barrier bound on the right hand side of \eqref{eq:BRZ-466}
\blue{utilizing \Cref{prop:BarrierSecGWPropUB} for $\KQ_{[\sfKth,n']}^{(\sfKth,n')}(j'',y)$}.
Summing the resulting bound 
over $\sfv \in \bF^\bb$ and $y \in \sfIp_{n'}$, 
we conclude,  analogously to \cite[(4.67)]{BRZ2},  that uniformly over $j'' \in \sfIp_{\sfKth}$, 
\begin{equation}\label{eq:BRZ-467}
e^{2 n''} \sup_{\bb} \bP\Big(\whUUp \cap \boxminus_{\bb,z; [\sfKth,n'],n''} \,\Big| \,  \whetap_{\bb,\sfKth} \in I_{j''}, 
\GG_{\sfKt}^{\sfv_\bb}  \Big) 
\le \hg(j'') 
e^{2 \sfKth} 
\sum_{y \in \sfIp_{n'}} \,  \frac{y^2}{\sfKf^{3/2}} e^{-\frac{y^2}{8 \sfKf}} = \hg(j'') e^{2 \sfKth} \,.
\end{equation}
To complete the proof of \eqref{pawfud0.1},  we utilize the following lemma (which is proved at the end of this section).
\begin{lemma}\label{lem:BRZ-CH4}
Fixing $z \in \bR$,  for any $\bb$,  $k \in [\sfKo,\sfKt]$,  $m \ge 0$ and integer $a$,  $j'$,  $j''$,  set 
\begin{align}\label{def:Lmpa} 
\LL'_{k,m,a} &:=
\{ \whetap_{\bb,k} \in I_{a} \} \cap
\WW^{\ff,\in m}_{\bb,k}(N_{k,a}) 
\\
\bB^{j',j''}_{\bb,k+3,\sfKt,\sfKth} &:= \bigcap_{j=k+3}^{\sfKt-1} \{ \whetajp \ge 2 j^\gak \}
\cap \{ \whetap_{\bb,\sfKt} \in I_{j'} \} \cap \{ \whetap_{\bb,\sfKth} \in I_{j''} \} \,,
\label{def:Bjp-jpp}
\end{align}
where $I_y=[y,y+1)$.
Then,  for some $c'>0$,  $C<\infty$,  all $\vep$ small enough,  any $k \in \sfJwp$,  $m \ge \log k$ and $\bb$, 
\begin{equation}\label{eq:BRZ-top}
\sum_{j''}  \hg(j'') \bP(\bigcup_{p,a,j'} \LL'_{k,m,a}  \cap \bB^{j',j''}_{\bb,k+3,\sfKt,\sfKth} ) \le C (1+z_+) e^{-2z} 
e^{-2 \sfKth} e^{-c' m^2}  \,,
\end{equation}
where the union is over 
$a \in \sfIp_k$,  $j' \in \sfIp_{\sfKt}$ and the sum 
over $j'' \in \sfIp_{\sfKth}$.
\end{lemma}
Indeed,  we consider which of the events in \eqref{def:Lmpa} and \eqref{def:Bjp-jpp} occur in
$\boxminus_{\bb,z; \sfJup^{(k)}} \cap  \{ \Lkp \notin \bB^\infty_{\bar \mu} (\boetak) \}$
(compare with \cite[(4.55) \& (4.69)]{BRZ2}),  and 
combine \eqref{eq:BRZ-467} with \eqref{eq:BRZ-top},  then 
sum over $m \ge \log k$ and $k \in \sfJwp$ to arrive at \eqref{pawfud0.1}.
\end{proof}

\begin{proof}[Proof of \Cref{lem:BRZ-CH4}] Even conditional upon 
$\whetap_{\bb,k}$,   the $\{\th_q\}$ in \eqref{eq:LBTruncatedSumzf} are \abbr{IID} 
from a law $\bar \mu$ whose density is bounded below (see \eqref{dfn:Poiss}).  
By \eqref{lbl:d-inf-bd} we thus have,  analogously to \cite[(4.54)]{BRZ2},  the estimate 
\begin{equation}\label{eq:BRZ-454}
\bP\Big( \WW^{\ff,\in m}_{\bb,k}(N_{k,a}) | \whetap_{\bb,k} \Big) \le 2 e^{-m^2/4}    
\end{equation}
(provided $c_0 \ge 6 \sqrt{2} \,  c$,   $c = \sqrt{1+\bar c^2}+\bar c$). 
Combining \eqref{eq:BRZ-454} and \cite[(4.56)]{BRZ2} (at $k$),  yields that 
\begin{equation}\label{eq:BRZ-455}
\bP(\LL'_{k,m,a}) \le c_1 e^{-2k} e^{-2(z-a)} e^{-(z-a)^2/(4k)} e^{-m^2/4} \,,
\end{equation}
for some $c_1$ finite,  fixed $z$ and all $k \in [\sfKo,\sfKt]$,  $a \in \sfIp_k$ and $m \ge \log k$.
By our choice of $\sfKth$, 
for some $C<\infty$ the bound \eqref{eq:BRZ-454} already implies \eqref{eq:BRZ-top} unless
$m \le C (\log n_\vep)^{1/(2 \gak)}$,  which we thus assume hereafter.  Moreover, 
as $\sfd_\sfw(\cdot,\cdot)  \le \sfd_\ff(\cdot,\cdot)$ (see \eqref{dfn:dw}),  our events 
$\WW^{\ff,\in \log k}_{\bb,k}(\sfs)$ imply the events $\WW_{\sfv_\bb,k}(\sfs)$ of \cite[(4.19)]{BRZ2}
and replacing $\log k$ by $m$ throughout \cite[Sec.  4.5]{BRZ2},  we see that 
\cite[Lemmas 4.7 and 4.8]{BRZ2} hold here as well\footnote{Throughout \cite{BRZ2} one has the barrier
$\gamma_1(\cdot)$ with parameter $\gak={1}/{4}$ (see \cite[(4.14)]{BRZ2}),  
whereas we must use $\gak<{1}/{12}$ (e.g.  part (d) relies on it).  However,    
one can check that all the proofs in \cite[Section  4.5]{BRZ2} work for such choice of $\gak$,
namely,  the height parameters $a,j',j''$ can start at a lower power exponent than in \cite{BRZ2},
at the cost of having a power $\wt{c}=\wt{c}(\gak) > {5}/{2}$ in \eqref{eq:BRZ-471}.}
with the event $\LL'_{k,m,a}$ of \eqref{def:Lmpa} instead of $\LL'_{k,m,p,a}$ of \cite[(4.53)]{BRZ2},
except that unlike
\cite[(4.59)]{BRZ2},  since here  $\sfKth-\sfKt \ne \sfKt - k$ we do not use the 
simplification step from \cite[(4.134)]{BRZ2} to \cite[(4.135)]{BRZ2}.  As  
a result the expressions in \cite[(4.71),(4.74)]{BRZ2} are respectively replaced by
\begin{align}\label{eq:BRZ-471}
G_{k,a,j',j''}(m) &=  (a+m^{\wt{c}}) e^{-2a} e^{-2 (\sfKth-k)}
\!\!\!\!\! \sum_{\stackrel{\{\bar j' : |\bar j' - j' | \le 2 M_0 m\}}{\{\bar j'' : |\bar j'' - j'' | \le 2 M_0 m}\}}\!\!\!\!\!
\frac{\bar j'}{(\sfKt-k)^{3/2}} e^{-\frac{(\bar j'-a)^2}{4(\sfKt-k)}} \frac{e^{2 \bar j''}}{\sqrt{\sfKth-\sfKt}}
 e^{-\frac{(\bar j''-\bar j')^2}{2(\sfKth-\sfKt)}} \,,
\\
H_{k,a,j',j''}(m) &= e^{-2a} e^{-2 (\sfKth-k)}
\!\!\!\!\! \sum_{\stackrel{\{\bar j' : |\bar j' - j' | \le 2 M_0 m\}}{\{\bar j'' : |\bar j'' - j'' | \le 2 M_0 m}\}}\!\!\!\!\!
\frac{1}{\sqrt{\sfKt-k}} e^{-\frac{(\bar j'-a)^2}{2(\sfKt-k)}} \frac{e^{2 \bar j''}}{\sqrt{\sfKth-\sfKt}}
 e^{-\frac{(\bar j''-\bar j')^2}{2(\sfKth-\sfKt)}} \,.
\label{eq:BRZ-474}
\end{align}
Here $k \in \sfJwp=[\sfKo,\sfKw]$,  hence $\sfKt - k  \ge \frac{1}{2} \sfKt \gg a \vee (\sfKth-\sfKt)$
and $a \wedge (\sfKth-\sfKt) \gg 1$.  It is then easy to check that for some finite $c$,  $C=C(M_0,\wt{c})$,
and independent $\bar j' \sim N(a,2(\sfKt-k))$,  $\bar j'' -\bar j' \sim N(0,\sfKth-\sfKt)$,
\begin{align}
\sum_{j',j''} G_{k,a,j',j''}(m)  \hg(j'') e^{2 \sfKth} & \le c \, e^{C m} \,  a e^{-2a}  e^{2 k} 
\frac{1}{\sfKt-k}  \bE [ (1+|\bar j'|)(1+ |\bar j''|) ] 
\nonumber \\ & 
\le c^2 e^{C m} \,  a \big(1+ \frac{a^2}{\sfKt}\big) e^{-2a} e^{2 k} \,, 
\label{eq:BRZ-472}\\
\label{eq:BRZ-473}
\sum_{j',j''} H_{k,a,j',j''} (m) \hg(j'') e^{2 \sfKth} & \le c e^{C m} e^{-2a} e^{2 k} 
 \bE [ (1+|\bar j''| )] \le c^2 e^{C m} (a + \sqrt{\sfKt}) e^{-2a} e^{2 k} \,.
\end{align}
Combining \eqref{eq:BRZ-455} and \eqref{eq:BRZ-472},  we find that for some 
$c<\infty$,  $z$ fixed,  any $k \in \sfJwp$ and $m \ge \log k$,  
\begin{align}
\sum_{a,j',j''} G_{k,a,j',j''} (m) \bP(\LL'_{k,m,a}) \hg(j'') e^{2 \sfKth} &\le c^3  e^{-2z} e^{-m^2/5}
\sum_a   a \big(1+ \frac{a^2}{\sfKt}\big) e^{-(z-a)^2/(4k)} \nonumber \\
& \le c^4 \sqrt{k} (\sqrt{k}+z_+) \big(1+\frac{k+z^2}{\sfKt} \big) e^{-2z} e^{-m^2/6} \le \hg(z) e^{-m^2/8} \,.
\label{eq:BRZ-452}
\end{align}
The same bound applies with $H_{k,a,j',j''}(m)$ instead of $G_{k,a,j',j''}(m)$,  provided $m \ge 4 (\log \sfKth)^{1/2}$.
Thus,  analogously to \cite[(4.52)]{BRZ2},  we conclude that \eqref{eq:BRZ-top} holds for 
for some $c'>0$,  $C<\infty$,  all $\vep$ small enough,  any $k \in \sfJwp$, $m \ge \log k$ and $\bb$.
\end{proof}
 
\subsection{Traversal times based on inter-city excursions: Proof of 
\Cref{Jay-4.3-4.3rev}}\label{subsec-clock-adjust}

Recall the traversal times $\tau_\sfx(m)$ of \eqref{eq-taux},  the traversal counts 
$T^{[\tau,\tau']}_{\sfv,j}$ of \Cref{def-trav}
(for Brownian stopping times $\tau$ and $\tau'$), 
and the $\sigma$-algebras $\GG^{\sfx,\tau}_j$
of \Cref{def:Fvjs}.  Equipped with these notations,  our next two lemmas establish \eqref{A.2} and \eqref{eq:Jay4.14},  respectively,  whereby completing the proof of \Cref{Jay-4.3-4.3rev}.
\begin{lemma}\label{Jay-Lemma4.3}
For $\delta \in (0,1)$,  some $c=c(\delta)>0$,  any $b>a$ and all $\sfx$, 
\begin{equation}\label{A.2}
\bP \big( T_j^{\sfx, b}  - T_j^{\sfx,a} \ge \ell \,|\,  \GG^{\tau_{\sfx}(a)} \,\big) \le e^{-c \, \ell/j},  \qquad 
\forall j \ge 1,  \;\; \ell \ge (1+2 \delta) (b-a) \,.
\end{equation}
\end{lemma}
\begin{proof} By the strong Markov property of the Brownian motion at the time $\tau_{\sfx}(a)$ it took to 
complete $a$ excursions from $\partial \bB_\sfd(\sfx,r_1)$ to $\partial \bB_\sfd(\sfx,r_0)$,  the law of
$T_j^{\sfx, b}  - T_j^{\sfx,a}$,  given $\GG^{\tau_{\sfx}(a)}$,  matches that of $T_j^{\sfx,b-a}$,  which
by radial symmetry,  is the sum of $b-a$ \abbr{iid} copies of $T_j^{\sfx,1}$.  Noting that  
$\bP(T_j^{\sfx,1}=i)=j^{-2} (1-\frac{1}{j})^{i-1}$ for $i \ge 1$,  we deduce that $\bE[T_j^{\sfx,1}]=1$ with
$\bE[e^{\lambda T_j^{\sfx,1}/j}] \le e^{\lambda(1+\delta)/j}$ for small $\lambda>0$ and any $j \ge 1$.
By Markov's inequality,  we arrive at \eqref{A.2} (for $c={\lambda \delta}/{(1+2\delta)}$).
\end{proof}

\begin{lemma}\label{Jay-Lemma4.3-rev}
There exist $c>0$,  so if $j \le  n_\vep$ and $T^{[\tau,\tau']}_{\sfx,j}$ is $\GG^{\sfx, \tau}_j$-measurable,
then
\begin{equation}\label{eq:Jay4.14}
\sup_{\sfv \in  \bB_\sfd(\sfx, r_j)}
\bP( T^{[\tau,\tau']}_{\sfv,n_\vep+1}=0 \;\big| \; \GG^{\sfx,\tau}_j,  T^{[\tau,\tau']}_{\sfx,j} \ge \ell )
 \le e^{-c \,  \ell/(n_\vep-j)} \,,   \qquad \forall \ell \ge 1 \,.
\end{equation}
\end{lemma}
\begin{proof} Given $\GG^{\sfx, \tau}_j$ we know the starting points $\sfw_i \in \partial \bB_\sfd (\sfx, r_{j-1})$ 
and ending points $\sfz_i \in \partial \bB_\sfd (\sfx, r_{j})$ of the outer part of all 
excursions at level $j$ around $\sfx$,  which occur during $[\tau,\infty)$.  Since
$T:= T^{[\tau,\tau']}_{\sfx,j}$ is $\GG^{\sfx, \tau}_j$-measurable,  we have no information
about the inner parts of the subset of those excursions which are completed during $[\tau,\tau']$,  
beyond the number $T$ of such excursions and their corresponding starting and ending points.  
Setting $\rho_j:=e^{-(n_\vep-j+1)}$,  the \abbr{lhs} of \eqref{eq:Jay4.14} is thus bounded above by the probability 
that $\ell$ independent Brownian path from $\partial \bB_\sfd(\sfx,r_j)$ to $\partial \bB_\sfd(\sfx,r_{j-1})$, 
with starting and ending points $(\sfz_i,\sfw_{i+1})$, 
miss $B_1 := \bB_\sfd(\sfv, \rho_j r_j)$.  Set $B_2 := \bB_\sfd(\sfv, e^{-1} r_j)$ and $B_3:= \bB_\sfd(\sfv, r_j)$,
noting that $B_1 \subseteq B_2 \subset B_3 \subset \bB_\sfd(\sfx,2 r_j)$.  
Let $q^{\sfx,\sfv}_{\sfz,\sfw}(j)$ denote the probability of the event $\Gamma$ 
that a Brownian path starting at $\sfz \in \partial \bB_\sfd (\sfx,r_j)$ and ending at 
$\sfw \in \partial \bB_\sfd (\sfx,r_{j-1})$ visits $\partial B_2$ before $\partial \bB_\sfd (\sfx, 2 r_j)$
and thereafter enters $B_1$ before exiting $B_3$,  using $q^{\sfx,\sfv}_{\sfz}(j)$ for such probability when 
$\sfw$ is not specified.  Since a Brownian path from $\sfy \in \partial B_2$ enters $B_1$ before exiting $B_3$ 
with probability
 \[
q (j) := \frac{1}{n_\vep-j+1}\,,
\]
we get \eqref{eq:Jay4.14} (with $c=\frac{1}{2} c_1 c_2$),  upon finding $c_i>0$ for $i=1,2$,  such that 
uniformly over $\sfz$,  $\sfw$,  $\sfx$,  $\sfv$ and $j$, 
\begin{equation}\label{eq:unif-bd}
q^{\sfx,\sfv}_{\sfz,\sfw} (j) \ge c_1 \,  q^{\sfx,\sfv}_{\sfz}(j) \,, \qquad 
q^{\sfx,\sfv}_{\sfz}(j) \ge c_2 \, q (j) \,.
\end{equation}
The quantities in \eqref{eq:unif-bd} are about planar Brownian motion,  so by 
translation invariance and Brownian scaling we can fix hereafter $\sfx=\npp$ and $r_j=1$.  
Doing so,  the probabilities $q^{\sfv}_{\sfz}$ that the
path from $\sfz$ hits $\partial \bB(\sfv,e^{-1})$ before exiting $\bB(\npp,2)$,  are continuous 
in $(\sfz,\sfv)$ and positive for $\|\sfz\|=1$ and $\|\sfv\| \le 1$.  Hence,
\[
c_2 := \inf_{\stackrel{\|\sfz\|=1}{\|\sfv\| \le 1}}  \{ q^{\sfv}_{\sfz} \} > 0 \,,
\] 
and the right-side of \eqref{eq:unif-bd} follows,  since by strong Markov
$q_{\sfz}^{\npp,\sfv}(j) \ge q_\sfz^\sfv \,  q (j)$.  Next,  the Poisson kernel $p(\sfy,\sfw)$  
for the exit measure of $\bB(\npp,e)$ is such that 
\[
c_1 := \inf_{\stackrel{\|\sfz_i\|=2}{\|\sfw\|=e}}  \;\;\; \Big\{ \frac{p(\sfz_1,\sfw)}{p(\sfz_2,\sfw)} \Big\} > 0 \,.
\]
With $\Gamma$ measurable on the Brownian path up to its first exit from $\bB(\npp,2)$,  
restricted to $\sigma(\Gamma)$,  the Radon-Nikodym derivative 
due to conditioning on exiting $\bB(\npp,e)$ at $\sfw$,  is at least $c_1$,  implying the
left-side of \eqref{eq:unif-bd}.
\end{proof}

\subsection{Towards the mean number of unvisited blocks}\label{sec:good-cover-end}

\begin{proof}[Proof of \Cref{lem:good-cover-end}] 
From \eqref{eq:g-def} we have that uniformly on $\sfI_{\sfKo} \times \sfIw$,  
\begin{equation}\label{eq:gstb-z}
\gstb(y,v)=(1+o_{\sfKz}(1)) \gstb(y-z,v) e^{-\sfco z} \,.
\end{equation}
Setting the non-negative
\begin{equation}\label{def:X-eps-pre}
\bar Q^{(1,\vep)}_{u,z} := 
\,  \alpha (\sfKfv,\sfKp) \bKz \,  e^{-\sfco z}
\avsum_{\bb  \in u} \bE \Big[ \gstb \big( \Ybz-z, \Vbz \big) ; \Ybz \in \sfI_{\sfKo},
\LKoz \in \bB^\infty_{\bar \mu} ( \boetasKo) 
\,\big|\,  \HHk  \Big] \,,
\end{equation}
we deduce,  in view of \eqref{eq:wtQ1-bd} and \eqref{eq:gstb-z}, 
that as claimed in \eqref{eq:eqv-wt-Q},  
\[
\wt{Q}^{(1,\vep)}_{u,z}
= 
(1+o_{\sfKfv,\sfKz,\vep}(1)) \bar Q^{(1,\vep)}_{u,z} \,.
\]
In view of \eqref{eq:Ball-to-angle},  only $\Vbz \in \sfIw$ can contribute to $ \bar Q^{(1,\vep)}_{u,z}$,
so comparing \eqref{def:X-eps} with \eqref{def:X-eps-pre} we deduce that
$\bar Q^{(1,\vep)}_{u,z} \le \sfX^{(\vep)}_{u,z} e^{-\sfco z}$ 
(since the indicator on $\sfI_{\sfKo} \times \sfIw$ is bounded above by $\gst(\cdot)$ of \eqref{eq:g-def-alt}).
The proof of \Cref{lem:good-cover-end} is thus concluded upon establishing \eqref{eq:wX-small}.  To this end,
note that $\gst(\cdot) \le \gstb(\cdot)$ with equality on $\sfI_{\sfKo} \times \sfIw$ and fixing 
$\kapom < \kapo$,  $\kappp > \kapp$,  $\cst$ and $z$,  for all large enough $\sfKo$, 
\[
\gst(y-z,v) >0 \quad \Longrightarrow \quad 
\gstb(y-z,v) e^{-\sfco z} \le 2 \gstb(y,v) \,,
\] 
while $\gst(\cdot)$ is zero outside 
\begin{equation}
\label{def:DI}
\tsfI \times \sfIwe := [\sfKo^{\kapom},\sfKo^{\kappp}] \times[-3\cst \log \sfKo,3\cst \log \sfKo] \,.
\end{equation}
For such $\sfKo$,  we thus have that 
\begin{equation}\label{eq:bd-wX}
\sfX^{(\vep)}_{u,z} e^{-\sfco z} - \bar Q^{(1,\vep)}_{u,z} 
\le 2 \alpha (\sfKfv,\sfKp) \bKz \Big[   
\avsum_{\bb  \in u} \wX^{(\vep,\sfwm)}_{\bb,z} + \avsum_{\bb \in u} \wX^{(\vep,\Cap)}_{\bb,z} \Big]\,,
\end{equation}
where for any $\bb \in u$,  
\begin{align}
\wX^{(\sfwm)}_{\bb,z} &:= 
\bE \Big[ \gstb \big( \Ybz, \Vbz \big); 
\Ybz \in \sfI_{\sfKo},
\Vbz \in \sfIwe,
\LKoz \notin \bB^\infty_{\bar \mu} ( \boetasKo) 
\,\big|\, 
 \HHu 
\Big] \,,
\nonumber
\\
\wX^{(\vep,\Cap)}_{\bb,z} &:= \bE \Big[ \gstb \big( \Ybz, \Vbz \big); \Ybz \in \tsfI \setminus \sfI_{\sfKo},
\Vbz \in \sfIwe
\,\big|\,  
 \HHu
 \Big] 
\end{align}
(and utilizing \Cref{rmk:z-a-vep},  we have changed to conditioning on the larger $\sigma$-algebra $\HHu$).
In view of \eqref{def:bar-bKz} \blue{and having $e^{-2 \sfKz} \sfM$ uniformly bounded bellow},
 it thus suffices to fix $\delta>0$ and show that 
\begin{equation}\label{eq:wX-small-alt}
\varlimsup_{\sfKz,\vep} 
\bP \Big[
\avsum_{\bb} \wX^{(\vep,\sfwm)}_{\bb,z} > \delta e^{-2 \sfKo} \Big] = 0 \,, 
\qquad 
\varlimsup_{\sfKz,\vep} \bP \Big[
\avsum_{\bb} \wX^{(\vep,\Cap)}_{\bb,z} 
  > \delta e^{-2 \sfKo} \Big] = 0 \,,
\end{equation}
where now we average over all blocks.
Moreover,  by \Cref{Thm4.1+4.2-Jay} we can and shall intersect the variables in \eqref{eq:wX-small-alt} 
with the events $\Akbz$ of \eqref{def:Jay-Thm4.1},  which are measurable on $\HHu$ when $\bb \in u$.
Upon doing so,  by Markov's inequality,  \eqref{eq:wX-small-alt} follows from having that for any fixed
$\sfk$,
\begin{align}\label{eq:wX-sfw-small}
\varlimsup_{\sfKz,\vep} 
 e^{2 \sfKo} \sup_{\bb} 
& \bE \Big[ \gstb\big(\Ybz,\Vbz\big); \Akbz,  \Ybz \in \sfI_{\sfKo},\Vbz \in \sfIwe,
\LKoz \notin \bB^\infty_{\bar \mu} ( \boetasKo)  
\Big] = 0 \,,\\
\varlimsup_{\sfKz,\vep} 
 e^{2 \sfKo} \sup_{\bb} 
& \bE \Big[ \gstb \big( \Ybz, \Vbz \big); \Akbz,  \Ybz \in \tsfI \setminus \sfI_{\sfKo},  \Vbz \in \sfIwe \Big] = 0 \,. 
 \label{eq:wX-Cap-small}
\end{align}
Recall \eqref{def:Jay-Thm4.1-cons} that $\Akbz$ implies having for $z'=z-\log \sfk$ and $z''=z+\log \sfk$,
\begin{equation}\label{eq:Akbz-cons}
T_{\sfKo}^{\sfv_\bb,\sfs_1(z')} \le \TKoz \le  T_{\sfKo}^{\sfv_\bb,\sfs_1(z'')} \quad \text{and} \quad \whetaKpp \le \Ybz \le \whetaKppp\,.
\end{equation}
In addition,  $2(\sfs_1(z'')-\sfs_1(z')) \le C_{\sfk} n_\vep$
for $C_\sfk$ finite,  all $\sfKz$ large and any $\vep$ (see \eqref{eq:b-a-bd}),  whereupon by \eqref{A.2-n}
\begin{equation}\label{A.2-alt}
\sup_{\bb} \bP\big( T_{\sfKo}^{\sfv_\bb,\sfs_1(z'')} - T_{\sfKo}^{\sfv_\bb,\sfs_1(z')} \ge C_\sfk n_\vep 
\big) \le 
e^{-c \, C_\sfk n_\vep/\sfKo} = o_\vep(1) \,.
\end{equation}
We set $\Vbp  =  \Vb(\LKopp)$ for $\Vb(\cdot)$ of \eqref{dfn:hat-L-alpha}.  Then,  noting that for some 
$c'=c'(\kapp,\cst)$ finite
\begin{equation}\label{eq:bd-gstb}
\gstb(y,v) \le \sfKo^{c'} e^{-2y} \le \sfKo^{c'}  \qquad \forall y \in \tsfI,  v \in \sfIwe \,,
\end{equation}
by \eqref{eq:Akbz-cons}--\eqref{eq:bd-gstb},  we can and will
replace the event $\Akbz$ on the \abbr{lhs} of \eqref{eq:wX-sfw-small}-\eqref{eq:wX-Cap-small} by 
\begin{equation}\label{eq:Akbz-alt}
n_\vep^{-1} \big( \TKoz - T_{\sfKo}^{\sfv_\bb,\sfs_1(z')} \big) \in [0, C_{\sfk}]\,,   \quad 
\Ybz - \whetaKpp \in [0, \frac{1}{2} C_{\sfk}]\,,
\quad 
|\Vbz  - \Vbp| 
\le 2 \bar c \, C_{\sfk}  \,.
\end{equation}
Under \eqref{eq:Akbz-alt},  for large $\sfKz$ it follows from \eqref{eq:dinf-comp} 
(similarly to our reasoning for \eqref{dfn:bb-extra}),  that also
\begin{align*}
\{ \Ybz \in \sfI_{\sfKo},  \LKoz \notin \bB^\infty_{\bar \mu} ( \boetasKo) \}  &\subset 
\{ \whetaKpp \in \tsfI,  \LKopp \notin \bB^\infty_{\bar \mu} (\boetaKo) \} \,.
\end{align*}
Further,  then having $\Ybz \in \tsfI \setminus \sfI_{\sfKo}$ implies the same for $\whetaKpp$
upon slightly reducing the parameters $\kapom$ and $\kapp$.
Finally,  with $y \mapsto \gstb(y,v)$ 
decreasing on $\sfI_{\sfKo}$ (for all large enough $\sfKz$),  we further 
replace $\gstb(\Ybz,\Vbz)$ in \eqref{eq:wX-Cap-small} by $\gstb(\whetaKpp,\Vbp)$ and in
\eqref{eq:wX-sfw-small} utilize \eqref{eq:bd-gstb} to replace it instead by 
$\sfKo^{c'} e^{-2 \whetaKpp}$.  

Implementing all those modifications,  our claims \eqref{eq:wX-sfw-small}-\eqref{eq:wX-Cap-small} are 
a direct consequence of 
\begin{align}\label{eq:wX-sfw-alt}
\varlimsup_{\sfKz,\vep} 
\sfKo^{c'} e^{2 \sfKo} 
& \bE \Big[ e^{-2  \whetaKpp};  \whetaKpp \in \tsfI,
\LKopp \notin \bB^\infty_{\bar \mu} (\boetaKo)  
\Big] = 0 \,,\\
\varlimsup_{\sfKz,\vep} 
 e^{2 \sfKo} 
& \bE \Big[  \whetaKpp e^{-2 \whetaKpp} e^{-2 \Vbp}; \whetaKpp \in \tsfI \setminus \sfI_{\sfKo} \Big] = 0 
 \label{eq:wX-Cap-alt}
\end{align}
(and we have dropped the supremum over $\bb$,  since $(\whetaKpp,\LKopp)$ have 
the same distribution at all blocks).  

Proceeding to confirm \eqref{eq:wX-sfw-alt}-\eqref{eq:wX-Cap-alt},  recall that 
conditional on completing precisely $\sfs' =\sfs_{\sfKo}(\whetaKpp)$ 
excursions at level $\sfKo$ during the first $\sfs_1(z')$ 
level one excursions,  $\LKopp$ is the empirical measure of $\sfs'$ \abbr{iid} 
samples from the law $\bar \mu$.  Thus,  similarly to \eqref{eq:BRZ-454} we have 
for some $c>0$ any $\cst \ge c_0$ and $\whetaKpp \in \tsfI$
\begin{equation}\label{eq:BRZ-454-new}
\bP\Big( 
\LKopp \notin \bB^\infty_{\bar \mu} (\boetaKo)  \,\big|\,  \whetaKpp \Big) \le 2 e^{-c (\log \sfKo)^2} \,.
\end{equation}
Likewise,  given $\whetaKpp$,  the variable $2 \Vbp$ is the sum of $\sfs'$ centered,  $[-1,1]$-valued,  \abbr{iid} 
variables multiplied by the factor $\bar c \,  \sfm_{\vep}/\sfs'$.  Hence,  with 
$\sup\{ \sfm_{\vep}^2/\sfs_{\sfKo}(y) : y \in \tsfI \}$
bounded uniformly in $\sfKz$ large and $\vep$ small,
\begin{equation}\label{eq:bd-Vbp}
\bE \Big[ e^{-2 \Vbp} \, \big| \,  \whetaKpp \Big] \le e^{(\bar c \, \sfm_{\vep})^2/(2 \sfs')} \le C
\end{equation}
for some $C$ finite.
Next,  applying
the Gaussian tail bounds of \cite[Proposition  1.4]{BRZ},  for $\whetaKpp$
and $\vep \to 0$,  we deduce that 
for $I_u=[u,u+1)$ and some $C=C(z')$ finite, 
\begin{equation}\label{eq:clt-bd-pKo}
\varlimsup_{\vep \to 0} \bP(\whetaKpp \in I_u) \le C \,  \sfKo^{-1/2} e^{-2(\sfKo-u)} e^{-u^2/(2 \sfKo)} \,,
\qquad \forall u \in \tsfI
\end{equation}
(as $z'=O(1)$ and $\sfKo \gg 1$).  Clearly,  \eqref{eq:clt-bd-pKo} implies having
for $2\kapo < \frac{1}{2} < \kapp$ and any $c>0$,  that 
\[
\sum_{u = \sfKo^{\kapom}}^{\sfKo^{\kappp}} e^{2(\sfKo-u)}
 \bP(\whetaKpp \in I_u)
 \Big(2 \sfKo^{c'}
e^{-c (\log \sfKo)^2} + u {\bf 1}_{u \le \sfKo^{\kapo}} + \sfKo^{\kappp} {\bf 1}_{u \ge \sfKo^{\kapp}} 
\Big) = o_{\sfKz,\vep}(1)  \,.
\]
We thus conclude by covering $\tsfI$ with intervals $I_u$ of integer end-points 
and plugging in \eqref{eq:wX-sfw-alt} and \eqref{eq:wX-Cap-alt}
our uniform bounds \eqref{eq:BRZ-454-new} and \eqref{eq:bd-Vbp},  respectively.
\end{proof}

\subsection{Tools for proving Propositions \ref{prop:2nd-moment} and  \ref{prop:tight}}\label{subsec:tools}

We start with a crude upper tail bound,  that we utilize in \Cref{sec:2nd-mom},
when proving \Cref{prop:2nd-moment}.
\begin{lemma}\label{lem:2mm-star}
Fix $\dlc$ and $\gak$,  $\okappa$ as in \eqref{dfn:aux-par}.  Then,  for any fixed $z$ and $\sfKfv$, 
\begin{equation}
\varlimsup_{\sfKz,\vep} \bP \Big( \sum_{u=1}^{\sfM} \wt{Q}^{(\star,\vep)}_{u,z} \ge 1 \Big) = 0 \,,
\label{eq:2mm-star}
\end{equation}
where 
\begin{equation}\label{dfn:Qst}
\wt{Q}^{(\star,\vep)}_{u,z} :=  \avsum_{\bb \in u}
\bE \big[ e^{-2 \YbKp}; \YbKp \in \sfI_{\sfKop}  | \FFuz \big] \qquad \hbox{with} \qquad \sfKop=\sfKo+3.
\end{equation}
\end{lemma}
\begin{proof} By \Cref{Thm4.1+4.2-Jay},  the event
$\Akz$ of \eqref{def:Jay-Thm4.1} holds up to decaying probability (of $C \sfk^{-4}$),  so it suffices to establish 
\eqref{eq:2mm-star} with $\Akz$ intersecting the event on its \abbr{lhs}.  Having done so,  since $\Akz$ is the intersection of events $\Akbz$ 
which are $\FFuz$-measurable whenever $\bb \in u$,  we can and will 
replace $\wt{Q}^{(\star,\vep)}_{u,z}$ in \eqref{eq:2mm-star} by
\[
\avsum_{\bb \in u} \sum_{i}
 e^{-2 (y_{(i)}-1)}
\bP \big[ \Akbz \cap \{ \YbKp \le y_{(i)} \} | \FFuz \big] \,,
\]
where $\{y_{(i)}\}$ forms a uniform partition of $\sfI_{\sfKop}$ to $\sfKo$ points.
With $\Akbz$ implying that 
$\YbKp \ge \Ybzp$ for $z'=z-\log \sfk$
(see \eqref{def:Jay-Thm4.1-cons}),  we further bound above this proxy for $\wt{Q}^{(\star,\vep)}_{u,z}$ by
\[
e^{2} \avsum_{\bb \in u} \sum_{i} e^{-2 y_{(i)}} \bP \big[  \Ybzp \le y_{(i)}  | \FFuz \big] \,.
\]
Plugging this bound on $\wt{Q}^{(\star,\vep)}_{u,z}$ into \eqref{eq:2mm-star},  we deduce 
by Markov's inequality that the probability on the \abbr{lhs} of \eqref{eq:2mm-star} is at most 
\begin{equation}\label{eq:418-alt}
e^{2} \,  \sfKo \,  \sfM\,  \max_{\bb,  y \in \sfI_{\sfKop}} \{ e^{-2 y} \bP (  \Ybzp \le y ) \} \,.
\end{equation}
Next,  recall  \eqref{eq-calgary},  that for any $z'$,  $y$ and $\bb$, 
\begin{equation}\label{eq:T-eta}
\bP \big( \Ybzp \le y \big) = \bP\big( T^{\sfv_\bb,\sfs_1(z')}_{\sfKop} \le \sfs_{\sfKop}(y) \big),
\end{equation}
where in view of \eqref{eq:Jay1.12},  the variable $T^{\sfv_\bb,m}_{\sfKop}$ is the sum $S_m$
of $m$ \abbr{iid} non-negative integer variables of the form \eqref{dfn:GJp}
for the parameter $q={1}/{\sfKop}$.  With $\sfKo \sfM  e^{-2 \sfKop} \to 0$ 
and the law of $S_m$ independent of $\bb$ and $\vep$,  we deduce in view of \eqref{eq:T-eta},
that the expression in \eqref{eq:418-alt} decays to zero with $\sfKz$,  if for some $c(z')$ finite
\begin{equation}\label{eq:amir-lem418-equiv}
\varlimsup_{\vep \to 0} \bP\big( S_{\sfs_1(z')} \le \sfs_{\sfKop}(y) \big)
\le c(z') e^{-2 \sfKop} e^{2y}\,,  \qquad \forall \sfKop,  \;\; y \in \sfI_{\sfKop} \,.
\end{equation}
Recall \eqref{eq:special-norm} and \eqref{eq:nominal} that $\sfs_1(z') \to \infty$ as $\vep \to 0$,  and
as $\bE[S_1]=1$,  ${\rm var}(S_1)=2(\sfKop-1)$ (see \eqref{eq:GJp-var}),  also an elementary calculation shows that
\[
\frac{ \sfs_{\sfKop}(y) - \bE[S_{\sfs_1(z')}]}{\sqrt{\hbox{var}(S_{\sfs_1(z')})}} \to
\frac{y-z' -2(\sfKop-1)}{\sqrt{\sfKop-1}} := \sft \,.
\]
Consequently,  by the \abbr{clt},  for any $z'$,  $y$ and $\sfKop$ (such that $\sft \le 0$),  
\[
\varlimsup_{\vep \to 0} \bP\big( S_{\sfs_1(z')} \le \sfs_{\sfKop}(y) \big) = 
\Phi(\sft) \le e^{-\frac{1}{2} \sft^2} \le e^{-2(\sfKop-1)} e^{2(y-z')}\,,
\]
as claimed in \eqref{eq:amir-lem418-equiv},  which thereby concludes our proof of \eqref{eq:2mm-star}.
\end{proof}

We proceed with key estimates for \eqref{ct2.25-new} and for the 
tightness of $\big( \sqrt{\smash[b]{\CCT}} - \osfm_\vep  \big)_+$
(see  \Cref{sec:tight}).
\begin{lemma}\label{lem:right-tight}
The  variables $\UUz$ of  \eqref{dfn:UUz} are such that 
\begin{equation}\label{eq:right-tight}
\lim_{z \to \infty} \varlimsup_{\sfKfv,\sfKz,\vep} \,\bP \big( \bigcup_{\bb} \UUz \big) = 0 \,.
\end{equation}
Further,  
the normalized $\vep$-cover time $\CCp$ of the square of edge-length $\srs$ centered at the Brownian motion 
starting point,  are such that 
\begin{equation}\label{eq:cov-p-lim}
\lim_{\srs \to 0} \varlimsup_{\vep \to 0} 
\Pbm \big( \sqrt{\CCp} > \osfm_\vep -  z \,  \big)  = 0  \qquad \forall z \in \bR \,.
\end{equation}
\end{lemma}
\begin{proof} 
By \eqref{dfn:UU-alt} and  \Cref{Thm4.1+4.2-Jay},  we get \eqref{eq:right-tight} upon showing that
for small $\wt{c}$,  $\srs=r_0$ and any fixed $\sfk \ge 1$,  
\begin{equation}\label{eq:81}
\lim_{z \to \infty} \varlimsup_{\sfKfv,\sfKz,\vep} \,\bP \big( \bigcup_{\sfv \in \bF} 
\{ \eta_{\sfv,n_\vep+1}^{[z]} = 0,  \Akz,  \Rkpp \} \big)  = 0 \,,
\end{equation}
where $z':=z-\log \sfk$ and the events $\Rkz$,  $\Akz$ are as
 in \eqref{def:Jay-Thm4.2} and \eqref{def:Jay-Thm4.1}.  To this end,  from \eqref{agar.67a},  \eqref{lem-first-moment} 
and the union bound,  it follows that 
\begin{equation}\label{eq:82}
\lim_{z \to \infty} \varlimsup_{\sfKfv,\vep} \,\bP \big( \bigcup_{\bb,\sfv \in \bF^\bb} \whUUvkp
 \big)  \le C_\sfk \lim_{z \to \infty} (1+z'_+) e^{-2 z'} = 0 \,.
\end{equation}
Further,  recall \eqref{def:Jay-Thm4.1-cons} and \eqref{def:UUvk} that 
\begin{equation}\label{eq:new-adj}
\{ \eta_{\sfv,n_\vep+1}^{[z]} = 0,  \Akz,  \Rkpp \} \subset \whUUvkp,  \qquad \forall \sfv \in \bF^\bb,
\end{equation}
hence \eqref{eq:81} is a consequence of \eqref{eq:82}.

Turning to \eqref{eq:cov-p-lim},  fix $\srs \le e^{-17}$,  setting 
$\dlc \le \frac{1}{18}$,  $\sfKz \ge \log (1/\srs)$,
$\sfxc:=\frac{1}{2} e^{-\sfKz}$ and $\sfxb:=\frac{1}{2} \vep e^{\sfKfv}$.  
\blue{Note that each point in 
$\bT^2$ is covered by the aggregate of the (non-extended) blocks within inner-cities 
over the sixteen tilings where,  first,  all city centers are simultaneously shifted by 
$(0,0)$,  $(\sfxc,0)$,  $(0,\sfxc)$,  $(\sfxc,\sfxc)$ and then all block centers in each of these city tilings 
are simultaneously shifted} by $(0,0)$,  $(\sfxb,0)$,  $(0,\sfxb)$ and $(\sfxb,\sfxb)$.
Let $\tau_0$ denote the first hitting time of a ball of radius $r_0=e^{-15}$ 
centered at the Brownian starting point,  with $\CCpp(i)$,  $1 \le i \le 16$ 
the normalized $\vep$-cover times by the Brownian path re-started at time $\tau_0$,
of the collection of blocks in the above $i$-th tiling which are in cities that
intersect the relevant square of edge length $\srs$ (around our starting point).  Obviously,  
\begin{equation}\label{eq:grid-bd}
 \CCp \le \tau_0 + \max_{i=1}^{16} \{ \CCpp(i) \}  \,.
\end{equation}
Fixing the $i$-th tiling,  by the strong Markov property and independence of the Brownian increments,  
$\CCpp(i)$ has the same law as with $\tau_0=0$,  except that now the path starts outside $\bB_{\sfd}(\sfv_\bb,r_1)$ for all blocks in cities $u \in (\sfM,\sfMo]$ (which intersected the square of 
edge length $\srs$ around the original starting point).  
Similarly to the explanation in \Cref{rmk:z-lz-adjust},  by moving to such a starting point we guarantee  
having the implication \eqref{def:Jay-Thm4.1-cons} for all the blocks in inner-cities $(\sfM,\sfMo]$.  
Taking now in $\Rkpp$ of \eqref{def:Jay-Thm4.2}
the intersection only over those blocks,  yields \eqref{eq:new-adj} and that the bound \eqref{eq:Jay-Thm-4.2} 
holds uniformly in $\dlc$ and $\srs \le r_2$.  Combining those facts at $z = z_\sfk = -\sfk^{\gak}/20$,
with \eqref{first-moment-patch} at $z' = z_\sfk - \log \sfk$,  
we deduce that for some finite $C$, 
$C'_{\sfk}
$,  any $\dlc$ and $\srs$  
\begin{align}\label{eq:cov-patch}
\varlimsup_{\sfKfv,\sfKz,\vep} & \Pbm \big( \CCpp > (\sfm_\vep+z)^2 \big)
\le \varlimsup_{\sfKfv,\sfKz,  \vep} 
\Pbm \Big(\bigcup_{u=\sfM+1}^{\sfMo} \bigcup_{\bb \in u} \bigcup_{\sfv \in \bF^\bb}
\{ \eta_{\sfv,n_\vep+1}^{[z]} = 0 \} \Big) \nonumber \\
&\le \varlimsup_{\sfKfv,\sfKz,  \vep} \Pbm \Big( \{ \Akz \cap \Rkpp \}^{\complement} \Big)
+ \varlimsup_{\sfKfv,\sfKz,\vep}
\Pbm \big( \bigcup_{u=\sfM+1}^{\sfMo}\bigcup_{\bb \in u} \bigcup_{\sfv \in \bF^\bb} \whUUvkp \big) 
 \le C \sfk^{-4} + C'_{\sfk} \srs^2 
\end{align}
(for the first inequality we also expressed $\CCpp$,  similarly to \eqref{cov-time-exc},
in terms of 
inter-city excursions from $\tau_0$ to $\tau_0+\CCpp$,  
utilizing the analog of \eqref{ct2.25cjj}).  Further,  
$\bP(\tau_0 \ge t_\sfk) \le \sfk^{-4}$ for some $t_\sfk$ finite,  while
$\osfm_\vep \ge \sfm_\vep$ and $\sfm_\vep+z_\sfk \to \infty$ as $\vep \to 0$.  It then follows 
from \eqref{eq:grid-bd},  \eqref{eq:cov-patch} and the union bound that for 
some finite $C'$,  $C'_{\sfk}$,  any $\sfk \ge 1$ and $\srs>0$ 
\[
\varlimsup_{\vep \to 0} 
\Pbm \big( \sqrt{\CCp} > \osfm_\vep -  \sfk^{\gak}/20  \,  \big)  \le  C'  \sfk^{-4} + C'_{\sfk} \srs^2  \,.
\]
Now taking $\srs \to 0$ we arrive at 
\[
\lim_{\srs \to 0} \varlimsup_{\vep \to 0} \Pbm \big( \sqrt{\CCp} > \osfm_\vep -  \sfk^{\gak}/20  \,  \big)  \le  C'  \sfk^{-4} \,,
\]
which amounts to \eqref{eq:cov-p-lim} since the \abbr{lhs} increases with $\sfk$,  
while the right hand side goes to zero as $\sfk \to \infty$.
\end{proof}

Our next result about $T_{\bb,\sfKo}^{[z]}$,  deals with $z=z(\sfKo) \ll 0$,  
which we utilize to control the left-tail of $\sqrt{\smash[b]{\CCd}}-\osfm_\vep$
when proving \Cref{prop:tight}.
\begin{lemma}\label{lem:clt-left-tight}
Taking  $z=z(\sfKo):=-2 \sfKo + \frac{1}{2} \sfKo^{\kapp}$, 
\begin{align}
\varlimsup_{\sfKo,\vep} \sup_\bb \bP ( 
\{ \LKoz \notin \bB^\infty_{\bar \mu} (2 \boetaKo) \} \cap \{ \wh{\eta}_{\bb,\sfKo}^{\,[z]}  \in \sfI_{\sfKo} \}) 
& = 0 \,,
\label{eq:VKo-good}
\\
\label{eq:YKo-good}
\varlimsup_{\sfKo,\vep} \sup_\beta \bP ( \wh{\eta}_{\bb,\sfKo}^{\,[z]}  \notin \sfI_{\sfKo} ) & = 0\,.
\end{align} 
\end{lemma}
\begin{proof} In view of the left-side of \eqref{eq:Jay-Thm-4.2},  \abbr{wlog} we intersect the events in 
\eqref{eq:VKo-good} and \eqref{eq:YKo-good} with $\Akz$ of \eqref{def:Jay-Thm4.1},  thereby having 
the properties of \eqref{def:Jay-Thm4.1-cons}.  With these in place,  as in the proof of \Cref{lem:adjust-AUz},
the event in \eqref{eq:VKo-good} implies that at $j=\sfKo$,  either \eqref{eq:Tb-Ta-W} holds,
or $\Ljpp \notin \bB^{\infty}_{\bar \mu} \big(\boetaj \big)$ and $\whetajpp \in [\frac{1}{2} j^{\kapo},j^{\kapp}]$.
Similarly,  the event of \eqref{eq:YKo-good} then implies 
that either $\wh{\eta}_{\bb,\sfKo}^{\,\langle z' \rangle} \le \sfKo^{\kapo}$
or  $\wh{\eta}_{\bb,\sfKo}^{\,\langle z'' \rangle} \ge \sfKo^{\kapp}$.

For $z=z(\sfKo)$ we still have \eqref{eq:b-a-bd} holding at any $\sfKo \ge \sfKo'(\sfk)$,  and thus 
deduce from \eqref{A.2-n},  similarly to \eqref{eq:step1-123},
that the probability of \eqref{eq:Tb-Ta-W} holding at $j=\sfKo$ remains $o_\vep(1)$.
Further,  conditional on $\whetajpp=\sqrt{2 \sfs}-\gamma_0(j)$,  we have that $\Ljpp$ is 
the empirical measure of $\sfs$ \abbr{iid} samples from $\bar \mu$. 
Hence,  by \eqref{lbl:d-inf-bd},  we have for $\cst \ge c_0$ 
and $j=\sfKo$,  that uniformly over $\whetajpp \in [\frac{1}{2} j^{\kapo},j^{\kapp}]$
\[
\bP \big( \Ljpp \notin \bB^{\infty}_{\bar \mu} \big(\boetaj \big) \, \big| \, 
 \whetajpp \big) \le 2 \sum_{m \ge \log j} e^{-m^2/4} 
\]
(see \eqref{eq:Lkp-WNa} and \eqref{eq:BRZ-454}).  With the latter bound being $o_j(1)$,  this 
concludes the proof of \eqref{eq:VKo-good}.

Next,  note that $\wh{\eta}_{\bb,\sfKo}^{\,\langle z' \rangle} \le \sfKo^{\kapo}$ 
implies having at most $\sfs'=\sfs_{\sfKo}(\sfKo^{\kapo})$ excursions at level $\sfKo$
within the first $\sfs=\sfs_1(z')$ excursions at level one.  Thus,  by the \abbr{clt} for 
$\big(\sqrt{2 T^{\sfv_\bb,\sfs}_{\sfKo}}-\sqrt{2\sfs}\, \big)/\sqrt{\sfKo}$ (as $\sfs \to \infty$),  see  
\eqref{eq-calgary} and \eqref{eq:Jay1.12},  or \cite[Proposition  1.4]{BRZ},
\begin{equation}\label{eq:YKo-clt-lb}
\bP\big( \wh{\eta}_{\bb,\sfKo}^{\,\langle z' \rangle} \le \sfKo^{\kapo} \big) \le C \int_{-\infty}^{w_1} e^{-u^2/2} du \,,
\end{equation}
where by our choice of $z(\sfKo)$
we have that for some $\alpha>0$ and any $\sfKo \ge \sfKo'(\sfk)$,  
\[
w_1:=\frac{\sqrt{2 \sfs'}-\sqrt{2 \sfs}}{\sqrt{\sfKo}}
=\frac{- \sfrho_\vep \sfKo + \sfKo^{\kapo}-z(\sfKo)+\log \sfk}{\sqrt{\sfKo}} \le - \sfKo^{\alpha} \,.
\]
The right hand side of \eqref{eq:YKo-clt-lb} is thus $o_{\sfKo}(1)$.  
Similarly,  taking now $\sfs=\sfs_1(z'')$ and $\sfs'=\sfs_{\sfKo}(\sfKo^{\kapp})$,
\begin{equation}\label{eq:YKo-clt-ub}
\bP\big( \wh{\eta}_{\bb,\sfKo}^{\,\langle z'' \rangle} \ge \sfKo^{\kapp} \big) \le C \int_{w_2}^{\infty} e^{-u^2/2} du \,,
\end{equation}
where by our choice of $z(\sfKo)$,  for any $\sfKo \ge \sfKo'(\sfk)$,  
\[
w_2=\frac{\sqrt{2 \sfs'}-\sqrt{2 \sfs}}{\sqrt{\sfKo}}
=\frac{- \sfrho_\vep \sfKo + \sfKo^{\kapp}-z(\sfKo)-\log \sfk}{\sqrt{\sfKo}} \ge \sfKo^{\alpha} \,.
\]
Consequently,  the right hand side of \eqref{eq:YKo-clt-ub} is also $o_{\sfKo}(1)$,  thereby concluding the proof 
of \eqref{eq:YKo-good}.
\end{proof}

\section{Essentially radial: Proof of \Cref{prop:angle-diff}}\label{sec:angle-diff} 

Throughout this section and \Cref{sec:good-cover}
we use the notation $n'=\nvf$, $n''=\nvfv$,  see \eqref{def.n-vep-i},  setting $\sfIs_i$ and $\ggi (\cdot)$ as in
\Cref{dfn:ggi-barsi}.  That is,  we restrict $\sfI_{\sfK_i}$, 
$i \le 3$ to $y_i$ with $\sfs_{\sfK_i}(y_i) \in \bN$,  while $\sfI_{n'}$ and $\sfI_{n''}$ 
are restricted to $y_i$ such that $\bar \sfs_{\sfK_i}(y_i) \in \bN$,  c.f.   \eqref{eq:nominal} 
and \eqref{eq:nominalno-vep}.  We also use within $[n',n_\vep]$
the $\vep$-free barrier and window heights of \Cref{rmk:no-vep-inJs} and 
as mentioned in \Cref{rmk:free-kappas},  slightly reduce $\kapo$ and increase $\kapp$
for $j=\sfK_3$ to make $\sfI_{\sfK_3}$ much wider than $\sfI_{\sfK_2}$.
Fixing a generic block $\beta$ of center $\sfv_\beta=\sfv$, 
we often use the $\sigma$-algebras $\GG^\sfv_j$ of \Cref{def:Fvjs} and the following simplified laws.
\begin{definition}\label{def:Fvj}
Denote by $\bar \bP^{(\ell)}_{y}$ the product law of $\sfs_{\ell}(y)$ \abbr{iid} 
Brownian excursions from $\partial \bB_{\sfd}(\sfv,r_\ell)$ to $\partial \bB_{\sfd}(\sfv,r_{\ell-1})$,  
starting independently,  uniformly on $\partial \bB_\sfd (\sfv,r_{\ell})$.  
The law it induces on the counts of  Brownian excursions at levels $j > \ell$ 
around $\sfv$ matches that of the occupation measure of a
reflected \abbr{srw} on $[\ell-1,\infty)$,  during its first $\sfs_{\ell}(y)$ moves from $\ell$ to $\ell-1$. 
\end{definition}


The proof of \Cref{prop:angle-diff}
consists of the following four steps:\\
{\bf I.} Having no constraints between levels $n'$ and $n''$ allows us,  conditional on the number of hits of level $n''$,
 to modify the law of the excursions inside the block
to be $\bar \bP^{(n')}$,  thereby contributing only factors independent of $\vep$,  $\underline{\sfx}$ 
and $\bb$; these ultimately lead to $\alpha(\sfKfv,\sfKp)$.\\
{\bf II.} Similarly,  thanks to the gap in constraints between levels $\sfKt$ and $\sfKth$,  
the probability of satisfying the constraints at levels $[\sfKth,n'']$ depends only 
on the number of $\sfKo$-level excursions that reach level $\sfKth$.  This  probability can be evaluated,  up to $o(1)$ relative error,
using $\sfKth$-level excursions under $\bar \bP^{(\sfKt)}$.\\
{\bf III.} \Cref{lem:uncond-comp} shows that removing the upper limit in $\sfI_j$,  $j \in (\sfKth,n')$ 
merely adds $o_{\sfKfv}(1)$ to $\alpha(\sfKfv,\sfKp)$,  and that once this is done,  adapting 
\cite[Sec.  6]{DRZ},  one recovers the right hand side of \eqref{eq:good-cover}. \\
{\bf IV.} \Cref{lem:alpha-unif-bd} shows that $\alpha(\sfKfv,\sfKp)$ is bounded above and away from zero 
by $\sfKp$-dependent positive constants,
uniformly in $\sfKfv$.

We now turn to the individual steps.

\noindent
{\bf I}.   Let $\CCe_\sfv(v,v')$ be the event of having precisely $\bar \sfs_{\sfKfv}(v')$ excursions 
from $\partial \bB_\sfd (\sfv, r_{n''})$ to $\partial \bB_\sfd(\sfv,r_{n'-1})$ during the first 
$\bar \sfs_{\sfKf}(v)$ excursions from $\partial \bB_\sfd(\sfv,r_{n'})$ to $\partial \bB_\sfd(\sfv,r_{n'-1})$.
We partition the 
events appearing in  $\UUpy(y)$,  see \eqref{eq:good-cover-y},  according to the disjoint events $\{\whetaynp=v,  \CCe_{\sfv_\bb}(v,v')\}$
(see \eqref{eq:norm-exc}-\eqref{eq:eta-centered} and note that then $\whetaynpp \ge v'$,  because each excursion from $\partial \bB_\sfd (\sfv, r_{n''})$ to $\partial \bB_\sfd(\sfv,r_{n'-1})$ generates at least one excursion from $\partial \bB_\sfd (\sfv, r_{n''})$ to $\partial \bB_\sfd (\sfv, r_{n''-1})$).  Suppressing the implicit dependence on $\beta$ and $y=y_1 \in \sfIs_1$,  we further represent each element of the partition 
as follows.  First,  let  $\BWb(v)$ denote the restriction to 
$\{\whetaynp=v\}$ of the event $\BBy \cap \WWy$ excluding the constraint $\{\whetaynpp \le \sfKfv\}$.  Next, let   $\QQe(v')$ denote the event of  having by
the end of the $\bar \sfs_{\sfKfv}(v')$-th excursion from $\partial \bB_\sfd (\sfv, r_{n''})$ to $\partial \bB_\sfd(\sfv,r_{n'-1})$, no more than $\bar \sfs_{\sfKfv} (\sfKfv)$  $n''$-level excursions 
around $\sfv$ and in addition there is some $\sfu \in \bF^\star$ with 
at least $\gamma_{-1} (j)^2/2$ excursions from 
$\partial \bB_\sfd (\sfv+\vep \sfu, r_{j})$ to $\partial \bB_\sfd(\sfv+ \vep \sfu,r_{j-1})$ for every 
$j \in \sfJs$ and some unhit $\vep$-ball centered within the square of edge-length 
$\sqrt{2} \wt{c} \, \vep$ around $\sfv+\vep \sfu$.  (With this definition, $\QQe (v')$ is the intersection of $\{\whetaynpp \le \sfKfv \}$ with $\RUy$ of \eqref{eq:barrier-event},
restricted to the event $\CCe_{\sfv_\bb}(v,v')$ occurring. Note that on this restriction, $y$ plays no role and is effectively replaced by $v'$.) Then, we have that
\begin{equation}\label{eq:decom-U}
 \UUpy = \bigcup_{y_i \in \sfIs_i} \Big\{ \BWb (y_4) \cap \CCe(y_4, y_5)
 \cap \QQe (y_5)   \Big\},
 \end{equation}
\blue{where} thanks to our choice of $\sfJs$, 
the event $\QQe(y_5)$ is measurable on the $\sigma$-algebra of those $\bar \sfs_{\sfKfv}(y_5)$ excursions 
at level $\sfKf$ around $\sfv$ which hit $\partial \bB_{\sfd}(\sfv,r_{n''})$. 

With our adjustment for $\gamma_b(j)$ in \Cref{rmk:no-vep-inJs},  and having 
$v$ and $v'$ fit \eqref{eq:nominalno-vep}, 
by Brownian scaling 
\begin{equation}\label{def:alpha56}
\alpha_{\sfKfv,\sfKp} (v') :=  \bar \bP^{(n')}_v \big( \QQe (v') \big| \, \CCe(v,v')  \big) \,,
\end{equation}
is independent of $\vep$,  of $v$,  and of the block center $\sfv=\sfv_\beta$ (see \Cref{def:Fvj} for $\bar \bP^{(\ell)}_y$). 
From Lemma \ref{recursionends2} we have 
 that for $\sfKf \ge \sfKfv + 5\log \sfKfv$, uniformly in $y_i \in \sfIs_i$, 
\begin{align}
 \bP \Big(\QQe (y_5) \big| \, \CCe (y_4,y_5), \GG^{\sfv}_{n'} \Big)
 &= (1+o_{\sfKfv} (1)) \alpha_{\sfKfv,\sfKp} (y_5) \,.
 \label{eq:K3-angles} 
 \end{align}
(See the text in the beginning of \Cref{probends} to match with our setting.)

The ordered collection $\underline{\sfx}$ of start point and endpoint of the first $\sfs$ of the $\sfKo$-level excursions 
of block $\beta$ and the event $\BWb(y_4)$ are both measurable on $\GG^\sfv_{n'}$.  It 
thus follows from \eqref{eq:decom-U} and \eqref{eq:K3-angles} that
\begin{align}\label{eq:K4-breakup}
\bP(\UUpy &| \underline{\sfx}) =
 \sum_{y_i \in  \sfIs_i}
\bE\Big[ {\bf 1}_{\BWb (y_4)} \bP(\CCe (y_4,y_5) | \GG^{\sfv}_{n'})
\bP\big(\QQe (y_5) \big| \, \CCe (y_4,y_5), \GG^{\sfv}_{n'} \big)
 \Big| \underline{\sfx} \Big]
\nn \\
&=
(1+o_{\sfKfv}(1))
\sum_{y_i \in \sfIs_i} \alpha_{\sfKfv,\sfKp}(y_5) 
\bP\Big( \BWb (y_4) \cap \CCe (y_4,y_5) \,\big|\, \underline{\sfx} \Big).
\end{align}
Recall that the event  $\BWb (y_4) \cap \CCe(y_4,y_5)$ is determined by 
the block excursion counts $\{ \Tjy \}$ (and with $r_{n'-1} \to r_{n''}$ excursions  
when $j=n''$),
and the associated \emph{un-ordered}
collections of angle-differences $\{\thbj, q \in[1,\Tjy]\}$,
for $j \in [\sfKo,n'']$
(see \eqref{eq:barrier-event-z} and \eqref{eq:was-event}).
These variables, which are measurable 
with respect to the first $\sfs$ of the $\sfKo$-level excursions of our block $\beta$, are invariant 
with respect to both the order of these $\sfKo$-level block excursions
and to rotations around $\sfv$ of the Brownian path within each of these $\sfKo$-level excursions. 
Consequently,  out of $\underline{\sfx}$,  
all that matters for $\BWb (y_4) \cap \CCe(y_4,y_5)$ is the 
un-ordered collection of angle-differences $\{\thbKo, \,  q \in [1,\sfs]\}$,  
captured by their empirical measure $\LKoy$.  That is,
\begin{equation}\label{eq:uux-to-emp-meas}
\bP\Big( \BWb (y_4) \cap \CCe(y_4,y_5) |\underline{\sfx} \Big) =
\bP\Big( \BWb (y_4) \cap \CCe(y_4,y_5) | \, \LKoy \,  \Big)\,,
\end{equation}
which in combination with \eqref{eq:K4-breakup} results 
with 
\begin{align}\label{eq:almost-radial}
\bP(\UUpy(y) | \underline{\sfx}) &= (1+o_{\sfKfv}(1))
\barP ( \UUpy(y) | \LKoy ) \,,
\end{align}
where we set hereafter for any possible value $L$ of $\LKoy$,
\begin{equation}\label{eq:pre-precise}
\barP ( \UUpy(y) |  \, L ) :=
\sum_{y_i \in \sfIs_i} \alpha_{\sfKfv,\sfKp}(y_5) 
\bP\Big( \BWb (y_4) \cap \CCe(y_4,y_5) \, \big| \,  L \Big).
\end{equation} 
Further,  
the non-random factors $\alpha_{\sfKfv,\sfKp}(\cdot)$ on the right hand side are independent 
of $\vep$, $\beta$, $y$ and 
the conditional probability on the right hand side of \eqref{eq:pre-precise}
depends on $\beta$ and $y$ only through $\sfs=\sfs_{\sfKo}(y)$ and the value $L$ of 
$\LKoy$.  We thus eliminate hereafter the value of $\beta$, 
setting the block center as $\sfv_\beta = \npp$ and using 
$\eta(j)$,   $\wh{\eta}(j)$ and $L_{\sfKo}(\sfs)$,  respectively,  for 
$\etajy$,  $\whetajy$ and  $\LKoy$.

\noindent
{\bf II.}  We define in terms of  $\WW^{(\vep)}$ of \eqref{eq:was-event}
and $\sfI_j$ of \eqref{dfn:windows},  the events
\begin{align}\label{eq:B-K1-K2}
\BW^{(\vep)}_{[\sfKo,\sfKt]}(y_2|y_1)  & :=  \WW^{(\vep)}(y_1) \bigcap \big \{\widehat{\eta} (\sfKt) = y_2 \big\} 
 \bigcap_{j \in [\sfKo,\sfKt]} \big\{ \wh{\eta} (j) \in \sfI_j  \big\} \,,    \quad\\
\BBKth (y_4,y_5) &:= \big\{ \wh{\eta} (n') = y_4, \; \wh{\eta}^{\sharp} (n'')=y_5 \big \} 
\!\!\!\!
\bigcap_{j \in [\sfKth,n']} \!\!\!\!\! \{ \wh{\eta} (j) \in \sfI_j \} \,, 
\label{eq:B-K3} \\
\wCe (y_2,y_3) & :=\{ \widehat{\eta} (\sfKth)=y_3 \,| \,  \widehat{\eta} (\sfKt)=y_2 \} \,,
\label{eq:hat-C}
\end{align}
where the latter event
amounts to precisely 
$\sfs_{\sfKth}(y_3)$ excursions from 
$\partial \bB_\sfd ({\bf 0}, r_{\sfKth-1})$ to $\partial \bB_\sfd({\bf 0},r_{\sfKth})$
during the first 
$\sfs_{\sfKt}(y_2)$ excursions from $\partial \bB_\sfd({\bf 0},r_{\sfKt})$ to 
$\partial \bB_\sfd({\bf 0},r_{\sfKt-1})$ (and $\wh{\eta}^{\sharp} (n'')$ denotes
the analog of $\wh{\eta}(n'')$ in terms of 
excursions from $\partial \bB_\sfd({\bf 0},r_{n'-1})$
to $\partial \bB_\sfd({\bf 0},r_{n''})$).  Now,  recall from \eqref{dfn:J-intervals} and the \abbr{lhs} of 
\eqref{eq:barrier-event-z} that 
\begin{align}\label{eq:decomp}
\BWb (y_4) \cap \CCe (y_4,y_5) = \bigcup_{y_2,y_3} \Big\{ \BB\WW^{(\vep)}_{[\sfKo,\sfKt]}(y_2|y_1) \cap 
\wCe(y_2,y_3) \cap \BBKth (y_4,y_5) \Big\}
\end{align}
with a disjoint union over $y_2 \in \sfIs_2$ and $y_3 \in \sfIs_3$.  
Recall $\GG_{\sfKt}:=\GG^{\bf 0}_{\sfKt}$ of \Cref{def:Fvjs}
(containing the end points of all excursions
from $\partial \bB_{\sfd}({\bf 0},r_{\sfKt})$ to  $\partial\bB_{\sfd}({\bf 0},r_{\sfKt-1})$).
With our choice of $\sfKth-\sfKt \ge 2 \log (5 n_\vep^2)$,
we deduce from Lemma \ref{recursionends},
analogously to \eqref{eq:K3-angles},  that uniformly in 
$y_i \in \sfIs_i$, $i=2,3$ and a.s.
\begin{align}\label{eq:g3}
 e^{2(n''-\sfKth)} \sum_{y_4,y_5}  \alpha_{\sfKfv,\sfKp}(y_5)
\bP \Big(\BBKth(y_4,y_5) \, \big| \,  \wCe (y_2,y_3), \GG_{\sfKt} \Big)
= (1+o_{\vep} (1)) \gKth(y_3|y_2) \,,
\end{align}
where for $\bar \bP^{(\ell)}_y$ of \Cref{def:Fvj},  
\begin{equation*}
\gKth(y_3|y_2) := e^{2(n''-\sfKth)}
\sum_{y_4,y_5} \alpha_{\sfKfv,\sfKp}(y_5) 
\bar{\bP}^{(\sfKt)}_{y_2} \Big( \BBKth (y_4,y_5)
  \big| \, \wCe (y_2,y_3) \Big) \,,
\end{equation*}
which by definition is zero for $y_3 \notin \sfIs_3$.
The event $\BBKth (y_4,y_5)$ concerns 
only excursion counts at levels $j \ge \sfKth$ within a tower of concentric circles,  
which under $\bar \bP^{(\sfKt)}_{y_2}(\cdot| \,  \wCe (y_2,y_3) \, )$ occur during 
prescribed $\sfs_{\sfKth}(y_3)$ excursions at the outer-most $\sfKth$-level,  whose 
angle-differences are \abbr{IID} samples from $\bar \mu$ of \eqref{dfn:Poiss}.
As such,  we may replace the aggregate 
of our $\sfKth$-level excursions by \abbr{iid} excursions,  thereafter
representing $\widehat{\eta}(\cdot)$ as the excess occupation measures for the law
$\bar \bP^{(\sfKth)}_{y_3}$ of a (reflected) \abbr{srw} on the interval $[\sfKth-1,n'']$, 
during its first $\sfs_{\sfKth}(y_3)$ moves from $\sfKth$ to $\sfKth-1$.  In particular, 
$\gKth (y_3|y_2)=\gKth(y_3)$ is independent of $y_2$,  with 
\begin{align}\label{eq:alt-g3}
\gKth(y_3) & := e^{2(n'-\sfKth)} 
\bar{\bE}^{(\sfKth)}_{y_3} \Big[ \gfo\big(\wh{\eta}(n')\big); \Bs \Big] \,,  \qquad 
\Bs := \!\!\! \!\! \bigcap_{j \in [\sfKth,n']} \{ \wh{\eta} (j) \in \sfI_j \}\, \\
\label{def:g4}
\gfo (y_4) &:= e^{2(\sfKf-\sfKfv)} {\bf 1}_{\sfIs_4}(y_4) \sum_{y_5 \in \sfIs_5} \bar{\bP}^{(n')}_{y_4} 
(\widehat{\eta}^\sharp (n'')=y_5)  \alpha_{\sfKfv,\sfKp}(y_5) \,.
\end{align}
Plugging \eqref{def:alpha56} into \eqref{def:g4},  unraveling the definitions of 
$\QQe (y_5)$ and $\CCe(y_4,y_5)$ and scaling space by $e^{-2}/r_{n'}$
we find that $\gfo(\cdot)$ of \eqref{def:g4} matches our definition in \eqref{dfn:g4-alt} 
and as such,  it is independent of $\vep$.

With $L=L_{\sfKo}(\sfs)$ and $\BW^{(\vep)}_{[\sfKo,\sfKt]}(y_2|y_1)$ measurable on $\GG_{\sfKt}$,  
we get from \eqref{eq:pre-precise}, \eqref{eq:decomp} and \eqref{eq:g3}, 
analogously to \eqref{eq:K4-breakup},  that in terms of $\gKth(\cdot)$ of \eqref{eq:alt-g3}, 
\begin{equation}\label{eq:step2}
e^{\sfco(n''-\sfKth)} \barP( \UUpy(y) |  \, L ) = (1+o_\vep(1)) 
\, \bE \big[ \gKth\big(\wh{\eta}(\sfKth) \big);  \BWl   \,| L \big]
\,,
\end{equation} 
where per fixed $y \in \sfIs_1$,  the event
\begin{equation}\label{eq:BW-K12-alt}
\BWl  := 
\bigcup_{y_2,y_3} \Big\{
\BW^{(\vep)}_{[\sfKo,\sfKt]}(y_2|y) \cap \wCe (y_2,y_3) \Big\}
\end{equation}
matches $\BW^{(\vep),(y)}_{\npp,[\sfKo,\sfKt]}$ of \eqref{eq:BW-K12}.

\noindent
{\bf III.} Recalling $\gfo(\cdot)$ of \eqref{dfn:g4-alt},  we set for $\ell \in [\sfKth,n']$, 
\begin{equation}\label{def:wh-g3-ell}
\gfl(y') := e^{2(n'-\ell)} \bar{\bE}^{(\ell)}_{y'} \Big[ \gfo\big(\wh{\eta}(n')\big); 
 \bigcap_{j \in (\ell,n')}  \{ \wh{\eta} (j) >  j_{n_\vep}^{\gak} \} \Big] \,,
\end{equation}
when $y' \in \sfI_\ell$ and $\sfs_\ell(y') \in \bN$,  while otherwise $\gfl(y')=0$.  Our next lemma shows that
\begin{equation}\label{def:wh-g3}
\whgth(y_3) = e^{2(n'-\sfKth)} \,  {\bf 1}_{\sfIs_3} (y_3)\, 
\bar{\bE}^{(\sfKth)}_{y_3} \Big[ \gfo \big(\wh{\eta}(n') \big); 
 \bigcap_{j \in (\sfKth,n')}  \{ \wh{\eta} (j) >  j_{n_\vep}^{\gak} \} \Big] \,,
\end{equation}
is a good proxy for $\gKth(\cdot)$ of \eqref{eq:alt-g3},  so that when combined with \eqref{eq:almost-radial} 
and \eqref{eq:step2},  it results with \eqref{eq:good-cover}.
\begin{lemma}\label{lem:uncond-comp}  
Recall $\gth(y)={\bf 1}_{\sfIs_3} (y) y e^{-\sfco y}$ and set
 \begin{equation}\label{eq:alpha-tot}
 \alpha(\sfKfv,\sfKp) := \frac{\sqrt{2} \, e^{-1}}{\sqrt{\pi \sfKf}} \int_0^\infty  
 \mathbb{E}\Big[\gfo\big(\sqrt{2Y_{\bar \sfs_{\sfKf}(y)}}- \sfco \sfKf\big)\Big] y e^{\sfco y} dy \,,
\end{equation} 
with $Y_\lambda:=$Poisson$(\lambda)$ and $\bar \sfs_{\sfKf}(\cdot)$ of \eqref{eq:nominalno-vep}.
Then,  uniformly over $y_3 \in \sfIs_3$,
\begin{align}\label{eq:g3-formula}
\whgth(y_3) 
& = (1+o_{\sfKfv,\vep} (1)) \alpha(\sfKfv,\sfKp) \gth(y_3) \,,\\
\label{eq:g3-to-wh-g3}
\gKth(y_3) &=  \whgth(y_3)+ o_{\sfKfv}(1) \gth(y_3) \,.
\end{align}
\end{lemma}
\begin{proof} Sharp asymptotics of rare event probabilities as in \eqref{eq:g3-formula}
are considered in \cite[(6.4)]{DRZ}.  To make the connection clearer we change below 
some of the notation in \cite{DRZ} to bring it in line to what we have here.  In particular, 
we write $m,n,\rho_n$ and $\ell$ of \cite{DRZ} as $\sfm_\vep$,  $n_\vep$,  $\sfrho_\vep$,  
and $\sfKf$,  respectively.  With these changes,  the asymptotics for \cite[(6.4)]{DRZ} 
would apply to \eqref{eq:g3-formula},  except for:\\
(a).  Replacing in \cite{DRZ} the correction term ${\log |\log \vep|}/{|\log \vep|}$ in the formula \eqref{def:mn} 
for the slope $\sfrho_\vep$ by $\log n_\vep/n_\vep$.\\
(b).  Replacing $\gfo(\cdot)$ of \eqref{dfn:g4-alt} by a function $\widehat{\gamma}_{\sfKf}(\cdot)$ 
supported on smaller interval $\sqrt{\sfKf}[(\sqrt{\log \sfKf})^{-1},\sqrt{\log \sfKf}]$.\\
(c).  Having in \cite{DRZ} a straight barrier (ie.  $\wh{\eta}(j)>0$),  while \eqref{def:wh-g3} has a curved 
one (that is,  $\wh{\eta}(j) > j_{n_\vep}^{\gak}$).\\
(d).  Changing in \cite{DRZ} our nominal slope value $\sfco$ to $c_\star = \sqrt{2 \log 2}$,
having there also $\sfKth=0$ and fixed $y_3$ (whose value does not grow as $n'-\sfKth \to \infty$).\\
Nevertheless,  it is not difficult to adapt the proof in \cite{DRZ} to the case in hand.  Specifically:\\
(a).  Since $\log n_\vep - n_\vep {\log |\log \vep|}/{|\log \vep|} \to 0$ 
as $\vep \to 0$ (see \eqref{def.nvep}),  the modification of the correction term in $\sfm_\vep$ is irrelevant. \\
(b).  The arguments in \cite[Sec. 6]{DRZ} apply verbatim for any positive nominal slope and
any function on the larger interval $\sfIs_4$ (the key is 
$\sup_{y \in \sfIs_4} \{ y/\sfKf + y^{-1} \log \sfKf \} = o_{\sfKf}(1)$, see \cite[(6.39)]{DRZ} and the text following it).\\
(c).  The straight barrier event probability exceeds the one for our curved barrier.
Moreover,  \cite[(6.8) and (6.18)]{DRZ} \emph{lower bound} the former probability,  up to relative 
error $o_{\sfKf}(1)$,  in terms of Brownian motion being above a barrier of the form 
$\gamma_0(j)+ 2 (h_{\sfKf} +  j_{n'}^{\delta})$,  \emph{for any} $h_\ell \ge (\log \ell)/2$.  
One then proceeds to show that the latter curving has a negligible effect when $\delta < {1}/{6}$ 
and $\sup_{y \in \sfIs_4} \{ \sfKf^{-1} y + y^{-1} h_{\sfKf}  \} = o_{\sfKf} (1)$ (see \cite[Lemma 6.3 and 
last display on Pg.  854]{DRZ}).  For some $c_\gak$ finite,  
$2 c_{\gak} \sfKf^{\gak} + 2 j_{n'}^{\gak} \ge j_{n_\vep}^\gak$.  Hence, 
upon taking $\delta=\gak$ and $h_\ell = c_{\gak} \ell^{\gak}$,  the above mentioned barrier of \cite{DRZ}
exceeds our barrier $j_{n_\vep}^\gak$,  while in view of \eqref{dfn:windows}, 
having $\kapo>\gak$ guarantees that the arguments of \cite{DRZ} apply for our
larger support of $\gfo(\cdot)$ and larger $h_{\sfKf}$ value.\\
(d).  The relative $o(1)$ error in \cite[Section 6]{DRZ},  is uniform over our range of
$h_{\sfKf}/{y_3} + {y_3}/{\sqrt{n'-\sfKth}} = o_\vep(1)$ (see \cite[(6.37) and derivation of (6.27)]{DRZ}).
Using the two displays above \cite[(6.26)]{DRZ} we see that the nominal slope of $c_\star$ results in
a leading exponential tail $e^{-(c_\star^2/2) (n'+1/2)}=2^{-(n'+1/2)}$ for the expectation on the 
right hand side of \cite[(6.4)]{DRZ}, 
whereas our slope $\sfco$ merely modifies it to be $e^{-(\sfco^2/2) (n'+1/2)} = e^{-2 n'-1}
$,  if we had $\sfKth=0$ as in \cite{DRZ} (this computation is also the reason for 
the additional factor $\sqrt{2}/e$ we have in \eqref{eq:alpha-tot}).
The only effect of having here $\sfKth \ne 0$ such that $\sfKth/n' \to 0$, 
is to modify the value of $n'$ to $n'-\sfKth$ in this leading factor.  We thus 
arrive at \eqref{eq:g3-formula} holding uniformly over $y_3 \in \sfIs_3$
(see \cite[(6.2)]{DRZ}, and recall that $\sfKth \to \infty$ as $\vep \to 0$).  

Turning to \eqref{eq:g3-to-wh-g3},  comparing \eqref{eq:alt-g3} with \eqref{def:wh-g3},
the bound \eqref{ubd2-g4} on $\gfo(\cdot)$ implies that for any $z \in \sfIs_3$,  
 \begin{align}\label{eq:ubd-whg3}
0 &\le \whgth(z) -  \gKth(z) 
\le  \frac{C(\sfKp)}{\sfKf}  e^{2(n'-\sfKth)}
\sum_{j=\sfKth+1}^{n'-1} \sum_{u =  j_{n_\vep}^{\okappa}}^\infty
\sum_{y} 
p_j(z,u,y) y  e^{-\sfco y}  e^{-\frac{y^2}{8 \sfKf}} 
\end{align}
where the inner-most sum is over all integer values of $y \in \sfI_{n'}$ and 
\begin{align*}
p_j(z,u,y):= 
\bar{\bP}^{(\sfKth)}_{z} \big( \wh{\eta} (j) \in I_u,  \wh{\eta}(n') \in I_y,
\wh{\eta}(\ell)>0 \;\; \text{for all} \;\;  \sfKth< \ell < n' 
\big) \,.
\end{align*}
Recall that $1 \ll \sfKth \le C (\log n_\vep)^{c}$ and as $\okappa > \kapp$,  here $u \gg z$ 
uniformly over $j \in (\sfKth, \frac{1}{2} n_\vep]$ and $z \in \sfI_{\sfKth}$,  while 
$u \gg y$ uniformly over $j \in [\frac{1}{2} n_\vep,n')$ and $y \in \sfI_{n'}$.
We thus get from the \blue{
\Cref{prop:BarrierSecGWPropUB}, 
similarly} to 
our derivation of \eqref{agar.45p},   
that for $f_{j,s}(u)$ of \eqref{def:fjs},
some constant $c$,  any such $z,y,u$ and $j \in (\sfKth, \frac{1}{2} n_\vep]$, 
\begin{equation}\label{eq:ubd-pjzuy}
e^{2(n'-\sfKth)} p_j(z,u,y) \blue{= e^{2(n'-\sfKth)} \KQ_{[\sfKth,n']}^{(\sfKth,j,n')}(z,u,y)}
\le 
c \,  \gth(z) \sfKf^{-1/2} y e^{\sfco y} f_{j_{n_\vep},\sfKth}(u)
\end{equation}
while for $j \in [\frac{1}{2} n_\vep,n')$ one merely replaces 
$\sfKth$ on the right hand side by $\sfKf$ 
(as done in \eqref{agar.45jj}).  Further,  with $\okappa \ge \gak$ and $a:=2\okappa -1 >0$, 
 if $j_{n_\vep}=s+\ell$ for $s,\ell \ge 1$,   then for some universal $c$, 
\[
\sum_{u = (s+\ell)^{\okappa}}^\infty f_{s+\ell,s}(u) \le 4 \sum_{u=(s+\ell)^{\okappa}}^\infty u^2 \, \ell^{-3/2} e^{-\frac{u^2}{16\ell }}
\le c \exp\big( -\frac{1}{18} (s+\ell)^{a} \big) 
\]
(see \eqref{def:fjs}),  and with $\sfKth \ge \sfKf$,   upon plugging the bounds 
of \eqref{eq:ubd-pjzuy} into \eqref{eq:ubd-whg3} we arrive at 
\[
\big|\,\whgth(z) -  \gKth(z)\big| \le 
C'(\sfKp) \gth(z) \Big[ \sum_{y} \sfKf^{-3/2} y^2 e^{-\frac{y^2}{8 \sfKf}} \Big]
\sum_{\ell \ge 1} \exp \big( -\frac{1}{18} (\sfKf+\ell)^{a} \big) \le C'(\sfKp) \,  \gth(z) e^{- c \,  \sfKf^{a}} 
\]
That is,  \eqref{eq:g3-to-wh-g3} holds uniformly over $y_3 \in \sfIs_3$. 
\end{proof}

\begin{remark}\label{rmk:g3-formula-gen}
Setting $\sfI_\ell$ of \eqref{dfn:windows} to have parameter $\gakp>\gak$,   our proof of  \eqref{eq:g3-formula} extends to show that for 
any $\ell=o(n')$,  uniformly over $y' \in \sfI_\ell$ with $\sfs_\ell(y') \in \bN$,  
\begin{align}\label{eq:g3-formula-ext}
\gfl(y') 
& = (1+o_{\sfKfv,\vep} (1)) \alpha(\sfKfv,\sfKp) \,  y' e^{-\sfco y'} 
\end{align}
and by a similar reasoning,  for some $c_1=c_1(\sfKfv)$ finite and any $\vep>0$,  
\begin{align}\label{eq:gfl-comp}
\sup \{\gfl(\wh{\eta}(\ell)),  \wh{\eta}(\ell) \in [b,b+2) \} \le c_1 \inf \{\gfl(\wh{\eta}(\ell)),  \wh{\eta}(\ell) \in [b,b+2) \} 
\quad \forall \ell \in [\sfKth, n'],  \;\; b \in \sfI_\ell \,.
\end{align}
\end{remark}

Before proceeding to step {\bf IV},  our next corollary will serve us when proving 
\Cref{lem:uncond-top}.
\begin{corollary}\label{cor:gK2-gK1-arg} 
Recall $\ggi(y)={\bf 1}_{\sfIs_i} (y) y e^{-2 y}$,  $i=1,2$.  Then,  uniformly in $y_i \in \sfIs_i$,
\begin{align}
e^{2 (\sfKt-\sfKo)} \,  
\bar \bE^{(\sfKo)}_{y_1} \Big[ \gtw \big( \wh{\eta}(\sfKt)\big);  \wh{\eta}(j) > j^{\gak},   j \in (\sfKo,\sfKt) \Big]
&=  (1+o_{\sfKz,\vep} (1)) \gon (y_1) \, .
\label{eq:K1-to-K2-u2}
\end{align}
\end{corollary} 
\begin{remark} One can adapt \cite[Section 6]{DRZ} to fit the 
simpler setting of \eqref{eq:K1-to-K2-u2} where the barrier length $\sfKt-\sfKo$ is much shorter 
than its ending nominal height $\gamma_0(\sfKt)$,  which is thus 
comparable to the starting barrier nominal height (and hence to $\sfm_\vep$).  But,  striving for 
a shorter proof,  we instead re-use our derivation of \eqref{eq:g3-formula} in part (a) below, 
at the cost of having to remove the obstacles of a non-integer 
$\log (r_{\sfKt}/r_{\sfKth})$ and $g_2(\cdot) \equiv 0$ outside $\sfI_{\sfKt}$
(see parts (b) and (c) below).
\end{remark}
\begin{proof}
The proof of \eqref{eq:K1-to-K2-u2} consists of the following three steps:\\
(a).  Our proof of \eqref{eq:g3-formula} requires a fixed ratio of radii between consecutive levels.
Specifically,  we could leverage it to verify \eqref{eq:K1-to-K2-u2},  if we had
 a fixed ratio of radii between consecutive levels,  and we could ignore the fact that $\gtw(\cdot)$ 
is zero outside $\sfIs_2$.  
To this end,  we consider \eqref{eq:K1-to-K2-u2} with $\vep$ replaced by $\vep':=\vep e^{\sfae}$ with $\sfae$ as in \eqref{def.nvep},
so that by \eqref{dr.1} and the following line we have $r_{j-1}/r_j = e$ for all $j$,
and assume that at $\vep'$,  uniformly over $\sfIs_1$,
\begin{equation}\label{def:wh-g1-excess}
e^{2(n'-\sfKo)} \,   \bar{\bE}^{(\sfKo)}_{y_1} \Big[ \gfo \big(\wh{\eta}(n') \big); 
 \bigcap_{j \in (\sfKo,n')}  \{ \wh{\eta} (j) >  j_{n_\vep}^{\gak} \},  \wh{\eta}(\sfKt) \notin \sfIs_2 
  \Big] = o_\vep(1) \,.
\end{equation}
(b).  The \abbr{lhs} of \eqref{eq:K1-to-K2-u2} at $\vep$ amounts to the same expression at $\vep'$
up to a certain non-random shift in $\wh{\eta}(\cdot)$.  This allows us to transfer
\eqref{eq:K1-to-K2-u2} from the special case of part (a),  to the general case.\\
(c).  By adapting the arguments in the proof of  \Cref{theo-win220}(d),  we establish \eqref{def:wh-g1-excess}.

\smallskip
\noindent
(a).  Proceeding to prove \eqref{eq:K1-to-K2-u2} at $\vep'$
instead of $\vep$,  note that $n_{\vep'}=n_\vep$ with the values of $\sfK_i \in \bN$,  $\sfIs_i$,  $i=1,2$ and 
$r_j$,  $j \le \sfKt$ also unchanged,  while the nominal barrier value $\gamma_0(j)$ for $\vep$
exceeds the one for $\vep'$ by 
\begin{equation}\label{eq:bar-gap}
\delta(j) := \sfrho_\vep \sfae + (\sfrho_{\vep} - \sfrho_{\vep'}) (n_\vep-j+1 )  = \sfco \sfae + o_\vep(1)
 \,.
\end{equation}
With $\log(r_{\sfK_i}/\vep') =  n_\vep-\sfK_i \in \bN$ we define 
$\wh{g}^{(\vep')}_{\sfK_i}(y_i)$,  $i=1,2$ analogously to \eqref{def:wh-g3} and subject to the 
obvious notational modifications,  our proof of \eqref{eq:g3-formula} also shows that
 for $i=1,2$,  uniformly in $y_i \in \sfIs_i$,
\begin{align}\label{eq:g12-formula}
\wh{g}^{(\vep')}_{\sfK_i}(y_i) & = (1+o_{\sfKfv,\sfKz,\vep} (1)) \alpha(\sfKfv,\sfKp) \ggi(y_i) 
\end{align}
(with the relative error $o_{\sfKz}(1)$ 
when $i=1$,  in lieu
of earlier relying on $\sfKth \to \infty$ when $\vep \to 0$).   Further,  having assumed
\eqref{def:wh-g1-excess},  we have \eqref{eq:g12-formula} holding for $i=1$ 
even with the additional constraint $\wh{\eta}(\sfKt) \in \sfIs_2$.  Thus,  
by the Markov property and \eqref{eq:g12-formula},  
\[ 
 \alpha(\sfKfv,\sfKp) \gon(y_1) = (1+o_{\sfKfv,\sfKz,\vep} (1))  \alpha(\sfKfv,\sfKp)
e^{2(\sfKt-\sfKo)} \,   \bar{\bE}^{(\sfKo)}_{y_1} \Big[ {\bf 1}_{\{ \wh{\eta} (j) >  j^{\gak},   j \in (\sfKo,\sfKt) \}} \, 
\gtw \big(\wh{\eta}(\sfKt) \big)
\Big] \,.
\]
Dividing both sides by $\alpha(\sfKfv,\sfKp)$ results with \eqref{eq:K1-to-K2-u2},  albeit at $\vep'$
(both $\gon(\cdot)$ and the \abbr{lhs} of \eqref{eq:K1-to-K2-u2} are independent of $\sfKfv$,  so 
after the fact we can remove $\sfKfv$ from the preceding $o_{\sfKfv,\sfKz,\vep} (1)$ error).  

\noindent
(b).  Now,  one gets from the \abbr{lhs} of \eqref{eq:K1-to-K2-u2} at $\vep'$ to the \abbr{lhs} expression at
$\vep$ by changing $\wh{\eta}(j) \to \wh{\eta}(j)-\delta(j)$ and replacing $y_1$ by $y_1 - \delta(\sfKo)$
(see \eqref{eq:eta-centered}).
In view of \eqref{eq:bar-gap},  for large $\sfKo$,  a slight change of the parameter $\kapo$ absorbs the shift in our 
constraints (barrier values),  whereas up to a $(1+o_{\sfKo}(1))$ factor,  such shift of the arguments of 
$\wh{g}_i(\cdot)$ merely multiplies 
both sides of \eqref{eq:K1-to-K2-u2} by the same factor $e^{2 \sfco \sfae}$ (as mentioned and used earlier, 
moving the borders of our windows $\sfIs_i$ by any constant,  or even a small power of $\sfK_i$, 
does not affect any of our statements).

\noindent
(c).  We proceed to establish \eqref{def:wh-g1-excess} at $\vep'$ and to this end,  
we will show that uniformly in $y \in \sfIs_1$,
\begin{align}\label{eq:bar-gap2}
e^{2 (n''-\sfKo)} \bP \big( \RR^{o,(y)}_{\sfKo} \cap \wh{\UU}^{(\vep',y)} \cap 
\{ \wh{\eta}(\sfKf) \in \sfIs_{4},  \wh{\eta}(\sfKfv) \in \sfIs_5,  \wh{\eta}(\sfKt) \notin \sfIs_2\} 
 \big) = o_\vep(1) \,,
\end{align}
where $\wh{\eta}(j)$ is defined using the nominal barrier $\gamma_0(\cdot)$ for $\vep'$,  while 
\begin{align*}
\RR^{o,(y)}_{\sfKo} :=  \!\!\!\! \bigcap_{j \in (\sfKo,n')} \!\!\!\! \{ \wh{\eta}(j) \ge  - j_{n_\vep}^{\gak} \} 
\, \quad\text{and}\quad 
\wh{\UU}^{(\vep',y)} &:= \bigcup_{\sfv\in \vep' \bF^\star} \Big\{ \bigcap_{j \in [\wt{n}_\vep^{(5)},n_\vep^{(6)}]}\!\!\!\!\!
\{\wh{\eta}_{\sfv} (j) \ge  - j_{n_\vep}^{\gak}\},  \eta_{\sfv} (n_\vep+1)= 0 \Big\} \,.
\end{align*}
Thus,  assume \eqref{eq:bar-gap2}.  Recall  \eqref{eq:K3-angles} and \eqref{def:g4},  that uniformly 
over $\sfIs_4$  and $\GG_{n'}^{\npp}$ the function $e^{2(n'-n'')} \gfo(\cdot)$ is within $1+o_1(\vep)$ 
of $\bP \big(\bigcup_{v} \{ \RU^{(\vep')} (v) \cap \, \CC^{(\vep')} (\cdot
, v) \} | \GG^{\npp}_{n'} \big)$.  As the event in this conditional probability
is contained within
$\wh{\UU}^{(\vep',y)} \cap  \{\wh{\eta}(\sfKfv) \in \sfIs_5\}$ (c.f.  \Cref{rem:alt-UU} or \eqref{eq:sharp-cont}),
and the upward curved barrier of \eqref{def:wh-g1-excess} 
is more restrictive than the downward curved barrier of $\RR^{o,(y)}_{\sfKo}$,  clearly
\eqref{def:wh-g1-excess} is a consequence of \eqref{eq:bar-gap2}.

It only remains to prove \eqref{eq:bar-gap2}.  Indeed,  our proof of \Cref{theo-win220}(d),  at $\sfk=\sfKp$,  
shows that the contribution to \eqref{lem-win1} from only $j=\sfKt$,  is $o_{\sfKt}(1) \hg(z)=o_\vep(1)$ 
(c.f.  \eqref{agar.66-new} and the summand for $j=\sfKt$ in \eqref{agar.66-ext}).  
This argument applies verbatim with window sizes $\sfIs_i$ at $i=2,4,5$,
instead of the smaller $\sfIp_j$ of \eqref{dfn:windows-alt}.  In contrast to
\eqref{lem-win1} where we consider also $o_{\sfKz}(1)$ contributions from levels near $\sfKo$,  here
with $\sfKo$ and $y_1$ independent of $\vep$,  our $o_\vep(1)$ smallness remains even when the rough barrier
$\Rkbp$ of \eqref{a216.0}  (which is part of $\whUUpk$ of \eqref{lem-win1}),  starts only
at $\sfKo$,  instead of at $\sfk=\sfKp$.  Having done so,  the preceding reasoning applies even 
when measuring the events of \eqref{lem-win1} during the $\sfs_{\sfKo}(y_1)$ first $\sfKo$-level excursions 
(instead of the first $\sfs_1(z)$ excursions at level one).  Fixing $\sfv_\bb=\npp$,  omitting the index $\bb$,
unraveling the notations  \eqref{a216.0},  \eqref{def:UUpk},  \eqref{14.4} and \eqref{a216.02} within \eqref{lem-win1}
and having now all events measured during the $\sfs_{\sfKo}(y)$ level $\sfKo$ excursions,  we arrive at \eqref{eq:bar-gap2}. 
\end{proof}
\noindent
{\bf IV.} We conclude the proof by showing the claimed bounds on $\alpha(\sfKfv,\sfKp)$.
\begin{proposition}\label{lem:alpha-unif-bd}
For $\sfKp$ large,  the $\alpha(\sfKfv,\sfKp)$ of \eqref{eq:alpha-tot} are bounded above and below 
by $\sfKp$-dependent positive constants,
uniformly in $\sfKfv$.
\end{proposition}
\begin{proof}  To show that $\alpha(\sfKfv,\sfKp) \le \bar \alpha(\sfKp)$,  we set
\begin{equation}\label{def:Zz}
Z:=\frac{\sqrt{2Y_{\la_z}}-\sqrt{2\la_z}}{\sqrt{\sfKf}} \quad \text{ for } \quad 
 \la_z:=\bar \sfs_{\sfKf}(\sqrt{\sfKf}\,  z)=\frac{1}{2}\big(\sfco \sfKf+\sqrt{\sfKf} \,  z\big)^{2}, 
\end{equation}
and note that 
\[
\sqrt{\sfKf}(Z+z) = \sqrt{2Y_{\la_z}}-\sqrt{2 \la_z}+\sqrt{\sfKf}\, z = \sqrt{2Y_{\la_z}}-\sfco \sfKf  \,.
\]
Thus,  setting $y=\sqrt{\sfKf} \, z$ in \eqref{eq:alpha-tot} results with 
\begin{equation}\label{alt:alpha-tot}
 \alpha(\sfKfv,\sfKp) = \frac{2e\sqrt{\sfKf}}{\sqrt{2\pi}} \int_0^\infty  
 \mathbb{E}\Big[\gfo\big(\sqrt{\sfKf}(Z+z)\big)\Big] z e^{\sfco \sqrt{\sfKf} \,  z} \, dz \,.
\end{equation}
Plugging 
\eqref{ubd2-g4} in \eqref{alt:alpha-tot},  after some algebra
it follows that for $h(x):=x_+ e^{-x^2/8}$,
\begin{align}\label{ubd-alpha}
\al(\sfKfv, \sfKp) &\leq C'( \sfKp)\int_{0}^{\ff} 
\mathbb{E}\Big[ h(Z+z) e^{-2\sqrt{\sfKf} Z } \Big] z\,dz \nonumber \\
& \le  C'(\sfKp) \int_0^\infty \bE \big[ h^2(Z+z) \big]^{1/2}  \bE \big[ e^{- 4 \sqrt{\sfKf} Z} \big]^{1/2}  z dz\,.
\end{align}
One clearly has that 
\[
 \bE \big[ e^{- 4 \sqrt{\sfKf} Z} \big] \le 1 + 4 \int_0^{\sqrt{2\lambda_z}} e^{4u} \,  \bP \big(\sqrt{\sfKf} Z \le -u \big) du \,.
\]
Hence,  by the Poisson tail estimate of \cite[(3.8)]{BRZ},  that for some $c>0$,  
\begin{equation}\label{eq:38-brz}
\bP(\sqrt{\sfKf} Z \le -u) \le e^{-c u^2} \,, \qquad \forall u \in [0,\sqrt{2\lambda_z}],   z \ge  0 \,,
\end{equation}
we deduce  that $\bE \big[ e^{- 4 \sqrt{\sfKf} Z} \big]$ is bounded above,  uniformly over $\sfKf$ and $z \ge 0$.
Further,  with $h^2(\cdot)$ uniformly bounded,  decreasing on $[2,\infty)$ and such that 
$h^2(x) \le (2+z)^{-6}$ for all $x \ge 5 \sqrt{\log (2+z)}$,  we get from \eqref{ubd-alpha} that 
$\sup_{\sfKfv}  \alpha(\sfKfv,\sfKp)$ is finite,  provided that 
\[
\sup_{\sfKf} \Big\{ \int_0^\infty z  \bP(Z \le -z + 5 \sqrt{\log (2+z)} )^{1/2} dz \Big\} < \infty \,,
\]
which in turn is also an easy consequence of the tail estimate \eqref{eq:38-brz}.

Next,  to show that $\alpha(\sfKfv,\sfKp) \ge \uu \alpha(\sfKp)>0$,  recall that 
$\inf_{\lambda  \ge 4} \Big\{ \bP\Big( \sqrt{2/\lambda} ($Poisson$(\lambda) -\lambda) \in [ -1,0]\Big) \Big\} > 0$.
Thus,  $Z$ of \eqref{def:Zz} is such that 
\[
\delta:= \inf_{z \in [3,4]} \{ \bP(0 \ge Z  \ge -1 ) \} \ge 0 \,.
\]
Restricting in the right hand side of \eqref{alt:alpha-tot} to $z \in [3,4]$ and $Z \in [-1,0]$,  we  
utilize the lower bound \eqref{eq:lbd-alpha} to arrive at
\[
 \alpha(\sfKfv,\sfKp) \ge \frac{4 e}{\sqrt{2\pi}} c(\sfKp)
 \inf_{z \in [3,4]} \bE \Big( e^{-2 \sqrt{\sfKf} Z} \,,  Z+z \in [2,4] \Big) \ge  \frac{4 e \delta}{\sqrt{2\pi}} c(\sfKp) \,,
\]
as claimed.  
\end{proof}
\section{Precise first moment: Proof of \Cref{prop:proxy-good-cover}}
\label{sec:good-cover}

Here we finish evaluating $\bP(\UUpy(y_1) | \underline{\sfx})$ up to $o(1)$ relative error,  using the notations 
$\sfs=\sfs_{\sfKo}(y_1)$ and $\btet = \btet(y_1) = \{\theta_q,  q \le \sfs \} := \{\theta_q^{(\sfKo)}, q \le \sfs\}$
for the angle differences,  denoting their empirical measures by $L_{\btet}$,  analogously to $\LKoy$ of 
\eqref{eq:almost-radial}.  For $\bar c$ and $\bar \mu$ of \eqref{dfn:Poiss} we set,  analogously to 
\eqref{dfn:hat-L-alpha},  
\begin{equation}\label{dfn:Vthet}
V_{\btet} := \frac{\bar c}{\sfm_\vep} \sum_{q=1}^{\sfs} \big[  \langle \cos(\cdot), \bar \mu \rangle 
- \cos (\theta_q) \big] 
\end{equation}
and consider the function 
\begin{align}\label{eq:K1-K3-alt}
\gKo(\btet(y_1)):= [r_{\sfKo}/r_{\sfKth}]^2\,  
\, \bE \big[ \gth\big( \wh{\eta}(\sfKth) \big) ;  \BWlo  \,| \btet(y_1) \big]  \,,
\end{align}
for the event $\BWlo$ of \eqref{eq:BW-K12-alt}.  As ${m_\vep}/{\sqrt{2 \sfs}} \to 1$ when $\vep \to 0$,
uniformly in $y_1 \in \sfIs_1$,  we have for $\Istar$ of \eqref{dfn:I-star-ext},
\begin{equation}\label{eq:bd-V-tet}
\{ (y_1,L_{\btet}) \in \Istar \} =
\{ \sfd_\infty (\bar \mu, L_{\btet}) \le \frac{2 \cst}{\sfm_\vep} \log \sfKo \} 
\subseteq \{ |V_{\btet}| \le  \, \cst \log \sfKo \} 
\end{equation}
(as $\bar c < 1$,   see also \eqref{dfn:dw-orig}-\eqref{dfn:dw}).  Further,  
$\gstb(y_1,V_{\btet}) \ge e^{-3 \sfKo^{\kapp}}$ for $\gstb (\cdot,\cdot)$ 
of \eqref{eq:g-def},  any large $\sfKz$ and all $(y_1,L_{\btet}) \in \Istar$.
\Cref{prop:proxy-good-cover} thus amounts to 
showing that for $\gstb$ and $\sfb$ of \eqref{eq:g-def},
\begin{align}\label{eq:V-tet}
\gKo(\btet(y_1)) = (1+o_{\sfKo,\vep} (1)) e^{-\sfb} \gstb (y_1,0) e^{-\sfco V_{\btet}} + o_\vep(1) \, ,
\end{align}
uniformly over $(y_1,L_{\btet}) \in \Istar$.
In proving \eqref{eq:V-tet},
we utilize the simpler proxy,  with neither  
higher barrier at levels $(\sfKo,\sfKt)$ nor any angle control
\begin{align}\label{dfn:wh-g1}
\whgo(\btet(y_1)) &:= [r_{\sfKo}/r_{\sfKth}]^2\,  
\, \bE \big[ \gth \big( \wh{\eta}(\sfKth) \big) ;  \Bl \,| \btet(y_1) \big], \\
\Bl & :=  \big\{ \wh{\eta}(\sfKt) \in \sfIs_2  \big\}
\bigcap_{j \in (\sfKo,\sfKt)} \big\{ \wh{\eta}(j) > j^{\gak}  \big\} 
\label{eq:BK12}
\end{align}
and our first step is the following unconditional result,
whose proof we defer to \Cref{sec:uncond-top}.
\begin{lemma}\label{lem:uncond-top}
Recalling $\bar \bP^{(\sfKo)}_{y_1}$ of \Cref{def:Fvj},  for some $c_0$ and any $\gak$,  $\kapo$,  $\kapp$
as in \eqref{dfn:aux-par},  uniformly over $\sfs$ such that $y_1 \in \sfIs_1$, 
\begin{align}\label{eq:K1-to-K3-u}
\bar \bE_{y_1}^{(\sfKo)} \Big[ \gth \big( \wh{\eta}(\sfKth) \big) ;  \Bl \,  \Big]
=  (1+o_{\sfKo,\vep} (1)) & [r_{\sfKth}/r_{\sfKo}]^2
\gstb (y_1,0) \, ,\\
\bar \bE_{y_1}^{(\sfKo)} \Big[ \gth \big( \wh{\eta}(\sfKth) \big) ;  \Bl \, ,
\sfd_\infty (\bar \mu, L_{\bpsi}) > c_0 \frac{\log \sfKo}{\sfm_\vep} 
\Big]
=  o_{\sfKo,\vep} (1) & [r_{\sfKth}/r_{\sfKo}]^2\,   \gstb (y_1,0) \, ,
\label{eq:K1-to-K3-trunc}
\end{align}
with $\bpsi = \{\psi_q, q \le \sfs\}$ the angle-differences of our $\sfKo$-level Brownian excursions.  
\end{lemma}

In the following we often abbreviate $\whgo(\btet(y_1))$ and $\gKo(\btet(y_1))$ as 
$\whgo(\btet)$ and $\gKo(\btet)$.  These functions
are
invariant under permutations of $\{\theta_q,  q \le \sfs\}$,   so \abbr{wlog} we make 
$q \mapsto \theta_q$ non-decreasing and proceed to couple the given excursions
with another set of $\sfs=\sfs_{\sfKo}(y_1)$ excursions where $\btet(y_1)$ is replaced by the non-decreasing,  
non-random $\sfs$-quantiles $\bmu = \bmu(y_1) := \{\bar \mu_q,  q \le \sfs\}$ of the law $\bar \mu$
(that is,  taking $\bar \mu_q:=\bar \mu^{-1}(q/\sfs)$).
The
variables in \eqref{eq:K1-K3-alt} and in \eqref{dfn:wh-g1}
involve only angle differences and excursion counts in concentric balls which are 
rotationally invariant,   so \abbr{wlog} the expectations there are with respect to $\sfs$ independent 
$\sfKo$-level Brownian excursions that start at $\sfz_q:=r_{\sfKo} e^{-i \theta_q}$ and 
are conditioned to end at $\sfu=(r_{\sfKo-1},0)$.  
Because of this we often write $\gKo(\btet)$ as $\gKo(\sfz)$. 
This notation has the advantage that,  since these variables are measurable on the
Brownian excursions starting from 
their first exit of the annulus $\sfB_\sfd(\npp, r_{\sfKo-1}) \setminus \sfB_{\sfd}(\npp,r_{\sfKo+1})$, 
the function $\gKo(\sfz)$ naturally extends
to 
any non-random starting points 
$\sfz=\{\sfz_q,  q \le \sfs\}$ inside this annulus (namely,  $|\sfz_q|\in (r_{\sfKo+1},r_{\sfKo-1})$).
We similarly write $\whgo(\btet)=\whgo(\sfz)$, 
$\gKo(\bmu)=\gKo(\sfz')$
and 
$\whgo(\bmu)=\whgo(\sfz')$,  where
$\sfz'_q := r_{\sfKo}e^{-i \bar \mu_q}$.  As
$\bar \mu$ has a $[c^{-1},c]$-valued density,  having
$\sfd_\infty(\bar \mu, L_{\btet}) \le 2 \boetaKo$ implies by \eqref{eq:ds-emp} and the triangle inequality,  that 
\begin{equation}\label{eq:quant}
| \theta_q - \bar \mu_q| \le \frac{3 \cst}{\sfm_\vep} \log \sfKo \,,   \qquad \forall q \le \sfs 
\end{equation}
(and all $\vep \le \vep_0$). 
Consequently,  $\sfz_q \in \bB_\sfd(\sfz'_q, \frac{1}{2} r_{\sfKo+1})$,  and in view of \Cref{lem-coupb} 
(with $r=r_{\sfKo}$),  setting 
\begin{equation}\label{dfn:bar-p-vep}
\bar p_\vep := C_{\sfKo}/\sfm_\vep,   \qquad C_{\sfKo} := 6 C \cst \log \sfKo, 
\end{equation}
the first hitting points $\zeta_{\sfz_q}$ and $\zeta_{\sfz'_q}$ of
$\partial \bB_\sfd(\sfz_q',r_{\sfKo+1})$ by the Brownian paths that started at
$\sfz_q$ and $\sfz'_q$,  respectively,  and are conditioned to end at $\sfu$,  are such that 
\begin{equation}\label{eq:zeta-tv-bd}
\sfd_{\rm TV}(\zeta_{\sfz_q},\zeta_{\sfz'_q}) \le C |\theta_q - \bar \mu_q| \le \frac{1}{2} \bar p_\vep \,.
\end{equation}
In particular, we have the following coupling 
\begin{equation}\label{dfn:coup}
\zeta_{\sfz_q} = (1-B_q) \zeta_{\sfz'_q}+B_q \xi_q 
\end{equation}
for some $\xi_q \in \partial \bB_\sfd(\sfz_q',r_{\sfKo+1})$ whose law is determined by 
$(\sfz_q,\sfz'_q)$, and \abbr{iid} Bernoulli($\bar p_\vep$) variables $\{B_q\}$ \emph{which 
are independent of everything else} (see \cite[Theorem 5.2]{L}).  In our setting
the balls $\bB_\sfd(\sfz_q',r_{\sfKo+1})$ are all 
strictly within the annulus $\bB_{\sfd}({\bf 0},r_{\sfKo-1}) \setminus \bB_{\sfd}({\bf 0},r_{\sfKo+1})$,
hence for $\zeta_{\sfz}:=\{\zeta_{\sfz_q},  q \le \sfs\}$ and $\zeta_{\sfz'}:=\{\zeta_{\sfz'_q},  q \le \sfs\}$ 
we have that 
\begin{equation}
\label{eq-011025a}
\whgo(\btet)=\whgo(\sfz)=\bE[\whgo(\zeta_{\sfz})|\sfz],  \qquad 
\whgo(\bmu)=\whgo(\sfz')=\bE[\whgo(\zeta_{\sfz'})|\sfz'],
\end{equation}
with the corresponding representations for $\gKo(\cdot)$.  Using \eqref{dfn:coup} and
\eqref{eq-011025a} yields the following coupling:
\newline
$\bullet$ Let $\whgo(\bmu)$ denote the expectation over a sample of
$\sfs$ independent $\sfKo$-level excursions of angle differences $\bmu$, 
taken 
by advancing each excursion to its first hitting time $\zeta_{\sfz'_q}$ and using 
the right hand side of \eqref{eq-011025a}.  
\newline
$\bullet$ Fix a sample $\uu{B}:=\{B_q, q \le \sfs\}$ of \abbr{iid} Bernoulli($\bar p_\vep$) variables,
independent of everything else.
\newline
$\bullet$ Denote by $\whes(\cdot)$ the analog of $\wh{\eta}(\cdot)$ due only to the 
$\sfKo$-excursions of angle differences $\bar \mu_q$  for which $B_q=0$.  
\newline
$\bullet$ Retain all those $\sfKo$-level excursions (with $B_q=0$), 
whereby $\zeta_{\sfz_q} = \zeta_{\sfz'_q}$,  while 
if $B_q=1$ replace the $q$-th excursion by a fresh 
one with $\zeta_{\sfz_q}$ drawn per the law of
$\xi_q$ (which depends on $\theta_q$ and $\bar \mu_q$).  
For $\ell \ge \sfKo$,  let $\Delj$ denote the number of $\ell$-level excursions within those fresh 
$\sfKo$-level excursions (that correspond to $B_q=1$),  and evaluate $\whgo(\btet)$ 
via the \abbr{lhs} of \eqref{eq-011025a} for the resulting total collection of $\sfs$ 
excursions,  with an extra (outer)
expectation over $\uu{B}$.

Utilizing this coupling,  we consider the following proxy for $\whgo(\cdot)$ of \eqref{dfn:wh-g1}
\begin{align}\label{dfn:wh-g1n}
\wtgo(\btet) &:= [r_{\sfKo}/r_{\sfKth}]^2\,  
\, \bE \big[  e^{-\frac{\sfco}{\sfm_\vep} \DelKth} \gth \big( \whes(\sfKth) \big) ;  \wtBl \big], \\
\wtBl & :=  \big\{ \whes(\sfKt) \in \sfIs_2  \big\}
\bigcap_{j \in (\sfKo,\sfKt)} \big\{ \whes(j) > j^{\gak}  \big\}\,, 
\label{eq:BK12n}
\end{align}
where the $\pm$ notation refers to our choice of $\gak$,  $\kapo$ and $\kapp$.  Specifically,  
in $\wtgm$ and $\wtBm$ we slightly decreased $\kapp$ 
while in $\wtBp$ we slightly decreased $\gak$ and $\kapo$.  
We prove in \Cref{sec:coupling} the following (sharp) comparison for the latter proxy.
\begin{lemma}\label{lem:quant-mu}
For any $\sfKo$ large enough and uniformly over  $(y_1,L_{\btet}) \in \Istar$,
\begin{align}
\wtgo(\btet) &= (1+o_{\vep}(1))  \wtgo(\bmu)\, e^{\sfco V_{\bmu}}  e^{-\sfco V_{\btet}} + o_\vep(1) \, ,
\label{eq:coupling-bdt}\\
\label{eq:coupling-bd1}
(1+o_{\vep}(1)) \wtgm(\btet) - o_\vep(1) \le \whgo(\btet) & \le (1+o_{\vep}(1)) \wtgp(\btet) + o_\vep(1) \, .
\end{align}
\end{lemma}

\noindent 
Having $\bpsi$ which are $\sfs$ \abbr{iid} samples 
from $\bar \mu$  amounts to $\sfKo$-excursions drawn from the (simpler) law $\bar \bP^{(\sfKo)}_{y_1}$.
Thus,  combining \eqref{eq:K1-to-K3-u} with \eqref{eq:K1-to-K3-trunc} we find in view 
of \eqref{dfn:wh-g1},  that for any $\cst \ge c_0$,  $\gak$,  $\kapo$,  $\kapp$
as in \eqref{dfn:aux-par},  
\begin{align}\label{eq:K1-to-K3n}
\bE \big[ \whgo(\bpsi);  \sfd_\infty(\bar \mu, L_{\bpsi}) \le 2 \boetaKo  \big]  &= (1+o_{\sfKo,\vep}(1)) \,  \gstb(y_1,0) \,,
\end{align}
uniformly over $\sfs$ such that $y_1 \in \sfIs_1$.  The right hand side of \eqref{eq:K1-to-K3n} is bounded away
from zero (as $\vep \to 0$),  so considering \eqref{eq:coupling-bd1} for $\btet=\bpsi$,  
with an indicator on $(y_1,L_{\bpsi}) \in \Istar$,  we find upon taking the expectation over $\bpsi$,  that 
\[
\bE \big[ \wtgm(\bpsi);  \sfd_\infty(\bar \mu, L_{\bpsi}) \le 2 \boetaKo  \big]  \le
 (1+o_{\sfKo,\vep}(1)) \,  \gstb(y_1,0) \le \bE \big[ \wtgp(\bpsi);  \sfd_\infty(\bar \mu, L_{\bpsi}) \le 2 \boetaKo  \big] \,,
\]
uniformly over $y_1 \in \sfIs_1$.  These bounds \emph{hold for any $\gak$,  $\kapo$,  $\kapp$
of \eqref{dfn:aux-par}},  hence uniformly in $y_1$ as above, 
\begin{align}\label{eq:K1-to-K3-uncond}
\bE
\big[ \wtgo(\bpsi);  \sfd_\infty(\bar \mu, L_{\bpsi}) \le 2 \boetaKo  \big]  
&= (1+o_{\sfKo,\vep}(1)) \,  \gstb (y_1,0) \,.
\end{align}
Moreover,  from \eqref{dfn:Vthet},  the variable $V_{\bpsi}$ is a normalized by ${\sfs}^{-1/2}$ sum 
of $\sfs$ \abbr{iid} $[-1,1]$-valued
variables of mean zero.
As $\vep \to 0$,  also $\sfs=\sfs_{\sfKo}(y_1) \to \infty$,  uniformly over $\sfIs_1$, 
so by the \abbr{clt},  for $\sfb$ of \eqref{eq:g-def} and any fixed $\sfc$,
\begin{equation}\label{eq:lim-V-psi}
\lim_{\vep\to 0} \bE[e^{-\sfc V_{\bpsi}}] = e^{\sfb(\sfc)} \,,  \qquad \sfb(\sfc):=(\sfc/\sfco)^2 \sfb \,.
\end{equation}
In addition,  we have 
for $y_1 \in \sfIs_1$ by our definitions of $\Istar$ and $\boetaKo$ that 
\[
\{ \bpsi :  (y_1,L_{\bpsi}) \notin \Istar \}  =
\{ \bpsi : \sfd_\infty (\bar \mu, L_{\bpsi}) > \frac{2 \cst}{\sfm_\vep} \log \sfKo \} \,.
\]
Applying  \eqref{lbl:d-inf-bd},  for the \abbr{IID} samples $\bpsi$ of \eqref{eq:K1-to-K3-uncond},
the probability of this event decays to zero 
(as $2 e^{-c' (\log \sfKo)^2}$,  some $c'>0$),  uniformly in $y_1 \in \sfI_1^{\star}$.   It then follows 
by Cauchy-Schwarz and \eqref{eq:lim-V-psi} that 
\begin{equation}\label{eq:lim-V-psi2}
\limsup_{\vep \to 0} \bE \big[ e^{-\sfco V_{\bpsi}} {\one}_{\{(y_1,L_{\bpsi}) \not\in \Istar\}} \big]
 = o_{\sfKo}(1) e^{2 \sfb} = o_{\sfKo}(1) \,.
\end{equation}
Consider now \eqref{eq:coupling-bdt} for $\btet=\bpsi$ as above,   with an indicator on $(y_1,L_{\bpsi}) \in \Istar$
and take the expectation over $\bpsi$.
In view of 
\eqref{eq:K1-to-K3-uncond} and \eqref{eq:lim-V-psi2},  we thereby find that uniformly over $y_1 \in \sfIs_1$,
\begin{align*}
 \gstb (y_1,0) = (1+o_{\sfKo,\vep}(1)) \, 
\bE \big[ \wtgo(\bpsi)  {\one}_{\{(y_1,L_{\bpsi}) \in \Istar\}} \big]
 & = (1+o_{\sfKo,\vep}(1))  \wtgo(\bmu)\, e^{\sfco V_{\bmu}}  e^{\sfb} + o_\vep(1) \,.
\end{align*}
In combination with \eqref{eq:coupling-bdt},  this implies in turn that uniformly over  $(y_1,L_{\btet}) \in \Istar$, 
\begin{equation}\label{eq:coupling-bd}
\wtgo(\btet) =  (1+o_{\vep}(1))  \wtgo(\bmu)\, e^{\sfco V_{\bmu}}  e^{-\sfco V_{\btet}} + o_\vep(1) =
(1+o_{\sfKo,\vep}(1)) e^{-\sfb} \,  \gst(y_1,0)  \, e^{-\sfco V_{\btet}}  + o_\vep(1)  \,.
\end{equation}
Combining \eqref{eq:coupling-bd1} and \eqref{eq:coupling-bd} we thus arrive at 
\begin{align}\label{eq:V-tet-alt}
\whgo(\btet) \,
&= (1+o_{\sfKo,\vep}(1)) \gstb (y_1,0)  e^{-\sfb}  e^{-\sfco V_{\btet}} + o_\vep(1) \,.
\end{align}
Our next lemma,  which is proved in \Cref{sec:pf-same-mean},  allows us to replace
$\whgo(\btet)$ of \eqref{dfn:wh-g1} 
in the \abbr{lhs} of \eqref{eq:V-tet-alt} 
by $\gKo(\btet)$ of \eqref{eq:K1-K3-alt},  thereby arriving at \eqref{eq:V-tet} and completing the proof 
of \Cref{prop:proxy-good-cover}. 
\bl \label{lem-} 
For any $\cst$ large enough,  uniformly in $(y_1,L_{\btet}) \in \Istar$
\bea
[r_{\sfKo}/r_{\sfKth}]^2\, 
\bE \big[ \gth \big( \wh{\eta}(\sfKth) \big) ;  \Bl \setminus \BWlo \,| \btet \big]  = 
o_{\sfKo,\vep}(1) \gstb (y_1,0) e^{-\sfco V_{\btet}} 
 \,.
\label{leftgoal.1}
\eea
\el

\subsection{Unconditional mean: Proof of \Cref{lem:uncond-top}}\label{sec:uncond-top}
Setting $L:=\log(r_{\sfKt}/r_{\sfKth})$,  recall $\ggi(y)={\bf 1}_{\sfIs_i} (y) y e^{-2 y}$ for $i=1,2,3$ and
\Cref{cor:gK2-gK1-arg}: 
\[
[r_{\sfKo}/r_{\sfKt}]^2  \,  
\bar \bE^{(\sfKo)}_{y_1} \Big[ \gtw \big( \wh{\eta}(\sfKt)\big);  \wh{\eta}(j) > j^{\gak},   j \in (\sfKo,\sfKt) \Big]
=  (1+o_{\sfKo,\vep} (1)) \gon (y_1) \,.
\]
By the Markov property of $j \mapsto \wh{\eta}(j)$ under the laws $\{ \bar \bP^{(\ell)} \}$,  it thus suffices for 
\eqref{eq:K1-to-K3-u} to show that,  
\begin{align}\label{eq:C23}
e^{2L} \,  \bar \bE_{y_2}^{(\sfKt)} \Big[ \gth \big( \wh{\eta}(\sfKth) \big) \Big] &= (1+o_{\vep}(1)) \gtw (y_2) \,,
\end{align}
uniformly over $y_2 \in \sfIs_2$.  To this end,  setting $m:=\sfs_{\sfKt}(y_2)$,
note that in \eqref{eq:C23} we have 
\begin{align}\label{eq:pre-clt}
\wh{\eta}(\sfKth) &= \sqrt{2 S_m} - \gamma_0(\sfKth) = \Gamma_m \wh{S}_m
+ y_2 + L \sfrho_\vep,
\end{align}
where $S_m$ denotes the sum of $m$ \abbr{IID} variables $\{\GJ'_i\}$ of the form 
\eqref{dfn:GJp} for 
$q:=1/(L+1)$
and 
\begin{align} 
\Gamma_m &:= \frac{2}{1+\sqrt{S_m/m}},  \qquad
\wh{S}_m := \frac{S_m-m}{\sqrt{2m}} \,.
\label{dfn:Gamma-m}
\end{align}
Here $L=o({n_\vep}/{\log n_\vep})$,  so $L (\sfrho_\vep-2) \to 0$ and \abbr{wlog} we replace $\sfrho_\vep$ by $2$ in \eqref{eq:pre-clt}.  Also,  $\sfKt \ge  \sfKth /3 \ge  L/3$,  being in $\sfIs_3$ amounts to 
$\wh{\eta}(\sfKth) \in [\sfKth^{\gak},\sfKth^{\okappa}]$ while having $y_2 \in \sfIs_2$ implies that
$y_2 \in [\sfKth^{\kapo},\sfKth^{\kapp}]$ for 
some $3/4>\okappa>\kapp>\kapo>\gak$.  Under these
constraints,  we can restrict in \eqref{eq:pre-clt} to 
\[
\Gamma_m |\wh{S}_m| = |\wh{\eta}(\sfKth) - y_2 - L \sfrho_\vep| \le O(\sfKth)\,.
\]
With $\sfKth = O((\log m)^{1/\gak})$ by an elementary manipulation of \eqref{dfn:Gamma-m},  
this implies in turn that 
\begin{equation}\label{eq:Gm-1}
|\Gamma_m -1 | |\wh{S}_m| = \frac{(\Gamma_m \wh{S}_m)^2}{\sqrt{8m}} =  O(\sfKth^2/\sqrt{m}) = o_m(1)\,,
\end{equation}
allowing us to omit the factor $\Gamma_m$ in  
\eqref{eq:pre-clt}.  Doing so and plugging also the explicit expressions for $\ggi(\cdot)$,  we re-write \eqref{eq:C23} 
in terms of $\wh{S}_m$,  as having uniformly in $y_2 \in \sfIs_2$, 
\begin{equation}\label{eq:C23-alt}
 \bE \Big[ \big(1+ y_2^{-1} \wh{S}_m \big) e^{-2 \wh{S}_m} ; \wh{S}_m + y_2 \in [\sfKth^{\gak},\sfKth^{\okappa}]
\Big] = (1+o_{\vep}(1)) \,  e^{2L} \,.
\end{equation}
As  $\bE \GJ'_1 =1$,  ${\sf var}(\GJ'_1)
=2L$ and $\sup_{L \ge L_o} \{ \bE\big[\exp(\GJ'_1/(2L))\big]\}$ is finite,  one has for $\delta_m:=\sqrt{2/m}$,  
\begin{equation}\label{eq:C23-mgf}
\bE\Big[ e^{-2 \wh{S}_m} \Big] = \bE\Big[ e^{-\delta_m (\GJ'_1-1)} \Big]^m = 
\Big[1+\frac{\delta_m^2}{2} {\sf var}(\GJ'_1) + O(L \delta_m^{3}) \Big]^m = 
(1+o_\vep(1)) e^{2 L} \,.
\end{equation}
In particular,  \eqref{eq:C23-alt} is thus a direct consequence of having 
\begin{equation}\label{eq:tilted}
 \wt{\bP}( |\wh{S}_m| \ge  \frac{1}{2} \sfKth^{\kapo}) \le 2 \sfKth^{-\kapo} \,  \wt{\bE}  |\wh{S}_m | = o_\vep(1) \,,
\end{equation}
for the tilted measure $\wt{\bP}$ corresponding to \abbr{IID} summands $\wt{\GJ}_i$ 
with $\bP(\wt{\GJ}_1=\ell) \propto  e^{-\delta_m \ell}  \bP(\GJ'_i=\ell)$.  It is not hard 
to check that $\bE\wt{\GJ}_1=\bE \GJ'_1 - \delta_m {\sf var}(\GJ'_1) + O(m^{-1})$ and
${\sf var}(\wt{\GJ}_1) = {\sf var}(\GJ'_1) + O(L^2 \delta_m)$.  We see that under $\wt{\bP}$ 
the variables $\wh{S}_m$ have the mean $-2L+o_\vep(1)$,  variance $L+o_\vep(1)$ 
and uniformly sub-Gaussian tails,  from which \eqref{eq:tilted} follows since $L = O(\sfKth^{\gak})$.

Turning to \eqref{eq:K1-to-K3-trunc},  we start with the following analog 
of \Cref{lem:BRZ-CH4} for $\sfKo$-level excursions from $\bar \bP_{y_1}^{\sfKo}$.
\begin{lemma}\label{lem:BRZ-CH4-bary}
For $\Bl$ of \eqref{eq:BK12} and integer $j'$, $j''$,  set  
\begin{align}
\bB^{j',j''}_{\sfKo,\sfKt,\sfKth} := \Bl \cap \{ \wh{\eta}(\sfKt) \in I_{j'} \} \cap \{ \wh{\eta}(\sfKth) \in I_{j''} \} \,.
\label{def:Bjp-jppn}
\end{align}
Then,  for $\hg(y):=y e^{-2y}$,  $\WW^{\in m}(\sfs)=\WW^{\ff,\in m}_{\npp,\sfKo}(\sfs)$ 
of \eqref{eq:LBTruncatedSumzf},  some finite $c_0$ and $C$,  any $y_1 \in \sfIs_1$ and $m \ge \log \sfKo$, 
\begin{equation}
\sum_{j''} \hg(j'')
\bar \bP_{y_1}^{(\sfKo)} \big( \bigcup_{j'}  \{ \WW^{\in m}(\sfs_{\sfKo}(y_1)) \cap
\bB^{j',j''}_{\sfKo,\sfKt,\sfKth} \} \big)
\le  C e^{-2 (\sfKth-\sfKo)} y_1 e^{-2 y_1} e^{-m^2/8}  \,,
\label{eq:K13-t}
\end{equation}
where the union is over $j' \in \sfI_{\sfKt}$ and the sum is over $j'' \in \sfI_{\sfKth}$.
\end{lemma}
\begin{proof} Using the short hand $\sfs:=\sfs_{\sfKo}(y_1)$,  recall similarly to \eqref{eq:BRZ-454},
that by \eqref{lbl:d-inf-bd} we have
for some $c_0$ finite,  the uniform bounds  
\begin{equation}\label{eq:BRZ-454nn}
\bar \bP^{(\sfKo)}_{y_1} \big( \WW^{\in m} (\sfs) \big) \le 2 e^{-m^2/4} \,,   \qquad \forall y_1 \in \sfIs_1,  m \ge 0 \,.
\end{equation}
With $\sfKth \le 2 (22 \log n_\vep)^{1/\gak}$,  it thus remains only to prove \eqref{eq:K13-t} for
$m \le C (\log n_\vep)^{1/(2 \gak)}$ and some universal $C<\infty$.  In contrast to \Cref{lem:BRZ-CH4},
here we have a fixed number $\sfs$ of excursions at level $\sfKo$ (which are further \abbr{iid}).
In such a setting we do not need the  counts of excursions at levels $k<\sfKo$ and
\cite[Lemmas 4.7 and  4.8]{BRZ2} hold with $\WW^{\in m} (\sfs)$ replacing $\LL'_{\sfKo,m,a}$ of
\eqref{def:Lmpa}.  That is, 
uniformly over our range of $y_1,j',j'',m$,
\begin{equation}\label{eq:BRZ-lem48}
\bar \bP_{y_1}^{(\sfKo)} \big( \WW^{\in m}(\sfs) \cap
\bB^{j',j''}_{\sfKo,\sfKt,\sfKth} \big)
\le c \,  H_{\sfKo,y_1,j',j''} (m) e^{-m^2/4} + e^{-\sqrt{n_\vep}} 
\end{equation}
for $H_{\sfKo,a,j',j''}(m)$ of \eqref{eq:BRZ-474},  while for $m \le 4 (\log \sfKth)^{1/2}$ 
we have a stronger bound,  with $H_{\sfKo,y_1,j',j''}(m)$ replaced in \eqref{eq:BRZ-lem48} by 
$G_{\sfKo,y_1,j',j''}(m)$ of \eqref{eq:BRZ-471}.  In either case,  multiplying both sides of this bound
by $g(j'')$ and summing over $j'$ and $j''$ yields \eqref{eq:K13-t} (note that $y_1^2 =o(\sfKt)$ and utilize \eqref{eq:BRZ-472} 
if $m \le 4 (\log \sfKth)^{1/2}$,  while if $m > 4 (\log \sfKth)^{1/2}$ use \eqref{eq:BRZ-473},  
for then $\sqrt{\sfKt} \le e^{m^2/32}$).
\end{proof}

As in our proof of \Cref{theo-win220}(e),  similarly to \eqref{eq:Lkp-WNa} we have 
for some $c_0$ finite and all $\vep$ small enough,
\begin{equation}\label{eq:LKo-WNa}
\Big \{ \sfd_\infty(\bar \mu,  L_{\bpsi}) >  c_0 \frac{\log \sfKo}{\sfm_\vep}  \Big\} 
\subset \bigcup_{m \ge  \log \sfKo} \WW^{\in m}(\sfs),   \qquad \forall y_1 \in \sfIs_1 \,.
\end{equation}
Further,  $\Bl \cap \{\wh{\eta}(\sfKth) \in \sfIs_3\}$ is contained in the union of 
$\bB^{j',j''}_{\sfKo,\sfKt,\sfKth}$ of \eqref{def:Bjp-jppn}  over integers 
$j' \in \sfI_{\sfKt}$ and $j'' \in \sfI_{\sfKth}$ (compare with \eqref{def:Bjp-jpp} and \cite[(4.69)]{BRZ2}). 
Recall that $\gth(y) \le y e^{-2y}$ and split the event on the \abbr{lhs} of \eqref{eq:K1-to-K3-trunc}  
according to which of the events $\bB^{j',j''}_{\sfKo,\sfKt,\sfKth} \cap \WW^{\in m}(\sfs)$ occur.
Upon doing so,  sum the bounds of \eqref{eq:K13-t} over $m \ge \log \sfKo$ to arrive at 
\eqref{eq:K1-to-K3-trunc}.

\subsection{Coupling with quantiles: Proof of \Cref{lem:quant-mu}}\label{sec:coupling}
Throughout this section $(y_1,L_{\btet}) \in \Istar$ (and in particular 
$L_{\btet} \in \bB^\infty_{\bar \mu}(2\boetaKo)$).  Our next lemma is the key to proving \eqref{eq:coupling-bdt}.
\begin{lemma}\label{lem:bd-hq}
Recall $\sfz_{q}=r_{\sfKo} e^{-i \theta_q}$,  $\sfz'_q=r_{\sfKo}e^{-i \bar \mu_q}$,
$\sfu=(r_{\sfKo-1},0)$ and $\bar p_\vep$ of \eqref{dfn:bar-p-vep}.   Let $T(\sfz_q)$ and $T(\sfz'_q)$
count the $\sfKth$-level excursions in an $\sfKo$ excursion
starting at $\sfz_q$,  or at $\sfz'_q$,  and ending at $\sfu$.   Then,   
\begin{align}\label{dfn:hq}
\htq  := \log \Big\{ 1+\frac{1}{\bar p_\vep} \Big(
 \frac{\bE \big[e^{-\sfco T(\sfz_q)/\sfm_\vep} \big]}{\bE \big[e^{-\sfco T(\sfz'_q)/\sfm_\vep} \big]} -1 \Big) \Big\} \,,
\end{align}
are such that uniformly over $(y_1,L_{\btet}) \in \Istar$ and all $\vep$ small enough,
\begin{equation}\label{eq:hq-bd}
\sup_{q} \Big | \sfm_\vep \htq - \frac{\sfco\,  \bar c}{\bar p_\vep} [ \cos(\theta_q)-\cos(\bar \mu_q) ] \Big| 
=o_\vep (1) \,. 
\end{equation}
\end{lemma}
\begin{remark}\label{rmk:hq} In view of \eqref{eq:quant} and \eqref{dfn:bar-p-vep},  it follows from 
\eqref{eq:hq-bd} that $\sup_q |\htq| \le C/\sfm_\vep$ for some $C$ finite.   Moreover,  with 
$\sfs \,  \bar p_\vep = O(\sfm_\vep)$,  it follows from \eqref{dfn:Vthet} and \eqref{eq:hq-bd} that 
\begin{equation}\label{eq:hq-sum}
| \bar p_\vep \sum_{q=1}^{\sfs} \htq - \sfco (V_{\bmu} - V_{\btet} ) | =o_\vep(1) \,.
\end{equation}
\end{remark}

\begin{proof} It is easy to check that $\bE [T(\sfz) ]$ is uniformly bounded above and below
over $|\sfz| = r_{\sfKo}$ and likewise $\sup_{\sfz} \bE [T(\sfz)^3] \le C \sfKth^2$.  
With $\bar p_\vep \sfm_\vep$ uniformly bounded in $\vep$,  by the inequalities 
$|e^{-x}-(1-x + \frac{1}{2} x^2) | \le |x|^3$ and  $|\log(1+x)-x| \le x^2$ whenever $|x| \le {1}/{2}$,  
we thus get \eqref{eq:hq-bd} upon showing that uniformly over $q \le \sfs$,  
\begin{align}\label{eq:bd-T-diff2}
\big| \bE[T(\sfz_q)^2]- \bE[T(\sfz'_q)^2] \big|  &= o_\vep(1) \,,  \\
\big| \bE[T(\sfz_q)]- \bE[T(\sfz'_q)] - \bar c [ \cos(\bar \mu_q)-\cos(\theta_q) ] \big| & = o_\vep(\bar p_\vep) \,.
\label{eq:bd-T-diff1}
\end{align}
Clearly,   $\sup_{\sfz} \bE [T(\sfz)^2] \le C \sfKth$ also for
$\sfz \in \Gamma_q := \partial \bB (\sfz'_q,r_{\sfKo+1})$,  with
$T(\sfz_q)=T(\zeta_{\sfz_q})$ and $T(\sfz'_q)=T(\zeta_{\sfz'_q})$ 
for the first hitting points $\zeta_{\sfz_q}$ and $\zeta_{\sfz'_q}$ 
of $\Gamma_q$ when starting at $\sfz_q$ and $\sfz'_q$ (see our justification of \eqref{eq-011025a}).
Further,  by the strong Markov property and \eqref{eq:zeta-tv-bd},  
\[
\big| \bE[T(\sfz_q)^2]- \bE[T(\sfz'_q)^2] \big| = \big| \bE[T(\zeta_{\sfz_q})^2]- \bE[T(\zeta_{\sfz'_q})^2] \big|  
 \le \sfd_{\rm TV}(\zeta_{\sfz_q},\zeta_{\sfz'_q})  \sup_{\sfz \in \Gamma_q} \bE[T(\sfz)^2] \le 
 C \sfKth \bar p_\vep 
\]
yielding \eqref{eq:bd-T-diff2}.  Turning to \eqref{eq:bd-T-diff1},  set $R:=r_{\sfKo-1}$,
$\rho:=r_{\sfKth-1}/r_{\sfKo-1}
$ and let $T'(\sfx)$ count the $\sfKth$-level 
excursions within a path starting at specific $\sfx \in  \partial \bB (\npp,\rho R)$ and ending at 
$\sfu=(R,0)$.  The unconditioned version of such $T'(\sfx)$,  
counting $\sfKth$-level excursions within a path from $\partial \bB (\npp, \rho R)$ to $\partial \bB (\npp, R)$,
is the Geometric variable $\GJo$ such that $\bP(\GJo=\ell)=(1-p) p^\ell$,  $\ell \ge 0$ 
for $p:=L/(L+1)$ and $L=\log (1/\rho)$,  see \cite[(2.1)]{BRZ2}.  Clearly,  $\bE[\GJo]=L$.  Moreover,  we have from 
\eqref{m1.5eu}
at $h_k=R$ and $h_{k'-1}= \rho R$,  the bound 
\[
\sup_{|\sfx| = \rho R} \Big\{ \big| \bE[T'(\sfx)]  - \bE[\GJo] \big| \Big\}  \le  
5 \sfKth \,  \rho \,  \bE[\GJo] \,.
\]
Starting at $\sfz \in \partial \bB_{\sfd}(\npp,R/e)$ and setting 
\[
A(\sfz):=\{\text{hitting} \;\; \partial \bB_{\sfd}(\npp,\rho R) \;\; \text{before exiting} \;\;
\bB_{\sfd}(\npp,R)  \;\;  \text{at}  \;\;  \sfu \}\,,
\]
we consequently have that 
\begin{equation}\label{eq:Tz-bd1}
\big| \bE[T(\sfz)] -  \bP(A(\sfz)) \bE[\GJo] \big| \le \bP(A(\sfz)) 
\sup_{|\sfx| = \rho R} \big\{ | \bE[T'(\sfx)] - \bE[\GJo] | \big\}
\le 5 \sfKth^2 e^{-(\sfKth-\sfKo)} = o_\vep(\bar p_\vep) \,.
\end{equation}
In computing $\bP(A(\sfz))$ our Brownian path is confined to $\bB(\npp,R)$,  
so \abbr{wlog} we set $R=1$ with a path starting at $\sfz=e^{-1} e^{-i \theta}$ and ending at $\sfu=(1,0)$
after hitting $\bB(\npp,\rho)$.  Then,  similarly to \eqref{10.11} we have by the strong Markov property,  that 
\begin{equation}\label{eq:Az-formula}
\bP(A(\sfz))= \frac{1}{p(\sfz,\sfu)} \bE_{\sfz} \big[ p(X_{\tau_\rho},\sfu); \tau_\rho < \tau_1 \big] \,,
\end{equation}
where $X_0=z$ with $\tau_r := \inf\{ t \ge 0 : |X_t |=r \}$ and the Poisson kernel
$p(\sfx,\sfu):=(1-|\sfx|^2)/|\sfx-\sfu|^2$  for the exit density of $\bB(\npp, 1)$ is such that 
\begin{equation}\label{eq:Pois-bds}
\sup_{|\sfx|=\rho} | p(\sfx,\sfu) -1 | \le 3 \rho \qquad \qquad \forall \rho \le \frac{1}{3} \,.
\end{equation}
Clearly,  $\inf \{ p(\sfz,\sfu) : |\sfz|=1/e \} = 1/c > 0$ and $\bP_{\sfz} (\tau_\rho < \tau_1 \big)= 1/L$
whenever $|\sfz|=1/e$ (see \cite[(2.1)]{BRZ2}),  so from \eqref{eq:Az-formula} and \eqref{eq:Pois-bds} we deduce that
\begin{equation}\label{eq:Tz-bd2}
\big| \bP(A(\sfz)) \bE[\GJo] - \frac{1}{p(\sfz,\sfu)} \big| \le 3 c \rho = o_\vep(\bar p_\vep)  \,.
\end{equation}
Finally,  note that
\[
\frac{1}{p(\sfz_q,\sfu)} - \frac{1}{p(\sfz'_q,\sfu)} =\frac{|\sfz_q-\sfu|^2-|\sfz'_q-\sfu|^2}{1-e^{-2}} =
 \bar c [ \cos(\bar \mu_q)-\cos(\theta_q) ] \,,
\]
which in view of \eqref{eq:Tz-bd1} and \eqref{eq:Tz-bd2} completes the proof of \eqref{eq:bd-T-diff1}.
\end{proof}

\begin{proof}[Proof of \eqref{eq:coupling-bdt}] Recall the first hitting points $\zeta_{\sfz_q}$ and $\zeta_{\sfz'_q}$ 
of $\Gamma_q := \partial \bB (\sfz'_q,r_{\sfKo+1})$ when starting at $\sfz_q$ and $\sfz'_q$,  see \eqref{dfn:coup}.
By the independence of our Brownian excursions,  we have for any $\ell \ge \sfKo$ the representation 
\begin{equation}\label{eq:Dlt-rep}
\Delj = \sum_{q=1}^{\sfs} B_q T_{q,\ell}
\end{equation}
for independent variables $\{T_{q,\ell},  B_q,  q \le \sfs\}$,  where $T_{q,\ell}$ counts the number of 
$\ell$-level excursions within a single $\sfKo$-level excursion that ends at $\sfu=(r_{\sfKo-1},0)$  
and starts at $\xi_q \in \Gamma_q$ such that 
\begin{equation}\label{eq:coup-st}
\bP(\xi_q \in A) =  \frac{1}{\bar p_\vep} \bP(\zeta_{\sfz_q} \in A) 
- \frac{(1-\bar p_\vep)}{\bar p_\vep} \bP(\zeta_{\sfz'_q} \in A) \qquad \forall A \subseteq \Gamma_q \,.
\end{equation}
With $\xi_q$ chosen according to the mixture \eqref{eq:coup-st},  
similarly to \eqref{eq-011025a},  we have for $T_q:=T_{q,\sfKth}$ and any $j \ge 0$, 
\begin{align*}
\bP(T_q=j)&=\bP(T(\zeta_{\sfz'_q})=j) + \frac{1}{\bar p_\vep} \big[
\bP(T(\zeta_{\sfz_q})=j) - \bP(T(\zeta_{\sfz'_q})=j) \big] \nonumber \\&=
\bP(T(\sfz'_q)=j) + \frac{1}{\bar p_\vep} \big[
\bP(T(\sfz_q)=j) - \bP(T(\sfz'_q)=j) \big] \,,
\end{align*}
from which we deduce by elementary algebra that
\begin{align}\label{dfn:hq-alt}
\bE \big[e^{-\sfco T_q/\sfm_\vep} \big] = e^{\htq}  \,  \bE \big[e^{-\sfco T(\sfz'_q)/\sfm_\vep} \big]  \,. 
\end{align}
In particular,  $h^{(\vep)}_q (\bar \mu_q)=0$.
Conditionally on $\uu{B}:=(B_1,\ldots,B_{\sfs})$,  the variable $\DelKth$ is independent 
of $\whes(\cdot)$ and $\wtgo$ depends on $\theta_q$ only via the law of $T_q$.  Thus,  setting 
\begin{align}\label{dfn:bar-g} 
\bgB (\uu{B}) &:= \bE \big[  \gth \big( \whes(\sfKth) \big) ; \wtBl \,| \uu{B} \big]
\prod_{\{q: B_q=1\}}  \bE \big[e^{-\sfco T(\sfz'_q)/\sfm_\vep} \Big] \, , \\
\bDl (\uu{B}) &:= \sum_{q=1}^{\sfs} (B_q - \bar p_\vep) \htq \,,
\label{dfn:bar-Dlt}
\end{align}
plug \eqref{eq:Dlt-rep} for $\ell=\sfKth$,  within \eqref{dfn:wh-g1n},  compute the conditional expectation 
given $\uu{B}$ and utilize \eqref{dfn:hq-alt},  to arrive at 
\begin{align}\label{dfn:g1n-alt}
\wtgo(\bmu) &= [r_{\sfKo}/r_{\sfKth}]^2\,   \bE \big[ \bgB  \big],\\
\label{dfn:g1t-alt}
\wtgo(\btet) &= [r_{\sfKo}/r_{\sfKth}]^2\,  \exp\big( \bar p_\vep \sum_{q=1}^\sfs \htq \big) 
\bE \big[ \bgB  e^{\bDl}  \big] 
\end{align}
(with both expectations over the vector $\uu{B}$ of \abbr{iid} Bernoulli$(\bar p_\vep)$ variables).
By \eqref{eq:hq-sum},  \eqref{dfn:g1n-alt},  \eqref{dfn:g1t-alt}
and the bound $|V_{\btet}|,  |V_{\bmu}| \le \cst \log \sfKo$ (as in \eqref{dfn:hat-L-alpha}),  
 our claim 
\eqref{eq:coupling-bdt} amounts to having uniformly over  $(y_1,L_{\btet}) \in \Istar$
\begin{align}\label{eq:cor-Dlt}
\bE \big[ \bgB e^{\bDl} \big] 
&= (1+o_{\vep}(1))  \bE \big[ \bgB \big] + o_\vep(e^{-2\sfKth}) \,.
\end{align}
Note that $\gth(\cdot) \in [0,1]$ for all small enough $\vep$ and thereby so is
$\bar g ^{(\vep)}$ of \eqref{dfn:bar-g}.  Setting $\delta_\vep := \bar p_\vep^{1/4} \to 0$,  
we lower bound the expected value of $\bgB e^{\bDl}$ by restricting it to $\bDl \ge -\delta_\vep$
and upper bounding it according to whether or not $\bDl \ge \delta_\vep$,  
arriving at \eqref{eq:cor-Dlt} upon showing that 
\begin{equation}\label{eq:Dlt-bdd-by-g}
\bE \Big[ e^{\bDl};  \bDl \ge \delta_\vep \Big]
+ \bP(\bDl \le -\delta_\vep) 
= o_{\vep}(e^{-2 \sfKth})  \,.
\end{equation}
Further,  the centered \abbr{cgf} of Bernoulli$(p)$ variables are such that for any $|\lambda| \le 1$ and $p \in [0,1]$, 
\[
\Lambda_p(\lambda) := \log \big( 1 + p (e^{\lambda}-1) \,  \big) - \lambda p \le \frac{e}{2} p  \lambda^2 \,,
\]
hence it follows from \eqref{dfn:bar-Dlt},  with $|h_q^{(\vep)}| \le C/\sfm_\vep$ (see  \Cref{rmk:hq}),  
$\sfs \le  \sfm_\vep^2/2$ and $\bar p_\vep \le 1/2$,  that 
\begin{equation}\label{eq:Dlt-mgf}
\log \bE \big [ e^{\lambda \bDl} \big] \le \sfs \Lambda_{\bar p_\vep} (C |\lambda|/\sfm_\vep) 
\le C^2 \bar p_\vep \lambda^2  \qquad \forall |\lambda| \le \sfm_\vep/C \,.
\end{equation}
By Markov's inequality and \eqref{eq:Dlt-mgf} at $\pm \lambda$, 
for $\lambda={\delta_\vep}/({2 C^2 \bar p_\vep})\in [1, \sfm_\vep/C]$, 
the \abbr{lhs} of \eqref{eq:Dlt-bdd-by-g} is at most
\[
e^{\delta_\vep} \bE \big[ e^{\lambda \bDl}\big]  e^{-\lambda \delta_\vep} + 
\bE \big[ e^{-\lambda \bDl}\big] e^{-\lambda \delta_\vep} \le 4 e^{-\delta_\vep^2/(4 C^2 \bar p_\vep)} =
o_{\vep}(e^{-2 \sfKth})
\]
(as $\sfKth \sqrt{\bar p_\vep} \to 0$).
This yields  \eqref{eq:Dlt-bdd-by-g} and thereby completes our proof of \eqref{eq:coupling-bdt}.
\end{proof}
Our next lemma,  which essentially replicates \cite[Lemmas 4.11 and  4.14]{BRZ2}, 
is the key to proving \eqref{eq:coupling-bd1}.
\begin{lemma}\label{lem:bd-Delj}
For some $C'$ finite,  any $\sfKo$ and all $\vep$, 
\begin{equation}\label{eq:bd-Delj}
\sup_{\btet} \{ \bP(\Delj \ge C' \,  \bar p_\vep \sfs) \} \le 2 \exp \big(-\frac{\bar p_\vep \sfs}{C' (\ell-\sfKo+1)} \big) \,, 
\qquad \forall y_1 \in \sfIs_1,  \;\;  \ell \ge \sfKo \,. 
\end{equation}
\end{lemma}
\begin{remark}\label{rmk:bd-Delj}
For $y_1 \in \sfIs_1$ and our choice of $\sfKth \ll n_\vep$ we have that 
$\sfs \le  \sfm_\vep^2/2 \le \gamma_0(\sfKth) \sfm_\vep$ and 
\[
2 \sfKth e^{-\bar p_\vep \sfs/(C' \sfKth)} = o_\vep(e^{-3 \sfKth}) \,.
\]
Recalling \eqref{dfn:bar-p-vep} that $\bar p_\vep = C_{\sfKo}/\sfm_\vep$,  we thus deduce 
from  \eqref{eq:bd-Delj}  and the union bound that 
\begin{equation}\label{eq:bd-Delu}
\sup_{y_1,\btet}  \bP \big( \Gamma_{\btet}^{(\vep)} \big) = o_\vep(e^{-3 \sfKth} ) \,,  \qquad
\Gamma_{\btet}^{(\vep)} := \bigcup_{\ell \in [\sfKo,\sfKth]} \{ \Delj \ge C' C_{\sfKo} \gamma_0(\ell) \} \,.
\end{equation}
\end{remark}
\begin{proof} Recall that $\DelKo = \sum_{q} B_q$ is a Binomial($\sfs,\bar p_\vep$),  whose \abbr{cgf}
is at most $\sfs \bar p_\vep (e^\lambda-1)$,  yielding 
by Markov's inequality at $\lambda=\log C_o$ that for large enough $C_o>1$,  
\begin{equation}\label{eq:bd-DelKo}
\bP\Big(\sum_{q=1}^{\sfs} B_q \ge C_o \bar p_\vep \sfs\Big) \le e^{-\bar p_\vep \sfs/C_o} 
\end{equation}
which is \eqref{eq:bd-Delj} at $\ell=\sfKo$.  Next,  setting $p_\ell:=1/(\ell-\sfKo+1)$,
recall that 
by Poisson kernel estimates 
$\bP(T_{q,\ell} = j ) \le {c} p_\ell^2 (1-p_\ell)^{j-1}/2$ for some finite $c$,  any $\ell > \sfKo$,  $q,j \ge 1$ 
and $\theta_q$ (see \cite[(4.124)]{BRZ2}).  By elementary algebra it then follows that 
at $\lambda_o>0$ such that $1-e^{-\lambda_o}= p_\ell/2$,  any $\ell > \sfKo$ and all $q$,  $\theta_q$,
\begin{equation}\label{eq:Tq-cgf}
\bE [e^{\lambda_o T_{q,\ell}}] \le 1 + \frac{c}{2}  \frac{p_\ell^2}{p_\ell - (1-e^{-\lambda_o})} = 1 + c \,  p_\ell 
\le e^{c \,  p_\ell} \;.
\end{equation}
With the representation \eqref{eq:Dlt-rep} for independent $\{T_{q,\ell},  B_q,  q \le \sfs\}$,  we deduce
from \eqref{eq:bd-DelKo} and Markov's inequality at $\lambda_o \ge 
 p_\ell/2$ 
of \eqref{eq:Tq-cgf},  that for any $C'>0$ and $\ell > \sfKo$,  
\begin{align*}
\sup_{\btet} \{ \bP(\Delj \ge C' \,  \bar p_\vep \sfs) \} &\le 
\bP\big(\sum_{q=1}^{\sfs} B_q \ge C_o \bar p_\vep \sfs\big) +
\sup_{\btet,\{q_i\}} \{ \bP(\sum_{i=1}^{C_o \bar p_\vep \sfs} T_{q_i,\ell} \ge C' \,  \bar p_\vep \sfs) \} \nonumber \\
& \le \exp(-\bar p_\vep \sfs/C_o) + \exp( (c \,  C_o - \frac{1}{2} C') p_\ell \,  \bar p_\vep \,  \sfs ) 
\end{align*}
and taking $C'=2 C_o (c+1)$ results with \eqref{eq:bd-Delj}.
\end{proof}
\begin{proof}[Proof of \eqref{eq:coupling-bd1}]  
Comparing \eqref{dfn:wh-g1} with \eqref{dfn:wh-g1n},  our claim is that for $\sfKo$ large,  
uniformly over $(y_1,L_{\btet}) \in \Istar$
\begin{align}
e^{-o_\vep(1)} \bE \big[ \whes(\sfKth) e^{-\sfco \whes(\sfKth) - \frac{\sfco}{\sfm_\vep} \DelKth} ;  \wtBm \big] 
\le \bE \big[ \wh{\eta}(\sfKth) e^{-\sfco \wh{\eta}(\sfKth)};  \wtBz \big] &+ o_\vep(e^{-2 \sfKth}) \nonumber \\
\le e^{o_\vep(1)} 
\, \bE \big[  \whes(\sfKth) e^{-\sfco \whes(\sfKth) - \frac{\sfco}{\sfm_\vep} \DelKth} ;  \wtBp \big]  &+ o_\vep(e^{-2 \sfKth}) \,,
\label{eq:bd-gopm}
\end{align}
where in view of \eqref{eq:BK12} and \eqref{eq:BK12n},
for some $\kappm<\kapp$,  $\kapom<\kapo$,  $\gakm <\gak$,  
\begin{align}
\wtBm & := \bigcap_{i=2,3} \big\{ \whes(\sfK_i) \in (\sfK_i^{\kapo},\sfK_i^{\kappm})  \big\}
\!\!\! \bigcap_{\ell \in (\sfKo,\sfKt)}\!\!\! \big\{ \whes(\ell) > \ell^{\gak}  \big\}\,, \nonumber \\
\wtBz &:=  \bigcap_{i=2,3} \big\{ \wh{\eta}(\sfK_i) \in (\sfK_i^{\kapo},\sfK_i^{\kapp})  \big\}
\!\!\! \bigcap_{\ell \in (\sfKo,\sfKt)} \!\!\! \big\{ \wh{\eta}(\ell) > \ell^{\gak}  \big\} \,,  \nonumber \\
\wtBp & := \bigcap_{i=2,3} \big\{ \whes(\sfK_i) \in (\sfK_i^{\kapom},\sfK_i^{\kapp})  \big\}
\!\!\! \bigcap_{\ell \in (\sfKo,\sfKt)} \!\!\! \big\{ \whes(\ell) > \ell^{\gakm}  \big\}\,.
\label{dfn:Bplus}
\end{align}
By definition $\wh{\eta}(\ell) \ge \whes(\ell)$,  $\gamma_0(\ell))+\whes(\ell) \ge 0$ and 
\begin{equation}\label{eq:Delj-alt}
\Delj = 
\big(\,  \gamma_0(\ell) +\frac{1}{2} \wh{\eta}(\ell) + \frac{1}{2} \whes(\ell) \big) \big( \wh{\eta}(\ell) -  \whes(\ell) \big)  \,.
\end{equation}
Hence,  it is easy to see that on the events $\wtBz$ and $\wtBl$
\begin{equation}\label{eq:bd-eta-etas}
0 \le \wh{\eta}(\ell) - \whes(\ell) \le 2 \Delj/\gamma_0(\ell)  \qquad \forall \ell \in [\sfKo,\sfKt] \cup \{\sfKth\} \,,
\end{equation}
from which we deduce that 
\begin{equation}\label{eq:Del-whes}
\begin{aligned}
2 C' C_{\sfKo} \le \sfKt^{\kapp} - \sfKt^{\kappm} \quad &\Longrightarrow \quad
\wtBm \cap \big( \Gamma_{\btet}^{(\vep)}\big)^{\complement} \subset \wtBz \,, \\
2 C' C_{\sfKo} \le (\sfKo^{\gak} - \sfKo^{\gakm}) \wedge ( \sfKt^{\kapo} - \sfKt^{\kapom})  \quad &\Longrightarrow \quad
\wtBz \cap \big( \Gamma_{\btet}^{(\vep)}\big)^{\complement} \subset \wtBp \,.
\end{aligned}
\end{equation}
In view of \eqref{dfn:bar-p-vep},  the \abbr{lhs} of \eqref{eq:Del-whes} holds for all $\sfKo$ large and 
$\vep$ small.  Further,  all expectations in \eqref{eq:bd-gopm} are over $[0,1]$-valued variables,  so by
\eqref{eq:bd-Delu} we establish \eqref{eq:bd-gopm} upon showing that 
\begin{equation}\label{eq:gth-star}
\begin{aligned}
\wtBm \cap \wtBz \cap \big( \Gamma_{\btet}^{(\vep)}\big)^{\complement}  \quad &\Longrightarrow \quad 
\wh{\eta}(\sfKth) - \whes(\sfKth) \le \frac{\DelKth}{\sfm_\vep} + o_\vep(1) \,,  \\
\wtBp \cap \wtBz \cap \big( \Gamma_{\btet}^{(\vep)}\big)^{\complement} \quad  &\Longrightarrow \quad 
o_\vep(\whes(\sfKth)) \ge \wh{\eta}(\sfKth) - \whes(\sfKth) \ge \frac{\DelKth}{\sfm_\vep} - o_\vep(1) \,.
\end{aligned}
\end{equation}
Note that  $0 \le \whes(\sfKth) \le \wh{\eta}(\sfKth) \le \sfKth$ on $\wtBz \cap \wtBl$.  Thus,  
substituting \eqref{eq:Delj-alt} at $\ell=\sfKth$ on the right hand side of \eqref{eq:gth-star}
yields the top inequality since $\sfKth (1- \frac{\gamma_0(\sfKth)}{\sfm_\vep}) =  o_\vep(1)$.
Moreover,  $\frac{\DelKth}{\gamma_0(\sfKth)} \le C' C_{\sfKo}=o_\vep(\whes(\sfKth))$
on $\wtBp \cap \big( \Gamma_{\btet}^{(\vep)}\big)^{\complement}$,  so with $\sfKth^2 = o_\vep(\sfm_\vep)$ 
the same reasoning yields also the bottom inequality of \eqref{eq:gth-star}.
\end{proof} 

\subsection{Removing higher barrier and angle control: Proof of \Cref{lem-}}\label{sec:pf-same-mean}

Fixing $(y_1,L_{\btet}) \in \Istar$,  for $j > \sfKo$ let $L_{\btet,j}$ denote the empirical measures of 
angle differences for the $j$-level excursions within the given $\sfs$ excursions of level $\sfKo$.
Comparing $\Bl$ of \eqref{eq:BK12} with $\BWlo=\BW^{(\vep),(y_1)}_{\npp,[\sfKo,\sfKt]}$ of \eqref{eq:BW-K12},   
upon plugging in \eqref{leftgoal.1}
the bound \eqref{eq:bd-V-tet} on $V_{\btet}$ and the values of $r_{\sfKo}$,  $r_{\sfKth}$ and
$\gst(y_1,0)$,  it follows by a union bound on the \abbr{lhs} of \eqref{leftgoal.1} that
\Cref{lem-} is a consequence of having for $\cst$ large enough,  some $c'>0$ and all
$j \in (\sfKo,\sfKt)$,  
the bounds
\begin{equation}\label{eq:Db-Wass}
\whgC(\btet) +  \whgw(\btet) {\bf 1}_{j \in \sfJw} \le  e^{-2 y_1} e^{-c' (\log  j)^2} \,,
\end{equation}
where $\sfJw=[\sfKo+3,\sfKw]$,  
\begin{equation}\label{dfn:Db-Wass}
\begin{aligned}
\whgC(\btet)&:= e^{2(\sfKth-\sfKo)} 
\bE \big[ \gth \big( \wh{\eta}(\sfKth) \big) ;  \Bl,   \wh{\eta}(j) \ge j^{\okappa}  
 \,| \btet \big]  \,, \\
\whgw(\btet) &:= e^{2(\sfKth-\sfKo)}
\bE \big[ \gth \big( \wh{\eta}(\sfKth) \big) ;  \Bl,   \wh{\eta}(j) < j^{\okappa},
L_{\btet,j} \notin \bB^{\infty}_{\bar \mu}\big( 2 \cst \boetaj \big)
 \,| \btet \big]  \,.
\end{aligned}
\end{equation}
Our proof then consists of the following three steps:\\
(a).  Reducing \eqref{eq:Db-Wass} to having for some $\gakm < \gak$,  $\kapom < \kapo$,  $\okappam <\okappa$ 
\begin{equation}\label{eq:Db-Wass-a}
\whgC + \whgw \,  {\bf 1}_{j \in \sfJw} \le  e^{-2 y_1} e^{-c' (\log  j)^2} \,,  \qquad \forall y_1 \in (\sfKo^{\kapom},\sfKo^{\kapp}),  \;\; \forall j \in (\sfKo,\sfKt), 
\end{equation}
now considering 
\begin{equation}\label{dfn:Db-Wass-a}
\begin{aligned}
\whgC &:= e^{2(\sfKth-\sfKo)} 
\bE \big[ \wh{\eta}(\sfKth) e^{-\sfco \wh{\eta} (\sfKth)} ;  \wtBp,   \wh{\eta} (j) > j^{\okappam}
\big]  \,, \\
\whgw &:= e^{2(\sfKth-\sfKo)}
\bE \big[ \wh{\eta} (\sfKth) e^{-\sfco \wh{\eta} (\sfKth)} ;  \wtBp,  \wh{\eta}(j) < j^{\okappa},
L_{\bpsi,j} \notin \bB^{\infty}_{\bar \mu}\big( \cst \boetaj \big)
 \big]  \,,
 \end{aligned}
\end{equation}
with expectations also over the vector $\bpsi$ of $\sfs$ \abbr{iid} samples from $\bar \mu$, 
independent of $\btet$,  and where now
\[ 
\wtBp := \bigcap_{i=2,3} \big\{ \wh{\eta}(\sfK_i) \in (\sfK_i^{\kapom},\sfK_i^{\kapp})  \big\}
\!\!\! \bigcap_{\ell \in (\sfKo,\sfKt)} \!\!\! \big\{ \wh{\eta}(\ell) > \ell^{\gakm}  \big\} 
\] 
(with a slight abuse of notation,  see \eqref{dfn:Bplus}).\\
(b).  Bounding $\whgC$ as stated,  by adapting the proof of \Cref{theo-win220}(d).\\
(c).  For $j \in \sfJw$ also bounding $\whgw$ as stated,  by adapting the proof of \Cref{theo-win220}(e).

\smallskip
\noindent
(a).
The probability of event $\AA_{\bpsi}:=\{ \sfd_\infty(\bar \mu, L_{\bpsi}) \le 2 \boetaKo \}$,  for a collection $\bpsi$ 
of $\sfs$ \abbr{iid} samples from $\bar \mu$,  which is independent of $\btet$,  is bounded below,  uniformly in 
$\cst \ge c_0$,  $y_1 \in \sfIs_1$,  $\sfKo$ and $\vep$.  \abbr{Wlog} we thus add $\AA_{\bpsi}$ as
part of the events of \eqref{dfn:Db-Wass}.  Now that
$(y_1,L_{\bpsi}) \in \Istar$,  upon ordering both $\btet$ and $\bpsi$,  we have that 
\eqref{eq:quant} also holds for $\psi_q$ instead of $\theta_q$.  This allows us to simultaneously 
couple the collections of $\sfs$ level $\sfKo$ excursions with angle differences $\bpsi$ 
and with angle differences $\btet$,  to such $\sfs$ excursions 
starting at $\sfz'_q=r_{\sfKo} e^{-i \bar \mu_q}$ and ending at $\sfu=(r_{\sfKo-1},0)$,
as we have  detailed after \eqref{eq-011025a}.  Moreover,  by \eqref{eq:zeta-tv-bd} and
the triangle inequality,  the first hitting points $\zeta_{\sfz_q}$ and $\zeta_{\sfz''_q}$ 
associated with the paths that started at $\sfz_q=r_{\sfKo} e^{-i \theta_q}$ 
and $\sfz''_q := r_{\sfKo} e^{-i \psi_q}$,  are then such that for $\bar p_\vep$ of \eqref{dfn:bar-p-vep},
\begin{equation}\label{eq:zeta-bd-n}
\sfd_{\rm TV}(\zeta_{\sfz_q},\zeta_{\sfz''_q}) \le 
\sfd_{\rm TV}(\zeta_{\sfz_q},\zeta_{\sfz'_q}) + \sfd_{\rm TV}(\zeta_{\sfz_q},\zeta_{\sfz'_q}) 
\le \bar p_\vep \,.
\end{equation}
Analogously to \eqref{dfn:coup},  for $\sfs$ 
\abbr{IID} Bernoulli($\bar p_\vep$) variables $B_q$ which are independent of $\btet$ and $\bpsi$,  the coupling 
\begin{equation}\label{dfn:coup-n}
\zeta_{\sfz_q} = (1-B_q) \zeta_{\sfz''_q}+B_q \xi_q 
\end{equation}
with some $\xi_q \in \partial \bB_\sfd(\sfz_q',r_{\sfKo+1})$ whose law is determined by 
$(\sfz_q,\sfz''_q)$,  allows us to replace in the expectations of  \eqref{dfn:Db-Wass} each starting point 
$\zeta_{\sfz_q}$ for which $B_q=0$,  by the corresponding $\zeta_{\sfz''_q}$.  
That is,  analogously to \eqref{dfn:wh-g1n} and \Cref{sec:coupling},
hereafter the quantities $\whes(\cdot)$ which are due only to such $\sfKo$-excursions 
be now formed by a Binomial($\sfs,1-\bar p_\vep$) number of $\sfKo$-excursions,  whose 
angle-differences are \abbr{iid} samples from $\bar \mu$,  independently of $\btet$, 
with $\Delj$ denoting the number of $\ell$-level excursions within the remaining 
$\sfKo$-level excursions (for which $B_q=1$).  While the starting points of the latter 
$\sfKo$-level excursions (or more precisely,  the laws of $\xi_q$),  depend on $\btet$ as well as 
on the realization of $\bpsi$ and $\uu{B}=(B_q)$,  the bounds of \Cref{lem:bd-Delj} hold
uniformly over all such starting points.  In particular,  see \eqref{eq:bd-Delu},  \abbr{wlog} we add 
\[
(\Gamma_{\btet}^{(\vep)})^{\complement} = \{ \Delj \le C' C_{\sfKo} \gamma_0(\ell),   \;\; \sfKo \le \ell \le \sfKth \} \,,
\]
as part of the events of \eqref{dfn:Db-Wass}.  In particular,  by the same reasoning as in \eqref{eq:bd-eta-etas},
we can assume \abbr{wlog} that for some $C$ finite, 
\begin{equation}\label{eq:wh-Ko-Kth}
0 \le \wh{\eta}(\sfKo) - \whes(\sfKo) \le C \log \sfKo,  \quad  \frac{\sfco}{\sfm_\vep} \DelKth \le C \log \sfKo \,.
\end{equation}
Moreover,  with 
$\gamma_j(0) \ge \frac{1}{2} \sfm_\vep$ for $j \le \sfKt$,  by 
\eqref{eq:dinf-comp} at $c:=\sqrt{1+\bar{c}^2} + \bar c$ and 
$\boeta = \cst \boetaj \le c$,  taking hereafter $\cst \ge 48 c \,  C' C$,  
we get that for $\wt{m} \ge \sfs_j(0)$ and
the angle-difference empirical measures $L_{\star,j}$ of our $m=\wt{m}-\Deljj$ 
excursions at level $j$ within those with $B_q=1$, 
\begin{equation}\label{eq:angle-cnt}
\{ \wh{\eta}(j) \ge 0,  \quad 
\Deljj \le C' C_{\sfKo} \gamma_0(j),  \quad 
L_{\btet,j} \notin \bB^{\infty}_{\bar \mu}\big(2 \cst \boetaj \big) \}
\;\; \Longrightarrow \;\; 
L_{\star,j}  \notin  \bB^{\infty}_{\bar \mu}(\cst \boetaj) \,.
\end{equation}
Combining \eqref{eq:angle-cnt} and the second implication of \eqref{eq:Del-whes} it then follows,  as in our proof of
the right-most inequality of \eqref{eq:bd-gopm},  that \eqref{eq:Db-Wass} is a consequence of 
having the bounds \eqref{eq:Db-Wass-a} for $\whgC$,  $\whgw$ of \eqref{dfn:Db-Wass-a}
and some $\gakm < \gak$,  $\kapom < \kapo$,  $\okappam <\okappa$.  Specifically,  this is because  
after using the right hand side of \eqref{eq:wh-Ko-Kth} to eliminate the contribution of $\DelKth$,  we can 
use its \abbr{lhs} to convert all $\whes(\cdot)$ counts to $\wh{\eta}(\cdot)$ as stated there 
(at the cost of slightly lowering our allowed range of $y_1$).

\noindent
(b).  In $\whgC$ we have \abbr{iid} excursions sampled from $\bar \bP^{(\sfKo)}$,  so replacing everywhere
$\kapo$ by $\kapom$ and combining 
the Markov property of $j \mapsto \wh{\eta}(j)$ under the laws $\{ \bar \bP^{(\ell)} \}$ with \eqref{eq:C23}, 
we find that for some $C<\infty$,  
\begin{align}
\whgC &= e^{2(\sfKth-\sfKo)} \sum_{y_2 \in \sfIs_2} 
\bar \bP_{y_1}^{(\sfKo)} \big[ \wh{\eta}(\ell) > \ell^{\gakm},  \ell \in (\sfKo,\sfKt),  \wh{\eta} (j) > j^{\okappam},  
\wh{\eta}(\sfKt) = y_2 \big] \bar \bE^{(\sfKt)}_{y_2} \big[ \gth(\wh{\eta}(\sfKth)) \big] 
\nonumber \\
&\le C e^{2(\sfKt-\sfKo)} \sum_{u=j^{\okappam}}^\infty \sum_y p_j(y_1,u,y) y e^{-2y} \,,
\label{eq:bd-whgC}
\end{align}
where $I_u:=[u,u+1)$,  the inner-most sum
is over all integer values of $y \in \sfI_{\sfKt}$ and 
\begin{align*}
p_j(z,u,y):= 
\bar{\bP}^{(\sfKo)}_{z} \big( \wh{\eta} (j) \in I_u,  \wh{\eta}(\sfKt) \in I_y,
\wh{\eta}(\ell)>0 \;\; \text{for all} \;\;  \sfKo< \ell < \sfKt 
\big) \,.
\end{align*}
With $1 \ll \sfKo \ll \sfKt \le C (\log n_\vep)^{c}$,  taking $\okappam > \kapp$,  here $u \gg z \gg 1$ 
and $y \gg 1$,  uniformly over $j \in (\sfKo, \sfKt)$,  $z \in \sfI_{\sfKo}$ and $y \in \sfI_{\sfKt}$.
Thus,  proceeding similarly to the derivation of \eqref{eq:ubd-pjzuy},  \blue{upon applying
the bound of \Cref{prop:BarrierSecGWPropUB} for $\KQ_{[\sfKo,\sfKt]}^{(\sfKo,j,\sfKt)}(z,u,y)$
with $\Delta_{\sfKt-\sfKo} \le 2/(\sfKt-\sfKo)$ 
(see \eqref{dfn:Delta}),
we}
find that for $f_{j,s}(u)$ of \eqref{def:fjs},
some constant $c$,  any such $z,y,u$ and $j \in (\sfKo, \sfKt)$, 
\begin{equation}\label{eq:ubd-pj}
e^{2(\sfKt-\sfKo)} p_j(z,u,y) y e^{-2y} \le c \,  \gon(z) f_{j,\sfKo}(u) \frac{y^2}{(\sfKt-j)^{3/2}} e^{-\frac{(y-u)^2}{16(\sfKt-j)}} \,.
\end{equation}
Bounding $\gon(z) \le u e^{-2z}$,  summing the right hand side of \eqref{eq:ubd-pj} 
over all integer $y \ge 0$ and plugging the 
value of $f_{j,\sfKo}(u)$,  we deduce that for some (other) constant $c$,  uniformly over $j \in (\sfKo,\sfKt)$ 
and our range of $z,u$,  
\begin{equation}\label{eq:ubd-pj-n}
e^{2(\sfKt-\sfKo)} \sum_y p_j(z,u,y) y e^{-2y} \le c \,  e^{-2z} \frac{u^5}{(j-\sfKo)^{3/2}} e^{- \frac{u^{2}}{16(j-\sfKo)}} \,. 
\end{equation}
With $a := 2 \okappam -1>0$,  we conclude this step by noting that for some $C<\infty$,  $c'>0$ and all $j \in (\sfKo,\sfKt)$,  
\[
 \sum_{u=j^{\okappam}}^\infty \frac{c \,  u^5}{(j-\sfKo)^{3/2}} e^{- \frac{u^{2}}{16(j-\sfKo)}} \le
C j^{3/2} e^{-c' j^{a}} \le e^{-2 c' (\log j)^2} \,.
\]

\noindent
(c).  Fix $\cst \ge c_0$ and $k \in \sfJw$.  Replacing $\kapo$ by $\kapom$,  
for any integer $a \in (k^{\gakm},k^{\okappa})$,  we let $N_{k,a}:=[\sfs_k(a+1)]$.  With
$L_{\bpsi,k}$ the empirical measure of $\{ \psi_q\}_{q=1}^{\sfs}$ for some   
$\sfs \in [\sfs_k(a),N_{k,a}]$ in case $\wh{\eta}(k) \in I_a$,  we consider for integer
$m \ge \log k$,  the events $\WW^{\ff,\in m}_{k}(\sfs)$ of \eqref{eq:LBTruncatedSumzf} and
$\LL'_{k,m,a}$ of  \eqref{def:Lmpa},  now at angles $\{\psi_q\}$.
Similarly to \eqref{def:Bjp-jpp} we also set for integers $j' \in (\sfKt^{\kapom},\sfKt^{\kapp})$,  
$j'' \in (\sfKth^{\kapom},\sfKth^{\kapp})$ the event
\begin{align}
\bB^{j',j''}_{k+3,\sfKt,\sfKth} &:= \bigcap_{\ell=k+3}^{\sfKt-1} \{ \wh{\eta}(\ell) > \ell^\gakm \} 
\cap \{ \wh{\eta}(\sfKt) \in I_{j'} \} \cap \{ \wh{\eta}(\sfKth) \in I_{j''} \} \,.
\label{def:jp-jpp}
\end{align}
Similarly to \eqref{eq:Lkp-WNa},  we then have that for all $\vep \le \vep_0$,  any $k \in \sfJw$ and $a$ as above, 
\begin{equation}\label{eq:Lk-WNa}
\Big\{ L_{\bpsi,k} \notin \bB^\infty_{\bar \mu} (\cst \boetak) \Big\} 
\subset \bigcup_{m \ge \log k} \WW^{\ff,\in m}_{k}(N_{k,a}) \,.
\end{equation}
Partition the event in the definition \eqref{dfn:Db-Wass-a} of $\whgk$,  according to the values of 
$a=\lfloor \wh{\eta}(k) \rfloor$,  $j' = \lfloor \wh{\eta}(\sfKt) \rfloor$,  using a union bound over
$j'' = \lfloor \wh{\eta}(\sfKth) \rfloor$ and $m$ from \eqref{eq:Lk-WNa},  to get
similarly to the end of our proof of \Cref{theo-win220}(e),  that for 
$\hg(j''):=j'' e^{-2 j''}$ and any $y_1 \in \sfIs_1$
\begin{equation}\label{eq:bd-whgk}
\whgk \le 
e^{2(\sfKth-\sfKo)} \sum_{m,j''} \hg(j'') 
\bar \bP^{(\sfKo)}_{y_1} \big( \bigcup_{a,j'} \LL'_{k,m,a} \cap \bB^{j',j''}_{k+3,\sfKt,\sfKth} \big)   \,.
\end{equation}
Even conditional upon 
$\wh{\eta}(\sfKo)$ and $\wh{\eta}(k)$,  our $\{\psi_q\}$ are \abbr{IID} samples from $\bar \mu$, 
hence as in \eqref{eq:BRZ-454}
\begin{align}\label{eq:BRZ-454n}
\bar \bP_{y_1}^{(\sfKo)} \big( \WW^{\ff,\in m}_{k}(N_{k,a}) \,  \big| \,  \wh{\eta}(k) \,  \big) &\le 2 e^{-m^2/4} \,, 
\end{align}
and consequently also
\begin{align}
\label{eq:BRZ-455n}
\bar \bP_{y_1}^{(\sfKo)}(\LL'_{k,m,a}) & \le 2 e^{-m^2/4} \bar \bP_{y_1}^{(\sfKo)} (\wh{\eta}(k) \in I_a) \,.
\end{align}
As in our proof of \Cref{lem:BRZ-CH4-bary},  by \eqref{eq:BRZ-454n} and our choice of $\sfKth$,  
it suffices to consider $m \le C (\log n_\vep)^{1/(2 \gak)}$ for some $C<\infty$.  Though 
$\wh{\eta}(k)$ is not fixed,  as before,  \cite[Lemmas 4.7 and  4.8]{BRZ2} hold 
under the law $\bar \bP_{y_1}^{(\sfKo)}$ and with our adjusted range of $a, j', j''$.  Thus, 
analogously to \eqref{eq:BRZ-lem48},  for $H_{k,a,j',j''}(m)$ of \eqref{eq:BRZ-474},  
we have now that uniformly over our range of $y_1,j',j'',k,m,a$,  
\begin{equation}\label{eq:BRZ-lem48n}
\bar \bP_{y_1}^{(\sfKo)} \big( \LL'_{k,m,a} \cap \bB^{j',j''}_{k+3,\sfKt,\sfKth} \big)
\le c \,  H_{k,a,j',j''} (m) \bar \bP_{y_1}^{(\sfKo)} \big( \LL'_{k,m,a} \big)  + e^{-\sqrt{n_\vep}} \,,
\end{equation}
with $G_{k,a,j',j''}(m)$ of \eqref{eq:BRZ-471} replacing $H_{k,a,j',j''}(m)$ when $m \le 4 (\log \sfKth)^{1/2}$.
In either case,  multiply both sides by $g(j'')$ and sum over $j'$,  $j''$,  utilizing 
\eqref{eq:BRZ-472} or \eqref{eq:BRZ-473},  according to $m \le 4 (\log \sfKth)^{1/2}$ or not.
Then,  employ \eqref{eq:BRZ-455n} on its right hand side and sum also 
on $m \ge \log k$,  to get for some $c'>0$,  
finite $C$,  any $k \in \sfJw$ and $y_1 \in \sfIs_1$, 
\begin{align*}
\whgk \le C e^{-3 c' (\log k)^2} 
e^{2(k- \sfKo)}
\sum_{a} \bar \bP^{(\sfKo)}_{y_1} ( \wh{\eta}(k) \in I_a) \hg(a)  \le 
e^{-2 y_1} e^{-2 c' (\log k)^2} 
\end{align*}
(the right inequality follows similarly to \eqref{eq:C23},  or from \cite[Prop. 1.4]{BRZ}). 
This completes the proof of  \eqref{eq:Db-Wass-a},  thereby also completing the proof of \Cref{lem-}.

\section{Second moment and tightness: Proof of Propositions \ref{prop:2nd-moment}
and  \ref{prop:tight}}\label{sec:2nd-mom-tight}  

\subsection{Second moment control: Proof of \Cref{prop:2nd-moment}}\label{sec:2nd-mom}
Fixing $z$,  $\dlc$,  the parameters of \eqref{dfn:aux-par} and $0<\gakm<\gak$, 
we set for each block $\bb$ in the inner city $\bA_u^o$, 
\begin{align}
\Pon_{\bb} &:=  \bP [ \, \UUpy(\Ybz) | \FFuz  ] \, ,  \qquad 
\Pm_{\bb,\sfKop} := \bP [ \,\UUmp(\YbKp) | \FFuz  ]  \,,
\label{dfn:P1-}
\end{align}
with $\sfKop \in \{\sfKo,\sfKo+3\}$ and $\UUpy(y)$ as in  \eqref{eq:good-cover-y},  whereas with $\RUy$ of \eqref{eq:barrier-event},
\begin{align}
\label{dfn:Umy}
\UUmy (y) &:=  \{ \LKoy \in \bB^\infty_{\bar \mu} (\boetasKo) \} 
 \!\!\!\!\!\!\!
\bigcap_{j \in \sfJm \cup \{\nvfv\}}   \!\!\!\!\!\!\!
\{ \whetajy \in \sfI_j  \} 
\bigcap_{j \in \sfJo} \{\whetajy > j_{n_\vep}^{\gakm} \} \bigcap \RUy \,,
\end{align}
similarly defining $\UUmp(y)$ at $\sfKop=\sfKo+3$,  apart from  
$\sfKo+3$ instead of $\sfKo$ (also in the definition of $\whetajy$ and at $\sfJm$,  so 
then $\sfJo = (\sfKo+3,\sfKt) \cup (\sfKth,n')$).
Note that  $\UUmy(y)$ requires $(y,\LKoy) \in \Istar$ and $\UUpy(y) \subset \UUmy(y)$ 
(compare \eqref{dfn:Umy} to $\BBy$ of \eqref{eq:barrier-event-z} and 
$\WWy$ \eqref{eq:was-event}).  Further,  in view of 
\eqref{eq:g3-to-wh-g3} and \eqref{leftgoal.1} we have that uniformly in $\bb$,
\begin{align}\label{eq:Po-Pm}
\Pon_{\bb} & \le \Pm_{\bb,\sfKo}  \le (1+o_{\sfKfv,\sfKz,\vep}(1)) \Pon_{\bb} 
\end{align}
with $o_{\sfKfv,\sfKz,\vep}(1)$ some non-random functions 
that go to zero as $\vep\to 0$,  then $\sfKz \to \infty$ followed by $\sfKfv \to \infty$.  
Hence,  setting $\gst (y,v) =  \gstb (y,v) \,  \sfIos (y,v)$ 
for the indicator $\sfIos(y,v)$ on $\sfI_{\sfKo} \times \sfIw$,  it follows from \eqref{UUpy-asymp} that 
\begin{align}
\Pm_{\bb,\sfKo} &= (1+o_{\sfKfv,\sfKz,\vep}(1)) 
\big[\vep e^{\sfKfv}/r_{\sfKo}\big]^2 \alpha(\sfKfv,\sfKp) e^{-\sfb} 
\bE [ \gst \big( \Ybz, \Vb (\LKoz) \big);
\LKoz \in \bB^\infty_{\bar \mu} (\boetasKo) | \FFuz ] 
\label{eq:Pm-form}
\end{align}
(again,  uniformly in $\bb$).  Recall \Cref{rmk:first-mom} that our proofs of Propositions \ref{prop:angle-diff} and \ref{prop:proxy-good-cover}  apply verbatim for any $\cst$ large,   any parameters satisfying \eqref{dfn:aux-par}
and also with $\sfKo+3$ instead of $\sfKo$.  In particular,  \eqref{eq:Pm-form} applies also with $\sfKo+3$ 
instead of $\sfKo$ and suitably adjusting the parameters for $\sfI_{\sfKop}$ and 
for the support 
of $\gst(\cdot,\cdot)$.

Next,  using $u=u(\bb)$ for the city within which the block $\bb$ is,  set for blocks
$\bb,\bb'$ in $\bA_u^o$,  
\begin{align}
\Ptw_{\bb'}&:= \!\!\!\!\! \sum_{\bb \ne \bb', \bb \in u(\bb')} \!\!\! \!\! \Ptw_{\bb,\bb'}\,  \qquad \text{where} \qquad
\label{dfn:P2+}
\Ptw_{\bb,\bb'} := \bP [ \,  \UUpp(\Ybpz) \cap \UUpy(\Ybz) | \FFuz  ] \,.
\end{align}
Writing $n'=\nvf$,  $n''=\nvfv$ with $\sfv_\bb$ the center of block $\bb$,
our next result is key to \Cref{prop:2nd-moment}.
\begin{proposition}    [Second moment bounds]
\label{prop:UpperBoundz}  
For any block $\bb'$,  partition the blocks in $u(\bb')$ according to
\begin{eqnarray}\label{eq-G0L}
\bG_{\ell}  :=  \left\{ \bb\in u(\bb') : 2 r_{\ell-1} \le \sfd(\sfv_{\bb'},\sfv_{\bb}) < 2 r_{\ell-2}\right\}.
\end{eqnarray}
Then,  for some $c=c(\sfKfv,\sfKp)$ finite,  any $\sfKp$,  
large enough $\sfKz$,  $\sfKfv$ and small enough $\vep>0$,
\begin{align}\label{eq:UpperBoundzr}
& \sup_{\bb \in \bG_\ell} \{ P^{(2)}_{\bb,\bb'} \} \leq  c \,  \Pm_{\bb',\sfKop}
e^{ -2(n''-\ell \vee \sfKo)}e^{- \frac{1}{8} (\ell \vee \sfKo)^{\gak}_{n_{\vep}}}\,, \qquad \forall \bb' \quad  
\forall \ell \le n_\vep - \sfKz + 2\,,
\end{align}
with $\sfKop=\sfKo$,  except for taking $\sfKop=\sfKo+3$ when $\ell \in \{\sfKo+1,\sfKo+2\}$.  Further,  for any block $\bb$,
\begin{align}\label{bd:Late-brch}
\Ptwl_{\bb} := \sum_{\{ \bb' : 0 < \sfd(\sfv_{\bb},\sfv_{\bb'}) \le 2 \vep e^{\sfKz-1} \} } \!\!\!\!\!\!\!\!
\Ptw_{\bb,\bb'} \le e^{-\frac{1}{3} \sfKfv^\gak} 
\Pon_{\bb} \,.
\end{align}
\end{proposition}
\begin{proof}[Proof of  \Cref{prop:2nd-moment}]
Recall our tiling of the inner city $\mathbb{A}^{o}_{u}$ 
of edge-length $(1-6\dlc) e^{-\sfKz}$ by blocks of edge-length $\vep e^{\sfKfv}$ each. 
It is thus easy to see that $\bG_\ell \neq \emptyset$ only when $\ell \ge \sfKz-\ello$,
in which case $|\bG_\ell| \le c' e^{2(n''-\ell)}$ for some universal finite constant $c'$.
Further,  if $\bb \in \bG_\ell$, $\ell > n_\vep - \sfKz + 2$ then necessarily
$\sfd(\sfv_{\bb'},\sfv_{\bb}) \le 2 \vep e^{\sfKz-1}$.   
Thus,  assuming \abbr{wlog} that $\sfKz \ge \sfKf + \ello$ and combining 
\eqref{eq:UpperBoundzr} and \eqref{bd:Late-brch} we find that 
\begin{align}
& \Ptw_{\bb'} =  \sum_{\bb\neq \bb', \bb \in u(\bb')} P^{(2)}_{\bb,\bb'}
\label{eq:SecondMomentDecomp}
\le \Ptwl_{\bb'} 
+ \sum_{\ell_{n_\vep} \ge \sfKz-\ello} \sum_{\bb \in \bG_{\ell}} P^{(2)}_{\bb,\bb'} 
\nn \\
& \le  e^{-\frac{1}{3} \sfKfv^\gak} \Pon_{\bb'} + 
 + c\, c' \Big [ \big(
 \sum_{\ell = \sfKz-\ello}^{\sfKo}  
 e^{2(\sfKo-\ell)- \frac{1}{8} \sfKo^{\gak}} 
+  \sum_{\ell = \sfKo+1}^{n_\vep-\sfKz+2}  e^{-\frac{1}{8} \ell_{n_\vep}^{\gak}} \Big) \Pm_{\bb',\sfKo} 
+ 2  e^{-\frac{1}{8} \sfKo^{\gak}} \,   \Pm_{\bb',\sfKo+3} \Big]
 \nn \\
 &
 \leq e^{-\frac{1}{3} \sfKfv^{\gak}} \Pm_{\bb',\sfKo} 
 + c  \,  \Big[ 
 c'' \big( e^{2(\sfKo-\sfKz)} e^{-\frac{1}{8} \sfKo^{\gak}} +
 e^{-\frac{1}{8}\sfKz^{\gak}}  \big) \Pm_{\bb',\sfKo} + 2 c' e^{-\frac{1}{8} \sfKo^{\gak}} \Pm_{\bb',\sfKo+3} \Big] 
 \end{align}
 for some $c''=c''(\ello)<\infty$.  Since $\sfKo \to \infty$ before $\sfKfv$,
 and recalling that $\sfKo-\sfKz = o(\sfKo^\gak)$,  we deduce that for any fixed $\sfKp$,  uniformly in $\bb'$,
\begin{equation}
\Ptw_{\bb'} \leq  o_{\sfKfv,\sfKz,\vep} (1)\,  \Pm_{\bb',\sfKo} + e^{-\frac{1}{9} \sfKo^{\gak}} \Pm_{\bb',\sfKo+3} \,.
\label{glb.6}
\end{equation}
Now,  recall \eqref{dfn:wt-Q},  \eqref{dfn:wt-Q2} that 
\[
\wt{Q}^{(i,\vep)}_{u,z} = \sum_{\bb' \in u} \Pii_{\bb'},  \qquad i=1,2 
\]
(where we also utilized \eqref{eq:good-cover-y}).  Recall that $\sup_{\sfKfv} \alpha(\sfKfv,\sfKp)<\infty$ 
(see \Cref{lem:alpha-unif-bd}),  and on the support of $\gst(y,v)$ we have that 
$y e^{-\sfco v} \le \sfKo^{\cstp}$ for some $\cstp(\cst)<\infty$.  So,  summing \eqref{glb.6} over 
the $[(1-6\dlc) e^{-\sfKz}/(\vep e^{\sfKfv})]^2$ blocks $\bb' \in u$,  using \eqref{eq:Po-Pm} for 
the term at $\sfKo$ and \eqref{eq:Pm-form} for the one at $\sfKo+3$, 
we have that for $\wt{Q}^{(\star,\vep)}_{u,z}$ of \eqref{dfn:Qst},
some $c=c(r_0,b,\sfKp)$ finite,  any large $\sfKfv$ and $\sfKz \ge \sfKz'(\sfKfv)$,  
$\vep$ small enough and all $u,z$, 
\begin{equation}\label{eq:Q2-Q1-Qs}
\wt{Q}^{(2,\vep)}_{u,z} \le  o_{\sfKfv,\sfKz,\vep} (1) \wt{Q}^{(1,\vep)}_{u,z}
+ c \,  e^{-\frac{1}{10} \sfKo^{\gak}} e^{2 (\sfKo-\sfKz)} \,  \wt{Q}^{(\star,\vep)}_{u,z} \,.
\end{equation}
In view of \Cref{lem:2mm-star},  summing this bound over $u$ results with  \eqref{eq:sec-mom-bd} for 
$\dst(\sfKfv,\sfKz)$ which is the maximum of $e^{-\frac{1}{11} \sfKo^{\gak}}$,  see \eqref{eq:2mm-star},
and twice the $\varlimsup_\vep o_{\sfKfv,\sfKz,\vep} (1)$ of \eqref{eq:Q2-Q1-Qs}.
\end{proof}

Fixing the city $u$ and using the short-hand $\FF:=\FFuz$ we devote the remainder of this section
to establishing \Cref{prop:UpperBoundz}.  To this end,  similarly to \eqref{eq:B-K3} denote by
$\wh{\eta}^{\sharp,(y)}_{\bb,n''}$ the analog of $\whetaynpp$ in terms of excursions 
from $\partial \bB_\sfd(\sfv_\bb,r_{n'-1})$ to $\partial \bB_\sfd(\sfv_\bb,r_{n''})$.  Then,  set the following 
proxies for $\Pm_{\bb,\sfKop}$ on the right hand side of \eqref{eq:UpperBoundzr},
\begin{equation}\label{dfn:Pml}
P_{\bb,\ell} := \sum_{y_4,y_5} \alpha_{\sfKfv,\sfKp}(y_5|y_4) \bP\big( \BBl( \YbKp, y_4, y_5) | \FF \big),  
\end{equation}
where the sum is over $y_i \in \sfIs_i$,  with $\sfJl=\sfJupl \setminus \{\ell-2,\ell-1\}$ and
\begin{equation}
\label{def:B-no-ell}
\BBl (y,y_4,y_5) :=  \{ L^{(y)}_{\bb,\sfKop} \in \bB^\infty_{\bar \mu} (\boetasKop) \} 
 \bigcap_{j \in \sfJ_\ell} \{ \whetajy \in \sfI_ j \} 
\bigcap \{ \whetaynp = y_4, \widehat{\eta}_{\bb,n''}^{\sharp,(y)} = y_5 \} \,,
\end{equation}
at $\sfKop=\sfKo$,   except for taking $\sfKop=\sfKo+3$ 
and $\sfJl=\sfJupl \cap [\sfKop,n'']$ when  $\ell \in  \{\sfKo+1,\sfKo+2\}$.  Indeed,   our next lemma
shows that dropping the barrier at two levels 
has a tolerable effect on $\Pm_{\bb,\sfKop}$.
\begin{lemma}\label{lem:2-levels} 
For some 
 $c(\sfKfv)$ finite,  large enough $\sfKfv$,  $\sfKz$ and any $\vep$ small enough
\begin{equation}\label{eq:drop2-levels}
P_{\bb,\ell}
\le c(\sfKfv) 
\,  \Pm_{\bb,\sfKop} \qquad \forall \bb,  \ell \in (\sfKo,n'] \,,
\end{equation} 
where $\sfKop=\sfKo$ except for $\sfKop=\sfKo+3$ when $\ell \in \{\sfKo+1,\sfKo+2\}$.
\end{lemma}
\begin{proof} For any block $\bb$ and $y=\YbKp$ we set,  analogously to \eqref{dfn:Umy},  the events
\begin{align}\label{dfn:UUpl}
\UUpl (y) := & 
 \{ \LKop \in \bB^\infty_{\bar \mu} ( \boetasKop) \}
\bigcap_{j \in \sfJl}  \{ \whetajy \in \sfI_j  \} \bigcap \RUy
\end{align}
and by the argument leading to \eqref{eq:pre-precise},  deduce that for all $\bb$ and $\ell$, 
\begin{equation}\label{dfn:Pml-alt}
P_{\bb,\ell} = (1+o_{\sfKfv}(1)) \bP \big(\UUpl (\YbKp)|\FF\big) \,.
\end{equation}
(a) Taking $\ell \in \{\sfKo+1,\sfKo+2\}$,  compare
\eqref{dfn:UUpl} with \eqref{dfn:Umy} at $\sfKop=\sfKo+3$,  to see 
that $\UUpl (y) \subset \UUmp (y)$ yielding \eqref{eq:drop2-levels} as an
immediate consequence of \eqref{dfn:Pml-alt} and the right hand side of \eqref{dfn:P1-}.

We proceed with $\ell \ge \sfKo+3$ where the starting level $\sfKop$ of $\UUpl(y)$ 
and the constraints imposed there,  match those of $\UUmy(y)$.  In view of \eqref{dfn:P1-} 
and \eqref{dfn:Pml-alt},  it suffices to 
compare the events $\UUpl(y)$ and $\UUmy(y)$ for fixed $\bb$ 
and \abbr{wlog} also $(y,\LKoy) \in \Istar$,  as we do next.\\
(b).  $\ell \in (\sfKt+2,\sfKth]$.  In this range,  comparing \eqref{dfn:Umy} with \eqref{dfn:UUpl} 
we see that \eqref{eq:drop2-levels} trivially holds since $\sfJl=\sfJupl$ and $\UUpl (y) \subset \UUmy (y)$. \\
(c).  $\ell \in \{\sfKt+1,\sfKt+2,\sfKth+1,\sfKth+2\}$.  For $\ell =\sfKt+i$,  $i=1,2$,  the removal 
of constraints at levels $\{\ell-2,\ell-1\}$ amounts to decreasing $\sfKt$ by $3-i$,  
while for $\ell=\sfKth+i$,  $i=1,2$,  it amounts to increasing $\sfKth$ by $i$ (where
to bound $P_{\bb,\ell}$ first replace $\gak$ by $\gakm$ and increase $\okappa$ 
outside its new boundary levels $\sfKt+i-3$,  or $\sfKth+i$,   then replace $\kapo$ by 
$\gak$ and $\kapp$ by $\okappa$ at the remainder of $\sfJm$).  
Recall \eqref{eq:flex-Ki} and \Cref{rmk:first-mom} that such adjustments of $\sfKt$,  $\sfKth$ 
and of our parameters do not affect the validity of \eqref{UUpy-asymp},  so in this case
$P_{\bb,\ell} \le (1+o_\vep(1)) \Pm_{\bb,\sfKo}$. \\
(d).  $\ell \in [\sfKth+3,n']$.  As in part (c),  first enlarge in $P_{\bb,\ell}$ the windows 
$\sfI_j$ to be $(j_{n_\vep}^{\gakm},\infty)$,  for all $j \in (\ell,n')$.  Setting $\gfl(\cdot)$ as in 
\eqref{def:wh-g3-ell},  except for replacing $\gak$ by $\gakm$ at levels $(\ell,n')$,
note that by the Markov property of $\wh{\eta}(j)$ under $\bar{\bP}^{(\sfKth)}_{y_3}$,  it suffices for
\eqref{eq:drop2-levels} to show that for large $\sfKfv$,  small $\vep>0$ and all allowable $x' \in \sfI_{\ell-3}$ 
\begin{equation}\label{eq:bd-drop-ell}
\bar{\bE}^{(\ell-3)}_{x'} \Big[ \gfl \big(\wh{\eta}(\ell)\big); \wh{\eta}(\ell) \in \sfI_\ell \Big] \le c(\sfKfv) 
\bar{\bE}^{(\ell-3)}_{x'} \Big[ \gfl \big(\wh{\eta}(\ell) \big); \wh{\eta}(\ell) \in \sfI_\ell,   
\wh{\eta} (\ell-i) > (\ell-i)_{n_\vep}^{\gakm},  i=1,2 \Big] \,.
\end{equation}
By \eqref{eq:gfl-comp} of \Cref{rmk:g3-formula-gen},  we get \eqref{eq:bd-drop-ell} 
upon having,  uniformly in $\ell$,  $x' \in \sfI_{\ell-3}$ and $y \in \gamma_0(\ell) + \sfI_\ell$,  that
\begin{equation}\label{eq:bd-d-ell-alt}
\bar{\bP}^{(\ell-3)}_{x'} \Big[ \eta(\ell) \in [y,y+2) \Big] \le c(\sfKfv) 
\bar{\bP}^{(\ell-3)}_{x'} \Big[ 
\wh{\eta} (\ell-i) > (\ell-i)_{n_\vep}^{\gakm},  i=1,2,  \eta(\ell) \in [y,y+2) \Big] \,.
\end{equation}
The \abbr{lhs} of \eqref{eq:bd-d-ell-alt} is $q_3(x,[y,y+2))$ of \Cref{dfn:Markov}  
at $x:=\gamma_0(\ell-3)+x'$.  Here
$x \wedge y \ge \gamma_0(\ell) \ge \sfrho_\vep \ell_{n_\vep}$ 
and $|x-y| \le (\sfrho_\vep \ell_{n_\vep})^{\okappa}$.  Further,  for 
$u_o$ of \Cref{lem:g4-cnc}(a),  the probability on the right hand side of \eqref{eq:bd-d-ell-alt} 
is at least $q^{(u)}_{3}(x,[y,y+2))$ with $u = \ell_{n_\vep}^{\gak} /2 \ge u_o$
once $\sfKfv$ is large enough,  whereby \eqref{eq:bd-d-ell-alt} follows from 
\eqref{eq:q3-bd}.
\newline
(e).  $\ell \in [\sfKo+3,\sfKt]$.  Here $\sfKop=\sfKo$ and for
$\btet=\btet(y_1)=\{\thbKo,  q \le \sfs_{\sfKo}(y_1)\}$ with $y_1 = \Ybz$ we set 
$\whgo(\cdot)$ as in \eqref{dfn:wh-g1},  except for replacing $\gak$ by $\gakm$.
It then follows from our derivation of \eqref{eq:step2}-\eqref{eq:g3-formula} that 
\begin{equation}\label{Pbm-hg3}
\Pm_{\bb,\sfKo} = (1+o_{\sfKfv,\vep}(1)) 
\big[\vep e^{\sfKfv}/r_{\sfKo}\big]^2
 \,  \alpha(\sfKfv,\sfKp) 
\bE \big[  \whgo(\btet); (y_1,L_{\btet}) \in \Istar  \,| \FF \big] \,.
\end{equation}
Setting analogously to \eqref{dfn:wh-g1}-\eqref{eq:BK12},
\begin{align}\label{dfn:wh-gl}
\whgl(\btet)& :=  [r_{\sfKo}/r_{\sfKth}]^2 \, \bE \big[ \gth \big( \wh{\eta}(\sfKth) \big) ;  \Bll  \,| \btet \big]
\quad \text{for} \quad \btet=\btet(y_1) \quad \text{and} \\
\Bll &:=  \!\!\!\! \bigcap_{j \in (\sfKo,\ell-3)} \!\!\!\! \big\{ \wh{\eta}(j) > j^{\gakm}  \big\} \!\!\!\!\! \bigcap_{j \in \{\ell-3,\ell\}}  
\!\!\!\!\! \big\{ \wh{\eta}(j) \in \sfI_j \big\} \!\!\!\! 
\bigcap_{j \in (\ell,\sfKt)} \!\!\! 
\big\{ \wh{\eta}(j) > j^{\gakm}  \big\} \bigcap \big\{ \wh{\eta}(\sfKt) \in \sfIs_2 \big\} \,,
\label{dfn:B-ell}
\end{align}
we enlarge $P_{\bb,\ell}$ by replacing at each $j \in (\sfKo,\ell-3) \cup (\ell,\sfKt) \cup (\sfKth,n')$,  the window
$\sfI_j$ by $(j_{n_\vep}^{\gakm},\infty)$,  to arrive,  analogously to \eqref{Pbm-hg3} at 
\begin{equation}\label{Pbl-hg3}
P_{\bb,\ell} \le 2 \big[\vep e^{\sfKfv}/r_{\sfKo}\big]^2 \,  \alpha(\sfKfv,\sfKp) 
\bE \big[  \whgl(\btet); (y_1,L_{\btet}) \in \Istar  \,| \FF \big] \,.
\end{equation}
Comparing \eqref{Pbm-hg3} with \eqref{Pbl-hg3},   it suffices for \eqref{eq:drop2-levels} to show 
that for some $C<\infty$,  
\begin{align*}
\whgl(\btet) \le C\,  
 \whgo(\btet)  \qquad  \qquad \forall (y_1,L_{\btet}) \in \Istar 
\end{align*}
and,  in view of \eqref{eq:K1-to-K3-u} and \eqref{eq:V-tet-alt},  for this we merely need to show that 
\begin{align}\label{eq:dr2-ell}
\whgl(\btet) \le C\,  e^{-\sfb} e^{-\sfco V_{\btet}} \,  [r_{\sfKo}/r_{\sfKth}]^2 \, 
\bar \bE_{y_1}^{(\sfKo)} \big[  \gth\big( \wh{\eta}(\sfKth) \big) ;  \BKo  \big]
\qquad  \forall (y_1,L_{\btet}) \in \Istar \,.
\end{align}
To this end,  we utilize our coupling with $\bmu$,  in analogy with \eqref{eq:BK12n}
to set $\wtBlp$ as in \eqref{dfn:B-ell} 
except for involving only the contributions $\whes(j)$ from the coupled 
$\sfKo$-level excursions (with $B_q=0$),  
while also slightly decreasing $\gakm$,  $\gak$ and $\kapo$ (so $\wtBKo=\wtBp$ of \eqref{eq:BK12n}),  
and consider the following proxy for $\whgl(\cdot)$,
\begin{align}\label{dfn:wt-g1n}
\wtglp (\btet) &:= [r_{\sfKo}/r_{\sfKth}]^2\,  
\, \bE \big[  e^{-\frac{\sfco}{\sfm_\vep} \DelKth} \gth \big( \whes(\sfKth) \big) ;  \wtBlp \big], 
\end{align}
(which at $\ell=\sfKo$ matches our earlier definition \eqref{dfn:wh-g1n}).  Now 
note that a minor change in \Cref{sec:coupling}
(of deriving \eqref{eq:bd-gopm} in terms of events $\Bll$ and $\wtBlp$ of \eqref{dfn:B-ell}
instead of those of \eqref{dfn:Bplus}),  
yields,  similarly to \eqref{eq:coupling-bdt}-\eqref{eq:coupling-bd1},  
\[
\whgl(\btet) \le (1+o_{\vep}(1))  \wtglp(\btet) + o_\vep(1) = (1+o_{\vep}(1)) \wtglp(\bmu)\, e^{\sfco V_{\bmu}}  
e^{-\sfco V_{\btet}} + o_\vep(1)  \quad  \forall (y_1,L_{\btet}) \in \Istar \,.
\]
Similarly to our derivation of \eqref{eq:V-tet-alt},  this,  
together with \eqref{eq:lim-V-psi}-\eqref{eq:lim-V-psi2},  yield that for ordered $\bpsi$ consisting 
of $
\sfs_{\sfKo}(y_1)$
\abbr{iid} samples from $\bar \mu$,   
\begin{equation}\label{eq:g-ell-bd}
\whgl(\btet) \le (1+o_{\sfKo,\vep}(1)) e^{-\sfb} e^{-\sfco V_{\btet}} 
\bE \big[ \wtglp(\bpsi)  {\one}_{\{(y_1,L_{\bpsi}) \in \Istar\}} \big]
+ o_\vep(1) 
\qquad  \forall (y_1,L_{\btet}) \in \Istar \,.
\end{equation}
Further,  denoting by $\Bsl$ the event $\Bll$ of \eqref{dfn:B-ell} but with  
the smaller parameters $\gakm$,  $\gak$,  $\kapo$ of $\wtBlp$ and with slightly larger 
values of $\kapp$,  $\okappa$,  our proof of the left-inequality of \eqref{eq:bd-gopm} can be adapted to show that 
\begin{equation}\label{eq:coupling-bd2}
\bE \big[ \wtglp(\bpsi)  {\one}_{\{(y_1,L_{\bpsi}) \in \Istar\}} \big] \le 2 \bE \big[ \whgl(\bpsi) \big] 
= 2  \,  [r_{\sfKo}/r_{\sfKth}]^2 \,  \bar \bE_{y_1}^{(\sfKo)} \big[  \gth\big( \wh{\eta}(\sfKth) \big) ;  \Bsl \big] \,,
\end{equation}
where in $\gth(\cdot)$ on the right hand side we have the same $\kapo$,  $\kapp$ values as in $\Bsl$.  In view of
\eqref{eq:g-ell-bd}-\eqref{eq:coupling-bd2},  the bound \eqref{eq:dr2-ell} is
a direct consequence of having for some $C<\infty$ and all $y_1 \in \sfIs_1$, 
\begin{equation}\label{eq:dr2-ell-alt}
\bar \bE_{y_1}^{(\sfKo)} \big[  \gth\big( \wh{\eta}(\sfKth) \big) ;  \Bsl \big] 
\le C\, 
\bar \bE_{y_1}^{(\sfKo)} \big[  \gth\big( \wh{\eta}(\sfKth) \big) ;  \BKo  \big]  \,.
\end{equation}
Now,  by the same reasoning as in part (d),  the bound \eqref{eq:bd-drop-ell} applies also for
\[
\gfl(y') = 
 \bar{\bE}^{(\ell)}_{y'} \Big[ \gth \big(\wh{\eta}(\sfKth)\big);
\{ \wh{\eta} (\sfKt) \in \sfIs_2 \} 
\!\! \bigcap_{j \in (\ell,\sfKt)} \!\! \{ \wh{\eta} (j) >  j_{n_\vep}^{\gakm} \} \Big] \,,
\]
yielding \eqref{eq:dr2-ell-alt} (as argued at
\eqref{eq:K1-to-K3n},  for 
$\sfKo \gg 1$ the right hand side 
of \eqref{eq:dr2-ell-alt} 
is independent of $\gakm$,  $\gak$,  $\kapo$, $\kapp$,  $\okappa$).
\end{proof}

Recalling that $\sfKw = (21 \log n_\vep)^{1/\gak}$,  we proceed to prove \Cref{prop:UpperBoundz},  
where similarly to \cite[Section 4]{BRZ2} we consider separately,  \emph{early branching},  namely 
\eqref{eq:UpperBoundzr} with $\ell \le \sfKw$,  \emph{branching in the bulk} where $\ell \in (\sfKw,n_\vep-\sfKz+2]$
and \emph{terminal branching},  namely \eqref{bd:Late-brch}. 
 
\subsection{Early branching}\label{sec-earlybrualt} 
Note that $\bb \in \bG_\ell$ implies that $\bB_\sfd(\sfv_{\bb},r_{\ell-1})$ is
inside the annulus  $\bB_\sfd(\sfv_{\bb'},r_{\ell-3}) \setminus \bB_\sfd(\sfv_{\bb'},r_{\ell-1})$. 
Thus, when $\ell \le \sfKo$ the events 
$\UUpp(\Ybpz)$ and $\UUpy(\Ybz)$
are independent given the ordered collections $\underline{\sfx}_\bb$ 
and $\underline{\sfx}_{\bb'}$ of start point and endpoint of the relevant 
$\sfKo$-level excursions of these two blocks.  Recall from \eqref{UUpy-asymp} and \eqref{eq:g-def-alt} that 
\[
\sup_{\bb,\underline{\sfx}_{\bb}} \{ \bP(\UUpy(\Ybz) | \underline{\sfx}_\bb) \} \le (1+o_{\sfKfv,\sfKz,\vep}(1)) 
e^{-2(n''-\sfKo)}
\bar \alpha(\sfKp)  \| \gst \|_\infty \,.
\] 
Thus,  for large $\sfKfv, \sfKz$ any $\vep \le \vep_0(\sfKfv,\sfKz)$ and $\bb \in G_\ell$ for some $\ell \le \sfKo$, 
\begin{align*}
P^{(2)}_{\bb,\bb'} & \le P^{(1)}_{\bb'} \sup_{\bb, \underline{\sfx}_{\bb}} \{
\Pbm(\UUpy | \underline{\sfx}_{\bb}) \}  
\le 2 e^{-2(n''-\sfKo)} \bar \alpha(\sfKp) \| \gst \|_\infty P^{(1)}_{\bb'} \,,
\end{align*}
so recalling 
that $\|\gst\|_\infty \le c e^{-\sfKo^{\kapo}}$,  we thus establish 
\eqref{eq:UpperBoundzr} for all $\ell \le \sfKo$. 

Turning to $\ell \in (\sfKo,\sfKw]$,  recall $\sfKop=\sfKo$ unless $\ell \in \{\sfKo+1,\sfKo+2\}$ in which case
$\sfKop=\sfKo+3$.  Setting $\ell':=\ell \vee (\sfKo+3)$ we keep from the event $\UUpy(\Ybz)$ 
only the angle control for the $\ell'$-level excursions and the barrier information starting  at that level,
while replacing the event $\UUpp(\Ybpz)$ by the larger event $\UUplp(\Ybpp)$ of \eqref{dfn:UUpl}.
Using the short-hand $y=\YbKp$,  $y':=\Ybpp$ with the corresponding $L^{(y)}_{\bb,j}$,  $L^{(y')}_{\bb',j}$ for 
$j \ge \sfKop$ as in \eqref{dfn:emp-ang},  except for taking the $j$-level 
excursions during the $\sfs_{\sfKop}(\cdot)$ first $\sfKop$-level excursions
(see \eqref{eq:nominal} for $\sfs_j(\cdot)$),
we deduce by the argument leading to \eqref{eq:pre-precise},  that 
\begin{align}\label{early-branch}
P^{(2)}_{\bb,\bb'} \le  & (1+o_{\sfKfv}(1)) 
\sum_{y'_4, y'_5} \alpha_{\sfKfv,\sfKp}(y'_5) 
\sum_{y_4, y_5} \alpha_{\sfKfv,\sfKp}(y_5) 
 \bP\big(  \whBBlp
 (y,y_4,y_5) 
\cap \BBpl (y',
y'_4,y'_5) \, \big|   \FF \big),
\end{align}
where $y,y' \in \sfI_{\sfKo'}$, with $\BBpl(\cdot)$ as in \eqref{def:B-no-ell} and 
\begin{align}\label{def:whBBl}
 \whBBl (y,y_4,y_5) &:=\{ 
L^{(y)}_{\bb,\ell} \in \bB^\infty_{\bar \mu} (\boetasl) \} 
 \!\!\!\! \bigcap_{j \in \sfJupl \cap [\ell,n'']} \!\!\!\! \{ \whetajy \in \sfI_j \} \bigcap
 \{ \whetaynp = y_4, \widehat{\eta}_{\bb}^\sharp(n'')  = y_5\}
\end{align}
while for all other values of $y,y'$ these events are empty
(compare with \eqref{eq:barrier-event-z} and \eqref{eq:was-event}).

We further claim that for some finite $c(\sfKp)$ and any $\ell \in (\sfKo,\sfKw]$,
\begin{align}
  \label{eq-decouplinearly}
 \sum_{y_4,y_5}  \alpha_{\sfKfv,\sfKp} (y_5) \bP(
 \whBBlp (y,y_4,y_5)
\cap & \BBpl (y',y'_4,y'_5)  | \FF) \nn \\ 
\leq  c(\sfKp) e^{-2(n''-\ell)-\ell^{\gak}}
 \bP(& \BBpl (y',y'_4,y'_5) | \FF ).  
\end{align}
Indeed,  note that for $\ell \ge \sfKo+3$ the event $\BBpl (y',y'_4,y'_5)$ is measurable 
on the $\sigma$-algebra $\GG^{\sfv_\bb}_{\ell'}$ since the path within
$\bB_\sfd(\sfv_{\bb},r_{\ell'-1}) \subset \bB_\sfd(\sfv_{\bb'},r_{\ell-3}) \setminus \bB_\sfd(\sfv_{\bb'},r_{\ell-1})$
impacts only $\bb'$ block excursions at levels $\ell-2$ and $\ell-1$.  For $\ell \in \{ \sfKo+1,\sfKo+2\}$ the
path within $\bB_\sfd(\sfv_{\bb},r_{\ell'-1})$ may impact the $\sfKo$ level excursions of block $\bb'$,  but 
by measuring for such $\ell$ both the angles and excursion counts within the $(\sfKo+3)$-level 
excursions of $\bb$ and $\bb'$ we retain this measurability property.  Now,  as the block $\bb$ is contained 
in $\bA_u^o$, also $\FF \subseteq \GG^{\sfv_\bb}_{\ell'}$, 
with \eqref{eq-decouplinearly} a direct consequence of the following lemma.
\begin{lemma}\label{lem:key-bd}
For some finite $c(\sfKp)$, uniformly over quantifiers as above, 
\begin{equation}
\label{eq-EB1}
\sum_{y_4,y_5} \alpha_{\sfKfv,\sfKp} (y_5)
\bP( \whBBlp (y,y_4,y_5)|\GG^{\sfv_\bb}_{\ell'}) \leq c(\sfKp) 
e^{-2(n''-\ell)} e^{-\frac{1}{7} \ell^{\gak}} \, .
\end{equation}
\end{lemma}
\begin{proof} For notational convenience we use hereafter $\ell \ge \sfKo+3$ instead of $\ell'$.  
Note that for $\sfv:=\sfv_\bb$,  the collection 
$\underline{\sfx}_{\bb}$ of start point and endpoint of all excursions from $\partial \bB_{\sfd}(\sfv,r_{\ell})$ 
to $\partial\bB_{\sfd}(\sfv,r_{\ell-1})$ is $\GG^{\sfv}_{\ell}$-measurable,  and with 
$\FF \subseteq \GG^\sfv_{\ell}$, the same applies for $y=\YbKp$ (as $\sfKop \le \ell$),   
$\widehat{\eta}^{(y)}_{\bb,\ell}$ and $L^{(y)}_{\bb,\ell}$.  It thus suffices to show that fixing 
any $a \in [\ell^{\gak},\ell^{\okappa}]$ such that $\sfs:=\sfs_{\ell}(a)$ of \eqref{eq:nominal} is integer,
the bound \eqref{eq-EB1} holds 
when replacing $\sfs_{\sfKop}(y)$ by the first $\sfs$ 
of the $\ell$-level excursions between concentric circles centered at $\sfv$.  Having done that,  
the conditional probability on the \abbr{lhs} of \eqref{eq-EB1} now only depends 
on the angular differences $\btet = \{\thbl, q \le \sfs\}$ between the start and end 
point of those $\sfs$ 
excursions.  We can further replace 
$\whetajy$ by $\wh{\eta}_{\bb,j}^{(a)}$ that correspond to our $\sfs$ excursions 
at level $\ell$ and for $y_i$ with $\sfs_{\sfK_i}(y_i) \in \bN$,  denote
by $\BBKth (y_4,y_5)$ and $\wCe (y_2,y_3)$ the events analogous to   
those of \eqref{eq:B-K3} and \eqref{eq:hat-C},  respectively.  In this setting 
\begin{align}\label{eq:B-U}
\whBBl 
(y,y_4,y_5)\subset 
& \{L_{\btet} \in \bB^\infty_{\bar \mu} (\boetasl)\} \bigcap \Bla 
\bigcap
\BBKth (y_4,y_5) 
 \,,
\end{align}
for the $\GG^{\sfv}_{\sfKt}$-measurable events
\begin{align*}
\Bla : =& 
 \!\!\!\! \bigcap_{j \in [\ell,\sfKt]} \!\!\!\! \{ \wh{\eta}^{(a)}_{\bb,j} \in \sfI_j \} 
\end{align*}
and we get \eqref{eq-EB1} upon showing that for any such $a$ and $L_{\btet} \in \bB^\infty_{\bar \mu} (\boetasl)$,  
\begin{equation}\label{eq-EBa}
\sum_{y_4,y_5} \alpha_{\sfKfv,\sfKp} (y_5)
\bP\big( \Bla \cap \BBKth (y_4,y_5) \,|\, \GG^{\sfv}_{\ell}\big) 
\leq c(\sfKp) 
e^{-2(n''-\ell)} e^{-\frac{1}{7} \ell^{\gak}} \, .
\end{equation}
Proceeding to show \eqref{eq-EBa},  with $\sfKth > \sfKt > \sfKw \ge \ell$,
applying Lemma \ref{recursionends}
allows us to decouple 
the events $\BBKth (y_4,y_5)$ (about concentric excursion counts starting at level $\sfKth$),
from $\GG^\sfv_{\sfKt}$-measurable events.
Specifically,  
we have as in the
derivation of \eqref{eq:g3},  that for some finite  $c'$, $c(\sfKp)$ and any $y_i \in \sfIs_i$, $i=2,3$, 
\begin{align}\label{eq:barrier}
& \sum_{y_4,y_5} \alpha_{\sfKfv,\sfKp} (y_5)
\bP \Big(\BBKth(y_4,y_5) \big| \, \wCe (y_2,y_3), \GG^{\sfv}_{\sfKt} \Big) \nn \\
& \le c' \sum_{y_4,y_5} \alpha_{\sfKfv,\sfKp} (y_5)
  \, \bar \bP_{y_2}^{(\sfKt)} \Big( \BBKth (y_4,y_5) \big| \, \wCe (y_2,y_3) \Big)
\le c(\sfKp) y_3 e^{-2y_3} e^{-2(n''-\sfKth)} \,,
\end{align}
where the last inequality follows as in the derivation of \eqref{eq:g3-formula} out of \eqref{eq:alt-g3}
(upon recalling that $\widehat{\eta}(\cdot)$ is distributed under $\bar{\bP}_{y_2}^{(\sfKt)}$ as the 
excess occupation measures for reflected \abbr{srw} on $[\sfKth-1,n'']$,  during its first 
$\sfs_{\sfKth}(y_3)$ moves from $\sfKth$ to $\sfKth-1$).  We further condition on  $\GG^{\sfv}_{\sfKt} \supset 
\GG^{\sfv}_\ell$
in the \abbr{lhs} of \eqref{eq-EBa},  splitting according to the value $y_2 \in \sfIs_2$ of $\wh{\eta}^{(a)}_{\bb,\sfKt}$,  
at which point we can combine \eqref{eq:barrier} with the elementary identity
\[
\bP \Big( \BBKth (y_4,y_5) | \GG^{\sfv}_{\sfKt} \Big) = \sum_{y_3 \in \sfIs_3} 
\bP (  \wCe (y_2,y_3) | \GG^{\sfv}_{\sfKt} ) \bP (
\BBKth (y_4,y_5) | \wCe (y_2,y_3),  \GG^{\sfv}_{\sfKt}), 
\]
to deduce that \eqref{eq-EBa} is a consequence of having that for any possible $a \in [\ell^{\gak},\ell^{\okappa}]$,  
\begin{equation}\label{eq-EB1-alt}
\bE \Big[ \gth (\wh{\eta}^{(a)}_{\bb,\sfKth}) ; \Bla
\, |\, \GG^{\sfv}_{\ell} \Big] {\one}_{\bB^\infty_{\bar \mu} (\boetasl)} (L_{\btet})
\leq  
e^{-2(\sfKth-\ell)} e^{-\frac{1}{7} \ell^{\gak}} \, .
\end{equation}
Our proof of \eqref{eq-EB1-alt} is an adaptation of the arguments in 
\cite[Section 4.4]{BRZ2}
and we only sketch the required changes.  
Specifically,  we recall that the \abbr{lhs} of \eqref{eq-EB1-alt} 
only depends on the angular differences $\btet$ between the start and end point of our $\sfs$ excursions, 
coupling those excursions,  as in \cite[Lemma 4.10]{BRZ2},  with $\sfs$ such excursions of \abbr{iid} angular 
differences $\bpsi = \{\psi_q, q \le \sfs\}$ from $\bar \mu$ of \eqref{dfn:Poiss}.  Hereafter we denote by $\wt{\bP}$
the law on the probability space enlarged with these coupled \abbr{iid} angles
(analogously to \cite[(4.110)]{BRZ2},  except that here we  also average over the 
independent realization of $\bpsi$).  Now,  in view of \eqref{dfn:dw},  the event 
$\{L_{\btet} \in \bB^\infty_{\bar \mu} (\boetasl)\}$ is contained (up to changing
certain universal constants),  in the $\GG_\ell^\sfv$-measurable event 
$\WW_{\sfv,\ell}(\sfs)$ of \cite[(4.19)]{BRZ2},  and the ``good coupling'' event $\Phi_{\ell,a}$ of \cite[(4.112)]{BRZ2},   
is such that uniformly over $\ell, a$ and $\sfv$, $\vep$,
\begin{align}\label{eq-EB2}
\widetilde \bP(  \Phi_{\ell,a}^{\complement}   | \GG^{\sfv}_{\ell},   \WW_{\sfv,\ell}(\sfs) ) & \leq e^{-10 n_\vep}
\end{align}
(see \cite[(4.108) and proof of Lemma 4.11]{BRZ2}).  
Similarly to  \cite[Lemma 4.9]{BRZ2}, 
utilizing the proof of \cite[Lemma 4.13]{BRZ2} at
$\sfKth-\sfKt \le \sfKt - \ell$,  yields that in terms of $G_{\ell,a,y_2,y_3}(m)$ of \eqref{eq:BRZ-471}
(which replaces the factor
 $F_{\ell,a,y_2,y_3}$ of \cite[(4.94)]{BRZ2}),  
\[
\wt{\bP} \Big(\Bla 
\cap \Phi_{\ell,a} \bigcap_{i=2,3} \{ \wh{\eta}^{(a)}_{\bb,\sfK_i}=y_i \}    
\, |\, \GG^{\sfv}_{\ell} \Big) \le G_{\ell,a,y_2,y_3}(m)\,, 
\]
where $m=\log \ell$ 
due to $
\Phi_{\ell,a}$.
Similarly to \eqref{eq:BRZ-472},  summing this bound over $y_2,y_3$ we arrive at 
\begin{align}\label{eq-EB3}
& \widetilde \bE\Big[\gth \big(\wh{\eta}^{(a)}_{\bb,\sfKth}\big) ;
\Bla 
\cap \Phi_{\ell,a}| \GG^{\sfv}_\ell \Big] 
\le a \big(1+ \frac{a^2}{\sfKt} \big) e^{-2a} e^{-2 (\sfKth-\ell)} e^{M \log \ell} 
\end{align}
for some universal $M$.  Now recall that $\gth(\cdot) \le 1$ and utilize \eqref{eq-EB2} to eliminate
the event $\Phi_{\ell,a}$ from \eqref{eq-EB3} in case $L_{\btet} \in \bB^\infty_{\bar \mu} (\boetasl)$.  
Thereby,  we can also change from $\wt{\bE}[\cdot]$ to $\bE[\cdot]$ in \eqref{eq-EB3}
and since $a \ge \ell^\gak \gg 1$,  conclude that 
\eqref{eq-EB1-alt} holds.
\end{proof}
Equipped with \Cref{lem:key-bd},  we combine \eqref{early-branch} with \eqref{eq-decouplinearly} to get
in terms of $P_{\bb',\ell}$ of \eqref{dfn:Pml},  that 
\begin{align*}
P^{(2)}_{\bb,\bb'} &\le  2 c(\sfKp) P_{\bb',\ell} \,  e^{-2(n''-\ell)} e^{-\ell^{\gak}}  
\end{align*}
whereby 
from  \Cref{lem:2-levels} 
we have that for some $c=c(\sfKfv,\sfKp)$ finite,  
\[
P^{(2)}_{\bb,\bb'} \le  c  \,  \Pm_{\bb',\sfKop} e^{-2(n''-\ell)} e^{-\ell^{\gak}}  \,,
\]
as claimed in \eqref{eq:UpperBoundzr}.

\subsection{Branching in the bulk}\label{sec-brbulku} 
In analogy with \eqref{eq:hat-C},  given $\wh{\eta}^{\,[z]\,}_{\bb,j}=a$ we set for $j' \ge j \ge \sfKo$ 
the events
\begin{equation}\label{def:Cjjp}
\wCe_{j,j'} (a,a') :=\{\wh{\eta}_{\bb,j'}^{(a)}=a'\}
\end{equation}
(namely,  that precisely $\sfs_{j'}(a')$ of the $j'$-level excursions around $\sfv:=\sfv_\bb$ have been completed 
during the first $\sfs_j(a)$ such level $j$ excursions),  so $\wCe_{j,j} (a,a)$ always holds 
and $\bigcup_a  \wCe_{j,j'} (a,a')$ amounts to $\wh{\eta}^{\,[z]\,}_{\bb,j'}=a'$.
Next,  note that for $\ell \in (\sfKw,n_\vep-\sfKz+2]$,  
the mapping $\ell \mapsto f(\ell):=\ell+\lceil 5 \log (n_\vep-\ell) \rceil$ is increasing 
with $f(\sfKt)=\sfKth$ (see \Cref{rem:more-int-effects}),  and $f(\ell) <  n'$ (since $\sfKz \gg \sfKf$).  Setting 
\begin{equation}\label{dfn:Ktp-Kthp}
(\sfKt',\sfKth') := \left\{\!\!\!\begin{array}{cll}
& (f(\ell),\sfKth),  & \qquad f(\ell) \in [f(\sfKw), \sfKt] \,,\\
& (\sfKt,\sfKth),  & \qquad  f(\ell) \in (\sfKt,\sfKth],   \\
& (\sfKth,f(\sfKth)), & \qquad \ell \in (\sfKt,\sfKth),  \\
& (\ell,f(\ell)), & \qquad \ell \in [\sfKth,n_\vep-\sfKz+2], 
\end{array}
\right.
\end{equation}
we have that $\sfKt' \ge \ell$ and $\sfKth' \ge f(\sfKt') \vee \sfKth$ are outside $(\sfKt,\sfKth)$.
 Aiming to show \eqref{eq:UpperBoundzr} for such $\ell$,  start as before with \eqref{early-branch}  
at $\ell'=\ell$ and $y,y' \in \sfI_{\sfKo}$,  replacing \eqref{def:whBBl} by the following 
larger event (which contains even the right hand side of \eqref{eq:B-U}),
\begin{align}\label{eq:B-U-ell}
 \whBBl (y,y_4,y_5) &:= 
\!\!\bigcup_{a,y_2,y_3}  
\!\!
\Big\{ \wCe_{\ell,\sfKt'} (a,y_2) \cap \wCe_{\sfKt',\sfKth'} (y_2,y_3) \cap \BBKthp (y_4,y_5)\Big\} \,,
\end{align} 
where we restrict to $y_2 \in \sfI_{\sfKt'}$,  $y_3 \in \sfI_{\sfKth'}$ and if $\ell \notin (\sfKt,\sfKth)$ then also 
to $a \in \sfI_{\ell}$.  
We show in the sequel that for some $c(\sfKfv,\sfKp)$ finite, 
all $\ell \in  (\sfKw,n_\vep-\sfKz+2]$ and relevant 
$y,y_4,y_5$,
\begin{equation}\label{eq:bulk-bd}
\bP( \whBBl
(y,y_4,y_5)\,|\,\GG^{\sfv}_{\ell}) \leq  
c(\sfKfv,\sfKp)
\, e^{-2(n'-\ell)} e^{- \frac{1}{6} \ell_{n_\vep}^{\gak}} \,.
\end{equation}
Since $\ell_{n_\vep} \ge \sfKz$ and $\sfKz \gg \sfKfv$,  upon recalling that $2(n''-n') = 10 \log \sfKfv$
we deduce from \eqref{eq:bulk-bd} that the analog of \Cref{lem:key-bd} holds here,  namely that 
uniformly over $\ell$ in our range, 
\begin{equation*}
\sum_{y_4,y_5} \alpha_{\sfKfv,\sfKp} (y_5)
\bP( \whBBl (y,y_4,y_5)|\GG^{\sfv}_{\ell}) \leq  
e^{-2(n''-\ell)} e^{-\frac{1}{7} \ell_{n_\vep}^{\gak}} \, .
\end{equation*}
From the latter bound we get \eqref{eq:UpperBoundzr} by applying \Cref{lem:2-levels},  
along the same path as in the end of \Cref{sec-earlybrualt}.

Turning to prove \eqref{eq:bulk-bd},  having $\sfKth' \ge f(\sfKt')$ and $y_2=\wh{\eta}^{[z]}_{\bb,\sfKt'}$ 
measurable on $\GG_{\sfKt'}^\sfv$,  yields 
for some $c'$, $c(\sfKfv)$ finite,  the following analog of \eqref{eq:barrier}, uniformly in 
$y_2 \in \sfI_{\sfKt'}$ and $y_3 \in \sfI_{\sfKth'}$, 
\begin{align}\label{eq:barrier-Ktp}
 \bP \Big(\BBKthp (y_4,y_5) \big| \, \wCe_{\sfKt',\sfKth'} (y_2,y_3), \GG^{\sfv}_{\sfKt'} \Big)
 &\le c' \, \bar \bP_{y_2}^{(\sfKt')} \Big( \BBKthp(y_4,y_5) \big| \, \wCe_{\sfKt',\sfKth'} (y_2,y_3) \Big) \nn \\
 & \le c(\sfKfv) y_3 e^{-2y_3} e^{-2(n'-\sfKth')} \,,
\end{align}
where the second inequality of \eqref{eq:barrier-Ktp}  is derived as in \cite[(8.3)]{BRZ2} 
after omitting the constraint at level $n''$ (taking there $L=n_\vep-\sfKth'$, $k=n'-\sfKth'$ 
and $z=y_3$; see also \blue{\eqref{18.20bendj1}}).
Setting $g_j(y):=ye^{-2y} {\bf 1}_{\sfI_j}(y)$,  
we  
utilize \eqref{def:Cjjp},  \eqref{eq:B-U-ell} and \eqref{eq:barrier-Ktp},  to deduce that 
\begin{align*}
\bP( 
\whBBl
 (y,y_4,y_5)\,|\,\GG^{\sfv}_{\ell}) &\le 
 c(\sfKfv) e^{-2(n'-\sfKth')} \bE\Big[ \hg_{\sfKth'} (\wh{\eta}^{\,[z]\,}_{\bb,\sfKth'})\, |\,  \GG^{\sfv}_{\ell} \Big]
 \,.
\end{align*}
In view of \eqref{dfn:windows},  the latter bound yields \eqref{eq:bulk-bd} provided that  
\begin{equation}\label{eq:bd-Kthp-ell}
(\sfKth')_{n_{\vep}}^{\gak} - \log (\sfKth')_{n_{\vep}} \ge \sfKth'-\ell +  \frac{1}{12} \ell_{n_\vep}^{\gak} \,.
\end{equation}
With $(\sfKth')_{n_\vep} \gg 1$ and $\frac{5}{6} \sfKth^{\gak} \ge 10 \log n_\vep = 2 (\sfKth-\sfKt)$ 
(see \Cref{rem:more-int-effects}),  it is not hard to verify from \eqref{dfn:Ktp-Kthp}
that \eqref{eq:bd-Kthp-ell},  and thereby \eqref{eq:bulk-bd},  hold whenever $\sfKt' \ge \sfKt$.  
Turning to the remaining case,  where $\sfKt' = f(\ell) < \sfKt$ and $\sfKth'=\sfKth$,  setting 
\begin{equation}\label{dfn:BKtp}
\BB^{(\vep)}_{[\sfKt',n'']}(y_2,y_4,y_5) := \bigcup_{y_3 \in \sfIs_3} \big\{ \wCe_{\sfKt',\sfKth} (y_2,y_3) \cap \BBKth(y_4,y_5) \big\}\,,
\end{equation}
we deduce similarly to \eqref{eq:barrier-Ktp} that 
for some finite $c'$,  $c(\sfKfv)$,  uniformly in $a \in \sfI_\ell$ and $y_2 \in \sfI_{\sfKt'}$,  
\begin{align}\label{eq:barrier-ell}
 \bP \Big( \BB^{(\vep)}_{[\sfKt',n'']}(y_2,y_4,y_5) \big| \, \wCe_{\ell,\sfKt'} (a,y_2), \GG^{\sfv}_{\ell} \Big)
 &\le c' \, \bar \bP_{a}^{(\ell)} \Big( \BB^{(\vep)}_{[\sfKt',n'']}(y_2,y_4,y_5) \big| \, \wCe_{\ell,\sfKt'} (a,y_2) \Big) \nn \\
 & \le c(\sfKfv) e^{-2(n'-\sfKth)} \bar \bE_a^{(\ell)} \Big[ \hg_{\sfKth} \big( \wh{\eta} (\sfKth) \big) \Big]
\end{align}
(by the Markov property of $j \mapsto  \wh{\eta}(j)$).  As $a=\wh{\eta}^{\, [z]}_{\bb,\ell}$ is
$\GG^{\sfv}_{\ell}$-measurable,  from \eqref{eq:B-U-ell} with \eqref{dfn:BKtp} one has that 
\begin{align*} 
\bP( \whBBl (y,y_4,y_5)\,|\,\GG^{\sfv}_{\ell}) 
&=  {\one}_{a \in \sfI_\ell} \bP\big(\!\! \bigcup_{y_2 \in \sfI_{\sfKt'}}\!\!
 \wCe_{\ell,\sfKt'} (a,y_2) \cap \BB^{(\vep)}_{[\sfKt',n'']} (y_2,y_4,y_5) \,|\,\GG^{\sfv}_{\ell} \big) \\
 &\le \sup_{a \in \sfI_\ell,y_2 \in \sfI_{\sfKt'}} \bP \big( \BB^{(\vep)}_{[\sfKt',n'']}(y_2,y_4,y_5) \big| \, \wCe_{\ell,\sfKt'} (a,y_2), \GG^{\sfv}_{\ell} \big) \,,
\end{align*}
which in view of \eqref{eq:barrier-ell} yields \eqref{eq:bulk-bd} upon showing that 
for some $C<\infty$ and all $\ell \in [\frac{1}{3} \sfKt,\sfKt]$,
\begin{align}\label{eq:rough-C23}
\sup_{a \in \sfI_\ell} \,  \bar \bE_{a}^{(\ell)} \Big[ \hg_{\sfKth} \big( \wh{\eta}(\sfKth) \big) \Big] \le 4 \,  \sfKth \, 
 e^{-2(\sfKth-\ell)}e^{-2 \ell^\gak} \,.
\end{align}
To this end,  first express $\wh{\eta}(\sfKth)$,  as in \eqref{eq:pre-clt},
in terms of $a$ and the normalized sum $\wh{S}_m$ of $m=\sfs_{\ell}(a)
$ 
\abbr{IID} variables $\{\GJ'_i\}$ of the form \eqref{dfn:GJp}.
Here $L=\log (r_{\ell}/r_{\sfKth}) \ge \sfKth-\ell$ may exceed $\sqrt{\sfKth}$,  so
\eqref{eq:tilted} and hence \eqref{eq:C23} do not apply.
But,  $\ell \ge \frac{1}{4} \sfKth$ and \eqref{eq:Gm-1} holds,  so in \eqref{eq:pre-clt}
we can still replace 
$\sfrho_\vep$ by $2$ and $\Gamma_m$ by $1$.  Having done so,  note that 
$\hg_{\sfKth}(y) \le \sfKth e^{-2y}$,  yielding,  as claimed,  that 
\[
\sup_{a \in \sfI_\ell} \,  \bar \bE_{a}^{(\ell)} \Big[ \hg_{\sfKth} \big( \wh{\eta}(\sfKth) \big) \Big] \le 
2 \sfKth \sup_{a \ge \ell^\gak} \big\{  e^{-2 a} e^{-4L}  \bE \big[ e^{- 2  \wh{S}_m} \big]  \big\} 
\le 4 \, \sfKth \,  e^{-2L} e^{-2 \ell^\gak}
\] 
(where the second inequality is due to \eqref{eq:C23-mgf}).

\subsection{Terminal branching}\label{sec-latebru}
Whereas in early and bulk branching we could completely ignore the dependence on $\sfKfv$, 
it matters in case of terminal branching,  where $\sfd\left(\sfv_{\bb},\sfv_{\bb'}\right)$ can be  
as close as the block size.  To deal with what happens within each block,  recall \eqref{dfn:UU-alt}
of \Cref{rem:alt-UU},  that under the event $\UUz$, 
the first $\sfk_\vep(z)$ inter-city excursions do not visit some ball $\bB_\sfd(\sfx,\frac{\vep}{e})$ 
centered at $\sfx$ in our $\widetilde{c} \, \vep$-net $\bF^\bb$ of block $\bb$.  
Thus,  replacing for $\sfx \in \bF^\bb$ the event $\RUy$ inside the block (see \eqref{eq:barrier-event}),  by 
not visiting $\bB_\sfd(\sfx,{\vep}/{e})$,  we set
in lieu of \eqref{eq:good-cover-y},  the events 
\begin{align}\label{dfn:UUpv}
\UUx (y) &:=  \{ \etaxny = 0 \} \bigcap \BBy \cap \WWy
\end{align}
noting that
\begin{equation}\label{eq:Ub-Uxb}
\UUpy(y) \subset \bigcup_{\sfx \in \bF^\bb} \UUx(y) \,.
\end{equation}
In preparation for taking such unions over pairs of centers,  fix $\sfx \in \bF^{\bb}$ and consider the subsets 
\begin{align}\label{dfn:tGj}
\wbG_{j}  &:=  \left\{ \sfx'  \in \bF_{n_\vep} :  \sfd(\sfx,\sfx') \in [e^{j-1} \vep,  e^j \vep) \right\},
\end{align}
where in view of \eqref{dfn:F-star}-\eqref{dfn:F-beta},
\begin{equation}\label{eq:bd-sfx}
\bigcup \,   \Big\{ \bF^{\bb'} :  \bb' \ne \bb,    \sfv_{\bb'} \in \bB_{\sfd} (\sfv_{\bb},  2 \vep e^{\sfKz-1} ) \Big\} 
\subset \bigcup_{j=\tsfKfv+2}^{\sfKz+1} \wbG_{j} \,,
\end{equation}
with $\tsfKfv=\sfKfv+\log \dlc$ as in \Cref{lem-xyN}.  Further,  setting 
\begin{align}
\label{dfn:UUxx}
\UUxx (y) &:=  \UUx(y) \bigcap \{ \etaxnyp = 0 \}\,,
\end{align}
in terms of the events of \eqref{dfn:UUpv},  note that for any $\sfx' \in \bB_{\sfd}(\sfv_{\bb},r_{\sfKo+1})$, 
\begin{align}\label{eq:Uxb-UUxx}
 \UUpx(\Ybz) \cap
 \UUxp(\Ybpz) 
 \subset \UUxx(\Ybz)  \,,
\end{align}
with our next lemma being key to the bound of \eqref{bd:Late-brch}.  
\begin{lemma}\label{lem:dec-2mm}
Fix $\dlc>0$ and large $\sfKp$.  Then,  for all $\sfKfv$ large enough,  
any $\bb$,  $y$ and vector $\uu{\sfx}_{\bb}$ of the start and end points of the $[\sfs_{\sfKo}(y)]$ 
first $\sfKo$-level excursions of $\bb$, 
\begin{equation}\label{eq:bd-2mm}
\sum_{\sfx \in \bF^{\bb}}  \sum_{j=\tsfKfv+2}^{\sfKz+1} \;\,
\sum_{\sfx' \in \wbG_j}
\bP \big[ \,  \UUxx(y) \big| \underline{\sfx}_\bb \big] \le e^{-\frac{1}{3} \sfKfv^\gak}  \bP(\UUpy(y) | \uu{\sfx}_\bb) \,.
\end{equation}
\end{lemma}
Indeed,  when $y=\Ybz$ and $\uu{\sfx}$ is the vector of start and end points of the first $\TKoz$ 
excursions of block $\bb$ at level  $\sfKo$,  the conditional probabilities on the \abbr{lhs} of \eqref{eq:bd-2mm} 
are precisely $\bP \big[ \,  \UUxx(\Ybz) \, |\, \GG^{\sfv_\bb}_{\sfKo} \big]$,  while the one on its right hand side is 
$\bP(\UUp\,|\, \GG^{\sfv_\bb}_{\sfKo})$.  As $\FFuz \subset \HHu \subset \GG^{\sfv_\beta}_{\sfKo}$
for any  $\bb \in u$ (compare \Cref{def:Fvjs} and \Cref{rmk:z-a-vep}),  we have then  
that 
\[
\sum_{\sfx \in \bF^{\bb}} \sum_{j=\tsfKfv+2}^{\sfKz+1} \;\,
\sum_{\sfx' \in \wbG_j}
\bP \big[ \,  \UUxx(\Ybz) \big| \FFuz \big] \le e^{-\frac{1}{3} \sfKfv^\gak}  \bP(\UUpy(\Ybz) | \FFuz) \,.
\]
Considering $\Ptwl_{\bb}$ of 
\eqref{bd:Late-brch} we thus get in view of \eqref{eq:Ub-Uxb},  \eqref{eq:bd-sfx},  \eqref{eq:Uxb-UUxx} 
and the union bound,  that 
\begin{align}
\Ptwl_{\bb}
 &\le
\sum_{\sfx \in \bF^{\bb}} \sum_{j=\tsfKfv+2}^{\sfKz+1} \;\,
\sum_{\sfx' \in \wbG_j}
\bP \big[ \,  \UUpx(\Ybz) \cap \UUxp(\Ybpz) \big| \FFuz \big]  \nonumber \\
& \le  \sum_{\sfx \in \bF^{\bb}} \sum_{j=\tsfKfv+2}^{\sfKz+1} \;\,
\sum_{\sfx' \in \wbG_j}  \bP \big[ \,  \UUxx(\Ybz) \big| \FFuz \big]  
\le  e^{-\frac{1}{3} \sfKfv^\gak}   \bP\big(\UUpy(\Ybz)  \big| \FFuz
 \big) \,.
\label{eq:pre-Late}
\end{align}
Since the conditional probability on the right hand side of \eqref{eq:pre-Late} 
 is $\Pon_\bb$ (see  \eqref{dfn:P1-}),  this is precisely our claim \eqref{bd:Late-brch}.

\begin{proof}[Proof of \Cref{lem:dec-2mm}] Fix $\bb$,  $\sfx \in \bF^\bb$,  $j \in [\sfKfv+\log \dlc+2,\sfKo+1]$, 
$\sfx' \in \wbG_j$,  $y$ and the vector $\uu{\sfx}_{\bb}$ of the start and end points
of the $[\sfs_{\sfKo}(y)]$ first $\sfKo$-level excursions of block $\bb$.  Set
$\wt{n}':=n_\vep-\ell \le n'$ for $\ell:=j+6\log j$,
$\sfv:=\sfv_{\bb}$ and define the following quantities in terms of the laws $\bar \bP^{(k)}_y$ of  \Cref{def:Fvj},
\begin{align}\label{dfn:gxxp}
g_{\sfx,\sfx'}(y_2,y_3) &:=
\sum_{y_4 \in \sfI_{\wt{n}'}
} \bar\bP^{(\sfKth)}_{y_3}  \Big(\BB^{(\vep)}_{[\sfKth,\wt{n}']} (y_3,y_4)  \Big) 
  \bar\bP^{(\sfKt)}_{y_2} \Big(
  A^{\ell,\sfv}_{\sfx,\sfx'} 
 (y_4) \,|\,  \BB^{(\vep)}_{[\sfKth,\wt{n}']} (y_3,y_4)  \Big) \,,\\
\BB^{(\vep)}_{[\sfKth,\wt{n}']} (y_3,y_4) &:=  \big\{ \wh{\eta}(\sfKth) = y_3 \big \} 
\!\!\!\!
\bigcap_{k \in [\sfKth,\wt{n}']} \!\!\!\!\! \{ \wh{\eta} (k) \in \sfI_k \} \bigcap \big\{ \wh{\eta}(\wt{n}') = y_4 \big \} \,,
\label{eq:B-K3n} 
\\
\label{eq:term2}
A^{\ell,\sfv}_{\sfx,\sfx'}(y_4) &:= \Big\{  T_{\sfx,n_\vep+1}^{\sfv,\wt{n}',m} =0, 
T_{\sfx',n_\vep+1}^{\sfv,\wt{n}',m} =0,    \;\;  m=\sfs_{\wt{n}'} (y_4) \in \bN \Big\},
\end{align}
where throughout this proof,  the traversal times of \eqref{eq-calgary}
and $\wh{\eta}(k)$ are for excursions around $\sfv$ during the 
first $\sfs_{\sfKo}(y)$ such $\sfKo$-level excursions.  A key step in our proof is to show that
\begin{align}
\bP [ \,  \UUxx(y) | \uu{\sfx}_\bb] \le 
4 & \sum_{y_2 \in \sfIs_2,  y_3 \in \sfIs_3} \bP \Big[ \BWlb \bigcap \{ \why_ {\bb,\sfKt} = y_2,
\why_{\bb,\sfKth}=y_3 \}  \, \big| \,  \uu{\sfx}_\bb \Big] g_{\sfx,\sfx'}(y_2,y_3) 
  \,,
\label{eq:dec-2mm}
\end{align}
where as in \eqref{eq:BW-K12}, 
\begin{align} \label{B-K1-K2n}
\BWlb & :=  \WWy 
\bigcap_{k \in [\sfKo,\sfKt]} \big\{ \whetaky \in \sfI_k  \big\}   \,.
\end{align}
To this end,  we 
omit from $\UUxx$ the excursion constraints at levels $(\wt{n}',n'']$ and decompose the resulting event 
as in Steps {\bf I} and {\bf II} of \Cref{sec:angle-diff},  now with $\wt{n}'$ and $\wt{n}'':=n_\vep-1-j \vee \sfKfv$
in lieu of $n'$ and $n''$.  Specifically,  
we have in analogy to \eqref{eq:decom-U} and \eqref{eq:B-K1-K2}-\eqref{eq:decomp}  that
\begin{align}\label{eq:decom-Un}
  \UUxx(y) &\subset \bigcup_{y_2,y_3,y_4} \Big\{  \wh{\BW}^{(\vep)}_{[\sfKo,\sfKt]} (y_2|y) 
  \bigcap \BB^{(\vep)}_{[\sfKth,\wt{n}']} (y_3,y_4) \bigcap 
  A^{\ell,\sfv}_{\sfx,\sfx'} (y_4)
 \Big\},
 \end{align}
 \begin{align}
 \label{eq:B-K1=K2}
\wh{\BW}^{(\vep)}_{[\sfKo,\sfKt]} (y_2|y) := \BWlb \bigcap \{ \why_ {\bb,\sfKt} = y_2 \} \,,
\end{align}
 and the union is over $y_2 \in \sfIs_2$,  $y_3 \in \sfIs_3$ and $y_4 \in \sfI_{\wt{n}'}$.
Since the number 
of $\wt{n}'$-level excursions around $\sfv$ and 
\[
 \wh{\BW}^{(\vep)}_{[\sfKo,\wt{n}']}(y_4|y) := 
  \bigcup_{y_2,y_3} \Big\{ \wh{\BW}^{(\vep)}_{[\sfKo,\sfKt]} (y_2|y) 
  \bigcap \BB^{(\vep)}_{[\sfKth,\wt{n}']} (y_3,y_4) \Big\},
\] 
are $\GG^{\sfv}_{\wt{n}'}$-measurable,  it follows from \eqref{eq:decom-Un} that 
\begin{equation}\label{eq:cnd-at-ell}
\bP [ \,  \UUxx(y) | \GG^{\sfv}_{\wt{n}'}] \le \sum_{y_4} 
{\one}_{\wh{\BW}^{(\vep)}_{[\sfKo,\wt{n}']}(y_4|y)} \,\,
\bP
\Big( A^{\ell,\sfv}_{\sfx,\sfx'} (y_4)   \, | \,  \GG^{\sfv}_{\wt{n}'} 
\Big) \,.
\end{equation}
Denoting by $\CCe_{\wt{n}',\wt{n}''} (y_4,y_5)$ having precisely $\sfs_{\wt{n}''} (y_5)$ 
of our $\sfs_{\wt{n}'} (y_4)$ level $\wt{n}'$ excursions enter $\bB_{\sfd}(\sfv,r_{\wt{n}''})$,  we set 
\begin{equation}\label{dfn:alpha-x}
\alpha_{\sfx,\sfx'} (y_5|y_2,y_3,y_4) := \bar \bP_{y_2}^{(\sfKt)} \Big( A^{\ell,\sfv}_{\sfx,\sfx'} (y_4)
 \,  | \,  \BB^{(\vep)}_{[\sfKth,\wt{n}']} (y_3,y_4) \bigcap 
\CCe_{\wt{n}',\wt{n}''} (y_4,y_5) \Big) \, .
\end{equation}
The $\GG^\sfv_{\wt{n}'}$-measurable 
event $\BB^{(\vep)}_{[\sfKth,\wt{n}']} (y_3,y_4)$ specifies the number $m=\sfs_{\wt{n}'}(y_4)$ of level $\wt{n}'$ excursions,  with $\alpha_{\sfx,\sfx'} (y_5|y_2,y_3,y_4)$ the conditional probability of event 
$A^{\ell,\sfv}_{\sfx,\sfx'} (y_4)$ when precisely $\sfs_{\wt{n}''}(y_5)$ of the  $m$ excursions reach level $\wt{n}''$
and the starting and end points of those $m$ excursions are chosen according to the law induced by 
$\bar \bP_{y_2}^{(\sfKt)}\big(\cdot|\BB^{(\vep)}_{[\sfKth,\wt{n}']} (y_3,y_4)\big)$.  
Here $m=\sfs_{\wt{n}'} (y_4)
 \le 5 \ell^2$,  while $\sfKfv$ large guarantees 
having $\wt{n}''-\wt{n}' \ge 2 \log (5 \ell^2)$.  Further,  
for $\sfx' \in \wbG_j$ and fixed $y_4$,  the event $A^{\ell,\sfv}_{\sfx,\sfx'}(y_4)$ 
is 
measurable on the path within $\bB_{\sfd}(\sfv,r_{\wt{n}''})$,  so by Lemma \ref{recursionends2},
we have similarly to \eqref{eq:K3-angles},  that for any $y_2 \in \sfIs_2$,  $y_3 \in \sfIs_3$
and relevant value of $y_4$,
\begin{align}
\bP \Big( A^{\ell,\sfv}_{\sfx,\sfx'} (y_4)   \, | \,  \GG^{\sfv}_{\wt{n}'} \Big) 
= \sum_{y_5}  \bP  \big( A^{\ell,\sfv}_{\sfx,\sfx'} & (y_4)   \, | \,  \GG^{\sfv}_{\wt{n}'},  \CCe_{\wt{n}',\wt{n}''} (y_4,y_5) 
\big) \bP \big( \CCe_{\wt{n}',\wt{n}''} (y_4,y_5)  \, | \,  \GG^{\sfv}_{\wt{n}'} ) 
\nonumber \\
&\le 2 \sum_{y_5} \alpha_{\sfx,\sfx'} (y_5|y_2,y_3,y_4) \bP \big( \CCe_{\wt{n}',\wt{n}''} (y_4,y_5)  \, | \,  \GG^{\sfv}_{\wt{n}'} ) \,.
\label{BRZ-93-2mt}
\end{align}
Here we used the fact that  the conditioning on $\GG^{\sfv}_{\wt{n}'}$ in $ \bP  \big( A^{\ell,\sfv}_{\sfx,\sfx'}  (y_4)   \, | \,  \GG^{\sfv}_{\wt{n}'},  \CCe_{\wt{n}',\wt{n}''} (y_4,y_5) 
\big)$ is the same as conditioning on the starting and end points of those $m$ excursions,  see the proof of 
Lemma \ref{recursionends2}.

Next,  combine \eqref{eq:cnd-at-ell} 
and \eqref{BRZ-93-2mt},  then evaluate the conditional expectation given 
$\GG^{\sfv}_{\sfKt} \subset \GG^{\sfv}_{\wt{n}'}$, 
on which $\wh{\BW}^{(\vep)}_{[\sfKo,\sfKt]}(y_2|y)$ 
is measurable.  Utilizing Lemma \ref{recursionends}
 similarly to \eqref{eq:g3},  this leads to 
\begin{equation}\label{eq:cnd-at-K2}
\bP [ \,  \UUxx(y) | \GG^{\sfv}_{\sfKt}] \le 4 \sum_{y_2,y_3} 
{\one}_{\wh{\BW}^{(\vep)}_{[\sfKo,\sfKt]}(y_2|y)} 
\bP \big( \wCe_{\sfKt,\sfKth} (y_2,y_3) \, | \,  \GG^{\sfv}_{\sfKt} \big) f_{\sfx,\sfx'} (y_2,y_3)\,,
\end{equation}
for events $\wCe_{\sfKt,\sfKth}(\cdot,\cdot)$ as in \eqref{def:Cjjp} (now in terms of our
$\sfs_{\sfKo}(y)$ first $\sfKo$-excursions around $\sfv$),  and 
\begin{align}\label{dfn:fxxp}
f_{\sfx,\sfx'}(y_2,y_3) := \sum_{y_4,y_5}  \alpha_{\sfx,\sfx'} (y_5|y_2,y_3,y_4)
 \bar\bP^{(\sfKt)}_{y_2} 
 \Big(\BB^{(\vep)}_{[\sfKth,\wt{n}']} (y_3,y_4) \bigcap \CCe_{\wt{n}',\wt{n}''} (y_4,y_5)
  \,| \,  \wCe_{\sfKt,\sfKth} (y_2,y_3)  \Big) \,.
 \end{align}
Note from \eqref{eq:B-K3n} that 
$\BB^{(\vep)}_{[\sfKth,\wt{n}']} (y_3,y_4) \subset  \wCe_{\sfKt,\sfKth} (y_2,y_3)$
under the law $\bar\bP^{(\sfKt)}_{y_2}$,  so plugging \eqref{dfn:alpha-x} into \eqref{dfn:fxxp},  we arrive at
  \begin{align*}
  f_{\sfx,\sfx'}(y_2,y_3)
&= \sum_{y_4,y_5}  
 \bar\bP^{(\sfKt)}_{y_2} 
 \Big(\BB^{(\vep)}_{[\sfKth,\wt{n}']} (y_3,y_4) \bigcap \CCe_{\wt{n}',\wt{n}''} (y_4,y_5) \bigcap A^{\ell,\sfv}_{\sfx,\sfx'} 
 (y_4)
  \,| \,  \wCe_{\sfKt,\sfKth} (y_2,y_3) \Big) \\
  &=\sum_{y_4}  
 \bar\bP^{(\sfKt)}_{y_2} 
 \Big(\BB^{(\vep)}_{[\sfKth,\wt{n}']} (y_3,y_4) \bigcap A^{\ell,\sfv}_{\sfx,\sfx'} 
 (y_4)
  \,| \,  \wCe_{\sfKt,\sfKth} (y_2,y_3) \Big) \\
 &= \sum_{y_4}  
 \bar\bP^{(\sfKt)}_{y_2} 
 \big(\BB^{(\vep)}_{[\sfKth,\wt{n}']} (y_3,y_4)  \,| \,  \wCe_{\sfKt,\sfKth} (y_2,y_3) \big) 
 \bar\bP^{(\sfKt)}_{y_2} \big( A^{\ell,\sfv}_{\sfx,\sfx'}  (y_4) \,|\,
 \BB^{(\vep)}_{[\sfKth,\wt{n}']} (y_3,y_4) \big) \,.
\end{align*}
This is precisely  $g_{\sfx,\sfx'}(y_2,y_3)$ of \eqref{dfn:gxxp},  since 
\begin{align*} 
 \bar\bP^{(\sfKt)}_{y_2} \big(
 \BB^{(\vep)}_{[\sfKth,\wt{n}']} (y_3,y_4)  \,| \,  \wCe_{\sfKt,\sfKth} (y_2,y_3) \big) 
 =
   \bar\bP^{(\sfKth)}_{y_3} 
 \big(\BB^{(\vep)}_{[\sfKth,\wt{n}']} (y_3,y_4)  \big) \,,
\end{align*}
thanks to the Markov property of the concentric excursion counts around $\sfv$.  

To proceed,  replace $f_{\sfx,\sfx'} (\cdot)$ by $g_{\sfx,\sfx'}(\cdot)$ in \eqref{eq:cnd-at-K2}
and evaluate the conditional expectation of both sides given 
$\GG^{\sfv}_{\sfKo} \subset \GG^{\sfv}_{\sfKt}$ to find that
\begin{equation}\label{eq:pre-dec-2mm}
\bP [ \,  \UUxx(y) | \GG^{\sfv}_{\sfKo}] \le 4 \sum_{y_2,y_3} \bP \Big[\wh{\BW}^{(\vep)}_{[\sfKo,\sfKt]}(y_2|y)
\bigcap  \wCe_{\sfKt,\sfKth} (y_2,y_3) \, | \,  \GG_{\sfKo}^{\sfv} \Big] g_{\sfx,\sfx'} (y_2,y_3)\,.
\end{equation}
This amounts to \eqref{eq:dec-2mm},  since $\uu{\sfx}_\bb$ is $\GG^{\sfv}_{\sfKo}$-measurable 
and one easily checks,  in view of \eqref{eq:B-K1=K2},  that the conditional probabilities 
on the right hand side of \eqref{eq:pre-dec-2mm} match those on 
\eqref{eq:dec-2mm}.

\blue{Now,  $\bar\bP^{(\sfKt)}_{y_2} \big(  \cdot
 \,|\,  \AA^{(m)} \big)$,  at $\AA^{(m)}=\BB^{(\vep)}_{[\sfKth,\wt{n}']} (y_3,y_4)$ 
which is determined by counts of excursions around $\sfv$
(for radii $R_{i+1-\sfKt}=r_i$,  $i \le \wt{n}'$),  
puts us
in the setting of \Cref{rmk:off-center-applies} for
$R=r_{\wt{n}'} = e^\ell \vep$ and precisely $m=\sfs_{\wt{n}'}(y_4)$ paths 
from $\bB_{\sfd}(\sfv,R)$ to $\bB_{\sfd}(\sfv,e R)$ completed during our 
$
\sfs:=
\sfs_{\sfKt}(y_2)$ \abbr{iid} outer Brownian excursions.   As
$A^{\ell,\sfv}_{\sfx,\sfx'} (y_4)$ requires avoiding 
$\bB_\sfd(\sfx,\vep/e) \cup \bB_\sfd(\sfx',\vep/e)$,  the 
conditional probability 
$\bar\bP^{(\sfKt)}_{y_2} \Big(  A^{\ell,\sfv}_{\sfx,\sfx'}  (y_4) \,|\,  \BB^{(\vep)}_{[\sfKth,\wt{n}']} (y_3,y_4)  \Big)$ 
is of the form of 
\Cref{lem-xyN},  with $\ell=j+6 \log j$ and $\sfd(\sfx,\sfx')$ 
as given in \eqref{dfn:tGj}.  Thus,  for $C'<\infty$ and large $\sfKfv$,   
any such $\sfx,\sfx',\sfv$ and $y_2,y_3,y_4$,   we deduce from \eqref{xyN.2} that 
\begin{equation}\label{eq:Axxp-bd}
 \,  \bar\bP^{(\sfKt)}_{y_2} \Big(
  A^{\ell,\sfv}_{\sfx,\sfx'} 
 (y_4) \,|\,  \BB^{(\vep)}_{[\sfKth,\wt{n}']} (y_3,y_4)  \Big) \le  C' e^{-4j} e^{-3y_4} \,.
\end{equation}}
In addition,  \blue{by \Cref{prop:BarrierSecGWPropUB}
we get that uniformly over the relevant range of $y_3$ and $b$}, 
\begin{align}\label{eq:barrier-K3-wtn}
\bar\bP^{(\sfKth)}_{y_3} 
 \Big( \big\{ \wh{\eta}(\wt{n}') \in [b,  b+1] \big \}  
\!\!\!\!
\bigcap_{k \in [\sfKth,\wt{n}']} \!\!\!\!\! \{ \wh{\eta} (k) \in \sfI_k \} 
\Big) \blue{ \le \KQ_{[\sfKth,\wt{n}']}^{(\sfKth,\wt{n'})}(y_3,b)}
 \leq 
\frac{C}{\sqrt{\ell}} \,  
b e^{2b} \, 
e^{-2(\wt{n}'-\sfKth)} \,  y_3 e^{-2y_3}  \,.
\end{align}
Since $ \sfI_{\wt{n}'} \subset [j^\gak,\infty)$ and $j \ge \tsfKfv \gg 1$,  plugging \eqref{eq:Axxp-bd}
into \eqref{dfn:gxxp} and then utilizing \eqref{eq:barrier-K3-wtn},  yields that 
\begin{align*}
 g_{\sfx,\sfx'}(y_2,y_3) &\le  C' e^{-4j} 
\sum_{b \in \sfI_{\wt{n}'} \cap \bN} e^{-3b} \, 
\bar\bP^{(\sfKth)}_{y_3}  \big(\bigcup_{y_4 \in [b,b+1)} \BB^{(\vep)}_{[\sfKth,\wt{n}']} (y_3,y_4)  \big) \\
&\le   
  C'  e^{-4j}  e^{-2(\wt{n}'-\sfKth)} \,  y_3 e^{-2y_3} 
 \sum_{b \in  \sfI_{\wt{n}'} \cap \bN} b e^{-b} \le \frac{1}{4} e^{-2(\wt{n}'-\sfKth)} e^{-4j} 
 e^{-\frac{1}{2} j^\gak}  \,  y_3 e^{-2y_3}  \,.
\end{align*}
Plugging this bound  in \eqref{eq:dec-2mm},  yields that for $\gth(y_3):={\bf 1}_{\sfIs_3}(y_3)y_3 e^{-2 y_3}$, 
any $\sfx \in \bF^\bb$ and $\sfx' \in \wbG_{j}$, 
\begin{align}\label{eq:bd-2mm-xxp}
\bP [ \,  \UUxx(y) |  \uu{\sfx}_\bb ] \le  e^{-2 (\wt{n}'-\sfKth)} e^{-4 j} e^{-\frac{1}{2} j^\gak}
\, \bE \big[ \gth(\wh{\eta}^{\,(y)}_{\bb,\sfKth});  \BWlb  \,|  \uu{\sfx}_\bb \big]  \,.
\end{align}
Now,  by the same reasoning that led to 
\eqref{eq:uux-to-emp-meas},  for any $\bb$,  $y$ and $\uu{\sfx}_\bb$,  
\[
\bE \big[ \gth(\wh{\eta}^{\,(y)}_{\bb,\sfKth});  \BWlb  \,| \uu{\sfx}_\bb \big] =
\bE \big[ \gth(\wh{\eta}^{\,(y)}_{\bb,\sfKth});  \BWlb  \,| \LKoy\big] \,,
\]
hence by \eqref{eq:good-cover},  for $\sfKp$ large,  some $\uu \alpha(\sfKp)>0$
all $\sfKfv$ large enough,  any $\beta$,  $y$ and $\uu{\sfx}_\bb$ 
\begin{equation}\label{bd:g-cover}
 \frac{1}{2} \uu \alpha(\sfKp) e^{-2 (n''-\sfKth)} 
\, \bE \big[ \gth(\wh{\eta}^{\,(y)}_{\bb,\sfKth});  \BWlb  \,| \uu{\sfx}_\bb \big]  \le \bP(\UUpy(y) | \uu{\sfx}_\bb) \,.
\end{equation}
Next,  recall \eqref{dfn:F-beta} and 
\eqref{dfn:tGj} that $|\bF^\bb| \le  c e^{2 \sfKfv}$ and 
$|\wbG_{j}| \le c e^{2j}$ for some $c(\wt{c})$ finite,  any $\bb$ and $j$.  Thus,   
with $\wt{n}'=n_\vep-j-6 \log j$,  we have for any $\sfKp$,  $\dlc>0$ and all $\sfKfv$ large enough,    
\begin{align}
|\bF^\bb| \sum_{j=\tsfKfv+2}^{\sfKz+1}  \,  |\wbG_j| \, 
e^{-2 (\wt{n}'-\sfKth)} e^{-4 j} e^{-\frac{1}{2} j^\gak}
&\le c^2 \,  e^{-2(n''-\sfKth)} 
 \sum_{j=\tsfKfv+2}^{\infty} j^{12} e^{-\frac{1}{2} j^\gak} \nn \\
 &\le \frac{1}{2} \uu \alpha(\sfKp) e^{-2 (n''-\sfKth)} 
 e^{-\frac{1}{3} \sfKfv^\gak} \,.
 \label{eq:sum-xxp}
\end{align}
Combining  \eqref{eq:bd-2mm-xxp} and \eqref{eq:sum-xxp} we deduce that 
\[
\sum_{\sfx \in \bF^{\bb}} \sum_{j=\tsfKfv+2}^{\sfKz+1} \;\,
\sum_{\sfx' \in \wbG_j}
\bP \big[ \,  \UUxx(y) \big| \underline{\sfx}_\bb \big] \le  \frac{1}{2} \uu \alpha(\sfKp) e^{-2 (n''-\sfKth)} 
 e^{-\frac{1}{3} \sfKfv^\gak} 
\, \bE \big[ \gth(\wh{\eta}^{\,(y)}_{\bb,\sfKth});  \BWlb  \,| \uu{\sfx}_\bb \big],
\]
which together with \eqref{bd:g-cover} implies that \eqref{eq:bd-2mm} holds as well.
\end{proof}

\subsection{Tightness:  Proof of \Cref{prop:tight}}\label{sec:tight}
We \blue{start with tightness of the left-tail of $\sqrt{\smash[b]{\CCT}}-\osfm_\vep$,  which in view of 
\eqref{eq:cov-p-lim} will also establish \eqref{ct2.25-new}.  Specifically, } fix small $\dlc$ and $\srs=r_0$,  
recalling \eqref{def:barQu},  \eqref{eq-good} and \eqref{def:barQub} that for any
$z \in \bR$,  $\sfKz$,  $\sfKfv$,  $u$ and $\vep$,
\[
1-Q^{(\vep)}_{u,z} \le 1-Q^{(+,\vep)}_{u,z} \le \exp(-Q^{(+,\vep)}_{u,z}) \,.
\]  
Hence,  in view of \eqref{eq:CC-Quz} we have that for $L'=|\log \delta|$ and any $\delta>0$,  
\begin{equation}\label{eq:lft-tail-equiv}
\bP(\CCex \le \sfk_\vep(z) ) = \bE \Big[   \prod_{u=1}^{\sfM}  (1- Q^{(\vep)}_{u,z} ) \Big] 
\le \delta + \bP \Big[  \sum_{u=1}^{\sfM} Q^{(+,\vep)}_{u,z} \le L'  \Big] \,.
\end{equation}
Further,  from
\eqref{dfn:wt-Q},  \eqref{dfn:wt-Q2} and \eqref{eq:Q2-Q1-Qs},  
for large  $\sfKfv$,  $\sfKz \gg \sfKfv$,  
all $z,u$,  small $\vep$ and $\FFuz$-a.e.  
\begin{equation}\label{eq:au-cu}
Q^{(+,\vep)}_{u,z} = \bP \big[    \Lambda^{(+,\vep)}_{u,z} \ge 1 |  \FFuz \big] \ge 
\frac{\bE \big[ \,\Lambda^{(+,\vep)}_{u,z} \,|\, \FFuz  \big]^2}
{ \bE\big[ (\Lambda^{(+,\vep)}_{u,z})^2  | \FFuz \big]} = \frac{\big(\wt{Q}^{(1,\vep)}_{u,z}\big)^2}
{ \wt Q^{(2,\vep)}_{u,z} + \wt{Q}^{(1,\vep)}_{u,z}} \ge \frac{\big(\wt{Q}^{(1,\vep)}_{u,z}\big)^2}
{2  \wt Q^{(1,\vep)}_{u,z} +  \wt Q^{(\star,\vep)}_{u,z} } \,.
\end{equation}
Proceeding to bound the right hand side of \eqref{eq:au-cu}, 
recall \eqref{eq:Ball-to-angle}
that $\Vb (\LKoz) \in \sfIw$
when $\LKoz \in \bB^\infty_{\bar \mu} (\boetasKo)$,
and by \eqref{eq:g-def},  if in addition $\Ybz  \in \sfI_{\sfKo}$ then for any $\sfKo$ large enough,
\begin{equation}\label{bd:e-3Ko}
\gstb\big( \Ybz,\Vb (\LKoz) \big) \ge  
\inf_{y \in \sfI_{\sfKo}} \{ y e^{-2 y} \} \inf_{v \in \sfIw} \{ e^{- 2 v} \}  
\ge e^{-3 \sfKo^{\kapp}} \,.
\end{equation}
We thus have for the bounded below $\alpha(\sfKfv,\sfKp)$ of  \Cref{prop:angle-diff}
and $\bKz$ of \eqref{def:bar-bKz},
for all large $\sfKo$,  
\begin{equation}\label{bd:e-4Ko}
\frac{1}{2} \alpha(\sfKfv,\sfKp) \bKz \gstb\big( \Ybz,\Vb (\LKoz) \big) \ge e^{-4 \sfKo^{\kapp}} \,.
\end{equation}
Further,  recall \Cref{rmk:z-a-vep},  for any $\bb \in u$ the law of $(\Ybz,\LKoz)$ is
unchanged by conditioning on $\HHk$ 
instead of only on $\FFuz$.  Hence,  by \eqref{eq:wtQ1-bd} and the preceding, 
for any such $\sfKz,\sfKfv,\sfKp$ and $\vep$ small enough, 
\begin{align*} 
\wt{Q}^{(1,\vep)}_{u,z} &\ge \frac{1}{2} \alpha (\sfKfv,\sfKp) \bKz
\,  \avsum_{\bb\in u} \bE \Big[ 
\gstb \big( \Ybz, \Vb (\LKoz)
\big); \Ybz \in \sfI_{\sfKo},  \LKoz \in \bB^\infty_{\bar \mu} (\boetasKo) 
 \,\big|\,  \HHk \Big] \\
&\ge  
 e^{-4 \sfKo^{\kapp}} \Big[ 1 - \avsum_{\bb\in u}
 \bP \big(  \{ \Ybz  \notin \sfI_{\sfKo} \} \cup
\{ \LKoz \notin \bB^\infty_{\bar \mu} (\boetasKo) \} 
\big| \,  \HHk \big) \Big]  \,,
\end{align*}
which upon using both \eqref{bd:e-3Ko} and \eqref{bd:e-4Ko},  leads to
\begin{align}\label{eq:av-Qub}
\Big( \wt{Q}^{(1,\vep)}_{u,z} \Big)^2 \ge  \frac{1}{2} \alpha (\sfKfv,\sfKp) \bKz  e^{-7 \sfKo^{\kapp}} 
\Big[ 1 - 2 \avsum_{\bb\in u}
 \bP \big(  \{ \Ybz  \notin \sfI_{\sfKo} \} \cup
\{ \LKoz \notin \bB^\infty_{\bar \mu} (\boetasKo) \} 
\big| \,  \HHk \big) \Big]  \,.
\end{align}
Moreover,  by \eqref{eq:wtQ1-bd} and \eqref{dfn:Qst} 
\begin{align*}
\wt{Q}^{(1,\vep)}_{u,z} \le  2 \alpha (\sfKfv,\sfKp) \bKz,  \qquad 
\wt{Q}^{(\star,\vep)}_{u,z} \le 1 \,,
\end{align*} 
which together with \eqref{eq:au-cu} and \eqref{eq:av-Qub} results with 
\begin{align}\label{eq:Qpv-bd}
Q^{(+,\vep)}_{u,z} \ge 
%
  e^{-8 \sfKz^{\kapp}} 
\Big[ 1 - 2 \avsum_{\bb\in u}
 \bP \big(  \{ \Ybz  \notin \sfI_{\sfKo} \} \cup
\{ \LKoz \notin \bB^\infty_{\bar \mu} (\boetasKo) \} 
\big| \,  \HHk \big) \Big]  \,.
\end{align}
Hence,   utilizing \eqref{eq:Qpv-bd},  we have that for $\sfKfv$,  $\sfKp$ and $\sfKz$ as above,    
\begin{align}\label{eq:wt-Q-lbd}
\frac{1}{\sfM} \sum_{u=1}^{\sfM} Q^{(+,\vep)}_{u,z} 
\ge e^{-8 \sfKz^{\kapp}} \Big( 1 - 2\bE\big[ \Gamma^{(\vep)}_{z} \,\big|\,  \HHk \big] \Big) \,,
\end{align}
where $\Gamma^{(\vep)}_{z}$ denotes the fraction of blocks $\bb$ for which either 
$\Ybz  \notin \sfI_{\sfKo}$ or $\LKoz \notin \bB^\infty_{\bar \mu} (\boetasKo)$.
Assuming \abbr{wlog} that
$\sfM = e^{2 \sfKz} \ge 3 e^{8 \sfKz^{\kapp}} L'$,  we deduce from \eqref{eq:lft-tail-equiv} and
\eqref{eq:wt-Q-lbd} that 
\begin{align}\label{eq:lft-tail1}
\varlimsup_{\vep \to 0} \bP\Big(\CCex \le \sfk_\vep(z) \Big) 
\le \delta + 
\varlimsup_{\vep \to 0} \bP \Big(  \bE[ \Gamma^{(\vep)}_{z} | \HHk]
\ge \frac{1}{3}  \Big) \,.
\end{align}
By definition,  
\[
\bE \big[ \Gamma^{(\vep)}_{z} ]  
= \avsum_{\bb} 
 \bP \Big(  \{ \Ybz  \notin \sfI_{\sfKo} \} \cup
\{ \LKoz \notin \bB^\infty_{\bar \mu} (\boetasKo) \}  \Big),
\]
so setting $z=z(\sfKo)=-2\sfKo+\frac{1}{2} \sfKo^{\kapp}$,  it follows from 
\Cref{lem:clt-left-tight} and Markov's inequality that for all 
$\sfKo$ large enough the right hand side of \eqref{eq:lft-tail1} is at most $4 \delta$.  
By \eqref{ct2.25cjj},  further increasing $\sfKz$ (and hence $\sfKo$),  yields that 
\begin{align}\label{eq:prelim-left-tail}
\varlimsup_{\vep\to 0} \bP \big( \sqrt{\CCd} \le  \sfm_\vep + z \big) \le 
\delta + \varlimsup_{\vep \to 0} \bP \big(  \, \CCex \le  \sfk_\vep(z)  \big) \le 5 \delta \,.
\end{align}
To recap,  \eqref{eq:prelim-left-tail} holds for any $\dlc$,  $\srs$,  
large $\sfKfv$ and all $\delta>0$,  provided $\sfKz \ge \sfKz(\delta,\srs,\sfKfv)$ and $z=z(\sfKo)$.
Consider the normalized $\vep$-cover time $\CCdz$ of the 
complement of the square (of edge-length $\srs$ around the Brownian motion starting point),
used for $\CCp$ of \Cref{lem:right-tight} (with both cover times independent of $\dlc$,  $\sfKz$ and $\sfKfv$). 
By definition
\begin{equation}\label{ord-cover}
\CCd 
 \le \CCdz \le  \max\{ \CCdz,\CCp\} = \CCT \,.
\end{equation}
Hence,  fixing $r_0$ and $\delta>0$ we get from \eqref{eq:prelim-left-tail} that for $\srs=r_0$
and some finite $z_\delta=z_\delta(r_0)$,
\begin{equation}\label{eq:left-tight}
\varlimsup_{\vep\to 0} \bP \big( \sqrt{\CCdz} \le  \osfm_\vep - z_\delta \big) \le 5 \delta \,,
\end{equation}
which,  since $\srs \mapsto \CCdz$ is non-increasing applies for all $\srs \le r_0$,  with
the same finite $z_\delta$.  Comparing the latter fact with \eqref{eq:cov-p-lim} at $z>z_\delta$ 
and taking $\delta \to 0$,  we conclude that   
\begin{equation}\label{eq:p-negl}
\lim_{\srs \to 0} \varliminf_{\vep \to 0} \Pbm \big( \CCT = \CCdz \big) \ge
\lim_{\srs \to 0} \varliminf_{\vep \to 0} \Pbm \big( \CCT = \CCdz > \CCp \big) = 1  \,.
\end{equation}
In particular,  \eqref{eq:left-tight} holds with $\CCT$ instead of 
$\CCdz$,  implying the tightness of $(\sqrt{\smash[b]{\CCT}}-\osfm_\vep)_-$.
Moreover,  by definition $\CCz \le \max\{\CCd,\CCp\}  \le \CCT$,  hence combining \eqref{ct2.25cj} 
with \eqref{eq:p-negl} implies that 
\[
\varlimsup_{\dlc,\srs,\sfKfv,\sfKz,\vep} \bP \big( \CCd \ne \CCT  \big) = 
\varlimsup_{\dlc,\srs,\sfKfv,\sfKz,\vep} \bP \big( \max\{\CCd,\CCp\} \ne \CCT  \big) \le
\varlimsup_{\dlc,\srs,\sfKfv,\sfKz,\vep} \bP \big( \CCz \ne \CCT  \big) = 0 \,,
\]
as claimed in \eqref{ct2.25-new}.  

Lastly,  given \eqref{ct2.25-new},  it suffices for the (uniform) tightness of 
$\big( \sqrt{\smash[b]{\CCT}} - \osfm_\vep  \big)_+$ to show that 
\[
\lim_{z \to \infty} \varlimsup_{\sfKfv,\sfKz,\vep} \bP \big( \CCd > \sfm_\vep + z  \big) = 0,   \qquad \forall \dlc,\srs=r_0 \,.
\] 
Fix $\dlc$,  $\srs=r_0$,  noting that by \eqref{ct2.25cjj} and \eqref{eq:right-tight},  
for $\CCex$,   $\sfk_\vep(z)$ and $\UUz$ of 
\eqref{cov-time-exc},  \eqref{dfn:bar-svep} and \eqref{dfn:UUz},  
\[
\lim_{z \to \infty} \varlimsup_{\sfKfv,\sfKz,\vep} \bP \big( \CCd > \sfm_\vep + z  \big) = 
\lim_{z \to \infty} \varlimsup_{\sfKfv,\sfKz,\vep} \,\bP \big(  \, \CCex  > \sfk_\vep(z)  \big) =
\lim_{z \to \infty} \varlimsup_{\sfKfv,\sfKz,\vep} \,\bP \big( \bigcup_{\bb} \UUz \big) = 0 \,,
\]
which thereby completes the proof of \Cref{prop:tight}.

\section{Limiting field: Proof of \Cref{prop:Gauss-field-continuity}}\label{sec:end-game}
Recall the decomposition of the Brownian path 
between its 
successive alternate hitting locations $\{\sfz_{j2},  j \ge 0\}$ of $\bTs
$
and its hitting locations $\{\sfz_{j1},  j  \ge 1\}$ of 
$\bTs^\sfe
$.  As explained in  \Cref{sec-exc-count},  
\abbr{wlog} we start our path at $\sfz_{02} \in \bTs$ 
according to the invariant measure $\mu_\sfe$ of 
the Doeblin aperiodic
chain $\{\sfz_{j2},j \ge 0\}$, 
whose uniformly positive transition density
depends only on $\dlc$ and $\sfKz$,  
with ${\sft}_1 = \Ebm_{\mu_{\sfe}} [T_{1}]$ the invariant expected
time between successive visits of $\bTs$.  By the same reasoning,  
the Markov chain of pairs 
$\{(\sfz_{(j-1)2},\sfz_{j1})\}$
has uniformly positive transition density on $\bTs \times \bTs^\sfe$,  thereby being 
also Doeblin aperiodic,  having the unique invariant measure $\bar \mu_{\sfe}$ 
(with $\mu_\sfe$ matching the $\bTs$-marginal of $\bar \mu_\sfe$).  

Recalling the filtration $\HH_k$ consisting of only the part of the 
path during its first $k$ successive inter-city excursions from $\bTs^\sfe$ to $\bTs$
(which start and end at 
$\sfz_{j1}$ and $\sfz_{j2}$,  respectively),   we
associate to each inner-city point $\sfv \in \bA^\sfo_{u}$
the counts $\Tvjz$
of excursions from 
$\partial \bB_\sfd\left(\sfv,r_{j}\right)$ to $ \partial \bB_\sfd\left(\sfv,r_{j-1}\right)$ 
during the first $\sfk_\vep (z)=(\sfm_\vep+z)^2/(2\pi\sft_1)$ 
inter-city excursions,  and set
\begin{equation}\label{17.01ad}
\Yvz:=\sqrt{2 N_{\sfv}} - \sfm_\vep + \sfrho_\vep \sfKo  \quad {\rm for} \quad 
N_{\sfv}:=\TvKoz.
\end{equation}
With $\{\th_{j,\sfv},  j\ge 1\}$ the angle increments of the consecutive $\sfKo$-level excursions around $\sfv$,  we also 
define analogously to \eqref{dfn:hat-L-alpha} (see also \eqref{dfn:emp-ang}),
\begin{equation}\label{def:V-eps-alt}
\Vvz :=  \frac{\sfm_\vep}{2 N_\sfv} \hat{N}^\th_\sfv\quad {\rm for} \quad 
\hat{N}^\th_\sfv := \sum_{j=1}^{N_\sfv}  \bar c [ \langle \cos(\cdot), \bar \mu \rangle - \cos(\theta_{j,\sfv})] \,.
\end{equation}
Similarly to \eqref{dr.2} we construct sub-lattices $\bar \bF_j$ of size $|\bar \bF_j| \le c \,  e^{2j}$ 
\emph{which are independent of $\vep$} and form $\wt{c} \,  e^{-j}$ nets for the union of all blocks in $\bT^2$.
 In \Cref{sec-contheight} 
we establish in the following uniform (in $\vep$) continuity of the fields $(\Yxz,\Vxz)$
over $\bF_{n_\vep}$ (where the block centers $\sfv_\bb$ are),  augmented by the intersection 
$\wt{\bF}_\vep$ of $\cup_u \bA_u^o$ and $\bar \bF_{\bar n_\vep}$,  for $\bar n_\vep:=2 \log (n_\vep \log n_\vep)$.

\begin{proposition}\label{th-heightcontinuity}
Fixing  $\dlc, \sfKz, \sfKfv,z$,  for any $\ga,\ze>0$,  there exists $\de>0$ so
\begin{align}
\limsup_{\vep \to 0} \Pbm \Big(\sup_{\stackrel{\sfx,\sfy \in \bF_{n_\vep} \cup \wt{\bF}_\vep}{\sfd (\sfx, \sfy) \le \de  }}
\|(\Yxz,\Vxz)
-(\Yyz,\Vyz)\|_\infty \ge \ga \Big) &\leq \ze \,.
\label{hc.2a}
\end{align}
\end{proposition}

We supplement this continuity result by proving in \Cref{sec-gauss} the following 
\abbr{clt} for the approximate finite-dimensional distributions (\abbr{fdd}) of these fields.
\begin{proposition}\label{theo-mclt}
Fix $\dlc,\sfKz,z$ and $\{\sfx_{i}: i \le p \} \subset  \wt{\bF}_\vep$,  independent of $\vep$ and $z$.
For $1 \le i \le p$,  let 
\begin{align}\label{def:qin}
q_{\rm in}^{(\vep)}(i) &:= \frac{1}{\sfm_\vep+z} \Big(N_{\sfx_i}-\bE\big[N_{\sfx_i}|\HHk \big]\Big),\\
q_{\rm in}^{(\vep)}(i+p) &:= \frac{1}{\sfm_\vep+z}
\Big(\hat{N}^\th_{\sfx_i}-\bE\big[\hat{N}^\th_{\sfx_i}|\HHk \big]\Big), 
\label{def:qin-th}
\\
q_{\rm out}^{(\vep)}(i) &:= \frac{1}{\sfm_\vep+z} \Big(\bE\big[N_{\sfx_i}|\HHk \big]-\bE_{\bar \mu_\sfe}\big[N_{\sfx_i}\big]\Big),
\label{def:qout}
\\
q_{\rm out}^{(\vep)}(i+p) &:=\frac{1}{\sfm_\vep+z}
 \Big(\bE\big[\hat{N}^\th_{\sfx_i}|\HHk \big]-\bE_{\bar \mu_\sfe}\big[\hat{N}^\th_{\sfx_i}\big]\Big)\,.
\label{def:qout-th}
\end{align}
There are 
non-negative definite matrices $\Sigma_1$,  $\Sigma_2$ independent of $z$,  such that 
\begin{align}
q_{\rm in}^{(\vep)} &\dconv  N\( 0, \Sigma_{1}\),  \quad \mbox{\rm conditionally on 
$\HHk$},\label{mCLT.03}\\
q_{\rm out}^{(\vep)} &\dconv  N\( 0, \Sigma_{2}\).
\label{mCLT.04}
\end{align}   
 \end{proposition}
 
We further evaluate in \Cref{sec-gauss} the following expected values.
\begin{lemma}\label{lem-expi}
For any $\dlc$,  large $\sfKo-\sfKz$ and all $\sfx \in \cup_u \bA^o_u$,
\begin{align}
 \Ebm_{\bar \mu_{\sfe}} \big[T_{\sfx,\sfKo}^{\bTs^\sfe(1)} \big]&=\pi {\bf  t_{1}}, \label{17.02b}\\
 \Ebm_{\bar \mu_{\sfe}}\Bigg[  \sum_{j=1}^{T_{\sfx,\sfKo}^{\bTs^\sfe(1)}}  \cos(\theta_{j,\sfx}) \Bigg]
 &=\pi {\bf  t_{1}} \langle \cos(\cdot), \bar \mu \rangle \,. \label{17.02c}
 \end{align}
 \end{lemma}
In particular,  from \eqref{17.02b} and \eqref{17.02c} it follows that for any $1 \le i \le p$, 
\begin{equation}\label{17.02-cnsq}
\Ebm_{\bar \mu_{\sfe}} [N_{\sfx_i}]=\pi \sft_1 \sfk_{\vep}(z) = \frac{1}{2} (\sfm_\vep+z)^2\, \quad {\rm and} 
\quad \Ebm_{\bar \mu_{\sfe}} [\hat{N}^\th_{\sfx_i}]=0\,.
\end{equation}
Conditional convergence in distribution (to a conditioning-independent limit) implies convergence in distribution.  Thus,  
by \Cref{theo-mclt},  both $q^{(\vep)}_{\rm in}$ 
and 
$q^{(\vep)}_{\rm out}$ 
converge in probability (on a suitably enlarged space),  so in view of the \abbr{lhs} of \eqref{17.02-cnsq},  as $\vep \to 0$,  
\begin{equation}\label{eq:N-conc}
\frac{N_{\sfx_i}}{\bE_{\bar \mu_{\sfe}}[N_{\sfx_i}]} - 1 = \frac{2}{\sfm_\vep+z} [q_{\rm in}^{(\vep)}(i)+q_{\rm out}^{(\vep)}(i)] 
\stackrel{p}{\to} 0 \,. 
\end{equation}
Further,  as $f_m(q):=m [\sqrt{1+2 q/m}-1] =q+O(q^2/m)$,  it follows upon 
plugging \eqref{def:qin},  \eqref{def:qout} and the \abbr{lhs} of \eqref{17.02-cnsq} in \eqref{17.01ad},  that uniformly over $i \le p$,
\begin{align}
\Yxiz -z &= \sfrho_\vep \sfKo + f_{\sfm_\vep+z} \big(
q_{\rm in}^{(\vep)}(i) + q_{\rm out}^{(\vep)}(i)\big) 
=\sfco \sfKo + q_{\rm in}^{(\vep)}(i)+q_{\rm out}^{(\vep)}(i)+o_{\vep} (1)\,.
\label{mCLT.01}
\end{align}
Similarly,  plugging \eqref{def:qin-th},  \eqref{def:qout-th} and the right hand side of \eqref{17.02-cnsq} into \eqref{def:V-eps-alt},
we deduce from \eqref{eq:N-conc} that uniformly over $i \le p$, 
\begin{align}
\Vxiz = \frac{\sfm_\vep}{\sfm_\vep+z} \frac{\Ebm_{\bar \mu_{\sfe}}[N_{\sfx_i}]}{N_{\sfx_i}} &
\big[ q_{\rm in}^{(\vep)}(i+p) + q_{\rm out}^{(\vep)}(i+p) \big] 
=  \big[ q_{\rm in}^{(\vep)}(i+p)+ q_{\rm out}^{(\vep)}(i+p) \big]
 \label{mCLT.02}
+o_{\vep}(1).
\end{align} 
We now derive \Cref{prop:Gauss-field-continuity} from \eqref{mCLT.01}-\eqref{mCLT.02} and
Propositions \ref{th-heightcontinuity},  \ref{theo-mclt}.  Indeed,  
by \Cref{th-heightcontinuity} we have for our bounded-Lipshitz,  compactly 
supported $\gst(\cdot)$ and $\{ \sfX^{(\vep)}_{u,z} \}_u$ of \eqref{def:X-eps},  the existence
of some $\{\sfx_i: i \le p\} \subset \wt{\bF}_\vep$ independent of $\vep$ and $z$,  for $p=p_\sfM$ large enough (to have in particular,  
some $\sfx_i$ within each city),  such that 
 \begin{equation}
\limsup_{\vep \to 0} 
\Pbm \Big(\sum_{u} \big| \sfX^{(\vep)}_{u,z} - \sfX^{(\vep,p)}_{u,z} \big| \ge \ga_p \Big) \le \ga_p \,,  \label{210.2}
 \end{equation}
for some $\ga_p \to 0$ as $p \to \infty$,  and
\begin{equation}\label{def:X-eps-pm}
\sfX^{(\vep,p)}_{u,z}  := \alpha(\sfKfv,\sfKp) \bKz
\,  \avsum_{\sfx_i \in \bA^o_u} \bE \Big[ \gst \big( \Yxiz-z,\Vxiz \big)  
\,\big|\,  \HHk \Big] \,.
\end{equation}
Next,  denoting by 
$\{\xi_{\rm in}(i)\}$ and $\{\xi_{\rm out}(i)\}$ the components of 
$N\(0,\Sigma_{1}\)$ and $N(0,\Sigma_2)$ vectors,  respectively,
it follows from \eqref{mCLT.03},  \eqref{mCLT.01}-\eqref{mCLT.02}, 
and the continuous mapping theorem,  that \eqref{210.2} holds also with
\begin{equation}\label{def:X-eps-pm-alt}
\sfX^{(\vep,p)}_{u,z} = \alpha(\sfKfv,\sfKp) \bKz \,  \avsum_{\sfx_i \in \bA^o_u} 
g^{(p)}_i\big(q_{\rm out}^{(\vep)}(i),q_{\rm out}^{(\vep)}(i+p)\big)\,,
\end{equation}
where 
\begin{align} \label{210.3}
g^{(p)}_i(y,v) := \bE \Big[ \gst \big( \sfco \sfKo +\xi_{\rm in}(i)+y, \xi_{\rm in}(i+p)+v \big)\Big]  \,.
\end{align}
Consider the Polish space $M_1(\bR^\sfM)$ equipped with the bounded-Lipshitz 
norm $\sfd_{\rm BL}(\cdot,\cdot)$ which is compatible with weak convergence (where $\sfd_{\rm BL}(Y,Y')$ 
denotes such distance between the laws of $Y,Y'$).  As $\ggi^{(p)}$ are bounded and 
continuous,  combining \eqref{mCLT.04} with \eqref{210.2} for $\{ \sfX^{(\vep,p)}_{u,z} \}$
of \eqref{def:X-eps-pm-alt},  
we deduce upon adjusting $\ga_p \to 0$,  that the collection 
\begin{equation}\label{def:210.2-alt}
\sfX^{(p)}_{u}  :=  \alpha(\sfKfv,\sfKp) \bKz
\avsum_{\sfx_i \in \bA^o_u} 
 g^{(p)}_i\big(\xi_{\rm out}(i),\xi_{\rm out}(i+p)\big)\,,
\end{equation}
is such that 
\begin{equation}\label{def:BL-conv}
\limsup_{\vep \to 0} \sfd_{\rm BL} (\{\sfX^{(\vep)}_{u,z}\}_u,\{\sfX^{(p)}_{u}\}_u) \le \ga_p \,.
\end{equation}
In particular,  \eqref{def:BL-conv}  implies that $p \mapsto \{\sfX^{(p)}_u\}_u$ is Cauchy 
in $(M_1(\bR^\sfM),\sfd_{\rm BL})$,  hence converging to some $\{\sfX_u\}_u$,  which 
is also the 
limit in law of $\{\sfX^{(\vep)}_{u,z}\}_u$.

\subsection{Continuity of the Field: Proof of Proposition \ref{th-heightcontinuity}}
\label{sec-contheight}
Set $\bar r:=r_{\sfKo-1}$,  $\wt r = r_{\sfKo-2}$ and $\wh{\bF}_\vep := \bF_{n_\vep} \cup \bar{\bF}_{\bar n_\vep}$.
With $\sfKo$ fixed,  upon adjustment of $\ga$ and $\ze$,  
it suffices for \eqref{hc.2a} to show that for fixed $\sfv$ and $\de=\de(\sfKo,\ga,\ze)$ small, 
 \begin{align}
\limsup_{\vep \to 0} \Pbm \Big(\sup_{\stackrel{\sfx,\sfy \in \wh{\bF}_\vep \cap \bB_\sfd(\sfv,\bar r)}{\sfd (\sfx, \sfy) \le \de  }}
\|(\Yxz,\Vxz)
-(\Yyz,\Vyz)\|_\infty \ge 6 \ga \Big) &\leq 5 \ze \,.
\label{hc.2ab}
\end{align}
For definiteness,  we fix $\sfv={\bf 0}$ and with $\dlc$ fixed and $\sfKz$ large,  
have $\bB_\sfd({\bf 0}, \wt r e)$ within a single city.  Further reduce $\de$ so 
$\bB_\sfd(\sfx,\bar r) \subseteq \bB_\sfd({\bf 0}, \wt r)$ for any
$\sfx$ we consider.  Setting
$N_{\sfx}:=\TxKoz$,  $\wt N_{\sfx}:=\TxKo2z$ and 
\begin{equation}
  D_{\vep}(c) : =\{ \wt \sfs (-c) \leq N_{\sfx}, \wt N_{\sfx}\leq \wt \sfs (c),\quad \forall 
  \sfx \in \wh{\bF}_{\vep}\}, \label{ff.0n}
\end{equation}
for $\wt \sfs (y):=n_\vep (2 n_\vep - \log n_\vep + y)$,  we then have the following a priori estimate.
\begin{lemma}\label{lem:D-c}
Fix $z$,  $\dlc$, $\sfKz$.  There exist $c_\ze<\infty$ such that for $\vep$ small enough,
\begin{equation}
\Pbm (D_\vep (c_\ze))\geq 1-\ze \,.
\label{ff.0}
\end{equation}
\end{lemma}
\begin{proof} Let $S_k$ denote the time to complete $k$ inter-city excursions.  In proving \Cref{lem:conc} 
we saw that $\frac{1}{\sqrt{k}} (S_k - k \sft_1)$ converges in law,  as $k \to \infty$, 
to a normal variable of zero mean and finite variance.
Recalling also \eqref{dfn:bar-svep} and \eqref{clt.2} that 
\[
\bar \sfs_\vep(z) = \pi \sfk_\vep(z) \sft_1 
 \in 
[\wt \sfs(2z-1),\wt \sfs(2z+2c+1)] \, ,
\]
we get that for some $\delta_b =\delta_b(z,\dlc,\sfKz) \to 0$ as $b\to\infty$ and all $\vep \le \vep_o(z,b)$,
\begin{equation}
 \Pbm \Big(  \wt{\sfs}(-b) \leq \pi S_{\sfk_{\vep}(z)}\leq \wt \sfs(b) \Big) \geq 1-\delta_b \,. \label{ff.0p}
\end{equation}
Further,  by the same reasoning as in \Cref{theo-etrt},  its conclusion applies for $k$-level excursions,
at any fixed $k$,  which in view of \eqref{ff.0p} results with \eqref{ff.0} holding for all $c_\ze$ large enough.
\end{proof}
By \eqref{17.01ad} we see that on $D_\vep:=D_\vep(c_\ze)$ of \Cref{lem:D-c},  we have that  
\begin{equation}
\label{field.9}
|\Yxz - \Yyz|
 = \sqrt{2} |\sqrt{N_\sfx} - \sqrt{N_\sfy}| = \sqrt{2} \frac{|N_\sfx-N_\sfy|}{\sqrt{N_\sfx}+\sqrt{N_\sfy}}
 \le \frac{1}{n_\vep} |N_\sfx-N_\sfy| \,.
\end{equation}
Recalling \eqref{def:V-eps-alt} and setting 
\begin{equation}\label{def:Nx-th}
N_{\sfv}^{\th} := \sum_{j=1}^{N_\sfv}  \cos(\theta_{j,\sfv}).
\end{equation}
we also have that 
\begin{align}
\label{17.01adp}
\Vxz - \Vyz = \frac{\bar c}{2} \sqrt{2 \bar \sfs_\vep(0)}
\Big[\frac{N^\th_\sfy}{N_\sfy}  - \frac{N_\sfx^\th}{N_\sfx}  \Big]\,.
\end{align}
Clearly,  $|N_\sfy^\th| \le N_\sfy$.  Further,  on $D_\vep$ 
we have that $N_\sfx \ge \bar c \sqrt{\bar \sfs_\vep(0)}\,  n_\vep$,  hence 
\begin{equation}
\label{field.9a}
|\Vxz - \Vyz| \le \frac{N_\sfx}{n_\vep}
\Big| \frac{N^\th_\sfy}{N_\sfy}  - \frac{N_\sfx^\th}{N_\sfx}  \Big| \le 
\frac{1}{n_\vep} | N^\th_\sfy  - N_\sfx^\th | + \frac{1}{n_\vep} |N_\sfy-N_\sfx| \,.
\end{equation}
It thus suffices for \eqref{hc.2ab} to show that 
 \begin{align}\label{hc.2aa}
\limsup_{\vep \to 0} \Pbm \Big(\sup_{\stackrel{\sfx,\sfy \in \wh{\bF}_{\vep} \cap \bB_\sfd({\bf 0},\bar r)}{\sfd(\sfx, \sfy) \le \de  }}
|N_\sfx-N_\sfy|  \ge 3 \ga n_\vep;  D_\vep \Big) &\leq 2\ze \,,\\
\limsup_{\vep \to 0} \Pbm \Big(\sup_{\stackrel{\sfx,\sfy \in \wh{\bF}_{\vep} \cap \bB_\sfd({\bf 0},\bar r)}{\sfd(\sfx, \sfy) \le \de  }}
|N^\th_\sfx-N^\th_\sfy|  \ge 3 \ga n_\vep;  D_\vep \Big) &\leq 2\ze \,.
\label{ac.2aa}
\end{align}
With our sub-lattices $\bF_j$ of size $|\bF_j| \le c \,  e^{2j}$ (see \eqref{dr.2}), 
being $\wt{c} \,  e^{-j}$ nets for the union of all blocks in $\bT^2$,  we get \eqref{hc.2aa} 
and \eqref{ac.2aa} out of the following.
\bl\label{lem-kolmo}  With $\bar n_\vep=2 \log (n_\vep \log n_\vep)$ and $\wh{\bF}_\vep=\bF_{n_\vep}
 \cup \bar{\bF}_{\bar{n}_\vep}$,
\begin{align}\label{height.70a}
\Pbm \Big(\sup_{\stackrel{\sfx,\sfy \in \wh{\bF}_{\vep} \cap \bB_\sfd({\bf 0},\bar r)}{2 \sfd(\sfx, \sfy) \le 
e^{-\bar n_\vep}}}
|N_\sfx-N_\sfy|  \ge \ga n_\vep;  D_\vep \Big) & = o_\vep(1) \,.
\end{align}
For some $\rho>0$,  $C_\ga<\infty$ any $\vep \le \vep_o$ and $\delta \le \delta_o$,  
\begin{equation}
\Pbm \Big(\sup_{\stackrel{\sfx,\sfy \in \bar \bF_{\bar n_\vep} \cap \bB_\sfd({\bf 0},\bar r)}{\sfd(\sfx, \sfy) \le 2 \de  }}
|N_\sfx-N_\sfy|  \ge \ga n_\vep;  D_\vep \Big) \leq C_\ga (2 \de)^\rho \,.
\label{height.74a}
\end{equation}
Moreover,  both \eqref{height.70a} and \eqref{height.74a} hold with $N_\sfx^\th-N_\sfy^\th$ instead of $N_\sfx-N_\sfy$.
\el
Indeed,  
set $\de<\de_o$ with $C_\ga (2 \de)^\rho \le \ze$ and reduce $\vep_o$ so 
$e^{-\bar n_\vep} \le \de$ 
with the $o_\vep(1)$ in \eqref{height.70a} being at most $\ze/2$.  For any $\sfx \in \wh{\bF}_{\vep} \cap \bB_\sfd({\bf 0},\bar r)$, 
there exists some $\sfx' \in \bar \bF_{\bar n_\vep} \cap \bB_\sfd({\bf 0},\bar r)$ with
$\sfd(\sfx,\sfx') \le  e^{-\bar n_\vep}/2$.  The same applies for the element $\sfy'$ near 
$\sfy$.  By the triangle inequality 
$\sfd(\sfx',\sfy') \le \de + e^{-\bar n_\vep} \le 2 \de$.  As
\[
|N_\sfx-N_\sfy|\leq |N_\sfx-N_{\sfx'}|+|N_{\sfx'}-N_{\sfy'}|+|N_{\sfy'}-N_{\sfy}|,
\]
we get \eqref{hc.2aa} by considering $\sfx',\sfy'$ in \eqref{height.74a} and pairs
$\sfx,\sfx'$ and $\sfy,\sfy'$ in \eqref{height.70a}.  Further,  applying the same reasoning 
for $N_\sfx^\th-N_\sfy^\th$,  yields \eqref{ac.2aa}.

\begin{proof}[Proof of Lemma \ref{lem-kolmo}] Recall that $\wt r =\bar r e$ and
$\bB_\sfd({\bf 0}, \wt  r e)$ lies inside one city,  with $\bB_\sfd(\sfx,\bar r) \subseteq \bB_\sfd({\bf 0}, \wt r)$ for any $\sfx$ we are to consider.  
As in \cite[Section 5]{BRZ2},  let $A_{\sfx,i}$ and $B_{\sfy,i}$ count traversals from 
$\partial \bB_\sfd(\sfx,\bar r)$ to $\partial \bB_\sfd(\sfx,\bar r/e)$,  respectively from
$\partial \bB_\sfd(\sfy,\bar r)$ to $\partial \bB_\sfd(\sfy,\bar r/e)$,
during the $i$-th excursion from $\partial \bB_\sfd({\bf 0},{\wt r})$ to
$\partial \bB_\sfd({\bf 0},  e {\wt r})$,  so
\begin{equation}
N_{\sfx}-N_{\sfy}=\sum_{i=1}^{\wt  N_{\bf 0}}\( A_{\sfx,i} -B_{\sfy,i}\)\,. \label{f1.2}
\end{equation}
By a union over pairs $\sfx,\sfy$ of distance at most $e^{-\bar n_\vep}$,  
the sign of $N_{\sfx}-N_{\sfy}$ and the $2 c_\ze n_\vep$ possible values of $\wt N_{\bf 0}$ 
under $D_\vep$,  it suffices for \eqref{height.70a},  to show that 
 \begin{align}\label{height.70a28}
n_\vep e^{4 n_\vep} e^{- 2 \bar n_\vep} 
\sup_{
\stackrel{\sfx,\sfy \in \bB_\sfd({\bf 0},\bar r),  2 \sfd(\sfx, \sfy) \le e^{-\bar n_\vep}}{\wt \sfs (-c_\ze) \le \ell \le 
\wt \sfs(c_\ze)}}
\Pbm \Big(  \sum_{i=1}^{\ell} (A_{\sfx,i} -B_{\sfy,i})   \ge \ga n_\vep \Big)  = o_\vep(1) \,.
\end{align}
Now,  if $2 \sfd(\sfx,\sfy) \le e^{-\bar n_\vep}$,  then by Chebychev's exponential inequality
and the bound \eqref{basicexp-ad} at $\lambda=5/\ga$,  we have for $\ell \le 2 n_\vep^2$
\[
\Pbm \Big(  \sum_{i=1}^{\ell} (A_{\sfx,i} -B_{\sfy,i})   \ge \ga n_\vep \Big)
\le e^{-5 n_\vep} e^{c_3 e^{4 \pi\lambda} e^{-\bar n_\vep} n_\vep^2} \,,
\]
from which \eqref{height.70a28} follows,  as $e^{-\bar n_\vep} n_\vep \to 0$.

As in the proof of Kolmogorov-Centsov continuity theorem (e.g.  \cite[p. 118,  solution of Problem 2.2.9]{KaS},  at
$d=2$,  $\alpha=14$,  $\beta=1$),  we get \eqref{height.74a} for some $\rho<1/14$ upon 
showing that for $C=C(\bar r)<\infty$,  all $\vep \le \vep_o$ and $\sfx,\sfy\in \bar \bF_{\bar n_\vep}\cap  \bB_\sfd(\npp,\bar r)$ with $\sfd(\sfx,\sfy) \le \delta_o(\bar r)$,  
\begin{equation}
\bE \Big[ |N_\sfx-N_\sfy|^{14};  D_\vep \Big] \le C n_\vep^{14}  \sfd(\sfx,\sfy)^{3} \,.
\label{kolmo.120-new}
\end{equation}
To this end,  by \eqref{f1.2} we have on $D_\vep$ that 
\begin{equation}
    \label{eq-rage1}
    N_\sfx - N_\sfy \leq 
    \sum_{i=1}^{\wt \sfs (-c_\ze)} \( A_{\sfx,i} -B_{\sfy,i}\)+
    \sum_{i=\wt \sfs(-c_\ze) +1}^{\wt \sfs(c_\ze)} |A_{\sfx,i}-B_{\sfy,i}|=: \AA_1+\AA_2 
  \end{equation}
(the absolute values in $\AA_2$ are needed because we have 
in \eqref{f1.2} some
$[\wt \sfs(-c_\ze),\wt \sfs(c_\ze)]$-valued random number of summands).  It thus suffices for \eqref{kolmo.120-new} to show that 
\[
\bE [ \AA_1^{14} ] \le C n_\vep^{14} \sfd(\sfx,\sfy)^{3}\,,  \quad 
\bE [ \AA_2^{14} ] \le C n_\vep^{14} \sfd(\sfx,\sfy)^{3}\,.
\]
For $\AA_1$ consider \eqref{eq-zav1-new} at $\ell=\wt \sfs(-c_\ze) \le 2 n_\vep^2$ and
$\lambda=n_\vep^{-1} \sfd(\sfx,\sfy)^{-3/14}$ which for
$\vep \le \vep_o$ and $\sfx \ne \sfy \in \bar \bF_{\bar n_\vep}$ is arbitrarily small,  to get as needed that 
\be
\label{eq-rage2}
\frac{\bE [\AA_1^{14}]}{14!} \le \frac{\bE [ \cosh(\lambda \AA_1)]}{\lambda^{14}} \le \lambda^{-14} c^{-1} e^{c \, \sfd(\sfx,\sfy) \lambda^2 n_\vep^2} \le 
C
\, n_\vep^{14} \,  \sfd(\sfx,\sfy)^3 \,.
\ee
For $\AA_2$ note that by  \eqref{basicexp-ad},  for $\lambda \in [0,1]$, 
\be\label{eq-zav2}
\bE[e^{\lambda \AA_2}] \le e^{2 c_\zeta c_3  \lambda \sfd(\sfx,\sfy) n_\vep} \,.
\ee
Considering as before $\lambda=n_\vep^{-1} \sfd(\sfx,\sfy)^{-3/14}$,  and having $\AA_2 \ge 0$,  we 
find that 
\[
\frac{\bE[\AA_2^{14}]}{14!} \le \frac{\bE[e^{\lambda \AA_2}]}{\lambda^{14}} \le  
e^{2 c_\zeta c_3  \delta_o^{11/14}} \,  n_\vep^{14} \sfd(\sfx,\sfy)^{3} \,,
\]
as needed to complete the proof of \eqref{kolmo.120-new} and thereby of \eqref{height.74a}.
Turning to the bounds for $N^\th_\sfx-N^\th_\sfy$,  we have analogously to \eqref{f1.2},  that
\begin{align}
N^\th_{\sfx}-N^\th_{\sfy}&=\sum_{i=1}^{\wt  N_{\bf 0}}( A^\th_{\sfx,i} -B^\th_{\sfy,i}),\label{f1.2-new}
\end{align}
for $A^{\th}_{\sfx,i}$ and $B^\th_{\sfy,i}$ of \eqref{angle.12}.  Reducing $\delta_o$ as needed,  
the same reasoning yields the bounds \eqref{height.70a} 
and \eqref{height.74a} with $N_\sfx^\th-N_\sfy^\th$ instead of $N_\sfx-N_\sfy$ 
(note that we made sure to have $e^{-{\bar n_\vep}/{2}} n_\vep \to 0$ and 
$\sqrt{\sfd(\sfx,\sfy)} \lambda^2 n_\vep^2$ uniformly bounded). 
\end{proof}

\subsection{Gaussian limits}\label{sec-gauss}

\begin{proof}[Proof of \Cref{theo-mclt}]
By the Cramer-Wold device,  it suffices for \Cref{theo-mclt} to fix 
$\{\sfx_i,  i \le p\}$,  $\{\la_i,  i \le 2p\} \in \bR^{2p}$ and prove the \abbr{clt}-s for 
\[
\OO^{(\vep)}_{\rm in} :=
\sum_{i=1}^{2p} \lambda_i q^{(\vep)}_{\rm in} (i) \quad {\rm and} \quad 
\OO^{(\vep)}_{\rm out}:=
\sum_{i=1}^{2p} \lambda_i q^{(\vep)}_{\rm out} (i) \,.
\]
To this end,  note that conditionally on $\HHk$ the $\bR$-valued $\OO^{(\vep)}_{\rm in}$ is of the form
\begin{equation}\label{def:O-vep-in}
\OO^{(\vep)}_{\rm in} = \frac{1}{\sqrt{\sfk_\vep(z)}} \sum_{j=1}^{\sfk_\vep(z)}
\big(  \OO_j  - \Ebm[ \OO |\sfz_{(j-1)2},\sfz_{j1}] \big)  \,,
\end{equation}
where the conditionally independent $\OO_j\deq \OO(\sfz,\sfz')$ count the contribution 
to $\OO^{(\vep)}_{\rm in}$ between the $(j-1)$'st and $j$'th inter-city excursions,  starting at $\sfz \in \bTs$ and 
ending at $\sfz' \in \bTs^\sfe$.  Indeed,  with $\sqrt{\sfk_\vep(z)}=\sfa (\sfm_\vep+z)$,
$\sfa:={1}/{\sqrt{2 \pi \sft_1}}$,  we have 
\begin{equation}\label{def:OO-z}
\OO(\sfz,\sfz') \deq
\sum_{i=1}^p \sum_{\ell=1}^{T_{\sfx_i,\sfKo}^{\bTs^\sfe(1)}}
\psi_i(\th_{\ell,\sfx_i})\,,
\end{equation}
aggregating during such single inter-city excursion,  over all $\sfKo$-level excursions of $\{\sfx_i,  i \le p \}$ 
with angle increments $\{\th_{\ell,\sfx_i}\}$,  the uniformly bounded scores 
\[
\psi_{i}(\th):=\sfa \Big\{\la_i + \la_{i+p}  \bar c[\langle \cos(\cdot),\bar \mu \rangle - \cos (\th)]\Big\}\,.
\]
The uniform separation between $\bigcup \{ \bB_{\sfd} (\sfv, r_{\sfKo}) :  \sfv \in \bA_u^o,  u \le \sfM \}$,
$\bTs$ and $\bTs^{\sfe}$ implies a uniform fourth moment control  
\[
\sup_{\sfz,\sfz',\sfv} \Ebm \Big[\big(T_{\sfv,\sfKo}^{\bTs^\sfe(1)}\big)^4\, \big|\,  \sfz,\sfz' \Big] < \infty \,,
\]
and thereby the same uniform control on $\sup_{\sfz,\sfz'} \Ebm [\OO(\sfz,\sfz')^4]$.  We thus deduce by
Lyapunov's \abbr{clt},  that as $\vep \to 0$,  the variables $\bsig_{\sfk_\vep(z)}^{-1} \OO^{(\vep)}_{\rm in}$ converge
in law conditionally on $\HHk$,  to a standard normal,  where 
\[
\bsig_k^2 := \frac{1}{k} \sum_{j=1}^{k} f(\sfz_{(j-1)2},\sfz_{j1}) \quad {\rm for} \quad 
f(\sfz,\sfz') := {\sf var}\big(\OO(\sfz,\sfz')|\sfz,\sfz'\big) \,.
\] 
Our fourth moment control implies that $f(\cdot,\cdot)$ is uniformly bounded.  Thus,  when $k \to \infty$, 
a.s.  $\bsig_k^2 \to \langle f,  \bar \mu_\sfe \rangle$ by the \abbr{lln} for integrable functions of such 
Doeblin chains (see \cite[Theorem 17.0.1]{Meyn-Tweedie}).  Replacing the random $\bsig_{\sfk_\vep(z)}$
by the non-random $\sqrt{\langle f,  \bar \mu_\sfe \rangle}$ thus concludes the proof of the conditional 
\abbr{clt} for $\OO^{(\vep)}_{\rm in}$ and consequently,  of \eqref{mCLT.03}.
Next,  note that similarly to \eqref{def:O-vep-in},   
\begin{equation}\label{def:O-vep-out}
\OO^{(\vep)}_{\rm out} = \wt{\OO}_{\sfk_\vep(z)} \quad {\rm for} \quad 
\wt{\OO}_k := \frac{1}{\sqrt{k}} \sum_{j=1}^{k}
\big(  \Ebm[ \OO |\sfz_{(j-1)2},\sfz_{j1}] - \Ebm_{\bar \mu_\sfe} [\OO] \big)  \,.
\end{equation}
As $\OO(\sfz,\sfz')$ of \eqref{def:OO-z} have uniformly bounded moments,  
$\Ebm[ \OO |\sfz,\sfz'] - \Ebm_{\bar \mu_\sfe} [\OO]$ is uniformly bounded,  
of zero $\bar \mu_\sfe$-mean.  Consequently,  $\{\wt{\OO}_k\}$ satisfies the \abbr{clt} (see 
\cite[Theorem 17.0.1]{Meyn-Tweedie}),  with some finite variance 
(of a form analogous to \eqref{eq:var-erg}).  In particular,  having established the \abbr{CLT}
for $\OO^{(\vep)}_{\rm out}$ yields \eqref{mCLT.04},  thereby completing the proof 
of the proposition.
\end{proof}
\begin{proof}[Proof of \Cref{lem-expi}] The tower property of conditional expectations allows us to
replace the stated expectations under the invariant law $\bar \mu_{\sfe}$ on $\bTs \times \bTs^\sfe$
by such expectations with respect to (only) the invariant law $\mu_\sfe$ on $\bTs$.  Having done so, 
for \eqref{17.02b} set $R=r_{\sfKo-1}$,  $r=r_{\sfKo}$ and denote by $D_1$ the time  
to go from $\partial \bB_\sfd (\sfx,R)$ to  $\partial \bB_\sfd (\sfx,r)$ and back to $\partial \bB_\sfd(\sfx,R)$.
Recall \cite[Lemma 8.1]{BK} that 
starting at the unique invariant measure $\mu_{\sfx,r}^R$ for such $\sfKo$-level excursions,  yields
$\Ebm_{\mu_{\sfx,r}^R}[D_1] ={1}/{\pi}$.
Further,  by an ergodic consideration (similar to the argument used in proving \Cref{lem-equalmean} below),  it takes time
${(k+o(k)) \sft_1}/{\Ebm_{\mu_{\sfe}} \big[T_{\sfx,\sfKo}^{\bTs^\sfe(1)} \big]}$ till 
completing the $k$-th $\sfKo$-level excursion around $\sfx$,  which on the other hand 
must also take $(k+o(k)) \Ebm_{\mu_{\sfx,r}^R}[D_1]$ time.  Equating the leading 
orders of both expressions yields \eqref{17.02b}.  Further,  recall that $\bar \mu$ 
of the paragraph containing \eqref{dfn:Poiss},  is the invariant law of the angle-differences $\{\th_{j,\sfx}\}$.  As it
must match the law induced on these angle-differences
by starting at the invariant measure $\mu_\sfe$ on $\bTs$,  we get \eqref{17.02c} out of \eqref{17.02b}.
\end{proof}
\section{Coupling and continuity estimates,  and proof of \Cref{lem:4.58m}}\label{sec:cont}

\Cref{sec:coupling-BM} provides a sharp coupling between the first (unconditional) 
inner Brownian excursions,  and such excursions conditioned upon certain specified 
number of outer excursions.  This coupling, which may be of independent interest, 
was used in the proof of \Cref{lem-sec3.3}.
It is supplemented in  \Cref{sec-continest} 
by general continuity estimates for excursion counts and angle-difference fields, 
which among other things we use in \Cref{sec:4.58m},  to prove \Cref{lem:4.58m}.

\subsection{Coupling conditional Brownian excursions}\label{sec:coupling-BM}
For $k \ge 1$ and the Brownian motion on $\bT^2$ starting at some 
point outside $\bB_{\sfd}({\bf 0}, e R)$, $R \le e^{-2}$, let 
$\{L^{(k)}_{i}, i \ge 1\}$ count the \textit{inner} excursions, 
each starting at $\partial \bB_\sfd({\bf 0},R e^{-k})$ and ending upon exiting 
$\bB_\sfd({\bf 0},R e^{-(k-1)})$, during the consecutive \textit{outer} (Brownian) excursions
from $\partial \bB_\sfd({\bf 0},R)$ till exiting $\bB_\sfd({\bf 0},e R)$. 

Having $S^{(k)}_\ell$ inner excursions completed during the first $\ell \ge 1$ outer ones, for any
integer $m \ge 1$, $b \ge a$ we denote by $\bP_{m,a,b}$ the joint law 
of the first $S^{(k)}_m$ inner excursions, conditional on the event 
\begin{equation}\label{eq:Amb-def}
\AA^{(k)}_{m,a,b} :=\{ S^{(k)}_m \in [m+a,m+b] \}\,.
\end{equation}
Let $\bP^\star_{m'}$ denote the unconditional joint law of the first $m'$ 
\textit{inner} excursions starting from $\partial \bB_\sfd({\bf 0},R e^{-k})$ to
$\partial \bB_\sfd({\bf 0},R e^{-(k-1)})$.

Recall that starting anywhere on $\partial \bB_\sfd({\bf 0},R)$, 
the order in which a Brownian path hits 
$\{\partial \bB_\sfd({\bf 0},R e^{-(l-1)}), l \ge 1 \}$, up to exiting $\bB_\sfd({\bf 0},e R)$, 
forms a \abbr{srw} on $\bZ_+$ which starts at $l=1$ and is absorbed at $l=0$.  In particular, 
$\{L^{(k)}_i, i \ge 1\}$ are \abbr{iid} copies of the law $L^{(k)}$ as in \Cref{lem:coupling}.
Utilizing that lemma,  we
next establish a
coupling of the un-ordered sets  of inner excursions under the laws $\bP^\star_{m'}$ and $\bP_{m,a,b}$, so that
 the number $\Delta_{m,a,b}$ of excursions that appear in one set but not the other is whp 
\blue{of order $(\log k \sqrt{m}) \vee |b|$.}
\begin{proposition}\label{prop:coupling} 
For some universal $\bar c>1$,
any $\ka>0$ small enough, integer $m' \in [m+a,m+b]$, 
$b-a \le {\sqrt{m}}/{\ka}$ and $(m,k,r,b) \in \Theta_\ka$ of
\eqref{dfn:Thetak}, there exists a coupling, possibly depending on $m',m,a,b,k,r,R$,
of the two (un-ordered) collections of inner excursions
distributed according to $\bP_{m,a,b}$ and $\bP^\star_{m'}$,
where the number $\Delta_{m,a,b}$ of differing inner excursions is such that 
\begin{align}\label{eq:coupling}
\sup_{\Theta_\ka,a,m'} \, \Big\{ \frac{e^r}{k} \bP \big( \Delta_{m,a,b} \ge 7  r {\bar c} \ka \sqrt{m} \big) \Big\} < \infty \,.
\end{align}
\end{proposition}
\begin{proof} 
We prove below Proposition \ref{prop:coupling} only for $a=b$ (and $m'=m+b$), adopting 
the notations $\AA_{m,b}$, $\bP_{m,b}$ and $\Delta_{m,b}$ for $\AA^{(k)}_{m,b,b}$,
$\bP_{m,b,b}$ and $\Delta_{m,b,b}$, respectively. Indeed, 
since \eqref{eq:coupling} trivially holds for $r \le r_o$ and any $r_o(\ka)$ finite,  we may replace
the stated constraint on $b-a$ by $b-a \le r \ka \sqrt{m}$. Further,  
$\bP_{m,a,b}$ is a mixture of the laws $\bP_{m,u}$, $u \in [a,b]$, and  also $\bP^\star_{m_1}$ and $\bP^\star_{m_2}$ can be trivially coupled so that they share $m_1\wedge m_2$ identical excursions.
It follows that weighting the coupling between $\bP_{m,u}$ and $\bP^\star_{m+u}$ 
as in Proposition \ref{prop:coupling}, according to these mixture coefficients, increases the 
number of different excursions with respect to $\bP^\star_{m'}$, for any fixed $m' \in [m+a,m+b]$,
by at most $b-a \le r \ka \sqrt{m}$. Thus, assuming the proposition for $a=b$ and $m'=m+b$, multiplying $\bar c$ by $8/7$ yields the general case in \eqref{eq:coupling}, as stated.

Setting hereafter $a=b$, $m'=m+b$,
we simplify the task of proving Proposition \ref{prop:coupling} upon noting that 
our inner Brownian excursions, under either  $\bP_{m,b}$ or $\bP^\star_{m}$, are independent conditional 
on the pairs $\sfy=(\sfy^{(1)},\sfy^{(2)})$, where
$\sfy^{(1)} \in \partial \bB_\sfd({\bf 0},R e^{-k})$  
and $\sfy^{(2)} \in \partial \bB_\sfd({\bf 0},R e^{-(k-1)})$ denote 
the starting and ending positions, respectively, of each inner excursion, 
with the conditional law of each such inner excursion determined,  in an identical way under both  $\bP_{m,b}$ or $\bP^\star_{m}$,
by the value of the corresponding 
pair $\sfy$. It thus suffices to provide a coupling between the joint laws $\bP_{m,b}$ and $\bP^\star_{m'}$ 
of the corresponding two collections $\{\sfy_{i} \}$, such that
the bound \eqref{eq:coupling} applies for the number $\Delta_{m,b}$ 
of differences between these two un-ordered lists. 

To this end, note that  
the consecutive starting and ending positions $\sfy^{(o)}_i=(\sfy_i^{(o,1)},\sfy_i^{(o,2)})$
of the outer excursions, form under the Brownian law a Markov chain, whose transition probabilities are composed 
of the outer part, yielding the conditional law of $\sfy_{i}^{(o,1)} \in  \partial \bB_\sfd({\bf 0},R)$
given $\sfy_{i-1}^{(o,2)} \in \partial \bB_\sfd({\bf 0}, e R)$ and the inner 
part, namely the conditional law of $\sfy_i^{(o,2)}$ given $\sfy_i^{(o,1)}$. By
the same reasoning as in the proof of  Lemma \ref{lem:conc},
the conditional law of $\sfy_i^{(o,1)}$ given $\sfy_{i-1}^{(o,2)}$ 
has a uniformly in $R \le e^{-2}$, bounded away from zero density with respect to 
Lebesgue's measure on the corresponding circle. The Markov chain $\{\sfy_i^{(o)}\}$
is thus Doeblin aperiodic (satisfying 
\cite[(16.9)]{Meyn-Tweedie}), admitting infinitely many renewal times $\{T_j, j \ge 1\}$,
where $\sfy^{(o,1)}$ is uniformly distributed, independently of the past,  
with \abbr{iid} positive Geometric ($1/c'$) inter-renewal distances 
\[
\chi_j := T_{j+1} - T_j,  \qquad \forall j \ge 0, 
\]
for some universal $c'>1$ (where $T_0=0$, so $\chi_0=T_1$).

In what follows, we refer to the law of the Brownian path as the \textit{unconditional} law.
Fixing $k \ge 1$, the unconditioned counts $\{L^{(k)}_i, i \ge 1\}$ of inner excursions completed 
during the consecutive outer excursions,
are \abbr{iid} and independent 
of the sequence $\{\chi_j\}_{j \ge 0}$.
For any $j \ge 0$, there are
\[
 \bar L^{(k)}_j := \sum_{i=0}^{\chi_j-1} L_{T_j+i}^{(k)}
\]
inner excursions between the $j$-th and $(j+1)$-st renewal events, to which correspond under $\bP^\star$ 
the mutually independent segments $\bar \YY_j^{(k)} := \{\sfy_i : i \in [T_j,T_{j+1})\}$
between the consecutive renewal events. Per $k \ge 1$,
the law of each $\bar \YY_j^{(k)}$, $j \ne 0$, depends only on the value of 
$\LL^{(k)}_j := \{\chi_j, L_{T_j+i}^{(k)}, i\in [0,\chi_j)\}$. Thus, we first construct 
the blocks $\LL^{(k)}_j$ out of  the
independent \abbr{iid} $\{L^{(k)}_i, i \ge 1\}$ and $\{\chi_j, j \ge 0\}$. Thereafter,
given these blocks
we construct the conditionally independent segments $\bar \YY^{(k)}_j$.
During the first $\ell$ outer excursions there are
\begin{equation}\label{eq:M-ell-def}
M_\ell := \sup \{ j \ge 0 : T_j \le \ell \} 
\end{equation}
renewal events, and we note that for any fixed $j \in [1,M_m)$, the law of $\bar \YY_j^{(k)}$ is 
affected by conditioning on $\AA_{m,b}$ only via the change 
in the law of $\LL^{(k)}_j$.

The measure $\bP_{m,b}$ on the collection $\{\sfy_i\}$ is then induced by the  
law $\bQ_{m,b}$ of $\{\wt{L}_i, i \le m\}$, namely the \abbr{iid} $\{L_i=L^{(k)}_i, i \le m\}$ 
conditional on the event $\AA_{m,b}=\AA^{(k)}_{m,b,b}$ 
of \eqref{eq:Amb-def} (and hereafter we often suppress the index $k$ in order to alleviate the notations). 
The stopping time $\tau_h$ of \eqref{def:tau-h} counts for such $L^{(k)}$
the number of outer excursions up to the completion of $h$ inner excursions, 
and the measure $\bP^\star_{m'}$ on the unconditional $\{\sfy_i, i \le m'\}$ is thus likewise
induced by the law
$\bQ^\star_{m'}$ of \Cref{lem:coupling}.

In view of the preceding discussion we note that for any given matched
block $\wt{\LL}_j=\LL^\star_j$,  we can and shall construct \emph{the same Brownian 
path between the corresponding pairs of renewal events} in our two sequences, 
thereby matching the relevant (ordered) segments $\bar \YY$ of the $\sfy_i$ values. Turning 
to construct (couple) the values of $\sfy_i$ in the collection $\bar \UU_{m,b}$ of 
non-matched blocks of outer excursions, note that by definition, $\bar \UU_{m,b}$ induces
a collection $\UU_{m,b}$ of 
$\bar \Delta_{m,b}$ positive summands $\{\wt{L}_{(\alpha)}, L^\star_{(\alpha^\star)}\}$, such that 
due to our choice of $m'=m+b$,
\[
\sum_\alpha \wt{L}_{(\alpha)} = \sum_{\alpha^\star} L^\star_{(\alpha^\star)} 
\]
(compare \eqref{eq:Amb-def} for $a=b$, with \eqref{eq:Amb-star-def}). 
We then further partition the positive integer entries of $\UU_{m,b}$ to yield, 
after re-ordering, at most $\bar \Delta_{m,b}$ parts satisfying 
\begin{equation}\label{eq:unmatch-piece}
\wh{L}_{(\alpha)}=\wh{L}^\star_{(\alpha^\star)} \ge 1 \,.
\end{equation}
Indeed, taking initially $\wh{L}_{(\alpha)}=\wt{L}_{(\alpha)}$ and $\wh{L}^\star_{(\alpha^\star)}=
L^\star_{(\alpha^\star)}$, we iteratively pair the $\arg\min_\alpha \{ \wh{L}_{(\alpha)} \}$ with 
$\arg\min_{\alpha^\star} \{ \wh{L}^\star_{(\alpha^\star)} \}$,
removing from consideration the common (matched) parts, after cutting as needed 
the larger of these two length.
With each step reducing the total number of yet unmatched parts of elements from $\UU_{m,b}$, 
we have the stated complete matching after at most $\bar \Delta_{m,b}$ operations. 
Now, let $R e^{-k} e^{i \sfz_j}$, 
for $0 \le j \le \ell$, denote the starting positions of the $\ell$ inner excursions 
which form such generic piece of an outer excursion, say $(\alpha)$, 
consisting of $\wh{L}_{(\alpha)}=\ell$ inner excursions, with 
$R e^{-k} e^{i \sfz^\star_j}$ for the starting positions in the corresponding paired
piece $(\alpha^\star)$ of some outer excursion in our other sequence. Conditional on the totality of 
$\{\wt{L}_i,L^\star_i,\chi_j\}$ and in particular on the identity and length of  
all matched parts appearing in \eqref{eq:unmatch-piece},  by radial symmetry of the 
Brownian motion within $\bB_\sfd({\bf 0},e R)$, the collections $\{ \sfz_{j} - \sfz_{j-1}, 1 \le j \le \ell \}$
and $\{\sfz^\star_j-\sfz^\star_{j-1}, 1 \le j \le \ell \}$, are each composed of 
\abbr{iid} $[2 \pi]$-valued variables, of density which is uniformly over $k \ge 1$ and $R \le e^{-2}$, bounded
below by some $1/c'' > 0$. Thus, even after the (paired) collection 
$\{\wh{L}_{(\alpha)}, \wh{L}^\star_{(\alpha^\star)} \}$ is specified, there exist \abbr{IID} Bernoulli($1/c''$) 
variables associated with the starting positions $\sfy_i^{(1)}$ of the unmatched inner excursions 
corresponding to $\UU_{m,b}$, 
except the first one in each piece,
whose occurrence signify the renewal event whereby the current 
$\sfy^{(1)}_i$ will be chosen uniformly, independently of the past. The corresponding 
distances $\wh{\chi}_{(\alpha)}$ and $\wh{\chi}^\star_{(\alpha^\star)}$ from
the start of a generic piece $(\alpha)$ or $(\alpha^\star)$, respectively, till the 
first such renewal, are \abbr{iid} positive Geometric($1/c''$). Further, after all such 
$\{\wh{\chi}_{(\alpha)}, \wh{\chi}^\star_{(\alpha^\star)}\}$ are chosen, the 
Brownian segments within the different unmatched pieces of both sequences, 
starting at the first renewal event in each piece, are independent of each other, with 
the law of any such segment determined solely by the remaining length within that 
piece, namely by $\wh{L}_{(\alpha)}-\wh{\chi}_{(\alpha)}$ and  
$\wh{L}^\star_{(\alpha^\star)}-\wh{\chi}^\star_{(\alpha^\star)}$, respectively. It thus follows
that the coupling $\wh{\chi}_{(\alpha)}=\wh{\chi}^\star_{(\alpha^\star)}$ for each pair
$(\alpha)$ and $(\alpha^\star)$ of \eqref{eq:unmatch-piece}, allows one  
to use the same Brownian path for each such pair, except for at most 
the first $\wh{\chi}_{(\alpha)}$ inner excursions within that piece. Thereby we arrive at 
\[
\Delta_{m,b} \le \sum_{(\alpha)} \wh{\chi}_{(\alpha)} \,,
\]
where the sum is over at most $\bar \Delta_{m,b}$ indices $(\alpha)$, conditional on which 
the variables $\{ \wh{\chi}_{(\alpha)} \}$ remain \abbr{iid} Geometric($1/c''$). 
In conclusion,  the bound of 
Lemma \ref{lem:coupling} on $\bar \Delta_{m,b}$, together with standard exponential tails for
sum of \abbr{iid} Geometric variables, complete the proof of Proposition \ref{prop:coupling}.
\end{proof}

\subsection{Continuity estimates}\label{sec-continest}

Supplementing \Cref{prop:coupling} by arguments adapted from \cite[Section 5]{BRZ2}
(where we now take throughout $h(r)=r$),  we 
establish in Lemmas \ref{lem-cont1} and \ref{lem-cont3} general continuity results for excursion counts,
building on the continuity results of \Cref{cor:HigherMomentBound} and  \Cref{lem:c3new},  which
also apply for the angle-difference fields in \Cref{th-heightcontinuity}.
Specifically,  fixing some point $\npp \in \bT^2$ which we identify with ${\bf 0} \in \bR^2$,  and 
constants $0<a<b<1$,  we assume analogously to \cite[(5.1)]{BRZ2} that 
\begin{equation}
  \label{eq-ratiobounds}
0<r<R<\tilde{r}<\tilde{R} < \frac{1}{2} \qquad {\rm with} \qquad  
  {r}/{R}, {R}/{\tilde r},{\tilde{r}}/{\tilde{R}}\in [a,b]\,,
\end{equation}
with all the results in the remainder of this section uniform in the choice of $\npp$ and 
such $r,R,\tilde{r},\tilde{R}$.  Similarly to \cite[(5.4)-(5.5)]{BRZ2},  we denote by 
$\wt R_{\npp,i}$ and $\wt D_{\npp,i}$ 
the successive arrival times of the Brownian motion $\BMS_{\cdot}$ to
$\partial \bB_{\sfd}(\npp,\tilde{r})$ and 
$\partial \bB_{\sfd}(\npp,\tilde{R})$ respectively, 
with 
\[
0 = \wt D_{\npp,0} < \wt R_{\npp,1}<\wt D_{\npp,1}<\wt R_{\npp,2}<\wt D_{\npp,2}<\ldots
\]
and $\wt X^{i}:=\BMS_{\left(\wt R_{\npp,i}+\cdot\right)\wedge \wt D_{\npp,i}}$
being the $i$-th excursion from $\partial \bB_{\sfd}(\npp,\tilde{r})$ to
$\partial \bB_{\sfd}(\npp,\tilde{R})$.  For $i\ge 1$ fixed and generic $\sfx \in \bB_\sfd (\npp,\tilde{r}-R)$ 
we define an $\sfx$-excursion to be an excursion from $\partial \bB_{\sfd}(\sfx,r)$ to $\partial \bB_{\sfd}(\sfx,R)$
during the excursion $\wt X^{i}$ and let $A_{\sfx,i}$ count the number of such $\sfx$-excursions.  This is also the 
number of traversals from $\partial \bB_{\sfd}(\sfx,R)$ to $\partial \bB_{\sfd}(\sfx,r)$ during $\wt{X}^{i}$,
as $A_{\sfx,i}$ is defined in \cite[Section 5]{BRZ2}.  Similarly,  for 
$\sfy \in \bB_\sfd (\npp,\tilde{r}-R)$  we let $B_{\sfy,i}$ count the number of $\sfy$-excursions during 
$\wt{X}^i$.   We also set the angle-difference counts
\begin{align}
A^{\th}_{\sfx,i}&=\sum_{j}
\cos \( \th_{j'}\), \hspace{.2 in} B^{\th}_{\sfy,i}=\sum_{j'}
\cos \( \th'_{j'}\),\label{angle.12}
\end{align} 
summing over  the $A_{\sfx,i}$ and $B_{\sfy,i}$,  respectively,   
$\Ipi$-valued angle-differences between the starting and the exit points of the 
$\sfx$-excursions,  and  respectively, $\sfy$-excursions,   
within 
the excursion $\wt X^i$.
Equipped with these notations and building on the proof of 
\cite[Lemma 5.2]{BRZ2},  we have the following exponential tail estimate 
(on which we rely in our continuity estimates).
\bl
\label{lem:HigherMomentBound}
Denote by $\bE^{\sfz,\sfv}(\cdot)$ the expectation over the excursion $\wt{X}^1$ which 
starts at $\sfz \in\partial \bB_\sfd(\npp,\tilde{r})$,  conditional on ending at 
$\sfv \in \partial \bB_{\sfd}(\npp,\tilde{R})$.  
For some $q_{1}=q_{1}(a,b)<\frac{1}{2}$,   
$c_3=c_3(a,b)<\ff$ and any $\sfx,\sfy \in \bB_\sfd (\npp,\tilde{r}-R)$ with $q:=\sfd(\sfx,\sfy)/R \leq q_{1}$,
if $\lambda>0$ and $e^{4 \pi \lambda} < 1/(q c_3)$,  then 
\be
\label{eq-zav1}
 \sup_{\sfz,\sfv} \bE^{\sfz,\sfv} \Big( \exp(\lambda |A_{\sfx,1}-B_{\sfy,1}|) \Big) \leq 
 \exp\big (c_3 q (e^{4 \pi \lambda}-1)\big) \,.
\ee 
With $q=\sqrt{\sfd(\sfx,\sfy)/R}$,  the same bound applies for $|A^\th_{\sfx,1}-B^\th_{\sfy,1}|$.
\el
\begin{proof} 
Fixing $\sfx$ and $\sfy$,  consider the following subsets of $\bB_\sfd(\npp,\tilde{r})$:
\begin{align*}
\Gamma_1& := \bB_\sfd(\sfx,r)\cap \bB_\sfd(\sfy,r)\,,\quad\;
\Gamma_2:= \bB_\sfd(\sfx,r)\cup \bB_\sfd(\sfy, r)\,, \\
\Gamma_3&:= \bB_\sfd(\sfx,R)\cap \bB_\sfd(\sfy,R)\,,\quad
\Gamma_4:= \bB_\sfd(\sfx,R)\cup \bB_\sfd(\sfy,R)\,,
\end{align*}
 \begin{figure}[h]
\begin{center}
\begin{tikzpicture}
\draw[red] (-.1,0) circle(1.5);
\draw[red]  (-.1,0) circle(3);
  \fill (-.1,0) circle (1pt) [red] node[red, left] {$\red{\sfx}$};
   \fill (-1.5,0)  node[blue, right] {$\Gamma_1$};
    \fill (-2.4,0)  node[green, right] {$\partial \Gamma_2$};
     \fill (-3.9,0)  node[purple, right] {$\partial \Gamma_4$};
      \fill (-3,0)  node[orange, right] {$\Gamma_3$};
  \draw (.1,0) circle(1.5);
\draw (.1,0) circle(3);
  \fill (.1,0) circle (1pt) node[right] {$\sfy$};
   \fill (0.8,0)  node[blue, right] {$\Gamma_1$};
    \fill (1.5,0)  node[green, right] {$\partial \Gamma_2$};
     \fill (2.3,0)  node[orange, right] {$\Gamma_3$};
      \fill (3,0)  node[purple, right] {$\partial \Gamma_4$};
\end{tikzpicture}
\end{center}
\caption{\small Boundary of $\textcolor{blue}{\Gamma_1},  \textcolor{green}{\Gamma_2}, 
\textcolor{orange}{\Gamma_3} ,  \textcolor{purple}{\Gamma_4}$ 
($R=2r$,  $q=\frac{1}{15}$,  \red{$\sfx$-balls in red}).}
\label{angle1}
\end{figure}
noting that $q_1<(1-b)/2$ guarantees that $\Gamma_2 \subseteq \Gamma_3$ (see Figure \ref{angle1}). 
Thus we have that 
\[
\Gamma_1 \subset \Gamma_2 \subset \Gamma_3 \subset \Gamma_4  \,.
\]
In addition to the
$\sfx$-excursions from $\partial \bB_\sfd(\sfx,r)$ to $\partial \bB_\sfd(\sfx,R)$ 
and the 
$\sfy$-excursions from $\partial \bB_\sfd(\sfy,r)$ to $\partial \bB_\sfd(\sfy,R)$ 
within our excursion around $\npp$,  we consider also the consecutive \emph{$(\sfx,\sfy)$-excursions} defined
to be the excursions of the Brownian path from $\partial \Gamma_2$ to $\partial \Gamma_4$,   and note that each 
$(\sfx,\sfy)$-excursion must contain at least one $\sfx$ or $\sfy$-excursion (see also Figure \ref{angle1}).
We call an $(\sfx,\sfy)$-excursion \emph{special} if the path enters $\partial \Gamma_1$ before $\partial \Gamma_3$ and do not visit $\Gamma_2$ after entering $\partial \Gamma_3$.   Note the following elementary relations:\\
$\bullet$ Any $\sfx$-excursion contains a path from $\partial \Gamma_2$ 
to $\partial \Gamma_3$;\\
$\bullet$ Any $\sfy$-excursion contains a path from $\partial \Gamma_2$ to $\partial \Gamma_3$;\\
$\bullet$ A path from $\partial \Gamma_1$ to $\partial \Gamma_4$ contains 
an  
$\sfx$-excursion and  $\sfy$-excursion;\\
$\bullet$ In a special $(\sfx,\sfy)$-excursion lies exactly one $\sfx$-excursion and 
one $\sfy$-excursion. \\ 

If the $j$-th $(\sfx,\sfy)$-excursion is special,  set $C_j:=0$ 
and $C^\th_j:=|\cos(\th)-\cos(\th')|$ for the angle-differences $\th$ and $\th'$ of the $\sfx$ and $\sfy$ 
excursions inside this $(\sfx,\sfy)$-excursion.  If the $j$-th $(\sfx,\sfy)$-excursion is non-special, 
set  $C^\th_j:=0$ and $C_j \ge 1$ be the number of excursions from $\partial \Gamma_2$ 
to $\partial \Gamma_3$ during this non-special $(\sfx,\sfy)$-excursion.  It is then easy to verify that
\begin{equation}\label{eq:geo-bd-ad}
|A_{\sfx,1}-B_{\sfy,1}| \le 2 \sum_j C_j \,,  \qquad |A^\th_{\sfx,1}-B^\th_{\sfy,1}| \le \sum_j (2 C_j + C^\th_j)  \,.
\end{equation}
Let $\GJo_p$ denote a \emph{non-negative} Geometric($1-p$) variable (so $\Pbm(\GJo_{p} \ge k)=  p^k$ for $k \ge 0$) 
and $\BB_p$ a Bernoulli($p$) variable independent of it.  We proceed to show that for some 
$q_1(a,b)>0$,  $c=c(a,b)<\ff$ and all $q=\sfd(\sfx,\sfy)/R \le q_1$,  conditionally on the 
path up to the start of its $j$-th $(\sfx,\sfy)$-excursion,  we have the stochastic domination
\begin{align}\label{eq:st-dom}
C_j & \stackrel{\rm st}{\preccurlyeq} 2 \GJo_{cq}\,,  \\ 
\label{eq:st-dom-th}
2 C_j + C^\th_j & \stackrel{\rm st}{\preccurlyeq} 4 \GJo_{cq} + c \sqrt{q} + 4 \pi \BB_{c \sqrt{q}} \,.
\end{align}
(The stochastic domination notation $X\stackrel{\rm st}{\preccurlyeq}Y$ means that the random variables $X$ and $Y$ can be coupled so that $X\leq Y$ a.s., or equivalently, that $\bP(X\geq x)\leq \bP(Y\geq x)$ for all real $x$.)
To this end,  set $\sfw:=\frac{1}{2}(\sfx+\sfy)$ and balls 
$\Gamma^\star_\ell:= \bB_{\sfd}(\sfw,\rho_\ell)$,  of radii 
\[
\rho_1:= r-\frac{1}{2} q R,  \quad \rho_2= r+\frac{1}{2}qR,  \quad \rho_3= R-\frac{1}{2}qR,  \quad \rho_4 = R+\frac{1}{2}qR \,.
\]
As $\sfd(\sfx,\sfy)=qR$,  it follows that 
\[
\Gamma^\star_1 \subset \Gamma_1 \subset \Gamma_2 \subset \Gamma^\star_2 \subset \Gamma^\star_3
\subset \Gamma_3 \subset \Gamma_4 \subset \Gamma^\star_4 \,.
\]
With
$\{C_j=1\}$ requiring that the $(\sfx,\sfy)$-excursion enter $\partial \Gamma_3$ before $\partial \Gamma_1$,
for any $\sfu \in \partial \Gamma_2$ we bound $\Pbm^\sfu(C_j=1)$ by the maximum over starting points in $\partial \Gamma_2^\star$ of the 
probability of entering $\partial \Gamma_3^\star$ before $\partial \Gamma^\star_1$,  which is 
\begin{equation}\label{eq:Cj-not-zero}
\gamma_1 := \frac{\log\Big(\frac{2 r+qR}{2r-qR}\Big)}{\log\Big(\frac{2R-qR}{2r-qR}\Big)} \,.
\end{equation}
Likewise,  the event $\{C_j \ge k\}$ requires at least $k-1 \ge 1$ returns 
from $\partial \Gamma_3$ to $\partial \Gamma_2$ prior to hitting $\partial \Gamma_4$.
Starting anywhere on $\partial \Gamma_3$ the probability of hitting $\partial \Gamma_2$ 
prior to $\partial \Gamma_4$ is at most the maximum over starting points in $\partial \Gamma_3^\star$ 
of the 
probability of entering $\partial \Gamma^\star_2$ prior to 
$\partial \Gamma^\star_4$,  which is 
\begin{equation}\label{eq:Cj-above-one}
\gamma_2:= \frac{\log\Big(\frac{2 R+qR}{2R-qR}\Big)}{\log\Big(\frac{2R+qR}{2r+qR}\Big)} \,.
\end{equation}
Clearly $\gamma_1 + \gamma_2 \le c q$ for any $q=\sfd(\sfx,\sfy)/R \le q_1<a$ and 
some $c=c(a,b,q_1)$ finite such that 
$c(a,b,q_1) q_1 \to 0$ when $q_1 \to 0$.  Assuming \abbr{wlog} that $c q_1 < 1$
and utilizing the strong Markov property of the Brownian path,  we thus arrive 
at $\bP(C_j \ge k) \le (cq)^{k-{\bf 1}_{k \ge 2}}$.  To get  \eqref{eq:st-dom} simply note
that $k-1 \ge k/2$ for $k \ge 2$.

Turning to \eqref{eq:st-dom-th},  recall that for any special $(\sfx,\sfy)$-excursion, 
$\th=\th_{ +}-\th_{-}$ where $\th_{-}$ the angle at which it first hits $\partial \bB_\sfd(\sfx,r)$ 
and $\th_{+}$ the angle 
at which it first hits $\partial \bB_\sfd(\sfx,R)$,  while analogously $\th'=\th'_{+}-\th'_{-}$ for 
the inner $\sfy$-excursion,  so with $\|\cos(\cdot)\|_{\rm Lip} \le 1$,  clearly 
\[
C^\th_j \le (|\th'_{-}-\th_{-}| + |\th'_+ - \th_+|) {\bf 1}_{\big\{\hbox{excursion}\; j \; \hbox{is special}\big\}} \,.
\]
Focusing first on $|\th_{-}-\th'_{-}|$,  we assume \abbr{wlog} that our special $(\sfx,\sfy)$-excursion 
starts at $\sfu \in \partial \bB_\sfd(\sfx,r)$ and construct a Cartesian coordinate system such that 
one axis contains the segment connecting $\sfu$ to $\sfy$ and the origin lies in $\partial \bB_\sfd(\sfy,r)$,  
see Figure \ref{angle},  with $\sfu$ represented 
as $(0, u)$ for some $0<u \le qR$ and $\sfy=(0, -r)$.  Note that
$\partial \bB_\sfd(\sfy,r)$ 
intersects the line $L$ of height $-u$ at $(\pm z  ,-u)$ for $z:=\sqrt{r^2-(r-u)^2}$.
 \begin{figure}[h]
\begin{center}
\begin{tikzpicture}
\draw (0,-2) circle(2);
\draw (0,-4.5) -- (0,1.3);
\draw (-2,0) -- (2,0);
  \fill (0,-2) circle (2pt) node[right] {$\sfy=(0,-r)$};
  \fill (0,.5) circle (2pt) node[right] {$\sfu=(0,u)$};
  \fill (1.33,-.5) circle (2pt) node[below] {$\!\!\!\!\!\!\!\!\!\!\!\! (z,-u)$};
  \fill (-1.33,-.5) circle (2pt) node[below] {$\;\;\;\;\;\;\;\;\;\;\; (-z,-u)$};
\draw (-3,-.5) -- (3,-.5);
\end{tikzpicture}
\end{center}
\caption{\small First hitting point}
\label{angle}
\end{figure}
Further measuring \abbr{wlog} the angles $\th_-$ and $\th'_{-}$ from $\sfy=(0, -r)$ we have that $\th_{-}=0$ 
with $\th'_{-}$ the angle from $(0,-r)$ of the point where the Brownian path,  started at $(0, u)$,  
first hits $\partial \bB_\sfd\left(\sfy,r\right)$.  Denote by $\tau$ the first hitting time 
of $L$ by a Brownian path $(X_t,  t \ge 0)$ started at $X_0=(0, u)$, 
with $X_{\tau}:=(x_{\tau}, -u)$.  Because we deal with a special excursion,  necessarily $\tau$ is smaller
than the hitting time of $\partial \Gamma_3$.  Note that if 
$|x_{\tau}|< z$ then the Brownian path $X_t$ must first hits $\partial \bB_\sfd(\sfy,r)$ 
on the upper arc of $\partial \bB_\sfd(\sfy,r)$ between $(z,-u)$ and $(-z,-u)$.  In that case 
$|\th'_{-}|$ is at most the angle of the point $(z,-u)$ measured from $\sfy=(0, -r)$.  That is, 
\[
|\th'_{-}| \le \tan^{ -1}\Big(\frac{z}{r-u}\Big) \le \tan^{-1} \Big(\sqrt{(1-qR/r)^{-2}-1}\Big) \le \frac{c}{2} \sqrt{q} \,,
\]
for some $c=c(a,q_1)<\ff$ and any $q \le q_1$.   
Consequently,  
\[
|\th'_{-}-\th_-| 
\leq \frac{c}{2} \sqrt{q} + 2\pi {\bf 1}_{\{   |x_{\tau}|\geq z\}}.
\]
 Further,  by the Poisson kernel for the upper half-plane
 \begin{align*}
 \mathbb{P}^{ \(   0,u \)} ( |x_{\tau}|\geq z ) 
 =  \frac{2}{\pi}\int_{z/(2u)}^{ \ff} \frac{1}{s^2+1}\,ds
  &= \frac{2}{\pi} \tan^{-1} \Big(\frac{2u}{z}\Big) 
  \le \frac{2}{\pi} \tan^{-1} \Big(\frac{2\sqrt{qR}}{\sqrt{2r-qR}}\Big)
\leq \frac{c}{2} \sqrt{q},  \nn
 \end{align*}
for any $q \le q_1$ and some $c=c(a,q_1) < \ff$ such that $c \sqrt{q_1} <1$.  
Omitting the similar analysis for $|\theta_+-\theta'_{+}|$,  we conclude that 
conditionally of the path up to the start of its $j$-th $(\sfx,\sfy)$-excursion,  the 
variable $C^\th_j$ is stochastically dominated by $c\sqrt{q}+4 \pi \BB_{c\sqrt{q}}$.  
Further,  $C^\th_j>0$ only for special excursions,  whereby $C_j=0$,  allowing us to take 
$\BB_{c \sqrt{q}}$ independent of $\GJo_{cq}$ in \eqref{eq:st-dom-th}.

Having proved \eqref{eq:st-dom} and \eqref{eq:st-dom-th} for each $(\sfx,\sfy)$-excursion,  
we still need to control the number $D$ of $(\sfx,\sfy)$-excursions within the excursion $\wt X^1$
in order to apply \eqref{eq:geo-bd-ad}.  To this end,  
note that to have $\{D \ge k\}$
requires at least $k$ visits from 
$\partial \Gamma_4 \subset \bB_\sfd(\npp,\tilde{r})$ to $\partial \Gamma_2$ before exiting 
$\bB_{\sfd}(\npp,\tilde{R}) \subset \bB_\sfd (\sfw,\rho_5):=\Gamma^\star_5$ where $\rho_5:=\tilde{R}+\tilde{r}$.
Similarly to our proof of \eqref{eq:st-dom},  each of these $k \ge 1$ visits requires hitting 
$\partial \Gamma^\star_2$ before $\partial \Gamma^\star_5$ and the
probability of doing so is at most the maximum of the latter probability 
over starting points in $\partial \Gamma^\star_3 \subset \Gamma_4$,  which is 
\begin{equation}\label{eq:-J-par}
\gamma_3 := \frac{\log\Big(\frac{2\tilde{R}+2\tilde{r}}{2R-qR}\Big)}{\log\Big(\frac{2\tilde{R}+2\tilde{r}}{2r+qR}\Big)}
\le p_o \,,
\end{equation}
for some $p_o=p_o(a,b) < 1$ and all $q \le q_1$.  
Consequently,  the number of $(\sfx,\sfy)$-excursions is stochastically dominated by $\GJo_{p_o}$,  
independent of the content of those excursions.  

Now,  for $\lambda>0$ and $3 c q e^{4 \lambda} \le 1$,
\begin{equation}\label{eq:bd-psi}
\bE\big[e^{4\lambda \GJo_{cq}}\big] = (1-c q) \sum_{k=0}^\infty (c q)^k e^{4 \lambda k} = 
\frac{1- c q}{1- c q e^{4 \lambda}} \le \exp\big(3 c q (e^{4\lambda}-1)\big)
\end{equation}
(as $e^{2t} (1-t)$ increase on $[0,\frac{1}{2}]$,   take here
$t=\frac{cq(e^{4\lambda}-1)}{1-cq}$).  Similarly,  setting $\psi_o:=\log(\frac{1+p_0}{2p_0})>0$
and $c_o:=\frac{1-p_0}{1+p_0}>0$,  we have that for any $\psi \in [0, \psi_o]$, 
\begin{equation}\label{eq:bd-J-po}
\bE[e^{\psi \GJo_{p_o}}] = \frac{1-p_0}{1-p_0 e^\psi} \le e^{\psi/c_o}\,.
\end{equation}
Further,  with $\log(1+t) \le t$,  we have for any $\lambda>0$,
\begin{equation}\label{eq:bd-psi-th}
\log \bE\big[e^{4\pi \lambda \BB_{c \sqrt{q}}}\big] + c \sqrt{q} \lambda 
 \le 2 c \sqrt{q} (e^{4 \pi \lambda} - 1)\,.
\end{equation}
We thus deduce from \eqref{eq:bd-psi} that for $\lambda>0$ with $e^{4 \lambda} \le \psi_o/(3 c q)$,
\[
\psi(\lambda) := \log \bE\big[e^{4\lambda \GJo_{cq}}\big] \le 3 c q (e^{4 \lambda}-1) \le \psi_o \,,
\]
and likewise we have in view of \eqref{eq:bd-psi-th} that when $e^{4\pi \lambda} \le \psi_o/(4 c \sqrt{q})$,
\[
\wt{\psi}(\lambda) := \psi(\lambda) + \log \bE\Big[e^{4\pi \lambda \BB_{c \sqrt{q}}}\Big] + c \sqrt{q} \lambda
\le 4 c \sqrt{q} (e^{4 \pi \lambda} - 1) \le \psi_o \,.
\]
Hence,  for such values of $\lambda$,  
 by \eqref{eq:st-dom}--\eqref{eq:st-dom-th} and the strong Markov property of the Brownian path,  
\begin{align}\label{eq:pre-zav1}
\bE^\sfz\big[\exp(\lambda |A_{\sfx,1}-B_{\sfy,1}|)\big] &\le \bE \big[e^{\psi(\lambda) \GJo_{p_o}}\big] \le 
e^{\psi(\lambda)/c_o} \le e^{(3c/c_o) q (e^{4 \lambda}-1)} \,,\\
\bE^\sfz\big[\exp(\lambda |A^\th_{\sfx,1}-B^\th_{\sfy,1}|)\big] &\le 
\bE \big[ e^{\wt \psi(\lambda) \GJo_{p_o}} \big] \le e^{\wt \psi(\lambda)/c_o} \le e^{(4c/c_o) \sqrt{q} (e^{4 \pi \lambda} - 1)}
\label{eq:pre-zav1-th}
\end{align}
(utilizing on the right hand side both \eqref{eq:bd-J-po} 
and our bounds on $\psi(\lambda)$ and $\wt \psi(\lambda)$).

On the \abbr{lhs} of \eqref{eq:pre-zav1} and \eqref{eq:pre-zav1-th} we have the \abbr{mgf} of some variable $Y \ge 0$ measurable 
on the part of $\wt X^1$ inside $\bB_{\sfd}(\npp,\tilde{r})$.
Similarly to \eqref{eq10.1},
restricted to $\sigma(Y)$,  the Radon-Nikodym derivative $d\Pbm^{\sfz,\sfv}/d\Pbm^{\sfz}$
due to conditioning on the end point of $\wt{X}^1$ is uniformly bounded by 
twice the maximum over $\sfz' \in \partial \bB(\npp,\sqrt{\tilde r \tilde R})$ of the Poisson kernel 
$p(\sfz',\sfv)$ for the exit measure of $\bB(\npp,\tilde{R})$, 
divided by its minimum over such $\sfz'$.  With $c_1=c_1(a,b)>1$ the maximum value of this factor 
for $\tilde{r}/\tilde{R} \in [a,b]$,  we deduce that for any such $Y$ and $\lambda>0$,
\[
\bE^{\sfz,\sfv}[e^{\lambda Y}] \le 1 + c_1 (\bE^{\sfz}[e^{\lambda Y}] -1) \le \Big(\bE^{\sfz}[e^{\lambda Y}]\Big)^{c_1}
\] 
(with the right hand side by the monotonicity in $x>0$ of $(e^{b x}-1)/x$,  when $b>0$).  Using this bound we absorb
the factor $c_1$ in the universal constant $c/c_o$ on the right hand side of \eqref{eq:pre-zav1}-\eqref{eq:pre-zav1-th},
arriving at \eqref{eq-zav1} and its analog for $|A^\th_{\sfx,1}-B^\th_{\sfy,1}|$.
\end{proof}

Next,  we state and prove a consequence of \Cref{lem:HigherMomentBound},  which will  take here the role of
\cite[Cor.  5.3]{BRZ2}.
\begin{corollary}\label{cor:HigherMomentBound}
Denote by $\HHs$ the $\sigma$-algebra of the Brownian motion excluding its
excursions $\{\wt X^{i},i \ge 1\}$ from $\partial \bB_{\sfd}(\npp,\tilde{r})$ to
$\partial \bB_{\sfd}(\npp,\tilde{R})$.  In the setting of \Cref{lem:HigherMomentBound},  for
any $N \ge 1$,  non-random $\{i_1,\ldots, i_N\}$ and $e^{4 \pi \lambda} \le 1/(q c_3)$,  a.s.
\be
\label{basicexp-ad}
\Ebm \Big( \exp(\lambda \sum_{\ell=1}^N |A_{\sfx,i_\ell}-B_{\sfy,i_\ell}|) \big| \HHs \Big) 
\leq e^{c_3 q (e^{4 \pi \lambda}-1) N} \,.
\ee 
Adjusting $c_3$, 
one can bound the right hand side
for $\lambda \in [0,1]$ 
by $e^{c_3 q \lambda N} \le (1+2 c_3 q \lambda)^N$ (reducing $q_1 \le 1/(2 c_3)$ if needed).
In particular,  for any $\varphi \ge 2 c_3 q N$,   a.s.
\begin{equation}\label{basicexp2}
\Pbm \Big(\sum_{\ell=1}^{N} |A_{\sfx,i_\ell}-B_{\sfy,i_\ell}| \ge \varphi \big| \HHs \Big) 
\leq  e^{-\frac{1}{2} (\varphi \wedge N) }  \, .
\end{equation}
If conditionally on $\HHs$,  the variables $\{\wt X^i\}$ are independent of a Geometric($p_2$) 
variable $\GJ \ge 1$,  then 
for any $\lambda\in [0, {p_{2}}/{3}]$,
\begin{equation}
\bE \Big(\exp\big(\lambda\sum_{i=1}^{\GJ} |A_{\sfx,i}-B_{\sfy,i}|\big)\Big)\leq e^{\frac{4 c_{3}}{p_2} q \lambda} \, . 
\label{basicexp1}
\end{equation}
Upon adjusting $q$ as in  \Cref{lem:HigherMomentBound},  
the bounds \eqref{basicexp-ad}-\eqref{basicexp1} apply for 
$|A^\th_{\sfx,i}-B^\th_{\sfy,i}|$.
\end{corollary}

\begin{proof} Conditioning on $\HHs$ is equivalent to knowing the starting and ending points of 
our $N$ excursions,  making them thereby independent copies of the single excursion of 
\Cref{lem:HigherMomentBound}.  Thus,  \eqref{basicexp-ad} is a direct 
consequence of \eqref{eq-zav1}.   
Using the simpler bound $e^{c_3 q N}$ at $\lambda=1$,  yields \eqref{basicexp2} 
by Markov's inequality,  as $\varphi- c_3 q N \ge \frac{1}{2} (\varphi \wedge N)$ for any
$\varphi \ge 2 c_3 q N$.  Similarly,  using the conditional on $\HHs$ 
independence between $\GJ$ and our excursions,  and summing 
the bound $(1+2 c_3 q \lambda)^{\GJ}$ in \eqref{basicexp-ad} 
over the values of $\GJ$,  the \abbr{lhs} of \eqref{basicexp1} is bounded by 
\[
p_2 \sum_{\ell=1}^\infty (1 + 2 c_3 q \lambda)^\ell (1-p_2)^{\ell-1} 
\le \Big(1-\frac{2 c_3 q \lambda
}{p_2}\Big)^{-1} \le e^{\frac{4 c_{3}}{p_2} q \lambda} \,,
\]
provided $2 c_3 q \lambda \le p_2/3$.
\end{proof}

\subsubsection{Interlude: zero-mean by ergodicity}\label{subsec:ergodic}
Recall that for any $\sfv \in \bT^2$ and $0<r<R$,  the starting points of the consecutive excursions 
from $\partial \bB_{\sfd}(\sfv,R)$ to $\partial \bB_{\sfd} (\sfv,r)$ and back to $\partial \bB_{\sfd}(\sfv,R)$
form a Doeblin aperiodic Markov chain on $\partial \bB_{\sfd}(\sfv,R)$,  which thereby admits the unique  
invariant measure,  denoted hereafter by $\mu_{\sfv,r}^{R}$ (see 
\cite[Lemma 2.1]{UenoOnRecurrentMarkovProcesses}).  Starting such 
excursion according to $\mu_{\sfv,r}^R$ induces a specific 
first hitting measure $\mu_{\sfv,R}^r$ on $\partial \bB_{\sfd} (\sfv,r)$ and we note in passing that
while $\sfv \mapsto  \mu_{\sfv,R}^r$ is invariant to translations,  
unlike the Brownian motion on $\bS^2$ of \cite[Section 5]{BRZ2},  here $\mu^r_{\sfv,R}$ is not
the uniform measure $\mu_{\sfv,r}$ on $\partial \bB_{\sfd}(\sfv,r)$.  However,  the density
with respect to $\mu_{\npp,\tilde{r}}$ of the Brownian hitting measure $\Pbm^{\sfu}$ 
of $\partial \bB_\sfd(\npp,\tilde{r})$ is bounded and bounded away from zero, 
uniformly over $\sfu \in \partial \bB_\sfd(\npp,\tilde{R})$.  Hence,  for any
$\mu \in M_1(\partial \bB_\sfd(\npp,\tilde{r}))$ of bounded density with respect to $\mu_{\npp,\tilde{r}}$ 
(and in particular for $\mu_{\npp,\tilde{R}}^{\tilde{r}}$), there exist
$\bar p(\mu,\tilde{r},\tilde{R})>0$ such that for any $p_3 \in (0,\bar p)$ one has
$\Pbm^{\sfu}=p_3 \mu + (1-p_3) \nu_{\sfu}$ with
$\nu_{\sfu} \in M_1(\partial \bB_\sfd(\npp,\tilde{r}))$. This yields in turn 
the corresponding 
renewal times construction of $\{\wt X^i\}$,
as in \cite[(5.28)]{BRZ2}.
\begin{cns}\label{con:renew}
Let $\{I_i\}$ be \abbr{iid} Bernoulli($p_3$) variables,  independent of the Brownian 
excursions $\{\wt{X}^i\}$.  Starting at $J_0=1$ with initial measure $\mu$ for $\wt{X}^1$,  at 
the renewal times $J_j:=\inf\{i>J_{j-1}: I_i=1\}$ set the initial measure for $\wt{X}^{J_j}$ to be $\mu$
independently of $\{\wt{X}^i,  i < J_j\}$,  while each $\wt{X}^{i+1}$ with $I_{i+1}=0$ 
starts at the measure $\nu_{\sfu_{i}}$ for $\sfu_{i} := \wt{X}_\infty^{i} \in \partial \bB_\sfd(\npp,\tilde{R})$. 
\end{cns}
Note that $J_j-J_{j-1}$,  $j \ge 1$  of Construction \ref{con:renew}
are \abbr{IID} strictly positive Geometric($p_3$) variables such that the bound \eqref{basicexp1}
holds for the excursions in each of the blocks $[J_{j-1},J_j)$.  That is,  the bound holds for the \abbr{mgf} of 
\begin{align}\label{eq:def-barV}
\bar V^\th_j :=\sum_{i=J_{j-1}}^{J_{j}-1} (A^\th_{\sfx,i}-B^\th_{\sfy,i}) \,,  \qquad
\bar V_j := \sum_{i=J_{j-1}}^{J_{j}-1} (A_{\sfx,i}-B_{\sfy,i}) 
\end{align}
(with the adjustment to $q$ in case of $\bar V^\th_j$).  The linear in $\lambda$ behavior of the exponent on the right side in \eqref{basicexp1}
is due to the positive mean of $|A_{\sfx,i}-B_{\sfy,i}|$.  As in \cite[(5.30)]{BRZ2},  also 
here it is crucial to get a quadratic in $\lambda$ behavior for the \abbr{cgf} of the variables 
in \eqref{eq:def-barV}.  This relies in turn on these variables having zero mean 
(as in \cite[(5.29)]{BRZ2}),
which we establish in our next lemma, 
without the spherical symmetry on which the proof of  \cite[Lemma 5.4]{BRZ2} relies.
\bl\label{lem-equalmean}
Fixing $\mu$,  $p_3>0$ and $\{J_j\}$ as in Construction \ref{con:renew},  one has that:\\
(a).  For any $\sfx,\sfy \in \bB_\sfd(\npp,\tilde{r}-R)$ and $j \ge 1$,  
\begin{align}
\Ebm^\mu [ \bar V^\th_j ] &= 0 \qquad
\hbox{and} \qquad
\Ebm^\mu [\bar V_j ]  
=0
 \,.\label{4.2}
\end{align}
(b).  Further,  in case $\mu=\mu_{\npp,\tilde{R}}^{\tilde{r}}$,  we have for any $\sfx,\sfy \in \bB_\sfd(\npp,\tilde{r}-R)$ and
$i \ge 1$, 
\begin{align}
\Ebm^\mu [ A^\th_{\sfx,i}-B^\th_{\sfy,i} | \{I_j\} ] &= 0 \qquad
\hbox{and} \qquad
\Ebm^\mu [A_{\sfx,i}-B_{\sfy,i} | \{I_j\} ] =0
 \,.\label{4.2-cond}
\end{align} 
\el
 \begin{proof} (a).  Our ergodic argument uses the Doeblin Markov chains of the starting points 
on $\partial \bB_\sfd(\sfx,r)$ for the excursions to $\partial \bB_\sfd(\sfx,R)$ and back,
the analogous chain for such excursions around $\sfy$ and the \abbr{IID} chain of $\mu$-distributed
starting points for blocks $[J_{j-1},J_j)$ of excursions 
from $\partial \bB_\sfd(\npp,\tilde{r})$ to $\partial \bB_\sfd(\npp,\tilde{R})$,  augmenting the state
of each chain with its excursion (block) length.\\
Since $\bT^2$ and the Brownian semi-group of $X_\cdot$ are translation invariant,  setting 
\begin{equation}\label{dfn:D-th}
D_{\sfx,k}:=\sum_{j=1}^{k} \cos(\th_{j,\sfx})\,, \qquad 
D_{\sfy,k}:=\sum_{j=1}^{k} \cos(\th_{j,\sfy})\,,
\end{equation}
the function $\sfu \mapsto \bE^{\sfx+\sfu}(D_{\sfx,1})$ is independent of $\sfx$ 
and translating $\partial \bB_\sfd(\sfx,\cdot)$ to $\partial \bB_\sfd(\sfy,\cdot)$ 
by adding $\sfy-\sfx$,  maps the invariant measure $\mu_{\sfx,r}^R$ into $\mu_{\sfy,r}^R$.
Thus,  applying Birkhoff's ergodic theorem to the bounded (random) functional $\cos(\th)$ 
and the chains of excursions around $\sfx$ and around $\sfy$,  we have that 
$k^{-1} D_{\sfx,k}$ and $k^{-1} D_{\sfy,k}$ converge a.s.  when 
$k \to \infty$,  to the same non-random limit (namely,  $\langle \cos(\cdot),\bar \mu \rangle$,  for $\bar \mu$ 
whose density is given in \eqref{dfn:Poiss}).  
Further,  with $N_\sfx(t)$ and $N_\sfy(t)$ counting the excursions completed by time $t$, 
from $\partial \bB_\sfd(\sfx,r)$ to $\partial \bB_\sfd(\sfx,R)$,  and from 
$\partial \bB_\sfd(\sfy,r)$ to $\partial \bB_\sfd(\sfy,R)$,  respectively,  
the length of such excursions possess finite exponential moments,  
uniform in the corresponding starting points,  and by the Doeblin property 
of the underlying chain the consecutive excursion lengths form an (exponentially) mixing sequence.  
Thus,  applying Birkhoff's ergodic theorem for the average of the $\sfx$ and $\sfy$ excursion length,  
using again the translation invariance of $\bT^2$ and $X_\cdot$,  we deduce that 
$t^{-1} N_{\sfx}(t)$ and $t^{-1} N_{\sfy}(t)$ have the same non-random limit when $t \to \infty$.  
It follows that 
\[
\frac{1}{t} \Delta_t:=\frac{1}{t} (D_{\sfx,N_\sfx(t)}-D_{\sfy,N_{\sfy(t)}}) =
\frac{D_{\sfx,N_\sfx(t)}}{N_{\sfx}(t)}\frac{N_\sfx(t)}{t}-\frac{D_{\sfy,N_{\sfy(t)}}}{N_{\sfy}(t)}\frac{N_{\sfy}(t)}{t}
\stackrel{a.s}{\to} 0 \,.
\] 
Next,  let  $T_j$ be the times of starting the $J_j$-th excursion from 
$\partial \bB_\sfd(\npp,\tilde{r})$ to $\partial \bB_\sfd(\npp,\tilde{R})$
and $\wt N_{\npp}(t):=\sup\{j \ge 1:T_j \le t\}$.  Then,  by the renewal \abbr{SLLN} for the integrable 
time $T_2-T_1$, 
there exists non-random $\infty>b>0$,  such that for any $\delta>0$
we have 
that $\Pbm(\AA^{\complement}_{t,\delta}) \to 0$ for $\AA_{t,\delta}:=\{|\wt N_{\npp}(t) - bt| \le \delta t-1\}$.  
On  the event $\AA_{t,\delta}$ we have 
that $t \in [T_{(b-\delta)t+1},T_{(b+\delta)t})$,  hence then  
\begin{equation}\label{eq:erg-ad}
|\Delta_{T_{bt}}| \le |\Delta_t| + \sum_{i = J_{(b-\delta) t}}^{J_{(b+\delta) t}-1} 
( A_{\sfx,i}+ B_{\sfy,i} )  
\end{equation}
(the last term counts the $\sfx$-excursions and $\sfy$-excursions in $[T_{(b-\delta)t},T_{(b+\delta)t})$).
With $\sup_{\sfu} \Ebm^{\sfu} A_{\sfx,1} + \sup_{\sfu'} \Ebm^{\sfu'} B_{\sfy,1} \le C$ and 
$\Ebm [J_j-J_{j-1}] = 1/p_3$,
by Markov's inequality the probability that the right-most term of \eqref{eq:erg-ad} 
exceeds $2 \delta t \ell C/p_3$ is at most $1/\ell$.   Clearly,  
\[
\Delta_{T_{bt}}=\sum_{j=1}^{bt} \bar V^\th_j \,,
\]
so dividing \eqref{eq:erg-ad} by $b t$,  taking $t \to \infty$,  then $\delta \to 0$ followed by $\ell \to \infty$,  we
conclude that a.s.
\begin{equation}\label{eq:erg-ad2}
\lim_{t \to \infty} \frac{1}{bt} \sum_{i=1}^{bt} \bar V^\th_j = 0 \,.
\end{equation}
By the \abbr{SLLN},  the limit in \eqref{eq:erg-ad2} must match the corresponding mean 
$\Ebm^\mu \bar V_j^\th$,  yielding the \abbr{lhs} of \eqref{4.2}. \\ 
Now,   instead of \eqref{dfn:D-th} set $D_{\sfx,k}=D_{\sfy,k}:=k$,  resulting with $\Delta_t = N_\sfx(t)-N_\sfy(t)$ 
and thereby with $\Delta_{T_{bt}} =\sum_{j \le bt} \bar V_j$.  Thus,  upon applying in this setting 
the same reasoning as above,  we arrive at the right hand side of \eqref{4.2}.\\
(b).  Fixing $i \ge 1$,  the variables in \eqref{4.2-cond} are measurable on $\wt X^i$,  hence their conditional 
expectations depend on $\{I_j\}$ only via the initial law of $\wt X^i$,  which in turn depends only on $\{I_j, j<i\}$.
Further,  as explained \cite[above (5.29)]{BRZ2},  in the invariant case $\mu=\mu_{\npp,\tilde{R}}^{\tilde{r}}$,  
one has that $\Ebm^\mu [ \nu_{\sfu_i} | \{I_j,j<i\} ] = \mu$,  so the conditional expectations 
in \eqref{4.2-cond} are independent of $\{I_j\}$ and $i \ge 1$.  In particular,  then 
\begin{align*}
\Ebm^\mu [\bar V^\th_1] &= (\Ebm J_1 -1)  \Ebm^\mu [ A^\th_{\sfx,1}-B^\th_{\sfy,1} ] \,,\\
\Ebm^\mu [\bar V_1] &= (\Ebm J_1-1)  \Ebm^\mu [ A_{\sfx,1}-B_{\sfy,1} ] \,,
\end{align*} 
with \eqref{4.2-cond} an immediate consequence of \eqref{4.2}.
\end{proof}

\subsubsection{Centered tails and sharper continuity bounds}\label{subsec:refined-cont}
It is not hard to verify that the Brownian hitting measures 
$\Pbm^\sfu$ of $\partial \bB_\sfd(\npp,\tilde{r})$ are such that for some $p_3=p_3(a,b)>0$, 
\begin{equation}\label{eq:unif-p3}
\inf_{\stackrel{\sfu,\sfu' \in \partial \bB_\sfd(\npp,\tilde{R})}{A \subset \partial \bB_\sfd(\npp,\tilde{r})}}
\Big\{ \Pbm^{\sfu}(A) - p_3 \Pbm^{\sfu'}(A) \Big\} \ge 0 \,,  \quad  \forall \tilde{r}/\tilde{R} \in [a,b],  \tilde{R} < \frac{1}{2}, 
\end{equation}
so Construction \ref{con:renew} applies for $\mu=\mu^{\tilde{r}}_{\npp,\tilde{R}}$ and 
$p_3=p_3(a,b)$.
As we show next,  using this construction as in \cite[Lemma 5.5]{BRZ2} and 
combining \Cref{lem-equalmean} with \eqref{basicexp1} improves the bound on the \abbr{mgf} of the 
variables in \eqref{eq:def-barV} and hence on the \abbr{mgf} of
\begin{equation}\label{eq:def-V}
V_\ell:= \sum_{i=1}^\ell \big(A_{\sfx,i}-B_{\sfy,i}\big)  \quad \hbox{and} \quad 
V^\th_\ell:= \sum_{i=1}^\ell \big(A^\th_{\sfx,i}-B^\th_{\sfy,i}\big) \,.
\end{equation}
We note in passing that the bounds of \eqref{eq-zav1-new} are utilized 
in  \Cref{lem-kolmo} 
(within a chaining argument),  to yield a uniform 
modulus of continuity both for $N_\sfx-N_\sfy$ of \eqref{17.01ad} and for the 
angle averages $N^\th_\sfx-N^\th_\sfy$ of \eqref{def:Nx-th}.
\bl\label{lem:c3new}
For some $q_{1}, p_3>0$ and  $c_4<\infty$,  depending only on $a, b$,  
and $\sfx,\sfy \in \bB_\sfd (\npp,\tilde{r}-R)$ with $q:=\sfd(\sfx,\sfy)/R \leq q_{1}$,  one  has that:\\
(a).  For $\mu=\mu_{\npp,\tilde{R}}^{\tilde{r}}$ and $\bar V_j$ of \eqref{eq:def-barV},
\begin{align}\label{basicexp1-imp}
\bE^\mu \Big[\exp\big(\lambda \sum_{j=1}^\ell \bar V_j \big)\Big] &\leq \exp\big(c_{4} q \lambda^2 \ell\big) \,,
\qquad \forall \ell \ge 1, \; |\lambda| \le \frac{p_3}{4}, \\
\label{eq-zav1-new}
\bE^{\mu}\big[\exp(\lambda V_\ell) \big] &\le \exp\big(2 c_4 q \lambda^2 \ell\big)  \,, \qquad \forall \ell \ge 1,   \;
|\lambda| \le \frac{p_3}{8}\,.
\end{align}
Up to the adjustment of $q$ these bounds apply for $\sum_j \bar V^\th_j$ and $V_\ell^\th$,  respectively.\\
(b).   Further,  for any $\tilde{R}/\wh R \in [a,b]$ and $\wh R < \frac{1}{2}$,   
\begin{equation}
\Ebm^{\mu_{\npp,\tilde{R}}}\big[ \exp(\la V_\ell) \,\big| \, 
\trav{1}{\npp}{\wh{R}}{\tilde{R}}{\npp}{\tilde{R}}{\tilde{r}}
 = \ell \big] \leq \exp(2 c_4 q \lambda^{2} \ell)\,,
\quad \forall \ell \ge 1,  |\la| \le \frac{p_3}{8} \,,
\label{ergtotal.3}
\end{equation}
and the same applies for $V_\ell^\th$,  upon adjusting $q$.
\el
\begin{proof} By symmetry,  it suffices to consider $\lambda>0$ as we do hereafter. \\
(a).  By  the construction of our renewal times,  $\{\bar V_j\}$ are \abbr{iid} 
(and so are $\{\bar V_j^\th\}$),  hence it suffices to establish  \eqref{basicexp1-imp}  for
 $\ell=1$.  In that case,  
from \eqref{basicexp1} at $p_2=p_3$
and $\lambda=\frac{p_3}{3}$,  we deduce that for any $m \ge 1$, 
\[
\frac{(p_3/3)^m}{m!} \Ebm\big( |\bar V_1|^m\big) \le \Ebm\big( e^{\frac{p_3}{3} |\bar V_1|}\big) - 1 \le 
e^{\frac{4}{3} c_3 q}-1 \le \frac{8}{3} c_3 q \,.
\]
Hence,   in view of the right hand side of \eqref{4.2},  for $c_4:=96 c_3/p_3^{2}$ and $|\lambda| \le \frac{p_3}{4}$,  as claimed
\begin{equation}\label{basicexp1-imp-ad}
\bE \big[e^{\lambda \bar V_1} \big] \le 1 + \sum_{m \ge 2} \frac{|\lambda|^m}{m!} \Ebm\big( |\bar V_1|^m\big)
\le 1 + \frac{c_4}{4}
 q  \lambda^2
\sum_{m \ge 0} \big(\frac{3}{4}\big)^m 
\le e^{c_4 q \lambda^2} \,.
\end{equation}
Similarly,  combine \eqref{basicexp1} for $\bar V_1^\th$ with the \abbr{lhs} of \eqref{4.2}
to get \eqref{basicexp1-imp} for $\bar V_1^\th$.

Turning to  \eqref{eq-zav1-new},  by Construction \ref{con:renew},  we have the renewals for 
$\mu=\mu^{\tilde{r}}_{\npp,\tilde{R}}$ 
at random times $1=J_0 < J_1 < \cdots < J_{U_{\ell}} \le \ell$
that break our excursions into $U_{\ell} < \ell$ segments,  each of $L_j:=J_j-J_{j-1}$ excursions, 
followed by a segment of $L_0:=\ell +1 -J_{U_{\ell}} \ge 1$ excursions.  
Conditioned on $\{J_j\}$ these segments are independent and therefore so are
$\bar V_0:=V_\ell-V_{J_{U_{\ell}}-1}$ and $\bar V_j := V_{J_j-1}-V_{J_{j-1}-1}$ for $1 \le j \le U_{\ell}$,
where by construction the conditional law of $\bar V_j$ depends only on the value of $L_j$.  
Further,  from the right hand side of \eqref{4.2-cond},  we deduce that 
$\bE^\mu [\bar V_j | L_j ]=0$.  Thus,  by Jensen's inequality 
$\bE^{\mu} [ e^{\la \bar V_j} | L_j ]  \ge 1$,  and thereby
\begin{align}
 \bE^{\mu}\big[e^{\la V_\ell}\,\big
| \,  \{J_j\} \big] 
&=
\prod_{j=0}^{U_{\ell}} \bE^{\mu} \big[ e^{\la \bar V_j} | L_j \big] \le 
 \prod_{j=0}^{\ell-1}  \bE^{\mu} \big[ e^{\la \bar V_j} | L_j \big] \,.
\label{2renew.6ad}
\end{align}
Taking the expectation 
on both sides of \eqref{2renew.6ad},  we deduce by 
Cauchy-Schwarz and the independence of $\{L_j,  j \ge 1\}$,  that 
\begin{align}
\Big( \bE^{\mu}\big[e^{\la V_\ell}\,& \big]\Big)^2  
 \le  \bE \Big(  \bE^{\mu} \big[ e^{\la \bar V_0} | L_0 \big]^2 \Big) 
\bE \Big( \prod_{j=1}^{\ell-1}  \bE^{\mu} \big[ e^{\la \bar V_j} | L_j \big]^2 \Big) 
 \nn \\ & \le
\bE \Big(  \bE^{\mu} \big[ e^{2 \la \bar V_0} | L_0 \big] \Big) 
\bE \Big( \prod_{j=1}^{\ell-1}  \bE^{\mu} \big[ e^{2 \la \bar V_j} | L_j \big] \Big)  
 =  \prod_{j=0}^{\ell-1}  \bE^{\mu} \big[ e^{2 \la \bar V_j} \big] 
\label{cs-Ls-ad}
\end{align}
By construction,  the  \abbr{IID} renewal epoch $\{L_j,  j\ge 1\}$ are positive
Geometric($p_3$),  so as in  \eqref{basicexp1-imp-ad} we have that for any $j \ge 1$, 
\begin{equation}\label{eq:bd-barVj}
\bE^{\mu} \big[ e^{2 \la \bar V_j} \big] \le e^{4 c_4 q \la^2} \,,  \qquad \forall |\lambda| \le \frac{p_3}{8} \,.
\end{equation}
The right hand side of \eqref{4.2-cond} implies that $\E^\mu[\bar V_0|L_0] = 0$.  Further,  $L_0 \le \ell$,  
being the last passage epoch 
before $\ell$,  
is stochastically dominated by such positive Geometric($p_3$) variable.   
The derivation of \eqref{basicexp1-imp-ad} thus applies for $\bar V_0$,
extending the validity of \eqref{eq:bd-barVj} to all $j \ge 0$.  
Substituting this in \eqref{cs-Ls-ad},  yields \eqref{eq-zav1-new}.

Using the \abbr{lhs} of \eqref{4.2-cond} to deduce that $\bE^\mu [\bar V^\th_j | L_j ]=0$ for all $j \ge 0$, 
and adjusting $q$ to have \eqref{eq:bd-barVj} holding for $\bar V^\th_j$ instead of $\bar V_j$,  we get
by the same reasoning that \eqref{eq-zav1-new} applies for $V^\th_\ell$ instead of $V_\ell$.\\
(b).  In  \eqref{ergtotal.3} we consider  
the sum of $\{A_{\sfx,i}-B_{\sfy,i}\}$ over the first $\ell$ excursions from $\sfv_i \in \partial \bB_\sfd(\npp,\tilde{r})$
to $\sfu_i \in \partial \bB_\sfd(\npp,\tilde{R})$ and back,  conditioned to stay inside $\bB_\sfd(\npp,\wh R)$.
Just like $\wt{X}^i$ of \cite[Section 5]{BRZ2},  these
are excursions of a \emph{planar Brownian motion} between concentric balls,  which here 
are further conditioned to remain inside a larger,  concentric ball.  Note that 
the Markov chain 
$\{\sfu_i,  i \le \ell\}$,  whose (uniform) invariant measure $\mu_{\npp,\tilde{R}}$ 
induces starting the $\partial \bB_{\sfd}(\npp,\tilde{r})$-valued chain $\{\sfv_i\}$ 
at its invariant measure $\mu_{\npp,\tilde{r}}$,  is independent of 
$\{T:=\trav{1}{\npp}{\wh{R}}{\tilde{R}}{\npp}{\tilde{R}}{\tilde{r}}=\ell\}$ (ie,
first hitting $\partial \bB_\sfd(\npp,\wh R)$ just before the inner part 
of the $(\ell+1)$'st excursion).  Thanks to our assumptions on $\wh R$,  $\tilde{R}$ and $\tilde{r}$,
the Brownian hitting measure of  $\partial \bB_\sfd(\npp,\tilde{r})$ 
\emph{conditioned to stay inside $\bB_\sfd(\npp,\wh R)$},  which 
with some abuse of notation we now denote also as $\Pbm^\sfu$,  satisfies
\eqref{eq:unif-p3}.  This allows us to apply the procedure of Construction \ref{con:renew}
now with $p_3=p_3(a,b)$ and $\mu:=\mu_{\npp,\tilde{r}}$,  to get renewal times 
$1=J_0 < J_1 < \cdots < J_{U_{\ell}} \le \ell$ for the starting points $\{\sfv_i\}$ of the 
excursions $\wt X^i$ which are now further conditioned to stay inside $\bB_\sfd(\npp,\wh R)$.
As before,  the renewal times break our excursions into $U_{\ell} < \ell$ segments,  each  
of $L_j$ excursions that are conditioned to stay inside $\bB_\sfd(\npp,\wh R)$,  
followed by a segment of $L_0 \ge 1$ excursions,  at the end of which the Brownian 
path reaches $\partial \bB_\sfd(\npp,\wh R)$.  Conditioned on $\{J_j\}$ and $\{T=\ell\}$,  
these segments are independent,  resulting with $\bar V_j$,  $0 \le j \le U_{\ell}$ being
conditionally mutually independent,  with the conditional law of $\bar V_j$ depending only on the value of $L_j$.
Starting at $\mu$ which is also the invariant measure of the chain $\{\sfv_i\}$ we further deduce,
as in the derivation of \eqref{4.2-cond} that  for $1 \le i \le \ell$, 
\begin{equation*}
\Ebm^\mu [ A^\th_{\sfx,i}-B^\th_{\sfy,i} \,|\,  \{I_j\},  T=\ell ] = 0 \quad
\hbox{and} \quad
\Ebm^\mu [A_{\sfx,i}-B_{\sfy,i} \,|\, \{I_j\}, T =\ell ] =0\,.
\end{equation*}
Therefore,  as in part (a),  here too $\bE^\mu [\bar V^\th_j | L_j ]=0$ and $\bE^\mu[\bar V_j | L_j]=0$.
Proceeding hereafter as in the proof of \eqref{eq-zav1-new},   but with extra conditioning on $\{T=\ell\}$ 
and using now $\mu=\mu_{\npp,\tilde{r}}$,  we thus get similarly to \eqref{2renew.6ad} that 
\begin{align}
 \bE^{\mu}\big[e^{\la V_\ell}\,\big
| \, T = \ell ,  \{J_j\} \big] 
&=
\prod_{j=0}^{U_{\ell}} \bE^{\mu} \big[ e^{\la \bar V_j} | L_j \big] \le 
 \prod_{j=0}^{\ell-1}  \bE^{\mu} \big[ e^{\la \bar V_j} | L_j \big] \,.
\label{2renew.6}
\end{align}
Take the expectation with respect to $\{J_j,j \ge 1\}$ on both sides of \eqref{2renew.6},  to 
likewise arrive by Cauchy-Schwarz and the independence of $\{L_j,  j \ge 1\}$,  at
\begin{align*}
\Big( \bE^{\mu}\big[e^{\la V_\ell}\,& \big| \, T = \ell \big]\Big)^2  
\le
  \prod_{j=0}^{\ell-1}  \bE^{\mu} \big[ e^{2 \la \bar V_j} \big] \,.
\end{align*}
As before,   \eqref{eq:bd-barVj} holds for $0 \le j < \ell$,  thereby resulting with 
\eqref{ergtotal.3} and upon adjusting $q$ as usual,  the same argument applies also for $V^\th_\ell$. 
\end{proof}

As in \cite[(5.6)]{BRZ2},  for  $\sfx,\sfy \in \bB(\npp,\tilde{r}-R)$,  the traversal counts 
$\trav{n}{\npp}{\tilde R}{\tilde r}{\sfx}{R}{r}$ and $\trav{n}{\npp}{\tilde R}{\tilde r}{\sfy}{R}{r}$
 of \Cref{def-trav} are such that 
\begin{equation}
\trav{n}{\npp}{\tilde R}{\tilde r}{\sfx}{R}{r}-\trav{n}{\npp}{\tilde R}{\tilde r}{\sfy}{R}{r}=
\travshort{n}{\npp}{\tilde{R}}{\tilde{r}}{\sfx}{r}{R} -  \travshort{n}{\npp}{\tilde{R}}{\tilde{r}}{\sfy}{r}{R}
= \sum_{i=1}^{n} A_{\sfx,i} -  \sum_{i=1}^{n}B_{\sfy,i} =  V_n \,.
 \label{1.5tot}
\end{equation}
Thus,  utilizing  \Cref{cor:HigherMomentBound} and \Cref{lem:c3new},
we establish the next lemma,  whose two parts are used for proving Lemmas 
\ref{rlem:Continuity lemma} and \ref{lem:4.58m},  respectively.
\bl\label{lem-cont1}
For any $C_{0}<\ff$ there exist $c_1<\infty$,  $c_{0},   q_{0} >0$,  depending only on $a,b,C_0$,
so for any $n \ge 1$ and 
$\sfx,\sfy \in \bB_{\sfd}(\npp,\wt{r}-R)$ with $q:=\sfd(\sfx,\sfy)/R\leq q_{0}$,
\begin{align}
\Pbm\left(\left| \trav{n}{\npp}{\tilde R}{\tilde r}{\sfx}{R}{r}
        - \trav{n}{\npp}{\tilde R}{\tilde r}{\sfy}{R}{r}
        \right|\ge\theta\sqrt{  qn}
\right) & \le  e^{-C_{0}\theta^{2}\sqrt{q}}\quad  \forall \theta \in [c_1, c_{0}\sqrt{n}\,] \,,
\label{cb.1} \\
\Pbm\left(\left| \trav{n}{\npp}{\tilde R}{\tilde r}{\sfx}{R}{r}
        - \trav{n}{\npp}{\tilde R}{\tilde r}{\sfy}{R}{r}
        \right|\ge\theta\sqrt{  qn}
\right)
& \le e^{-C_{0}\theta} \,, \qquad 
\forall \theta \in [c_{1}, (nq)^{2} \wedge c_0 n] \,. 
\label{cb2.0}
\end{align}
\el
\begin{proof} 
In view of \eqref{1.5tot},  this is about tail estimates for $V_n$ of \eqref{eq:def-V} and the proof of
\eqref{cb.1} parallels that of \cite[(5.7)]{BRZ2}, 
except that here we do not condition on their events $Y_m \in \mathcal{I}_m$.
Indeed,  we merely follow the bound of  \cite[(5.35)-(5.36)]{BRZ2},  utilizing in \cite[(5.37)]{BRZ2}
our bound \eqref{basicexp1-imp} for their ``main" term $A_1$,  while using 
\eqref{basicexp1} to control the ``reminder'' term $A_2$ there.
Armed with the above,  the reader can check that the proof of \cite[(5.7)]{BRZ2} carries over 
here to get \eqref{cb.1}.  

Turning to larger values of $\theta$ (that exhibit
the transition from Gaussian moderate deviations to exponential large deviations),  the 
proof of \eqref{cb2.0} proceeds similarly,  by adapting \cite[proof of Lemma 5.6]{BRZ2}.
\end{proof}
Our next result is a conditional version of \eqref{cb2.0},  used for proving \eqref{eq-thanks4lem}
(see \eqref{eq-calgary} for $T_{k'}^{\npp,\bar k,m}$).
\bl\label{lem-cont3}  
For any $C_{0}<\ff$ there exist  $c_1<\infty$,  $q_0>0$,  depending only on $a,b,C_0$,
so if  \eqref{eq-ratiobounds} holds for $\tilde{r}= r_{k'}$ and $\tilde{R}=r_{k'-1}$,  then for 
any $n \ge 1$ and $\sfx,\sfy \in \bB_{\sfd}(\npp,\wt{r}-R)$ with $q:=\sfd(\sfx,\sfy)/R\in [\frac{c_1^2}{n}, q_{0}]$,
\begin{equation}
\Pbm\Big(
\big| \trav{n}{\npp}{ r_{k'-1}}{r_{k'}}{\sfx}{R}{r}
        - \trav{n}{\npp}{r_{k'-1}}{r_{k'}}{\sfy}{R}{r}
        \big|
\ge\theta\sqrt{  qn}\,\Big| 
\, T_{k'}^{\npp,\bar k,m} = n, \, {\mathcal G}_{\bar k}^\npp \Big)
 \le e^{-C_{0}\theta}\,,
\label{cb2.1}
\end{equation}
for all $\theta \ge c_{1}$,  $m \ge 1$ and $k' - \bar k \ge 2 \vee \log (n \log n)$.
\el
\begin{proof} We will show for some 
$p_4,q_0>0$,  $c_5, C<\infty$ depending only on $a$,  $b$,
\bea
 \bE\big(e^{\la V_n} 
 \,\big| \, 
T_{k'}^{\npp,\bar k,m} = n, \, {\mathcal G}_{\bar k}^\npp \big)\leq C e^{ c_5 q  \lambda^{2} n}, 
\quad \forall n \ge 1,  \la \in [0,p_4]\,,  q \le q_0 \,.
\label{ergtotal.1}
\eea
By Chebyshev's inequality,  such a bound implies that
\begin{align}\label{2renew3a}
\mathbb{P}\big( V_n \ge\theta\sqrt{q n}\,\Big
| \, T_{k'}^{\npp,\bar k,m}
 = n, \, {\mathcal G}_{\bar k}^\npp \big) 
& \leq  C e^{ c_5 q  \lambda^{2} n} e^{ -\la \theta\sqrt{qn}}.
\end{align}
By symmetry,  \eqref{2renew3a} holds also for $-V_n$.   Recalling \eqref{1.5tot} and taking 
$\lambda = \min\{ \theta/(2c_5\sqrt{qn}), p_4\}$,  we bound the 
logarithm of the \abbr{lhs} of \eqref{cb2.1} by
\[
\log (2C) - \frac{\theta}{4c_5} \min\{\theta, 2 p_4 c_5 \sqrt{qn}\}  \le - C_0 \theta \,,
\]
provided $\theta \wedge \sqrt{qn} \geq c_1$, 
with $c_1$ large enough (depending on $p_4$,  $c_5$,  $C$ and $C_0$). 

Turning to prove \eqref{ergtotal.1},  let $Z_{k'}^{i'}$ count the excursions from $\partial \bB_{\sfd}(\npp,r_{k'-1})$ to $\partial \bB_{\sfd}(\npp,r_{k'})$ during the $i'$'th excursion from $\partial \bB_{\sfd}(\npp,r_{\bar k})$ to $\partial \bB_{\sfd}(\npp,r_{\bar k-1})$ and set 
\[
C := \sup_{\ell \ge 2 \vee \log (n\log n)} \{ (1+ 28.2 (\ell+1) e^{-\ell})^n \} < \infty \,.
\]
By the law of total probability,  \eqref{ergtotal.1} follows upon showing for 
$s_m := \displaystyle{\sum_{j=1}^m} \ell_{j}$,  $S :=1+ 28.2 (k'-\bar k+1) e^{-(k'-\bar k)}$
and any fixed $\{ \ell_{j} \in \bZ_+,  j \ge 1\}$, 
\be
\mathbb{E}\big( e^{\la V_{s_m}} \,\big
| \, \{Z^{i'}_{k'}=\ell_{i'}, i'\le m\} , \, {\mathcal G}_{\bar k}^{\npp} \big)\leq (S e^{  c_5 q \lambda^{2}})^{s_m} 
\,,\quad \forall \la\in (0,p_4]\,,
\label{ergtotal.1z}
\ee
Setting $V_0:=0$  and
$T := \trav{1}{\npp}{r_{\bar k-1}}{r_{k'-1}}{\npp}{r_{k'-1}}{r_{k'}}$ counting
the $k'$-level excursions around $\npp$ till reaching $\partial \bB_\sfd(\npp,r_{\bar k-1})$,  note that 
by 
\eqref{m1.5eucc},
\begin{align}
 \mathbb{E}\big(e^{\la V_{s_m}} \,\big|  \, \{Z^{i'}_{k'}=\ell_{i'}\} ,  \, {\mathcal G}_{\bar k}^{\npp} \big)
& =\mathbb{E}\big(\prod_{j=1}^{m} e^{\la (V_{s_j}-V_{s_{j-1}})}
\,\big| \, \{Z^{i'}_{k'}=\ell_{i'}\} , \, {\mathcal G}_{\bar k}^{\npp} \big) 
\nn \\ & 
\leq \prod_{j=1}^{ m} \, S^{\ell_{j}} \,
\bE^{\mu_{\npp,r_{k'-1}}} \big(e^{\la V_{\ell_{j}}}\,\big| T = \ell_{j}\big).
 \label{2renew.4}
\end{align}
Thus,  applying \eqref{ergtotal.3} for
$\tilde{r}=r_{k'}$,  $\tilde{R}=r_{k'-1}$,  $\bar R = r_{\bar k-1}$ and $\ell=\ell_{j}$ for $1 \le j \le m$,  
yields \eqref{ergtotal.1z},  thereby completing the proof of the lemma.  
\end{proof}

\subsection{Proof of \Cref{lem:4.58m}}\label{sec:4.58m}


We set $\bar k := l - 4 \log (\alpha(l)+j+1) \ge 1$.  Decomposing the event in \eqref{eq-thanks4lem} 
according to the value of $T_{\bar k}^{{\bf 0},t} \in \GG^{\bf 0}_{\bar k}$,  it suffices to show that 
\eqref{eq-thanks4lem} holds under the (bigger) conditioning of \eqref{cb2.1},  
where $n=\frac{u^2}{2}$.  To this end,  we utilize \Cref{lem-cont3} for $C_0=8/d_0$,
$k'= l-2$,  $R= \wt{r}_{l-1}^+$ and $r=\wt r_l$,  noting that \eqref{eq-ratiobounds} 
then holds (with $a=1/(3e)$,  $b=3/(2e)$, as $\beta \le 1/2$).  Since  
$\wt{c} \le \frac{1}{2} \le e^2 - \frac{3}{2}e$,  clearly $
\bB_{\sfd}\left({\bf 0},\wt{c} r_{l}\right) \subset 
\bB_{\sfd}({\bf 0},r_{k'}-R)$ and taking $\wt c \le q_0$,  also 
\be\label{eq:q-q0}
q:=\frac{\sfd(\wt \sfy,\sfy)}{R}
= \frac{\tilde q}{(1+\beta) e} \le \frac{2 \wt{c}}{e}  \le q_0 \,.
\ee
For $j_o$ large enough 
$n q \ge u^2 \tilde{q}/(3e) \ge u \wt{c}/(12e) \ge j_o \wt{c}/(12 e) \ge c_1^2$
and since $j \ge j_o$ yields $n \ge n_o(j_o)$,  then also 
\[
k' - \bar k = 4 \log (\al(l)+j+1) - 2 \ge \log (4 n^2/e^2) \ge 2 \vee \log (n\log n) \,.
\]
Having verified all the conditions of \Cref{lem-cont3},  we proceed to apply \eqref{cb2.1}
at $\theta := d_0 j m \sqrt{(1+\beta) e/2}$.
Indeed,  such $\th \ge d_0 j_o \ge c_1$ for large $j_o$,
and  
$C_0 \theta \ge 8 
j m \ge 4 (j+m)$,  as claimed.

We get the unconditional version of \eqref{eq-thanks4lem} by applying instead 
\eqref{cb2.0} for $\tilde{r}=r_{l-2}$,  $\tilde{R}=r_{l-3}$
and $R,r$ as before (in particular,  \eqref{eq-ratiobounds} holds and $\wt \sfy$,  $\sfy$ 
satisfy the conditions of \Cref{lem-cont1},  with $q$ given again by \eqref{eq:q-q0}).
Setting again $n=\frac{u^2}{2}$,  our previous choices of $C_0$ and $\theta$ yield 
$C_0 \theta \ge 4(j+m)$ and $\theta \sqrt{qn} = \frac{d_0}{2} j m u \sqrt{\wt q}$.
By definition $j \beta \le 1$,  so our assumption that $m \beta \le 1/j_o$ implies  
that $j m \beta^2 \le 1/j_o$.  Hence,  $\theta/n=2 \theta/u^2$ is arbitrarily small for $j_o$ large enough
(which as we saw before,  also guarantees that $\theta \ge c_1$).  Lastly,  since 
\[
n q^2 \ge 2 \frac{(u \tilde{q})^2}{(3 e)^2} \ge \frac{\wt{c}^2}{8 (3 e)^2} := c_2 \,,
\]
also $(n q)^2 \ge c_2 n \ge \theta$ when $j_o$ is large enough,  thereby allowing us to use
\Cref{lem-cont1} and complete the proof of \Cref{lem:4.58m}.

\appendix

\section{Barrier estimates}\label{sec:BoundaryCrossing}

Recall our notations $T_{\ell}^{\sfx,k,m} =\trav{m}{\sfx}{r_{k-1}}{r_{k}}{\sfx}{r_{\ell}}{r_{\ell-1}}$ 
of \eqref{eq-calgary},  $\al(\cdot):=\gamma_{-1}(\cdot)$ of \eqref{eq:norm-exc},
$t_z:=\sfs_1(z)$ of \eqref{eq:nominal},   and $I_{y }=[y,y+1)$ (denoted $H_y$ in \cite{BRZ}).

Since our barrier estimates are translation invariant,  \abbr{wlog} we fix $\sfx=\npp$.  We then set
\[
\eta^{\langle z \rangle,k}_\ell:=\sqrt{2T_{\ell}^{\npp,k, \sfs_k(z)}} \quad {\rm with}  \quad  
\wh{\eta}^{\langle z \rangle,k}_\ell := \eta^{\langle z \rangle,k}_\ell - \gamma_0(\ell),  \quad 
\etap_\ell:=\eta^{\langle z \rangle,1}_\ell  \quad
{\rm and} \quad 
\whetap_\ell := \wh{\eta}^{\langle z \rangle,1}_\ell, 
\]
building in this appendix on \cite[Theorem 1.1]{BRZ} and \Cref{cor:bend1},  to bound the barrier probabilities
\begin{align}\label{dfn:KQ}
 \KQ^{(s_-,k_+)}_{[s,k]}(u,v) 
\,&:=\Pbm\left(\etapsm_\ell \ge \al(\ell) , \, s \le \ell \le k;\,  \whetapsm_{k_+} \in I_{v} \right)\,,  \\
\KQ^{(s_-,j,k)}_{[s,k]}(u,w,v) &:=
\Pbm\left(\etapsm_\ell \ge \al(\ell) , \, s \le \ell \le k;\,  
\whetapsm_{j} \in I_{w}\,, 
\whetapsm_{k} \in I_{v} \right)\,,
\label{dfn:KQj}
\end{align}
for various values of $0<s_- \le s < j \le k \le k_+ < n_\vep$.  We first deal in \Cref{prop:BarrierSecGWPropUB}
with $s_-=s$ and $k_+=k$ (generalizing \cite[Lemma 8.1]{BRZ2},  where $j=k$ and $s=0$).  We then bound
the probability in \eqref{dfn:KQ} for $1=s_-<s=o(n_\vep)<\frac{1}{2} (k-s)$,  first in \Cref{lem-bendc} 
where $k_+=k$,  followed by \Cref{lem-bendd} where also $k_+-k =o(n_\vep)>0$.
\bl\label{prop:BarrierSecGWPropUB}
For some $C$ finite,  all small $\vep>0$ and any $0<s<j<k < n_{\vep}$,
\[
u \in [-s_{n_\vep}^\gak,100 (n_\vep-s)],  \quad w \in [-j_{n_\vep}^\gak,100(n_\vep-j)],  \quad 
v  \geq -k_{n_\vep}^\gak,
\] 
we have that 
\begin{align}
\KQ^{(s,k)}_{[s,k]}(u,v)  & \leq C \,  \bar{\KQ}_{s,k,u,v}  \,  e^{-\frac{\left( u-v \right)^{2}}{4(k-s)}} 
  \,,  \label{18.20bendj1} \\
\KQ^{(s,j,k)}_{[s,k]}(u,w,v)  & \le C \, \bar{\KQ}_{s,k,u,v} \,
\frac{(1+w+j_{n_\vep}^\gak)^2} {((k-j) \wedge (j -s))^{3/2}} \, e^{-\frac{(u-w)^2}{4 (j-s)}} e^{-\frac{(w-v)^2}{4 (k-j)}}
 \,, 
 \label{eq:bendA}
\end{align} 
where 
\begin{align}\label{dfn:bar-Q}
\bar{\KQ}_{s,k,u,v} & := \Delta_{k-s} \,
 \frac{\left(1+u+ s_{n_\vep}^{\gak} \right)\left(1+ v +   k_{n_{\vep}}^{\gak}\right)}{\sqrt{(n_\vep-k) \wedge (k-s)}} \, 
 e^{-2(k-s)-2(u-v)}\,\\
\label{dfn:Delta}
\Delta_{k-s} &:= \frac{1}{k-s} \exp\Big( \frac{(k-s)}{n_\vep} \log n_\vep\Big) \le 2 \,.
\end{align}
\el
\begin{proof} Starting with \eqref{18.20bendj1},  with $\gak \le \frac{2}{3}$ we utilize \Cref{cor:bend1} to
enlarge the event considered in \eqref{dfn:KQ} by replacing the requirements
of the form $\etaps_\ell \ge \al(\ell)$ with $\etaps_\ell \ge \gamma_{-1}^{[s,k]}(\ell)$ of \eqref{dfn:gamma-J}.
This allows us to apply \cite[(1.1)]{BRZ} with the parameters  
\begin{equation}\label{eq:xy-BRZ}
L=k-s,  \quad a=\al(s),  \quad b=\al(k),  \quad 
x=\sfrho_{\vep} (n_\vep-s) +u \qquad {\rm  and} \qquad  y= \sfrho_\vep (n_\vep-k)+v\,,
\end{equation}
yielding the bound
\begin{equation}\label{eq:bd11-BRZ}
\KQ^{(s,k)}_{[s,k]}(u,v) \le c \,  \frac{(1+x-a)(1+y-b)}{L} \sqrt{\frac{x}{y L}} e^{-\frac{(x-y)^2}{2L}} 
\end{equation}
(see \cite[above (5.10)]{DRZ} on removing the requirement that $\max(x, y) \le c' L$).
We note in passing that our shift of $\sfae \in [0,1)$ at level $\sfKtp$,  see \eqref{dr.1},  
is easily accommodated by increasing the constant $c$ in \eqref{eq:bd11-BRZ}. 
Turning to simplify \eqref{eq:bd11-BRZ} towards \eqref{18.20bendj1},  
first recall from  \eqref{def:mn} and \eqref{def.nvep},  that 
\[
2-\sfrho_\vep = \frac{\log (n_\vep)}{2n_\vep} \Big(1 + O(n_\vep^{-1}) \Big).
\]
Hence, 
\begin{align}\label{eq:lin-expan}
\frac{(u-v- (2-\sfrho_\vep) (k-s))^2}{2(k-s)} &= \frac{(u-v)^2}{2(k-s)} - (2-\sfrho_\vep) (u-v) + o_\vep(1)
\ge \frac{(u-v)^2}{4(k-s)} - o_\vep(1) 
\end{align}
(using $(A-B)^{2}\geq A^{2}/2-B^{2}$),  and we consequently have that
\begin{align}
\frac{(x-y)^2}{2L} 
= \frac{(\sfrho_\vep (k-s) + u-v)^{2}}{2(k-s)} 
&= 2 (k-s) + 2(u-v) - 2 (2-\sfrho_\vep) (k-s) + \frac{(u-v+ (\sfrho_\vep-2) (k-s))^2}{2(k-s)} \nonumber \\
&\ge  2 (k-s) + 2(u-v) - \frac{(k-s)}{n_\vep} \log n_\vep + \frac{(u-v)^2}{4(k-s)} - o_\vep(1) \,.
\label{eq:ld-xyL}
\end{align}
Moreover,  with $L=k-s$,  in view of \eqref{eq:xy-BRZ},  our conditions on $u$ and $v$ guarantee that 
 \begin{equation}\label{eq:xyL-bd}
 x \le 
 102 (n_\vep-s) = 
 102 (n_\vep-k+L) \le 204 \max\{n_\vep-k,L\} \le \frac{204\,  y L}{(n_\vep-k)\wedge (k-s)}  \,.
 \end{equation}
Plugging \eqref{eq:ld-xyL} and \eqref{eq:xyL-bd} in \eqref{eq:bd11-BRZ} yields \eqref{18.20bendj1}, 
with $\Delta_{j} \le 2$ by the concavity of the logarithm.

Turning to \eqref{eq:bendA},  note that by the Markov property of $\ell \mapsto \etaps_{\ell}$,
\[
\KQ^{(s,j,k)}_{[s,k]}(u,w,v) \le \KQ^{(s,j)}_{[s,j]}(u,w) \sup_{u' \in I_w} \{ \KQ^{(j,k)}_{[j,k]}(u',v) \} \,.
\]
Hence,  following our derivation of \eqref{18.20bendj1} and embedding the worst case effect 
of  $u'-w \in [0,1]$ within our universal constant,  we get by twice applying  \cite[(1.1)]{BRZ},  that 
\begin{align}\label{eq:bd111-BRZ}
\KQ^{(s,j,k)}_{[s,k]}(u,w,v) &\le c \,  \frac{(1+z-b')^2}{(L_1 \wedge L_2)^{3/2}}\,  \frac{(1+x-a)(1+y-b)}{L} 
\sqrt{\frac{x}{y L}} e^{-\frac{(x-z)^2}{2L_1}} e^{-\frac{(z-y)^2}{2L_2}}  \; 
\end{align}
for $x,y,a,b,L$ of \eqref{eq:bd11-BRZ} and 
\[
L_1=j-s,  \quad L_2=k-j,  \quad z = \sfrho_\vep (n_\vep-j) + w  = b' + w + j_{n_\vep}^\gak  \,.
\]
Expanding the quadratic terms in the exponent of \eqref{eq:bd111-BRZ}
as in \eqref{eq:lin-expan} and \eqref{eq:ld-xyL},  yields
\begin{equation}\label{eq:quad-exp2}
\frac{(x-z)^2}{2L_1}+\frac{(z-y)^2}{2L_2}  \ge 2 (k-s) + 2(u-v) -  \frac{(k-s)}{n_\vep} \log n_\vep 
+ \frac{(u-w)^2}{4 L_1} + \frac{(w-v)^2}{4L_2} - o_\vep(1) \,.
\end{equation}
Now,  plugging \eqref{eq:quad-exp2} in \eqref{eq:bd111-BRZ} and having bounded already 
by $\bar\KQ_{s,k,v,u}$ those terms which are common with the right hand side 
of \eqref{eq:bd11-BRZ},  results with \eqref{eq:bendA}.
\end{proof}

Hereafter,  we ignore the effects 
of our $\sfae$ fractional shift and if needed also replace $\sfs_0(z)$ by $\sfs_1(z)$, 
since the same rational as in the proceeding proof applies to all our barrier estimates.

As promised,  we next generalize \Cref{prop:BarrierSecGWPropUB},  imposing a barrier from level $s>1$ 
on the counts of excursions completed during a specified number of level one excursions.
\bl\label{lem-bendc}
For some $C$ finite,  all small $\vep>0$,  any $s \log n_{\vep}/n_{\vep}\leq c$,  $3 s<  k < n_{\vep}$,   
$|z| \leq \log n_{\vep}$ and $v  \in [ -k_{n_\vep}^\gak,20\sqrt{n_\vep (k-s)}]$, 
\begin{align}
\KQ^{(1,k)}_{[s,k]} (z,v) 
\label{18.20srgbend}
\leq & C \,  \Delta_{k-s} \,  e^{-2k-2z+2v}\,
 \frac{\left(1+\sqrt{s}+ z_+ \right) \left(1+ v +   k_{n_{\vep}}^{\gak}\right)}{\sqrt{(n_\vep-k) \wedge (k-s)}} \, 
 e^{-\frac{\left( z-v \right)^{2}}{8(k-s)}}  \,.
\end{align}
\el

\begin{proof} For small $\vep$,  by our assumptions $s \le n_\vep/2$ and by the Markov property,  
\begin{equation}
\KQ^{(1,k)}_{[s,k]} (z,v) \leq \Pbm ( \whetap_s \ge 20 n_\vep )
\label{18.34b}
+\sum_{u=-s^{\gak}}^{20 n_{\vep}} \Pbm \big( \whetap_s \in I_{u} \big) \sup_{y \in I_u} \{ \KQ^{(s,k)}_{[s,k]}(y, v) \} \,.
\end{equation}
Further,  by \cite[Prop. 1.4]{BRZ}
and using that $s\log n_{\vep}/n_{\vep}\leq c$,  in a computation similar to that of \eqref{eq:ld-xyL},
\bea
\Pbm \left( \whetap_s \in I_{u} \right)
\leq  C' \sqrt{\frac{\sfrho_\vep n_{\vep}+z}{(\sfrho_\vep (n_\vep-s)+ u) s}}\,\,e^{-(\sfrho_\vep s+z-u)^{2}/2s} 
\le  
 \frac{c_1}{\sqrt{s}}\,\,e^{-2s-2 (z-u)} 
e^{-\frac{(z-u)^{2}}{4s}}\,.   
\label{180.20sr2}
\eea
In particular,  using \eqref{180.20sr2} we see that 
\begin{equation}\label{eq:neglig-A}
\Pbm ( \whetap_s \ge 20 n_\vep ) \le 
c_1 e^{2s} 
\sum_{u = 20 n_\vep}^{\infty} \frac{1}{\sqrt{s}} e^{-\frac{(u-z-4s)^{2}}{4s}} \le c_2 e^{-80 n_\vep} \,,
\end{equation}
which is negligible in comparison to the right hand side of \eqref{18.20srgbend}.  As before,  the effect of having 
$y-u \in [0,1]$ in \eqref{18.34b} is negligible,  so plugging 
\eqref{18.20bendj1} and \eqref{180.20sr2} on the right hand side of \eqref{18.34b}
and setting $u=w-s^\gak$,  yields the stated bound of the lemma,  once we show that 
\begin{equation}\label{eq:jay-quad}
\exp\Big(\frac{( z-v)^{2}}{8(k-s)}\Big) \,
\sum_{w=0}^{\infty} 
\frac{(1+w)}{\sqrt{s}}\, e^{-\frac{\left(w-s^\gak-v\right)^{2}}{4(k-s)}} e^{-\frac{(z+s^\gak-w)^{2}}{4s}} \, 
\le C (1+\sqrt{s}+ z_+)  \, 
\end{equation}
Now,  using $(A-B)^{2}\geq A^{2}/2-B^{2}$ for $A=z-v$ and $B=z+s^\gak-w$,  as well as $2s< k-s$,   we have that 
\begin{equation}\label{eq:A-B-bd}
- \frac{\left(w  - s^{ \gak} - v \right)^{2}}{4(k-s)} \leq - \frac{\left(z-v \right)^{2}}{8(k-s)} + \frac{(z+ s^{ \gak}-w)^{2}}{8s}.
\end{equation}
Applying the latter inequality,  we bound the \abbr{lhs} of \eqref{eq:jay-quad} by 
\[
\sum_{w=0}^{\infty} 
\frac{(1+w)}{\sqrt{s}}\,  \exp\Big(-\frac{(w-z-s^\gak)^{2}}{8s} \Big) \le  C (1+\sqrt{s}+ z_+) 
\]
(recall that $\gak \le \frac{1}{2}$),  whereupon the proof of \eqref{18.20srgbend} is completed.
\end{proof}

Lastly,  we generalize the result of \Cref{lem-bendc},  
by specifying the excursion count only at level $k_+$,  while the barrier stopped already at level $k<k_+$.
\bl\label{lem-bendd}
For some $C$ finite,  small $\vep>0$ any 
\[
(k_+-k) \log n_{\vep}/n_{\vep}\leq c,  \quad k_+-k \le k - s,  \quad v \in [-(k_+)_{n_\vep}^{\gak}, 20\sqrt{n_\vep(k-s)}],
\]
and all $s$,  $z$,  $k < k_+ < n_\vep$ as in \Cref{lem-bendc},
\begin{align}
\KQ^{(1,k_+)}_{[s,k]}(z,v) 
\label{18.20srbend}
&\leq  C \,  \Delta_{k-s} \,  e^{-2k_+-2z+2v}\,
 \frac{\big(1+\sqrt{s}+ z_+ \big) \big(1+  \sqrt{k_+-k} +
 v +   (k_+)_{n_{\vep}}^{\gak} \big)}{\sqrt{(n_\vep-k_+) \wedge (k-s)}} \, 
 e^{-\frac{\left( z-v \right)^{2}}{16(k-s)}}  \,.
\end{align}
\el
\begin{proof} Setting $b=15 \sqrt{n_\vep k} \le 20 \sqrt{n_\vep (k-s)}$,  
we have as in \eqref{18.34b},  that 
\begin{equation}\label{18.34A}
\KQ^{(1,k_+)}_{[s,k]}(z,v)  \le \Pbm( \whetap_k \ge b) + 
 \sum_{u=-k_{n_\vep}^\gak}^{b} \KQ^{(1,k)}_{[s,k]}(z,u) \sup_{w \in I_u} \{
\Pbm (\wh{\eta}^{\langle w \rangle,k}_{k_+} \in I_v) \} \,.
\end{equation}
Similarly to \eqref{eq:neglig-A},  the first term on the right hand side is at most $e^{-40 n_\vep}$,  hence negligible.
Setting $r=k_+-k$,  we find by applying \cite[Prop. 1.4]{BRZ}
and having $r \log n_{\vep}/n_{\vep}\leq c$,  that similarly to \eqref{180.20sr2},
\begin{align}
\sup_{w \in I_u} \{ \Pbm (\wh{\eta}^{\langle w \rangle,k}_{k_+} \in I_v) \} &\le
\sqrt{\frac{\gamma_0(k)+u}{\gamma_0(k_+)+ v}}
 \frac{c_1}{\sqrt{r}}\,  e^{-2r-2 u+2v} e^{-\frac{(u-v)^{2}}{4r}} \,.
\label{180.20t2}
\end{align}
The first factor on the right hand side of \eqref{180.20t2} may be large,  but we absorb it
by adjusting the term $1/\sqrt{y}$ on the right hand side of \eqref{eq:bd11-BRZ},  en-route to bounding 
$\KQ^{(1,k)}_{[s,k]}(z,u)$ analogously to  \eqref{18.20srgbend}.
Doing so,  plugging in \eqref{18.34A} the other part of the right hand side of \eqref{180.20t2}
and writing $u=w-k_{n_\vep}^\gak$,   bounds the sum 
in \eqref{18.34A}  by  
\[
C \,  \Delta_{k-s} \,  e^{-2k_+-2z+2v}\,
 \frac{\left(1+\sqrt{s}+ z_+ \right)}{\sqrt{(n_\vep-k_+) \wedge (k-s)}} \,  \sum_{w=0}^\infty 
 \frac{(1+ w)}{\sqrt{r}}
 e^{-\frac{\left( z+k_{n_\vep}^\gak-w \right)^{2}}{8(k-s)}}  \,  e^{-\frac{(w-k_{n_\vep}^\gak -v)^{2}}{4r}} \,.
\]
By our assumption that $r<k-s$,  we get similarly to \eqref{eq:A-B-bd} that  
\[
- \frac{(z+ k_{n_\vep}^{ \gak}-w)^{2}}{8(k-s)} \leq - \frac{\left(z-v \right)^{2}}{16(k-s)} 
+ \frac{\left(w  - k_{n_\vep}^{ \gak} - v \right)^{2}}{8r} \,,
\]
so after verifying that 
\[
 \sum_{w=0}^\infty 
 \frac{(1+ w)}{\sqrt{r}} \exp\Big(-\frac{(w-k_{n_\vep}^\gak -v)^{2}}{8r}\Big)
\le C (1 + \sqrt{r} + v+ (k_+)_{n_\vep}^\gak) \,.
\]
we establish the bound  \eqref{18.20srbend} and complete the proof.
\end{proof}

\begin{remark}\label{rmk:straight-line}
While we stated the bounds of \Cref{lem-bendc} and \Cref{lem-bendd} only for the probabilities in 
\eqref{dfn:KQ},  under exactly the same modifications we get the corresponding bounds on
$\KQ_{[s,k]}^{(1,j,k_+)}(u,w,v)$.  Also,  by the preceding derivations,  if in \eqref{dfn:KQ} or \eqref{dfn:KQj} 
we replace the curved barrier $\al(\cdot)$ by the straight one $\gamma_0(\cdot)$,  
then the terms $s_{n_\vep}^\gak$,  $j_{n_\vep}^\gak$,  $k_{n_\vep}^\gak$ and $(k_+)_{n_\vep}^\gak$
can be omitted in the right hand side of our bounds.
 \end{remark}

\section{Conditional excursion probabilities}\label{probends}
In \cite[Appendix II]{BRZ2}, the  claim that proving \cite[(9.24)]{BRZ2} is sufficient for 
\cite[Lemma 9.1]{BRZ2} is false, leading to an incorrect statement in \cite[(9.1)]{BRZ2}. The corrected statement of \cite[(9.1)]{BRZ2} is Lemma \ref{recursionends2} below. The statement 
\cite[(9.2)]{BRZ2} is correct, and appears here as Lemma \ref{recursionends}. Both statements rely on \cite[Lemma 9.2]{BRZ2}, whose proof there is correct. For completeness, we repeat \cite[Lemma 9.2]{BRZ2} here as Lemma \ref{lemprobends}, and provide the proof, using the notation of 
\cite{BRZ2} and working on the sphere $\bS^2$, as follows. For  $ k \geq 1$ and $y\in \bS^2$,  let ${\mathcal G}_k^y$ denote the $\sigma$-algebra generated by the excursions
from $\partial B_{d}(y,h_{k-1})$ to $\partial B_{d}(y,h_{k})$ (and if we start outside $\partial B_{d}(y,h_{k-1})$ we include the initial excursion to $\partial B_{d}(y,h_{k-1})$). 
Note that ${\mathcal G}_k^y$ includes the information on the end points of the excursions from  $\partial B_{d}(y,h_{k})$ to  $\partial B_{d}(y,h_{k-1})$.

Note that  \cite[Appendix II]{BRZ2} and this appendix are  stated for a Brownian motion on $\bS^2$, however a suitable
one-to-one mapping of $\bB_\sfd (\sfv,r_k)$, $k \ge n'-1$ to concentric 
caps on $\bS^2$ sends our Brownian path, up to time rescale,  to the  path of a Brownian motion
on $\bS^2$.
Everywhere  $h_j$ should be read as $r_j$ and $B_d$ should be read as $\bB_\sfd$ (and everywhere $j\leq n_\vep$).
We let  $\trav{n}{y}{k-1}{k}{y}{l-1}{l}$ denote the number of excursions    from
	$\partial B_{d}(y,h_{l-1})\to\partial B_{d}(y,h_{l})$ during the first $n$
 excursions    from
	$\partial B_{d}(y,h_{k })\to\partial B_{d}(y,h_{k-1})$, which is consistent with  \eqref{eq-calgary}.

\begin{lemma}\label{recursionends}
For any $k '>k+10\geq 11$ and $n \ge 1$,
%
  all $m_l ;l=k',\ldots L$, and all $y \in \bS^2$, 
\bea
&& \Pbm(\trav{n}{y}{k-1}{k}{y}{l-1}{l} = m_l ;l=k'+1,\ldots L \,
|\, \trav{n}{y}{k-1}{k}{y}{k'-1}{k'} = m_{k'}, \, {\mathcal G}_k^y) \nonumber \\ 
&&=
\left(1+ O\( (k'-k+1)\frac{h_{k'-1}}{h_{k-1}}\)\right)^{m_{k'}} 
\label{m1.2e}
\Pbm(\trav{n}{y}{k-1}{k}{y}{l-1}{l} = m_l ;l=k'+1,\ldots L \,
|\, \trav{n}{y}{k-1}{k}{y}{k'-1}{k'} = m_{k'} ).
\eea
\end{lemma}
\noindent
In Lemma \ref{recursionends}
we have
\begin{equation}
| O\( (k'-k+1)\frac{h_{k'-1}}{h_{k-1}}\)|\leq 50\( (k'-k+1)\frac{h_{k'-1}}{h_{k-1}}\).\label{}
\end{equation}

The key to the proof of Lemma \ref{recursionends} is the following Lemma. When we write $|\frac{A}{B}-1|\leq C$ we mean $|A-B|\leq CB$, even when $B=0.$
Also, we make the usual abuse of notation when writing conditioning on $X_{\bar \tau}=v$, using that all random variables for which 
conditionings of this form are made possess a continuous density.

\begin{lemma}\label{lemprobends}
For $ k \geq 0$,   $k' \geq k+10$
and a Brownian path $X_\cdot$ 
starting at $z\in \partial B_{d}(y,h_{k'-1})$,
let $\bar{\tau}=\inf \{t>0\, :\,X_t\in \partial B_{d}(y,h_{k})\}$ and let 
$F$ denote an event measurable with respect to the path of the Brownian motion
inside $B_d(y,h_{k'})$, prior to $\bar \tau$. Let $\wh{\tau}_{k'}=\inf \{t>0\, :\,X_t\in \partial B_{d}(   0, h_{k'})\}$ and $D_{k'}=\{ \wh{\tau}_{k'}<\bar{\tau}\}$. 
Then,  
uniformly in $z\in \partial B_{d}(y,h_{k'-1})$, $v\in \partial B_{d}(y,h_{k})$ and $y$, 
\begin{equation}
\label{m1.5eua}
  \left|\frac{\Ebm^z (F, D_{k'}\big |\,X_{\bar{\tau}}=v)}{\Ebm^z 
      (F, D_{k'}   )}-1\right|
\leq  4.8 (k'-k)\frac{h_{k'-1}}{h_k}.
\end{equation}
and
\begin{equation}
\label{m1.5euaa}
  \left|\frac{\Ebm^z (F\big |\, D_{k'},X_{\bar{\tau}}=v)}{\Ebm^z 
      (F\big |\, D_{k'}   )}-1\right|
\leq  9.8 (k'-k)\frac{h_{k'-1}}{h_k}.
\end{equation}

\noindent
Further, let
$Z_{k'}=\trav{1}{y}{k}{k'-1}{y}{k'-1}{k'}$ 
denote the number of excursions of
the path  from 
$\partial B_{d}(y,h_{k'-1})\to\partial B_{d}(y,h_{k'})$, prior to
$\bar{\tau}=\inf \{t>0\, :\,X_t\in \partial B_{d}(y,h_{k})\}$, then,  
uniformly in $z\in \partial B_{d}(y,h_{k'-1})$, $v\in \partial B_{d}(y,h_{k})$, $j\geq 0$ and $y$, 
\begin{equation}
\label{m1.5eu}
  \left|\frac{\Ebm^z (F ;Z_{k'}=j \,\big |\,X_{\bar{\tau}}=v)}{\Ebm^z 
      (F;Z_{k'}=j )}-1\right|
\leq  4.8 (k'-k)\frac{h_{k'-1}}{h_k},
\end{equation}
and
\begin{equation}
\label{m1.5euc}
  \left|\frac{\Ebm^z (F \,\big |\,Z_{k'}=j, X_{\bar{\tau}}=v)}{\Ebm^z 
      (F\,\big |\,Z_{k'}=j )}-1\right|
\leq  9.8 (k'-k)\frac{h_{k'-1}}{h_k}.
\end{equation}

\noindent 
Finally, uniformly in $x\in \partial B_{d}(y,h_{k+1})$, $v\in \partial B_{d}(y,h_{k})$, $j\geq 1$, and $y$,  
\begin{equation}
\label{m1.5eucc}
  \left|\frac{\Ebm^x (F \,\big |\,Z_{k'}=j, X_{\bar{\tau}}=v)}{\Ebm^{\la_{k'-1}} 
      (F\,\big |\,Z_{k'}=j )}-1\right|
\leq  28.2 (k'-k)\frac{h_{k'-1}}{h_k},
\end{equation}
where $\la_{k'-1}$ denotes uniform measure on $\partial B_{d}(y,h_{k'-1})$,
and
\begin{equation}
\label{m1.5eucco}
\left|\frac{\Ebm^x (F\,\big |\,D_{k'}, X_{\bar{\tau}}=v)}{\Ebm^{\la_{k'-1}} 
      (F\big |\,D_{k'})}-1\right|
\leq  28.2 (k'-k)\frac{h_{k'-1}}{h_k}.
\end{equation}
\end{lemma}

In words, conditioning by the endpoint of the excursion at level $k$ has only minor influence on the probability
of events involving pieces of excursions inside the ball of radius $h_{k'}$.

\noindent{\bf Proof of Lemma \ref{lemprobends}:}
We first prove \eqref{m1.5eu} and  \eqref{m1.5euc}. Toward the proof of \eqref{m1.5eu}, it is enough to prove that
    \begin{equation}
\Big| \Ebm^z (F ;Z_{k'}=j \,\big |\,X_{\bar{\tau}}=v)-
\Ebm^z (F; Z_{k'}=j )\Big|\leq 
     4.8 (k'-k)\frac{h_{k'-1}}{h_k}\Ebm^z (F;
    Z_{k'}=j ).\label{92.1}
\end{equation}
    Without loss of generality we can take $y=0$.  
Fixing $k \geq 0$ and
$z \in \partial B_{d}(0, h_{k'-1})$,  
consider  a positive continuous function $g$ on $\partial B_{d}(0,h_{k})$.
Let $\bar{\tau} = \inf \{ t:\, w_t \in \partial B_{d}(0,h_{k})\}$,
$\tau_0=0$ and for $i=0,1,\ldots$ define
\begin{eqnarray*}
\tau_{2i+1} & = & \inf \{ t \geq \tau_{2i}  :\; X_t \in \partial
B_{d}(0,h_{k'}) \cup \partial B_{d}(0,h_{k}) \} \\
\tau_{2i+2} & = &
 \inf \{ t \ge \tau_{2i+1} :\; X_t \in \partial B_{d}(0,h_{k'-1})\}\,.
\end{eqnarray*}
 Then, by the strong Markov property at $\tau_{2j}$,
\[
\Ebm^z (g(X_{\bar\tau}); F
    ;Z_{k'}=j  )
=
\Ebm^z \big[ \Ebm^{X_{\tau_{2j}}} (g(X_{\bar\tau}); Z_{k'}=0) ;F,
Z_{k'}=j, {\bar\tau} \geq \tau_{2j} \big]
\,.
\]
and, substituting $g=1$,
\bea
&&
\Pbm^z \(F,Z_{k'}=j \)\nn
=
\Ebm^z \left[ \Pbm^{X_{\tau_{2j}}} ( Z_{k'}=0) ;F,
Z_{k'}=j, {\bar\tau} \geq \tau_{2j} \right]\nn
 \;.
\eea
Consequently,
\begin{eqnarray*}
&&\Pbm^z \(F, Z_{k'}=j \)
\inf_{x\in \partial B_{d}(0,h_{k'-1})}
\frac{\Ebm^{x} \( g(X_{\bar\tau});\, Z_{k'}=0 \)} {\Pbm^x \(Z_{k'}=0\)}
\leq \Ebm^z [ g(X_{\bar\tau}); F, Z_{k'}=j ]\\ 
&& \leq 
\Pbm^z \(F, Z_{k'}=j \)
\sup_{x\in \partial B_{d}(0,h_{k'-1})}
\frac{\Ebm^{x} \( g(X_{\bar\tau});\, Z_{k'}=0 \)} {\Pbm^x \(Z_{k'}=0\)}
\,,
\end{eqnarray*}
and, using again the strong Markov property at time $\tau_{2}$,
\begin{eqnarray*}
\Ebm^{x} \( g(X_{\bar\tau});\, Z_{k'}=0 \) &=&
\Ebm^{x} \( g(X_{\bar\tau}) \) -
\Ebm^{x} \( \Ebm^{X_{\tau_2}} (g(X_{\bar\tau}));\, Z_{k'} \geq 1 \) \\
&\leq& \Ebm^{x} \( g(X_{\bar\tau}) \) -
\Pbm^{x} (Z_{k'} \geq 1) \inf_{y\in \partial B_{d}(0,h_{k'-1})} \Ebm^{y} (g(X_{\bar\tau})),
\end{eqnarray*}
with the reversed inequality if the $\inf$ is replaced by $\sup$.
Writing $p=\Pbm^{x} (Z_{k'} \geq 1) =1-1/(k'-k)$
whenever $x\in \partial B_{d}(0,h_{k'-1})$, c.f. \cite[(2.7)]{BRZ2}.
it thus
follows that
\begin{equation}
\label{eq10.1}
 \Ebm^z \left[  g(X_{\bar\tau}); F, Z_{k'}=j  \right] 
\leq
\Pbm^z\(F, Z_{k'}=j  \)
\Ebm^z\( g(X_{\bar\tau})\)
(1-p)^{-1}
\Big\{ \frac{\sup_{x\in \partial B_{d}(0,h_{k'-1})} \Ebm^x( g(X_{\bar\tau}))} 
{\inf_{y\in \partial B_{d}(0,h_{k'-1})} \Ebm^y( g(X_{\bar\tau}))} -p\Big\}, 
\end{equation}
with the reversed inequality if the $\inf$ and $\sup$ are interchanged.

Let $p_{B_{d}(   0, h_{k})}(   z,x)$ denote 
 the Poisson kernel  for $B_{d}(   0, h_{k})\subseteq \bS^{ 2}$,  see \cite[(2.8)]{BRZ2}. Then,
$$
\Ebm^{z'} g (X_{\bar{\tau} })
= \int_{\partial B_{d}(0,h_{k})}  p_{B_{d}(0,h_{k})}(   z',u) g(u) du.
$$
Therefore,
we get the Harnack inequality
\bea
 \frac{\sup_{x\in \partial B_{d}(0,h_{k'-1})} \Ebm^x( g(X_{\bar\tau}))}
{\inf_{y\in \partial B_{d}(0,h_{k'-1})} \Ebm^y( g(X_{\bar\tau}))}& \leq &
\frac{\max_{x\in \partial B_{d}(0,h_{k'-1}),u\in \partial B_{d}(0,h_{k })}p_{B_{d}(0,h_{k})}( x,u)}
{\min_{y\in \partial B_{d}(0,h_{k'-1}),u\in \partial B_{d}(0,h_{k })} p_{B_{d}(0,h_{k})}( y,u)}\nn \\
&=&\frac{\max_{y\in \partial B_{d}(0,h_{k'-1}),u\in \partial B_{d}(0,h_{k })}\sin^{ 2} ( d(  u,y)/2)}
{\min_{x\in \partial B_{d}(0,h_{k'-1}),u\in \partial B_{d}(0,h_{k })} \sin^{ 2} ( d(  u,x)/2)}
\nn \\
&\leq&
\frac{\sin^{ 2}((h_k+h_{k'-1})/2) }{\sin^{ 2}((h_k-h_{k'-1})/2) }=:B_{k,k'-1}^2.\label{eq10.2}
\eea
Writing $\al=(h_k-h_{k'-1})/2$, $\bb=h_{k'-1}$  we have
\be
B_{k,k'-1}=\frac{\sin (\al+\bb) }{\sin (\al) }=\cos(\bb)+\frac{\sin ( \bb)\cos (\al) }{\sin (\al) }.
\label{eq10.2a}
\ee
Using the bounds  $\sin (\alpha)\geq .9\alpha $, $\sin (\bb)\leq \bb $, $\cos (\al) , \cos \beta\leq 1$, which are valid for  $0<\alpha,\beta <0.01$, we obtain that 
$B_{k,k'-1}$ is bounded above by 
\begin{equation}
1+ \bb/.9\al=1+2 h_{k'-1}/.9(h_k-h_{k'-1}).\label{jr,b}
\end{equation}
Since $k' \geq k+10$ we have that  $.99 h_k \leq h_k-h_{k'-1}$, so that 
(\ref{jr,b}) is bounded by 
$1+2.3h_{k'-1}/h_k$, and therefore $B_{k,k'-1}^2$ 
is bounded above by 
$1+4.8h_{k'-1}/h_k$ while its reciprocal is bounded below by
$1-4.8h_{k'-1}/h_k$.
Substituting this bound into 
 (\ref{eq10.1}) (and its reversed version with inequality and sup/inf interchanged), and using the value of $p$, yields 
 (\ref{92.1}). 
 
 For (\ref{m1.5euc}) we have by (\ref{92.1}), first for $F$ and then replacing $F$ by $1$
     \begin{equation}
 \Ebm^z (F ;Z_{k'}=j \,\big |\,X_{\bar{\tau}}=v) \leq 
     \(1+4.8 (k'-k)\frac{h_{k'-1}}{h_k}\)\Ebm^z (F; Z_{k'}=j ).
     \label{92.1m}
\end{equation}
and
     \begin{equation}
 \Ebm^z (Z_{k'}=j \,\big |\,X_{\bar{\tau}}=v) \geq 
     \(1-4.8 (k'-k)\frac{h_{k'-1}}{h_k}\)\Ebm^z (Z_{k'}=j ).
\label{92.1mm}
\end{equation}
This gives the upper bound 
\begin{equation}
\Ebm^z (F \,\big |\,Z_{k'}=j, X_{\bar{\tau}}=v)= { \Ebm^z (F ;Z_{k'}=j \,\big |\,w_{\bar{\tau}}=v) \over  \Ebm^z (Z_{k'}=j \,\big |\,w_{\bar{\tau}}=v)}
\leq       \(1+9.8 (k'-k)\frac{h_{k'-1}}{h_k}\)\Ebm^z (F \,\big |\,Z_{k'}=j). \label{9.15jay}
\end{equation}
The lower bound is obtained similarly.

 \eqref{m1.5eua}
follows from \eqref{m1.5eu} by multiplying both sides by
$\Ebm^z (F;
    Z_{k'}=j )$ and summing over $j\geq 1$, and \eqref{m1.5euaa} will follow from \eqref{m1.5eua} exactly  in the way  that  \eqref{m1.5euc} followed from  \eqref{m1.5eu}.

We turn to the proof of  \eqref{m1.5eucc}. Let $\wh{\tau}_{k'-1}=\inf \{t>0\, :\,X_t\in \partial B_{d}(   0, h_{k'-1})\}$. We will need to show that
the law of $X_{\wh\tau}$, started at $x\in \partial B_d(0, h_{k+1})$ and
conditioned on the event $D_{k'-1}=\{\wh \tau_{k'-1}< \bar \tau\}$, is close to uniform. Toward this end, we begin  with an unconditional statement. 
Let $0^{\ast}$ denote the antipode of $0\in \bS^{ 2}$. Note that $\wh{\tau}_{k'-1}=\inf \{t>0\, :\,X_t\in \partial B_{d}(   0, h_{k'-1})=\partial B_{d}(  0^{\ast}, \pi-h_{k'-1})\}$.
 Then for any positive  function $h$ on  $\partial B_{d}(   0, h_{k'-1})=\partial B_{d}(  0^{\ast}, \pi-h_{k'-1})$,  and $x\in \partial B_{d}(   0, h_{k+1 })=\partial B_{d}(  0^{\ast}, \pi-h_{k+1 })$, writing   $\bar h=h(X_{\wh \tau_{k'-1}})$,
\be
\Ebm^{x}\bar h
= \int_{\partial B_{d}(   0, h_{k'-1})}  p_{B_{d}(0^{\ast}, \pi-h_{k'-1})}( x,u) h(u) d\la_{k'-1}(u),\label{poissonk}
\ee
where $\la_{k'-1}$ is uniform measure on $\partial B_{d}(   0, h_{k'-1})$. We have,   see \cite[(2.8)]{BRZ2},
  \bea
p_{B_{d}(0^{\ast}, \pi-h_{k'-1})}( x,u)&=&\frac{\sin^{ 2} (\pi/2-h_{k'-1}/2)-\sin^{ 2} (\pi/2-h_{k+1}/2)}{\sin^{ 2} ( d( u,x)/2)} \nn\\
&=&\frac{\sin^{ 2} (h_{k+1}/2)-\sin^{ 2} (h_{k'-1}/2)}{\sin^{ 2} ( d( u,x)/2)}.\label{poissonj}
  \eea
For the upper bound, using $\sin^{ 2}(a)-\sin^{ 2}(b)=\sin (a+b)\sin (a-b)$
  \begin{equation}
 p_{B_{d}(0^{\ast}, \pi-h_{k'-1})}( x,u)\leq \frac{\sin  (  ( h_{k+1}+h_{k'-1})/2)}{\sin  (  ( h_{k+1}-h_{k'-1})/2)}, \label{poissonj1}
  \end{equation}
which is   $B_{k+1, k'-1}$, see (\ref{eq10.2a}),   so that 
   \begin{equation}
 p_{B_{d}(0^{\ast}, \pi-h_{k'-1})}( x,u)\leq   1+2.3 \(\frac{   h_{k'-1}  }{ h_{k+1}  }\). \label{poissonj3}
  \end{equation} 
  
  For the lower bound, using $\sin^{ 2}(a)-\sin^{ 2}(b)=\sin (a+b)\sin (a-b)$ again we  note that
    \begin{equation}
p_{B_{d}(0^{\ast}, \pi-h_{k'-1})}( x,u)\geq \frac{\sin  (  ( h_{k+1}-h_{k'-1})/2)}{\sin  (  ( h_{k+1}+h_{k'-1})/2)}={1 \over B_{k+1, k'-1}}.  \label{poissonj5}
  \end{equation}
It follows  that
     \begin{equation}
p_{B_{d}(0^{\ast}, \pi-h_{k'-1})}( x,u)\geq 1-2.4 \(\frac{   h_{k'-1}  }{ h_{k+1}  }\).  \label{poissonj8}
  \end{equation}
   Combining this with (\ref{poissonk}), we conclude that with $h^{\lambda_{k'-1}}=\int h(u) d\lambda_{k'-1}(u)$,
 %
  \begin{equation}
  \label{eq-monday1}
  \big|\Ebm^x(\bar h)-h^{\lambda_{k'-1}}\big| \leq 2.4 \(\frac{   h_{k'-1}  }{ h_{k+1}  }\) h^{\lambda_{k'-1}}=:\delta_{k'-1,k+1}h^{\lambda_{k'-1}}.
  \end{equation}
  
  We next turn to proving the analogous  conditional statement. Recall the event $D_{k'-1}=\{\wh \tau_{k'-1}< \bar \tau\}$.
By the Markov property $\Ebm^x(\bar h, \{\bar \tau<\wh \tau_{k'-1}\})=\Ebm^x(\Ebm^{X_{\bar \tau}}(\bar h), \{\bar \tau<\wh \tau_{k'-1}\})$, and hence by  \eqref{eq-monday1}  with $k+1$ replaced by $k$
  \begin{equation}
  \label{eq-monday2}
  \big|\Ebm^x(\bar h|D_{k'-1}^{c})- h^{\lambda_{k'-1}}\big| \leq\delta_{k'-1,k}h^{\lambda_{k'-1}}.
  \end{equation}
  Note that  by \cite[(2.7)]{BRZ2} we have that
$\Pbm^x(D_{k'-1})= 1/(k'-k-1)=:\delta'_{k',k}$ and is independent of $x$. We have
  \[ \Ebm^x(\bar h)=\Ebm^x(\bar h|D)\Pbm^x(D_{k'-1})+\Ebm^x(\bar h|D_{k'-1}^c)\Pbm^x(D_{k'-1}^c)\]
  and therefore, using the last display and \eqref{eq-monday1} and \eqref{eq-monday2} in the inequality,
  \bea
  \big|\Ebm^x(\bar h|D_{k'-1})-h^{\lambda_{k'-1}}\big|\label{eq-monday3}
  &&=\big|\frac{1}{\Pbm^{x}(D_{k'-1})} \left(\Ebm^x(\bar h)-\Ebm^x(\bar h|D_{k'-1}^c)\right)+\Ebm^x(\bar h|D_{k'-1}^c) -h^{\lambda_{k'-1}}\big|\\
  && \leq \left(\frac{(\delta_{k'-1,k+1}+ \delta_{k'-1,k}  )}{\delta'_{k',k}} +\delta_{k'-1,k} \right) h^{\lambda_{k'-1}}\leq 9.6(k'-k)\(\frac{   h_{k'-1}  }{ h_{k}  }\)h^{\lambda_{k'-1}}.\nn
  \eea
This gives the desired conditional estimate.

To obtain \eqref{m1.5eucc} we first  consider
$H=F1_{\{Z_{k'}=j\}} g(X_{\bar{\tau}})$. For all $j\geq 1$ we have that $H=H\circ\th_{\wh{\tau}_{k'-1}}\,\,\,1_{\{\wh{\tau}_{k'-1}<\bar{\tau}\}}$.          Hence
 by the strong Markov property, 
$$\Ebm^x (F1_{\{Z_{k'}=j\}}g(X_{\bar{\tau}}))=\Ebm^x \(\Ebm^{X_{\wh{\tau}_{k'-1}}  } \(F1_{\{Z_{k'}=j\}} g(X_{\bar{\tau}})\), \,D_{k'-1}\). $$
 Let $h(y)= \Ebm^{y  } \(F1_{\{Z_{k'}=j\}} g(X_{\bar{\tau}})\)$. (\ref{eq-monday3}) then shows that
\begin{eqnarray}
&&|\Ebm^x (F1_{\{Z_{k'}=j\}}g(X_{\bar{\tau}}))-\Ebm^{\la_{k'-1}} (F1_{\{Z_{k'}=j\}}g(X_{\bar{\tau}}))\Pbm^x(D_{k'-1})|
\label{1922}\\
&&   \leq 9.6(k'-k)\(\frac{   h_{k'-1}  }{ h_{k}  }\)\Ebm^{\la_{k'-1}} (F1_{\{Z_{k'}=j\}}g(X_{\bar{\tau}}))\Pbm^x(D_{k'-1}).\nonumber
\end{eqnarray}
Using this also with $F=1$, proceeding as in the proof of (\ref{m1.5euc}),  and then letting $g\to \de_{v}$ gives  
\begin{eqnarray}
&&|\Ebm^x (F\,\big |\,Z_{k'}=j, X_{\bar{\tau}}=v)-\Ebm^{\la_{k'-1}} (F\,\big |\,Z_{k'}=j, X_{\bar{\tau}}=v)|
\label{1923}\\
&&   \leq 19.3(k'-k)\(\frac{   h_{k'-1}  }{ h_{k}  }\)\Ebm^{\la_{k'-1}} (F\,\big |\,Z_{k'}=j, X_{\bar{\tau}}=v),\nonumber
\end{eqnarray}
and (\ref{m1.5eucc}) then follows from (\ref{m1.5euc}).

To obtain \eqref{m1.5eucco} let $D_{k'}=\{\wh \tau_{k'}< \bar \tau\}$ where $\wh \tau_{k'} =\inf \{t>0\, :\,X_t\in \partial B_{d}(   0, h_{k'})\}$. Summing (\ref{1922}) over $j\geq 1$
 shows that
\begin{eqnarray}
&&|\Ebm^x (F1_{D_{k'}}g(X_{\bar{\tau}}))-\Ebm^{\la_{k'-1}} (F1_{D_{k'}}g(X_{\bar{\tau}}))\Pbm^x(D_{k'-1})|
\label{k1922}\\
&&   \leq 9.6(k'-k)\(\frac{   h_{k'-1}  }{ h_{k}  }\)\Ebm^{\la_{k'-1}} (F1_{D_{k'}}g(X_{\bar{\tau}}))\Pbm^x(D_{k'-1}),\nonumber
\end{eqnarray}
and as in (\ref{1923})  we then obtain
\begin{eqnarray}
&&|\Ebm^x (F\,\big |\,D_{k'}, X_{\bar{\tau}}=v)-\Ebm^{\la_{k'-1}} (F\,\big |\,D_{k'}, X_{\bar{\tau}}=v)
\label{k1923}\\
&&   \leq 19.3(k'-k)\(\frac{   h_{k'-1}  }{ h_{k}  }\)\Ebm^{\la_{k'-1}} (F\,\big |\,D_{k'}, X_{\bar{\tau}}=v).\nonumber
\end{eqnarray}
Using (\ref{m1.5euaa}) then gives the upper bound in \eqref{m1.5eucco} and the lower bound is proven similarly.
\qed

 \noindent{\bf Proof of Lemma \ref{recursionends}:} This follows from Lemma \ref{lemprobends}
in the same manner as \cite[Lemma 7.3]{DPRZthick} was derived from  \cite[Lemma 7.4]{DPRZthick}: using the strong Markov property, each excursion 
from $\partial B_{d}(y,h_{k'})$ to $\partial B_{d}(y,h_{k-1})$ will contribute a multiplicative
factor $(1+50 (k'-k+1)h_{k'-1}/{h_{k-1}})$ to the probability.

In more detail, for the $i$'th excursion from $\partial B_{d}(y,h_{k})$ to $\partial B_{d}(y,h_{k-1})$
let $x_{i}, v_{i}$ denote the the starting and endpoint of the excursion (note that these are random variables!)
and 
let $Z_{k'}^{i}$ denote the number of excursions from $\partial B_{d}(y,h_{k'-1})$ to $\partial B_{d}(y,h_{k'})$. Set
\begin{equation}
S=\left(1+O\( (k'-k+1)\frac{h_{k'-1}}{h_{k-1}}\)\right)^{m_{k'}}\label{}
\end{equation}
with $| O\( (k'-k+1)\frac{h_{k'-1}}{h_{k-1}}\)|\leq 50 (k'-k+1)\frac{h_{k'-1}}{h_{k-1}}$.

To prove (\ref{m1.2e}) we first show that for all $j_{i}; i\in [1, n]$ with $\sum_{i=1}^{n}j_{i}= m_{k'}$,
\bea
&&\Pbm \(\trav{n}{y}{k-1}{k}{y}{l-1}{l} = m_l ;l=k'+1,\ldots L \, \,\big |\,x_{i}, v_{i}, Z^{i}_{k'}=j_{i}; i\in [1, n]\)\nn\\
&&=S\Pbm(\trav{n}{y}{k-1}{k}{y}{l-1}{l} = m_l ;l=k'+1,\ldots L \,
|\, \trav{n}{y}{k-1}{k}{y}{k'-1}{k'} = m_{k'} ).\label{indf.01}
\eea

Using the fact that $\sum_{i=1}^{n}j_{i}= m_{k'}$ and all excursion counts are centered at $y$ we have
\begin{eqnarray}
&&\Pbm \(\trav{n}{y}{k-1}{k}{y}{l-1}{l} = m_l ;l=k'+1,\ldots L \,\,\big |\,Z^{i}_{k'}=j_{i}; i\in [1, n]\)
\label{indf.1n}
\\
&&=\Pbm(\trav{n}{y}{k-1}{k}{y}{l-1}{l} = m_l ;l=k'+1,\ldots L \,
|\, \trav{n}{y}{k-1}{k}{y}{k'-1}{k'} = m_{k'} ),
\nonumber
\end{eqnarray}
so to prove (\ref{indf.01}) it suffices to show that for all $j_{i}; i\in [1, n]$ with $\sum_{i=1}^{n}j_{i}= m_{k'}$,
\bea
&&\Pbm \(\trav{n}{y}{k-1}{k}{y}{l-1}{l} = m_l ;l=k'+1,\ldots L \, \,\big |\,x_{i}, v_{i}, Z^{i}_{k'}=j_{i}; i\in [1, n]\)\nn\\
&&=S \Pbm \(\trav{n}{y}{k-1}{k}{y}{l-1}{l} = m_l ;l=k'+1,\ldots L \,\,\big |\,Z^{i}_{k'}=j_{i}; i\in [1, n]\).\label{indf.1}
\eea

 Let $F_{i}=1$ if $Z^{i}_{k'}=0$. Otherwise let $\trav{(i)}{y}{k-1}{k}{y}{l-1}{l}$ be the number of excursions    from
	$\partial B_{d}(y,h_{l-1})\to\partial B_{d}(y,h_{l})$ during the  $i$'th
 excursion    from
	$\partial B_{d}(y,h_{k })\to\partial B_{d}(y,h_{k-1})$, so that $\trav{(i)}{y}{k-1}{k}{y}{k'-1}{k'}=Z^{i}_{k'}=j_{i}$,  and set
 \begin{equation}
 F_{i}=\{\trav{(i)}{y}{k-1}{k}{y}{l-1}{l} = m_{l,i} ;l=k'+1,\ldots L\}\label{}
 \end{equation}
with  $\sum_{i=1}^{n}m_{l,i}= m_{l};\,l=k'+1,\ldots L$.

To prove (\ref{indf.1}) it suffices to show that
\bea
&&\Pbm \(\prod_{i=1}^{n} F_{i} \, \,\big |\,x_{i}, v_{i}, Z^{i}_{k'}=j_{i}; i\in [1, n]\)\label{indf.1l}\\
&&=\left(1+O\( (k'-k+1)\frac{h_{k'-1}}{h_{k-1}}\)\right)^{\sum_{i=1}^{n}j_{i}}\Pbm \(\prod_{i=1}^{n} F_{i} \,\,\big |\,Z^{i}_{k'}=j_{i}; i\in [1, n]\),\nn
\eea
for all choices of $m_{l,i}$ satisfying $\sum_{i=1}^{n}m_{l,i}= m_{l};\,l=k'+1,\ldots L$, since summing over all choices gives (\ref{indf.1})

By the Markov property and (\ref{m1.5eucc}) with $k$ replaced by $k-1$ we have 
\bea
&&\Ebm \(\prod_{i=1}^{n}F_{i} \,\big |\,x_{i}, v_{i}, Z^{i}_{k'}=j_{i}; i\in [1, n]\)\label{indf.3}\\
&&=\Ebm \(\prod_{i=1}^{n-1}F_{i} \,\big |\,x_{i}, v_{i}, Z^{i}_{k'}=j_{i}; i\in [1, n-1]\)\Ebm^{x_{n}} \(F_{n} \,\big |\, v_{n}, Z^{n}_{k'}=j_{n}\)\nn\\
&&\leq \left(1+28.2\( (k'-k+1)\frac{h_{k'-1}}{h_{k-1}}\)\right)^{1_{\{j_{n}>0\}}}\nn\\
&&\hspace{.3 in}\times \Ebm \(\prod_{i=1}^{n-1}F_{i} \,\big |\,x_{i}, v_{i}, Z^{i}_{k'}=j_{i}; i\in [1, n-1]\)\Ebm^{\la_{k'-1}} \( F_{n} \,\big |\, Z^{n}_{k'}=j_{n}\),\nn
\eea
and by induction
\begin{equation}
\Ebm \(\prod_{i=1}^{n}F_{i} \,\big |\,x_{i}, v_{i}, Z^{i}_{k'}=j_{i}; i\in [1, n]\)\label{indf.4}
\!\leq \!\left(1+28.2\( (k'-k+1)\frac{h_{k'-1}}{h_{k-1}}\)\right)^{\sum_{i=1}^{n}j_{i}}\prod_{i=1}^{n}\Ebm^{\la_{k'-1}} \( F_{i} \,\big |\, Z^{i}_{k'}=j_{i}\).
\end{equation}
Similarly, by the analogue of (\ref{1923}), but without conditioning on the endpoint, which follows from (\ref{1922}) with $g\equiv1$, we have
\bea
&&\hspace{.3 in}\Ebm \(\prod_{i=1}^{n}F_{i} \,\big |\,  Z^{i}_{k'}=j_{i}; i\in [1, n]\)\label{indf.5}\\ 
&&\geq \left(1-19.3\( (k'-k+1)\frac{h_{k'-1}}{h_{k-1}}\)\right)^{\sum_{i=1}^{n}j_{i}}\prod_{i=1}^{n}\Ebm^{\la_{k'-1}} \( F_{i} \,\big |\, Z^{i}_{k'}=j_{i}\).\nn
\eea
Together this gives the upper bound in  (\ref{indf.1l}) and the lower bound is proven similarly.

We now explain why (\ref{indf.01}) implies (\ref{m1.2e}).
Let $\mathcal{G}$ denote the sigma-algebra generated by the 
starting and endpoints $(x_i,v_i)$  of the the $i$'th excursion from $\partial B_{d}(y,h_{k})$ to $\partial B_{d}(y,h_{k-1})$ and set
\begin{equation}
\wt{\mathbb{P}}(\cdot)=\mathbb{P}\left(\cdot\,\Big
| \, \trav{n}{0}{  k-1}{  k}{0}{k'-1}{k'} = m, \, \mathcal{G}\right ). \label{}
\end{equation}
By the law of total probability, with $A=\{\trav{n}{y}{k-1}{k}{y}{l-1}{l} = m_l ;l=k'+1,\ldots L\}$,
\begin{equation}
\wt{\mathbb{P}}(A)=\sum_{\{(j_{1},\ldots, j_{n})\,|\, \sum_{i'=1}^{n}j_{i'}=m\}}\wt{\mathbb{P}}(A\,|\,Z^{i'}_{k'}=j_{i'}; i'\in [1, n])
 \label{m1.2e1}
\wt{\mathbb{P}}(Z^{i'}_{k'}=j_{i'}; i'\in [1, n]).
\end{equation}
 But 
 \begin{eqnarray}
 &&\wt{\mathbb{P}}(A\,|\,Z^{i'}_{k'}=j_{i'}; i'\in [1, n])
=\frac{\wt{\mathbb{P}}(A,\,Z^{i'}_{k'}=j_{i'}; i'\in [1, n])}
{ \wt{\mathbb{P}}(Z^{i'}_{k'}=j_{i'}; i'\in [1, n)}
 \label{m1.2e2}
 \\
 &&=\frac{{\mathbb{P}}(A,\,Z^{i'}_{k'}=j_{i'}; i'\in [1, n] \,| \, \mathcal{G})}{ {\mathbb{P}}(Z^{i'}_{k'}=j_{i'}; i'\in [1, n]\,| \, \mathcal{G} )}={\mathbb{P}}(A\,|\,Z^{i'}_{k'}=j_{i'}; i'\in [1, n], \, \, \mathcal{G})
= S \mathbb{P}\left(A\,| \, \trav{n}{0}{  k-1}{  k}{0}{k'-1}{k'} = m\right)   \nn
 \end{eqnarray}
 by (\ref{indf.01}). Since this is independent of $Z^{i'}_{k'}=j_{i'}; i'\in [1, n]$, 
  (\ref{m1.2e}) now follows by plugging (\ref{m1.2e2}) into (\ref{m1.2e1}).
  \qed

The following lemma is 
a small variant of what we already proved. Let $\trav{n}{y}{k-1}{k}{y}{k'}{k-1}$ be the number of excursions    from
	$\partial B_{d}(y,h_{k'})\to\partial B_{d}(y,h_{k-1})$ during the first $n$
 excursions    from
	$\partial B_{d}(y,h_{k })\to\partial B_{d}(y,h_{k-1})$.

\begin{lemma}\label{recursionends2}
For any $k '>k+10\geq 11$ and $n \ge 1$,
let $\mathcal{A}_{k'}$ denote an event, measurable on the excursions of the Brownian motion 
inside $ B_{d}(y,h_{k'})$  during the first $n$
 excursions    from
	$\partial B_{d}(y,h_{k })\to\partial B_{d}(y,h_{k-1})$. Then, for any $w\in B^c_{d}(y,h_{k-1})$
\be\Pbm^w(\mathcal{A}_{k'} \,
| \, \trav{n}{y}{k-1}{k}{y}{k'}{k-1} = m_{k'}, \, {\mathcal G}_k^y)
=
\left(1+ O\( (k'-k+1)\frac{h_{k'-1}}{h_{k-1}}\)\right)^{m_{k'}} \Pbm^w(\mathcal{A}_{k'}\,
| \, \trav{n}{y}{k-1}{k}{y}{k'}{k-1} = m_{k'})  \,.
\label{m1a.2e0}
\ee
In particular, 
%
for  all $m_l ;l=k',\ldots L$, and all $y \in \bS^2$, 
\bea
&& \Pbm^w(\trav{n}{y}{k-1}{k}{y}{l-1}{l} = m_l ;l=k'+1,\ldots L \,
|\, \trav{n}{y}{k-1}{k}{y}{k'}{k-1} = m_{k'}, \, {\mathcal G}_k^y) 
\nonumber \\ 
&&=
\left(\!1+ O\( (k'-k+1)\frac{h_{k'-1}}{h_{k-1}}\)\!\right)^{m_{k'}} 
\label{m1a.2e}
\!\Pbm^w(\trav{n}{y}{k-1}{k}{y}{l-1}{l} = m_l ;l=k'+1,\ldots L \,
|\, \trav{n}{y}{k-1}{k}{y}{k'}{k-1} = m_{k'} ).
\eea
\end{lemma}
\noindent
In Lemma \ref{recursionends2}
we have  $| O\( (k'-k+1)\frac{h_{k'-1}}{h_{k-1}}\)|\leq 50\( (k'-k+1)\frac{h_{k'-1}}{h_{k-1}}\)$.

 \noindent{\bf Proof of Lemma \ref{recursionends2}:} This is similar to the proof of  Lemma \ref{recursionends}, except that
 no specific form, and in particular no rotation invariance, of the event 
${\mathcal A}_{k'}$ is assumed.
 (Recall that  no specific form   was imposed on $F$
in Lemma \ref{lemprobends}.)

In more detail, for the $j$'th excursion $X_{j}$ from $\partial B_{d}(y,h_{k'})$ to $\partial B_{d}(y,h_{k-1})$
let $x_{j}, v_{j}$ denote the the starting and endpoint of the excursion from $\partial B_{d}(y,h_{k})$ to $\partial B_{d}(y,h_{k-1})$ which contains $X_{j}$. 
It suffices to show that
\bea
&&
\Pbm^w \({\mathcal A}_{k'} \,\big |\,\trav{n}{y}{k-1}{k}{y}{k'}{k-1} = m_{k'};\,\, x_{j}, v_{j},  ; j\in [1, m_{k'}]\)\label{indfa.1}\\
&&=\left(1+O\( (k'-k+1)\frac{h_{k'-1}}{h_{k-1}}\)\right)^{m_{k'}}\Pbm^w \({\mathcal A}_{k'} \,\big |\,\trav{n}{y}{k-1}{k}{y}{k'}{k-1} = m_{k'}\),\nn
\eea
with $| O\( (k'-k+1)\frac{h_{k'-1}}{h_{k-1}}\)|\leq 50 (k'-k+1)\frac{h_{k'-1}}{h_{k-1}}$. 

Let $F=\prod_{i=1}^{m_{k'}}F_{i}$ where $F_{i}$ depends on the excursions inside $B_{d}(y,h_{k'})$ during the $i$'th excursion from $\partial B_{d}(y,h_{k'})$ to $\partial B_{d}(y,h_{k-1})$. 
Since such events generate ${\mathcal A}_{k'}$, to prove (\ref{indfa.1}) it suffices to show that
\bea
&&\Ebm^w \(F \,\big |\,\trav{n}{y}{k-1}{k}{y}{k'}{k-1} = m_{k'};\,\, x_{i}, v_{i},  ; i\in [1, m_{k'}]\)\label{indfa.2}\\
&&=\left(1+O\( (k'-k+1)\frac{h_{k'-1}}{h_{k-1}}\)\right)^{m_{k'}}\Ebm^w \(F \,\big |\,\trav{n}{y}{k-1}{k}{y}{k'}{k-1} = m_{k'}\).\nn
\eea
Let $\bar{\tau_i}$ denote the end of the $i$th excursion
from $\partial B_d(y,h_k)$ to $\partial B_d(y, h_{k-1})$ that reached
$\partial B_d(y,h_{k'})$, let $\wh{\tau}_{k',i}$ be defined analogously to
$\wh{\tau}_{k'}$ for that excursion, let $D_{k',i}=\{\wh{\tau}_{k',i}<\bar{\tau}_i<\ff\}$, and
let $x_i,v_i$ denote the starting and end points of that excursion.
By the Markov property and \eqref{m1.5eucco}
with $k$ replaced by $k-1$ we have 
\bea
&&\Ebm^w \(\prod_{i=1}^{m_{k'}}F_{i} \,\big |\,\trav{n}{y}{k-1}{k}{y}{k'}{k-1} = m_{k'};\,\, x_{j}, v_{j},  ; j\in [1, m_{k'}]\)\label{indfa.3}\\
&&=\Ebm^w \(\prod_{i=1}^{m_{k'}-1}F_{i} \,\big |\,x_{i}, v_{i}, D_{k',i}
; i\in [1,  m_{k'}]
\)\nn\\
&&=\Ebm^w \(\prod_{i=1}^{m_{k'}}F_{i} \,\big |\,x_{i}, v_{i}, D_{k',i}
; i\in [1,  m_{k'}-1]\)\Ebm^{ x_{m_{k'}}}
\(F_{ m_{k'}} \,\big |\, v_{ m_{k'}},  D_{k',m_{k'}}
\)\nn\\
&&\leq \left(1+28.2\( (k'-k+1)\frac{h_{k'-1}}{h_{k-1}}\)\right)
 \Ebm^w \(\prod_{i=1}^{m_{k'}-1}F_{i} \,\big |\,x_{i}, v_{i},
D_{k',i}
; i\in [1, m_{k'}-1]\)\Ebm^{\la_{k'-1}} \( F_{m_{k'}}\,\big |\, D_{k',m_{k'}} \)
\nn
\eea
and by induction
\bea
&&\Ebm^w \(\prod_{i=1}^{m_{k'}}F_{i} \,\big |\,\trav{n}{y}{k-1}{k}{y}{k'}{k-1} = m_{k'};\,\, x_{j}, v_{j},  ; j\in [1, m_{k'}]\)\label{indfa.4}\\ 
&&\leq \left(1+28.2\( (k'-k+1)\frac{h_{k'-1}}{h_{k-1}}\)\right)^{ m_{k'}}
\prod_{i=1}^{ m_{k'}}\Ebm^{\la_{k'-1}} \( F_{i} \big |\,D_{k', i}\).\nn
\eea
The same proof gives also that the LHS of \eqref{indfa.4} is bounded below by the right hand side, with $+28.2$ replaced by $-28.2$. Integrating the later  with respect to
$\Pbm^w(\cdot|\trav{n}{y}{k-1}{k}{y}{k'}{k-1} = m_{k'})$ gives
\be
\Ebm^w \(\prod_{i=1}^{m_{k'}}F_{i} \,\big |\,\trav{n}{y}{k-1}{k}{y}{k'}{k-1} = m_{k'}\)\label{indfa.41}
\geq \left(1-28.2\( (k'-k+1)\frac{h_{k'-1}}{h_{k-1}}\)\right)^{ m_{k'}}
\prod_{i=1}^{ m_{k'}}\Ebm^{\la_{k'-1}} \( F_{i} \big |\,D_{k', i}\).
\ee
Together with (\ref{indfa.4}) this gives 
 the upper bound in  (\ref{indfa.2}) and the lower bound is proven similarly.\qed

\bibliographystyle{plain}

\newpage
\subsection*{Index of Notation}
The following notation is used frequently; the second column points to its definition or first systematic use.

\begingroup
\small
\setlength{\tabcolsep}{0pt}
\renewcommand{\arraystretch}{1.10}
\setlength{\emergencystretch}{2em}
\noindent\begin{tabular}{@{}p{0.315\textwidth}@{\hspace{0.025\textwidth}}p{0.62\textwidth}@{}}
  $\CC_{\sfx,\vep,\bM},\ \widetilde{\CC}^{\star}_{\vep,\bM}$ $\CC^{\star}_{\vep,\bM},\ \sfA_{\bM}$ &
  point-cover time, cover time, normalized cover time and Riemannian area, \eqref{eq:dprz}\\[2pt]
  $\CCT,\ \osfm_\vep,\ \sfrho_\vep\ $ &
  torus cover time, centering and slope, \eqref{def:mn}\\[2pt]
  $W_\star,\ X_\star,\ Y_\star$ &
  standard Gumbel variable, random shift and limiting cover-time variable, \eqref{eq-Gumble}, Theorem \ref{theo-covertime}\\[2pt]
  $\dlc,\ \sfMo$ &
  city-boundary width and number of extended cities, beginning of Section \ref{subsec-cities}\\[2pt]
  $\bA_u^\sfe,\ \bA_u,\ \bA_u^\sfo,\ \bTs^\sfe,\ \bTs$ &
  extended city, city, inner city and the two boundary unions, Section \ref{subsec-cities}\\[2pt]
  $\sfz_{01},\ \sfz_{02},\ \sfz_{j1},\ \sfz_{j2},\ T_j$ &
  inter-city hitting points and excursion duration, Section \ref{subsec-cities}\\[2pt]
  $\HH_k,\ \HHu,\ \HH_\infty^\star,\ \tau_u(k)$, $\sfx_i^{(u)},\ \FF_{u,s}$ &
  inter-city filtration and city-excursion data, Section \ref{subsec-cities}, \eqref{eq:FF-to-HH}\\[2pt]
  $\bb,\ \sfv_\bb,\ \bB_\sfd(\sfv,r)$, $\sfd,\ \npp$ &
  blocks, block centers, metric balls, torus distance and the origin, Section \ref{subsec-cities}\\[2pt]
  $\CCd,\ \srs,\ \sfM$ &
  modified cover time, retained-city radius and number of retained cities, Definition \ref{dfn:CCd}, \eqref{ct2.25c}\\[2pt]
  $\varlimsup_{\dlc,\sfKfv,\sfKz,\vep}$, $\varliminf_{\dlc,\sfKfv,\sfKz,\vep}$ &
  iterated-limit convention, Section \ref{subsec-cities} and the notation paragraph preceding Section \ref{sec:outline}\\[2pt]
  $\mu_\sfe,\ S_k,\ \sft_1,\ \sfv_\sfe$ &
  invariant law, inter-city clock, mean and asymptotic variance, \eqref{dfn:Sk}, Lemma \ref{lem:conc}\\[2pt]
  $\UUk,\ \CCex$ &
  uncovered-block event and cover time measured in inter-city excursions, \eqref{cov-time-exc}\\[2pt]
  $r_0,\ \ello,\ \sfm_\vep$ &
  reference radius, logarithmic offset and shifted centering, \eqref{def:mn-alt}\\[2pt]
  $\sfk_\vep(z),\ \bar\sfs_\vep(z)$ &
  inter-city excursion count and its nominal clock, \eqref{dfn:bar-svep}\\[2pt]
  $\UUz,\ \FFuz$, $\Lambda^{(\vep)}_{u,z},\ Q^{(\vep)}_{u,z}$ &
  uncovered-block event, city conditioning, count and conditional probability, \eqref{dfn:UUz}, \eqref{def:barQu}\\[2pt]
  $n_\vep,\ \sfae,\ n_\vep^{(i)},\ \wnvfv$ &
  radial depth, fractional correction and shifted depths, \eqref{def.nvep}, \eqref{def.n-vep-i}\\[2pt]
  $\sfKz,\ \sfKo$, $\sfKt,\ \sfKth$ &
  city, radial-tower, coupling and barrier-gap levels, Section \ref{subsec-circle}\\[2pt]
  $\sfKf,\ \sfKfv$, $\sfKp,\ \sfKw,\ \sfKtp$ &
  terminal, block, rough-barrier, angle-control and bridge levels, Section \ref{subsec-circle}, Remark \ref{rem:more-int-effects}\\[2pt]
  $r_j,\ T_j^{\sfv,\sfKo,s}$ &
  concentric radii and traversal counts, \eqref{dr.1}, \eqref{Tjsdef}\\[2pt]
  $\wt c,\ \gak,\ \kapo,\ \kapp,\ \okappa,\ \cst$ &
  net, lower-barrier, window and angle-control parameters, \eqref{dfn:aux-par}\\[2pt]
  $\gakm,\ \kapom,\ \kappm,\ \okappam$ &
  slightly smaller exponents used in relaxed barriers, Section \ref{sec:pf-same-mean} and Section \ref{sec:2nd-mom}\\[2pt]
  $\UUp,\ \RUb,\ \BBb,\ \WWb$ &
  truncated uncovering event and its restricted-uncovering, barrier and angle components, \eqref{eq-good}\\[2pt]
  $\Lambda^{(+,\vep)}_{u,z},\ Q^{(+,\vep)}_{u,z}$ &
  truncated city count and conditional probability, \eqref{def:barQub}\\[2pt]
  $\bF^\star,\ \bF^\bb,\ \bF_{n_\vep},\ \bF_j$ &
  reference, block and multiscale nets, \eqref{dfn:F-star}, \eqref{dfn:F-beta}\\[2pt]
  $\gamma_b(j),\ j_\ell$ &
  curved excursion-count barriers and distance from the nearer endpoint, \eqref{eq:norm-exc}\\[2pt]
 $\sfs_j(y)$ &
  nominal excursion count at level $j$, \eqref{eq:nominal}\\[2pt]
  $\Tvjz,\ \Tjz,\ \etavjz,\ \etajz$, $\whetavjz,\ \whetajz$ &
  $z$-clock counts, square-root counts and centered excess counts, Definition \ref{def:ex-cnt}, \eqref{eq:eta-centered}\\[2pt]
$\Tvjy,\ \Tjy,\ \etavjy,\ \etajy$, $\whetavjy,\ \whetajy$ &
  $y$-clock versions of the preceding quantities, Definition \ref{def:ex-cnt}\\[2pt]
  $\whUUz$ &
  net-based uncovered proxy, \eqref{dfn:UU-alt}\\[2pt]
\end{tabular}
\endgroup
\newpage
\subsection*{Index of Notation (continued)}\leavevmode\par

\begingroup
\small
\setlength{\tabcolsep}{0pt}
\renewcommand{\arraystretch}{1.10}
\setlength{\emergencystretch}{2em}
\noindent\begin{tabular}{@{}p{0.315\textwidth}@{\hspace{0.025\textwidth}}p{0.62\textwidth}@{}}
   $\sfJupl,\ \sfJo,\ \sfJm,\ \sfI_j$ &
  barrier level sets and excursion-count windows, \eqref{dfn:J-intervals}, \eqref{dfn:windows}\\[2pt]
  $\BBy,\ \BBb$ &
  $y$- and $z$-clock double-barrier events, \eqref{eq:barrier-event-z}\\[2pt]
  $\Ybz,\ \UUy(y),\ \UUv$ &
  block top-level excess and $y$-parametrized uncovering events, \eqref{eq:time-ident}\\[2pt]
  $\sfJs,\ \RUb,\ \RRv,\ \RUy,\ \RRy$ &
  terminal levels and restricted-uncovering events, \eqref{def:res-u-z}, \eqref{eq:barrier-event}\\[2pt]
  $\thbj,\ \Ljy,\ \Ljz,\ \LKoy,\ \LKoz$ &
  angle differences and their empirical measures, \eqref{dfn:emp-ang}\\[2pt]
  $\Ipi,\ \bar\mu,\ \bar c$ &
  angle space and reference Poisson-kernel law, \eqref{dfn:Poiss}\\[2pt]
  $\sfd_\infty,\ \bB_\mu^\infty(c)$ &
  sup-quantile distance and its balls, \eqref{dfn:dsup}\\[2pt]
  $\boetaj,\ \boetasj,\ \sfJwp,\ \sfJw$ &
  angle tolerances and controlled levels, \eqref{dfn:angle-width}, \eqref{eq:was-event-z}\\[2pt]
  $\WWy,\ \WWb,\ \UUpy(y)$ &
  angle-control events and the complete $y$-clock truncation, \eqref{eq:was-event}, \eqref{eq:was-event-z}, \eqref{eq:good-cover-y}\\[2pt]
  $\wt Q^{(1,\vep)}_{u,z},\ \wt Q^{(2,\vep)}_{u,z}$ &
  conditional first and second factorial moments, \eqref{dfn:wt-Q}, \eqref{dfn:wt-Q2}\\[2pt]
  $\sfIs_i,\ \ggi(y),\ \bar\sfs_{\sfK_i}(y)$ &
  integer-compatible windows, first-moment weights and terminal clocks, \eqref{dfn:ggi}, \eqref{eq:nominalno-vep}\\[2pt]
  $\uu{\sfx},\ \BWlb$ &
  ordered endpoint data and the large-scale good event, \eqref{eq:BW-K12}\\[2pt]
  $\uu\alpha(\sfKp),\ \alpha(\sfKfv,\sfKp),\ \bar\alpha(\sfKp)$ &
  within-block contribution and its uniform bounds, Proposition \ref{prop:angle-diff}\\[2pt]
  $\Istar,\ \Vb(\nu)$ &
  admissible top-level data and angular-defect functional, \eqref{dfn:I-star-ext}, \eqref{dfn:hat-L-alpha}\\[2pt]
  $\gstb(y,v),\ \sfb$ &
  first-moment kernel and Gaussian correction, \eqref{eq:g-def}\\[2pt]
  $\GG^\tau,\ \GG_j^{\sfv,\tau},\ \GG_j^\sfv$ &
  radial excursion sigma-algebras, Definition \ref{def:Fvjs}\\[2pt]
  $\sfIw,\ \bKz,\ \Vbz$ &
  angular-defect range, normalization and block defect, \eqref{eq:Ball-to-angle}, \eqref{def:bar-bKz}\\[2pt]
  $\gst(y,v),\ \sfIos(y,v),\ \sfX^{(\vep)}_{u,z},\ \bar Q^{(1,\vep)}_{u,z}$ &
  compactly supported kernel, its cutoff, field proxy and first-moment proxy, \eqref{eq:g-def-alt}, \eqref{def:X-eps}, \eqref{def:X-eps-pre}\\[2pt]
  $\wX^{(\sfwm)}_{\bb,z},\ \wX^{(\vep,\Cap)}_{\bb,z}$ &
  angle and window error contributions, \eqref{eq:bd-wX}, Section \ref{sec:good-cover-end}\\[2pt]
  $\sfX_u,\ \overline X,\ X_\star$ &
  limiting city field, aggregate variable and final random shift, Proposition \ref{prop:Gauss-field-continuity}, \eqref{dfn:Xstar-K}, Theorem \ref{theo-covertime}\\[2pt]
  $\dst$ &
  second-moment error parameter, \eqref{eq:C-star-good}\\[2pt]
  $\gamma_b^J(j),\ \sfrho_\vep^{J,b}$ &
  barriers bent relative to an interval $J$, \eqref{dfn:gamma-J}\\[2pt]
  $\NN_{\sfv,r,\ell'},\ \NN_{\sfv,r,m,k}$, $\PRzu,\ \PRsu,\ (\bP_R^1)^m$ &
  off-centered excursion-count sigma-algebras and laws, beginning of Section \ref{subsec:off-center}, \eqref{eq:off-center-bd}\\[2pt]
  $\bar\bP_z^{(\star)},\ \eta_{\sfv,j}$, $\eta^\sharp_{\sfv,\sfKf},\ \sfLw$ &
  within-block proxy law, counts and terminal depth, Section \ref{sec:g4-bds}, \eqref{def:why}, \eqref{def:eta-sharp}\\[2pt]
 $\gfo(z),\ \QQs$ &
  within-block contribution and restricted uncovering event, \eqref{dfn:g4-alt}, \eqref{dfn:QQs}\\[2pt]
  $\Delta,\ \tsfKfv,\ \sfLo,\ \sfLn,\ \sfJus$ &
  terminal re-scaling and level interval, \eqref{dfn:tsfKfv}\\[2pt]
$\gamma_{-1},\ \sfIp_j,\ \gamma_{-1}^{\sfJus}$ &
  terminal barriers and windows, \eqref{dfn:gam}, \eqref{dfn:Ip-sfJs}, \eqref{dfn:gam-Js}\\[2pt]
  $\wQs_{\sfv},\ \whcI_{\sfv},\ \whcI_{\sfv,\ell}$ &
  terminal barrier events and their relaxations, \eqref{def:QQs}, \eqref{def:whcI}, \eqref{def:Aell}\\[2pt]
$q^{(u)}_{\ell',\ell},\ q^{(\ge)}_{\sfJus},\ q^\Cap_{\sfJus}$ &
  terminal transition kernels, Definition \ref{dfn:Markov}\\[2pt]
  $\GJ_i',\ S_m^{(j)},\ \Lambda_q,\ \widehat S_m^{(j)}$ &
  offspring variables, sums, cumulant generating function and normalization, \eqref{dfn:GJp}--\eqref{eq:GJp-var}\\[2pt]
  $\sfd_\sfw,\ \sfd_{\rm KS}$ &
  Wasserstein and Kolmogorov--Smirnov distances, \eqref{dfn:dw-orig}, \eqref{eq:quantiles}\\[2pt]
  $L_m^\mu,\ L_\mu(\sfs)$ &
  discrete quantile and empirical measures, \eqref{eq:ds-emp}, \eqref{lbl:d-inf-bd}\\[2pt]
$\zeta_\sfz,\ \sfd_{\rm TV}$ &
  conditional hitting law and total-variation distance, Lemma \ref{lem-coupb}\\[2pt]
  $L^{(k)},\ q_i^{(k)},\ \sfv_k$ &
  conditional-sum marginals and variances, \eqref{eq:mean-var}--\eqref{eq:qi-spread}\\[2pt]
\end{tabular}
\endgroup
\newpage
\subsection*{Index of Notation (continued)}\leavevmode\par

\begingroup
\small
\setlength{\tabcolsep}{0pt}
\renewcommand{\arraystretch}{1.10}
\setlength{\emergencystretch}{2em}
\noindent\begin{tabular}{@{}p{0.405\textwidth}@{\hspace{0.025\textwidth}}p{0.56\textwidth}@{}}
  $S_\ell^{(k)},\ \bQ_{m,b}^{(k)}$ &
  partial sums and laws conditioned on their total, \eqref{dfn:Sl}\\[2pt]
  $\Theta_\ka,\ \Gamma$ &
  admissible coupling parameters and comparison set, \eqref{dfn:Thetak}, \eqref{dfn:Gamma}\\[2pt]
  $\bar\Delta_{(m,b),(\widehat m,\widehat b)}$, $N_i^{(m,b)}$ &
  coupling discrepancy and count profile, \eqref{dfn:bar-Del}, \eqref{dfn:Nmb}\\[2pt]
  $\tau_h,\ \bQ_{m'}^\star,\ \AA_{m'}^\star$ &
  stopping index, stopped law and exact-sum event, \eqref{def:tau-h}, \eqref{eq:Amb-star-def}\\[2pt]
  $\trav{m}{\sfx}{\widetilde R}{\widetilde r}{\sfv}{r}{R}$ &
  Brownian traversal count, Definition \ref{def-trav}\\[2pt]
  $T_{\sfv,j}^{\sfx,m},\ T_j^{\sfx,k,m}$, $T_{\sfv,j}^{[\tau,\tau']}$ &
  standard traversal abbreviations, \eqref{eq-calgary}\\[2pt]
  $\wt\sfs_{n_\vep}(y)$, $\tau_{\sfv,r,r'}(m),\ \tau_\sfv(m)$ &
  approximate nominal count and excursion clocks, \eqref{clt.2}, \eqref{eq-tauab}, \eqref{eq-taux}\\[2pt]
  $T_{\sfv,j}^{\partial\bC^\sfe(k)}$ &
  traversals accumulated during inter-city excursions, \eqref{dfn:T-Cek}\\[2pt]
  $\BWz,\ \BWs$ &
  enlarged barrier--angle good events, \eqref{dfn:bw}, \eqref{dfn:bws}\\[2pt]
  $\Rkz,\ \Rkbp,\ \Rkvpp$ &
  rough lower-barrier events, \eqref{def:Jay-Thm4.2}, \eqref{a216.0}\\[2pt]
  $\Akz,\ \Akbz$ &
  clock-comparison events, \eqref{def:Jay-Thm4.1}\\[2pt]
  $\whUUvk,\ \AUvz$ &
  effective-cover and exceptional events, \eqref{def:UUvk}, \eqref{def:AUvz}\\[2pt]
  $\whUUpk,\ \whUUp,\ \BWp$ &
  block effective-cover proxies and their good event, \eqref{def:UUpk}, \eqref{dfn:bw-Jay}\\[2pt]
  $\sfJup,\ \sfIp_j$ &
  narrower effective-cover levels and windows, \eqref{dfn:windows-alt}\\[2pt]
  $\bar\bP_y^{(\ell)}$ &
  concentric excursion product law, Definition \ref{def:Fvj}\\[2pt]
  $\CCe_\sfv(v,v'),\ \QQe(v')$, $\alpha_{\sfKfv,\sfKp}(v')$ &
  inner-block decomposition and terminal contribution, \eqref{eq:decom-U}, \eqref{def:alpha56}\\[2pt]
  $\BWb(y_4),\ \BBKth(y_4,y_5), \ \wCe(y_2,y_3)$ &
  radial pieces of the good event, \eqref{eq:B-K1-K2}--\eqref{eq:hat-C}\\[2pt]
  $\gth,\ \gfo,\ \gKth,\ \gfl,\ \whgth$ &
  basic barrier weights and radial first-moment functions, Section \ref{sec:angle-diff}, \eqref{eq:alt-g3}, \eqref{def:wh-g3-ell}, \eqref{def:wh-g3}\\[2pt]
  $\alpha(\sfKfv,\sfKp)$ &
  integrated within-block constant, \eqref{eq:alpha-tot}\\[2pt]
$\btet,\ L_{\btet},\ V_{\btet}$ &
  top-level angle data and angular defect, \eqref{dfn:Vthet}\\[2pt]
  $\gKo,\ \whgo,\ \Bl$ &
  conditional and lower-barrier first-moment proxies, \eqref{eq:K1-K3-alt}--\eqref{eq:BK12}\\[2pt]
  $\bmu,\ \bar p_\vep,\ \zeta_{\sfz_q},\ B_q$ &
  quantile coupling and Bernoulli replacement variables, \eqref{eq:quant}--\eqref{dfn:coup}\\[2pt]
  $\whes(j),\ \Delj,\ \wtgo$ &
  retained-count process, coupling defect and replacement proxy, \eqref{dfn:wh-g1n}\\[2pt]
 $\htq,\ \bgB(\uu B),\ \bDl(\uu B)$ &
  log-Laplace correction, conditional first-moment factor and centered coupling correction, \eqref{dfn:hq-alt}, \eqref{dfn:bar-g}, \eqref{dfn:bar-Dlt}\\[2pt]
  $\sfKop,\ \Pon_\bb,\ \Pm_{\bb,\sfKop}$ &
  auxiliary top scale and single-block probabilities, \eqref{dfn:P1-}\\[2pt]
  $\UUmy(y),\ \UUmp(y)$ &
  relaxed good events at the two top levels, \eqref{dfn:Umy}\\[2pt]
  $\Ptw_{\bb,\bb'},\ \Ptw_{\bb'}$ &
    pair probabilities and their block sums, \eqref{dfn:P2+}\\[2pt]
 $\bG_\ell,\ \Ptwl_\bb$ &
  distance classes and terminal-pair contribution, \eqref{eq-G0L}, \eqref{bd:Late-brch}\\[2pt]
  $\bar\mu_\sfe$ &
  invariant law of inter-city endpoint pairs, beginning of Section \ref{sec:end-game}\\[2pt]
  $N_\sfv,\ \Yvz$ &
  top-level excursion count and centered square-root field, \eqref{17.01ad}\\[2pt]
  $\widehat N^\th_\sfv,\ \Vvz$ &
  angular-defect count and field, \eqref{def:V-eps-alt}\\[2pt]
 $\bar\bF_j,\ \wh\bF_\vep,\ \bar \bF_{\bar n_\vep},\ \bar n_\vep$ &
  epsilon-independent nets and continuity mesh, Section \ref{sec:end-game}\\[2pt]
  $q_{\rm in}^{(\vep)},\ q_{\rm out}^{(\vep)}$ &
  inner and outer fluctuation vectors, \eqref{def:qin}--\eqref{def:qout-th}\\[2pt]
  $\Sigma_1,\ \Sigma_2$ &
  limiting covariance matrices, \eqref{mCLT.03}, \eqref{mCLT.04}\\[2pt]
\end{tabular}
\endgroup
\newpage
\subsection*{Index of Notation (continued)}\leavevmode\par

\begingroup
\small
\setlength{\tabcolsep}{0pt}
\renewcommand{\arraystretch}{1.10}
\setlength{\emergencystretch}{2em}
\noindent\begin{tabular}{@{}p{0.315\textwidth}@{\hspace{0.025\textwidth}}p{0.62\textwidth}@{}}
  $\sfX^{(\vep,p)}_{u,z},\ g_i^{(p)},\ \sfX_u^{(p)}$ &
  finite-dimensional field approximants, \eqref{def:X-eps-pm}--\eqref{def:210.2-alt}\\[2pt]
  $\sfd_{\rm BL}$ &
  bounded-Lipschitz metric for weak convergence, \eqref{def:BL-conv}\\[2pt]
  $N_\sfv^\th$ &
  cosine sum used in the continuity estimates, \eqref{def:Nx-th}\\[2pt]
  $\AA_{m,a,b}^{(k)},\ \bP_{m,a,b}$, $\bP_{m'}^\star$ &
  conditioned and free collections of inner excursions, \eqref{eq:Amb-def}\\[2pt]
  $\Delta_{m,a,b}$ &
  number of unmatched excursions in the coupling, Proposition \ref{prop:coupling}\\[2pt]
  $\al(\ell),\ t_z,\ I_y$ &
  appendix barrier-estimate notation, beginning of \Cref{sec:BoundaryCrossing}\\[2pt]
\end{tabular}
\endgroup

\end{document}